\documentclass[a4paper,12pt,amsart,frenchb]{article}

\usepackage{amsmath,amsbsy,amsfonts,amssymb}
\usepackage[french]{babel}
\newcommand{\ass}[2]{\vskip0.3cm\noindent
{\bf {#1}}. { \sl {#2}}\vskip0.3cm\noindent
}
  
\begin{document}

    \title{Stabilisation de la formule des traces tordue VI:   la partie g\'eom\'etrique de cette formule}
\author{C. Moeglin, J.-L. Waldspurger}
\date{9 juin 2014}
\maketitle

{\bf Introduction}

\bigskip
Le but  ultime de notre travail est la stabilisation de la formule des traces tordue. Le pr\'esent article \'enonce les r\'esultats que l'on a en vue concernant la partie g\'eom\'etrique de cette formule. Ce sont les g\'en\'eralisations au cas tordu des th\'eor\`emes    \'enonc\'es par Arthur dans son premier article sur la stabilisation ([A1]), du moins de ceux qui concernent cette partie g\'eom\'etrique. Nous ne d\'emontrons pas ici les r\'esultats en question. Ils seront d\'emontr\'es plus tard en reprenant les m\'ethodes des deux autres articles d'Arthur sur le sujet.  On doit  dire que  g\'en\'eraliser au cas tordu les constructions d'Arthur pose certains probl\`emes techniques, mais aucun probl\`eme conceptuel. C'est-\`a-dire que l'essentiel est d\^u \`a Arthur lui-m\^eme.

La premi\`ere section pr\'esente le cadre global dans lequel on se place. On consid\`ere un corps de nombres $F$, un groupe r\'eductif connexe $G$ d\'efini sur $F$, un espace tordu $\tilde{G}$ sous $G$ et un caract\`ere $\omega$ de $G({\mathbb A})$ trivial sur $G(F)$, o\`u ${\mathbb A}$ est l'anneau des ad\`eles de $F$.   On  d\'efinit les int\'egrales orbitales pond\'er\'ees globales, certaines de leurs variantes et on \'enonce les formules de descente qui les relient \`a leurs avatars locaux. La section 2 \'enonce la partie g\'eom\'etrique de la formule des traces tordue $\omega$-\'equivariante. On \'enonce cette formule d'une fa\c{c}on un peu plus abstraite qu'Arthur. Il est traditionnel et naturel  de l'\'ecrire comme une somme avec coefficients d'int\'egrales orbitales $I^{\tilde{G}}(\gamma,\omega,f)$ ou plus g\'en\'eralement d'int\'egrales orbitales pond\'er\'ees $\omega$-\'equivariantes $I_{\tilde{M}}^{\tilde{G}}(\gamma,\omega,f)$. La pr\'esence du caract\`ere $\omega$ perturbe d\'ej\`a la situation: les avatars locaux des int\'egrales $I^{\tilde{G}}(\gamma,\omega,f)$ ne d\'ependent  pas seulement de la classe de conjugaison de $\gamma$ mais bien du point base $\gamma$ choisi. Surtout, comme on le sait, les coefficients dont sont affect\'es ces int\'egrales ne sont \`a ce jour pas connus explicitement (ils sont connus si $\gamma$ est fortement r\'egulier, mais pas si $\gamma$ contient une partie unipotente). On a choisi de regrouper les int\'egrales $I^{\tilde{G}}(\gamma,\omega,f)$, affect\'ees de leurs coefficients, selon la classe de conjugaison \`a laquelle appartient la partie semi-simple de $\gamma$. On obtient ainsi des distributions  not\'ees $A^{\tilde{G}}(V,{\cal O},\omega)$ d\'ependant d'un ensemble fini assez grand $V$ de places du corps de base $F$  et d'une classe de conjugaison semi-simple ${\cal O}$ dans $\tilde{G}(F)$. On \'enonce en 2.3 leur d\'efinition. Ces distributions  sont les ingr\'edients globaux de la partie g\'eom\'etrique de la formule des traces.    Elles v\'erifient une formule de descente qui les ram\`ene au cas basique o\`u $\tilde{G}=G$ n'est pas tordu et o\`u ${\cal O}$ est simplement la classe $\{1\}$ (dans ce cas, $A^{\tilde{G}}(V,{\cal O},\omega)$ est exactement la "partie unipotente" de la formule des traces). Dans la section 3, on pr\'esente la th\'eorie de l'endoscopie dans le cadre global. La diff\'erence essentielle avec le cas local d\'evelopp\'e en [I] est que, dans le cas global, pour une donn\'ee endoscopique ${\bf G}'$ relevante, il y a un facteur de transfert canonique. La situation tordue pose ici un probl\`eme technique. Usuellement, la d\'efinition d'un tel facteur utilise un point $\delta\in \tilde{G}'(F)$, assez r\'egulier, qui se transf\`ere en un \'el\'ement $\gamma_{v}\in \tilde{G}(F_{v})$ pour toute place $v$ de $F$. Dans le cas non tordu, un tel point existe si ${\bf G}'$ est relevante. Ce n'est plus vrai dans le cas tordu. On dit que ${\bf G}'$ est relevante si et seulement $\tilde{G}'(F)\not=\emptyset$ et la donn\'ee locale ${\bf G}'_{v}$ est relevante pour toute place $v$. On a pos\'e cette d\'efinition parce que c'est la seule notion que l'on sache contr\^oler. Or cela n'assure pas l'existence d'un $\delta\in \tilde{G}'(F)$ v\'erifiant les propri\'et\'es ci-dessus. Il est possible que les donn\'ees ${\bf G}'$ pour lesquelles il n'existe pas de tel $\delta$ puissent \^etre \'elimin\'ees du processus de stabilisation mais cela ne nous para\^{\i}t pas clair. On a plut\^ot  choisi de donner une d\'efinition du facteur de transfert global sous la seule hypoth\`ese de relevance telle que d\'efinie ci-dessus. D'abord, en inspectant les d\'efinitions de Kottwitz-Shelstad ou de Labesse, on voit que l'on n'a pas vraiment besoin d'un $\delta$ comme ci-dessus. Il suffit qu'il existe un sous-tore tordu maximal $\tilde{T}'$ de $\tilde{G}'$, d\'efini sur $F$, de sorte que, pour toute place $v$, il existe un \'el\'ement $\delta_{v}\in \tilde{T}'(F_{v})$ assez r\'egulier qui se transf\`ere en un \'el\'ement $\gamma_{v}\in \tilde{G}(F_{v})$. M\^eme cette propri\'et\'e moins forte n'est pas assur\'ee par notre hypoth\`ese de relevance. Mais on peut plonger $\tilde{G}$ et $\tilde{G}'$ dans des espaces plus gros qui satisfont cette propri\'et\'e. On d\'efinit alors le facteur de transfert comme la restriction \`a notre couple $(\tilde{G}',\tilde{G})$ du facteur de transfert d\'efini sur ces espaces plus gros. Cela est fait au paragraphe 3.9. Dans la section 4, on d\'efinit les avatars stables et endoscopiques des int\'egrales orbitales pond\'er\'ees invariantes. On \'enonce le r\'esultat principal en 4.5, \`a savoir qu'une int\'egrale orbitale pond\'er\'ee endoscopique est en fait une int\'egrale orbitale pond\'er\'ee ($\omega$-\'equivariante) tout court. Ce r\'esultat se d\'eduit des analogues locaux \'enonc\'es en [II] et [V], qui restent \`a d\'emontrer. La section 5 \'enonce la version stable de la partie g\'eom\'etrique de la formule des traces. Les distributions $A^{\tilde{G}}(V,{\cal O},\omega)$ sont remplac\'ees par des distributions $SA^{\tilde{G}}(V,{\cal O})$, o\`u cette fois, ${\cal O}$ est une classe de conjugaison stable dans $\tilde{G}(F)$. Ces distributions doivent \^etre stables. Le th\'eor\`eme 5.4 exprime que les coefficients initiaux $A^{\tilde{G}}(V,{\cal O},\omega)$ se r\'ecup\`erent comme somme de transferts de tels coefficients $SA^{{\bf G}'}(V,{\cal O}^{{\bf G}'})$, o\`u ${\bf G}'$ d\'ecrit les donn\'ees endoscopiques de $(G,\tilde{G},{\bf a})$ qui sont elliptiques, relevantes et non ramifi\'ees hors de $V$. Le th\'eor\`eme principal 5.10 exprime que la partie g\'eom\'etrique de la formule des traces $\omega$-\'equivariante se r\'ecup\`ere de m\^eme comme somme de parties g\'eom\'etriques de formules des traces stables associ\'ees \`a ces m\^emes donn\'ees ${\bf G}'$. Dans la section 6, on montre que ce  dernier th\'eor\`eme r\'esulte des autres. Il s'agit ici d'une reprise du paragraphe 10 de l'article [A1]. Il y a deux ingr\'edients. D'une part, la proposition combinatoire 6.5 qui fournit deux expressions a priori tr\`es diff\'erentes d'une double somme sur des donn\'ees endoscopiques ${\bf G}'$ et sur des "Levi" ${\bf M}'$ de ${\bf G}'$. D'autre part, une proposition d'annulation 6.6 qui permet dans le paragraphe suivant de faire dispara\^{\i}tre les termes apparaissant dans les formules de traces stables des donn\'ees endoscopiques qui ne correspondent \`a rien du c\^ot\'e de l'espace initial $\tilde{G}$. La d\'emonstration de cette proposition 6.6 est un amusant exercice bas\'e sur les propri\'et\'es des facteurs de transfert globaux. 

\bigskip

\section{Les d\'efinitions}

\bigskip

\subsection{Groupes et espaces tordus}
Soit $F$ un corps de nombres. On note $Val(F)$ l'ensemble de ses places, $Val_{\infty}(F)$ le sous-ensemble des places archim\'ediennes, $Val_{f}(F)$ celui des places finies et  ${\mathbb A}_{F}$ son anneau d'ad\`eles. On fixe une cl\^oture alg\'ebrique $\bar{F}$ de $F$. Il est commode de supposer tacitement que, pour toute place $v\in Val(F)$, on a choisi un prolongement $\bar{v}$ de $v$ \`a $\bar{F}$. Ainsi, on peut d\'efinir le sous-groupe de d\'ecomposition $\Gamma_{F_{v}}\subset \Gamma_{F}$ comme le fixateur de $\bar{v}$. De m\^eme, pour une vari\'et\'e $X$ d\'efinie sur $F$, on a une application $X(\bar{F})\to X(\bar{F}_{v})$ obtenue en identifiant $\bar{F}$ \`a un sous-corps de $\bar{F}_{v}$ gr\^ace \`a la place $\bar{v}$.

Soient $G$ un groupe r\'eductif connexe d\'efini sur $F$ et $\tilde{G}$ un espace tordu sous $G$. On utilise les d\'efinitions des quatre premiers paragraphes de [I]. Soit ${\bf a}$ un \'el\'ement de $H^1(W_{F};Z(\hat{G}))/ker^1(W_{F};Z(\hat{G}))$, o\`u $ker^1(W_{F};Z(\hat{G}))$ est le noyau (fini) de l'homomorphisme de localisation
$$H^1(W_{F};Z(\hat{G}))\to \prod_{v\in Val(F)}H^1(W_{F_{v}};Z(\hat{G})).$$
Cet \'el\'ement ${\bf a}$ d\'etermine un caract\`ere $\omega$ de $G({\mathbb A})$, trivial sur $G(F)$. L'application ${\bf a}\mapsto \omega$ est bijective. On impose les hypoth\`eses analogues \`a celles de [I] 1.5:

$\bullet$ $\tilde{G}(F)\not=\emptyset$;

$\bullet$ $\theta^*$ est d'ordre fini, o\`u ici $\theta^*$ est la restriction de $ad_{\gamma}$ \`a $Z(G)$ pour n'importe quel \'el\'ement $\gamma\in \tilde{G}$.

On impose de plus

$\bullet$ $\omega$ est unitaire.

On pourrait ajouter la condition

$\bullet$ $\omega$ est trivial sur $Z(G;{\mathbb A}_{F})^{\theta}$,

\noindent faute de laquelle la th\'eorie devient vide. Mais, pour des raisons de r\'ecurrence, il vaut mieux ne pas l'imposer d\`es le d\'epart.

Pour prouver une assertion concernant $(G,\tilde{G},{\bf a})$, on raisonne par r\'ecurrence sur $dim(G_{SC})$, o\`u $G_{SC}$ est le rev\^etement simplement connexe du groupe d\'eriv\'e de $G$, cf. [II] 1.1. Si $(G,\tilde{G},{\bf a})$ est quasi-d\'eploy\'e et \`a torsion int\'erieure, on suppose connues toutes les assertions concernant des triplets $(G',\tilde{G}',{\bf a}')$ quasi-d\'eploy\'es et \`a torsion int\'erieure tels que $dim(G'_{SC})< dim(G_{SC})$. Si $(G,\tilde{G},{\bf a})$ n'est pas  quasi-d\'eploy\'e et \`a torsion int\'erieure, on suppose connues toutes les assertions concernant des triplets $(G',\tilde{G}',{\bf a}')$ quasi-d\'eploy\'es et \`a torsion int\'erieure tels que $dim(G'_{SC})\leq dim(G_{SC})$ et on suppose connues toutes les assertions concernant des triplets $(G',\tilde{G}',{\bf a}')$  quelconques tels que $dim(G'_{SC})< dim(G_{SC})$.  Quand une assertion est relative \`a un espace de Levi $\tilde{M}$ de $\tilde{G}$, on suppose connues toutes les assertions concernant le triplet $(G,\tilde{G},{\bf a})$ relatives \`a un espace de Levi $\tilde{L}$ contenant strictement $\tilde{M}$. Au paragraphe 1.16, nous introduirons la notion de $K$-triplet $(KG,K\tilde{G},{\bf a})$. Quand on travaille avec un tel triplet, on suppose connues toutes les assertions concernant des triplets $(G',\tilde{G}',{\bf a}')$ quasi-d\'eploy\'es et \`a torsion int\'erieure tels que $dim(G'_{SC})\leq dim(G_{SC})$ et on suppose connues toutes les assertions concernant des $K$-triplets $(KG',K\tilde{G}',{\bf a}')$ tels que $dim(G'_{SC})<dim(G_{SC})$. 

Revenons \`a notre triplet $(G,\tilde{G},{\bf a})$. Pour toute place $v\in Val(F)$, on d\'eduit par localisation de $(G,\tilde{G},{\bf a})$ un triplet $(G_{v},\tilde{G}_{v},{\bf a}_{v})$ sur $F_{v}$, qui v\'erifie les hypoth\`eses de [I]. On note $V_{ram}$, ou plus pr\'ecis\'ement $V_{ram}(\tilde{G},{\bf a})$, le plus petit ensemble de places $v$ contenant les places archim\'ediennes et tel que, pour $v\not\in V_{ram}$, on ait:

- $G$ et ${\bf a}$ sont non ramifi\'es en $v$ et $\tilde{G}(F_{v})$ poss\`ede un sous-espace hypersp\'ecial;

 - en notant $p$ la caract\'eristique r\'esiduelle de $F_{v}$ et $e_{v}=[F_{v}:{\mathbb Q}_{p}]$, on a $p>5$  et $p> N(G)e_{v}+1$, o\`u $N(G)$ est l'entier d\'efini en [W2] 4.3.

On fixe une paire parabolique $(P_{0},M_{0})$ de $G$ d\'efinie sur $F$ et minimale. On en d\'eduit une paire parabolique $(\tilde{P}_{0},\tilde{M}_{0})$ de $\tilde{G}$. 

 Soit $V$ un ensemble fini de places contenant $V_{ram}$ et tel que $G$ et $\tilde{G}$ puissent \^etre d\'efinis sur l'anneau $\mathfrak{o}^{V}$ des \'el\'ements de $F$ qui sont entiers hors de $V$. Fixer des structures de $G$ et $\tilde{G}$ sur $\mathfrak{o}^V$ permet de d\'efinir les ensembles $G(\mathfrak{o}_{v})$ et $\tilde{G}(\mathfrak{o}_{v})$ pour tout $v\not\in V$, o\`u $\mathfrak{o}_{v}$ est l'anneau d'entiers de $F_{v}$. Pour tout $v\not\in V_{ram}$, on fixe un sous-groupe compact hypersp\'ecial $K_{v}$ de $G(F_{v})$ et un sous-espace hypersp\'ecial $\tilde{K}_{v}$ de $\tilde{G}(F_{v})$ associ\'e \`a $K_{v}$. On impose que $K_{v}=G(\mathfrak{o}_{v})$ et $\tilde{K}_{v}=\tilde{G}(\mathfrak{o}_{v})$ pour presque tout $v$. Cette condition ne d\'epend ni du choix de $V$, ni de celui  des structures sur $\mathfrak{o}^V$. Elle implique que, pour $\gamma\in \tilde{G}(F)$, on a $\gamma\in \tilde{K}_{v}$ pour presque tout $v$. Pour $v\in V_{ram}$, on fixe un sous-groupe compact maximal $K_{v}$ de $G(F_{v})$, sp\'ecial si $v$ est non archim\'edienne. On impose
 
 - pour toute place $v\in Val(F)$, $K_{v}$ est en bonne position relativement \`a $M_{0}$.
 
 C'est possible puisque, pour $v\not\in V_{ram}$, tout sous-groupe compact hypersp\'ecial est conjugu\'e par un \'el\'ement de $G(F_{v})$ \`a un tel sous-groupe en bonne position relativement \`a $M_{0}$. Pour $v\not\in V_{ram}$, on note ${\bf 1}_{\tilde{K}_{v}}$ la fonction caract\'eristique de $\tilde{K}_{v}$ dans $\tilde{G}(F_{v})$ et on appelle mesure canonique sur $G(F_{v})$ la mesure de Haar sur ce groupe telle que $mes(K_{v})=1$ (elle est en effet canonique car tous les sous-groupes compacts hypersp\'eciaux ont m\^eme mesure).

 Pour un espace de Levi $\tilde{M}$ de $\tilde{G}$, on utilise les notations d'Arthur ${\cal P}(\tilde{M})$, ${\cal L}(\tilde{M})$ etc... Les \'el\'ements de ces ensembles sont des espaces paraboliques, resp. des espaces de Levi etc... d\'efinis sur $F$. Pour $\tilde{M}\supset \tilde{M}_{0}$ et $v\in Val(F)$, on pose $\tilde{K}_{v}^{\tilde{M}}=\tilde{K}_{v}\cap \tilde{M}(F_{v})$. Ces donn\'ees v\'erifient les m\^emes conditions que celles pour $\tilde{G}$.
 
 Rappelons que l'on note $A_{\tilde{G}}$ le plus grand sous-tore d\'eploy\'e contenu dans le sous-groupe d'invariants $Z(G)^{\theta^*}$. On note ${\cal A}_{\tilde{G}}=X_{*}(A_{\tilde{G}})\otimes_{{\mathbb Z}}{\mathbb R}$. On dispose de l'homomorphisme habituel $H_{\tilde{G}}:G({\mathbb A}_{F})\to {\cal A}_{\tilde{G}}$. On d\'efinit une application $\tilde{H}_{\tilde{G}}:\tilde{G}({\mathbb A}_{F})\to {\cal A}_{\tilde{G}}$ par les conditions suivantes:
 
 $\tilde{H}_{\tilde{G}}(\dot{\gamma})=0$ pour tout $\dot{\gamma}\in \tilde{G}(F)$;
 
 $\tilde{H}_{\tilde{G}}(x\gamma)=H_{\tilde{G}}(x)+\tilde{H}_{\tilde{G}}(\gamma)$ pour tous $x\in G({\mathbb A}_{F})$ et $\gamma\in \tilde{G}({\mathbb A}_{F})$.
 
 Soit $V$ un ensemble fini de places de $F$. On pose $F_{V}=\prod_{v\in V}F_{v}$ et on note ${\mathbb A}_{F}^V$ le sous-anneau des ad\`eles dont les composantes sur $F_{v}$ sont nulles pour tout $v\in V$.  On pose
$$C_{c}^{\infty}(\tilde{G}(F_{V}))=\otimes_{v\in V}C_{c}^{\infty}(\tilde{G}(F_{v})).$$
On d\'efinit de la m\^eme fa\c{c}on l'espace $I(\tilde{G}(F_{V}),\omega)$.

{\bf Remarque.} On adopte cette d\'efinition par commodit\'e. On pourrait aussi bien utiliser un espace un peu plus gros en regroupant les places archim\'ediennes. C'est-\`a-dire en rempla\c{c}ant la partie archim\'edienne 
$$\otimes_{v\in V\cap V_{\infty}(F)}C_{c}^{\infty}(\tilde{G}(F_{v}))$$
du produit tensoriel ci-dessus par
$$C_{c}^{\infty}(\prod_{v\in V\cap V_{\infty}(F)}\tilde{G}(F_{v})).$$
\bigskip

Consid\'erons le cas o\`u $V\supset V_{ram}$. Dans ce cas, on peut identifier $C_{c}^{\infty}(\tilde{G}(F_{V}))$ \`a un sous-espace de $C_{c}^{\infty}(\tilde{G}({\mathbb A}_{F}))$ en compl\'etant un produit $\otimes_{v\in V}f_{v}$ en le produit $ {\bf 1}_{\tilde{K}^V}\otimes (\otimes_{v\in V}f_{v})$, o\`u ${\bf 1}_{\tilde{K}^V}=\otimes_{v\not\in V}{\bf 1}_{\tilde{K}_{v}}$.  . On peut aussi identifier $Mes(G(F_{V}))$ \`a $Mes(G({\mathbb A}_{F}))$ en prolongeant toute mesure sur $G(F_{V})$ par les produit sur $v\not\in V$ des mesures canoniques sur $G(F_{v})$. Pour tout $v\in Val(F)$, notons $\tilde{G}_{v}$ l'espace $\tilde{G}$ vu comme un espace sur $F_{v}$.  En [II] 1.6, on a d\'efini l'ensemble $\tilde{{\cal A}}_{\tilde{G}_{v}}$ de la fa\c{c}on suivante. On note $G(F_{v})^1$ le noyau et ${\cal A}_{\tilde{G}_{v},F_{v}}$ l'image de l'homomorphisme $H_{\tilde{G}_{v}}:G(F_{v})\to {\cal A}_{\tilde{G}_{v}}$. Alors le sous-groupe ${\cal A}_{\tilde{G}_{v},F_{v}}$ agit naturellement sur $G(F_{v})^1\backslash \tilde{G}(F_{v})$. On pose 
 $$\tilde{{\cal A}}_{\tilde{G}_{v}}=(G(F_{v})^1\backslash \tilde{G}(F_{v}))\otimes_{{\cal A}_{\tilde{G}_{v},F_{v}}}{\cal A}_{\tilde{G}_{v}}.$$
 C'est un espace affine sous ${\cal A}_{\tilde{G}_{v}}$. On note $\tilde{H}_{\tilde{G}_{v}}:\tilde{G}(F_{v})\to \tilde{{\cal A}}_{\tilde{G}_{v}}$ l'application naturelle.   Pour $v\not\in V_{ram}$, l'image de $\tilde{K}_{v}$ par $\tilde{H}_{\tilde{G}_{v}}$ est r\'eduite \`a un point et on identifie $\tilde{{\cal A}}_{\tilde{G}_{v}}$ \`a ${\cal A}_{\tilde{G}_{v}}$ en identifiant ce point \`a $0$. Posons
 $$\tilde{{\cal A}}_{\tilde{G}_{V}}=\prod_{v\in V}\tilde{{\cal A}}_{\tilde{G}_{v}}.$$
 On d\'efinit une application
 $$\tilde{p}_{V}:\tilde{{\cal A}}_{\tilde{G}_{V}}\to {\cal A}_{\tilde{G}}$$
 de la fa\c{c}on suivante. Fixons $\dot{\gamma}\in \tilde{G}(F)$. Tout \'el\'ement $\tilde{X}\in \tilde{{\cal A}}_{\tilde{G}_{V}}$ s'\'ecrit $\tilde{X}=(\tilde{H}_{\tilde{G}_{v}}(\dot{\gamma})+X_{v})_{v\in V}$, o\`u $X_{v}\in {\cal A}_{\tilde{G}_{v}}$. On pose
 $$\tilde{p}_{V }(\tilde{X})=(\sum_{v\in V}X_{v,\tilde{G}})-(\sum_{v\not\in V}(\tilde{H}_{\tilde{G}_{v}}(\dot{\gamma}))_{\tilde{G}}),$$
 o\`u les indices $\tilde{G}$ d\'esignent les projections orthogonales sur l'espace ${\cal A}_{\tilde{G}}$.  Cette d\'efinition ne d\'epend pas du point $\dot{\gamma}$ choisi. On d\'efinit une application
 $$\tilde{H}_{\tilde{G}_{V}}:\tilde{G}(F_{V})\to {\cal A}_{\tilde{G}}$$
 par $\tilde{H}_{\tilde{G}_{V}}((\gamma_{v})_{v\in V})=\tilde{p}_{V}((\tilde{H}_{\tilde{G}_{v}}(\gamma_{v}))_{v\in V})$. 

On introduit le sous-espace  $C_{c}^{\infty}(\tilde{G}(F_{V}),K)$ de $C_{c}^{\infty}(\tilde{G}(F_{V}))$ form\'e des fonctions qui, en toute place archim\'edienne $v\in V$, sont $K_{v}$-finies \`a droite et \`a gauche. On note $I(\tilde{G}(F_{V}),K,\omega)$ son image dans $I(\tilde{G}(F_{V}),\omega)$.

 Consid\'erons le cas de deux ensembles finis de places $S\supset V\supset V_{ram}$. On pose $F_{S}^V=\prod_{v\in S-V}F_{v}$. On peut identifier $C_{c}^{\infty}(\tilde{G}(F_{V}))$ \`a un sous-espace de $C_{c}^{\infty}(\tilde{G}(F_{S}))$ en compl\'etant un produit $\otimes_{v\in V}f_{v}$ en le produit ${\bf 1}_{\tilde{K}_{S}^V}\otimes (\otimes_{v\in V}f_{v})$. De m\^eme, on   peut identifier $I(\tilde{G}(F_{V}),\omega)$ \`a un sous-espace de $I(\tilde{G}(F_{S}),\omega)$. On peut aussi identifier $Mes(G(F_{V}))$ \`a $Mes(G(F_{S}))$. 
 
 \bigskip
 
 \subsection{Remarque sur les hypoth\`eses}
 On va montrer
 
 (1) il existe un groupe alg\'ebrique non connexe $G^+$  d\'efini sur $F$ et r\'eductif, de composante neutre $G$, de sorte que $\tilde{G}$ s'identifie \`a une composante connexe de $G^+$ munie des actions de $G$ par multiplication \`a droite et \`a gauche.
 
 Preuve. Fixons $\gamma\in \tilde{G}(F)$, posons $\theta=ad_{\gamma}$. Parce que l'on suppose $\theta^*$ d'ordre fini, il existe un entier $n\geq1$ tel que $\theta^n$ soit un automorphisme int\'erieur de $G$. Il existe donc $x\in G_{SC}(\bar{F})$ tel que $\theta^n=ad_{x}$. Parce que $\theta$ est d\'efini sur $F$, $ad_{x}$ l'est aussi. Donc $\sigma(x)\in xZ(G_{SC})$ pour tout $\sigma\in \Gamma_{F}$.  Parce que $\theta$ commute \`a $\theta^n$, donc \`a $ad_{x}$, on a aussi $\theta(x)\in xZ(G_{SC})$. Notons $m$ le nombre d'\'el\'ements de $Z(G_{SC})$. Posons $N=mn$ et $y=x^m$. Alors $\theta^N=ad_{y}$ et on a $y\in G_{SC}(F)$ et $\theta(y)=y$. Notons $G^+$ l'ensemble des \'el\'ements $(g,\theta^{i})$ avec $g\in G$ et $i\in \{0,...,N-1\}$. On d\'efinit la multiplication par
 $$(g,\theta^{i})(g',\theta^j)=\left\lbrace\begin{array}{cc}(g\theta^{i}(g'),\theta^{i+j}),&\text{ si }i+j\leq N-1,\\ (g\theta^{i}(g')y,\theta^k),&\text{ si }i+j=N+k \text{ avec }k\geq0.\\ \end{array}\right.$$
 On obtient un groupe r\'eductif non connexe d\'efini sur $F$. L'application $(g,\theta)\mapsto g\gamma$ identifie la composante $G\theta$ de $G^+$ \`a $\tilde{G}$. $\square$
 
 Cette remarque ne nous servira pas directement. Mais elle nous permet d'appliquer \`a notre espace $\tilde{G}$ les r\'esultats d\'emontr\'es dans la litt\'erature pour les groupes non connexes. 
 
 \bigskip
 
 \subsection{Mesures sur les espaces ${\cal A}_{\tilde{M}}$}
  
 Soit $\tilde{M}$ un espace de Levi de $\tilde{G}$. On aura besoin d'une mesure sur l'espace ${\cal A}_{\tilde{M}}$. On a choisi d'\'eviter autant que possible de normaliser les mesures.  Logiquement, on devrait faire de m\^eme pour la mesure sur cet espace. Toutefois, cela conduirait \`a des formulations par trop inhabituelles des formules de descente. On fixe donc sur tout espace ${\cal A}_{\tilde{M}}$ une mesure de Haar, \`a laquelle on impose quelques conditions \'evidentes; par exemple, si $\tilde{M}$ et $\tilde{M}'$ sont conjugu\'es par un \'el\'ement $g\in G(F)$, on suppose que les mesures sur ${\cal A}_{\tilde{M}}$ et ${\cal A}_{\tilde{M}'}$ se correspondent par cette conjugaison. Si $\tilde{M}\subset \tilde{L}$ sont deux espaces de Levi, on munit ${\cal A}_{\tilde{M}}^{\tilde{L}}$ de la mesure pour laquelle la d\'ecomposition ${\cal A}_{\tilde{M}}={\cal A}_{\tilde{M}}^{\tilde{L}}\oplus {\cal A}_{\tilde{L}}$ est compatible aux mesures.
 
 Il y a au moins  deux fa\c{c}ons de d\'efinir ces mesures. Identifions ${\cal A}_{\tilde{M}}$ \`a $Hom(X^*(M)^{\Gamma_{F},\theta},{\mathbb R})$, o\`u $X^*(M)$ est le groupe des caract\`eres alg\'ebriques de $M$. Notons ${\cal A}_{\tilde{M},{\mathbb Z}}$ le r\'eseau 
 
 \noindent $Hom(X^*(M)^{\Gamma_{F},\theta},{\mathbb Z})$. On peut imposer que ce r\'eseau est de covolume $1$. C'est la normalisation de la th\'eorie des mesures de Tamagawa. Elle a l'inconv\'enient de se comporter assez mal vis-\`a-vis des suites exactes. Une autre m\'ethode est la suivante. Consid\'erons $\underline{la}$ paire de Borel \'epingl\'ee ${\cal E}^*=(B^*,T^*,(E^*_{\alpha})_{\alpha\in \Delta})$ de $G$, munie de l'action galoisienne quasi-d\'eploy\'ee. Posons $\bar{{\cal A}}_{T^*}=X_{*}(T^*)\otimes_{{\mathbb Z}}{\mathbb R}$. On munit cet espace vectoriel r\'eel d'une forme quadratique d\'efinie positive invariante par l'action galoisienne quasi-d\'eploy\'ee, par celle du groupe de Weyl $W$ et par l'automorphisme $\theta$ associ\'e \`a ${\cal E}^*$. C'est possible puisque le groupe d'automorphismes de $\bar{{\cal A}}_{T^*}$ engendr\'e par ces actions est fini.   Soit $\tilde{M}$ un espace de Levi de $\tilde{G}$. En choisissant un \'el\'ement $\tilde{P}\in {\cal P}(\tilde{M})$, on peut identifier ${\cal A}_{\tilde{M}}$ \`a un sous-espace de $\bar{{\cal A}}_{T^*}$. Alors ${\cal A}_{\tilde{M}}$ se retrouve muni de la restriction de la forme quadratique pr\'ec\'edente. Les invariances impos\'ees \`a cette derni\`ere impliquent que cette restriction ne d\'epend pas du choix de $\tilde{P}$. On munit ${\cal A}_{\tilde{M}}$ et plus g\'en\'eralement tout sous-espace de ${\cal A}_{\tilde{M}}$ de la mesure euclidienne associ\'ee \`a cette forme quadratique. 
 
 En tout cas, on suppose fix\'ees les mesures sur ces espaces, d'une fa\c{c}on ou d'une autre. De m\^eme, si $v$ est une place de $F$ et $\tilde{M}^v$ est un espace de Levi de $\tilde{G}_{v}$, on suppose fix\'ee une mesure de Haar sur ${\cal A}_{\tilde{M}^v}$. Notons qu'on n'impose pas de relation entre les mesures "locales" et les mesures "globales". Par exemple, si $G$ est d\'eploy\'e et si $\tilde{M}$ est un espace de Levi de $\tilde{G}$, on a l'\'egalit\'e ${\cal A}_{\tilde{M}}={\cal A}_{\tilde{M}_{v}}$, mais on ne demande pas que les mesures sur ces espaces soient les m\^emes. 
 
 Notons $G_{{\mathbb Q}}$ le groupe sur ${\mathbb Q}$ d\'eduit de $G$ par restriction des scalaires et, comme toujours, $A_{G_{{\mathbb Q}}}$ le plus grand tore d\'eploy\'e central dans $G_{{\mathbb Q}}$.  On note $\mathfrak{A}_{G}$ la composante neutre topologique de $A_{G_{{\mathbb Q}}}({\mathbb R})$. Notons que $A_{G_{{\mathbb Q}}}$ est  aussi le plus grand tore d\'eploy\'e dans le groupe sur ${\mathbb Q}$ d\'eduit de $A_{G}$ par restriction des scalaires. On en d\'eduit des inclusions
$$\mathfrak{A}_{G}\subset A_{G}(F_{\infty})\subset A_{G}({\mathbb A}_{F})\subset G({\mathbb A}_{F}),$$
o\`u $F_{\infty}=\prod_{v\in Val_{\infty}(F)}F_{v}$. L'espace $\mathfrak{A}_{G}$ est invariant par $\theta$.  On note $\mathfrak{A}_{\tilde{G}}$ le sous-espace des invariants par $\theta$.
La restriction \`a $\mathfrak{A}_{G}$ de l'homomorphisme $H_{G}:G_({\mathbb A}_{F})\to {\cal A}_{G}$ est un isomorphisme qui permet d'identifier $\mathfrak{A}_{G}$ \`a ${\cal A}_{G}$ et $\mathfrak{A}_{\tilde{G}}$ \`a ${\cal A}_{\tilde{G}}$.  On munit l'espace $\mathfrak{A}_{\tilde{G}}$ de la mesure telle que ce dernier isomorphisme pr\'eserve les mesures.  
  
  \bigskip
 
 \subsection{Formule de descente des $(\tilde{G},\tilde{M})$-familles}
 Soient $\tilde{M}$ un espace de Levi de $\tilde{G}$ et $V$ un ensemble fini non vide de places de $F$. Rappelons que, pour $v\in V$, on note  $\tilde{G}_{v}$, $\tilde{M}_{v}$ etc... les espaces $\tilde{G}$, $\tilde{M}$ etc... vus comme des espaces sur $F_{v}$. Pour tout $v$, soit $\tilde{R}^v$ un espace de Levi de $\tilde{M}_{v}$ d\'efini sur $F_{v}$. On a mis $v$ en exposant pour \'eviter une possible  confusion: $\tilde{R}^v$ n'a pas de raison d'\^etre issu d'un espace $\tilde{R}$ d\'efini sur $F$. Soit $(c_{v}(\tilde{S}^v;\Lambda))_{\tilde{S}^v\in {\cal P}(\tilde{R}^v)}$ une $(\tilde{G}_{v},\tilde{R}^v)$-famille (la variable $\Lambda$ appartient \`a $i{\cal A}_{\tilde{R}^v}^*$). On en d\'eduit une $(\tilde{G},\tilde{M})$-famille de la fa\c{c}on suivante. Pour $\tilde{P}\in {\cal P}(\tilde{M})$ et $v\in V$, on choisit $\tilde{S}^{v}\in {\cal P}(\tilde{R}^v)$ tel que $\tilde{S}^v\subset \tilde{P}_{v}$. L'espace $i{\cal A}_{\tilde{M}}^*$ se plonge dans $i{\cal A}_{\tilde{R}^v}^*$. Pour $\Lambda\in i{\cal A}_{\tilde{M}}^*$, on pose
 $$c(\tilde{P};\Lambda)=\prod_{v\in V}c_{v}(\tilde{S}^v;\Lambda).$$
 Cela ne d\'epend pas du choix des $\tilde{S}^v$ et la famille $(c(\tilde{P};\Lambda))_{\tilde{P}\in {\cal P}(\tilde{M})}$ est une $(\tilde{G},\tilde{M})$-famille. Pour $\tilde{P}\in {\cal P}(\tilde{M})$, on note $\Delta_{\tilde{P}}$ l'ensemble des restrictions \`a ${\cal A}_{\tilde{M}}$ de racines simples relativement \`a $P$. A toute racine $\alpha\in \Delta_{\tilde{P}}$, on associe une coracine $\check{\alpha}\in {\cal A}_{\tilde{M}}$ (la normalisation pr\'ecise de $\check{\alpha}$ n'importe pas, la demi-droite qu'elle porte \'etant d\'efinie sans ambigu\"{\i}t\'e). On note ${\mathbb Z}(\check{\Delta}_{\tilde{P}})$ le r\'eseau de ${\cal A}_{\tilde{M}}^{\tilde{G}}$ engendr\'e par ces coracines. On d\'efinit la fonction m\'eromorphe
 $$\epsilon_{\tilde{P}}^{\tilde{G}}(\Lambda)=mes({\cal A}_{\tilde{M}}^{\tilde{G}}/{\mathbb Z}[\check{\Delta}_{\tilde{P}}])\prod_{\alpha\in \Delta_{\tilde{P}}}<\Lambda,\check{\alpha}>^{-1}$$
 sur le complexifi\'e ${\cal A}_{\tilde{M},{\mathbb C}}^*$. On d\'eduit de la $(\tilde{G},\tilde{M})$-famille une fonction
 $$c_{\tilde{M}}^{\tilde{G}}(\Lambda)=\sum_{\tilde{P}\in {\cal P}(\tilde{M})}c(\tilde{P};\Lambda)\epsilon_{\tilde{P}}^{\tilde{G}}(\Lambda).$$
 
 D'autre part, pour $v\in V$ et $\tilde{Q}^v=\tilde{L}^v\tilde{U}_{Q^v}\in {\cal F}(\tilde{R}^v)$, on d\'eduit de $(c_{v}(\tilde{S}^v;\Lambda))_{\tilde{S}^v\in {\cal P}(\tilde{R}^v)}$ une $(\tilde{L}^v,\tilde{R}^v)$-famille $(c_{v}(\tilde{S}^v;\Lambda))_{\tilde{S}^v\in {\cal P}(\tilde{R}^v); \tilde{S}^v\subset \tilde{Q}^v}$, puis  une fonction $c_{v,\tilde{R}^v}^{\tilde{Q}^v}(\Lambda)$. Posons $\tilde{R}^V=(\tilde{R}^v)_{v\in V}$ et notons ${\cal L}(\tilde{R}^V)$ l'ensemble des familles $\tilde{L}^V=(\tilde{L}^v)_{v\in V}$ o\`u $\tilde{L}^v\in {\cal L}(\tilde{R}^v)$ pour tout $v\in V$. Pour une telle famille, on pose
 $${\cal A}_{\tilde{L}^V}=\oplus_{v\in V}{\cal A}_{\tilde{L}^v},\,\,{\cal A}_{\tilde{L}^V}^{\tilde{G}}=\oplus_{v\in V}{\cal A}_{\tilde{L}^v}^{\tilde{G}}$$
 o\`u ${\cal A}_{\tilde{L}^v}^{\tilde{G}}$ est l'orthogonal dans ${\cal A}_{\tilde{L}^v}$ de l'image naturelle de ${\cal A}_{\tilde{G}}$ (notons que cette image est incluse dans ${\cal A}_{\tilde{G}_{v}}$ mais l'inclusion peut \^etre stricte: $\tilde{G}$ peut \^etre plus d\'eploy\'e sur $F_{v}$ que sur $F$). L'espace ${\cal A}_{\tilde{R}^V}^{\tilde{G}}$ contient ${\cal A}_{\tilde{L}^V}^{\tilde{G}}$ comme sous-espace. Il contient aussi l'espace $\Delta_{V}({\cal A}_{\tilde{M}}^{\tilde{G}})$  o\`u $\Delta_{V}$ est le plongement diagonal.  On d\'efinit le coefficient $d_{\tilde{R}^V}^{\tilde{G}}(\tilde{M},\tilde{L}^V)$. Il est nul sauf si
 $${\cal A}_{\tilde{R}^V}^{\tilde{G}}=\Delta_{V}({\cal A}_{\tilde{M}}^{\tilde{G}})\oplus {\cal A}_{\tilde{L}^V}^{\tilde{G}}.$$
 Si cette \'egalit\'e est v\'erifi\'ee, c'est le rapport entre la mesure sur le membre de droite et celle sur le membre de gauche. Pour tout $\tilde{L}^V$ tel que ce nombre soit non nul, Arthur d\'efinit, au moyen d'une donn\'ee auxiliaire, une famille $(\tilde{Q}^v)_{v\in V}$ telle que $\tilde{Q}^v\in {\cal P}(\tilde{L}^v)$ pour tout $v\in V$. On a alors la formule
 
 $$(1) \qquad c_{\tilde{M}}^{\tilde{G}}(\Lambda)=\sum_{\tilde{L}^V\in {\cal L}(\tilde{R}^V)}d_{\tilde{R}^V}^{\tilde{G}}(\tilde{M},\tilde{L}^V)\prod_{v\in V}c_{v,\tilde{R}^v}^{\tilde{Q}^v}(\Lambda).$$
 Cf. [A2] proposition 7.1.
 
 On appliquera souvent cette formule \`a la famille $\tilde{R}^V=(\tilde{M}_{v})_{v\in V}$. On note $\tilde{M}_{V}$ cette famille.
 
 \bigskip
 
 \subsection{Caract\`eres pond\'er\'es}
 Soient $\tilde{M} \in {\cal L}(\tilde{M}_{0})$ et $V$ un ensemble fini de places de $F$. Rappelons que, pour une place $v\in Val(F)$, on a d\'efini en [W1] 2.5 la notion de $\omega$-repr\'esentation de $\tilde{G}(F_{v})$.  A toute telle $\omega$-repr\'esentation $\tilde{\pi}$ est associ\'ee une repr\'esentation sous-jacente $\pi$ de $G(F_{v})$.  On note $\tilde{{\cal A}}_{\tilde{G}_{v}}^*$ l'espace des applications affines $\tilde{\lambda}:\tilde{{\cal A}}_{\tilde{G}_{v}}\to {\mathbb R}$. A tout $\tilde{\lambda}\in \tilde{{\cal A}}_{\tilde{G}_{v}}^*$ est associ\'ee une forme lin\'eaire $\lambda$ sur ${\cal A}_{\tilde{G}_{v}}$. Si $\tilde{\pi}$ est une $\omega$-repr\'esentation de $\tilde{G}(F_{v})$ et $\tilde{\lambda}\in \tilde{{\cal A}}_{\tilde{G}_{v}}^*$, on d\'efinit la $\omega$-repr\'esentation $\tilde{\pi}_{\tilde{\lambda}}$ par $\tilde{\pi}_{\tilde{\lambda}}(\gamma)=e^{<\tilde{\lambda},\tilde{H}_{\tilde{G}_{v}}(\gamma)>}\tilde{\pi}(\gamma)$, o\`u on note $<\tilde{\lambda},\tilde{H}_{\tilde{G}_{v}}(\gamma)>$ l'\'evaluation de $\tilde{\lambda}$ au point $\tilde{H}_{\tilde{G}_{v}}(\gamma)$. Pour une $\omega$-repr\'esentation $\tilde{\pi}$ admissible et de longueur finie de $\tilde{M}(F_{v})$, on a d\'efini, \`a la suite d'Arthur, le caract\`ere pond\'er\'e $f\mapsto J_{\tilde{M}_{v}}^{\tilde{G}_{v}}(\tilde{\pi},f)$. Dans la suite de ce paragraphe, on va globaliser cette d\'efinition.
 
 Pour tout $v\in V$, soit $\tilde{\pi}_{v}$ une $\omega$-repr\'esentation de $\tilde{M}(F_{v})$, admissible et de longueur finie. Posons $\tilde{\pi}_{V}=\otimes_{v\in V} \tilde{\pi}_{v}$. Fixons $\tilde{P}\in {\cal P}(\tilde{M})$, introduisons la repr\'esentation induite $Ind_{\tilde{P}}^{\tilde{G}}(\tilde{\pi})$ de $\tilde{G}(F_{V})$, que l'on r\'ealise dans son espace habituel que l'on note $V_{\pi,P}$. Supposons dans un premier temps que $\tilde{\pi}$ est en position g\'en\'erale de sorte que les op\'erateurs d'entrelacement qui vont appara\^{\i}tre soient bien d\'efinis et inversibles. Pour $\tilde{Q}\in {\cal P}(\tilde{M})$, l'op\'erateur $J_{P\vert Q}(\pi)J_{Q\vert P}(\pi)$ est un automorphisme de $V_{\pi,P}$. Notons $\mu_{Q\vert P}(\pi)$ son inverse. Pour $\Lambda\in i{\cal A}_{\tilde{M}}^*$, posons
 $${\cal M}(\pi;\Lambda,\tilde{Q})=\mu_{Q\vert P}(\pi)^{-1}\mu_{Q\vert P}(\pi_{\Lambda/2})J_{Q\vert P}(\pi)^{-1}J_{Q\vert P}(\pi_{\Lambda}).$$
 La famille $({\cal M}(\pi;\Lambda,\tilde{Q}))_{\tilde{Q}\in {\cal P}(\tilde{M})}$ est une $(\tilde{G},\tilde{M})$-famille \`a valeurs op\'erateurs.  On en d\'eduit un op\'erateur ${\cal M}_{\tilde{M}}^{\tilde{G}}(\pi;\Lambda)$. On pose ${\cal M}_{\tilde{M}}^{\tilde{G}}(\pi)={\cal M}_{\tilde{M}}^{\tilde{G}}(\pi;0)$. Le caract\`ere pond\'er\'e est la forme lin\'eaire sur $C_{c}^{\infty}(\tilde{G}(F_{V}),K)$ d\'efinie par
 $$J_{\tilde{M}}^{\tilde{G}}(\tilde{\pi},f)=trace({\cal M}_{\tilde{M}}^{\tilde{G}}(\pi)Ind_{\tilde{P}}^{\tilde{G}}(\tilde{\pi},f)).$$
 On v\'erifie que cette d\'efinition ne d\'epend pas de l'espace parabolique $\tilde{P}$ choisi. Au moins si $\tilde{\pi}$ est temp\'er\'ee, on peut supprimer la condition de $K$-finitude et d\'efinir $J_{\tilde{M}}^{\tilde{G}}(\tilde{\pi},f)$ pour tout $f\in C_{c}^{\infty}(\tilde{G}(F_{V}))$. 
 
 En utilisant la formule de descente du paragraphe pr\'ec\'edent  appliqu\'ee \`a $\tilde{M}_{V} $ et l'ind\'ependance de l'espace parabolique que l'on vient d'indiquer, on montre que l'on a l'\'egalit\'e
 
 $$(1) \qquad J_{\tilde{M}}^{\tilde{G}}(\tilde{\pi},f)=\sum_{\tilde{L}^V\in {\cal L}(\tilde{M}_{V})}d_{\tilde{M}_{V}}^{\tilde{G}}(\tilde{M},\tilde{L}^V)\prod_{v\in V}J_{\tilde{M}_{v}}^{\tilde{L}^{v}}(\tilde{\pi}_{v},f_{v,\tilde{Q}^v,\omega}) .$$

 Levons l'hypoth\`ese que $\tilde{\pi}$ est en position g\'en\'erale.  L'ensemble $\oplus_{v\in V}\tilde{{\cal A}}_{\tilde{M}_{v}}$ est un espace affine sous $\oplus_{v\in V}{\cal A}_{\tilde{M}_{v}}$, lequel se projette naturellement sur ${\cal A}_{\tilde{M}}$. On note $\tilde{{\cal A}}_{\tilde{M}}$ le quotient de $\oplus_{v\in V}\tilde{{\cal A}}_{\tilde{M}_{v}}$ par le noyau de cette projection. C'est un espace affine sous ${\cal A}_{\tilde{M}}$. On note encore $\tilde{{\cal A}}_{\tilde{M}}^*$ l'espace des fonctions affines sur cet espace, et $\tilde{{\cal A}}_{\tilde{M},{\mathbb C}}^*$ son complexifi\'e. Pour $\tilde{\lambda}\in \tilde{{\cal A}}_{\tilde{M},{\mathbb C}}^*$, on note $\lambda\in {\cal A}_{\tilde{M},{\mathbb C}}^*$ la forme lin\'eaire sous-jacente \`a $\tilde{\lambda}$. 
  Pour $\tilde{\pi}$ quelconque et pour $\tilde{\lambda}\in \tilde{{\cal A}}_{\tilde{M},{\mathbb C}}^*$ en position g\'en\'erale, $\tilde{\pi}_{\tilde{\lambda}}$ est en position g\'en\'erale et l'op\'erateur ${\cal M}_{\tilde{M}}^{\tilde{G}}(\pi_{\lambda})$ comme la forme lin\'eaire $f\mapsto J_{\tilde{M}}^{\tilde{G}}(\tilde{\pi}_{\tilde{\lambda}},f)$ sont bien d\'efinis. Ces termes sont m\'eromorphes en $\tilde{\lambda}$. S'ils sont r\'eguliers en $\tilde{\lambda}=0$, on note ${\cal M}_{\tilde{M}}^{\tilde{G}}(\pi)$ et $f\mapsto J_{\tilde{M}}^{\tilde{G}}(\tilde{\pi},f)$ leurs valeurs en $\tilde{\lambda}=0$.

  \ass{Proposition}{Si $\tilde{\pi}$ est unitaire, les termes ci-dessus sont r\'eguliers en $\tilde{\lambda}=0$.}

Preuve. La formule (1) nous ramene au  cas local, qui est trait\'e par Arthur ([A4], proposition 2.3). $\square$

Evidemment, la formule (1) s'\'etend au cas o\`u tous les termes de la formule sont d\'efinis.

 \subsection{L'application $\phi_{\tilde{M}}$}
 Pour toute place $v\in Val_{f}(F)$, on note  $C_{ac}^{\infty}(\tilde{G}(F_{v}))$ l'espace des fonctions  $f:\tilde{G}(F_{v})\to {\mathbb C}$ telles que:
 
 (1)  pour toute fonction $b\in C_{c}^{\infty}(\tilde{{\cal A}}_{\tilde{G}})$, la fonction $\gamma\mapsto b(\tilde{H}_{\tilde{G}_{v}}(\gamma))f(\gamma)$ appartient \`a $C_{c}^{\infty}(\tilde{G}(F_{v}))$;
 
 - il existe un sous-groupe ouvert compact $H$ de $G(F_{v})$ tel que $f$ soit biinvariante par $H$.
 
 Pour toute place $v\in Val_{\infty}(F)$, on note $C_{ac}^{\infty}(\tilde{G}(F_{v}))$ l'espace des fonctions $f:\tilde{G}(F_{v})\to {\mathbb C}$ qui v\'erifient (1). On note $C_{ac}^{\infty}(\tilde{G}(F_{v}),K_{v})$ le sous-espace des \'el\'ements $K_{v}$-finis \`a droite et \`a gauche de $C_{ac}^{\infty}(\tilde{G}(F_{v}))$.
 
 On note $I_{ac}(\tilde{G}(F_{v}),\omega)$ le quotient de $C_{ac}^{\infty}(\tilde{G}(F_{v}))$ par le sous-espace des $f\in C_{ac}^{\infty}(\tilde{G}(F_{v}))$ telles que $I^{\tilde{G}_{v}}(\gamma,\omega,f)=0$ pour tout \'el\'ement $\gamma\in \tilde{G}_{reg}(F_{v})$ (on rappelle que l'on note ainsi l'ensemble des \'el\'ements semi-simples et  fortement r\'eguliers de $\tilde{G}(F_{v})$). Si $v$ est archim\'edienne, on d\'efinit de m\^eme la variante $I_{ac}(\tilde{G}(F_{v}),\omega,K_{v})$.

 Soient $\tilde{M}$ un espace de Levi de $\tilde{G}$ et $V$ un ensemble fini de places de $F$. Pour $v\in V$, on d\'efinit une application lin\'eaire $\phi_{\tilde{M}_{v}}:C^{\infty}_{ac}(\tilde{G}(F_{v}))\to I_{ac}(\tilde{M}(F_{v}),\omega)$. Elle est d\'efinie en [W1] 6.4 dans le cas o\`u $v$ est non-archim\'edienne, en [V] 1.2  dans le cas o\`u $v$ est archim\'edienne. Dans ce dernier cas, l'application se restreint en une application lin\'eaire $C_{ac}^{\infty}(\tilde{G}(F_{v}),K_{v})\to I_{ac}(\tilde{M}(F_{v}),\omega,K_{v})$, cf. [W1] 6.4.
 On definit comme en 1.1 les espaces $C^{\infty}_{ac}(\tilde{G}(F_{V}))$ et $I_{ac}(\tilde{G}(F_{V}),\omega)$ (avec la variante $I_{ac}(\tilde{G}(F_{V}),\omega,K)$ si $V$ contient des places archim\'ediennes). Comme dans le paragraphe pr\'ec\'edent, on applique les d\'efinitions   du paragraphe  1.2 \`a la famille $\tilde{M}_{V}=(\tilde{M}_{v})_{v\in V}$.    On d\'efinit une application $\phi_{\tilde{M}}:C^{\infty}_{ac}(\tilde{G}(F_{V}))\to I_{ac}(\tilde{M}(F_{V}),\omega)$ par
 $$(2) \qquad \phi_{\tilde{M}}(f)=\sum_{\tilde{L}^V\in {\cal L}(\tilde{M}_{V})}d_{\tilde{M}_{V}}^{\tilde{G}}(\tilde{M},\tilde{L}^V)\left(\otimes_{v\in V}\phi_{\tilde{M}_{v}}^{\tilde{L}^v}(f_{v,\tilde{Q}^v,\omega})\right)$$
 pour une fonction $f=\otimes_{v\in V}f_{v}\in C^{\infty}_{ac}(\tilde{G}(F_{V}))$. 
 
 Cette d\'efinition d\'epend a priori d'une donn\'ee auxiliaire puisque les $\tilde{Q}^v$ en d\'ependent. Pour qu'elle soit loisible, on doit montrer qu'en fait, elle n'en d\'epend pas. Pour cela, on utilise la caract\'erisation de [W1] 6.4(5).  Fixons pour tout $v\in V$ une $\omega$-repr\'esentation temp\'er\'ee $M$-irr\'eductible $\tilde{\pi}_{v}$ de $\tilde{M}(F_{v})$. Fixons aussi $X_{v}\in \tilde{{\cal A}}_{\tilde{M}_{v}}$. On a d\'efini une forme lin\'eaire $\varphi_{v}\mapsto I^{\tilde{M}_{v}}(\tilde{\pi}_{v},X_{v},\varphi_{v})$ sur $I_{ac}(\tilde{M}(F_{v}),\omega)$, cf. [W1] 6.4 et [V] 1.2. Quand $\tilde{\pi}_{v}$ et $X_{v}$ varient, ces formes lin\'eaires s\'eparent les \'el\'ements de $I_{ac}(\tilde{M}(F_{v}),\omega)$. Il suffit donc de prouver que la valeur sur le membre de droite de (2) de la forme lin\'eaire 
$$\otimes_{v\in V}\varphi_{v}\mapsto \prod_{v\in V}I^{\tilde{M}_{v}}(\tilde{\pi}_{v},X_{v},\varphi_{v})$$
est bien d\'efinie. Cette valeur est
 $$\sum_{\tilde{L}^V\in {\cal L}(\tilde{M}_{V})}d_{\tilde{M}_{V}}^{\tilde{G}}(\tilde{M},\tilde{L}^V)\prod_{v\in V}I^{\tilde{M}_{v}}(\tilde{\pi}_{v},X_{v},\phi_{\tilde{M}_{v}}^{\tilde{L}^v}(f_{v,\tilde{Q}^v,\omega})),$$
 ou encore
 $$\sum_{\tilde{L}^V\in {\cal L}(\tilde{M}_{V})}d_{\tilde{M}_{V}}^{\tilde{G}}(\tilde{M},\tilde{L}^V)\prod_{v\in V}J_{\tilde{M}_{v}}^{\tilde{L}^{v}}(\tilde{\pi}_{v},X_{v},f_{v,\tilde{Q}^v,\omega}).$$
 Par d\'efinition,  on a pour tout $v\in V$ une \'egalit\'e
 $$J_{\tilde{M}_{v}}^{\tilde{L}^{v}}(\tilde{\pi}_{v},X_{v},f_{v,\tilde{Q}^v,\omega})=c_{v}\int_{i{\cal A}_{\tilde{M}_{v},F_{v}}^*}J_{\tilde{M}_{v}}^{\tilde{L}^{v}}(\tilde{\pi}_{v,\tilde{\lambda}_{v}},f_{v,\tilde{Q}^v,\omega})e^{-<\tilde{\lambda}_{v},X_{v}>}d\lambda_{v}$$
 (les deux fonctions que l'on int\`egre d\'ependent de $\tilde{\lambda}_{v}\in i\tilde{{\cal A}}_{\tilde{M}_{v}}$; leur produit se descend en une fonction sur $i{\cal A}_{\tilde{M}_{v},F_{v}}^*$; $c_{v}$ est une constante d\'ependant seulement des mesures de Haar). Il suffit donc de prouver que, pour toute famille $(\tilde{\lambda}_{v})_{v\in V}\in \oplus i\tilde{{\cal A}}_{\tilde{M}_{v}}^*$, le terme
  $$\sum_{\tilde{L}^V\in {\cal L}(\tilde{M}_{V})}d_{\tilde{M}_{V}}^{\tilde{G}}(\tilde{M},\tilde{L}^V)\prod_{v\in V}J_{\tilde{M}_{v}}^{\tilde{L}^{v}}(\tilde{\pi}_{v,\tilde{\lambda}_{v}},f_{v,\tilde{Q}^v,\omega}) $$
  est bien d\'efini. Quitte \`a remplacer $\tilde{\pi}_{v}$ par $\tilde{\pi}_{v,\tilde{\lambda}_{v}}$, il suffit de consid\'erer
  $$ \sum_{\tilde{L}^V\in {\cal L}(\tilde{M}_{V})}d_{\tilde{M}_{V}}^{\tilde{G}}(\tilde{M},\tilde{L}^V)\prod_{v\in V}J_{\tilde{M}_{v}}^{\tilde{L}^{v}}(\tilde{\pi}_{v},f_{v,\tilde{Q}^v,\omega}) .$$
  La formule 1.5(1) dit que cette expression est \'egale \`a $J_{\tilde{M}}^{\tilde{G}}(\tilde{\pi},f)$, o\`u $\tilde{\pi}=\otimes_{v\in V}\tilde{\pi}_{v}$. Ce terme ne d\'ependant d'aucun param\`etre auxiliaire, cela d\'emontre l'assertion.
  
  \bigskip
  
  \subsection{Une propri\'et\'e globale de l'application $\phi_{\tilde{M}}$}
  La situation est la m\^eme que dans le paragraphe pr\'ec\'edent, mais on suppose que $V$ contient $V_{ram}$. On se rappelle l'application $\tilde{H}_{\tilde{G}_{V}}:\tilde{G}(F_{V})\to {\cal A}_{\tilde{G}}$ de 1.1. L'espace $ C^{\infty}({\cal A}_{\tilde{G}})$ des fonctions $C^{\infty}$ sur ${\cal A}_{\tilde{G}}$ op\`ere sur $C_{c}^{\infty}(\tilde{G}(F_{V}))$: \`a $b\in C^{\infty}({\cal A}_{\tilde{G}})$ et $f\in C_{c}^{\infty}(\tilde{G}(F_{V}))$, on associe la fonction produit $f(b\circ \tilde{H}_{\tilde{G}_{V}})$. L'espace $C^{\infty}({\cal A}_{\tilde{G}})$ op\`ere de m\^eme sur $C^{\infty}_{ac}(\tilde{G}(F_{V}))$. Ces actions se descendent en des actions sur $I(\tilde{G}(F_{V}),\omega)$ et $I_{ac}(\tilde{G}(F_{V}),\omega)$. Notons $C^{\infty}_{ac,glob}(\tilde{G}(F_{V}))$ le sous-espace des $f\in C^{\infty}_{ac}(\tilde{G}(F_{V}))$ tels que $f(b\circ \tilde{H}_{\tilde{G}_{V}})\in C_{c}^{\infty}(\tilde{G}(F_{V}))$ pour tout $b\in C_{c}^{\infty}({\cal A}_{\tilde{G}})$. Notons $I_{ac,glob}(\tilde{G}(F_{V}),\omega)$ l'image de cet espace dans $I_{ac}(\tilde{G}(F_{V}),\omega)$.
  
  \ass{Lemme}{L'homomorphisme $\phi_{\tilde{M}}$ envoie $C^{\infty}_{ac,glob}(\tilde{G}(F_{V}))$ dans $I_{ac,glob}(\tilde{M}(F_{V}),\omega)$.}
  
  Preuve. Soient $f\in C^{\infty}_{ac,glob}(\tilde{G}(F_{V}))$ et $b\in C_{c}^{\infty}(\tilde{A}_{\tilde{M}})$. On doit montrer que $\phi_{\tilde{M}}(f)(b\circ\tilde{H}_{\tilde{M}_{V}})$ appartient \`a $I(\tilde{M}(F_{V}),\omega)$. Soit $b'\in C_{c}^{\infty}(\tilde{A}_{\tilde{G}})$ qui vaut $1$ sur la projection dans ${\cal A}_{\tilde{G}}$ du support de $b$. Alors $b\circ\tilde{H}_{\tilde{M}_{V}}=(b\circ\tilde{H}_{\tilde{M}_{V}})(b'\circ\tilde{H}_{\tilde{G}_{V}})$. Il r\'esulte de la d\'efinition de $\phi_{\tilde{M}}$ que $\phi_{\tilde{M}}(f)(b'\circ\tilde{H}_{\tilde{G}_{V}})=\phi_{\tilde{M}}(f(b'\circ\tilde{H}_{\tilde{G}_{V}}))$. D'o\`u
  $$\phi_{\tilde{M}}(f)(b\circ\tilde{H}_{\tilde{M}_{V}})=\phi_{\tilde{M}}(f(b'\circ\tilde{H}_{\tilde{G}_{V}}))(b\circ\tilde{H}_{\tilde{M}_{V}}).$$
  Quitte \`a remplacer $f$ par $f(b'\circ\tilde{H}_{\tilde{G}_{V}})$, on est ramen\'e au cas o\`u $f\in C_{c}^{\infty}(\tilde{G}(F_{V}))$. En consid\'erant la formule (1) du paragraphe pr\'ec\'edent, on voit qu'il nous suffit de fixer $\tilde{L}^V\in {\cal L}(\tilde{M}_{V})$ tel que $d_{\tilde{M}_{V}}^{\tilde{G}}(\tilde{M},\tilde{L}^V)\not=0$ et de prouver que la fonction
  $$(1) \qquad \left(\otimes_{v\in V}\phi_{\tilde{M}_{v}}^{\tilde{L}^{v}}(f_{v,\tilde{Q}_{v},\omega})\right)(b\circ\tilde{H}_{\tilde{M}_{V}})$$
  appartient \`a $I_{ac,glob}(\tilde{M}(F_{V}),\omega)$. On rel\`eve chaque terme $\phi_{\tilde{M}_{v}}^{\tilde{L}^{v}}(f_{v,\tilde{Q}_{v},\omega})$ en un \'el\'ement de $C^{\infty}_{ac}(\tilde{M}(F_{v}))$. Notons $\Xi_{v}$ l'image de son support  dans $\tilde{{\cal A}}_{\tilde{M}_{v}}$ par l'application $\tilde{H}_{\tilde{M}_{v}}$. Parce que $f_{v}$ est maintenant \`a support compact, on sait que l'on peut supposer que la projection naturelle  de $\Xi_{v}$ dans $\tilde{{\cal A}}_{\tilde{L}^v}$ est compacte.  Notons $\Xi$ l'image du support de la fonction (1) dans $\tilde{{\cal A}}_{\tilde{M}_{V}}$ par l'application $\prod_{v\in V}\tilde{H}_{\tilde{M}_{v}}$. Pour que la fonction (1) appartienne \`a  $I_{ac,glob}(\tilde{M}(F_{V}),\omega)$, il suffit que $\Xi$ soit compact. Or $\Xi$ est le sous-ensemble des \'el\'ements de $\prod_{v\in V}\Xi_{v}$ dont l'image par $\tilde{p}_{V}$ appartient au support compact de la fonction $b$. En fixant des points bases de nos espaces affines, qui permettent d'identifier $\tilde{{\cal A}}_{\tilde{M}_{V}}$ \`a ${\cal A}_{\tilde{M}_{V}}$, on est ramen\'e \`a la situation suivante. On a un sous-ensemble ferm\'e $\Xi\subset {\cal A}_{\tilde{M}_{V}} $ dont la projection dans ${\cal A}_{\tilde{L}^V}$ est compacte et dont l'image par l'application
  $$\begin{array}{cccc}p_{V}:&{\cal A}_{\tilde{M}_{V}}&\to&{\cal A}_{\tilde{M}}\\ &(H_{v})_{v\in V}&\mapsto& \sum_{v\in V}H_{v,\tilde{M}}\\ \end{array}$$
  est compacte. On doit montrer que ce sous-ensemble est compact. Il suffit que l'intersection des noyaux de $p_{V}$ et de la projection dans ${\cal A}_{\tilde{L}^V}$ soit r\'eduite \`a $0$. Ou encore que la somme de ${\cal A}_{\tilde{L}^V}$ et de l'orthogonal du noyau de $p_{V}$ soit l'espace ${\cal A}_{\tilde{M}_{V}}$ tout entier. Cela r\'esulte de la condition $d_{\tilde{M}_{V}}^{\tilde{G}}(\tilde{M},\tilde{L}^V)\not=0$ puisque l'orthogonal du noyau de $p_{V}$ est l'espace ${\cal A}_{\tilde{M}}$ plong\'e diagonalement dans ${\cal A}_{\tilde{M}_{V}}$. $\square$
  
  On a la variante suivante: $\phi_{\tilde{M}}$ envoie $C_{ac,glob}^{\infty}(\tilde{G}(F_{V}),K)$ dans $I_{ac,glob}(\tilde{M}(F_{V}),\omega,K)$.
  
  \bigskip
  
  \subsection{Espaces de distributions}
 Soient $\tilde{M}$ un espace de Levi de $\tilde{G}$ et $V$ un ensemble fini de places de $F$.   Pour tout $v\in V$, on a d\'efini l'espace de distributions $D_{g\acute{e}om}(\tilde{M}(F_{v}),\omega)$ en [I] 5.1 et 5.2. C'est celui des distributions $\omega$-\'equivariantes support\'ees par un nombre fini de classes de conjugaison. Supposons $v$ archim\'edienne. On a d\'efini en [V] 1.3 et 2.1 les sous-espaces
  $$D_{orb}(\tilde{M}(F_{v}),\omega)\subset D_{tr-orb}(\tilde{M}(F_{v}),\omega) \subset D_{g\acute{e}om}(\tilde{M}(F_{v}),\omega)\supset D_{g\acute{e}om,\tilde{G}-\acute{e}qui}(\tilde{M}(F_{v}),\omega).$$
  L'espace $D_{orb}(\tilde{M}(F_{v}),\omega)$ est engendr\'e par les int\'egrales orbitales ordinaires ($\omega$-\'equivariantes).   L'espace $D_{g\acute{e}om,\tilde{G}-\acute{e}qui}(\tilde{M}(F_{v}),\omega)$ est le sous-espace des \'el\'ements de $D_{g\acute{e}om}(\tilde{M}(F_{v}),\omega)$ dont le support est form\'e d'\'el\'ements de $\tilde{M}(F_{v})$ qui sont $\tilde{G}$-\'equisinguliers. L'espace $D_{tr-orb}(\tilde{M}(F_{v}),\omega)$ est 
 d\'efini par r\'ecurrence.   Il est engendr\'e par $D_{orb}(\tilde{M}(F_{v}),\omega)$ et par les images par transfert des espaces $D_{tr-orb}({\bf M}'_{v})$ pour les donn\'ees endoscopiques elliptiques ${\bf M}_{v}'$ de $(M_{v},\tilde{M}_{v})$, avec la restriction  $M_{v}'\not=M_{v}$ si $(M_{v},\tilde{M}_{v},{\bf a}_{v})$ est quasi-d\'eploy\'e et \`a torsion int\'erieure.
    
{\bf Remarque.} Si on applique les m\^emes d\'efinitions pour une place $v$ non-archim\'edienne, on a l'\'egalit\'e  $D_{orb}(\tilde{M}(F_{v}),\omega) =D_{tr-orb}(\tilde{M}(F_{v}),\omega)=D_{g\acute{e}om}(\tilde{M}(F_{v}),\omega)$.
  
    On pose $D_{g\acute{e}om}(\tilde{M}(F_{V}),\omega)=\otimes_{v\in V}D_{g\acute{e}om}(\tilde{M}(F_{v}),\omega)$. On d\'efinit les sous-espaces $D_{orb}(\tilde{M}(F_{V}),\omega)$, $D_{tr-orb}(\tilde{M}(F_{V}),\omega)$  et $D_{g\acute{e}om,\tilde{G}-\acute{e}qui}(\tilde{M}(F_{V}),\omega)$ en rempla\c{c}ant, pour toute place  archim\'edienne $v\in V$, l'espace  $D_{g\acute{e}om}(\tilde{M}(F_{v}),\omega)$ par $D_{orb}(\tilde{M}(F_{v}),\omega)$ etc... 
  
  \bigskip
  
  \subsection{Int\'egrales orbitales pond\'er\'ees}
  Soient $\tilde{M}\in {\cal L}(\tilde{M}_{0})$ et $V$ un ensemble fini de places de $F$. Pour tout $g=(g_{v})_{v\in V}\in G(F_{V})$, Arthur d\'efinit une $(\tilde{G},\tilde{M})$-famille $(v_{\tilde{P}}(g;\Lambda))_{\tilde{P}\in {\cal P}(\tilde{M})}$, d'o\`u une fonction $v_{\tilde{M}}^{\tilde{G}}(g;\Lambda)$. On pose $v_{\tilde{M}}^{\tilde{G}}(g)=v_{\tilde{M}}^{\tilde{G}}(g,0)$.
 
  Soit $\gamma=(\gamma_{v})_{v\in V}\in \tilde{M}(F_{V})$.  On  pose $M_{\gamma}=\prod_{v\in V}M_{v,\gamma_{v}}$ que l'on peut consid\'erer comme un groupe d\'efini sur l'anneau $F_{V}$. On d\'efinit de m\^eme $G_{\gamma}$, $Z_{M}(\gamma)$ etc... Si $\omega$ n'est pas trivial sur $M_{\gamma}(F_{V})$, on pose $J_{\tilde{M}}^{\tilde{G}}(\gamma,\omega,f)=0$ pour tout $f\in C_{c}^{\infty}(\tilde{G}(F_{V}))$. Supposons dans la suite que $\omega$ est trivial sur $M_{\gamma}(F_{V})$. On fixe des mesures de Haar sur $G(F_{V})$ et $M_{\gamma}(F_{V})$. Supposons d'abord que $\gamma$ soit $\tilde{G}$-\'equisingulier, c'est-\`a-dire que $M_{\gamma}=G_{\gamma}$. Pour $f\in C_{c}^{\infty}(\tilde{G}(F_{V}))$, on pose
  $$J_{\tilde{M}}^{\tilde{G}}(\gamma,\omega,f)=D^{\tilde{G}}(\gamma)^{1/2}\int_{M_{\gamma}(F_{V})\backslash G(F_{V})}\omega(g)f(g^{-1}\gamma g)v_{\tilde{M}}(g)\,dg.$$
  Pour $\gamma$ quelconque, Arthur d\'efinit $J_{\tilde{M}}^{\tilde{G}}(\gamma,\omega,f)$ par un proc\'ed\'e de limite. Nous allons le rappeler bri\`evement, tout en le modifiant. Soit $v\in V$ et $a_{v}\in A_{\tilde{M}_{v}}(F_{v})$. On a d\'efini en [II] 1.5 une $(\tilde{G}_{v},\tilde{M}_{v})$-famille $(r_{\tilde{P}}(\gamma_{v},a_{v};\Lambda))_{\tilde{P}\in {\cal P}(\tilde{M}_{v})}$, pourvu que $a_{v}$ soit en position g\'en\'erale.  Plus pr\'ecis\'ement, notons $\eta_{v}$ la partie semi-simple de $\gamma_{v}$.  Il suffit que $a_{v}$ v\'erifie $\alpha(a_{v})\not=\pm 1$ pour toute racine $\alpha$ de $A_{M_{v,\eta_{v}}}$ dans $G_{v,\eta_{v}}$ pour que les fonctions pr\'ec\'edentes soient d\'efinies. Signalons que la d\'efinition de   ces fonctions est  l\'eg\`erement diff\'erente de celle d'Arthur. Soit maintenant $a=(a_{v})_{v\in V}\in A_{\tilde{M}}(F_{V})$. Si $a$ est en position g\'en\'erale, les $a_{v}$ ne sont pas v\'eritablement "en position g\'en\'erale" car le tore $A_{M_{v},\eta_{v}}$ est en g\'en\'eral plus gros que le localis\'e de $A_{\tilde{M}}$, mais ils v\'erifient la condition pr\'ecise ci-dessus. Comme en 1.3, on d\'eduit alors des $(\tilde{G}_{v},\tilde{M}_{v})$-familles ci-dessus une famille produit $(r_{\tilde{P}}(\gamma,a;\Lambda))_{\tilde{P}\in {\cal P}(\tilde{M})}$. On en d\'eduit une fonction $r_{\tilde{M}}^{\tilde{G}}(\gamma,a;\Lambda)$ et on pose $r_{\tilde{M}}^{\tilde{G}}(\gamma,a)=r_{\tilde{M}}^{\tilde{G}}(\gamma,a;0)$. Consid\'erons l'expression
  $$\sum_{\tilde{L}\in {\cal L}(\tilde{M})}r_{\tilde{M}}^{\tilde{L}}(\gamma,a)J_{\tilde{L}}^{\tilde{G}}(a\gamma,\omega,f).$$
  Tous les termes sont bien d\'efinis puisque $M_{a\gamma}=G_{a\gamma}$ pour $a$ en position g\'en\'erale. Arthur montre que cette expression a une limite quand $a$ tend vers $1$ ([A5] th\'eor\`eme 5.2).   La modification que l'on a apport\'ee aux d\'efinitions n'affecte  pas cette propri\'et\'e. On note $J_{\tilde{M}}^{\tilde{G}}(\gamma,\omega,f)$  cette limite.
  
  {\bf Remarque.} La d\'efinition est globale: les $(\tilde{G},\tilde{M})$-familles intervenant sont index\'ees par des espaces paraboliques d\'efinis sur $F$. M\^eme dans le cas o\`u $V$ est r\'eduit \`a une seule place $v$, les int\'egrales orbitales pond\'er\'ees ci-dessus ne co\"{\i}ncident pas en g\'en\'eral avec leurs similaires locales relatives au corps de base $F_{v}$. La relation entre les deux objets est donn\'ee par la formule (1) suivante. La m\^eme remarque s'appliquera aux int\'egrales orbitales pond\'er\'ees $\omega$-\'equivariantes d\'efinies au paragraphe suivant.
  
  \bigskip

  En utilisant plusieurs fois la formule 1.4(1), on montre que, pour $f=\otimes_{v\in V}f_{v}$, on a l'\'egalit\'e
  $$(1) \qquad J_{\tilde{M}}^{\tilde{G}}(\gamma,\omega,f)=\sum_{\tilde{L}^V\in {\cal L}(\tilde{M}_{V})}d_{\tilde{M}_{V}}^{\tilde{G}}(\tilde{M},\tilde{L}^V)\prod_{v\in V}J_{\tilde{M}_{v}}^{\tilde{L}^v}(\gamma_{v},\omega,f_{v,\tilde{Q}^v,\omega}).$$
  
    Comme dans le cas local, on peut formaliser les d\'efinitions ci-dessus et les rendre ind\'ependantes de tout choix de mesures en d\'efinissant $J_{\tilde{M}}^{\tilde{G}}(\boldsymbol{\gamma},{\bf f})$ pour $\boldsymbol{\gamma}\in D_{orb}(\tilde{M}(F_{V}),\omega)\otimes Mes(M(F_{V}))^*$ et ${\bf f}\in C_{c}^{\infty}(\tilde{G}(F_{V}))\otimes Mes(G(F_{V}))$. Signalons que ces int\'egrales d\'ependent tout de m\^eme de la mesure fix\'ee sur ${\cal A}_{\tilde{M}}^{\tilde{G}}$. 
  
  On aura besoin d'une variante de la formule (1). Supposons $V$ r\'eunion  disjointe de deux sous-ensembles $V_{1}$ et $V_{2}$. Pour $i=1,2$, soient $\boldsymbol{\gamma}_{i}\in D_{orb}(\tilde{M}(F_{V_{i}}),\omega)\otimes Mes(M(F_{V_{i}}))$ et ${\bf f}_{i}\in C_{c}^{\infty}(\tilde{G}(F_{V_{i}}))\otimes Mes(G(F_{V_{i}}))$. Posons $\boldsymbol{\gamma}=\boldsymbol{\gamma}_{1}\otimes \boldsymbol{\gamma}_{2}$ et ${\bf f}={\bf f}_{1}\otimes {\bf f}_{2}$. 
  Supposons que, pour tout $\tilde{Q}=\tilde{L}\tilde{U}_{Q}\in {\cal F}(\tilde{M})$, l'int\'egrale orbitale pond\'er\'ee  $J_{\tilde{M}}^{\tilde{L}}(\boldsymbol{\gamma}_{1},{\bf f}_{1,\tilde{Q},\omega})$ ne d\'epende que de $\tilde{L}$. On a alors l'\'egalit\'e
 $$(2) \qquad J_{\tilde{M}}^{\tilde{G}}(\boldsymbol{\gamma},{\bf f})=\sum_{\tilde{L}\in {\cal L}(\tilde{M})} J_{\tilde{M}}^{\tilde{L}}(\boldsymbol{\gamma}_{1},{\bf f}_{1,\tilde{Q},\omega})J_{\tilde{L}}^{\tilde{G}}(\boldsymbol{\gamma}^{\tilde{L}}_{2},{\bf f}_{2}),$$
 o\`u $\tilde{Q}$ est un \'el\'ement quelconque de ${\cal P}(\tilde{L})$ et o\`u $\boldsymbol{\gamma}_{2}^{\tilde{L}}$ est l'induite de $\boldsymbol{\gamma}_{2}$ \`a $\tilde{L}(F_{V_{2}})$. 
 
 \bigskip
 
 \subsection{Syst\`eme de fonctions $B$}
 Supposons $(G,\tilde{G},{\bf a})$ quasi-d\'eploy\'e et \`a torsion int\'erieure. Rappelons que dans ce cas, on supprime le caract\`ere   trivial  $\omega$ des notations. Pour une place $v$ de $F$, on a d\'efini en [II] 1.9 la notion de syst\`eme de fonctions $B$ sur $\tilde{G}(F_{v})$. Pour globaliser cette notion, on va en donner une d\'efinition un peu diff\'erente.   Fixons une paire de Borel $(B^*,T^*)$ de $G$ d\'efinie sur $F$. Notons $\tilde{T}^*$ l'ensemble des $\gamma\in \tilde{G}$ tels que $ad_{\gamma}$ conserve cette paire. Pour  $\eta\in \tilde{T}^*$,  notons $\Sigma^{G_{\eta}}(T^*)$ l'ensemble des racines de $T^*$ dans $\mathfrak{g}_{\eta}$.  L'ensemble $\Sigma^{G_{\eta}}(T^*)$ est un sous-ensemble de l'ensemble $\Sigma(T^*)$ des racines de $T^*$ dans $\mathfrak{g}$. Pour tout sous-ensemble $\Sigma'\subset \Sigma(T^*)$, consid\'erons l'ensemble des $\eta\in \tilde{T}^*$ tels que $\Sigma^{G_{\eta}}(T^*)=\Sigma'$. C'est un sous-ensemble alg\'ebrique de $\tilde{T}^*$ que l'on d\'ecompose en composantes connexes. En faisant varier $\Sigma'$, on obtient une d\'ecomposition de $\tilde{T}^*$ en r\'eunion disjointe finie de sous-ensembles alg\'ebriques connexes. On note $\underline{\Omega}$ cet ensemble de sous-ensembles alg\'ebriques. Pour $\Omega\in \underline{\Omega}$, on note $\Sigma(\Omega)$ l'ensemble $\Sigma^{G_{\eta}}(T^*)$ pour un \'el\'ement quelconque $\eta\in \Omega$. Le groupe de Weyl $W$ agit sur $\tilde{T}^*$. Pour $w\in W$ et $\eta\in \tilde{T}^*$, l'\'el\'ement $w$ d\'efinit une bijection $w:\Sigma^{G_{\eta}}(T^*)\to \Sigma^{G_{w(\eta)}}(T^*)$. Le groupe de Galois $\Gamma_{F}$ agit aussi et, pour $\sigma\in \Gamma_{F}$ et $\eta\in \tilde{T}^*$, on a aussi une bijection $\sigma:\Sigma^{G_{\eta}}(T^*)\to \Sigma^{G_{\sigma(\eta)}}(T^*)$. Il en r\'esulte que les actions de $W$ et $\Gamma_{F}$ sur $\tilde{T}^*$ permutent les \'el\'ements de $\underline{\Omega}$. 
 
 On se donne pour tout $\Omega\in \underline{\Omega}$ une fonction $B_{\Omega}:\Sigma(\Omega)\to {\mathbb Q}_{>0}$. On suppose v\'erifi\'ees les conditions (1) et (2) suivantes pour tout $\Omega\in \underline{\Omega}$:
 
 (1) pour toute composante irr\'eductible $\Sigma'$ du syst\`eme de racines $\Sigma(\Omega)$, ou bien $B_{\Omega}$ est constante sur $\Sigma'$, ou bien la fonction $\beta\mapsto \frac{B_{\Omega}(\beta)}{(\beta,\beta)}$ est constante sur $\Sigma'$, o\`u $(.,.)$ est une forme quadratique d\'efinie positive et invariante par le groupe de Weyl sur l'espace $X^*(T^*)\otimes_{{\mathbb Z}}{\mathbb R}$;

 (2) pour $w\in W$, $\sigma\in \Gamma_{F}$ et $\beta\in \Sigma(\Omega)$, on a les \'egalit\'es $B_{w(\Omega)}(w(\beta))=B_{\sigma(\Omega)}(\sigma(\beta))=B_{\Omega}(\beta)$.
 
 A ces conditions, on dit que les fonctions $B_{\Omega}$ forment un "syst\`eme de fonctions $B$" sur $\tilde{G}$.  Fixons un tel syst\`eme.  Il convient d'\'elargir l'ensemble $V_{ram}$ de 1.1 de sorte que
 
 (3) soit $v\in Val(F)-V_{ram}$, notons $p$ la caract\'eristique r\'esiduelle de $v$; alors, pour tout \'el\'ement $\Omega\in \underline{\Omega}$, les valeurs de $B_{\Omega}$ sont premi\`eres \`a $p$. 
 
 C'est possible puisque l'ensemble $\underline{\Omega}$ est  fini. 
 
 Soit $\eta\in \tilde{T}^*(\bar{F})$. Il existe un unique $\Omega\in \underline{\Omega}$ tel que $\eta\in \Omega(\bar{F})$. On pose $B_{\eta}=B_{\Omega}$. Plus g\'en\'eralement, pour un \'el\'ement semi-simple $\eta\in \tilde{G}(\bar{F})$, fixons une paire de Borel $(B,T)$ de $G$ conserv\'ee par $ad_{\eta}$. On d\'efinit comme ci-dessus le syst\`eme de racines $\Sigma^{G_{\eta}}(T)$. On fixe un \'el\'ement $x\in G(\bar{F})$ tel que $ad_{x}(B,T)=(B^*,T^*)$. Alors $ad_{x}$ identifie $\Sigma^{G_{\eta}}(T)$ \`a $\Sigma^{G_{ad_{x}(\eta)}}(T^*)$. En transportant la fonction $B_{ad_{x}(\eta)}$ par cet isomorphisme, on obtient une fonction $B_{\eta}$ sur $\Sigma^{G_{\eta}}(T)$, qui ne d\'epend pas de l'\'el\'ement $x$ choisi. 
 
 Soit $v$ une place de $F$ et soit $\eta$ un \'el\'ement semi-simple de $\tilde{G}(\bar{F}_{v})$. Le m\^eme proc\'ed\'e permet de d\'efinir une fonction $B_{\eta}$ sur le syst\`eme de racines de $G_{\eta}$. La restriction de ces fonctions aux \'el\'ements $\eta\in \tilde{G}(F_{v})$ est un syst\`eme de fonctions $B$ sur $\tilde{G}(F_{v})$, au sens de [II] 1.9.

 Soient $\tilde{M}$ un espace de Levi de $\tilde{G}$, $V$ un ensemble fini de places de $F$ et  $\gamma=(\gamma_{v})_{v\in V}\in \tilde{M}(F_{V})$. Pour $v\in V$, on a d\'efini en [II] 1.9 des   $(\tilde{G}_{v},\tilde{M}_{v})$-familles $(r_{\tilde{P}}(\gamma_{v},a_{v},B;\Lambda))_{\tilde{P}\in {\cal P}(\tilde{M}_{v})}$. En utilisant ces familles dans les constructions du paragraphe pr\'ec\'edent, on d\'efinit l'int\'egale orbitale pond\'er\'ee $J_{\tilde{M}}^{\tilde{G}}(\gamma,B,f)$.  Elle co\"{\i}ncide avec $J_{\tilde{M}}^{\tilde{G}}(\gamma,f)$  dans le cas o\`u $M_{\gamma}=G_{\gamma}$.
  
  \bigskip
  
  \subsection{Int\'egrales orbitales pond\'er\'ees $\omega$-\'equivariantes}
  Soient $\tilde{M}\in {\cal L}(\tilde{M}_{0})$ et $V$ un ensemble fini de places de $F$. Soit $\gamma\in \tilde{M}(F_{V})$. On fixe encore des mesures de Haar sur $G(F_{V})$ et $M_{\gamma}(F_{V})$. Pour $f\in C_{c}^{\infty}(\tilde{G}(F_{V}))$, on d\'efinit l'int\'egrale orbitale pond\'er\'ee $\omega$-\'equivariante par la formule de r\'ecurrence
  $$I_{\tilde{M}}^{\tilde{G}}(\gamma,\omega,f) =J_{\tilde{M}}^{\tilde{G}}(\gamma,\omega,f)-\sum_{\tilde{L}\in {\cal L}(\tilde{M}),\tilde{L}\not=\tilde{G}}I_{\tilde{M}}^{\tilde{L}}(\gamma,\omega,\phi_{\tilde{L}}(f)).$$
  
  {\bf Remarque.} Comme souvent, certaines propri\'et\'es de ces termes sont suppos\'ees connues pour les donn\'ees $(G',\tilde{G}',{\bf a}')$ similaires \`a $(G,\tilde{G},{\bf a})$ telles que  $dim(G'_{SC})<dim(G_{SC})$. Les propri\'et\'es utilis\'ees ici est que $I_{\tilde{M}}^{\tilde{G}}(\gamma,\omega,f)$ ne d\'epend que de l'image de $f$ dans $I(\tilde{G}(F_{V}),\omega)$ et que la d\'efinition s'\'etend \`a $f\in C^{\infty}_{ac}(\tilde{G}(F_{V}))$. Gr\^ace \`a ces hypoth\`eses, les termes $I_{\tilde{M}}^{\tilde{L}}(\gamma,\omega,\phi_{\tilde{L}}(f))$ sont bien d\'efinis pour $\tilde{L}\not=\tilde{G}$. La formule (1) ci-dessous, qui se d\'eduit de la simple formule de d\'efinition ci-dessus, ram\`ene la v\'erification des hypoth\`eses de r\'ecurrence aux propri\'et\'es des int\'egrales analogues locales, pour lesquelles on renvoie \`a [II] et [V].
  
 Pour $f=\otimes_{v\in V}f_{v}$, on a l'\'egalit\'e
  $$(1) \qquad I_{\tilde{M}}^{\tilde{G}}(\gamma,\omega,f)=\sum_{\tilde{L}^V\in {\cal L}(\tilde{M}_{V})}d_{\tilde{M}_{V}}^{\tilde{G}}(\tilde{M},\tilde{L}^V)\prod_{v\in V}I_{\tilde{M}_{v}}^{\tilde{L}^v}(\gamma_{v},\omega,f_{v,\tilde{L}^v,\omega}).$$
  
      Comme en 1.9, on peut formaliser les d\'efinitions ci-dessus et en particulier les rendre ind\'ependantes de tout choix de mesures en d\'efinissant $I_{\tilde{M}}^{\tilde{G}}(\boldsymbol{\gamma},{\bf f})$ pour $\boldsymbol{\gamma}\in D_{orb}(\tilde{M}(F_{V}),\omega)\otimes Mes(M(F_{V}))^*$ et ${\bf f}\in C_{c}^{\infty}(\tilde{G}(F_{V}))\otimes Mes(G(F_{V}))$, ou ${\bf f}\in I(\tilde{G}(F_{V}),\omega)\otimes Mes(G(F_{V}))$. Ces int\'egrales d\'ependent tout de m\^eme   de la mesure  fix\'ee sur ${\cal A}_{\tilde{M}}^{\tilde{G}}$.

Supposons $V=V_{1}\sqcup V_{2}$. Pour $i=1,2$, soient $\boldsymbol{\gamma}_{i}\in D_{orb}(\tilde{M}(F_{V_{i}}),\omega)\otimes Mes(M(F_{V_{i}}))^*$ et ${\bf f}_{i}\in I(\tilde{G}(F_{V_{i}}),\omega)\otimes Mes(G(F_{V_{i}}))$. Posons $\boldsymbol{\gamma}=\boldsymbol{\gamma}_{1}\otimes \boldsymbol{\gamma}_{2}$ et ${\bf f}={\bf f}_{1}\otimes {\bf f}_{2}$. On a la formule de scindage
$$(2) \qquad I_{\tilde{M}}^{\tilde{G}}(\boldsymbol{\gamma},{\bf f})=\sum_{\tilde{L}_{1},\tilde{L}_{2}\in {\cal L}(\tilde{M})}d_{\tilde{M}}^{\tilde{G}}(\tilde{L}_{1},\tilde{L}_{2})I_{\tilde{M}}^{\tilde{L}_{1}}(\boldsymbol{\gamma}_{1},(f_{1})_{\tilde{L}_{1},\omega})I_{\tilde{M}}^{\tilde{L}_{2}}(\boldsymbol{\gamma}_{2},(f_{2})_{\tilde{L}_{2},\omega}).$$

Conform\'ement aux r\'esultats de [V], on peut d\'efinir le terme $I_{\tilde{M}}^{\tilde{G}}(\boldsymbol{\gamma},{\bf f})$  dans le cas o\`u $\boldsymbol{\gamma}$ appartient \`a $D_{g\acute{e}om,\tilde{G}-\acute{e}qui}(\tilde{M}(F_{V}),\omega)\otimes Mes(M(F_{V}))^*$: on le d\'efinit par la formule (1). Si  $(G,\tilde{G},{\bf a})$ est quasi-d\'eploy\'e et \`a torsion int\'erieure, ou si l'on suppose v\'erifi\'ee l'hypoth\`ese (Hyp) de [V] 2.5 en toute place archim\'edienne de $V$,  on peut aussi le d\'efinir dans le  cas o\`u $\boldsymbol{\gamma}$ appartient \`a $D_{tr-orb}(\tilde{M}(F_{V}),\omega)\otimes Mes(M(F_{V}))^*$. Les propri\'et\'es ci-dessus s'\'etendent \`a tous les cas o\`u les termes sont d\'efinis. 

\bigskip

\subsection{Une propri\'et\'e de support}
Soient $\tilde{M}\in {\cal L}(\tilde{M}_{0})$ et $V$ un ensemble fini de places de $F$ contenant $V_{ram}$. 

\ass{Lemme}{Soit $\Xi\subset {\cal A}_{\tilde{M}}$ un ensemble compact et soit $f\in C_{ac,glob}^{\infty}(\tilde{G}(F_{V}))$. Alors il existe un sous-ensemble compact $\tilde{C}_{V}$ de $\tilde{M}(F_{V})$ tel que, pour tout $\gamma\in \tilde{M}(F_{V})$ v\'erifiant les deux conditions:

- $\tilde{H}_{\tilde{M}_{V}}(\gamma)\in \Xi$,

- $I_{\tilde{M}}^{\tilde{G}}(\gamma,\omega,f)\not=0$,

\noindent $\gamma$ soit conjugu\'e \`a un \'el\'ement de $\tilde{C}_{V}$ par un \'el\'ement de $M(F_{V})$.}

Preuve. On choisit une fonction $b\in C_{c}^{\infty}({\cal A}_{\tilde{G}})$ qui vaut $1$ sur la projection de $\Xi$ dans ${\cal A}_{\tilde{G}}$. On a alors l'\'egalit\'e $I_{\tilde{M}}^{\tilde{G}}(\gamma,\omega,f)=I_{\tilde{M}}^{\tilde{G}}(\gamma,\omega,f(b\circ \tilde{H}_{\tilde{G}_{V}}))$ pour tout $\gamma\in \tilde{M}(F_{V})$ tel que $\tilde{H}_{\tilde{M}_{V}}(\gamma)\in \Xi$. Cela nous permet de remplacer $f$ par $f(b\circ\tilde{H}_{\tilde{G}_{V}})$. En oubliant cela, on peut supposer $f$ \`a support compact. On utilise la d\'efinition donn\'ee en 1.9. Pour que $I_{\tilde{M}}^{\tilde{G}}(\gamma,\omega,f)$ soit non nul, il faut que $J_{\tilde{M}}^{\tilde{G}}(\gamma,\omega,f)$ soit non nul ou qu'il existe $\tilde{L}\in {\cal L}(\tilde{M})$ avec $\tilde{L}\not=\tilde{G}$ tel que $I_{\tilde{M}}^{\tilde{L}}(\gamma,\omega,\phi_{\tilde{L}}(f))$ soit non nul. Dans le premier cas, $\gamma$ est conjugu\'e par un \'el\'ement de $M(F_{V})$ \`a un \'el\'ement du support de $f$ et la conclusion s'ensuit. Dans le deuxi\`eme cas, le lemme 1.7 nous dit que $\phi_{\tilde{L}}(f)$ appartient \`a $I_{ac,glob}(\tilde{M}(F_{V}),\omega)$. Puisque $\tilde{L}\not=\tilde{G}$, on peut appliquer le lemme par r\'ecurrence, d'o\`u encore la conclusion. $\square$

\bigskip

\subsection{Le cas non ramifi\'e}

Soit $V$ un ensemble fini de places de $F$. Contrairement \`a l'habitude, on suppose ici $V\cap V_{ram}=\emptyset$. En particulier, les places dans $V$ sont finies. Soit $\tilde{M}\in {\cal L}(\tilde{M}_{0})$. On se d\'ebarrasse des espaces de mesures en fixant sur $G(F_{V})$ et $M(F_{V})$ les mesures canoniques. On d\'efinit une forme lin\'eaire $r_{\tilde{M}}^{\tilde{G}}(.,\tilde{K}_{V})$ sur $D_{g\acute{e}om}(\tilde{M}(F_{V}),\omega)$ par 
$$r_{\tilde{M}}^{\tilde{G}}(\boldsymbol{\gamma},\tilde{K}_{V})=J_{\tilde{M}}^{\tilde{G}}(\boldsymbol{\gamma},{\bf 1}_{\tilde{K}_{V}})$$
pour tout $\boldsymbol{\gamma}\in D_{g\acute{e}om}(\tilde{M}(F_{V}),\omega)$. Remarquons que, pour tout espace parabolique $\tilde{Q}=\tilde{L}U_{Q}\in {\cal L}(\tilde{M}_{0})$  et tout $v\in V$, on a l'\'egalit\'e $({\bf 1}_{\tilde{K}_{v}})_{\tilde{Q},\omega}={\bf 1}_{\tilde{K}_{v}^{\tilde{L}}}$.   Pour $\boldsymbol{\gamma}=\otimes_{v\in V}\boldsymbol{\gamma}_{v}$, la formule de descente 1.8(1) donne donc
$$r_{\tilde{M}}^{\tilde{G}}(\boldsymbol{\gamma},\tilde{K}_{V})=\sum_{\tilde{L}^V\in {\cal L}(\tilde{M}_{V})}d_{\tilde{M}_{V}}^{\tilde{G}}(\tilde{M},\tilde{L}^V)\prod_{v\in V}r_{\tilde{M}_{v}}^{\tilde{L}^{v}}(\boldsymbol{\gamma}_{v},\tilde{K}_{v}^{\tilde{L}^v}),$$
o\`u les derniers facteurs sont les termes locaux  d\'efinis en [II] 4.1.

\bigskip

\subsection{Int\'egrales orbitales pond\'er\'ees invariantes et syst\`emes de fonctions $B$}
On suppose $(G,\tilde{G},{\bf a})$ quasi-d\'eploy\'e et \`a torsion int\'erieure. Fixons un syst\`eme de  fonctions $B$ comme en 1.9. Soient $\tilde{M}\in {\cal L}(\tilde{M}_{0})$, $V$ un ensemble fini de places et $\gamma\in \tilde{M}(F_{V})$. De la m\^eme fa\c{c}on que dans le paragraphe  1.10, et modulo des choix de mesures de Haar, on d\'efinit l'int\'egrale orbitale pond\'er\'ee invariante $I_{\tilde{M}}^{\tilde{G}}(\gamma,B,f)$ pour $f\in C_{c}^{\infty}(\tilde{G}(F_{V}))$. Elle v\'erifie les m\^emes propri\'et\'es qu'en  1.10 et 1.11.

\ass{Lemme}{Supposons que $V$ contienne $V_{ram}$. Soit $\dot{\gamma}\in \tilde{M}(F)$, notons $\gamma$ son image naturelle dans $\tilde{M}(F_{V})$. 
 Alors on a l'\'egalit\'e
$$I_{\tilde{M}}^{\tilde{G}}(\gamma,B,f)=I_{\tilde{M}}^{\tilde{G}}(\gamma,f)$$
pour tout $f\in C_{c}^{\infty}(\tilde{G}(F_{V}))$.}

Preuve.  On v\'erifie que le proc\'ed\'e de limite utilis\'e pour d\'efinir les int\'egrales orbitales pond\'er\'ees s'\'etend aux int\'egrales invariantes. On a donc
$$I_{\tilde{M}}^{\tilde{G}}(\gamma,f)=lim_{a\to 1}\sum_{\tilde{L}\in {\cal L}(\tilde{M})}r_{\tilde{M}}^{\tilde{L}}(\gamma,a)I_{\tilde{L}}^{\tilde{G}}(a\gamma,f)$$
pour tout $f\in C_{c}^{\infty}(\tilde{G}(F_{V}))$, o\`u $a$ parcourt les \'el\'ements de $A_{\tilde{M}}(F_{V})$ en position g\'en\'erale. De m\^eme,
$$I_{\tilde{M}}^{\tilde{G}}(\gamma,B,f)=lim_{a\to 1}\sum_{\tilde{L}\in {\cal L}(\tilde{M})}r_{\tilde{M}}^{\tilde{L}}(\gamma,a,B)I_{\tilde{L}}^{\tilde{G}}(a\gamma,f).$$
Rappelons que, pour $a$ en position g\'en\'erale, l'\'el\'ement $a\gamma$ est $\tilde{G}$-\'equisingulier donc 
$I_{\tilde{L}}^{\tilde{G}}(a\gamma,B,f)=I_{\tilde{L}}^{\tilde{G}}(a\gamma,f)$  pour tout $\tilde{L}$.  Rappelons que les termes  $r_{\tilde{M}}^{\tilde{L}}(\gamma,a)$ et $r_{\tilde{M}}^{\tilde{L}}(\gamma,a,B)$ sont issus de $(\tilde{G},\tilde{M})$-familles $(r_{\tilde{P}}(\gamma,a;\Lambda))_{\tilde{P}\in {\cal P}(\tilde{M})}$ et $(r_{\tilde{P}}(\gamma,a,B;\Lambda))_{\tilde{P}\in {\cal P}(\tilde{M})}$. D\'efinissons une $(\tilde{G},\tilde{M})$-famille $(c_{\tilde{P}}(\gamma,a,B;\Lambda))_{\tilde{P}\in {\cal P}(\tilde{M})}$ par
$$c_{\tilde{P}}(\gamma,a,B;\Lambda)=r_{\tilde{P}}(\gamma,a,B;\Lambda)r_{\tilde{P}}(\gamma,a;\Lambda)^{-1}.$$
Pour tout $\tilde{L}\in {\cal L}(\tilde{M})$, on a la formule de d\'ecomposition
$$r_{\tilde{M}}^{\tilde{L}}(\gamma,a,B)=\sum_{\tilde{R}\in {\cal L}(\tilde{M}),\tilde{R}\subset \tilde{L}}c_{\tilde{M}}^{\tilde{R}}(\gamma,a,B)r_{\tilde{R}}^{\tilde{L}}(\gamma,a).$$
D'o\`u
$$I_{\tilde{M}}^{\tilde{G}}(\gamma,B,f)=lim_{a\to 1}\sum_{\tilde{R}\in {\cal L}(\tilde{M})}c_{\tilde{M}}^{\tilde{R}}(\gamma,a,B)\sum_{\tilde{L}\in {\cal L}(\tilde{R})}r_{\tilde{R}}^{\tilde{L}}(\gamma,a)I_{\tilde{L}}^{\tilde{G}}(a\gamma,f).$$
En utilisant la formule de descente 1.4(1)  et [II] 1.7(12), on voit que pour tout $\tilde{R}$, la somme int\'erieure a une limite quand $a$ tend vers $1$. Cette limite est $I_{\tilde{R}}^{\tilde{G}}(\boldsymbol{\gamma}^{\tilde{R}},f)$, o\`u $\boldsymbol{\gamma}^{\tilde{R}}$ est la distribution induite \`a $\tilde{R}(F_{V})$ de l'int\'egrale orbitale dans $\tilde{M}(F_{V})$ associ\'ee \`a $\gamma$.  Pour prouver le lemme, il suffit de prouver la relation
$$ lim_{a\to 1}c_{\tilde{M}}^{\tilde{R}}(\gamma,a,B)=\left\lbrace\begin{array}{cc}1,&\text{ si }\tilde{R}=\tilde{M},\\ 0,& \text{ si }\tilde{R}\not=\tilde{M}\\ \end{array}\right.$$
pour tout $\tilde{R}\in {\cal L}(\tilde{M})$. Un tel $\tilde{R}$ \'etant fix\'e, on peut remplacer l'espace ambiant $\tilde{G}$ par $\tilde{R}$. Cela nous ram\`ene \`a prouver la relation
$$(1) \qquad  lim_{a\to 1}c_{\tilde{M}}^{\tilde{G}}(\gamma,a,B)=\left\lbrace\begin{array}{cc}1,&\text{ si }\tilde{G}=\tilde{M},\\ 0,& \text{ si }\tilde{G}\not=\tilde{M}\\ \end{array}\right.$$
Le cas $\tilde{G}=\tilde{M}$ est \'evident.  On  suppose d\'esormais $\tilde{M}\not=\tilde{G}$. 

Rappelons la construction de nos $(\tilde{G},\tilde{M})$-familles. Ecrivons $\dot{\gamma}=u\eta$, o\`u $\eta\in \tilde{M}(F)$ est semi-simple et $u\in M_{\eta}(F)$ est unipotent.  On introduit les ensembles $\Sigma^{G_{\eta}}(Z(M_{\eta})^0)$ et $\Sigma(A_{\tilde{M}})$ de racines de $Z(M_{\eta})^0$ dans $\mathfrak{g}_{\eta}$, resp. de $A_{\tilde{M}}$ dans $\mathfrak{g}$. On a une application de restriction
$$\begin{array}{ccc}\Sigma^{G_{\eta}}(Z(M_{\eta})^0)&\to&\Sigma(A_{\tilde{M}})\\ \alpha&\mapsto& \alpha_{\tilde{M}}\\ \end{array}$$
Posons ${\cal Z}=X_{*}(Z(M)^0)\otimes _{{\mathbb Z}}{\mathbb R}$. On a ${\cal A}_{\tilde{M}}\subset {\cal Z}$. Puisque ${\cal Z}$ est muni d'une forme quadratique d\'efinie positive (cf. 1.3), on a aussi une inclusion d'espaces duaux ${\cal A}_{\tilde{M}}^*\subset {\cal Z}^*$.    Fixons $v\in V$ et $\alpha\in \Sigma^{G_{\eta}}(Z(M_{\eta})^0)$. En [II] 1.4, on a d\'efini un \'el\'ement de ${\cal Z}$, not\'e alors $\rho(\alpha,u)$. Sa d\'efinition d\'epend a priori de la place $v$, notons-le plut\^ot $\rho_{v}(\alpha,u)$.  Soient $\tilde{P}\in {\cal P}(\tilde{M})$,  $a=(a_{v})_{v\in V}\in A_{\tilde{M}}(F_{V})$ et $\Lambda\in i{\cal A}_{\tilde{M}}^*$. Il r\'esulte des d\'efinitions que l'on a l'\'egalit\'e
$$(2) \qquad r_{\tilde{P}}(\gamma,a;\Lambda)=\prod_{v\in V}\prod_{\alpha\in \Sigma^{G_{\eta}}(Z(M_{\eta})^0); \alpha_{\tilde{M}}>_{P}0}\vert \alpha(a_{v})-\alpha(a_{v})^{-1}\vert _{F_{v}}^{<\Lambda,\rho_{v}(\alpha,u)>/2},$$
o\`u le symbole $>_{P}$ d\'esigne la positivit\'e relative \`a $P$. 

 Fixons une paire de Borel $(B,T)$ de $G$ conserv\'ee par $ad_{\eta}$ et telle que $M$ soit standard pour cette paire. Introduisons l'ensemble $\Sigma^{G_{\eta}}(T)$ des racines de $T$ dans $\mathfrak{g}_{\eta}$. On a introduit en [II] 1.8 l'ensemble $\Sigma^{G_{\eta}}(T,B_{\eta})$ form\'e des $B_{\eta}(\alpha)^{-1}\alpha$ pour $\alpha\in \Sigma^{G_{\eta}}(T,B_{\eta})$ (on consid\`ere ces \'el\'ements comme des formes lin\'eaires sur $\mathfrak{t}$). On note $\Sigma^{G_{\eta}}(Z(M_{\eta})^0,B)$ l'ensemble des restrictions \`a $\mathfrak{z}(M_{\eta})$ d'\'el\'ements de $\Sigma^{G_{\eta}}(T,B_{\eta})$.    Fixons $v\in V$ et $\alpha'\in \Sigma^{G_{\eta}}(Z(M_{\eta})^0,B)$. En [II] 1.4, on a d\'efini un \'el\'ement de ${\cal Z}$, not\'e alors $\rho(\alpha',u,B)$, qu'il convient de noter plut\^ot $\rho_{v}(\alpha',u,B)$. Soient $\tilde{P}\in {\cal P}(\tilde{M})$,  $a=(a_{v})_{v\in V}\in A_{\tilde{M}}(F_{V})$ et $\Lambda\in i{\cal A}_{\tilde{M}}^*$. 
Un \'el\'ement  $\alpha'\in \Sigma^{G_{\eta}}(Z(M_{\eta})^0,B)$  se restreint \`a $\mathfrak{a}_{\tilde{M}}$ en un   \'el\'ement  $\alpha'_{\tilde{M}}=q\beta$, o\`u $q\in {\mathbb Q}_{>0}$ et $\beta\in \Sigma(A_{\tilde{M}})$.  On dit que $\alpha'_{\tilde{M}}>_{P}0$ si et seulement si $\beta>_{P}0$. D'autre part, l'\'el\'ement $a$ \'etant suppos\'e proche de $1$, on peut \'ecrire $a_{v}=exp(H_{v})$ pour tout $v\in V$, o\`u $H_{v}\in \mathfrak{a}_{\tilde{M}}(F_{v})$ est proche de $0$. On pose $\alpha'(a_{v})=exp(q\beta(H_{v}))$. Il r\'esulte  alors des d\'efinitions que l'on a l'\'egalit\'e
$$(3) \qquad r_{\tilde{P}}(\gamma,a,B;\Lambda)=\prod_{v\in V}\prod_{\alpha'\in \Sigma^{G_{\eta}}(Z(M_{\eta})^0,B); \alpha'_{\tilde{M}}>_{P}0}\vert \alpha'(a_{v})-\alpha'(a_{v})^{-1}\vert _{F_{v}}^{<\Lambda,\rho_{v}(\alpha',u,B)>/2}. $$

Notons $\Sigma_{ind}(A_{\tilde{M}})$  l'ensemble des \'el\'ements indivisibles de $\Sigma(A_{\tilde{M}})$. Un \'el\'ement $\alpha$ intervenant dans (2) se restreint en un \'el\'ement $\alpha_{\tilde{M}}$ qui est un multiple entier positif d'un unique \'el\'ement $\beta\in \Sigma_{ind}(A_{\tilde{M}})$. On regroupe les $\alpha$ selon cet \'el\'ement $\beta$ et on obtient une d\'ecomposition en produit
$$r_{\tilde{P}}(\gamma,a;\Lambda)=\prod_{\beta\in \Sigma_{ind}(A_{\tilde{M}}), \beta_{\tilde{M}}>_{P}0}r_{\beta}(\gamma,a;\Lambda).$$
Un \'el\'ement $\alpha'$ intervenant dans (3) se restreint en un \'el\'ement $\alpha'_{\tilde{M}}$ qui est un multiple rationnel positif d'un unique \'el\'ement $\beta\in \Sigma_{ind}(A_{\tilde{M}})$. On regroupe les $\alpha'$ selon cet \'el\'ement $\beta$ et on obtient une d\'ecomposition en produit
$$r_{\tilde{P}}(\gamma,a,B;\Lambda)=\prod_{\beta\in \Sigma_{ind}(A_{\tilde{M}}), \beta_{\tilde{M}}>_{P}0}r_{\beta}(\gamma,a,B;\Lambda).$$
Fixons  $\beta\in \Sigma_{ind}(A_{\tilde{M}})$. On peut introduire une coracine $\check{\beta}\in {\cal A}_{\tilde{M}}$, normalis\'ee par la condition $<\beta,\check{\beta}>=2$. Il resulte des constructions de [II] 1.4  que, si $\alpha$ est un \'el\'ement de  $\Sigma^{G_{\eta}}(Z(M_{\eta})^0)$ tel que $\alpha_{\tilde{M}}$ est un multiple entier de $\beta$, alors, pour tout $v\in V$, la projection orthogonale de $\rho_{v}(\alpha,u)$ sur ${\cal A}_{\tilde{M}}$ est colin\'eaire \`a $\check{\beta}$. Il en r\'esulte que
$$<\Lambda, \rho_{v}(\alpha,u)>/2=<\beta,\rho_{v}(\alpha,u)><\Lambda,\check{\beta}>/2.$$
De m\^eme, il r\'esulte des constructions de [II] 1.8 que, si $\alpha'$ est un \'el\'ement de  $\Sigma^{G_{\eta}}(Z(M_{\eta})^0,B)$ tel que $\alpha'_{\tilde{M}}$ est un multiple rationnel  de $\beta$, alors, pour tout $v\in V$, la projection orthogonale de $\rho_{v}(\alpha',u,B)$ sur ${\cal A}_{\tilde{M}}$ est colin\'eaire \`a $\check{\beta}$. Il en r\'esulte que
$$<\Lambda, \rho_{v}(\alpha',u,B)>/2=<\beta,\rho_{v}(\alpha',u,B)><\Lambda,\check{\beta}>/2.$$
D\'efinissons une fonction $c_{\beta}(\gamma,a,B;x)$ d'une variable r\'eelle $x$ par l'\'egalit\'e
$$(4)  \qquad c_{\beta}(\gamma,a,B;x)=\prod_{v\in V}\left(\prod_{\alpha'\in \Sigma^{G_{\eta}}(Z(M_{\eta})^0,B); \alpha'_{\tilde{M}}\in {\mathbb Q}_{>0}\beta}\vert \alpha'(a_{v})-\alpha'(a_{v})^{-1}\vert _{F_{v}}^{ix<\beta,\rho_{v}(\alpha',u,B)>/4}\right)$$
$$\left(\prod_{\alpha\in \Sigma^{G_{\eta}}(Z(M_{\eta})^0); \alpha_{\tilde{M}}\in {\mathbb Z}_{>0}\beta}\vert \alpha(a_{v})-\alpha(a_{v})^{-1}\vert _{F_{v}}^{-ix<\beta,\rho_{v}(\alpha,u)>/4}\right).$$
On obtient alors l'\'egalit\'e
$$c_{\beta}(\gamma,a,B;-i<\Lambda,\check{\beta}>)=r_{\beta}(\gamma,a,B;\Lambda)r_{\beta}(\gamma,a;\Lambda)^{-1}.$$
D'o\`u l'\'egalit\'e 
$$c_{\tilde{P}}(\gamma,a,B;\Lambda)=\prod_{\beta\in \Sigma_{ind}(A_{\tilde{M}}), \beta>_{P}0}c_{\beta}(\gamma,a,B;-i<\Lambda,\check{\beta}>).$$
 Cela montre  que la famille $(c_{\tilde{P}}(\gamma,a,B;\Lambda))_{\tilde{P}\in {\cal P}(\tilde{M})}$  est de la forme particuli\`ere \'etudi\'ee par Arthur en [A6].  
Le lemme 7.1 de cette r\'ef\'erence  calcule explicitement $c_{\tilde{M}}^{\tilde{G}}(\gamma,a,B)$:  c'est une combinaison lin\'eaire de produits des  d\'eriv\'ees $c'_{\beta}(\gamma,a,B;x)$ des fonctions $c_{\beta}(\gamma,a,B;x)$ \'evalu\'ees en $x=0$. Pour prouver (1), il suffit de fixer $\beta\in \Sigma_{ind}(A_{\tilde{M}})$ et de prouver la relation
$$lim_{a\to 1}c'_{\beta}(\gamma,a,B;0)=0.$$
Fixons $\beta$ et consid\'erons la formule (4). Notons $\Sigma_{ind}^{G_{\eta}}(Z(M_{\eta})^0)$ l'ensemble des \'el\'ements indivisibles de $\Sigma^{G_{\eta}}(Z(M_{\eta})^0)$. Pour $\alpha_{ind}\in \Sigma_{ind}^{G_{\eta}}(Z(M_{\eta})^0)$, d\'efinissons une fonction $c_{\alpha}(\gamma,a,B;x)$  par une formule analogue \`a (4), o\`u on se restreint aux $\alpha$ qui sont multiples entiers positifs de $\alpha_{ind}$ et aux $\alpha'$ qui sont multiples rationnels positifs de $\alpha_{ind}$. 
 On obtient une d\'ecomposition
$$c_{\beta}(\gamma,a,B;x)=\prod_{\alpha_{ind}\in \Sigma_{ind}^{G_{\eta}}(Z(M_{\eta})^0), \alpha_{ind,\tilde{M}}\in {\mathbb Z}_{>0}\beta}c_{\alpha_{ind}}(\gamma,a,B;x).$$
Il nous suffit de fixer $\alpha_{ind}\in \Sigma_{ind}^{G_{\eta}}(Z(M_{\eta})^0)$ et de prouver la relation 
$$lim_{a\to 1}c'_{\alpha_{ind}}(\gamma,a,B;0)=0.$$
Fixons donc un \'el\'ement de $\Sigma_{ind}^{G_{\eta}}(Z(M_{\eta})^0)$. Pour la commodit\'e de l'\'ecriture, notons-le simplement $\alpha$. On a vu en [II] 1.8 que l'ensemble $\Sigma_{ind}^{G_{\eta}}(Z(M_{\eta})^0,B)$  poss\'edait beaucoup des propri\'et\'es  des  syst\`emes de racines. En particulier, l'ensemble des $\alpha'\in \Sigma_{ind}^{G_{\eta}}(Z(M_{\eta})^0,B)$ qui sont des multiples rationnels positifs de $\alpha$ poss\`ede un unique \'el\'ement minimal, notons-le simplement $\alpha'$. Les autres \'el\'ements de cet ensemble sont des multiples entiers positifs de $\alpha'$. Posons alors
$$(5) \qquad X_{\alpha}(\gamma,a)=\sum_{v\in V}\sum_{n\geq1}log(\vert \alpha(a_{v})^n-\alpha(a_{v})^{-n}\vert _{F_{v}})\rho_{v}(n\alpha,u),$$
$$(6) \qquad X_{\alpha'}(\gamma,a,B)=\sum_{v\in V}\sum_{n\geq1}log(\vert \alpha'(a_{v})^n-\alpha'(a_{v})^{-n}\vert _{F_{v}})\rho_{v}(n\alpha',u,B),$$
o\`u, par convention, les termes $\rho_{v}(n\alpha,u)$ et $\rho_{v}(n\alpha',u,B)$ sont nuls si $n\alpha\not\in \Sigma^{G_{\eta}}(Z(M_{\eta})^0)$, resp.  $n\alpha'\not\in \Sigma^{G_{\eta}}(Z(M_{\eta})^0,B)$.  On calcule
$$c'_{\alpha}(\gamma,a,B;0)=\frac{i}{4}<\beta,X_{\alpha'}(\gamma,a,B)-X_{\alpha}(\gamma,a)>,$$
o\`u $\beta$ est l'unique \'el\'ement de $\Sigma_{ind}(A_{\tilde{M}})$ tel que $\alpha_{\tilde{M}}$ soit un multiple entier positif de $\beta$. 
Il nous suffit de prouver l'\'egalit\'e
$$(7)\qquad lim_{a\to 1}(X_{\alpha}(\gamma,a)-X_{\alpha'}(\gamma,a,B))=0.$$

On a besoin de deux r\'esultats pr\'eliminaires. D'abord

(8) pour tout $n\geq1$, les termes $\rho_{v}(n\alpha,u)$ et $\rho_{v}(n\alpha',u,B)$ sont ind\'ependants de $v\in V$. 

Preuve. C'est clair si $n\alpha\not\in \Sigma^{G_{\eta}}(Z(M_{\eta})^0)$, resp. $n\alpha'\not\in \Sigma^{G_{\eta}}(Z(M_{\eta})^0,B)$. Soit $n\geq2$, supposons que $n\alpha\in \Sigma^{G_{\eta}}(Z(M_{\eta})^0)$.  Soit $v\in V$.   Notre terme $\rho_{v}(n\alpha,u)$ d\'epend du groupe ambiant $G_{\eta}$, notons-le ici $\rho_{v}^{G_{\eta}}(n\alpha,u)$. On a introduit en [II] 1.4 un sous-groupe $G_{\eta,n\alpha}$ de $G_{\eta}$. On a l'\'egalit\'e $\rho_{v}^{G_{\eta}}(n\alpha,u)=\rho_{v}^{G_{\eta,n\alpha}}(n\alpha,u)$. L'hypoth\`ese $n\geq2$ implique que $dim(G_{\eta,n\alpha,SC})<dim(G_{\eta})$. En raisonnant par r\'ecurrence sur cette dimension, on peut supposer que $\rho_{v}^{G_{\eta,n\alpha}}(n\alpha,u)$ est ind\'ependant de $v$. Donc $\rho_{v}^{G_{\eta}}(n\alpha,u)$ aussi. Le m\^eme raisonnement vaut pour $\rho_{v}^{G_{\eta}}(n\alpha',u,B)$, en utilisant le groupe $G_{\eta,n\alpha'}$ de [II] 1.8 (ce n'est plus un sous-groupe de $G_{\eta}$ mais peu importe). L'assertion \'etant d\'emontr\'ee pour $n\geq2$, il nous suffit pour conclure de prouver que les sommes
$$\sum_{n\geq1}\rho_{v}(n\alpha,u)\, \text{ et }\sum_{n\geq1}\rho_{v}(n\alpha',u,B)$$
sont ind\'ependantes de $v$. D'apr\`es [II] 1.8(6), ces deux sommes sont \'egales. D'apr\`es la d\'efinition de [II] 1.4, elles valent  le terme $\rho_{v}^{Art}(\alpha,u)\check{\alpha}$ d\'efini  primitivement par Arthur  en [A5] paragraphe 3. Rappelons bri\`evement sa d\'efinition. Fixons une extension finie $F'_{v}$ de $F_{v}$ sur laquelle $G_{\eta}$ est d\'eploy\'ee. A $\alpha$ est associ\'e un Levi $M_{\eta,\alpha}$ de $G_{\eta}$ qui contient strictement $M_{\eta}$ et qui est minimal parmi les Levi v\'erifiant cette propri\'et\'e. Le terme $\rho_{v}^{Art}(\alpha,u)\check{\alpha}$  relatif au groupe ambiant $G_{\eta}$ est \'egal \`a celui relatif \`a $M_{\eta,\alpha}$. On ne perd rien \`a supposer $G_{\eta}=M_{\eta,\alpha}$.  Fixons $P=M_{\eta}U_{P}\in {\cal P}(M_{\eta})$, notons $\bar{P}=M_{\eta}U_{\bar{P}}$ le parabolique oppos\'e et notons ${\cal U}$ l'orbite de $u$ pour la conjugaison par $M_{\eta}$. Fixons un poids $\omega$ de $A_{M_{\eta}}$ qui est dominant pour $P$, fixons une repr\'esentation alg\'ebrique irr\'eductible $\Lambda_{\omega}$ de $G_{\eta}$ dans un espace $V_{\omega}$, de plus haut poids $\omega$. Fixons un vecteur extr\'emal $\phi_{\omega} \in V_{\omega}$, de poids $\omega$. Pour $\underline{a}\in A_{M_{\eta}}$ en position g\'en\'erale et pour $\pi=n\nu \in {\cal U}U_{\bar{P}}$, avec $n\in {\cal U}$ et $\nu\in U_{\bar{P}}$, introduisons l'\'el\'ement $\nu(\underline{a},\pi)\in U_{\bar{P}}$ tel que $\underline{a}\pi=\nu(\underline{a},\pi)^{-1}an\nu(\underline{a},\pi)$. On pose $\varphi(\underline{a},\pi)=\Lambda_{\omega}(\nu(\underline{a},\pi)^{-1})\phi_{\omega}$. Arthur montre en [A4] page 238 que $\varphi$ est une application rationnelle sur $A_{M_{\eta}}\times {\cal U}U_{\bar{P}}$, \`a valeurs dans $V_{\omega}$, et qu'il existe un unique entier $k\in {\mathbb Z}$ tel que la fonction
$$(\underline{a},\pi)\mapsto (\alpha(\underline{a})-\alpha(\underline{a})^{-1})^k\varphi(\underline{a},\pi)$$
soit r\'eguli\`ere et non nulle sur $\{1\}\times {\cal U}U_{\bar{P}}$. Une coracine $\check{\alpha}$ \'etant fix\'ee, $\rho_{v}^{Art}(\alpha,u)$ est l'unique r\'eel tel que $k=\omega(\check{\alpha})\rho_{v}^{Art}(\alpha,u)$. 
Fixons une extension finie $F'$ de $F$ telle que $G_{\eta}$ soit d\'eploy\'e sur $F'$.   La construction ci-dessus \'etant de nature alg\'ebrique, on peut la refaire sur le corps de base $F'$. On obtient un nombre r\'eel $\rho^{Art}(\alpha,u)$ ind\'ependant de la place $v$ et il est clair que $\rho_{v}^{Art}(\alpha,u)=\rho^{Art}(\alpha,u)$ pour tout $v$. Cela prouve (8).

Supprimons d\'esormais les indices $v$ des termes $\rho_{v}(n\alpha,u)$ et $\rho_{v}(n\alpha',u,B)$. Comme on l'a dit dans la preuve ci-dessus, il r\'esulte de [II] 1.8(6) que
$$(9)\qquad \sum_{n\geq1}\rho(n\alpha,u)=\sum_{n\geq1}\rho(n\alpha',u,B).$$

Soit $\beta\in \Sigma^{G_{\eta}}(T)$ telle que $\alpha'$ soit la restriction de $B_{\eta}(\beta)^{-1}\beta$ \`a $Z(M_{\eta})^0$. Posons $b=B_{\eta}(\beta)$ et soit $m\geq1$ tel que la restriction de $\beta$ soit $m\alpha$. Alors $\alpha'=\frac{m}{b}\alpha$. L'\'el\'ement $a$ \'etant suppos\'e proche de $1$, on \'ecrit $a_{v}=exp(H_{v})$ pour tout $v\in V$, o\`u $H_{v}\in \mathfrak{a}_{\tilde{M}}(F_{v})$ est proche de $0$. Pour $n\geq1$, on a $\alpha'(a_{v})^n=exp(\frac{nm}{b}<\alpha,H_{v}>)$. 
Donc  l'expression $log(\vert \alpha'(a_{v})^n-\alpha'(a_{v})^{-n}\vert _{F_{v}})-log(\vert \alpha(a_{v})-\alpha(a_{v})^{-1}\vert _{F_{v}})-log (\vert \frac{nm}{b}\vert _{F_{v}})$ tend vers $0$ quand $a$ tend vers $1$.    Posons
$$Y_{\alpha'}(\gamma,a,B)=\left(\sum_{v\in V}log(\vert \alpha(a_{v})-\alpha(a_{v})^{-1}\vert _{F_{v}})\sum_{n\geq1}\rho(n\alpha',\gamma,B)\right)$$
$$+\left(\sum_{n\geq1}\rho(n\alpha',\gamma,B)\sum_{v\in V}log(\vert \frac{nm}{b}\vert _{F_{v}}) \right).$$
En se reportant \`a l'expression (6), on voit que
 $X_{\alpha'}(\gamma,a,B)-Y_{\alpha'}(\gamma,a,B)$ tend vers $0$ quand $a$ tend vers $1$. 
 On a suppos\'e que $V$ contenait $V_{ram}$. Donc $b$ est une unit\'e en toute place $v\in Val(F)-V$. Le m\^eme calcul que dans la preuve du lemme [II] 1.9 montre que les seuls premiers pouvant diviser les nombres $n$ et $m$ intervenant ci-dessus  sont $2$, $3$ et $5$. Ces nombres sont donc eux-aussi des unit\'es hors de $V$. La formule du produit entra\^{\i}ne alors 
 $$\sum_{v\in V}log(\vert \frac{nm}{b}\vert _{F_{v}}) =0.$$
La d\'efinition de $Y_{\alpha'}(\gamma,a,B)$ se simplifie en
$$Y_{\alpha'}(\gamma,a,B)=\sum_{v\in V}log(\vert \alpha(a_{v})-\alpha(a_{v})^{-1}\vert _{F_{v}})\sum_{n\geq1}\rho(n\alpha',\gamma,B).$$
Un calcul analogue vaut pour $X_{\alpha}(\gamma,B)$. Si on pose
$$Y_{\alpha}(\gamma,a)=\sum_{v\in V}log(\vert \alpha(a_{v})-\alpha(a_{v})^{-1}\vert _{F_{v}})\sum_{n\geq1}\rho(n\alpha,\gamma),$$
on a $lim_{a\to 1}(X_{\alpha}(\gamma,a)-Y_{\alpha}(\gamma,a))=0$. Mais  (9) entra\^{\i}ne que $Y_{\alpha}(\gamma,a)=Y_{\alpha'}(\gamma,a,B)$. On en d\'eduit la relation (7), ce qui ach\`eve la preuve. $\square$

  .

 \bigskip

\subsection{Variante avec caract\`ere central}
On suppose $(G,\tilde{G},{\bf a})$ quasi-d\'eploy\'e et \`a torsion int\'erieure. On suppose donn\'ee une extension
$$1\to C_{1}\to G_{1}\to G\to 1$$
o\`u $C_{1}$ est un tore induit central, et une extension compatible $\tilde{G}_{1}\to \tilde{G}$, encore \`a torsion int\'erieure. On fixe un caract\`ere $\lambda_{1}$ de $C_{1}({\mathbb A}_{F})$, automorphe c'est-\`a-dire trivial sur $C_{1}(F)$.   Notons $V_{1,ram}$ (ou plus pr\'ecis\'ement $V_{ram}(\tilde{G}_{1},\lambda_{1})$) le plus petit ensemble de places de $F$ contenant $V_{ram}(\tilde{G}_{1})$ et tel que $\lambda_{1,v}$ soit non ramifi\'e pour $v\not\in V_{1,ram}$.

 On doit fixer pour tout $v\not\in V_{1,ram}$ un sous-espace hypersp\'ecial $\tilde{K}_{1,v}$ de $\tilde{G}_{1}(F_{v})$, soumis aux conditions de 1.1. On suppose que $\tilde{K}_{1,v}$ se projette sur $\tilde{K}_{v}$. Le groupe $K_{C_{1},v}=C_{1}(F_{v})\cap K_{1,v}$ est le plus grand sous-groupe compact de $C_{1}(F_{v})$.

Pour toute place $v$, on a d\'efini en [I] 2.4 l'espace $C_{c,\lambda_{1}}^{\infty}(\tilde{G}_{1}(F_{v}))$. Pour $v\not\in V_{1,ram}$, on note ${\bf 1}_{\tilde{K}_{1,v},\lambda_{1}}$ l'unique \'el\'ement de cet espace \`a support dans $C_{1}(F_{v})\tilde{K}_{1,v}$ qui vaut $1$ sur $\tilde{K}_{1,v}$. On note $C_{c,\lambda_{1}}^{\infty}(\tilde{G}_{1}({\mathbb A}_{F}))$ le produit tensoriel restreint des $C_{c,\lambda_{1}}^{\infty}(\tilde{G}_{1}(F_{v}))$ relativement \`a ces \'el\'ements ${\bf 1}_{\tilde{K}_{1,v},\lambda_{1}}$. Pour un ensemble fini $V$ de places de $F$, on d\'efinit
$$C_{c,\lambda_{1}}^{\infty}(\tilde{G}_{1}(F_{V}))=\otimes_{v\in V}C_{c,\lambda_{1}}^{\infty}(\tilde{G}_{1}(F_{v})).$$
Dualement, on d\'efinit de m\^eme l'espace $D_{g\acute{e}om,\lambda_{1}}(\tilde{G}_{1}(F_{V}))$ et ses sous-espaces $D_{orb,\lambda_{1}}(\tilde{G}_{1}(F_{V}))$ etc... Les constructions des paragraphes pr\'ec\'edents s'\'etendent \`a cette situation. En particulier, soit $\tilde{M}\in {\cal L}(\tilde{M}_{0})$. Notons $\tilde{M}_{1}$ son image r\'eciproque dans $\tilde{G}_{1}$. On a l'\'egalit\'e ${\cal A}_{\tilde{M}_{1}}^{\tilde{G}_{1}}={\cal A}_{\tilde{M}}^{\tilde{G}}$ et on munit le premier espace de la mesure pour laquelle cette \'egalit\'e pr\'eserve les mesures. Pour $\boldsymbol{\gamma}\in D_{orb,\lambda_{1}}(\tilde{M}_{1}(F_{V}))\otimes Mes(M(F_{V}))^*$ et ${\bf f}\in C_{c,\lambda_{1}}^{\infty}(\tilde{G}_{1}(F_{V}))\otimes Mes(G(F_{V}))$, on d\'efinit l'int\'egrale orbitale pond\'er\'ee $J_{\tilde{M}_{1},\lambda_{1}}^{\tilde{G}_{1}}(\boldsymbol{\gamma},{\bf f})$ et son avatar invariant $I_{\tilde{M}_{1},\lambda_{1}}^{\tilde{G}_{1}}(\boldsymbol{\gamma},{\bf f})$.

Supposons que $V\cap V_{ram}=\emptyset$. On d\'efinit une forme lin\'eaire $r_{\tilde{M}_{1},\lambda_{1}}^{\tilde{G}_{1}}(.,\tilde{K}_{1,V})$ sur

\noindent $D_{g\acute{e}om,\lambda_{1}}(\tilde{M}_{1}(F_{V}))$ par
$$r_{\tilde{M},\lambda_{1}}^{\tilde{G}}(\delta,\tilde{K}_{1,V})=J_{\tilde{M}}^{\tilde{G}}(\delta,{\bf 1}_{\tilde{K}_{V},\lambda_{1}}).$$
Elle v\'erifie une formule de descente analogue \`a celle de 1.13.

 Consid\'erons  d'autres donn\'ees
$$1\to C_{2}\to G_{2}\to G\to 1,\,\,\tilde{G}_{2}\to \tilde{G}, \, \lambda_{2}$$
 et des sous-espaces hypersp\'eciaux $\tilde{K}_{2,v}$ de $\tilde{G}_{2}(F_{v})$ pour $v\not\in V_{2,ram}$, v\'erifiant les m\^emes hypoth\`eses. On note $G_{12}$ le produit fibr\'e de $G_{1}$ et $G_{2}$ au-dessus de $G$ et $\tilde{G}_{12}$ le produit fibr\'e de $\tilde{G}_{1}$ et $\tilde{G}_{2}$ au-dessus de $\tilde{G}$. On suppose donn\'e un caract\`ere automorphe $\lambda_{12}$ de $G_{12}({\mathbb A}_{F})$ tel que la restriction de $\lambda_{12}$ \`a $C_{1}({\mathbb A}_{F})\times C_{2}({\mathbb A}_{F})$ soit $\lambda_{1}\times \lambda_{2}^{-1}$.  Notons $V_{12,ram}$ le plus petit ensemble de places de $F$  contenant $V_{1,ram}$ et $V_{2,ram}$ et tel que $\lambda_{12}$ soit non ramifi\'e hors de $V_{12,ram}$.
 
  Pour un ensemble fini $V$ de places de $F$, notons $\lambda_{12,V}$ la restriction de $\lambda_{12}$ \`a $G_{12}(F_{V})$. Fixons une fonction $\tilde{\lambda}_{12,V}$ sur $\tilde{G}_{12}(F_{V})$ telle que
   $$(1) \qquad \tilde{\lambda}_{12,V}(x_{V}\gamma_{V})=\lambda_{12,V}(x_{V})\tilde{\lambda}_{12,V}(\gamma_{V})$$
pour $x_{V}\in G_{12}(F_{V})$ et $\gamma_{V}\in \tilde{G}_{12}(F_{V})$. On d\'efinit un isomorphisme
$$\begin{array}{ccc}C_{c,\lambda_{1}}^{\infty}(\tilde{G}_{1}(F_{V}))&\to& C_{c,\lambda_{2}}^{\infty}(\tilde{G}_{2}(F_{V}))\\ f_{1}&\mapsto&f_{2}\\ \end{array}$$
par la formule
$$f_{2}(\gamma_{2})=\tilde{\lambda}_{12}(\gamma_{1},\gamma_{2})f_{1}(\gamma_{1})$$
pour tous $(\gamma_{1},\gamma_{2})\in \tilde{G}_{12}(F_{V})$. Cet isomorphisme se dualise en un  isomorphisme
$$D_{g\acute{e}om,\lambda_{1}}(\tilde{G}_{1}(F_{V}))\simeq D_{g\acute{e}om,\lambda_{2}}(\tilde{G}_{2}(F_{V})).$$
 On v\'erifie que les int\'egrales orbitales pond\'er\'ees et leurs avatars invariants se recollent selon ces isomorphismes. C'est-\`a-dire, soit $\tilde{M}\in {\cal L}(\tilde{M}_{0})$ et soient, pour $i=1,2$, $\boldsymbol{\gamma}_{i}\in D_{orb,\lambda_{i}}(\tilde{M}_{i}(F_{V}))\otimes Mes(M(F_{V}))^*$ et ${\bf f}_{i}\in C_{c,\lambda_{i}}^{\infty}(\tilde{G}_{i}(F_{V}))\otimes Mes(G(F_{V}))$. Supposons que, par les isomorphismes pr\'ec\'edents, $\boldsymbol{\gamma}_{1}$ et $\boldsymbol{\gamma}_{2}$ se correspondent, ainsi que ${\bf f}_{1}$ et ${\bf f}_{2}$. Alors on a l'\'egalit\'e
 $$J_{\tilde{M}_{1},\lambda_{1}}^{\tilde{G}_{1}}(\boldsymbol{\gamma}_{1},{\bf f}_{1})=J_{\tilde{M}_{2},\lambda_{2}}^{\tilde{G}_{2}}(\boldsymbol{\gamma}_{2},{\bf f}_{2}).$$
 
 Si $V$ contient $V_{12,ram}$, il y a une fonction $\tilde{\lambda}_{12,V}$ canonique construite de la fa\c{c}on suivante.  Pour $v\not\in V_{12,ram}$, on d\'efinit une fonction $\tilde{\lambda}_{12,v}$ sur $\tilde{G}_{12}(F_{v})$ par les conditions:

- $\tilde{\lambda}_{12,v}$ vaut $1$ sur $\tilde{G}_{12}(F_{v})\cap (\tilde{K}_{1,v}\times \tilde{K}_{2,v})$;

- pour $x\in G_{12}(F_{v})$ et $\gamma\in \tilde{G}_{12}(F_{v})$, $\tilde{\lambda}_{12,v}(x\gamma)=\lambda_{12,v}(x)\tilde{\lambda}_{12,v}(\gamma)$.

On d\'efinit une fonction $\tilde{\lambda}_{12}$ sur $\tilde{G}({\mathbb A}_{F})$ par les conditions:

- $\tilde{\lambda}_{12}$ vaut $1$ sur $\tilde{G}(F)$;

- pour $x\in G_{12}({\mathbb A}_{F})$ et $\gamma\in \tilde{G}_{12}({\mathbb A}_{F})$,  $\tilde{\lambda}_{12}(x\gamma)=\lambda_{12}(x)\tilde{\lambda}_{12}(\gamma)$.

 Il existe alors une unique fonction $\tilde{\lambda}_{12,V}$ sur $\tilde{G}_{12}(F_{V})$ de sorte que, pour $\gamma=\gamma_{V}\prod_{v\not\in V}\gamma_{v}\in \tilde{G}_{12}({\mathbb A}_{F})$, on ait l'\'egalit\'e
$$\tilde{\lambda}_{12}(\gamma)=\tilde{\lambda}_{12,V}(\gamma_{V})\prod_{v\not\in V}\tilde{\lambda}_{12,v}(\gamma_{v}).$$
Evidemment, elle v\'erifie (1).

Si $V\cap V_{12,ram}=\emptyset$, il y a aussi une fonction $\tilde{\lambda}_{12,V}$ canonique: le produit des $\tilde{\lambda}_{12,v}$ ci-dessus pour $v\in V$. Les formes lin\'eaires $r_{\tilde{M}_{1},\lambda_{1}}^{\tilde{G}_{1}}(.,\tilde{K}_{1,V})$ et $r_{\tilde{M}_{2},\lambda_{2}}^{\tilde{G}_{2}}(.,\tilde{K}_{2,V})$ se correspondent par l'isomorphisme 
$$D_{g\acute{e}om,\lambda_{1}}(\tilde{M}_{1}(F_{V}))\simeq D_{g\acute{e}om,\lambda_{2}}(\tilde{M}_{2}(F_{V}))$$
d\'eduit de cette fonction $\tilde{\lambda}_{12,V}$.

   \bigskip
  
  \subsection{$K$-espaces}
  On utilise dans ce paragraphe les notations usuelles pour divers ensembles de cohomologie: $H^1(F;G)$, $H^1(F_{v};G)$, $H^1({\mathbb A}_{F};G)$ etc... Par exemple $H^1(F;G)=H^1(\Gamma_{F};G(\bar{F}))$. On renvoie \`a [Lab] pour des d\'efinitions compl\`etes. 
  
On rappelle que l'on note $\pi:G_{SC}\to G$ la projection naturelle ainsi que les applications qui s'en d\'eduisent fonctoriellement. Ainsi, on a une application $\pi:H^1(F;G_{SC})\to H^1(F,G)$. D'autre part, pour $r\in \tilde{G}(F)$, l'application $ad_{r}$ d\'efinit naturellement un automorphisme de $H^1(F;G)$ qui ne d\'epend pas de $r$. On le note $\theta$. On note $Val_{{\mathbb R}}(F)$ l'ensemble des places r\'eelles de $F$. 
\ass{Lemme}{L'application naturelle
$$\pi(H^1(F;G_{SC}))\cap H^1(F;G)^{\theta}\to \prod_{v\in Val_{{\mathbb R}}(F)}\pi(H^1(F_v;G_{SC}))\cap H^1(F_v;G)^{\theta}$$
est bijective.}

Peuve. On commence par prouver que l'application
$$(1) \qquad \pi(H^1(F;G_{SC}))\to  \prod_{v\in Val_{{\mathbb R}}(F)}\pi(H^1(F_v;G_{SC}))$$
est bijective. La surjectivit\'e r\'esulte de celle de l'application
$$H^1(F;G_{SC})\to  \prod_{v\in Val_{{\mathbb R}}(F)}H^1(F_v;G_{SC}),$$
cf. [Lab] th\'eor\`eme 1.6.9. Soient $p,p'\in H^1(F;G_{SC})$ tels que $\pi(p)$ et $\pi(p')$ aient m\^eme image par l'application (1). Alors $\pi(p)$ et $\pi(p')$ ont m\^eme image dans $H^1({\mathbb A}_{F};G)$ (les images aux places non r\'eelles sont triviales). Parce qu'ils proviennent de $G_{SC}$, les \'el\'ements $\pi(p)$ et $\pi(p')$ ont aussi une image nulle dans $H^1_{ab}(F;G)$, cf. [Lab] 1.6 pour la d\'efinition de ce groupe. D'apr\`es le th\'eor\`eme 1.6.10 de [Lab], ces deux propri\'et\'es entra\^{\i}nent $\pi(p)=\pi(p')$. D'o\`u l'injectivit\'e de (1).

Il est imm\'ediat que l'application (1) est \'equivariante pour les actions naturelles de $\theta$. La bijectivit\'e de (1) entra\^{\i}ne donc celle de l'application obtenue en rempla\c{c}ant les espaces de d\'epart et d'arriv\'ee par leurs sous-espaces d'invariants par $\theta$. Cette application n'est autre que celle de l'\'enonc\'e. $\square$

On devra \`a diverses occasions travailler non pas avec les donn\'ees $(G,\tilde{G},{\bf a})$ de 1.1, mais avec une collection finie de telles donn\'ees, que l'on appellera un $K$-espace. On a d\'efini  (apr\`es Arthur) cette notion en [I] 1.11 dans le cas local. La d\'efinition s'\'etend au cas global. Rappelons-la. On consid\`ere une famille finie $(G_{p},\tilde{G}_{p})_{p\in \Pi}$ o\`u, pour tout $p$, $G_{p}$ est un groupe r\'eductif connexe d\'efini sur $F$ et $\tilde{G}_{p}$ est un espace tordu sur $G_{p}$.  On suppose donn\'ees des familles $(\phi_{p,q})_{p,q\in \Pi}$, $(\tilde{\phi}_{p,q})_{p,q\in \Pi}$ et $(\nabla_{p,q})_{p,q\in \Pi}$. Pour $p,q\in \Pi$, $\phi_{p,q}:G_{q}\to G_{p}$ et $\tilde{\phi}_{p,q}:\tilde{G}_{q}\to \tilde{G}_{p}$ sont des isomorphismes compatibles d\'efinis sur $\bar{F}$  et $\nabla_{p,q}:\Gamma_{F}\to G_{p,SC}(\bar{F})$ est un cocycle. On suppose les hypoth\`eses (1) \`a (5) v\'erifi\'ees pour tous $p,q,r\in \Pi$ et $\sigma\in \Gamma_{{\mathbb R}}$:
 
 (1) $\phi_{p,q}\circ\sigma(\phi_{p,q})^{-1}=ad_{\nabla_{p,q}(\sigma)}$ et 
 $\tilde{\phi}_{p,q}\circ\sigma(\tilde{\phi}_{p,q})^{-1}=ad_{\nabla_{p,q}(\sigma)}$ (ce dernier automorphisme est l'action par conjugaison de $\nabla_{p,q}(\sigma)$ sur $\tilde{G}_{p}$);
 
 (2) $\phi_{p,q}\circ\phi_{q,r}=\phi_{p,r}$ et $\tilde{\phi}_{p,q}\circ\tilde{\phi}_{q,r}=\tilde{\phi}_{p,r}$;
 
 (3) $\nabla_{p,r}(\sigma)=\phi_{p,q}(\nabla_{q,r}(\sigma))\nabla_{p,q}(\sigma)$;
 
 (4) $\tilde{G}_{p}(F)\not=\emptyset$;

 (5) la famille $(\nabla_{p,q})_{q\in \Pi}$ s'envoie bijectivement sur $\pi(H^1(F;G_{p,SC}))\cap H^1(F;G_{p})^{\theta}$.
 
 Sous ces hypoth\`eses, on d\'efinit le $K$-groupe $KG$ comme la r\'eunion disjointe des $G_{p}$ pour $p\in \Pi$ et le $K$-espace $K\tilde{G}$ comme la r\'eunion disjointe des $\tilde{G}_{p}$. 
 
 Comme dans le cas local, $\underline{les}$ paires de Borel \'epingl\'ees des diff\'erents $G_{p}$ s'identifient et les $G_{p}$ ont un $L$-groupe commun, que l'on note $^LG$. La donn\'ee suppl\'ementaire d'un \'el\'ement ${\bf a}\in H^1(W_{F};Z(\hat{G}))/ker^1(W_{F};Z(\hat{G}))$ d\'etermine un caract\`ere de chaque $G_{p}({\mathbb A}_{F})$, que l'on note simplement $\omega$. 
 
 Pour tout $p\in \Pi$, on fixe une paire parabolique $(P_{p,0},M_{p,0})$ de $G_{p}$ d\'efinie sur $F$ et minimale. Pour toute place $v\in Val(F)$, on fixe une paire parabolique $(P_{p,v,0},M_{p,v,0})$ de $G_{p}$ d\'efinie sur $F_{v}$ et minimale, de sorte que $P_{p,v,0}\subset P_{p,0}$ et $M_{p,v,0}\subset M_{p,0}$. Soulignons que ces inclusions peuvent \^etre strictes: le groupe $G_{p}$ peut \^etre plus d\'eploy\'e sur $F_{v}$ que sur $F$. Il se d\'eduit de ces paires des paires d'espaces $(\tilde{P}_{p,0},\tilde{M}_{p,0})$ et $(\tilde{P}_{p,v,0},\tilde{M}_{p,v,0})$. 
 
 Soit $v$ une place de $F$ finie ou complexe. Alors $H^1(F_v;G_{p,SC})=\{1\}$ pour tout $p$. Fixons un \'el\'ement $p'\in \Pi$. On peut  fixer pour tout $p\in \Pi$ un \'el\'ement $x_{p}\in G_{p',SC}(\bar{F}_{v})$ tel que $\nabla_{p',p}(\sigma)=x_{p}^{-1}\sigma(x_{p})$ pour tout $\sigma\in \Gamma_{F_{v}}$.   D\'efinissons $\phi'_{p',p}=ad_{x_{p}}\circ\phi_{p',p}$ et $\tilde{\phi}'_{p',p}=ad_{x_{p}}\circ\tilde{\phi}_{p',p}$.   
 Alors $\phi'_{p',p}:G_{p}\to G_{p'}$ et $\tilde{\phi}'_{p',p}:\tilde{G}_{p}\to \tilde{G}_{p'}$ sont des isomorphismes d\'efinis sur $F_{v}$.   Quitte \`a multiplier $x_{p}$ \`a gauche par un \'el\'ement de $G_{p',SC}(F_{v})$, on peut supposer que $\phi'_{p',p}$ envoie $(P_{p,v,0},M_{p,v,0})$ sur $(P_{p',v,0},M_{p',v,0})$ et que $\tilde{\phi}_{p',p}$ envoie $(\tilde{P}_{p,v,0},\tilde{M}_{p,v,0})$ sur $(\tilde{P}_{p',v,0},\tilde{M}_{p',v,0})$.  Les diff\'erents espaces $I(\tilde{G}_{p}(F_{v}),\omega)$ s'identifient gr\^ace \`a ces isomorphismes $\tilde{\phi}'_{p',p}$ \`a l'espace $I(\tilde{G}_{p'}(F_{v}),\omega)$, que l'on peut noter simplement $I(\tilde{G}(F_{v}),\omega)$. Notons que, bien que les $\tilde{\phi}'_{p',p}$ d\'ependent du choix des $x_{p}$, les isomorphismes qui s'en d\'eduisent entre les espaces  $I(\tilde{G}_{p}(F_{v}),\omega)$ ne d\'ependent pas de ce choix. 
  
  Puisque les diff\'erents espaces $\tilde{G}_{p}$ s'identifient sur $F_{v}$ pour toute place finie, les ensembles $V_{ram}(\tilde{G}_{p},{\bf a})$ sont les m\^emes, on les note $V_{ram}(\tilde{G},{\bf a})$, ou simplement  $V_{ram}$. Pour l'\'el\'ement $p'$ fix\'e ci-dessus et pour  toute place $v\not\in V_{ram}$, on choisit un espace hypersp\'ecial $\tilde{K}_{p',v}$ de $\tilde{G}_{p'}(F_{v})$ de sorte que les conditions de 1.1 soient v\'erifi\'ees. Pour $p\in \Pi$, on note $\tilde{K}_{p,v}$ l'image de $\tilde{K}_{p',v}$ par $(\tilde{\phi}'_{p',p})^{-1}$.
    Alors les fonctions caract\'eristiques ${\bf 1}_{\tilde{K}_{p,v}}$ de $\tilde{K}_{p,v}$, pour $p\in \Pi$,  ont m\^eme image dans $I(\tilde{G}(F_{v}),\omega)$. On peut noter ${\bf 1}_{\tilde{K}_{v}}$ cette image. 
    
  Les choses sont plus compliqu\'ees en une place r\'eelle. Soit $v$ une telle place. Pour $p,q\in \Pi$, notons $\nabla_{p,q,v}$ la restriction de $\nabla_{p,q}$ \`a $\Gamma_{F_{v}}$. Disons que $p$ et $q$ sont $v$-\'equivalents si et seulement si $\pi(\nabla_{p,q,v})$ est cohomologiquement trivial. On v\'erifie que c'est une relation d'\'equivalence.  Fixons un ensemble $\Pi_{v}\subset \Pi$ de repr\'esentants des classes de $v$-\'equivalence. On voit que la famille $(G_{p,v},\tilde{G}_{p,v})_{p\in \Pi_{v}}$, munie des familles $(\phi_{p,q})_{p,q\in \Pi_{v}}$, $(\tilde{\phi}_{p,q})_{p,q\in \Pi_{v}}$ et $(\nabla_{p,q,v})_{p,q\in \Pi_{v}}$, d\'efinit un $K$-espace sur $F_{v}$, au sens de [I] 1.11. Soit $p\in \Pi$. Notons $p'$ l'unique \'el\'ement de $\Pi_{v}$ tel que $p$ soit $v$-\'equivalent \`a $p'$. On peut fixer $x_{p}\in G_{p'}(\bar{F}_{v})$ tel que $\pi(\nabla_{p',p,v})=x_{p}^{-1}\sigma(x_{p})$. Les applications $\phi'_{p',p}=ad_{x_{p}}\circ\phi_{p',p}$ et $\tilde{\phi}'_{p',p}=ad_{x_{p}}\circ \tilde{\phi}_{p',p}$ sont encore des isomorphismes d\'efinis sur $F_{v}$. Quitte \`a multiplier $x_{p}$ \`a gauche par un \'el\'ement de $G_{p'}(F_{v})$, on peut supposer que $\phi'_{p',p}$ envoie $(P_{p,v,0},M_{p,v,0})$ sur $(P_{p',v,0},M_{p',v,0})$ et que $\tilde{\phi}_{p',p}$ envoie $(\tilde{P}_{p,v,0},\tilde{M}_{p,v,0})$ sur $(\tilde{P}_{p',v,0},\tilde{M}_{p',v,0})$. D\'ecomposons $x_{p}$ en $\pi(x_{p,sc})z_{p}$, o\`u $x_{p,sc}\in G_{p',SC}(\bar{F}_{v})$ et $z_{p}\in Z(G_{p'};\bar{F}_{v})$. D\'efinissons la cocha\^{\i}ne $\nabla_{p}:\Gamma_{F_{v}}\to G_{p',SC}$ par $\nabla_{p}(\sigma)=x_{p,sc}\nabla_{p',p}(\sigma)\sigma(x_{p,sc})^{-1}$. C'est un cocycle \`a valeurs dans $Z(G_{p',SC})$ et le couple $(\nabla_{p},z_{p})$ d\'efinit un \'el\'ement de $H^{1,0}(F_{v}; Z(G_{p',SC})\to Z(G_{p'}))=G_{p',ab}(F_{v})$, que l'on note $h_{p}$. On sait que ${\bf a}$ d\'efinit un caract\`ere de ce groupe. On note $\omega(h_{p})$ sa valeur en $h_{p}$. On d\'efinit un homomorphisme
  $$\begin{array}{ccc}C_{c}^{\infty}(\tilde{G}_{p}(F_{v}))&\to& C_{c}^{\infty}(\tilde{G}_{p'}(F_{v}))\\ f&\mapsto &f'\\ \end{array}$$
  par $f'\circ\tilde{\phi}'_{p',p}(r)=\omega(h_{p})f(r)$ pour tout $r \in \tilde{G}_{p}(F_{v})$. C'est un isomorphisme qui se quotiente en un isomorphisme de $I(\tilde{G}_{p}(F_{v}),\omega)$ sur $I(\tilde{G}_{p'}(F_{v}),\omega)$. On v\'erifie que ce dernier isomorphisme ne d\'epend pas des choix de $x_{p}$ et $x_{p,sc}$.
  
Soit $V$ un ensemble fini de places contenant les places r\'eelles.  Gr\^ace au lemme ci-dessus et aux diff\'erents isomorphismes que l'on vient de construire, on obtient  un isomorphisme
$$\oplus_{p\in \Pi}I(\tilde{G}_{p}(F_{V}),\omega)\simeq \left(\otimes_{v\in V-Val_{{\mathbb R}}(F)}I(\tilde{G}(F_{v}),\omega)\right)\otimes\left(\otimes_{v\in Val_{{\mathbb R}}(F)}\left(\oplus_{p'\in \Pi_{v}}I(\tilde{G}_{p'}(F_{v}),\omega)\right)\right).$$
On note $I(K\tilde{G}(F_{V}),\omega)$ le membre de gauche. On a \'evidemment des isomorphismes duaux pour les espaces de distributions.

Les diff\'erents groupes $G_{p}$ \'etant formes int\'erieures l'un de l'autre, les espaces de mesures de Haar $Mes(G_{p}(F_{v}))$ s'identifient pour toute place $v$ \`a un espace commun que l'on note $Mes(G(F_{v}))$. De m\^eme, les espaces  $\mathfrak{A}_{\tilde{G}_{p}}$ s'identifient \`a un espace commun $\mathfrak{A}_{\tilde{G}}$.  Comme dans le cas local, pour un triplet $(G,\tilde{G},{\bf a})$ comme en 1.1, on peut construire un $K$-espace dont $\tilde{G}$ soit une composante connexe, cf. [I] 1.11.

\bigskip

\subsection{$K$-espaces de Levi}
On poursuit avec les m\^emes donn\'ees que dans le paragraphe pr\'ec\'edent. Notons $(B^*,T^*)$ $\underline{la}$ paire de Borel commune des groupes $G_{p}$ et $\Delta$ l'ensemble de racines simples associ\'e. Le groupe $\Gamma_{F}$ et l'automorphisme $\theta$ agissent sur $\Delta$. Toute paire parabolique $(P_{p},M_{p})$ de $G_{p}$ d\'etermine un sous-ensemble $\Delta^{M_{p}}\subset \Delta$. Ainsi, aux paires $(P_{p,0},M_{p,0})$ et $(P_{p,v,0},M_{p,v,0})$, pour $v\in Val(F)$, sont associ\'es des sous-ensembles $\Delta^{M_{p,0}}$ et $\Delta^{M_{p,v,0}}$ de $\Delta$.  Ils sont invariants par l'action de $\theta$. L'ensemble $\Delta^{M_{p,0}}$ est invariant par l'action de $\Gamma_{F}$ et l'ensemble $\Delta^{M_{p,v,0}}$ est invariant par l'action de $\Gamma_{F_{v}}$. Pour une place $v$ r\'eelle, on a vu que notre $K$ espace d\'eterminait un $K$-espace sur $F_{v}={\mathbb R}$. En [I] 3.5, on lui a associ\'e un sous-ensemble $\Delta_{min}\subset \Delta$, qu'il convient de noter $\Delta_{min,v}$. Il est invariant par $\theta$. Le lemme de cette r\'ef\'erence dit que $\Delta_{min,v}\subset \Delta^{M_{p,v,0}}$ pour tout $p$ et que cette inclusion est une \'egalit\'e pour au moins un $p$. La m\^eme construction vaut  pour une place complexe ou non-archim\'edienne. On r\'ef\`ere pour ce cas \`a [A 10] lemme 2.1. Pour une telle place, on d\'efinit donc un sous-ensemble $\Delta_{min,v}\subset \Delta$. Puisque, pour une telle place,  les groupes $G_{p}$ sont tous isomorphes sur $F_{v}$, on a simplement $\Delta_{min,v}=\Delta^{M_{p,v,0}}$ pour tout $p$. On vient de voir que, pour une place $v$ quelconque, il existait $p\in \Pi$ tel que $\Delta_{min,v}=\Delta^{M_{p,v,0}}$. Le lemme du paragraphe pr\'ec\'edent renforce ce r\'esultat en \'echangeant les quantificateurs: 

(1) il existe $p\in \Pi$ tel que  $\Delta_{min,v}=\Delta^{M_{p,v,0}}$ pour tout $v$. 

Notons $\Delta'_{min}=\cup_{v\in Val(F)}\Delta_{min,v}$ et $\Delta_{min}$ le plus petit sous-ensemble de $\Delta$ qui contient 
$ \Delta'_{min}$ et qui est invariant par l'action de $\Gamma_{F}$. Puisque les actions de $\Gamma_{F}$ et $\theta$ commutent, $\Delta_{min}$ est invariant par $\theta$. 

\ass{Lemme}{On a l'inclusion $\Delta_{min}\subset \Delta^{M_{p,0}}$ pour tout $p$. Cette inclusion est une \'egalit\'e pour $p$ v\'erifiant (1).}

Preuve.  Pour tout $p\in \Pi$ et toute place $v$, on a l'inclusion $\Delta^{M_{p,v,0}}\subset \Delta^{M_{p,0}}$. Donc $\Delta^{M_{p,0}}$ contient $\Delta'_{min}$. Puisqu'il est invariant par $\Gamma_{F}$, il contient $\Delta_{min}$. Fixons $p$ v\'erifiant (1), posons $G=G_{p}$, $M_{0}=M_{p,0}$ et $M_{v,0}=M_{p,v,0}$ pour tout $v$. Introduisons une forme quasi-d\'eploy\'ee $G^*$ de $G$ et un torseur int\'erieur $\psi:G\to G^*$  de sorte que $(B^*,T^*)$ s'identifie \`a une paire de Borel de $G^*$ d\'efinie sur $F$. Puisque $\Delta_{min}$ est invariant par $\Gamma_{F}$, il existe une paire parabolique $(P^*,M^*)$ de $G^*$, d\'efinie sur $F$, contenant $(B^*,T^*)$, de sorte que $\Delta_{min}=\Delta^{M^*}$. Soit $v$ une place de $F$. Puisque $\Delta_{min}$ est invariant par $\Gamma_{F_{v}}$ et que $\Delta^{M_{v,0}}\subset \Delta_{min}$, la paire $(P^*,M^*)$ se transf\`ere en une paire parabolique de $G$ d\'efinie sur $F_{v}$. Comme on le sait, \`a la forme int\'erieure $G$ de $G^*$ est associ\'e un \'el\'ement $\kappa\in H^1(F;G_{AD}^*)$. Cet \'el\'ement se localise en un \'el\'ement $\kappa_{v}\in H^1(F_v;G^*_{AD})$. La condition pr\'ec\'edente signifie que $\kappa_{v}$ appartient \`a l'image de l'application $H^1(F_v;M^*_{ad})$. D'apr\`es [Lab] proposition 1.6.12, on a un diagramme commutatif
$$(2)\qquad \begin{array}{ccccc}H^1(F;G^*_{AD})&\to& \oplus_{v\in Val(F)}H^1(F_v;G^*_{AD})&\to&H^{2,1}({\mathbb A}_{F}/F;T^*_{sc}\to T^*_{ad})\\ \,\,\uparrow \iota_{F}&&\,\,\uparrow \iota_{{\mathbb A}_{F}}&&\,\,\uparrow \iota_{ab}\\ H^1(F;M^*_{ad})&\to& \oplus_{v\in Val(F)}H^1(F_v;M^*_{ad})&\to&H^{2,1}({\mathbb A}_{F}/F;T^*_{M^*,sc}\to T^*_{ad})\\ \end{array}$$
Conform\'ement \`a la convention de [I], on a not\'e $H^{2,1}$ les groupes de cohomologie que Labesse note $H^1$ et que Kottwitz et Shelstad notent $H^2$. On renvoie \`a ces auteurs (ou \`a 3.6 ci-dessous) pour la d\'efinition de la cohomologie $H^{2,1}({\mathbb A}_{F}/F;.)$.  On a not\'e $T^*_{sc}$ et  $T^*_{M^*,sc}$ les images r\'eciproques de $T^*$ dans les rev\^etements simplement connexes $G^*_{SC}$ et $M^*_{SC}$ des groupes d\'eriv\'es de $G^*$ et $M^*$. Posons $\kappa_{{\mathbb A}_{F}}=\oplus_{v\in Val(F)}\kappa_{v}$. D'apr\`es ce que l'on a dit ci-dessus, on peut fixer $\kappa^{M^*}_{{\mathbb A}_{F}}\in \oplus_{v\in Val(F)}H^1(F_v;M^*_{ad})$ tel que $\iota_{{\mathbb A}_{F}}(\kappa^{M^*}_{{\mathbb A}_{F}})=\kappa_{{\mathbb A}_{F}}$. Montrons que

(3) $\iota_{ab}$ est injectif. 

Le groupe $X_{*}(T^*_{sc})$ a pour base l'ensemble $\check{\Delta}$ de coracines associ\'e \`a $\Delta$. Le groupe $X_{*}(T^*_{M^*,sc})$ a pour base le sous-ensemble $\check{\Delta}^{M^*}$ associ\'e \`a $\Delta^{M^*}$. On a donc une suite exacte
$$1\to T^*_{M^*,sc}\to T^*_{sc}\to T_{1}\to 1,$$
o\`u $T_{1}$ est un tore tel que $X_{*}(T_{1})$ a pour base l'ensemble $\check{\Delta}-\check{\Delta}^{M^*}$. Puisque $\Gamma_{F}$ agit sur $\Delta$ par permutations laissant stable $\Delta^{M^*}$, les trois tores sont induits. On a le diagramme commutatif suivant 
$$\begin{array}{ccccccccc}1&\to&T^*_{M^*,sc}&\to &T^*_{sc}&\to&T_{1}&\to&1\\ &&\downarrow&&\downarrow&&\downarrow&&\\ 1&\to&T^*_{ad}&\to&T^*_{ad}&\to &1&&\\ \end{array}$$
 dont les suites horizontales sont exactes. Il s'en d\'eduit une suite exacte
$$H^1({\mathbb A}_{F}/F;T_{1})\to H^{2,1}({\mathbb A}_{F}/F;T^*_{M^*,sc}\to T^*_{ad})\stackrel{\iota_{ab}}{\to}H^{2,1}({\mathbb A}_{F}/F;T^*_{sc}\to T^*_{ad}).$$
Puisque $T_{1}$ est induit, le premier groupe est  nul, ce qui prouve (3).

Puisque $\kappa_{{\mathbb A}_{F}}$ provient de $\kappa$, son image dans $H^{2,1}({\mathbb A}_{F}/F;T^*_{sc}\to T^*_{ad})$ est nulle. La commutativit\'e du diagramme (2) et l'assertion (3) entra\^{\i}nent que l'image de $\kappa^{M^*}_{{\mathbb A}_{F}}$ dans $H^{2,1}({\mathbb A}_{F}/F;T^*_{M^*,sc}\to T^*_{ad})$ est nulle. Donc $\kappa^{M^*}_{{\mathbb A}_{F}}$ provient d'un \'el\'ement $\kappa^{M^*}\in H^1(F;M^*_{ad})$. Notons $\kappa'$ l'image de cet \'el\'ement dans $H^1(F;G_{AD}^*)$.  Les \'el\'ements $\kappa$ et $\kappa'$ ont m\^eme image $\kappa_{{\mathbb A}_{F}}$ dans $\oplus_{v\in Val(F)}H^1(F_v;G^*_{AD})$. D'apr\`es [Lab] proposition 1.6.12, la fibre au-dessus de $\kappa_{{\mathbb A}_{F}}$ de l'application 
$$H^1(F;G^*_{AD})\to \oplus_{v\in Val(F)}H^1(F_v;G^*_{AD})$$
est isomorphe au noyau $ker^1(F;G_{AD})$ de l'application similaire
$$H^1(F;G_{AD})\to \oplus_{v\in Val(F)}H^1(F_v;G_{AD}).$$
 Or, $G_{AD}$ \'etant adjoint, cet ensemble est r\'eduit \`a $\{1\}$ ([S] corollaire 5.4). Donc $\kappa=\kappa'$ et $\kappa$ provient d'un \'el\'ement de $H^1(F:M^*_{ad})$. Cela \'equivaut \`a dire que la paire parabolique $(P^*,M^*)$ se transf\`ere \`a $G$, c'est-\`a-dire qu'il existe une paire parabolique de $G$ qui est d\'efinie sur $F$ et qui est conjugu\'ee \`a $\psi^{-1}(P^*,M^*)$. Cela implique $\Delta^{M_{0}}\subset \Delta^{M^*}$. D'o\`u le lemme. $\square$

Notons $\Pi^{M_{0}}$ le sous-ensemble des $p\in \Pi$ tels que $\Delta^{M_{p,0}}=\Delta_{min}$. Il est non vide d'apr\`es le lemme. Pour $p\in \Pi$ et $p'\in \Pi^{M_{0}}$, fixons un \'el\'ement $x_{p',p}\in G_{p'}$ tel que $ad_{x_{p',p}}\circ \tilde{\phi}_{p',p}(\tilde{P}_{p,0},\tilde{M}_{p,0})$ contienne $(\tilde{P}_{p',0},\tilde{M}_{p',0})$. Le lemme implique l'existence d'un tel \'el\'ement. On note $K\tilde{M}_{0}$ la famille $(\tilde{M}_{p,0})_{p\in \Pi^{M_{0}}}$. C'est un $K$-espace de Levi de $K\tilde{G}$, la d\'efinition de cette notion \'etant similaire \`a celle du cas local, cf. [I] 3.5. On note ${\cal L}(K\tilde{M}_{0})$  l'ensemble des $K$-espaces de Levi $K\tilde{L}=(\tilde{L}_{p})_{p\in \Pi^L}$ de $K\tilde{G}$ v\'erifiant les deux conditions suivantes:

- $\tilde{L}_{p}\supset \tilde{M}_{p,0}$ pour tout $p\in \Pi^L$;

- $ad_{x_{p',p}}\circ \tilde{\phi}_{p',p}(\tilde{L}_{p})=\tilde{L}_{p'}$ pour tous $p\in \Pi^L$, $p'\in \Pi^{M_{0}}$.

Plus g\'en\'eralement, pour $K\tilde{M}=(\tilde{M}_{p})_{p\in \Pi^M}\in {\cal L}(\tilde{M}_{0})$, on note ${\cal L}(K\tilde{M})$ l'ensemble des $K\tilde{L}=(\tilde{L}_{p})_{p\in \Pi^L}\in {\cal L}(K\tilde{M}_{0})$ tels que $\Pi^M\subset \Pi^L$ et $\tilde{M}_{p}\subset \tilde{L}_{p}$ pour tout $p\in \Pi^M$. On d\'efinit de fa\c{c}on similaire les ensembles ${\cal P}(K\tilde{M})$ et ${\cal F}(K\tilde{M})$. 

Dans la suite, on travaillera avec un triplet $(G,\tilde{G},{\bf a})$ comme en 1.1, en n'introduisant les $K$-espaces que lorsque cela sera indispensable.  Dans ce cas, on supposera fix\'e un espace de Levi minimal $K\tilde{M}_{0}$. 

  \bigskip

\section{La partie g\'eom\'etrique de la formule des traces}

\bigskip

\subsection{La partie g\'eom\'etrique de la formule des traces non invariante}
 On a d\'efini en 1.1 le noyau $G({\mathbb A}_{F})^1$ de l'homomorphisme $H_{\tilde{G}}$ (contrairement \`a ce que sugg\`ere la notation, ce groupe d\'epend de $\tilde{G}$ et pas seulement de $G$). On a la d\'ecomposition en produit direct $G({\mathbb A}_{F})= \mathfrak{A}_{\tilde{G}}\times G({\mathbb A}_{F})^1$.   On note $\tilde{G}({\mathbb A}_{F})^1$ l'ensemble des $\gamma\in\tilde{G}({\mathbb A}_{F})$ tels que $\tilde{H}_{\tilde{G}}(\gamma)=0$ (c'est-\`a-dire $\tilde{G}({\mathbb A}_{F})^1=G({\mathbb A}_{F})^1 \tilde{G}(F) $). On a encore la d\'ecomposition en produit direct $\tilde{G}({\mathbb A}_{F})= \mathfrak{A}_{\tilde{G}}\times \tilde{G}({\mathbb A}_{F})^1$.

Soit $f\in C_{c}^{\infty}(\tilde{G}({\mathbb A}_{F}))$. C'est-\`a-dire que $f$ est combinaison lin\'eaire de fonctions $\otimes_{v\in Val(F)}f_{v}$ o\`u, pour tout $v$, $f_{v}\in C_{c}^{\infty}(\tilde{G}(F_{v}))$ et, pour presque tout $v$, $f_{v}={\bf 1}_{\tilde{K}_{v}}$. Notons $f^{1} $ sa restriction \`a $\tilde{G}({\mathbb A}_{F})^1$. D\'efinissons une fonction $f_{*}$ sur $\tilde{G}({\mathbb A}_{F})$ par $f_{*}(a\gamma)=f^1(\gamma)$ pour tous $a\in \mathfrak{A}_{\tilde{G}}$ et $\gamma\in \tilde{G}({\mathbb A}_{F})^1$. Pour un instant, d\'efinissons l'espace  $\bar{C}_{c}^{\infty}(\tilde{G}({\mathbb A}_{F}))$ un peu plus gros que $C_{c}^{\infty}(\tilde{G}({\mathbb A}_{F}))$, obtenu en rempla\c{c}ant la composante archim\'edienne $\otimes_{v\in Val_{\infty}(F)}C_{c}^{\infty}(\tilde{G}(F_{v}))$ par $C_{c}^{\infty}(\prod_{v\in Val_{\infty}(F)}\tilde{G}(F_{v}))$. Introduisons une fonction $\varphi\in \bar{C}_{c}^{\infty}(\tilde{G}({\mathbb A}_{F}))$ telle que
$$f_{*}(\gamma)=\int_{\mathfrak{A}_{\tilde{G}}}\varphi^a(\gamma) \,da$$
pour tout $\gamma\in \tilde{G}({\mathbb A}_{F})$, 
o\`u $\varphi^a$ est la fonction $\varphi^a(\gamma)=\varphi(a\gamma)$.
 Fixons une mesure de Haar $dg$ sur $G({\mathbb A}_{F})$. Posons ${\bf f}=f\otimes dg\in C_{c}^{\infty}(\tilde{G}({\mathbb A}_{F}))\otimes Mes(G({\mathbb A}_{F}))$. Nous notons $J^{\tilde{G}}_{g\acute{e}om}({\bf f},\omega)$ le membre de gauche de l'\'egalit\'e du th\'eor\`eme 11.3.2 de [LabW] appliqu\'e \`a la fonction $\varphi$, multipli\'e par  $\vert det((1-\theta)_{\vert {\cal A}_{G}/{\cal A}_{\tilde{G}}})\vert ^{-1}$. Dans le cas o\`u $\omega=1$, c'est aussi le terme qu'Arthur note $J_{g\acute{e}om}(f^{1})$. 

Rappelons la d\'efinition pour m\'emoire.  On renvoie \`a [LabW] pour les notations, qui sont d'ailleurs les notations standard. Pour un espace parabolique standard $\tilde{P}=\tilde{M}U_{P}$, on d\'efinit une fonction $K_{\tilde{P}}(f)$ sur $G({\mathbb A}_{F})$ par
$$K_{\tilde{P}}(f,g)=\int_{U_{P}(F)\backslash U_{P}({\mathbb A}_{F})}\sum_{\gamma\in \tilde{P}(F)}f(g^{-1}\gamma ug)\,du.$$
Soit $T\in {\cal A}_{\tilde{M}_{0}}$ assez positif relativement \`a $\tilde{P}_{0}$. On pose
$$k^T_{g\acute{e}om}(f,g)=\sum_{\tilde{P}\supset \tilde{P}_{0}}(-1)^{a_{\tilde{P}}-a_{\tilde{G}}}\sum_{\xi\in P(F)\backslash G(F)}\hat{\tau}_{\tilde{P}}(H_{\tilde{P}}(\xi g)-T)K_{\tilde{P}}(f,\xi g),$$
puis
$$J^T_{g\acute{e}om}(f,\omega)=\int_{\mathfrak{A}_{\tilde{G}}G(F)\backslash G({\mathbb A}_{F})}k^T_{g\acute{e}om}(f,g)\omega(g)\,dg.$$
Quand $T$ tend vers l'infini dans la chambre positive d\'etermin\'ee par $\tilde{P}_{0}$, cette expression est asymptote \`a un polyn\^ome en $T$. Alors $J_{g\acute{e}om}^{\tilde{G}}({\bf f},\omega)$ est la valeur de ce polyn\^ome au point $T_{0}$ d\'efini en [LabW] lemme 3.3.3.

  On note $\tilde{G}_{ss}(F)$, resp. $\tilde{G}(F)_{ell}$, l'ensemble des \'el\'ements semi-simples  de $\tilde{G}(F)$, resp. celui des \'el\'ements semi-simples et elliptiques. On note  $\tilde{G}_{ss}(F)/conj$, resp. $\tilde{G}(F)_{ell}/conj$, l'ensemble des classes de conjugaison par $G(F)$ dans $\tilde{G}_{ss}(F)$, resp. $\tilde{G}(F)_{ell}$. 

{\bf Remarque.} Puisque $F$ est un corps de nombres, les deux d\'efinitions possibles de l'ellipticit\'e sont \'equivalentes. C'est-\`a-dire que, pour un \'el\'ement semi-simple $\gamma\in \tilde{G}(F)$, $\gamma$ est elliptique si et seulement s'il v\'erifie les conditions \'equivalentes:

(1)  il existe un sous-tore tordu maximal $\tilde{T}$ de $\tilde{G}$ d\'efini sur ${\mathbb Q}$ et elliptique tel que $\gamma\in \tilde{ T}(F)$;

(2)  $A_{G_{\gamma}}=A_{\tilde{G}}$.

\bigskip

Soit ${\cal O}\in \tilde{G}_{ss}(F)/conj$. Pour un espace parabolique standard $\tilde{P}$, on d\'efinit une fonction $K_{\tilde{P},{\cal O}}(f)$ en restreignant dans la d\'efinition ci-dessus la somme sur $\gamma\in \tilde{P}(F)$ aux $\gamma$ dont la partie semi-simple appartient \`a ${\cal O}$. En poursuivant comme ci-dessus, on construit des termes $k^T_{{\cal O}}(f,g)$, $J^T_{{\cal O}}(f,\omega)$, puis  $J_{{\cal O}}({\bf f},\omega)$ (ou $J_{{\cal O}}^{\tilde{G}}({\bf f},\omega)$ s'il est bon de pr\'eciser l'espace). Tous ces termes sont nuls 	sauf si la classe de conjugaison par $G({\mathbb A}_{F})$ engendr\'ee par ${\cal O}$ coupe le support de $f$. Il n'y a qu'un nombre fini de ${\cal O}$ v\'erifiant cette condition, en vertu du lemme suivant. En vue d'une application ult\'erieure, on \'enonce celui-ci sous une forme plus forte qu'il ne nous est n\'ecessaire \`a pr\'esent.

\ass{Lemme}{Soit  $V$ un ensemble fini de places contenant $V_{ram}$. Soit $\tilde{U}_{V}$ un sous-ensemble compact de $\tilde{G}(F_{V})$. Consid\'erons l'ensemble $\Xi$ des  classes de conjugaison ${\cal O}$ par $G(F)$ d'\'el\'ements semi-simples de $\tilde{G}(F)$ telles que  la classe de conjugaison ${\cal O}_{V}$ par $G(F_{V})$ engendr\'ee par ${\cal O}$ 
 coupe $\mathfrak{A}_{\tilde{G}} \tilde{U}_{V}$ et, pour tout $v\not\in V$, la classe de conjugaison   par $G(F_{v})$ engendr\'ee par ${\cal O}$ coupe $\tilde{K}_{v}$. Alors
 
 (i) $\Xi$ est fini;
 
 (ii) il existe un sous-ensemble compact $C\subset \mathfrak{A}_{\tilde{G}}$ tel que, pour tout ${\cal O}\in \Xi$, ${\cal O}_{V}\cap \mathfrak{A}_{\tilde{G}} \tilde{U}_{V}$ soit contenu dans $C\tilde{U}_{V}$.}

Preuve. Soient ${\cal O}\in \Xi$ et $\gamma\in {\cal O}$. Les conditions impos\'ees hors de $V$ et la d\'efinition de $\tilde{H}_{\tilde{G}_{V}}$ entra\^{\i}nent que $\tilde{H}_{\tilde{G}_{V}}(\gamma)=0$. Soit $a\in \mathfrak{A}_{\tilde{G}}$ et $\gamma_{V}\in \tilde{U}_{V}$ tels que $a\gamma_{V}\in {\cal O}_{V}$. Alors $H_{\tilde{G}}(a)=-\tilde{H}_{\tilde{G}_{V}}(\gamma_{V})$ et, puisque $\tilde{U}_{V}$ est compact, $a$ reste dans un compact. Cela prouve (ii). Quitte \`a agrandir $\tilde{U}_{V}$, on peut donc  supposer que, pour ${\cal O}$ et $\gamma$ comme ci-dessus,  la classe de conjugaison de $\gamma$ par $G(F_{V})$ coupe $\tilde{U}_{V}$. On peut fixer $x_{V}\in G(F_{V})$ tel que $x_{V}^{-1}\gamma x_{V}\in \tilde{U}_{V}$ et, pour tout $v\not\in V$, un \'el\'ement $x_{v}\in G(F_{v})$ tel que $x_{v}^{-1}\gamma x_{v}\in \tilde{K}_{v}$. Mais $\gamma\in \tilde{K}_{v}$ pour presque tout $v$. On peut donc supposer $x_{v}=1$ pour presque tout $v$. Les \'el\'ements $x_{V}$ et $(x_{v})_{v\not\in V}$ se regroupent alors en un \'el\'ement $x\in G({\mathbb A}_{F})$ tel que $x^{-1}\gamma x\in \tilde{U}$, o\`u $\tilde{U}=\tilde{U}_{V}\times (\prod_{v\not\in V}\tilde{K}_{v})$. La classe de conjugaison par $G({\mathbb A}_{F})$ de $\gamma$ coupe donc $\tilde{U}$ et il reste \`a appliquer le lemme 3.7.4 de [LabW]. $\square$

On a
$$(3) \qquad J^{\tilde{G}}_{g\acute{e}om}({\bf f},\omega)=\sum_{{\cal O}\in \tilde{G}_{ss}(F)/conj}J_{{\cal O}}({\bf f},\omega).$$

Notons que tous ces termes sont nuls si $\omega$ n'est pas trivial sur $Z(G;{\mathbb A}_{F})^{\theta}$.

{\bf Cas particulier.} Dans le cas o\`u $\tilde{G}=G$, $\tilde{K}_{v}=K_{v}$ pour tout $v\not\in V_{ram}$ et ${\cal O}=\{1\}$, on remplace l'indice ${\cal O}$ par $unip$ dans les termes d\'efinis ci-dessus et on ajoute l'exposant $G$. Par exemple $J^{G}_{unip}({\bf f},\omega)=J_{\{1\}}({\bf f},\omega)$.

\bigskip

\subsection{ Le terme unipotent de la formule des traces non invariante}
Soient $V$ un ensemble fini de places contenant $V_{ram}$ et ${\cal O}\in \tilde{G}_{ss}(F)/conj$. 
On a dit en 1.1 que l'on pouvait identifier $C_{c}^{\infty}(\tilde{G}(F_{V}))$ \`a un sous-espace de $C_{c}^{\infty}(\tilde{G}({\mathbb A}_{F}))$ et identifier $Mes(G(F_{V}))$ \`a $Mes(G({\mathbb A}_{F}))$. Pour ${\bf  f}_{V}\in C_{c}^{\infty}(\tilde{G}(F_{V}))\otimes Mes(G(F_{V}))$, le terme $J_{{\cal O}}({\bf f}_{V},\omega)$ est donc bien d\'efini. 

On consid\`ere dans ce paragraphe le cas o\`u $\tilde{G}=G$, $\tilde{K}_{v}=K_{v}$ pour tout $v\not\in V_{ram}$ et o\`u ${\cal O}=1$.   En [A7] th\'eor\`eme 8.1, Arthur prouve l'existence d'une distribution
$A^G_{unip}(V,\omega)\in D_{orb}(G(F_{V}))\otimes Mes(G(F_{V}))^*$ qui v\'erifie les propri\'et\'es suivantes:

(1) pour tout ${\bf f}_{V}\in  C_{c}^{\infty}(G(F_{V}))\otimes Mes(G(F_{V}))$, on a l'\'egalit\'e:
$$I^{G}(A^{G}_{unip}(V,\omega),{\bf f}_{V})=J^G_{unip}({\bf f}_{V},\omega)-\sum_{M\in {\cal L}(M_{0}), M\not=G}\vert W^{M}\vert \vert W^{G}\vert ^{-1} $$
$$J_{M}^{G}(A^{M}_{unip}(V,\omega),{\bf f}_{V});$$

(2) $A^{G}_{unip}(V,\omega)$ est combinaison lin\'eaire d'int\'egrales orbitales associ\'ees \`a des \'el\'ements unipotents $u\in G(F)$, ou plus exactement aux projections dans $G(F_{V})$ de tels \'el\'ements;

(3) $A^G_{unip}(V,\omega)=0$  si $\omega$ n'est pas trivial sur $Z(G,{\mathbb A}_{F})$.

Supposons $\omega$ trivial sur $Z(G,{\mathbb A}_{F})$. En g\'en\'eral, on ne sait pas expliciter la distribution $A^G_{unip}(V,\omega)$. On peut toutefois la calculer dans le cas particulier o\`u $G$ est un tore. Dans ce cas, soit $f_{V}\in C_{c}^{\infty}(G(F_{V}))$ et $dg_{V}$ une mesure sur $G(F_{V})$, qui s'identifie \`a une mesure sur $G({\mathbb A}_{F})$ par produit avec les mesures canoniques  sur $G(F_{v})$ pour $v\not\in V$. Posons ${\bf f}_{V}=f_{V}\otimes dg_{V}$. Alors
$$I^{G}(A_{unip}^G(V,\omega),{\bf f}_{V})=mes(\mathfrak{A}_{G}G(F)\backslash G({\mathbb A}_{F}))f_{V}(1).$$

La distribution $A^G_{unip}(V,\omega)$ d\'epend, comme la partie g\'eom\'etrique de la formule des traces non-invariante, des choix du Levi minimal $M_{0}$ et des compacts $K_{v}$ pour tout $v$. Rappelons que ceux-ci sont suppos\'es en bonne position relativement \`a $M_{0}$. Arthur montre que $A^G_{unip}(V,\omega)$ ne d\'epend pas des $K_{v}$ pour $v\in V$ ([A9] proposition 13.2).   Notons plus pr\'ecis\'ement $A^G_{unip}(V,\omega,M_{0},K^V)$ notre distribution. Soit $x\in G(F)M_{0}({\mathbb A}_{F})$, rempla\c{c}ons $M_{0}$ et les $K_{v}$ par $M'_{0}=ad_{x}(M_{0})$ et les $K'_{v}=ad_{x}(K_{v})$. Un raisonnement par transport de structure montre que 
$$(4) \qquad A^G_{unip}(V,\omega,M'_{0},K^{'V})=\omega(x^V)^{-1}A^G_{unip}(V,\omega,M_{0},K^V),$$
 o\`u $x^V$ est la projection de $x$ dans $G({\mathbb A}_{F}^V)$.  On peut alors  modifier formellement la d\'efinition de la fa\c{c}on suivante. On oublie pour un temps les objets fix\'es en 1.1. Pour tout $v\not\in V$, soit $K_{v}$ un sous-groupe compact hypersp\'ecial de $G(F_{v})$. On suppose comme en 1.1 que $K_{v}=G(\mathfrak{o}_{v})$ pour presque tout $v$. Fixons un Levi minimal $M_{0}$. On peut choisir $x^V=(x_{v})_{v\not\in V}\in G({\mathbb A}_{F}^V)$ de sorte que, pour tout $v\not\in V$,  $K'_{v}=ad_{x_{v}}(K_{v})$ soit en bonne position relativement \`a $M_{0}$.    On pose
$$A^G_{unip}(V,\omega,K^V)=\omega(x^V)A^G_{unip}(V,\omega,M_{0},K^{'V}).$$
On doit prouver

(5) cette d\'efinition ne d\'epend pas des choix de $x^V$ et $M_{0}$.

Preuve. Montrons d'abord que, pour $M_{0}$ fix\'e, la d\'efinition ne d\'epend pas de $x^V$. Choisissons $y^V=(y_{v})_{v\not\in V}\in G({\mathbb A}_{F}^V)$ tel que, pour tout $v\not\in V$, $K''_{v}=ad_{y_{v}}(K_{v})$ soit en bonne position relativement \`a $M_{0}$. Soit $v\not\in V$. Notons $o'_{v}$ et $o''_{v}$ les points hypersp\'eciaux de l'immeuble de $G_{AD}$ sur $F_{v}$ qui sont fix\'es respectivement par $K'_{v}$ et  $K''_{v}$. Notons ${\cal S}_{v}$ l'ensemble des sous-tores de $M_{0}$ d\'efinis et d\'eploy\'es sur $F_{v}$, qui sont maximaux parmi les sous-tores v\'erifiant ces conditions. Pour $S\in {\cal S}_{v}$, notons $A(S)$ l'appartement de l'immeuble associ\'e \`a $S$. 
 Les deux groupes $K'_{v}$ et $K''_{v}$ sont en bonne position relativement \`a $M_{0}$. Cela signifie que l'on peut fixer des sous-tores $S',S''\in {\cal S}_{v}$  de sorte que $o'_{v}$ appartienne \`a l'appartement  $A(S')$    et que $o''_{v}$ appartienne \`a l'appartement $A(S'')$. 
Posons $g_{v}=y_{v}x_{v}^{-1}$. On a $K''_{v}=ad_{g_{v}}(K'_{v})$. Donc $o''_{v}$ appartient \`a $A(ad_{g_{v}}(S'))$. Puisque $o''_{v}$ est hypersp\'ecial, le fait que $o''_{v}$ appartienne \`a la fois \`a $A(S'')$ et \`a $A(ad_{g_{v}}(S'))$ entra\^{\i}ne que les deux tores $S''$ et $ad_{g_{v}}(S')$ sont conjugu\'es par un \'el\'ement de $K''_{v}$. Fixons un \'el\'ement $h_{v}\in K''_{v}$ tel que $S''=ad_{h_{v}g_{v}}(S')$.  L'ensemble ${\cal S}_{v}$ forme une seule classe de conjugaison par $M_{0}(F_{v})$. Fixons donc $r_{v}\in M_{0}(F_{v})$ de sorte que $S'=ad_{r_{v}}(S'')$. On obtient que   $h_{v}g_{v}r_{v}$ appartient au normalisateur de $S''$ dans $G(F_{v})$. Comme on le sait, tout \'el\'ement du groupe de Weyl de $S''$ a un repr\'esentant dans $K''_{v}$. Cela entra\^{\i}ne que le normalisateur de $S''$ dans $G(F_{v})$ est contenu dans le produit de $K''_{v}$ et du centralisateur de $S''$ dans $G(F_{v})$. Ce dernier groupe est contenu dans $M_{0}(F_{v})$. On obtient que $h_{v}g_{v}r_{v}\in K''_{v}M_{0}(F_{v})$. Puisque $h_{v}\in K''_{v}$ et $r_{v}\in M_{0}(F_{v})$, on a aussi $g_{v}\in K''_{v}M_{0}(F_{v})$. Ecrivons donc $g_{v}=k_{v}m_{v}$, avec $k_{v}\in K''_{v}$ et $m_{v}\in M_{0}(F_{v})$. Remarquons que $K''_{v}=ad_{k_{v}^{-1}}(K''_{v})=ad_{k_{v}^{-1}g_{v}}(K'_{v})=ad_{m_{v}}(K'_{v})$. Dans la relation (4), rempla\c{c}ons $K^V$ par $K^{_{'}V}$ et $x$ par un \'el\'ement $m\in M_{0}({\mathbb A}_{F})$ dont les composantes hors de $V$ sont les $m_{v}$. On obtient l'\'egalit\'e
$$A^G_{unip}(V,\omega,M_{0},K^{_{''}V})=\omega(m^V)^{-1}A^G_{unip}(V,\omega,M_{0},K^{_{'}V}).$$
Les \'egalit\'es $y_{v}x_{v}^{-1}=g_{v}=k_{v}m_{v}$ pour tout $v\not\in V$ entra\^{\i}nent $\omega(m^V)=\omega(y^V)\omega(x^V)^{-1}$. L'\'egalit\'e ci-dessus devient  
 $$\omega(y^V)A^G_{unip}(V,\omega,M_{0},K^{''V})=\omega(x^V)A^G_{unip}(V,\omega,M_{0},K^{'V}).$$
 Cela d\'emontre que notre d\'efinition ne d\'epend pas du choix de $x^V$. 
 Changeons maintenant de Levi $M_{0}$. Cela signifie qu'on le remplace par $ad_{g}(M_{0})$ pour un $g\in G(F)$. D'apr\`es l'ind\'ependance d\'ej\`a prouv\'ee de $x^V$, on peut remplacer $x^V$ par $g^Vx^{V}$, o\`u $g^V$ est la projection de $g$ dans $G({\mathbb A}_{F}^V)$. La relation \`a prouver devient
  $$\omega(x^V)A^G_{unip}(V,\omega,M_{0},K^{'V})=\omega(g^Vx^V)A^G_{unip}(V,\omega,ad_{g}(M_{0}),ad_{g^V}(K^{'V})).$$
  Cela r\'esulte de (4). $\square$
\bigskip

\subsection{Les distributions associ\'ees \`a une classe rationnelle semi-simple}
Revenons au cas o\`u $(G,\tilde{G},{\bf a})$ est quelconque.  Soit ${\cal O}\in \tilde{G}_{ss}(F)/conj$. Fixons $\dot{\gamma}\in {\cal O}$ et une paire de Borel \'epingl\'ee ${\cal E}=(B,T,(E_{\alpha})_{\alpha\in \Delta})$ (d\'efinie sur $\bar{F})$ de sorte que $\dot{\gamma}$ conserve $(B,T)$. Fixons $e\in Z(\tilde{G},{\cal E})$ et posons $\theta=ad_{e}$. En notant $\Sigma(T)$ l'ensemble des racines de $T$ dans $G$, on sait d\'efinir $N\alpha$ pour tout $\alpha\in \Sigma(T)$: c'est la somme des \'el\'ements de l'orbite de $\alpha$ sous l'action du groupe d'automorphismes engendr\'e par $\theta$. On peut \'ecrire $\dot{\gamma}=te$, avec $t\in T$. En fait, on peut fixer une extension finie $E$ de $F$ telle que $t\in T(E)$. On a alors $(N\alpha)(t)\in E^{\times}$ pour tout $\alpha$. On note $S({\cal O})$ le plus petit ensemble de places de $F$ contenant $V_{ram}$ et tel que les propri\'et\'es suivantes soient v\'erifi\'ees. Soit $w$ une place de $E$ au-dessus d'une place $v\not\in S({\cal O})$. Notons $\mathfrak{o}^{\times}_{w}$ le groupe des unit\'es de $E_{w}$ et ${\mathbb E}_{w}$ le corps r\'esiduel. Alors

(1) pour tout $\alpha\in \Sigma(T)$, on a $(N\alpha)(t)\in \mathfrak{o}^{\times}_{w}$;

(2) pour tous $\alpha\in \Sigma(T)$, si $(N\alpha)(t)\not= \pm 1$, alors la r\'eduction dans ${\mathbb E}_{w}$ de $(N\alpha)(t)$ n'est pas \'egale \`a $\pm 1$. 

 On voit que cette d\'efinition ne d\'epend pas des choix auxiliaires, en particulier ne d\'epend pas du choix de $\dot{\gamma}$. Soit $v\not\in S({\cal O})$. On a
 
 (3) supposons que la classe de conjugaison par $G(F_{v})$ engendr\'ee par ${\cal O}$ coupe $\tilde{K}_{v}$; alors l'ensemble des $x\in G(F_{v})$ tels que $x^{-1}\dot{\gamma}x\in \tilde{K}_{v}$ forme une unique double classe dans $G_{\dot{\gamma}}(F_{v})\backslash G(F_{v})/K_{v}$; pour un \'el\'ement $x$ de cet ensemble, le groupe $xK_{v}x^{-1}\cap G_{\dot{\gamma}}(F_{v})$ est un sous-groupe compact hypersp\'ecial de $G_{\dot{\gamma}}(F_{v})$.
 
 Preuve. Il s'agit essentiellement de la proposition 7.1 de [K3], que l'on a reprise dans le cadre tordu en [W2]. L'hypoth\`ese $v\not\in V_{ram}$ permet d'appliquer les r\'esultats de [W2] chapitres 4 et 5. En particulier, on d\'ecrit $\tilde{G}(F_{v})$ et $\tilde{K}_{v}$ comme en [W2] 4.4. Notons simplement $\gamma$ l'image de $\dot{\gamma}$ dans $\tilde{G}(F_{v})$. L'hypoth\`ese que la classe de conjugaison de $\gamma$ coupe $\tilde{K}_{v}$ signifie que $\gamma$ est un \'el\'ement compact dans la terminologie de cette r\'ef\'erence. On peut le d\'ecomposer en $\gamma=\gamma_{tu}\gamma_{p'}$, o\`u $\gamma_{p'}$ est un \'el\'ement d'ordre fini premier \`a la caract\'eristique r\'esiduelle $p$ et $\gamma_{tu}$ est un \'el\'ement topologiquement unipotent de $G(F_{v})$ qui commute \`a $\gamma_{p'}$. Les conditions (1) et (2) ci-dessus entra\^{\i}nent l'\'egalit\'e $G_{\gamma}=G_{\gamma_{p'}}$. Par ailleurs, la $p'$-composante d'un \'el\'ement de $\tilde{K}_{v}$ appartenant aussi \`a $\tilde{K}_{v}$, la condition $x^{-1}\gamma x\in \tilde{K}_{v}$ implique $x^{-1}\gamma_{p'}x\in \tilde{K}_{v}$. Il reste \`a appliquer les lemmes [W2] 5.6(ii) et 5.4(ii). $\square$

 Soit $V$ un ensemble fini de places de $F$ contenant $S({\cal O})$. On va d\'efinir une distribution $A^{\tilde{G}}(V,{\cal O},\omega)\in D_{orb}(\tilde{G}(F_{V}),\omega)\otimes Mes(G(F_{V}))^*$. Elle est nulle sauf si les trois conditions suivantes sont v\'erifi\'ees:
 
 (4) $\omega$ est trivial sur $Z(G,{\mathbb A}_{F})^{\theta}$ et sur $Z(G_{\dot{\gamma}},{\mathbb A}_{F})$ pour $\dot{\gamma}\in {\cal O}$;
 
 (5) ${\cal O}$ est form\'e d'\'el\'ements elliptiques;
 
 (6) pour tout $v\not\in V$, la classe de conjugaison par $G(F_{v})$ engendr\'ee par ${\cal O}$ coupe $\tilde{K}_{v}$.
 
 Supposons ces conditions v\'erifi\'ees. On fixe $\dot{\gamma}\in {\cal O}$. Pour tout $v\not\in V$, on fixe $x_{v}\in G(F_{v})$ tel que $x_{v}^{-1}\dot{\gamma}x_{v}\in \tilde{K}_{v}$. On pose $K_{\dot{\gamma},v}=x_{v}K_{v}x_{v}^{-1}\cap G_{\dot{\gamma}}(F_{v})$. Puisque $\dot{\gamma}\in \tilde{K}_{v}$ pour presque tout $v$, on peut supposer $x_{v}=1$ pour presque tout $v$. Alors l'\'el\'ement $x^V=(x_{v})_{v\not\in V}$ appartient \`a $G({\mathbb A}_{F}^V)$ et la famille $(K_{\dot{\gamma},v})_{v\not\in V}$ de compacts hypersp\'eciaux v\'erifie la condition de compatibilit\'e globale de 1.1. La distribution $A^{G_{\dot{\gamma}}}_{unip}(V,\omega,K_{\dot{\gamma}}^V)$ est bien d\'efinie.  Fixons des mesures $dg_{V}$ sur $G(F_{V})$ et $dh_{V}$ sur $G_{\dot{\gamma}}(F_{V})$. Soit $f_{V}\in C_{c}^{\infty}(\tilde{G}(F_{V}))$, posons ${\bf f}_{V}=f_{V}\otimes dg_{V}$. Pour $y\in G(F_{V})$, notons $^yf_{V\vert G_{\dot{\gamma}}}$ la fonction $h\mapsto f(y^{-1}h\gamma y)$ sur $G_{\dot{\gamma}}(F_{V})$. Posons $^y{\bf f}_{V\vert  G_{\dot{\gamma}}}={^yf}_{V\vert G_{\dot{\gamma}}}\otimes dh_{V}$.  On d\'efinit  la distribution $A^{\tilde{G}}(V,{\cal O},\omega)$  par la formule suivante:
 $$(7) \qquad I^{\tilde{G}}(A^{\tilde{G}}(V,{\cal O},\omega),{\bf f}_{V})=[Z_{G}(\dot{\gamma},F):G_{\dot{\gamma}}(F)]^{-1}\omega(x^V)$$
 $$\int_{G_{\dot{\gamma}}(F_{V})\backslash G(F_{V})}I^{G_{\dot{\gamma}}}(A_{unip}^{G_{\dot{\gamma}}}(V,\omega,K_{\dot{\gamma}}^V),{^y{\bf f}_{V\vert G_{\dot{\gamma}}}})\omega(y)\,dy,$$
o\`u   $dy$ est la mesure $\frac{dg_{V}}{dh_{V}}$ sur $G_{\dot{\gamma}}(F_{V})\backslash G(F_{V})$. On v\'erifie que cela ne d\'epend ni du choix de $\dot{\gamma}$, ni de celui de l'\'el\'ement $x^V$, ni de celui des mesures. Par contre, la distribution $A^{\tilde{G}}(V,{\cal O},\omega)$ d\'epend des $\tilde{K}_{v}$ pour $v\not\in V$.

Cas particulier: supposons ${\cal O}$ form\'ee d'\'el\'ements elliptiques et fortement r\'eguliers dans $\tilde{G}(F)$. Alors la formule (7) se simplifie en
$$I^{\tilde{G}}(A^{\tilde{G}}(V,{\cal O},\omega),{\bf f}_{V})=[Z_{G}(\gamma,F):G_{\gamma}(F)]^{-1}\omega(x^V)mes(\mathfrak{A}_{G}G_{\dot{\gamma}}(F)\backslash G_{\dot{\gamma}}({\mathbb A}_{F}))$$
$$\int_{G_{\dot{\gamma}}(F_{V})\backslash G(F_{V})}f(y^{-1}\dot{\gamma}y)\omega(y)\,dy .$$

Consid\'erons maintenant un ensemble fini $V$ de places, contenant $V_{ram}$ mais pas forc\'ement $S({\cal O})$. On fixe un ensemble fini $S$ de places contenant $V$ et $S({\cal O})$. Rappelons que, pour tout espace de Levi $\tilde{M}$, l'intersection $\tilde{M}(F)\cap {\cal O}$ se d\'ecompose en un nombre fini de classes de conjugaison semi-simples par $M(F)$. 
Pour une telle classe ${\cal O}_{\tilde{M}}$, on a $S({\cal O}_{\tilde{M}})\subset S({\cal O})$. La distribution $A^{\tilde{M}}(S,{\cal O}_{\tilde{M}},\omega)$ (relative \`a $\tilde{K}^{\tilde{M},S}$) est bien d\'efinie. On peut \'ecrire
$$(8) \qquad A^{\tilde{M}}(S,{\cal O}_{\tilde{M}},\omega)=\sum_{i=1,...,n({\cal O}_{\tilde{M}})}k_{i}^{\tilde{M}}({\cal O}_{\tilde{M}},\omega)_{S}^V\otimes A_{i}^{\tilde{M}}({\cal O}_{\tilde{M}},\omega)_{V},$$
o\`u $k_{i}^{\tilde{M}}({\cal O}_{\tilde{M}},\omega)_{S}^V\in D_{g\acute{e}om}(\tilde{M}(F_{S}^V),\omega)\otimes Mes(M(F_{S}^V))^*$ et $A_{i}^{\tilde{M}}({\cal O}_{\tilde{M}},\omega)_{V}\in D_{orb}(\tilde{M}(F_{V}),\omega)\otimes Mes(M(F_{V}))^*$. On note $A_{i}^{\tilde{G}}({\cal O}_{\tilde{M}},\omega)_{V}\in D_{orb}(\tilde{G}(F_{V}),\omega)\otimes Mes(G(F_{V}))^*$ la distribution induite de $A_{i}^{\tilde{M}}({\cal O}_{\tilde{M}},\omega)_{V}$ \`a $\tilde{G}(F_{V})$. On se rappelle l'application lin\'eaire $r_{\tilde{M}}^{\tilde{G}}(.,\tilde{K}_{S}^V)$ d\'efinie sur $D_{g\acute{e}om}(\tilde{M}(F_{S}^V),\omega)\simeq D_{g\acute{e}om}(\tilde{M}(F_{S}^V),\omega)\otimes Mes(M(F_{S}^V))^*$ en 1.13. 
On d\'efinit la distribution $A^{\tilde{G}}(V,{\cal O},\omega) $ par la formule
$$(9) \qquad A^{\tilde{G}}(V,{\cal O},\omega)=\sum_{\tilde{M}\in {\cal L}(\tilde{M}_{0})}\vert W^{\tilde{M}}\vert \vert W^{\tilde{G}}\vert ^{-1}\sum_{{\cal O}_{\tilde{M}}\in \tilde{M}_{ss}(F)/conj, {\cal O}_{\tilde{M}}\subset {\cal O}}$$
$$\sum_{i=1,...,n({\cal O}_{\tilde{M}})}r_{\tilde{M}}^{\tilde{G}}(k_{i}^{\tilde{M}}({\cal O}_{\tilde{M}},\omega)_{S}^V,\tilde{K}_{S}^V)A_{i}^{\tilde{G}}({\cal O}_{\tilde{M}},\omega)_{V}.$$

\ass{Proposition}{La distribution $A^{\tilde{G}}(V,{\cal O},\omega) $ v\'erifie les propri\'et\'es suivantes: 

(i) pour tout ${\bf f}_{V}\in  C_{c}^{\infty}(\tilde{G}(F_{V}))\otimes Mes(G(F_{V}))$, on a l'\'egalit\'e:
$$ J_{{\cal O}}({\bf f}_{V},\omega)=\sum_{\tilde{M}\in {\cal L}(\tilde{M}_{0})}\vert W^{\tilde{M}}\vert \vert W^{\tilde{G}}\vert ^{-1}\sum_{{\cal O}_{\tilde{M}}\in \tilde{M}_{ss}(F)/conj, {\cal O}_{\tilde{M}}\subset {\cal O}}$$
$$J_{\tilde{M}}^{\tilde{G}}(A^{\tilde{M}}(V,{\cal O}_{\tilde{M}},\omega),{\bf f}_{V});$$

(ii) $A^{\tilde{G}}(V,{\cal O},\omega)$ est combinaison lin\'eaire d'int\'egrales orbitales associ\'ees \`a des \'el\'ements $\dot{\gamma}$ (ou plus exactement aux projections dans $\tilde{G}(F_{V})$ de tels \'el\'ements), o\`u $\dot{\gamma}\in \tilde{G}(F) $ est un \'el\'ement dont la partie semi-simple appartient \`a ${\cal O}$;  

(iii) $A^{\tilde{G}}(V,{\cal O},\omega)=0$ sauf si, pour tout $v\not\in V$, la classe de conjugaison par $G(F_{v})$ engendr\'ee par ${\cal O}$ coupe $\tilde{K}_{v}$;

(iv) $A^{\tilde{G}}(V,{\cal O},\omega)=0$ sauf si $\omega$ est trivial sur $Z(G;{\mathbb A}_{F})^{\theta}$ et sur $Z(G_{\dot{\gamma}},{\mathbb A}_{F})$ pour $\dot{\gamma}\in {\cal O}$.}

{\bf Remarque.}    L'\'egalit\'e (i) suffit  \`a caract\'eriser la distribution $A^{\tilde{G}}(V,{\cal O},\omega)$ par r\'ecurrence. Cela prouve en particulier que la d\'efinition (9) ne d\'epend pas de l'ensemble $S$ choisi.

\bigskip

Preuve. On suppose d'abord que $V$ contient $S({\cal O})$ et que $A^{\tilde{G}}(V,{\cal O},\omega) $ est d\'efinie par la formule (7). Les propri\'et\'es (ii), (iii) et (iv) sont imm\'ediates d'apr\`es la d\'efinition et l'analogue de la propri\'et\'e (ii) pour la distribution $A^{G_{\dot{\gamma}}}_{unip}(V,\omega,K_{\dot{\gamma}}^V)$. L'assertion principale (i) est \`a peu pr\`es celle que prouve Arthur dans [A10] th\'eor\`eme 8.1. La seule diff\'erence avec le r\'esultat d'Arthur est  la d\'efinition de l'ensemble $S({\cal O})$. Expliquons ce point. Parmi les espaces de Levi  $\tilde{M}\in {\cal L}(\tilde{M}_{0})$ tels que ${\cal O}\cap \tilde{M}(F)\not=\emptyset$, fixons un \'el\'ement minimal $\tilde{M}_{1}$. Au lieu de la condition (6), Arthur impose qu'il existe un \'el\'ement $\dot{\gamma}\in {\cal O}\cap \tilde{M}_{1}(F)$ tel que $\dot{\gamma}\in \tilde{K}_{v}$ pour tout $v\not\in V$. L'existence d'un tel \'el\'ement d\'efini sur $F$ ne r\'esulte pas de (6). Mais on peut modifier tr\`es l\'eg\`erement la preuve d'Arthur pour obtenir notre assertion. En effet, fixons $\dot{\gamma}\in  {\cal O}\cap \tilde{M}_{1}(F)$. On a

(10) l'hypoth\`ese (6) entra\^{\i}ne que, pour tout $v\not\in V$, la classe de conjugaison par $M_{1}(F_{v})$ engendr\'ee par $\dot{\gamma}$ coupe $\tilde{K}^{\tilde{M}_{1}}_{v}$.

Il y a en tout cas un $x\in G(F_{v})$ tel que $x^{-1}\dot{\gamma}x\in \tilde{K}_{v}$. Fixons $\tilde{P}_{1}=\tilde{M}_{1}U_{P_{1}}\in {\cal P}(\tilde{M}_{1})$. Gr\^ace \`a la d\'ecomposition d'Iwasawa, on \'ecrit $x=muk$, avec $m\in M_{1}(F_{v})$, $u\in U_{P_{1}}(F_{v})$, $k\in K_{v}$. On a encore $(mu)^{-1}\dot{\gamma}mu\in \tilde{K}_{v}$. Mais $(mu)^{-1}\dot{\gamma}mu=u'\gamma'$, avec $u'\in U_{P_{1}}(F_{v})$ et $\gamma'=m^{-1}\dot{\gamma}m\in \tilde{M}_{1}(F_{v})$. Parce que $K_{v}$ est en bonne position relativement \`a $M_{1}$, cela entra\^{\i}ne $\gamma'\in \tilde{K}_{v}^{\tilde{M}_{1}}$. D'o\`u (10).

Gr\^ace \`a (10), on peut construire $x^V\in M_{1}({\mathbb A}_{F}^V)$ tel que $(x^V)^{-1}\dot{\gamma}x^V\in \tilde{K}^{\tilde{M}_{1},V}$. Dans la preuve de [A10], l'\'el\'ement $y'$ des derni\`eres lignes de la page 203 appartient alors, non pas \`a $K_{\dot{\gamma}}^V\backslash K^V$, mais \`a $K_{\dot{\gamma}}^V\backslash x^VK^V$. Mais un \'el\'ement de $x^VK^V$ ne contribue pas aux fonctions $v'_{R}$ de la page 204 parce que ces fonctions sont invariantes \`a droite par $K^V$ et \`a gauche par $M_{1}({\mathbb A}_{F})$, puisque les Levi $R$ intervenant contiennent $M_{1}$. Ces \'el\'ements disparaissent et la suite de la d\'emonstration d'Arthur s'applique sans changement. Il reste toutefois le terme $\omega(x^V)$ qui n'apparaissait pas dans Arthur simplement parce que celui-ci traitait le cas $\omega=1$. Cela prouve la proposition dans le cas o\`u $V$ contient $S({\cal O})$.

Levons cette hypoth\`ese et d\'efinissons $A^{\tilde{G}}(V,{\cal O},\omega)$ par la formule (9). Les propri\'et\'es (ii), (iii) et (iv) sont claires d'apr\`es  cette d\'efinition et les m\^emes propri\'et\'es de  la distribution $A^{\tilde{G}}(S,{\cal O},\omega)$.  Il faut v\'erifier (i). Soit ${\bf f}_{V}\in C_{c}^{\infty}(\tilde{G}(F_{V}))\otimes Mes(G(F_{V}))$. Notons ${\bf f}_{S}\in C_{c}^{\infty}(\tilde{G}(F_{S}))\otimes Mes(G(F_{S}))$ l'\'el\'ement auquel ${\bf f}_{V}$ s'identifie, c'est-\`a-dire ${\bf f}_{S}={\bf f}_{V}\otimes {\bf 1}_{\tilde{K}_{S}^V}$, o\`u le dernier terme est tacitement tensoris\'e avec la mesure canonique sur $G(F_{S}^V)$. On conna\^{\i}t le d\'eveloppement (i) pour $S$:
$$(11)\qquad  J_{{\cal O}}({\bf f}_{S},\omega)=\sum_{\tilde{M}\in {\cal L}(\tilde{M}_{0})}\vert W^{\tilde{M}}\vert \vert W^{\tilde{G}}\vert ^{-1}\sum_{{\cal O}_{\tilde{M}}\in \tilde{M}_{ss}(F)/conj, {\cal O}_{\tilde{M}}\subset {\cal O}}$$
$$J_{\tilde{M}}^{\tilde{G}}(A^{\tilde{M}}(S,{\cal O}_{\tilde{M}},\omega),{\bf f}_{S}).$$
Notons que $ J_{{\cal O}}({\bf f}_{S},\omega)= J_{{\cal O}}({\bf f}_{V},\omega)$ par d\'efinition. Soient $\tilde{M}$ et ${\cal O}_{\tilde{M}}$ intervenant ci-dessus. Pour tout espace $\tilde{Q}=\tilde{L}U_{Q}\in {\cal F}(\tilde{M})$, la fonction $({\bf 1}_{\tilde{K}_{S}^V})_{\tilde{Q},\omega}$ est \'egale \`a ${\bf 1}_{\tilde{K}^{\tilde{L},V}_{S}}$. Elle ne d\'epend donc que de $\tilde{L}$. En utilisant (8), la   relation 1.9(2) et la d\'efinition de $r_{\tilde{M}}^{\tilde{G}}(.,\tilde{K}_{S}^V)$, on obtient l'\'egalit\'e
$$J_{\tilde{M}}^{\tilde{G}}(A^{\tilde{M}}(S,{\cal O}_{\tilde{M}},\omega),{\bf f}_{S})=\sum_{i=1,...,n({\cal O}_{\tilde{M}})}\sum_{\tilde{L}\in {\cal L}(\tilde{M})}r_{\tilde{M}}^{\tilde{L}}(k_{i}^{\tilde{M}}({\cal O}_{\tilde{M}},\omega)_{S}^V,\tilde{K}_{S}^V)J_{\tilde{L}}^{\tilde{G}}(A_{i}^{\tilde{L}}({\cal O}_{\tilde{M}},\omega)_{V},{\bf f}_{V}).$$
La formule (11) se r\'ecrit
$$(12) \qquad J_{{\cal O}}({\bf f}_{V},\omega)=\sum_{\tilde{L}\in {\cal L}(\tilde{M}_{0})}\vert W^{\tilde{L}}\vert \vert W^{\tilde{G}}\vert ^{-1}\sum_{\tilde{M}\in {\cal L}^{\tilde{L}}(\tilde{M}_{0})}\vert W^{\tilde{M}}\vert \vert W^{\tilde{L}}\vert ^{-1}\sum_{{\cal O}_{\tilde{M}}\in \tilde{M}_{ss}(F)/conj, {\cal O}_{\tilde{M}}\subset {\cal O}}$$
$$\sum_{i=1,...,n({\cal O}_{\tilde{M}})}r_{\tilde{M}}^{\tilde{L}}(k_{i}^{\tilde{M}}({\cal O}_{\tilde{M}},\omega)_{S}^V,\tilde{K}_{S}^V)J_{\tilde{L}}^{\tilde{G}}(A_{i}^{\tilde{L}}({\cal O}_{\tilde{M}},\omega)_{V},{\bf f}_{V}).$$
On peut d\'ecomposer la somme en ${\cal O}_{\tilde{M}}$ en une somme sur ${\cal O}_{\tilde{L}}\in \tilde{L}_{ss}(F)/conj$ tel que ${\cal O}_{\tilde{L}}\subset {\cal O}$ et une somme sur ${\cal O}_{\tilde{M}}\in \tilde{M}_{ss}(F)/conj$ tel que ${\cal O}_{\tilde{M}}\subset {\cal O}_{\tilde{L}}$. La contribution d'un couple $(\tilde{L},{\cal O}_{\tilde{L}})$ est alors
$$\sum_{\tilde{M}\in {\cal L}^{\tilde{L}}(\tilde{M}_{0})}\vert W^{\tilde{M}}\vert \vert W^{\tilde{L}}\vert ^{-1}\sum_{{\cal O}_{\tilde{M}}\in \tilde{M}_{ss}(F)/conj, {\cal O}_{\tilde{M}}\subset {\cal O}_{\tilde{L}}}$$
$$\sum_{i=1,...,n({\cal O}_{\tilde{M}})}r_{\tilde{M}}^{\tilde{L}}(k_{i}^{\tilde{M}}({\cal O}_{\tilde{M}},\omega)_{S}^V,\tilde{K}_{S}^V)J_{\tilde{L}}^{\tilde{G}}(A_{i}^{\tilde{L}}({\cal O}_{\tilde{M}},\omega)_{V},{\bf f}_{V}).$$
Ceci qui n'est autre que 
$J_{\tilde{L}}^{\tilde{G}}(A^{\tilde{L}}(V,{\cal O}_{\tilde{L}},\omega),{\bf f}_{V})$
par d\'efinition de $A^{\tilde{L}}(V,{\cal O}_{\tilde{L}},\omega)$. Alors la formule (12) devient l'\'egalit\'e (i) de l'\'enonc\'e. $\square$

{\bf Remarque.}  La m\^eme preuve montre que l'\'egalit\'e (9) est vraie pour tous $S$, $V$ tels que $V_{ram}\subset V\subset S$.

\bigskip

\subsection{D\'eveloppement de la partie g\'eom\'etrique de la formule des traces non invariante}

Soit ${\bf f}\in C_{c}^{\infty}(\tilde{G}({\mathbb A}_{F})\otimes Mes(G({\mathbb A}_{F})$. Si $V$ est un ensemble fini de places assez grand, on peut \'ecrire ${\bf f}={\bf f}_{V}\otimes {\bf 1}_{\tilde{K}^V}$, o\`u ${\bf f}_{V}\in C_{c}^{\infty}(\tilde{G}(F_{V}))\otimes Mes(G(F_{V}))$ et ${\bf 1}_{\tilde{K}^V}$ est tacitement tensoris\'e avec la mesure canonique sur $G({\mathbb A}_{F}^V)$.  En vertu de 2.1(3) et de la proposition 2.3(i),  on a l'\'egalit\'e
$$(1) \qquad J^{\tilde{G}}_{g\acute{e}om}({\bf f},\omega)=\sum_{\tilde{M}\in {\cal L}(\tilde{M}_{0})}\vert W^{\tilde{M}}\vert \vert W^{\tilde{G}}\vert ^{-1}\sum_{{\cal O}\in \tilde{M}_{ss}(F)/conj} J_{\tilde{M}}^{\tilde{G}}(A^{\tilde{M}}(V,{\cal O},\omega),{\bf f}_{V}).$$
Le lemme 2.1 et sa preuve montrent que cette somme est finie ind\'ependamment de $V$. Cette finitude et les d\'efinitions  entra\^{\i}nent que, pour ${\bf f}$ fix\'e, il existe un ensemble $S({\bf f})$ tel que, pour $V\supset S({\bf f})$, les seuls couples $(\tilde{M},{\cal O})$ qui contribuent \`a la formule  v\'erifient  $S({\cal O})\subset V$ et ${\cal O}\subset \tilde{M}(F)_{ell}$. Notons que l'ensemble fini des indices contribuant \`a la formule (1) ainsi que l'ensemble $S({\bf f})$ peuvent \^etre choisis ind\'ependants de ${\bf f}$ si on se limite \`a des fonctions dont le support est contenu dans un compact fix\'e.

   \bigskip

\subsection{Variante avec caract\`ere central}
On suppose $(G,\tilde{G},{\bf a})$ quasi-d\'eploy\'e et \`a torsion int\'erieure. On suppose donn\'ee une extension
$$1\to C_{1}\to G_{1}\to G\to 1,$$
 une extension compatible $\tilde{G}_{1}\to \tilde{G}$, o\`u $\tilde{G}_{1}$ est \`a torsion int\'erieure,  et un caract\`ere automorphe $\lambda_{1}$ de $C_{1}({\mathbb A}_{F})$. On introduit des donn\'ees compl\'ementaires comme en 1.15. On a une suite exacte
 $$0\to \mathfrak{A}_{C_{1}}\to \mathfrak{A}_{G_{1}}\to \mathfrak{A}_{G}\to 0.$$
 On suppose des mesures choisies sur ces espaces de fa\c{c}on compatible \`a cette suite. On fait de m\^eme quand on remplace $\tilde{G}$ par un espace de Levi $\tilde{M}$. Notons que cette hypoth\`ese ne serait pas v\'erifi\'ee si on normalisait les mesures "\`a la Tamagawa", cf. 1.3.

 Fixons des mesures $dg$ sur $G({\mathbb A}_{F})$ et $dc$ sur $C_{1}({\mathbb A}_{F})$. On en d\'eduit une mesure $dg_{1}$ sur $G_{1}({\mathbb A}_{F})$ compatible \`a la suite exacte
 $$1\to C_{1}({\mathbb A}_{F})\to G_{1}({\mathbb A}_{F})\to G({\mathbb A}_{F})\to 1.$$
  Soit $f\in C_{c,\lambda_{1}}^{\infty}(\tilde{G}_{1}({\mathbb A}_{F}))$. On peut fixer une fonction $\phi\in C_{c}^{\infty}(\tilde{G}_{1}({\mathbb A}_{F}))$ telle que
$$f(\gamma_{1})=\int_{C_{1}({\mathbb A}_{F})}\phi^c(\gamma_{1})\lambda_{1}(c)\,dc,$$
o\`u $\phi^c(\gamma_{1})=\phi(c \gamma_{1})$.  On pose
$$(1) \qquad J_{g\acute{e}om,\lambda_{1}}^{\tilde{G}_{1}}(f\otimes dg)=mes(\mathfrak{A}_{C_{1}}C_{1}(F)\backslash C_{1}({\mathbb A}_{F}))^{-1}\int_{C_{1}(F)\backslash C_{1}({\mathbb A}_{F})}J^{\tilde{G}_{1}}_{g\acute{e}om}(\phi^c\otimes dg_{1})\lambda_{1}(c)\,dc.$$
On v\'erifie en effet, d'une part que la fonction \`a int\'egrer est invariante par $C_{1}(F)$ et \`a support compact modulo ce groupe, d'autre part que l'expression ci-dessus ne d\'epend  que de $f$ et $dg$ et  pas des choix de $\phi$ et  de la mesure $dc$. De m\^eme, elle ne d\'epend que de la mesure sur $\mathfrak{A}_{\tilde{G}}$ et pas de celle sur $\mathfrak{A}_{C_{1}}$. Ces v\'erifications se font en exprimant comme en 2.1 les diverses expressions comme les valeurs en $T_{0}$ d'un polyn\^ome asymptote \`a des expressions qui, elles, v\'erifient \'evidemment ces propri\'et\'es.  Remarquons qu'il y a une surjection naturelle
$$(2) \qquad \tilde{G}_{1,ss}(F)/conj\to \tilde{G}_{ss}(F)/conj.$$
Pour ${\cal O}\in \tilde{G}_{ss}(F)/conj$, on pose
$$(3) \qquad J_{{\cal O},\lambda_{1}}^{\tilde{G}_{1}}(f\otimes dg)=mes(\mathfrak{A}_{C_{1}}C_{1}(F)\backslash C_{1}({\mathbb A}_{F}))^{-1}\int_{C_{1}(F)\backslash C_{1}({\mathbb A}_{F})}$$
$$\sum_{{\cal O}_{1}\in Fib^{\tilde{G}_{1}}({\cal O})}J^{\tilde{G}_{1}}_{{\cal O}_{1}}(\phi^c\otimes dg_{1})\lambda_{1}(c)\,dc,$$
o\`u  $Fib^{\tilde{G}_{1}}({\cal O})$ est la fibre de l'application (2) au-dessus de ${\cal O}$. De nouveau, ce terme ne d\'epend  que de $f$ et $dg$ et  pas des choix de $\phi$ et  de la mesure $dc$.

Pour un ensemble fini $V$ de places contenant $V_{1,ram}$, on peut encore identifier $C_{c,\lambda_{1}}^{\infty}(\tilde{G}_{1}(F_{V}))\otimes Mes(G(F_{V}))$ \`a un sous-espace de $C_{c,\lambda_{1}}^{\infty}(\tilde{G}_{1}({\mathbb A}_{F}))\otimes Mes(G({\mathbb A}_{F}))$. Ainsi, $J_{g\acute{e}om,\lambda_{1}}^{\tilde{G}_{1}}({\bf f}_{V})$ est d\'efini pour ${\bf f}_{V}\in C_{c,\lambda_{1}}^{\infty}(\tilde{G}_{1}(F_{V}))\otimes Mes(G(F_{V}))$. De m\^eme, $J_{{\cal O},\lambda_{1}}^{\tilde{G}_{1}}({\bf f}_{V})$ est d\'efini. Pour ${\cal O}\in \tilde{G}_{ss}(F)/conj$, on va d\'efinir une distribution $A^{\tilde{G}_{1}}_{\lambda_{1}}(V,{\cal O})\in D_{orb,\lambda_{1}}(\tilde{G}_{1}(F_{V}))\otimes Mes( G(F_{V}))^*$ qui v\'erifie les propri\'et\'es suivantes:

(4) pour ${\bf f}_{V}\in C_{c,\lambda_{1}}^{\infty}(\tilde{G}_{1}(F_{V}))\otimes Mes(G(F_{V}))$,
$$I^{\tilde{G}_{1}}_{\lambda_{1}}(A^{\tilde{G}_{1}}_{\lambda_{1}}(V,{\cal O}),{\bf f}_{V})=J_{{\cal O},\lambda_{1}}^{\tilde{G}_{1}}({\bf f}_{V})$$
$$-\sum_{M\in {\cal L}(M_{0}), M\not=G}\vert W^M\vert \vert W^G\vert ^{-1}\sum_{{\cal O}_{\tilde{M}}\in \tilde{M}_{ss}(F)/conj,{\cal O}_{\tilde{M}}\subset {\cal O}}J_{\tilde{M}_{1},\lambda_{1}}^{\tilde{G}_{1}}(A^{\tilde{M}}_{\lambda_{1}}(V,{\cal O}_{\tilde{M}}),{\bf f}_{V});$$

(5) $A^{\tilde{G}_{1}}_{\lambda_{1}}(V,{\cal O})$ est combinaison lin\'eaire d'int\'egrales orbitales associ\'ees \`a des \'el\'ements $\gamma_{1}$ (ou plus exactement aux projections dans $\tilde{G}_{1}(F_{V})$ de tels \'el\'ements), o\`u $\gamma_{1}\in \tilde{G}_{1}(F)$ se projette en un \'el\'ement  de $\tilde{G}(F)$ dont la partie semi-simple appartient \`a ${\cal O}$;

(6) $A^{\tilde{G}_{1}}_{\lambda_{1}}(V,{\cal O})=0$ sauf si, pour tout $v\not\in V$, la classe de conjugaison par $G(F_{v})$ engendr\'ee par ${\cal O}$ coupe $\tilde{K}_{v}$.

Comme en 2.2, la propri\'et\'e (4) caract\'erise $A^{\tilde{G}_{1}}_{\lambda_{1}}(V,{\cal O})$. Signalons que cette distribution d\'epend de $\tilde{K}_{1}^V$.

Soit ${\bf f}_{V}\in C_{c,\lambda_{1}}^{\infty}(\tilde{G}_{1}(F_{V}))\otimes Mes(G(F_{V}))$, que l'on \'ecrit ${\bf f}_{V}=f_{V}\otimes dg$. On fixe $\phi_{V}\in C_{c}^{\infty}(\tilde{G}_{1}(F_{V}))$ de sorte que
$$f_{V}=\int_{C_{1}(F_{V})}\phi_{V}^c\,dc.$$
On se rappelle que les distributions $A^{\tilde{G}_{1}}(V,{\cal O}_{1})$ d\'ependent du choix de $\tilde{K}_{1}^V$. Notons-les $A^{\tilde{G}_{1}}(V,{\cal O}_{1},\tilde{K}_{1}^V)$. Soit $c\in C_{1}({\mathbb A}_{F})$, que l'on d\'ecompose en $c=c_{V}c^V$.  Consid\'erons le terme
$$ I^{\tilde{G}_{1}}(A^{\tilde{G}_{1}}(V,{\cal O}_{1},c^V\tilde{K}_{1}^V),\phi_{V}^{c_{V}}\otimes dg_{1})$$
pour ${\cal O}_{1}\in Fib^{\tilde{G}_{1}}({\cal O})$. Soit $\Gamma$ un sous-ensemble compact de $C_{1}({\mathbb A}_{F})$. Le lemme 2.1(ii)  entra\^{\i}ne qu'il existe un sous-ensemble fini $\Delta$ de $Fib^{\tilde{G}_{1}}({\cal O})$ tel que, pour tout ${\cal O}_{1}\in Fib^{\tilde{G}_{1}}({\cal O})-\Delta$ et tout $c\in \mathfrak{A}_{C_{1}}\Gamma$, le terme ci-dessus soit nul. Le (i) du m\^eme lemme montre que, pour ${\cal O}_{1}\in \Delta$, ce terme est une fonction de $c$ \`a support compact dans $\mathfrak{A}_{C_{1}}\Gamma$. C'est aussi une fonction lisse de $c$. Les m\^emes propri\'et\'es valent donc pour la somme
$$\sum_{{\cal O}_{1}\in Fib^{\tilde{G}_{1}}({\cal O})} I^{\tilde{G}_{1}}(A^{\tilde{G}_{1}}(V,{\cal O}_{1},c^V\tilde{K}_{1}^V),\phi_{V}^{c_{V}}\otimes dg_{1}).$$
Comme fonction de $c$, celle-ci est invariante par $C_{1}(F)$. En effet, on voit facilement que, pour $\xi\in C_{1}(F)$ et ${\cal O}_{1}\in Fib^{\tilde{G}_{1}}({\cal O})$, on a l'\'egalit\'e
$$(7)\qquad  I^{\tilde{G}_{1}}(A^{\tilde{G}_{1}}(V,\xi{\cal O}_{1},\xi^Vc^V\tilde{K}_{1}^V),\phi_{V}^{\xi_{V}c_{V}}\otimes dg_{1})= I^{\tilde{G}_{1}}(A^{\tilde{G}_{1}}(V,{\cal O}_{1},c^V\tilde{K}_{1}^V),\phi_{V}^{c_{V}}\otimes dg_{1}).$$
L'invariance requise en r\'esulte. En se rappelant que l'on peut choisir $\Gamma$ de sorte que $C_{1}({\mathbb A}_{F})=C_{1}(F)\mathfrak{A}_{C_{1}}\Gamma$, on voit que l'int\'egrale
$$(8) \qquad \int_{C_{1}(F)\backslash C_{1}({\mathbb A}_{F})}\sum_{{\cal O}_{1}\in Fib^{\tilde{G}_{1}}({\cal O})} I^{\tilde{G}_{1}}(A^{\tilde{G}_{1}}(V,{\cal O}_{1},c^V\tilde{K}_{1}^V),\phi_{V}^{c_{V}}\otimes dg_{1})\lambda_{1}(c)\,dc$$
est convergente. On va montrer qu'elle ne d\'epend pas du choix de $\phi_{V}$ et d\'efinit donc une forme lin\'eaire en ${\bf f}_{V}$ qui appartient \`a $D_{g\acute{e}om,\lambda_{1}}(\tilde{G}_{1}(F_{V}))\otimes Mes( G(F_{V}))^*$. Tout d'abord une distribution $A^{\tilde{G}_{1}}(V,{\cal O}_{1},c^V\tilde{K}_{1}^V)$ n'est non nulle que si la condition suivante est v\'erifi\'ee:

(9)  pour tout $v\not\in V$, la classe de conjugaison par $G_{1}(F_{v})$ engendr\'ee par ${\cal O}_{1}$ coupe $c_{v}\tilde{K}_{1,v}$. 

Par projection dans $\tilde{G}(F_{v})$ cela entra\^{\i}ne que la classe de conjugaison par $G(F_{v})$ engendr\'ee par ${\cal O}$ coupe $\tilde{K}_{v}$. On obtient

(10) l'expression (8) est nulle s'il existe $v\not\in V$ telle que la classe de conjugaison par $G(F_{v})$ engendr\'ee par ${\cal O}$ ne coupe pas $\tilde{K}_{v}$.

Supposons que, pour tout $v\not\in V$,  la classe de conjugaison par $G(F_{v})$ engendr\'ee par ${\cal O}$ coupe $\tilde{K}_{v}$.
 Posons $\Xi=C_{1}(F)\cap (K_{C_{1}}^VC_{1}(F_{V}))$. On va montrer

(11) l'ensemble des $c\in C_{1}({\mathbb A}_{F})$ tels qu'il existe ${\cal O}_{1}\in Fib^{\tilde{G}_{1}}({\cal O})$ v\'erifiant (9) est une unique classe modulo $C_{1}(F)K_{C_{1}}^VC_{1}(F_{V})$;

(12) pour $c$ dans cette classe, l'ensemble des ${\cal O}_{1}\in Fib^{\tilde{G}_{1}}({\cal O})$ v\'erifiant (9) est une unique classe sous l'action de $\Xi$ par multiplication.

On d\'emontre d'abord

(13) soient $v\not\in V_{ram}(\tilde{G}_{1})$, $\gamma_{1}\in \tilde{K}_{1,v}$ et $\xi\in C_{1}(F_{v})$; alors  $\xi\gamma_{1}$ est conjugu\'e \`a un \'el\'ement de $\tilde{K}_{1,v}$ si et seulement si $\xi\in K_{C_{1},v}$. 

Preuve de (13). L'assertion "si" est \'evidente. Dans l'autre sens, soient $\gamma_{2}\in \tilde{K}_{1,v}$ et $g\in G_{1}(F_{v})$, supposons $\xi\gamma_{1}=g^{-1}\gamma_{2}g$. Ecrivons $\gamma_{2}=k\gamma_{1}$, avec $k\in K_{1,v}$, et posons $\theta=ad_{\gamma_{1}}$. Alors $\xi=g^{-1}k\theta(g)$. Notons $X^*(G_{1})$ et $X^*(C_{1})$ les groupes des caract\`eres alg\'ebriques  de $G_{1}$ et $C_{1}$. Pour $\chi\in X^*(G_{1})^{\Gamma_{F_{v}}}$, le caract\`ere continu $x\mapsto \vert \chi(x)\vert _{F_{v}}$ de $G_{1}(F_{v})$ est trivial sur $K_{1,v}$ puisque ce groupe est compact. Il est invariant par $G_{AD}(F_{v})$ parce que $\chi$ l'est. La torsion \'etant int\'erieure, $\theta$ est l'automorphisme associ\'e \`a un \'el\'ement de ce groupe $G_{AD}(F_{v})$. Donc le caract\`ere pr\'ec\'edent est invariant par $\theta$. On en d\'eduit qu'il vaut $1$ sur $\xi$. On sait que l'application de restriction
$$X^*(G_{1})^{\Gamma_{F_{v}}}\to X^*(C_{1})^{\Gamma_{F_{v}}}$$
a un conoyau fini. Puisque les caract\`eres pr\'ec\'edents sont \`a valeurs dans ${\mathbb R}_{>0}$, on en d\'eduit que $\vert \chi(\xi)\vert _{F_{v}}=1$ pour tout $\chi\in X^*(C_{1})^{\Gamma_{F_{v}}}$.  Cette relation entra\^{\i}ne que que $\xi$ appartient au sous-groupe compact maximal de $C_{1}(F_{v})$, qui n'est autre que $K_{C_{1},v}$. Cela prouve (13).

D\'emontrons (11). Il est clair que l'ensemble des $c$ en question est invariant par $C_{1}(F)K_{C_{1}}^VC_{1}(F_{V})$. Il est non vide. En effet, fixons une classe  ${\cal O}_{1}\in Fib^{\tilde{G}_{1}}({\cal O})$. L'hypoth\`ese sur ${\cal O}$ implique que, pour tout $v\not\in V$,  la classe de conjugaison par $G_{1}(F_{v})$ engendr\'ee par ${\cal O}_{1}$ coupe l'image r\'eciproque de $\tilde{K}_{v}$ dans $\tilde{G}_{1}(F_{v})$, c'est-\`a-dire $C_{1}(F_{v})\tilde{K}_{1,v}$. On peut donc choisir $c_{v}\in C_{1}(F_{v})$ tel que cette classe coupe $c_{v}\tilde{K}_{1,v}$. Il est clair que l'on peut choisir $c_{v}=1$ pour presque tout $v$. Les $c_{v}$ se regroupent alors en un \'el\'ement $c\in C_{1}({\mathbb A}_{F})$ qui v\'erifie la condition requise. Soient maintenant $c,c'$ dans l'ensemble en question. On choisit ${\cal O}_{1}$ v\'erifiant (9) et ${\cal O}'_{1}$ v\'erifiant son analogue pour $c'$. Il existe $\xi\in C_{1}(F)$ tel que $\xi{\cal O}'_{1}={\cal O}_{1}$. Alors ${\cal O}_{1}$ v\'erifie (9) pour $c$ comme pour $\xi c'$. L'assertion (13) entra\^{\i}ne que $c^V\in (\xi c')^VK_{C_{1}}^V$. Cela entra\^{\i}ne que $c\in c'C_{1}(F)K_{C_{1}}^VC_{1}(F_{V})$. D'o\`u (11).

Le m\^eme calcul prouve (12): pour $c=c'$, la condition finale devient $\xi^V\in K_{C_{1}}^V$, c'est-\`a-dire $\xi\in \Xi$.

 Fixons un repr\'esentant $\zeta$ de la classe d\'efinie par (11), que l'on peut choisir dans $C_{1}({\mathbb A}_{F}^V)$. Dans la formule (8), on peut remplacer l'int\'egrale sur $C_{1}(F)\backslash C_{1}({\mathbb A}_{F})$ par une int\'egrale sur $C_{1}(F)\backslash \zeta C_{1}(F)K_{C_{1}}^VC_{1}(F_{V})$. La fonction que l'on int\`egre est invariante par $K_{C_{1}}^V$. Cela permet, modulo translation par $\zeta$, de remplacer cette derni\`ere int\'egrale par une int\'egrale sur $\Xi_{V}\backslash C_{1}(F_{V})$, o\`u $\Xi_{V}$ est la projection de $\Xi$ dans $C_{1}(F_{V})$.  
  L'int\'egrale (8) est donc \'egale \`a
 $$\int_{\Xi_{V}\backslash C_{1}(F_{V})}\sum_{{\cal O}_{1}\in Fib^{\tilde{G}_{1}}({\cal O})} I^{\tilde{G}_{1}}(A^{\tilde{G}_{1}}(V,{\cal O}_{1},\zeta\tilde{K}_{1}^V),\phi_{V}^{c_{V}}\otimes dg_{1})\lambda_{1}(\zeta c_{V})\,dc_{V}.$$
 D'apr\`es (12), l'ensemble des ${\cal O}_{1}$ qui contribuent \`a cette formule est une unique classe sous $\Xi$. Fixons un \'el\'ement ${\cal O}_{1}$ de cette classe. 
  L'action de $\Xi$ n'est pas libre en g\'en\'eral, il y a un noyau fini. Notons $d$ le nombre d'\'el\'ements de ce noyau. Alors la formule pr\'ec\'edente devient
 $$d^{-1}\int_{\Xi_{V}\backslash C_{1}(F_{V})}\sum_{\xi\in \Xi}I^{\tilde{G}_{1}}(A^{\tilde{G}_{1}}(V,\xi{\cal O}_{1},\zeta\tilde{K}_{1}^V),\phi_{V}^{c_{V}}\otimes dg_{1})\lambda_{1}(\zeta c_{V})\,dc_{V}.$$
 En utilisant (7), on obtient
 $$d^{-1}\int_{\Xi_{V}\backslash C_{1}(F_{V})}\sum_{\xi\in \Xi_{V}}I^{\tilde{G}_{1}}(A^{\tilde{G}_{1}}(V,{\cal O}_{1},\zeta\tilde{K}_{1}^V),\phi_{V}^{\xi c_{V}}\otimes dg_{1})\lambda_{1}(\zeta c_{V})\,dc_{V}.$$
On a aussi $\lambda_{1}(\xi)=1$ pour tout $\xi\in \Xi_{V}$ et l'\'egalit\'e pr\'ec\'edente se r\'ecrit
$$(14) \qquad d^{-1}\int_{C_{1}(F_{V})}I^{\tilde{G}_{1}}(A^{\tilde{G}_{1}}(V,{\cal O}_{1},\zeta\tilde{K}_{1}^V),\phi_{V}^{c_{V}}\otimes dg_{1})\lambda_{1}(\zeta c_{V})\,dc_{V}$$
$$=d^{-1}\lambda_{1}(\zeta)I^{\tilde{G}_{1}}_{\lambda_{1}}(A^{\tilde{G}_{1}}(V,{\cal O}_{1},\zeta\tilde{K}_{1}^V),f_{V}\otimes dg_{1}),$$
o\`u, dans cette derni\`ere \'egalit\'e, $A^{\tilde{G}_{1}}(V,{\cal O}_{1},\zeta\tilde{K}_{1}^V)$ d\'esigne l'image de cette distribution dans $D_{orb,\lambda_{1}}(\tilde{G}_{1}(F_{V}))\otimes Mes(G_{1}(F_{V}))^*$.  Cela d\'emontre que l'int\'egrale (8) a les propri\'et\'es voulues.  Par ailleurs, on voit ais\'ement qu'en multipliant cette int\'egrale par $mes(\mathfrak{A}_{C_{1}}C_{1}(F)\backslash Z({\mathbb A}_{F}))^{-1}$, elle ne d\'epend que de la mesure sur $G(F_{V})$ et pas de celle sur $C_{1}({\mathbb A}_{F})$. On peut donc d\'efinir  une distribution $A^{\tilde{G}}_{\lambda_{1}}(V,{\cal O})\in D_{orb,\lambda_{1}}(\tilde{G}_{1}(F_{V}))\otimes Mes( G(F_{V}))^*$ par l'\'egalit\'e
$$(15) \qquad I^{\tilde{G}_{1}}_{\lambda_{1}}(A^{\tilde{G}}_{\lambda_{1}}(V,{\cal O}),{\bf f}_{V})=mes(\mathfrak{A}_{C_{1}}C_{1}(F)\backslash C_{1}({\mathbb A}_{F}))^{-1} \int_{C_{1}(F)\backslash C_{1}({\mathbb A}_{F})}$$
$$\sum_{{\cal O}_{1}\in Fib^{\tilde{G}_{1}}({\cal O})} I^{\tilde{G}_{1}}(A^{\tilde{G}_{1}}(V,{\cal O}_{1},c^V\tilde{K}_{1}^V),\phi_{V}^{c_{V}}\otimes dg_{1})\lambda_{1}(c)\,dc.$$
Elle v\'erifie la propri\'et\'e (5). La propri\'et\'e (6) r\'esulte de (10). Il reste \`a v\'erifier (4).  On conserve les termes ${\bf f}_{V}$, $f_{V}$ et $\phi_{V}$ ci-dessus. On compl\`ete $f_{V}$ en $f=f_{V}\otimes {\bf 1}_{\tilde{K}_{1}^V,\lambda_{1}}$ et $\phi_{V}$ en $\phi=\phi_{V}\otimes {\bf 1}_{\tilde{K}_{1}^V}$. On a alors
$$f=\int_{C_{1}({\mathbb A}_{F})}\phi^c\lambda_{1}(c)\,dc$$
et l'\'egalit\'e (3). Pour $c\in C_{1}({\mathbb A}_{F})$ et ${\cal O}_{1}\in Fib^{\tilde{G}_{1}}({\cal O})$, on peut d\'evelopper $J_{{\cal O}_{1}}^{\tilde{G}_{1}}(\phi^c\otimes dg_{1})$ selon l'\'egalit\'e de la proposition 2.3(i) relative \`a l'espace hypersp\'ecial $c^V\tilde{K}_{1}^V$. On obtient
$$J_{{\cal O},\lambda_{1}}^{\tilde{G}_{1}}({\bf f}_{V})=m^{-1}\int_{C_{1}(F)\backslash C_{1}({\mathbb A}_{F})}\sum_{{\cal O}_{1}\in Fib^{\tilde{G}_{1}}({\cal O})}\sum_{M\in{\cal L}(M_{0})}\vert W^M\vert \vert W^G\vert ^{-1}$$
$$\sum_{{\cal O}_{\tilde{M},1}\in \tilde{M}_{1,ss}(F)/conj, {\cal O}_{\tilde{M},1}\subset {\cal O}_{1}}J_{\tilde{M}_{1}}^{\tilde{G}_{1}}(A^{\tilde{M}_{1}}(V,{\cal O}_{\tilde{M},1},c^V\tilde{K}_{1}^{\tilde{M}_{1},V}),\phi_{V}^{c_{V}}\otimes dg_{1})\lambda_{1}(c)\,dc,$$
o\`u $m=mes(\mathfrak{A}_{C_{1}}C_{1}(F)\backslash C_{1}({\mathbb A}_{F}))$. Pour un Levi $M$ fix\'e, sommer en ${\cal O}_{1}\in Fib^{\tilde{G}_{1}}({\cal O})$ puis en ${\cal O}_{\tilde{M},1}\in \tilde{M}_{1,ss}(F)/conj$, ${\cal O}_{\tilde{M},1}\subset {\cal O}_{1}$ revient \`a sommer en ${\cal O}_{\tilde{M}}\in \tilde{M}_{ss}(F)/conj$, ${\cal O}_{\tilde{M}}\subset {\cal O}$, puis en ${\cal O}_{\tilde{M},1}\in Fib^{\tilde{M}_{1}}({\cal O}_{\tilde{M}})$.   On obtient
$$(16) \qquad J_{{\cal O},\lambda_{1}}^{\tilde{G}_{1}}({\bf f}_{V})=m^{-1}\int_{C_{1}(F)\backslash C_{1}({\mathbb A}_{F})}\sum_{M\in{\cal L}(M_{0})}\vert W^M\vert \vert W^G\vert ^{-1}\sum_{{\cal O}_{\tilde{M}}\in \tilde{M}_{ss}(F)/conj, {\cal O}_{\tilde{M}}\subset {\cal O}}$$
$$\sum_{{\cal O}_{\tilde{M},1}\in Fib^{\tilde{M}_{1}}({\cal O}_{\tilde{M}})}J_{\tilde{M}_{1}}^{\tilde{G}_{1}}(A^{\tilde{M}_{1}}(V,{\cal O}_{\tilde{M},1},c^V\tilde{K}_{1}^{\tilde{M}_{1},V}),\phi_{V}^{c_{V}}\otimes dg_{1})\lambda_{1}(c)\,dc.$$
Fixons  $M\in {\cal L}(M_{0})$ et ${\cal O}_{\tilde{M}}\in \tilde{M}_{ss}(F)/conj$. Il n'est pas difficile de voir en d\'evissant les d\'efinitions que la d\'efinition (15) pour la distribution $A^{\tilde{M}}_{\lambda_{1}}(V,{\cal O}_{\tilde{M}})$ s'\'etend aux int\'egrales orbitales pond\'er\'ees. C'est-\`a-dire que l'on a
$$J_{\tilde{M}_{1},\lambda_{1}}^{\tilde{G}_{1}}(A^{\tilde{M}}_{\lambda_{1}}(V,{\cal O}_{\tilde{M}}),f_{V}\otimes dg)= m^{-1} \int_{C_{1}(F)\backslash C_{1}({\mathbb A}_{F})}$$
 $$\sum_{{\cal O}_{\tilde{M},1}\in Fib^{\tilde{M}_{1}}({\cal O}_{\tilde{M}})} J_{\tilde{M}_{1}}^{\tilde{G}_{1}}(A^{\tilde{M}_{1}}(V,{\cal O}_{\tilde{M},1},c^V\tilde{K}_{1}^{\tilde{M}_{1},V}),\phi_{V}^{c_{V}}\otimes dg_{1})\lambda_{1}(c)\,dc.$$
Alors, la relation (16) devient (5).

 Avec ces d\'efinitions, on obtient la formule
$$(17) \qquad J_{g\acute{e}om,\lambda_{1}}^{\tilde{G}_{1}}({\bf f}_{V})=\sum_{\tilde{M}\in {\cal L}(\tilde{M}_{0})}\vert W^{\tilde{M}}\vert \vert W^{\tilde{G}}\vert ^{-1}\sum_{{\cal O}_{\tilde{M}}\in \tilde{M}_{ss}(F)/conj}J_{\tilde{M}_{1},\lambda_{1}}^{\tilde{G}_{1}}(A_{\lambda_{1}}^{\tilde{M}_{1}}(V,{\cal O}_{\tilde{M}}),{\bf f}_{V}).$$

Si $S$ est un ensemble fini de places contenant $V$, on a une formule analogue \`a 2.2(7). Plus pr\'ecis\'ement,   pour $\tilde{M}\in {\cal L}(\tilde{M}_{0})$ et ${\cal O}_{\tilde{M}}\in \tilde{M}_{ss}(F)/conj$ tel que ${\cal O}_{\tilde{M}}\subset {\cal O}$, on peut \'ecrire
$$A_{\lambda_{1}}^{\tilde{M}_{1}}(S,{\cal O}_{\tilde{M}})=\sum_{i=1,....,n({\cal O}_{\tilde{M}})}k_{\lambda_{1},i}^{\tilde{M}_{1}}({\cal O}_{\tilde{M}})_{S}^V\otimes A_{\lambda_{1},i}^{\tilde{M}_{1}}({\cal O}_{\tilde{M}})_{V},$$
o\`u, cette fois, $k_{\lambda_{1},i}^{\tilde{M}_{1}}({\cal O}_{\tilde{M}})_{S}^V\in D_{g\acute{e}om,\lambda_{1}}(\tilde{M}_{1}(F_{S}^V))\otimes Mes(M(F_{S}^V))^*$ et
$A_{\lambda_{1},i}^{\tilde{M}_{1}}({\cal O}_{\tilde{M}})_{V}\in $

\noindent $D_{orb,\lambda_{1}}(\tilde{M}_{1}(F_{V}))\otimes Mes(M(F_{V}))^*$.  On note $A_{\lambda_{1},i}^{\tilde{G}_{1}}({\cal O}_{\tilde{M}})_{V}$ la distribution induite \`a $\tilde{G}_{1}(F_{V})$. On se rappelle la forme lin\'eaire $r_{\tilde{M},\lambda_{1}}^{\tilde{G}}(.,\tilde{K}_{1,S}^V)$ sur $D_{g\acute{e}om,\lambda_{1}}(\tilde{M}_{1}(F_{S}^V))\otimes Mes(M(F_{V}^V))^*$ de 1.15.  On a alors l'\'egalit\'e
$$ A_{\lambda_{1}}^{\tilde{G}_{1}}(V,{\cal O})=\sum_{\tilde{M}\in {\cal L}(\tilde{M}_{0})}\vert W^{\tilde{M}}\vert \vert W^{\tilde{G}}\vert ^{-1}\sum_{{\cal O}_{\tilde{M}}\in \tilde{M}_{ss}(F)/conj, {\cal O}_{\tilde{M}}\subset {\cal O}}$$
$$\sum_{i=1,...,n({\cal O}_{\tilde{M}})}r_{\tilde{M},\lambda_{1}}^{\tilde{G}}(k_{\lambda_{1},i}^{\tilde{M}_{1}}({\cal O}_{\tilde{M}} )_{S}^V,\tilde{K}_{1,S}^V)A_{\lambda_{1},i}^{\tilde{G}_{1}}({\cal O}_{\tilde{M}})_{V}.$$

 Il est utile d'obtenir une expression plus explicite pour la distribution $A_{\tilde{\lambda}_{1}}^{\tilde{G}_{1}}(V,{\cal O})$ dans le cas o\`u $V$ contient $S({\cal O})$, ce que l'on suppose dans ce qui suit.  On suppose aussi que, pour tout $v\not\in V$, la classe de conjugaison par $G(F_{V})$ engendr\'ee par ${\cal O}$ coupe $\tilde{K}_{v}$ (sinon, $A_{\tilde{\lambda}_{1}}^{\tilde{G}_{1}}(V,{\cal O})$ est assez explicite d'apr\`es (6)). Fixons $\dot{\gamma}\in {\cal O}$, fixons  une classe ${\cal O}_{1}\in \tilde{G}_{1}(F)$ d'image ${\cal O}$ et $\dot{\gamma}_{1}\in  {\cal O}_{1}$ d'image $\gamma$. On fixe $\zeta\in C_{1}({\mathbb A}_{F}^V)$ tel que, pour tout $v\not\in V$, la classe de conjugaison par $G_{1}(F_{v})$ engendr\'ee par ${\cal O}_{1}$ coupe $\zeta_{v}\tilde{K}_{1,v}$. On construit comme en 2.3 un \'el\'ement $x^V\in G({\mathbb A}_{F}^V)$ tel que $(x^V)^{-1}\dot{\gamma}x^{V}\in \tilde{K}^{V}$ et on d\'efinit le groupe compact hypersp\'ecial $K_{\dot{\gamma}}^V=x^VK^{V}(x^V)^{-1}\cap G_{\dot{\gamma}}({\mathbb A}_{F}^V)$ de $G_{\dot{\gamma}}({\mathbb A}_{F}^V)$. Il s'en d\'eduit un groupe hypersp\'ecial  $K_{\dot{\gamma}_{1}}^V$ de $G_{1,\dot{\gamma}_{1}}({\mathbb A}_{F}^V)$, qui n'est autre que $x^VK_{1}^{V}(x^V)^{-1}\cap G_{1,\dot{\gamma}_{1}}({\mathbb A}_{F}^V)$. La distribution $A^{G_{\dot{\gamma}}}(V,K_{\dot{\gamma}}^V)$ est bien d\'efinie.  Consid\'erons la suite exacte
$$0\to \mathfrak{c}_{1}(F_{V})\to \mathfrak{g}_{1,\dot{\gamma}_{1}}(F_{V})\to \mathfrak{g}_{\dot{\gamma}}(F_{V})\to 0.$$
On peut la scinder: le scindage est canonique au-dessus de la partie semi-simple de $\mathfrak{g}_{\dot{\gamma}}(F_{V})$, il suffit de la scinder au-dessus du centre de cette alg\`ebre. Notons $\iota_{V}:\mathfrak{g}_{\dot{\gamma}}(F_{V})\to \mathfrak{g}_{1,\dot{\gamma}_{1}}(F_{V})$ un tel scindage, qui est un homomorphisme d'alg\`ebres de Lie. Soit $U_{V}$ un voisinage invariant de l'unit\'e dans $G_{\dot{\gamma}}(F_{V})$, assez petit pour que l'exponentielle y soit d\'efinie. On d\'efinit un plongement $i_{V}:U_{V}\to G_{1,\dot{\gamma}_{1}}(F_{V})$ de la fa\c{c}on suivante: pour $x\in U_{V}$, on \'ecrit $x=exp(X)$ avec $X\in \mathfrak{g}_{\dot{\gamma}}(F_{V})$ et on pose $i_{V}(x)=exp(\iota_{V}(X))$. Ce plongement est \'equivariant pour l'action par conjugaison de $G_{\dot{\gamma}}(F_{V})$. Il est canonique sur l'intersection de $U_{V}$ avec la partie semi-simple de $G_{\dot{\gamma}}(F_{V})$, a fortiori sur les \'el\'ements unipotents. On fixe des mesures $dg_{V}$ sur $G(F_{V})$ et $dh_{V}$ sur $G_{\dot{\gamma}}(F_{V})$. Soit $f_{V}\in C_{c,\lambda_{1}}^{\infty}(\tilde{G}_{1}(F_{V}))$, posons ${\bf f}_{V}=f_{V}\otimes dg_{V}$. Pour $y\in G(F_{V})$, notons $^yf_{V\vert G_{\dot{\gamma}}}$ un \'el\'ement de $ C_{c}^{\infty}(G_{\dot{\gamma}}(F_{V}))$ qui co\^{\i}ncide sur $U_{V}$ avec la fonction $x\mapsto f_{V}(y^{-1}i_{V}(x)\dot{\gamma}_{1}y)$. Posons $^y{\bf f}_{V\vert  G_{\dot{\gamma}}}={^yf}_{V\vert  G_{\dot{\gamma}}}\otimes dh_{V}$. On d\'efinit une distribution $\underline{A}_{\lambda_{1}}^{\tilde{G}_{1}}(V,{\cal O})\in D_{orb,\lambda_{1}}(\tilde{G}_{1}(F_{V}),\omega)\otimes Mes(G(F_{V}))^*$ par les formules:

- $\underline{A}^{\tilde{G}_{1}}_{\lambda_{1}}(V,{\cal O})=0$ si ${\cal O}\not\subset \tilde{G}(F)_{ell}$;

- si ${\cal O}\subset \tilde{G}(F)_{ell}$,
$$I^{\tilde{G}_{1}}_{\lambda_{1}}(\underline{A}_{\tilde{\lambda}_{1}}^{\tilde{G}_{1}}(V,{\cal O}),{\bf f}_{V})=[Z_{G}(\dot{\gamma},F):G_{\dot{\gamma}}(F)]^{-1}\lambda_{1}(\zeta)\int_{G_{\dot{\gamma}}(F_{V})\backslash G(F_{V})}I^{G_{\dot{\gamma}}}(A^{G_{\dot{\gamma}}}_{unip}(V, K_{\dot{\gamma}}^V),{^y{\bf f}}_{V\vert  G_{\dot{\gamma}}})\,dy$$
o\`u la mesure $dy$ est d\'eduite des mesures $dg_{V}$ et $dh_{V}$. 

Alors

(18) on a l'\'egalit\'e $A_{\tilde{\lambda}_{1}}^{\tilde{G}_{1}}(V,{\cal O})=\underline{A}_{\tilde{\lambda}_{1}}^{\tilde{G}_{1}}(V,{\cal O})$.

Preuve. La formule (14) donne 
$$(19) \qquad I^{\tilde{G}_{1}}_{\lambda_{1}}(A_{\tilde{\lambda}_{1}}^{\tilde{G}_{1}}(V,{\cal O}),{\bf f}_{V})=m^{-1}d^{-1}\lambda_{1}(\zeta)I^{\tilde{G}_{1}}(A^{\tilde{G}_{1}}(V,{\cal O}_{1},\zeta\tilde{K}_{1}^V),f_{V}\otimes dg_{1}).$$
 Rappelons que $d$ est le nombre d'\'el\'ements du noyau de l'action de $C_{1}(F)$ sur $Fib^{\tilde{G}_{1}}({\cal O})$. On v\'erifie que
$$d=[Z_{G}(\gamma,F):G_{\gamma}(F)][Z_{G_{1}}(\gamma_{1},F):G_{1,\gamma_{1}}(F)]^{-1}.$$
On a $S({\cal O}_{1})=S({\cal O})$ par d\'efinition. La distribution $A^{\tilde{G}_{1}}(V,,{\cal O}_{1},\zeta\tilde{K}_{1}^V)$ est d\'efinie par la formule (7) de 2.3. En particulier, elle est nulle si ${\cal O}_{1}$ n'est pas elliptique, ce qui \'equivaut \`a ce que ${\cal O}$ ne le soit pas. Supposons ${\cal O}$ elliptique. On fixe sur $G_{1,\dot{\gamma}_{1}}(F_{V})$ la mesure $dh_{1,V}$ d\'eduite des mesures $dh_{V}$ sur $G_{\dot{\gamma}}(F_{V})$ et $dz$ sur $C_{1}({\mathbb A}_{F})$ (cette derni\`ere \'etant identifi\'ee comme toujours \`a une mesure sur $C_{1}(F_{V})$). Notons que $G_{1,\dot{\gamma}_{1}}(F_{V})\backslash G_{1}(F_{V})=G_{\dot{\gamma}}(F_{V})\backslash G(F_{V})$ et que la mesure sur ce quotient d\'eduite  de $dg_{1,V}$ et $dh_{1,V}$ co\"{\i}ncide avec celle d\'eduite de $dg_{V}$ et $dh_{V}$. Notons aussi que le terme $[Z_{G_{1}}(\dot{\gamma}_{1},F):G_{1,\dot{\gamma}_{1}}(F)]$ intervenant dans $d$ va compenser le  facteur intervenant dans 2.3(7). La conjonction de la formule (19) ci-dessus et de 2.3(7) donne
$$I^{\tilde{G}_{1}}_{\lambda_{1}}(A_{\tilde{\lambda}_{1}}^{\tilde{G}_{1}}(V,{\cal O}),{\bf f}_{V})=m^{-1}[Z_{G}(\dot{\gamma},F):G_{\dot{\gamma}}(F)]^{-1} \lambda_{1}(\zeta)$$
$$\int_{G_{\dot{\gamma}}(F_{V})\backslash G(F_{V})}I^{G_{1,\dot{\gamma}_{1}}}(A^{G_{1,\dot{\gamma}_{1}}}_{unip}(V,K_{\dot{\gamma}_{1}}^V),{^yf}_{V\vert  G_{1,\dot{\gamma}_{1}}}\otimes dh_{1,V})\,dy.$$
Pour d\'emontrer (18), il suffit de prouver que, pour tout $y$, on a l'\'egalit\'e
$$m^{-1}I^{G_{1,\gamma_{1}}}(A_{unip}^{G_{1,\dot{\gamma}_{1}}}(V,K_{\dot{\gamma}_{1}}^V),{^yf}_{V\vert  G_{1,\dot{\gamma}_{1}}}\otimes dh_{1,V})=I^{G_{\dot{\gamma}}}(A^{G_{\dot{\gamma}}}_{unip}(V,K_{\dot{\gamma}}^V),{^yf}_{V\vert  G_{\gamma}}\otimes dh_{V}).$$
Remarquons que la fonction ${^yf}_{V\vert  G_{\gamma}}$ co\"{\i}ncide dans $U_{V}$ avec ${^yf}_{V\vert  G_{1,\gamma_{1}}}\circ i_{V}$. L'\'egalit\'e \`a prouver r\'esulte alors de (17) ci-dessous. $\square$

On va simplement reformuler cette \'egalit\'e sous une forme plus g\'en\'erale. On consid\`ere $C_{1}$, $G_{1}$ et $G$ comme pr\'ec\'edemment, en oubliant $\tilde{G}$ et $\tilde{G}_{1}$. On fixe comme plus haut des compacts qui se correspondent $K_{1,v}$ et $K_{v}$ pour $v\not\in V_{ram}(G_{1})$ et des mesures $dg$, $dc$ et $dg_{1}$. Soit $V$ un ensemble fini de places contenant $V_{ram}(G_{1})$. Notons $G_{1,unip}(F_{V})$ et $G_{unip}(F_{V})$ les ensembles d'\'el\'ements unipotents de $G_{1}(F_{V})$ et $G(F_{V})$. La projection $G_{1,unip}(F_{V})\to G_{unip}(F_{V})$ est un isomorphisme. Notons $i_{V}$ son inverse. Soit $f_{1,V}\in C_{c}^{\infty}(G_{1}(F_{V}))$ et  $f_{V}\in C_{c}^{\infty}(G(F_{V}))$. On suppose que   $f_{V}=f_{1,V}\circ i_{V}$ sur $G_{unip}(F_{V})$. Alors on a l'\'egalit\'e
$$(20)\qquad mes(\mathfrak{A}_{C_{1}}C_{1}(F)\backslash C_{1}({\mathbb A}_{F}))^{-1}I^{G_{1}}(A^{G_{1}}_{unip}(V),f_{1,V}\otimes dg_{1,V})=I^{G}(A^{G}_{unip}(V),f_{V}\otimes dg_{V}),$$
o\`u les distributions unipotentes sont  relatives aux compacts fix\'es.

Preuve. Fixons un ensemble ${\cal U}_{V}$ de repr\'esentants des classes de conjugaison par $G(F_{V})$ dans l'ensemble des \'el\'ements unipotents de ce groupe. Son image $i_{V}({\cal U}_{V})$ dans $G_{1}(F_{V})$ est un ensemble analogue pour $G_{1}(F_{V})$. Pour tout $u\in {\cal U}_{V}$, fixons une mesure sur $G_{u}(F_{V})$, dont on d\'eduit une mesure sur $G_{1,i_{V}(u)}(F_{V})$. Le membre de gauche ci-dessus est combinaison lin\'eaire de termes
$$\int_{G_{1,i_{V}(u)}(F_{V})\backslash G_{1}(F_{V})}f_{1,V}(x_{1}^{-1}i_{V}(u)x_{1})\,dx_{1}$$
pour $u\in {\cal U}_{V}$, o\`u la mesure $dx_{1}$ sur le quotient se d\'eduit des mesures sur les deux groupes,  tandis que celui de droite est combinaison lin\'eaire de termes analogues
$$\int_{G_{u}(F_{V})\backslash G(F_{V})}f_{V}(x^{-1}ux)\,dx.$$
De la d\'efinition de $f_{V}$ se d\'eduit que les deux int\'egrales ci-dessus sont \'egales. L'assertion (20) revient \`a dire que les coefficients des deux combinaisons lin\'eaires sont \'egaux. On d\'eduit ais\'ement de cela que, pour un Levi $M\not=G$, si l'on suppose par r\'ecurrence  (20) vrai pour $M$, cette assertion se g\'en\'eralise aux int\'egrales orbitales pond\'er\'ees, c'est-\`a-dire que l'on a l'\'egalit\'e
$$mes(\mathfrak{A}_{C_{1}}C_{1}(F)\backslash C_{1}({\mathbb A}_{F}))^{-1}J_{M_{1}}^{G_{1}}(A^{M_{1}}_{unip}(V),f_{1,V}\otimes dg_{1,V})=J_{M}^G(A_{unip}^M(V),f_{V}\otimes dg_{V}).$$
  En vertu de 2.2(1), l'\'egalit\'e (20) r\'esulte par r\'ecurrence de l'\'egalit\'e que l'on prouvera ci-dessous
$$(21)\qquad mes(\mathfrak{A}_{C_{1}}C_{1}(F)\backslash C_{1}({\mathbb A}_{F}))^{-1}J^{G_{1}}_{unip}(f_{1,V}\otimes dg_{1})=J^G_{unip}(f_{V}\otimes dg).$$
Soit $v\not\in V$. Alors ${\bf 1}_{K_{1,v}}$ et ${\bf 1}_{K_{v}}$ se correspondant comme $f_{V}$ et $f_{1,V}$, c'est-\`a-dire que leurs restrictions aux ensembles d'unipotents $G_{1,unip}(F_{v})$ et $G_{unip}(F_{v})$ se correspondent par la projection bijective du premier ensemble sur le second. En effet, pour $g_{1}\in G_{1}(F_{v})$, $g_{1}$ appartient \`a $K_{1,v}$ si et seulement s'il v\'erifie les deux conditions suivantes:

- sa projection $g$ dans $G(F_{v})$ appartient \`a $K_{v}$;

- pour tout $\chi\in X^*(G_{1})^{\Gamma_{F_{v}}}$ (cf. preuve de (13)), on a $\vert \chi(g_{1})\vert _{F_{v}}=1$.

La deuxi\`eme condition \'etant toujours r\'ealis\'ee pour $g_{1}$ unipotent, l'assertion s'ensuit.
Notons $f_{1}=f_{1,V}\otimes {\bf 1}_{K_{1}^V}$, $f=f_{V}\otimes {\bf 1}_{K^V}$.  Alors les restrictions des fonctions $f_{1}$ et $f$ aux \'el\'ements unipotents de $G_{1}({\mathbb A}_{F})$ et $G({\mathbb A}_{F})$ s'identifient encore. On voit alors que, pour tout $g_{1}\in G_{1}({\mathbb A}_{F})$, on a l'\'egalit\'e
$$k_{unip}^{G_{1},T}(f_{1},g_{1})=k^{G,T}_{unip}(f,g),$$
o\`u $g$ est la projection de $g_{1}$. On construit $J^{G_{1},T}_{unip}(f_{1})$, resp. $J^{G,T}_{unip}(f)$, en int\'egrant cette fonction sur $\mathfrak{A}_{G_{1}}G_{1}(F)\backslash G_{1}({\mathbb A}_{F})$, resp. $\mathfrak{A}_{G}G(F)\backslash G({\mathbb A}_{F})$. Ces int\'egrales sont \'egales, au rapport des mesures pr\`es, et ce rapport est bien s\^ur $mes(\mathfrak{A}_{C_{1}}C_{1}(F)\backslash C_{1}({\mathbb A}_{F}))$. En prenant ensuite la valeur de ces termes au point $T=T_{0}$, on obtient (17). Cela ach\`eve les preuves de (21) et (19). $\qquad$

\bigskip

\subsection{Variante avec caract\`ere central, suite}
On continue avec les m\^emes donn\'ees que dans le paragraphe pr\'ec\'edent. On se donne d'autres donn\'ees
$$1\to C_{2}\to G_{2}\to G\to 1,\,\,\tilde{G}_{2}\to \tilde{G}$$
 $\lambda_{2}$  v\'erifiant les m\^emes hypoth\`eses, munies de donn\'ees compl\'ementaires comme en 1.15.  On suppose donn\'e un caract\`ere automorphe $\lambda_{12}$ de $G_{12}({\mathbb A}_{F})$ tel que la restriction de $\lambda_{12}$ \`a $C_{1}({\mathbb A}_{F})\times C_{2}({\mathbb A}_{F})$ soit $\lambda_{1}\times \lambda_{2}^{-1}$.    Soit $V$ un ensemble fini de places de $F$ contenant $V_{12,ram}$. On a construit en 1.15 une fonction canonique $\tilde{\lambda}_{12,V}$ sur $\tilde{G}_{12}(F_{V})$, dont on d\'eduit  un isomorphisme
$$D_{g\acute{e}om,\lambda_{1}}(\tilde{G}_{1}(F_{V}))\otimes Mes(G(F_{V}))^*\simeq D_{g\acute{e}om,\lambda_{2}}(\tilde{G}_{2}(F_{V}))\otimes Mes(G(F_{V}))^*$$
qui se restreint en un isomorphisme
$$D_{orb,\lambda_{1}}(\tilde{G}_{1}(F_{V}))\otimes Mes(G(F_{V}))^*\simeq D_{orb,\lambda_{2}}(\tilde{G}_{2}(F_{V}))\otimes Mes(G(F_{V}))^*.$$

\ass{Lemme}{Pour tout ${\cal O}\in \tilde{G}_{ss}(F)/conj$, les distributions $A^{\tilde{G}_{1}}_{\lambda_{1}}(V,{\cal O})$ et $A^{\tilde{G}_{2}}_{\lambda_{2}}(V,{\cal O})$ se correspondent par cet isomorphisme.}

Preuve. Consid\'erons d'abord la situation du paragraphe pr\'ec\'edent dans le cas particulier o\`u $\lambda_{1}=1$. Alors $C_{c,\lambda_{1}}^{\infty}(\tilde{G}_{1}({\mathbb A}_{F}))$ s'identifie \`a $C_{c}^{\infty}(\tilde{G}({\mathbb A}_{F}))$. Pour un \'el\'ement $f$ de cet espace et une mesure $dg$ sur $G({\mathbb A}_{F})$, on a d\'efini deux termes $J_{g\acute{e}om,\lambda_{1}}^{\tilde{G}_{1}}(f\otimes dg)$ et $J_{g\acute{e}om}^{\tilde{G}}(f\otimes dg)$. Le premier en 2.3 en  consid\'erant $f$ comme un \'el\'ement du premier espace, le second en 2.1 en consid\'erant $f$ comme un \'el\'ement du second espace. En d\'evissant les d\'efinitions, on v\'erifie que ces deux termes sont \'egaux. Plus finement, pour ${\cal O}\in \tilde{G}_{ss}(F)/conj$, les deux termes $J_{{\cal O},\lambda_{1}}^{\tilde{G}_{1}}(f\otimes dg)$ et $J_{{\cal O}}^{\tilde{G}}(f\otimes dg)$ d\'efinis en 2.3 et 2.1 sont \'egaux.

D'autre part, soit $f\in C_{c}^{\infty}(\tilde{G}({\mathbb A}_{F}))$ et $\lambda$ un caract\`ere automorphe de $G({\mathbb A}_{F})$. On en d\'eduit comme en 1.15 une fonction $\tilde{\lambda}$ sur $\tilde{G}({\mathbb A}_{F})$ qui vaut $1$ sur $\tilde{G}(F)$. On v\'erifie imm\'ediatement que l'on a l'\'egalit\'e
$$J_{g\acute{e}om}^{\tilde{G}}((f\tilde{\lambda})\otimes dg)=J_{g\acute{e}om}^{\tilde{G}}(f\otimes dg).$$
Plus finement, pour tout ${\cal O}\in \tilde{G}_{ss}(F)/conj$, on a l'\'egalit\'e
$$J_{{\cal O}}^{\tilde{G}}((f\tilde{\lambda})\otimes dg)=J_{{\cal O}}^{\tilde{G}}(f\otimes dg).$$

Revenons maintenant \`a l'\'enonc\'e et fixons des mesures sur tous nos groupes. En consid\'erant la relation 2.3(4), on voit que, par r\'ecurrence, il nous suffit de prouver que, pour $f_{1,V}\in C_{c,\lambda_{1}}^{\infty}(\tilde{G}_{1}(F_{V}))$ et $f_{2,V}\in C_{c,\lambda_{2}}^{\infty}(\tilde{G}_{2}(F_{V}))$ se correspondant par l'isomorphisme ci-dessus, on a l'\'egalit\'e
$$J_{{\cal O},\lambda_{1}}^{\tilde{G}_{1}}(f_{1,V}\otimes dg)=J_{{\cal O},\lambda_{2}}^{\tilde{G}_{2}}(f_{2,V}\otimes dg).$$
Posons $f_{i}=f_{i,V}\otimes{\bf 1}_{\tilde{K}_{i}^V,\lambda_{i}}$ pour $i=1,2$. En vertu des d\'efinitions, ces fonctions se correspondent par la relation
$$f_{2}(\gamma_{2})=\tilde{\lambda}_{12}(\gamma_{1},\gamma_{2})f_{1}(\gamma_{1})$$
pour tous $(\gamma_{1},\gamma_{2})\in \tilde{G}_{12}({\mathbb A}_{F})$. Et la relation \`a prouver est
$$J_{{\cal O},\lambda_{1}}^{\tilde{G}_{1}}(f_{1}\otimes dg)=J_{{\cal O},\lambda_{2}}^{\tilde{G}_{2}}(f_{2}\otimes dg).$$
On s'est ainsi d\'ebarrass\'e de l'ensemble $V$. Introduisons une fonction $\phi_{1}\in C_{c}^{\infty}(\tilde{G}_{1}({\mathbb A}_{F}))$ telle que
$$f_{1}=\int_{C_{1}({\mathbb A}_{F})}\phi_{1}^c\lambda_{1}(c)\,dc.$$
On peut consid\'erer $\phi_{1}$ comme une fonction sur $\tilde{G}_{12}({\mathbb A}_{F})$, invariante par le sous-groupe $C_{2}({\mathbb A}_{F})$. En la consid\'erant ainsi, on peut introduire une fonction $\phi_{12}\in C_{c}^{\infty}(\tilde{G}_{12}({\mathbb A}_{F}))$ telle que
$$\phi_{1}=\int_{C_{2}({\mathbb A}_{F})}\phi_{12}^c\,dc.$$
Notons $\phi_{21}$ le produit de $\phi_{12}$ et de la fonction $\tilde{\lambda}_{12}$ sur $\tilde{G}_{12}({\mathbb A}_{F})$. D\'efinissons une fonction $\phi_{2}$ invariante par $C_{1}({\mathbb A}_{F})$ par
$$\phi_{2}=\int_{C_{1}({\mathbb A}_{F})}\phi_{21}^c\,dc.$$
On peut la consid\'erer comme un \'el\'ement de $C_{c}^{\infty}(\tilde{G}_{2}({\mathbb A}_{F}))$. On v\'erifie qu'on a alors l'\'egalit\'e
$$f_{2}=\int_{C_{2}({\mathbb A}_{F})}\phi_{2}^c\lambda_{2}(c)\,dc.$$
Utilisons la formule 2.3(3):
$$J_{{\cal O},\lambda_{1}}^{\tilde{G}_{1}}(f_{1}\otimes dg)=m_{1}^{-1}\int_{C_{1}(F)\backslash C_{1}({\mathbb A}_{F})}\sum_{{\cal O}_{1}}J^{\tilde{G}_{1}}_{{\cal O}_{1}}(\phi_{1}^{c_{1}}\otimes dg_{1})\lambda_{1}(c_{1})\,dc_{1},$$
o\`u $m_{1}=mes(\mathfrak{A}_{C_{1}}C_{1}(F)\backslash C_{1}({\mathbb A}_{F}))$ et o\`u
la somme porte sur les ${\cal O}_{1}$ au-dessus de ${\cal O}$. Pour chaque ${\cal O}_{1}$, on utilise la premi\`ere remarque ci-dessus. Elle entra\^{\i}ne que, pour tout $c_{1}$, on a l'\'egalit\'e
$$J^{\tilde{G}_{1}}_{{\cal O}_{1}}(\phi_{1}^{c_{1}}\otimes dg_{1})=m_{2}^{-1}\int_{C_{2}(F)\backslash C_{2}({\mathbb A}_{F})}\sum_{{\cal O}_{12}}J^{\tilde{G}_{1,2}}_{{\cal O}_{12}}(\phi_{12}^{c_{1}c_{2}}\otimes dg_{12})\,dc_{2},$$
o\`u la signification de $m_{2}$ est claire et o\`u la somme porte sur les ${\cal O}_{12}\in \tilde{G}_{12,ss}(F)/conj$ au-dessus de ${\cal O}_{1}$. On obtient
$$J_{{\cal O},\lambda_{1}}^{\tilde{G}_{1}}(f_{1}\otimes dg)=m_{1}^{-1}m_{2}^{-1}\int_{C_{1}(F)\backslash C_{1}({\mathbb A}_{F})\times C_{2}(F)\backslash C_{2}({\mathbb A}_{F})}\sum_{{\cal O}_{12}}J^{\tilde{G}_{1,2}}_{{\cal O}_{12}}(\phi_{12}^{c_{1}c_{2}}\otimes dg_{12})\lambda_{1}(c_{1})\,dc_{1}\,dc_{2},$$
o\`u cette fois, la somme porte sur les ${\cal O}_{12}$ au-dessus de ${\cal O}$.
On a une formule analogue en permutant les indices $1$ et $2$.
Il nous suffit de fixer ${\cal O}_{12}$ au-dessus de ${\cal O}$ et de d\'emontrer que, pour tous $c_{1},c_{2}$, on a l'\'egalit\'e
$$J^{\tilde{G}_{1,2}}_{{\cal O}_{12}}(\phi_{12}^{c_{1}c_{2}}\otimes dg_{12})\lambda_{1}(c_{1})=J^{\tilde{G}_{1,2}}_{{\cal O}_{12}}(\phi_{21}^{c_{1}c_{2}}\otimes dg_{12})\lambda_{2}(c_{2}).$$
Mais on voit que la fonction $\lambda_{2}(c_{2})\phi_{21}^{c_{1}c_{2}}$ est le produit de $\lambda_{1}(c_{1})\phi_{12}^{c_{1}c_{2}}$ et du caract\`ere automorphe $\tilde{\lambda}_{12}$ de $\tilde{G}_{12}({\mathbb A}_{F})$. L'\'egalit\'e pr\'ec\'edente r\'esulte de la deuxi\`eme remarque du d\'ebut de la preuve. $\square$

 \bigskip
 
\subsection{La partie g\'eom\'etrique de la formule des traces $\omega$-\'equivariante}
 On revient \`a la situation de 2.1. On fixe un ensemble fini $V$ de places contenant $V_{ram}$. On note $\tilde{G}_{ss}(F_{V})/conj$ l'ensemble des  classes de conjugaison semi-simples par $G(F_{V})$ dans $\tilde{G}(F_{V})$, autrement dit l'ensemble des familles ${\cal O}_{V}=({\cal O}_{v})_{v\in V}$, o\`u,  pour tout $v\in V$, ${\cal O}_{v}$ est  une classe de conjugaison semi-simple par $G(F_{v})$ dans $\tilde{G}(F_{v})$. Soit ${\cal O}_{V}$ une telle famille.  Pour ${\cal O}\in \tilde{G}_{ss}(F)/conj$, la relation ${\cal O}\subset {\cal O}_{V}$ signifie que la classe de conjugaison par $G(F_{V})$ engendr\'ee par ${\cal O}$ est \'egale \`a ${\cal O}_{V}$. On d\'efinit  un \'el\'ement $A^{\tilde{G}}({\cal O}_{V},\omega)\in D_{orb}(\tilde{G}(F_{V}),\omega)\otimes Mes(G(F_{V}))^*$ par la formule
 $$A^{\tilde{G}}({\cal O}_{V},\omega)=\sum_{{\cal O}\in \tilde{G}_{ss}(F)/conj, {\cal O}\subset {\cal O}_{V}}A^{\tilde{G}}(V,{\cal O},\omega).$$
 Notons que les ${\cal O}$ qui contribuent v\'erifient non seulement ${\cal O}\subset {\cal O}_{V}$, mais aussi que, pour tout $v\not\in V$, la classe de conjugaison par $G(F_{v})$ engendr\'ee par ${\cal O}$ coupe $\tilde{K}_{v}$, cf. 2.2(3). Il n'y en a qu'un nombre fini d'apr\`es le lemme 2.1. Pour ${\bf f}\in C_{c}^{\infty}(\tilde{G}(F_{V}))\otimes Mes(G(F_{V}))$, la formule 2.2(8) se r\'ecrit
 $$(1) \qquad J_{g\acute{e}om}^{\tilde{G}}({\bf f},\omega)=\sum_{\tilde{M}\in {\cal L}(\tilde{M}_{0})}\vert W^{\tilde{M}}\vert \vert W^{\tilde{G}}\vert ^{-1}\sum_{{\cal O}_{V}\in \tilde{M}_{ss}(F_{V})/conj} J_{\tilde{M}}^{\tilde{G}}(A^{\tilde{M}}({\cal O}_{V},\omega),{\bf f}_{V}).$$
 Toutes les distributions intervenant sont \`a support dans l'ensemble des $\gamma\in \tilde{G}(F_{V})$ tels que $\tilde{H}_{\tilde{G}}(\gamma)=0$. Elles se prolongent donc \`a tout $C^{\infty}_{ac,glob}(\tilde{G}(F_{V}))\otimes Mes(G(F_{V}))$ et l'\'egalit\'e ci-dessus reste vraie  sur cet espace.

 En utilisant le lemme 1.7, on d\'efinit une distribution ${\bf f}\mapsto I^{\tilde{G}}_{g\acute{e}om}({\bf f},\omega)$ sur $C^{\infty}_{ac,glob}(\tilde{G}(F_{V}))\otimes Mes(G(F_{V}))$ par la formule de r\'ecurrence
$$(2) \qquad I^{\tilde{G}}_{g\acute{e}om,}({\bf f},\omega)=J^{\tilde{G}}_{g\acute{e}om}({\bf f},\omega)-\sum_{\tilde{L}\in {\cal L}(\tilde{M}_{0}),\tilde{L}\not=\tilde{G}}\vert W^{\tilde{L}}\vert \vert W^{\tilde{G}}\vert ^{-1}I^{\tilde{L}}_{g\acute{e}om}(\phi_{\tilde{L}}({\bf f}),\omega).$$
En [A3] paragraphe 2, Arthur montre que c'est une distribution $\omega$-\'equivariante (cette propri\'et\'e est comme toujours suppos\'ee par r\'ecurrence pour que la d\'efinition ci-dessus ait un sens).  Elle se factorise donc en une forme lin\'eaire sur $I_{ac,glob}(\tilde{G}(F_{V}))\otimes Mes(G(F_{V}))$. 
Soulignons que cette distribution d\'epend implicitement des $\tilde{K}_{v}$ pour $v\not\in V$.

 \ass{Proposition}{ Pour ${\bf f}\in C_{c}^{\infty}(\tilde{G}(F_{V}))\otimes Mes(G(F_{V}))$, on a l'\'egalit\'e
$$I_{g\acute{e}om}^{\tilde{G}}({\bf f},\omega)=\sum_{\tilde{M}\in {\cal L}(\tilde{M}_{0})}\vert W^{\tilde{M}}\vert \vert W^{\tilde{G}}\vert ^{-1}\sum_{{\cal O}_{V}\in \tilde{M}_{ss}(F_{V})/conj}I_{\tilde{M}}^{\tilde{G}}(A^{\tilde{M}}({\cal O}_{V},\omega),{\bf f}).$$
  La somme est en fait finie.}

Preuve. On d\'emontre plus g\'en\'eralement la formule de l'\'enonc\'e pour ${\bf f}\in C^{\infty}_{ac,glob}(\tilde{G}(F_{V}))\otimes Mes(G(F_{V}))$. 
En raisonnant par r\'ecurrence, on exprime les termes $I^{\tilde{L}}_{g\acute{e}om}(\phi_{\tilde{L}}({\bf f}),\omega)$ qui interviennent dans la relation (2) par la formule  de l'\'enonc\'e.  On exprime le terme $J_{g\acute{e}om}^{\tilde{G}}({\bf f},\omega)$ par la formule (1). En regroupant les termes, on  obtient
$$I^{\tilde{G}}_{g\acute{e}om}({\bf f},\omega)=\sum_{\tilde{M}\in {\cal L}(\tilde{M}_{0})}\vert W^{\tilde{M}}\vert \vert W^{\tilde{G}}\vert ^{-1}\sum_{{\cal O}_{V}\in \tilde{M}_{ss}(F_{V})/conj}$$
$$\left(J_{\tilde{M}}^{\tilde{G}}(A^{\tilde{M}}({\cal O}_{V},\omega),{\bf f})-\sum_{\tilde{L}\in {\cal L}(\tilde{M}), \tilde{L}\not=\tilde{G}}I_{\tilde{M}}^{\tilde{L}}(A^{\tilde{M}}({\cal O}_{V},\omega),\phi_{\tilde{L}}({\bf f}))\right).$$
La derni\`ere somme int\'erieure est par d\'efinition $I_{\tilde{M}}^{\tilde{G}}(A^{\tilde{M}}({\cal O}_{V},\omega),{\bf f})$ et on obtient la formule de l'\'enonc\'e. $\square$

 \bigskip
 
\subsection{La partie g\'eom\'etrique de la formule des traces invariante, variante avec caract\`ere central}
Dans la situation de 2.3, pour un ensemble fini $V$ de places contenant $V_{1,ram}$, on peut rendre invariante la distribution ${\bf f}_{V}\mapsto J_{g\acute{e}om,\lambda_{1}}^{\tilde{G}_{1}}({\bf f}_{V})$ sur $C_{c,\lambda_{1}}^{\infty}(\tilde{G}_{1}(F_{V}))\otimes Mes(G(F_{V}))^*$ par le m\^eme proc\'ed\'e qu'au paragraphe pr\'ec\'edent. En adaptant les d\'efinitions de fa\c{c}on plus ou moins \'evidente, on obtient la formule 
$$(1) \qquad I^{\tilde{G}_{1}}_{g\acute{e}om,\lambda_{1}}({\bf f}_{V})=\sum_{\tilde{M}\in {\cal L}(\tilde{M}_{0})}\vert W^{\tilde{M}}\vert \vert W^{\tilde{G}}\vert ^{-1}\sum_{{\cal O}\in \tilde{M}_{ss}(F_{V})/conj}I_{\tilde{M}_{1},\lambda_{1}}^{\tilde{G}_{1}}(A_{\lambda_{1}}^{\tilde{M}_{1}}(V,{\cal O}),{\bf f}_{V}).$$

Remarquons que, pour ${\cal O}\in \tilde{G}_{ss}(F_{V})/conj$, la   distribution $A^{\tilde{G}_{1}}_{\lambda_{1}}(V,{\cal O})$ est donn\'ee  par une formule similaire \`a 2.3(11), \`a savoir avec les notations de ce paragraphe:
 $$(2) \qquad I^{\tilde{G}_{1}}_{\lambda_{1}}(A^{\tilde{G}}_{\lambda_{1}}(V,{\cal O}),{\bf f}_{V})=mes(\mathfrak{A}_{C_{1}}C_{1}(F)\backslash C_{1}({\mathbb A}_{F}))^{-1} \int_{C_{1}(F)\backslash C_{1}({\mathbb A}_{F})}$$
$$\sum_{{\cal O}_{1}\in Fib^{\tilde{G}_{1}}({\cal O})} I^{\tilde{G}_{1}}(A^{\tilde{G}_{1}}(V,{\cal O}_{1},c^V\tilde{K}_{1}^V),\phi_{V}^{c_{V}}\otimes dg_{1})\lambda_{1}(c)\,dc.$$
La somme int\'erieure porte sur les \'el\'ements de $\tilde{G}_{1,ss}(F_{V})/conj$ qui se projettent sur ${\cal O}$. Cet ensemble n'est pas d\'enombrable, mais seuls un nombre fini de termes sont non nuls. 

 Dans la situation de 2.6, on a un r\'esultat similaire \`a celui du lemme de ce paragraphe:
 
 (3) les distributions $I^{\tilde{G}_{1}}_{g\acute{e}om,\lambda_{1}}$ et  $I^{\tilde{G}_{2}}_{g\acute{e}om,\lambda_{2}}$ d\'efinies sur $C_{c,\lambda_{1}}^{\infty}(\tilde{G}_{1}(F_{V}))\otimes Mes(G(F_{V}))^*$ et $C_{c,\lambda_{2}}^{\infty}(\tilde{G}_{2}(F_{V}))\otimes Mes(G(F_{V}))^*$ se correspondent par l'isomorphisme d\'efini en 2.6 entre ces espaces.
 
 \bigskip
 
 \subsection{Variante pour les $K$-espaces}
 
Consid\'erons un $K$-espace comme en 1.16. Soit $V$ un ensemble fini de places contenant $V_{ram}$. Pour ${\bf f}=\oplus_{p\in \Pi}{\bf f}_{p}\in I(K\tilde{G}(F_{V}),\omega)$, on pose 
$$I^{K\tilde{G}}_{geom}({\bf f},\omega)=\sum_{p\in \Pi}I^{\tilde{G}_{p}}_{geom}({\bf f}_{p},\omega).$$
Par simple sommation sur $p$ du d\'eveloppement de  $I^{\tilde{G}_{p}}_{geom}({\bf f}_{p},\omega)$ fourni par la proposition 2.7, on obtient
$$I_{g\acute{e}om}^{K\tilde{G}}({\bf f},\omega)=\sum_{K\tilde{M}\in {\cal L}(K\tilde{M}_{0})}\vert W^{K\tilde{M}}\vert \vert W^{K\tilde{G}}\vert ^{-1}\sum_{{\cal O}_{V}\in K\tilde{M}_{ss}(F_{V})/conj}I_{K\tilde{M}}^{K\tilde{G}}(A^{K\tilde{M}}({\cal O}_{V},\omega),{\bf f}).$$
Donnons quelques explications sur les notations. Un $K$-espace de Levi $K\tilde{M}$ intervenant ci-dessus est une collection $(\tilde{M}_{p})_{p\in \Pi^M}$ o\`u $\Pi^M$ est un sous-ensemble de $\Pi$ et, pour tout $p\in \Pi^M$, $\tilde{M}_{p}$ est un espace de Levi de $\tilde{G}_{p}$ (cf. 1.17). Un \'el\'ement ${\cal O}_{V}\in K\tilde{M}_{ss}(F_{V})/conj$ est un \'el\'ement de $  \tilde{M}_{p,ss}(F_{V})/conj$ pour un certain $p\in \Pi^M$. Par d\'efinition, $A^{K\tilde{M}}({\cal O}_{V},\omega)=A^{\tilde{M}_{p}}({\cal O}_{V},\omega)$ et 
 $$I_{K\tilde{M}}^{K\tilde{G}}(A^{K\tilde{M}}({\cal O}_{V},\omega),{\bf f})=I_{\tilde{M}_{p}}^{\tilde{G}_{p}}(A^{\tilde{M}_{p}}({\cal O}_{V},\omega),{\bf f}_{p}).$$

\bigskip

\section{Endoscopie}
 
\bigskip

\subsection{Donn\'ees endoscopiques}
Une donn\'ee endoscopique de $(G,\tilde{G},{\bf a})$ est un triplet ${\bf G}'=(G',{\cal G}',\tilde{s})$ v\'erifiant les m\^emes conditions que dans le cas local (cf. [I] 1.5), \`a l'unique diff\'erence suivante pr\`es. On suppose qu'il existe un  cocycle $a:W_{F}\to Z(\hat{G})$ tel que:

- pour tout $(g,w)\in {\cal G}'$, on ait l'\'egalit\'e $ad_{\tilde{s}}(g,w)=(a(w)g,w)$;

- l'image de $a$ dans $H^1(W_{F},,Z(\hat{G}))/ker^1(W_{F},Z(\hat{G}))$ est \'egal \`a ${\bf a}$.

La notion d'\'equivalence de donn\'ees endoscopiques est la m\^eme que dans le cas local. Pour une donn\'ee endoscopique ${\bf G}'$ comme ci-dessus, on d\'efinit l'espace endoscopique $\tilde{G}'=G'\times_{{\cal Z}(G)}{\cal Z}(\tilde{G})$ comme en [I] 1.7 (l'ensemble ${\cal Z}(\tilde{G})$ a une structure sur $F$). On note $V_{ram}({\bf G}')$ la r\'eunion de $V_{ram}$ et de l'ensemble des places $v$ o\`u ${\bf G}'_{v}$ est ramifi\'e (on rappelle que, pour $v\not\in V_{ram}$, on dit que ${\bf G}'_{v}$ est non ramifi\'e si le groupe d'inertie $I_{v}\subset W_{F_{v}}$ est contenu dans ${\cal G}'_{v}$).  
On dit que ${\bf G}'$ est relevante si elle v\'erifie les deux conditions

$\bullet$ $\tilde{G}'(F)\not=\emptyset$;

$\bullet$ pour tout $v\in Val(F)$, la donn\'ee ${\bf G}'_{v}$ est relevante.

{\bf Attention.} Il ne s'ensuit pas de ces conditions qu'il existe $\delta\in \tilde{G}'(F)$ qui soit $\tilde{G}$-r\'egulier et tel qu'en toute place $v$, $\delta$ corresponde \`a un \'el\'ement $\gamma_{v}\in \tilde{G}(F_{v})$. La situation se simplifie toutefois dans le cas o\`u $(G,\tilde{G},{\bf a})$ est  quasi-d\'eploy\'e et \`a torsion int\'erieure, d'apr\`es le lemme suivant.

\ass{Lemme}{Supposons $(G,\tilde{G},{\bf a})$ quasi-d\'eploy\'e et \`a torsion int\'erieure. Supposons $\tilde{G}'(F)\not=\emptyset$. Alors l'ensemble des \'el\'ements $\tilde{G}$-r\'eguliers de $\tilde{G}'(F)$ n'est pas vide et, pour tout  \'el\'ement $\delta$ de cet ensemble, il existe $\gamma\in \tilde{G}_{reg}(F)$ qui correspond \`a $\delta$.}

 La preuve est la m\^eme que dans le cas local, cf. [I] lemme 1.9.
 
 {\bf Remarque.} Il se peut que $\tilde{G}'(F)$ soit vide. Par exemple, soient $d\in F^{\times}$, $G=SL(2)$, $\tilde{G}=\{\gamma\in GL(2); det(\gamma)=d\}$ et ${\bf a}=1$. Pour toute extension quadratique  $E$ de $F$, il y a une donn\'ee endoscopique ${\bf G}'$ telle que $G'(F)$ est le groupe des \'el\'ements de $E$ de norme $1$. Alors $\tilde{G}'(F)$ est l'ensemble des \'el\'ements de $E$ de norme $d$. On peut trouver choisir $d$ et $E$ de sorte que cet ensemble soit vide.

\bigskip

Comme dans le cas local, on dit que ${\bf G}'$ est elliptique si $Z(\hat{G}')^{\Gamma_{F},0}=Z(\hat{G})^{\Gamma_{F},\hat{\theta},0}$.

On fixe une paire parabolique  $(P'_{0},M'_{0})$ de $G'$, d\'efinie sur $F$ et minimale. Pour tout $v\not\in V_{ram}({\bf G}')$, supposons fix\'e un sous-groupe compact hypersp\'ecial $K'_{v}$ de $G'(F_{v})$ v\'erifiant  des conditions analogues \`a celles de 1.1.  D'apr\`es [I] 6.2, l'espace $\tilde{K}_{v}$ d\'etermine un espace hypersp\'ecial $\tilde{K}'_{v}$ dans $\tilde{G}'(F_{v})$ associ\'e \`a $K'_{v}$. Ces ensembles v\'erifient encore la condition de compatibilit\'e globale pr\'ec\'edente.

\bigskip

\subsection{Plongements de tores et ramification}

Soient $\hat{T}$ un tore complexe et $E$ une extension galoisienne finie de $F$. On suppose $\hat{T}$  muni d'une action alg\'ebrique de $Gal(E/F)$. Introduisons le groupe de Weil relatif $W_{E/F}$. C'est un quotient de $W_{F}$ et il s'ins\`ere dans une suite exacte
$$1\to E^{\times}\backslash {\mathbb A}_{E}^{\times}\to W_{E/F}\to Gal(E/F)\to 1,$$
cf. [T] paragraphe 1. Notons $I_{F}^{E}$ le sous-groupe distingu\'e de $W_{F}$ engendr\'e par les images des groupes d'inertie $I_{v}\subset W_{F}$ pour toutes les places finies $v$ de $F$ non ramifi\'ees dans $E$ et notons $I_{E/F}$ son image dans $W_{E/F}$. Le groupe $I_{E/F}$ est aussi l'image dans $W_{E/F}$ du sous-groupe $\prod_{w}\mathfrak{o}_{w}^{\times}$, o\`u le produit est pris sur les places finies  $w$ de $E$ au-dessus d'une telle place $v$ de $F$ et o\`u $\mathfrak{o}_{w}^{\times}$ est le groupe des unit\'es de $E_{w}^{\times}$. Soit $a$ un $2$-cocycle de $Gal(E/F)$ dans $\hat{T}$, que l'on rel\`eve en un $2$-cocycle de $W_{E/F}$ dans $\hat{T}$. Alors

(1) il existe une cocha\^{\i}ne continue $b:W_{E/F}\to \hat{T}$, biinvariante par $I_{E/F}$, telle que $a$ soit le bord de $b$. 

Preuve. Langlands prouve qu'il existe une cocha\^{\i}ne continue $b':W_{E/F}\to \hat{T}$ telle que $a$ soit le bord de $b'$ ([Lan] lemme 4). Munissons le groupe compact $I_{E/F}$ de la mesure de Haar de masse totale $1$. Pour $w\in W_{E/F}$, posons 
$$b(w)=\int_{I_{E/F}}b'(iw)\,di.$$ 
C'est une cocha\^{\i}ne continue et invariante \`a gauche par $I_{E/F}$, donc aussi \`a droite puisque $I_{E/F}$ est distingu\'e. Pour $w_{1},w_{2}\in W_{E/F}$, on a l'\'egalit\'e
$$a(w_{1},w_{2})b'(w_{1}w_{2})=b'(w_{1})w_{1}(b'(w_{2})) .$$
Soient $i_{1},i_{2}\in I_{E/F}$. Rempla\c{c}ons $w_{1}$ par $i_{1}w_{1}$ et $w_{2}$ par $i_{2}w_{2}$. Puisque les images de $i_{1}$ et $i_{2}$ dans $Gal(E/F)$ sont triviales, on obtient
$$a(w_{1},w_{2})b'(i_{3}w_{1}w_{2})=b'(i_{1}w_{1})w_{1}(b'(i_{2}w_{2})),$$
o\`u $i_{3}=i_{1}w_{1}i_{2}w_{1}^{-1}$. En int\'egrant cette \'egalit\'e en $i_{1}$ et $i_{2}$, on obtient
$$a(w_{1},w_{2})b(w_{1}w_{2})=b(w_{1})w_{1}(b(w_{2})) ,$$
 c'est-\`a-dire que $a$ est le bord de $b$. $\square$
 
Soient $E$ une extension galoisienne finie de $F$ et  $\hat{T}_{i}$, pour $i=1,2,3$, des tores complexes munis d'actions de $\Gamma_{E/F}$. On suppose donn\'ee une suite exacte
 $$1\to \hat{T}_{1}\to \hat{T}_{2}\to \hat{T}_{3}\to 1$$
 \'equivariante pour ces actions.  Soit $U$ un sous-groupe de $I_{F}^{E}$. On a
 
 (2) pour tout cocycle continu $b_{3}:W_{F}\to \hat{T}_{3}$ biinvariant par $U$, il existe un cocycle continu $b_{2}:W_{F}\to \hat{T}_{2}$, biinvariant par $U$, dont l'image par composition avec l'homomorphisme $\hat{T}_{2}\to \hat{T}_{3}$ soit \'egale \`a $b_{3}$.
 
 Preuve. D'apr\`es la preuve de [Lan] lemme 4, il existe un cocycle $b'_{2}:W_{F}\to \hat{T}_{2}$ dont l'image par composition avec l'homomorphisme $\hat{T}_{2}\to \hat{T}_{3}$ soit \'egale \`a $b_{3}$. On peut fixer une extension galoisienne finie $E'$ de $E$ telle que $b'_{2}$ et $b_{3}$ se factorisent par $W_{E'/F}$. Notons $U'$ la cl\^oture de l'image de $U$ dans $W_{E'/F}$ et munissons ce groupe compact de la mesure de Haar de masse totale $1$. Le cocycle $b_{3}$ est encore biinvariant par $U'$. Pour $w\in W_{E'/F}$, posons 
 $$b_{2}(w)=\int_{U'}b'_{2}(iw)\,di.$$
 Le m\^eme calcul que ci-dessus montre que $b_{2}$ est encore un cocycle. Son rel\`evement en un cocycle sur $W_{F}$ est biinvariant par $U$. Puisque $b_{3}$ est   biinvariant par $U'$, $b_{2}$ v\'erifie comme $b'_{2}$ la derni\`ere condition de l'assertion. $\square$
 
 {\bf Remarque.} Les deux propri\'et\'es ci-dessus ont \'et\'e prouv\'ees de fa\c{c}on diff\'erente par Arthur ([A1], preuve du lemme 7.1).

Soient ${\bf G}'=(G',{\cal G}',\tilde{s})$ une donn\'ee endoscopique de $(G,\tilde{G},{\bf a})$ et $S'$ un sous-tore maximal de $G'$ d\'efini sur $F$. Fixons comme en [I] 1.5 une paire de Borel \'epingl\'ee $\hat{{\cal E}}=(\hat{B},\hat{T},(\hat{E}_{\alpha})_{\alpha\in \Delta})$ de $\hat{G}$ conserv\'ee par $ad_{\tilde{s}}$. On note $\hat{\theta}$ l'automorphisme d\'etermin\'e par cette paire et on adapte l'action galoisienne de sorte que $\hat{{\cal E}}$ soit conserv\'ee par la nouvelle action. La paire de Borel $(\hat{B}',\hat{T}')=(\hat{B}\cap \hat{G}', \hat{T}^{\hat{\theta},0})$ se compl\`ete en une paire de Borel \'epingl\'ee de $\hat{H}$.   Le tore dual de $S'$ s'identifie \`a $\hat{T}^{\hat{\theta},0}$ muni de l'action galoisienne $\sigma\mapsto \omega_{S',G}(\sigma)\circ \sigma$, o\`u $\omega_{S',G}$ est un cocycle \`a valeurs dans $W^{\theta}$. On note $\hat{S}'$ ce tore muni de cette action galoisienne.

\ass{Lemme}{Le plongement $\hat{S}'\subset \hat{G'}\subset {\cal G}'$ se prolonge en un plongement $^LS'\to {\cal G}'\subset {^LG}$ de sorte que, pour toute place $v\in Val(F)$ telle que $v\not\in V_{ram}({\bf G}')$ et $S'$ soit non ramifi\'e en $v$, l'image d'un \'el\'ement $w$ du groupe d'inertie $I_{v}$ soit $(1,w)\in {^LG}$.} 

Preuve. Notons $V$ la r\'eunion de $V_{ram}({\bf G}') $ et de l'ensemble des places au-dessus desquelles $S'$ est ramifi\'e. Posons $\hat{G}^1=\hat{G}^{\hat{\theta},0}$ et $^LG^1=\hat{G}^1\rtimes W_{F}\subset {^LG}$. On commence par prouver que 

(3) le plongement $\hat{S}'\subset \hat{G}^{1}$ se prolonge en un plongement $^LS'\to {^LG^1}$ qui poss\`ede les propri\'et\'es de l'\'enonc\'e.

Notons $E$  l'extension galoisienne finie  de $F$ telle que $\Gamma_{E}$ soit l'intersection des noyaux des actions de $\Gamma_{F}$ sur $\hat{G}$ et sur $\hat{S}'$.  Cette extension n'est pas ramifi\'ee au-dessus des places hors de $V$. Le cocycle $\omega_{S',G}$ se factorise  par $Gal(E/F)$. Tout \'el\'ement de $W^{\theta}$ se relevant dans $\hat{G}^1$, on peut fixer une application $b:Gal(E/F)\to \hat{G}^1$ de sorte que $b(\sigma)$ soit un rel\`evement de $\omega_{S',G}(\sigma)$ pour tout $\sigma\in Gal(E/F)$. Le cobord
$$db(\sigma,\sigma') =b(\sigma)\sigma(b(\sigma'))b(\sigma\sigma')^{-1}$$
est un $2$-cocycle de $Gal(E/F) $ dans $\hat{S}'$. D'apr\`es  (1), on peut fixer une cocha\^{\i}ne continue $\mu:W_{F}\to \hat{S}'$ telle que $\mu(I_{v})=1$ pour tout $v\not\in V$ et $d\mu=db$. Alors le plongement 
$$\begin{array}{ccc}{^LS'}&\to &{^LG^1}\\ (t,w)&\mapsto &(t\mu(w)^{-1}b(w),w)\\ \end{array}$$
satisfait les conditions requises (on a relev\'e $b$ en une application d\'efinie sur $W_{F}$ via la projection $W_{F}\to Gal(E/F)$).

On fixe un plongement $(t,w)\mapsto (tg^1(w),w)$ v\'erifiant (3).  Rappelons que, pour tout $w\in W_{F}$, il existe $g_{w}=(g(w),w)\in {\cal G}'$ tel que $ad_{g_{w}}$ pr\'eserve la paire \'epingl\'ee de $\hat{G}'$ que l'on a fix\'ee plus haut. L'\'el\'ement $g(w)$ normalise $\hat{T}$, son image $\omega_{G'}(w)$ dans $W$ est fixe par $\theta$. Pour tout $w\in W_{F}$, on a une \'egalit\'e $\omega_{S',G}(w)=\omega_{S',G'}(w)\omega_{G'}(w)$, o\`u $\omega_{S',G'}(w)\in W^{G'}$. Il en r\'esulte que, pour tout $w$, on peut choisir $(x(w),w)\in {\cal G}'$ tel que $x(w)$ normalise $\hat{T}$ et que son image dans $W$ soit $\omega_{S',G}(w)$. Notons que $x(w)$ est  bien d\'etermin\'e \`a multiplication \`a gauche pr\`es par un \'el\'ement de $\hat{S}'$. Il existe $t(w)\in \hat{T}$ tel que $x(w)=t(w)g^1(w)$. L'image $d(w)$ de $t(w)$ dans $\hat{T}/\hat{S}'$ est uniquement d\'etermin\'ee. Montrons que l'application $d:W_{F}\to \hat{T}/\hat{S}'$ est continue. Il suffit de montrer que, pour tout $w_{0}$, on peut choisir $x(w)$ continu au voisinage de $w_{0}$. Par hypoth\`ese, on peut choisir une section continue $w\mapsto s_{w}$ de la projection ${\cal G}'\to W_{F}$. En g\'en\'eral, l'\'el\'ement $s_{w}$ ne conserve pas la paire $(\hat{B}',\hat{T}')$. Mais la paire $s_{w}(\hat{B}',\hat{T}')$ varie continuement en $w$.  Dans un voisinage de $w_{0}$, on peut donc fixer une application continue $w\mapsto h(w)$ \`a valeurs dans $\hat{G}'$ de sorte que $ad_{h(w)}ad_{s_{w}}(\hat{B}',\hat{T}')=\omega_{S',G'}(w)(\hat{B}',\hat{T}')$. En posant $h(w)s_{w}=(x(w),w)$, l'application $x$ ainsi d\'efinie au voisinage de $w_{0}$ convient. Cela prouve la continuit\'e de $d$. On peut munir le tore $\hat{T}$ de l'action galoisienne $\sigma\mapsto \omega_{S',G}(\sigma)\circ\sigma$. Notons $\hat{S}$ ce tore muni de cette action. On peut consid\'erer $d$ comme une application \`a valeurs dans $\hat{S}/\hat{S}'$. On v\'erifie qu'alors, c'est un cocycle. De plus, pour $v\not\in V$, $d$ est triviale sur $I_{v}$. En effet, pour $w\in I_{v}$, on a $g^1(w)=1$ et l'hypoth\`ese $v\not\in V_{ram}({\bf G}')$ nous permet de choisir $x(w)=1$, d'o\`u $d(w)=1$. On consid\`ere la suite exacte
$$1\to \hat{S}'\to \hat{S}\to \hat{S}/\hat{S}'\to 1$$
D'apr\`es (2), on peut relever $d$ en un cocycle $e:W_{F}\to \hat{S}$ tel que, pour tout $v\not\in V$, $e$ soit trivial sur $I_{v}$.  Fixons un tel cocycle. Pour tout $w\in W_{F}$, posons $y(w)=e(w)g^1(w)$. Alors $(y(w),w)$ appartient \`a ${\cal G}'$. L'application
$$\begin{array}{ccc}{^LS'}&\to &{\cal G}'\\ (t,w)&\mapsto& (te(w),w)\\ \end{array}$$
v\'erifie les conditions du lemme. $\square$

\bigskip

\subsection{Donn\'ees auxiliaires}
Soit ${\bf G}'=(G',{\cal G}',\tilde{s})$ une donn\'ee endoscopique de $(G,\tilde{G},{\bf a})$. On suppose $\tilde{G}'(F)\not=\emptyset$. La notion de donn\'ees auxiliaires $G'_{1}$, $\tilde{G}'_{1}$, $C_{1}$, $\hat{\xi}_{1}$ se d\'efinit comme dans le cas local ([I] 2.1). Toujours comme dans le cas local, \`a de telles donn\'ees est associ\'e un caract\`ere $\lambda_{1}$ qui est  cette fois  automorphe, c'est-\`a-dire un caract\`ere du groupe $C_{1}({\mathbb A})/C_{1}(F)$.

\ass{Lemme}{On peut choisir de telles donn\'ees auxiliaires de sorte que, pour $v\not\in V_{ram}({\bf G}')$, le groupe $G'_{1,v}$ soit non ramifi\'e, le plongement $\hat{\xi}_{1,v}$ soit l'identit\'e sur $I_{v}$ et $\lambda_{1}$ soit non ramifi\'e en $v$.}

{\bf Remarque.}  C'est le lemme 7.1 de [A1]. Nous reprenons la d\'emonstration car Arthur formule les conditions de non-ramification un peu diff\'eremment de nous.

Preuve. Fixons une paire de Borel \'epingl\'ee $\hat{{\cal E}}$ de $\hat{G}$ conserv\'ee par l'action galoisienne. Il s'en d\'eduit un automorphisme $\hat{\theta}$. Notons $(\hat{B},\hat{T})$ la paire de Borel sous-jacente \`a $\hat{{\cal E}}$. On peut supposer que $ad_{\tilde{s}}$ conserve cette paire. Posons $\hat{B}'=\hat{B}\cap \hat{G}'$, $\hat{T}'=\hat{T}^{\hat{\theta},0}$. C'est une paire de Borel de $\hat{G}'$ et on peut supposer qu'elle est conserv\'ee par l'action galoisienne $\sigma\mapsto \sigma_{G'}$. Posons 
$$\hat{G}'_{1}=(\hat{G}'\times \hat{T}')/diag(Z(\hat{G}')).$$
Ce groupe est \`a centre connexe (isomorphe \`a $\hat{T}'$) et s'ins\`ere dans une extension
$$1\to \hat{G}'\to \hat{G}'_{1}\to \hat{T}'_{ad}\to 1.$$
Le tore $\hat{T}'_{ad}$ est induit. Alors $\hat{G}'_{1}$ v\'erifie les conditions requises pour appliquer le lemme 2.2.A de [KS]: l'inclusion $\hat{G}'\to \hat{G}'_{1}$ se prolonge en un plongement $\hat{\tau}:{\cal G}'\to {^LG'_{1}}$. Remarquons que l'action galoisienne sur $\hat{G}'_{1}$ est non ramifi\'ee hors de $V_{ram}({\bf G}')$. Appliquons le lemme 3.2 au tore $\hat{T}'$. Soit $\hat{\iota}:{^LT'}\to {\cal G}'$ un plongement v\'erifiant les propri\'et\'es de ce lemme. Pour $w\in W_{F}$, posons $\hat{\tau}\circ \hat{\iota}(w)=(x(w),w)$.  Cet \'el\'ement agit comme $w$ sur $\hat{T}'$, donc $x(w)$ appartient au commutant $\hat{T}_{1}$ de $\hat{T}'$ dans $\hat{G}'_{1}$. 
L'application $w\mapsto x(w)$ est un cocycle de $W_{F}$ dans $\hat{T}_{1}$. Pour $v\not\in V_{ram}({\bf G}')$ et $w\in I_{v}$, on a $\hat{\iota}(w)=(1,w)\in {^LG}$, donc $\hat{\iota}(w)$ commute \`a $\hat{G}'$. Donc $(x(w),w)$ puis $x(w)$ commutent aussi \`a $\hat{G}'$. Il en r\'esulte que $x(w)\in Z(\hat{G}'_{1})$. Notons $t_{3}$ le compos\'e de $x$ et de la projection $\hat{T}_{1}\to \hat{T}_{1}/Z(\hat{G}'_{1})$. Alors $t_{3}$ est un cocycle trivial sur $I_{v}$ pour tout $v\not\in V_{ram}({\bf G}')$. Puisque $Z(\hat{G}'_{1})$ est connexe, on peut appliquer 3.2(2) (on prend pour $E$ la plus petite extension sur laquelle $\hat{T}'$ est d\'eploy\'ee et pour $U$ le sous-groupe de $W_{F}$ engendr\'e par les $I_{v}$ pour $v\not\in V_{ram}({\bf G}')$). Il existe donc un cocycle $t_{2}:W_{F}\to \hat{T}_{1}$ relevant $t_{3}$ et non ramifi\'e hors de $V_{ram}({\bf G}')$. Pour $w\in W_{F}$, posons $\zeta(w)=x(w)^{-1}t_{2}(w)$. Alors $\zeta$ est un cocycle \`a valeurs dans $Z(\hat{G}'_{1})$. Rempla\c{c}ons le plongement $\hat{\tau}$ par $\hat{\xi}_{1}$ d\'efini par
$$\hat{\xi}_{1}(g,w)=\zeta(w)\hat{\tau}(g,w)$$
pour tout $(g,w)\in {\cal G}'$. On a alors $\hat{\xi}_{1}\circ \hat{\iota}(w)=(t_{2}(w),w)$ pour tout $w\in W_{F}$. Pour $v\not\in V_{ram}({\bf G}')$ et $w\in I_{v}$, on a $\hat{\iota}(w)=w$ et $t_{2}(w)=1$, donc $\hat{\xi}_{1}(w)=w$. C'est-\`a-dire que $\hat{\xi}_{1}$  est l'identit\'e sur $I_{v}$. On prend pour $G'_{1}$ le groupe quasi-d\'eploy\'e sur $F$ dont $\hat{G}'_{1}$ est le groupe dual. Alors les deux premi\`eres conditions   du lemme sont satisfaites. La condition de non-ramification de $\lambda_{1}$ r\'esulte formellement de ces deux premi\`eres conditions.

La d\'efinition d'un espace tordu $\tilde{G}'_{1}$ ne pose pas de probl\`eme. On fixe $\gamma\in \tilde{G}'(F)$. L'automorphisme $ad_{\gamma}$ se rel\`eve en un automorphisme $\theta_{1}$ de $G'_{1}$ qui est d\'efini sur $F$. On prend pour $\tilde{G}'_{1}$ le sous-ensemble $G'_{1}\times \{\theta_{1}\}$ du produit semi-direct de $G'_{1}$ avec son groupe d'automorphismes. La projection $\tilde{G}'_{1}\to \tilde{G}'$ est d\'efinie en envoyant $\theta_{1}$ sur $\gamma$. $\square$

Soit $V$ un ensemble fini de places de $F$ contenant $V_{ram}({\bf G}')$.  Soient $G'_{1}$ etc... des donn\'ees auxiliaires telles que  les conditions de l'\'enonc\'e soient satisfaites pour $v\not\in V$.  Pour $v\not\in V$, le   sous-groupe hypersp\'ecial $K'_{v}$ de $G'(F_{v})$ d\'etermine un tel sous-groupe $K'_{1,v}$ de $G'_{1}(F_{v})$. On fixe   un sous-espace hypersp\'ecial $\tilde{K}'_{1,v}$ de $\tilde{G}'_{1}(F_{v})$ au-dessus de $\tilde{K}'_{v}$. On suppose que ces  sous-espaces v\'erifient la m\^eme condition de compatibilit\'e globale qu'en 1.1.  On adjoint cette famille de sous-espaces hypersp\'eciaux aux donn\'ees auxiliaires et on appelle donn\'ees auxiliaires non ramifi\'ees hors de $V$ les donn\'ees $G'_{1}$, $\tilde{G}'_{1}$, $C_{1}$, $\hat{\xi}_{1}$, $(\tilde{K}'_{1,v})_{v\not\in V}$.  Si $V'$ est un autre ensemble fini de places contenant $V$, des donn\'ees auxiliaires non ramifi\'ees hors de $V$ se restreignent en des donn\'ees non ramifi\'ees hors de $V'$ en oubliant les $\tilde{K}'_{1,v}$ pour $v\in V'-V$.  

Consid\'erons deux s\'eries de donn\'ees auxiliaires $G'_{1}$ etc... et $G'_{2}$ etc... non ramifi\'ees hors de $V$. On d\'efinit comme en [I] 2.5 le produit fibr\'e $G'_{12}$ de $G'_{1}$ et $G'_{2}$ au-dessus de $G'$ et le produit fibr\'e analogue $\tilde{G}'_{12}$. Toujours comme en [I] 2.5, on d\'efinit un caract\`ere $\lambda_{12}$ qui est cette fois un caract\`ere de $G'_{12}({\mathbb A})$ trivial sur $G'_{12}(F)$. Il est non ramifi\'e hors de $V$, donc trivial sur $K'_{1,v}$ pour $v\not\in V$. 
 
Comme en 1.15, on  d\'eduit de $\lambda_{12}$ des fonctions $\tilde{\lambda}_{12}$ sur $\tilde{G}'_{12}({\mathbb A}_{F})$, $\tilde{\lambda}_{12,v}$ sur $\tilde{G}'_{12}(F_{v})$ pour $v\not\in V$ et $\tilde{\lambda}_{12,V}$ sur $\tilde{G}'_{12}(F_{V})$. Cette derni\`ere fonction permet de recoller les espaces $C_{c,\lambda_{1}}^{\infty}(\tilde{G}'_{1}(F_{V}))$ et $C_{c,\lambda_{1}}^{\infty}(\tilde{G}'_{2}(F_{V}))$: \`a $f_{1}\in C_{c,\lambda_{1}}^{\infty}(\tilde{G}'_{1}(F_{V}))$, on associe $f_{2}\in C_{c,\lambda_{1}}^{\infty}(\tilde{G}'_{2}(F_{V}))$ telle que $f_{2}(\delta_{2})=\tilde{\lambda}_{12,V}(\delta_{1},\delta_{2})f_{1}(\delta_{1})$ o\`u $\delta_{1}$ est n'importe quel \'el\'ement de $\tilde{G}'_{1}(F_{V})$ tel que $(\delta_{1},\delta_{2})\in \tilde{G}'_{12}(F_{V})$. Ce recollement v\'erifie une propri\'et\'e de transitivit\'e \'evidente qui permet de d\'efinir un espace $C_{c}^{\infty}({\bf G}'_{V})$ comme la limite inductive des $C_{c,\lambda_{1}}^{\infty}(\tilde{G}'_{1}(F_{V}))$ sur les donn\'ees $G'_{1}$,...,$(\tilde{K}'_{1,v})_{v\not\in V}$ non ramifi\'ees hors de $V$, les applications de transition \'etant celles que l'on vient de d\'efinir. Cette d\'efinition pose le m\^eme probl\`eme logique que dans le cas local, que l'on peut lever comme dans ce cas, cf. [I] 2.5. On d\'efinit de m\^eme les espaces $I({\bf G}'_{V})$ et $SI({\bf G}'_{V})$. De nouveau, si $V'$ est un ensemble fini de places contenant $V$, les espaces $C_{c}^{\infty}({\bf G}'_{V})$ etc... s'identifient \`a des sous-espaces de $C_{c}^{\infty}({\bf G}'_{V'})$ etc...

\bigskip

\subsection{Levi}
Les relations entre donn\'ees endoscopiques d'espaces de Levi de $\tilde{G}$ et groupes de Levi de donn\'ees endoscopiques de $(G,\tilde{G},{\bf a})$ sont essentiellement les m\^emes dans le cas global que dans le cas local, cf. [I] paragraphes 3.2, 3.3 et 3.4. Signalons l'analogue global de la relation 3.2(2) de [I]. Soit $\tilde{M}$ un espace de Levi de $\tilde{G}$. On r\'ealise $\hat{M}$ comme Levi standard de $\hat{G}$ comme dans cette r\'ef\'erence. Alors l'homomorphisme
$$H^1(W_{F},Z(\hat{G}))/ker^1(W_{F},Z(\hat{G}))\to H^1(W_{F},Z(\hat{M}))/ker^1(W_{F},Z(\hat{M}))$$
est injectif ([A11],  lemme 2).

Soit ${\bf G}'=(G',{\cal G}',\tilde{s})$ une donn\'ee endoscopique  de $(G,\tilde{G},{\bf a})$.  On suppose $\tilde{G}'(F)\not=\emptyset$. Soit $M'$ un Levi de $G'$ contenant $M'_{0}$, dont on d\'eduit un espace de Levi $\tilde{M}'$. Puisque $\tilde{G}'$ est \`a torsion int\'erieure, $\tilde{M}'$ est l'ensemble des $\gamma\in \tilde{G}$ tels qu'il existe $m\in M$ de sorte que $ad_{\gamma}=ad_{m}$.  Pour une place $v\not\in V_{ram}({\bf G}')$, on d\'efinit le sous-groupe hypersp\'ecial $K^{M'}_{v}=M'(F_{v})\cap K'_{v}$ et l'espace hypersp\'ecial $\tilde{K}^{M'}_{v}=\tilde{M}'(F_{v})\cap \tilde{K}'_{v}$. On construit comme en [I] 3.4 un sous-groupe ${\cal M}'\subset {\cal G}'$ qui est une extension de $\hat{M}'$ par $W_{F}$. On pose ${\bf M}'=(M',{\cal M}',\tilde{s})$.   Il peut exister un espace de Levi $\tilde{M}$ de $\tilde{G}$ de sorte que ${\bf M}'$ s'identifie \`a une donn\'ee endoscopique de $(M,\tilde{M},{\bf a}_{M})$. Mais, comme dans le cas local, un tel espace de Levi n'existe pas toujours. Toutefois, ind\'ependamment de l'existence d'un tel espace, pour un ensemble fini de places $V$ contenant $V_{ram}({\bf G}')$, on peut d\'efinir des espaces $C_{c}^{\infty}({\bf M}'_{V})$, $I({\bf M}'_{V})$ etc... En effet, on n'a pas besoin qu'il existe un espace $\tilde{M}$ pour d\'efinir la notion de donn\'ees auxiliaires de ${\bf M}'$ non ramifi\'ees hors de $V$ ni pour d\'efinir des fonctions de recollement $\tilde{\lambda}_{12}$. Et cela suffit pour d\'efinir les espaces pr\'ec\'edents.

 Soit $\tilde{M}'$ un espace de Levi de $\tilde{G}'$. On doit fixer une mesure sur ${\cal A}_{\tilde{M}'}$ (qui est d'ailleurs \'egal \`a ${\cal A}_{M'}$ puisque la torsion sur $\tilde{G}'$ est int\'erieure).  Comme on vient de le dire,  $\tilde{M}'$ peut correspondre ou non \`a un espace de Levi $\tilde{M}$. Dans le cas o\`u $\tilde{M}'$ ne correspond \`a aucun espace de Levi $\tilde{M}$, la mesure sur ${\cal A}_{\tilde{M}'}$ ne nous importera pas, on la choisit arbitrairement. Dans le cas o\`u $\tilde{M}'$ correspond \`a un espace de Levi $\tilde{M}$, on a un isomorphisme naturel ${\cal A}_{\tilde{M}'}\simeq {\cal A}_{\tilde{M}}$ et on choisit la mesure sur le premier espace qui correspond par cet isomorphisme \`a celle fix\'ee sur le second. L'espace de Levi $\tilde{M}$ n'est bien d\'efini qu'\`a conjugaison pr\`es, mais notre d\'efinition est insensible \`a une telle conjugaison.

\bigskip

\subsection{La partie g\'eom\'etrique de la formule des traces invariante pour une donn\'ee endoscopique}

Soit ${\bf G}'=(G',{\cal G}',\tilde{s})$ une donn\'ee endoscopique  de $(G,\tilde{G},{\bf a})$.  On suppose $\tilde{G}'(F)\not=\emptyset$. Soit $V$ un ensemble fini de places de $F$ contenant $V_{ram}({\bf G}')$. Fixons des donn\'ees auxiliaires $G'_{1}$, $\tilde{G}'_{1}$, $C_{1}$, $\hat{\xi}_{1}$, $(\tilde{K}'_{1,v})_{v\not\in V}$ non ramifi\'ees hors de $V$. On construit comme en 2.8 une distribution ${\bf f}_{V}\mapsto I^{\tilde{G}'_{1}}_{g\acute{e}om,\lambda_{1}}({\bf f}_{V})$ sur $C_{c,\lambda_{1}}^{\infty}(\tilde{G}'_{1}(F_{V}))\otimes Mes(G'(F_{V}))$. Pour ${\cal O}\in \tilde{G}'_{ss}(F)/conj$, on construit une distribution $A^{\tilde{G}'_{1}}_{\lambda_{1}}(V,{\cal O})\in D_{orb,\lambda_{1}}(\tilde{G}'_{1}(F_{V}))\otimes Mes(G'(F_{V}))^*$. Si on change de donn\'ees auxiliaires, le lemme 2.6 et la relation 2.8(3) disent que ces distributions se recollent selon l'isomorphisme d\'efini au paragraphe pr\'ec\'edent. On peut donc les consid\'erer comme des distributions sur les espaces $C_{c}^{\infty}({\bf G}'_{V})\otimes Mes(G'(F_{V}))$. On les note alors $I^{{\bf G}'}_{g\acute{e}om}$ et $A^{{\bf G}'}(V,{\cal O})$.  Soit $M'\in {\cal L}(M'_{0})$. Notons  ${\bf M}'=(M',{\cal M}',\tilde{s})$ la donn\'ee introduite dans le paragraphe pr\'ec\'edent. On a une forme bilin\'eaire sur
$$(D_{orb,\lambda_{1}}(\tilde{M}'_{1}(F_{V}))\otimes Mes(M'(F_{V}))^*)\times (C_{c,\lambda_{1}}^{\infty}(\tilde{G}'_{1}(F_{V}))\otimes Mes(G'(F_{V})))$$
qui, \`a $(\boldsymbol{\gamma}_{V},{\bf f}_{V})$, associe $I_{\tilde{M}'_{1},\lambda_{1}}^{\tilde{G}'_{1}}(\boldsymbol{\gamma}_{V},{\bf f}_{V})$. Comme ci-dessus, quand on change de donn\'ees auxiliaires, ces formes lin\'eaires se recollent. On obtient une forme bilin\'eaire
$$(\boldsymbol{\gamma}_{V},{\bf f}_{V})\mapsto I_{{\bf M}'}^{{\bf G}'}(\boldsymbol{\gamma}_{V},{\bf f}_{V})$$
sur
$$(D_{orb}({\bf M}'_{V})\otimes Mes(M'(F_{V}))^*)\times(C_{c}^{\infty}({\bf G}'_{V})\otimes Mes(G'(F_{V}))).$$

Il r\'esulte de 2.8(1) que l'on a l'\'egalit\'e
$$I^{{\bf G}'}_{g\acute{e}om}({\bf f}_{V})=\sum_{\tilde{M}'\in {\cal L}(\tilde{M}'_{0})}\vert W^{\tilde{M}'}\vert \vert W^{\tilde{G}'}\vert ^{-1}\sum_{{\cal O}\in \tilde{M}'_{ss}(F_{V})/conj}I_{{\bf M}'}^{{\bf G}'}(A^{{\bf M}'}(V,{\cal O}),{\bf f}_{V})$$
pour tout ${\bf f}_{V}\in C_{c}^{\infty}({\bf G}'_{V})\otimes Mes(G'(F_{V}))$.

 \bigskip

\subsection{Facteur de transfert global, cas particulier}
Soit ${\bf G}'=(G',{\cal G}',\tilde{s})$ une donn\'ee endoscopique relevante de $(G,\tilde{G},{\bf a})$.  Rappelons que, puisque $\tilde{G}'$ est \`a torsion int\'erieure,  \`a tout tore maximal $T'$  de $G'$ d\'efini sur $F$  est associ\'e un unique tore tordu $\tilde{T}'$, \`a savoir l'ensemble des $\delta\in \tilde{G}'$ qui commutent \`a tout \'el\'ement de $T'$. On peut toutefois avoir $\tilde{T}'(F)=\emptyset$. On impose dans ce paragraphe l'hypoth\`ese

{\bf (Hyp)} {\it  il existe un sous-tore maximal $T'$ de $G'$, d\'efini sur $F$,  de sorte que, pour toute place $v$ de $F$, il existe un couple $(\delta_{v},\gamma_{v})\in {\cal D}({\bf G}'_{v})$ avec $\delta_{v}\in \tilde{T}'(F_{v})$. } 

Fixons un tel tore $T'$. On fixe pour tout $v\in Val(F)$ un couple  $(\delta_{v},\gamma_{v})\in {\cal D}({\bf G}'_{v})$ avec $\delta_{v}\in \tilde{T}'(F_{v})$. Nous allons imposer \`a ces couples des conditions de "non-ramification". On impose d'abord

(1) $\delta_{v}\in \tilde{K}'_{v}$ et $\gamma_{v}\in \tilde{K}_{v}$ pour presque tout $v$.

Consid\'erons une place finie $v$. Fixons une paire de Borel \'epingl\'ee ${\cal E}_{v}=(B_{v},T_{v},(E_{\alpha})_{\alpha\in \Delta})$ de $G$ telle que $ (B_{v},T_{v})$ soit conserv\'ee par $ad_{\gamma_{v}}$. Ecrivons $\gamma_{v}=t_{v}e_{v}$ avec $e_{v}\in Z(\tilde{G},{\cal E})$ et $t_{v}\in T_{v}$. Rappelons que l'on note $\Sigma(T_{v})$ l'ensemble des racines de $T_{v}$ dans $\mathfrak{g}$, $F_{v}^{nr}$ la plus grande extension non ramifi\'ee de $F_{v}$, $\mathfrak{o}_{v}^{nr}$ son anneau d'entiers et $\mathfrak{o}_{v}^{nr,\times}$ le groupe des unit\'es. Notons aussi ${\mathbb F}_{v}$ le corps r\'esiduel de $\mathfrak{o}_{v}$ et $\bar{{\mathbb F}}_{v}$ sa cl\^oture alg\'ebrique, qui est aussi le corps r\'esiduel de $\mathfrak{o}_{v}^{nr}$. Avec ces notations, on impose

(2) pour presque toute place finie $v$ et pour tout $\alpha\in \Sigma(T_{v})$, on a $(N\alpha)(t_{v})\in \mathfrak{o}_{v}^{nr,\times}$ et la r\'eduction de $(N\alpha)(t_{v}) $ dans $\bar{{\mathbb F}}_{v}$ est diff\'erente de $\pm 1$. 

Cette condition ne d\'epend pas des choix de ${\cal E}_{v}$ et de $e_{v}$. En effet, on ne peut changer $e_{v}$ qu'en le multipliant par un \'el\'ement de $Z(G)$. Cela multiplie $t_{v}$ par l'inverse de cet \'el\'ement, ce qui ne change pas les $(N\alpha)(t_{v})$. Rempla\c{c}ons ${\cal E}_{v}$ par une autre paire de Borel \'epingl\'ee ${\cal E}'_{v}$  dont la paire de Borel sous-jacente est conserv\'ee par $ad_{\gamma_{v}}$. Le tore $T_{v}$ reste le m\^eme: puisque $\gamma_{v}$ est fortement r\'egulier, c'est le commutant de $G_{\gamma_{v}}$ dans $G$.  Soit $x\in G$ tel que $ad_{x}({\cal E}_{v})={\cal E}'_{v}$. Alors $ad_{x}$ conserve $T_{v}$.  D'apr\`es [I] 1.3(2),  l'image de $x$ dans $W$  (identifi\'e au groupe de Weyl de $G$ pour $T_{v}$) est invariante par $\theta=ad_{e_{v}}$. Un tel \'el\'ement se rel\`eve dans le groupe  $G_{e_{v}}$. On peut donc \'ecrire $x=\tau y$ avec $\tau \in T_{v}$ et $y\in G_{e_{v}}$. Alors $e'_{v}=ad_{\tau}(e_{v})=(1-\theta)(\tau)e_{v}\in Z(\tilde{G},{\cal E}'_{v})$. On peut prendre pour d\'ecomposition relative \`a ${\cal E}'_{v}$ l'\'egalit\'e $\gamma_{v}= t'_{v}e'_{v}$ o\`u $t'_{v}=(\theta-1)(\tau)t_{v}$. Puisque $(N\alpha)((\theta-1)(\tau))=1$, on voit que la condition (2) ne change pas.

Nous montrerons plus loin qu'il existe des familles $(\delta_{v})_{v\in Val(F)}$ et $(\gamma_{v})_{v\in Val(F)}$ v\'erifiant les deux conditions (1) et (2) ci-dessus.  

Fixons un ensemble fini $V$ de places de $F$ contenant $V_{ram}({\bf G}')$, ainsi que des donn\'ees auxiliaires $G'_{1}$, $\tilde{G}'_{1}$, $C_{1}$, $\hat{\xi}_{1}$, $(\tilde{K}'_{1,v})_{v\not\in V}$, non ramifi\'ees hors de $V$.   On fixe pour tout $v$ un rel\`evement $\delta_{1,v}\in \tilde{G}'_{1}(F_{v})$ de $\delta_{v}$. La condition (1) permet d'imposer
que $\delta_{1,v}\in \tilde{K}'_{1,v}$ pour presque tout $v$. 

Posons $\delta_{1}=(\delta_{1,v})_{v\in Val(F)}$, $\delta=(\delta_{v})_{v\in Val(F)}$, $\gamma=(\gamma_{v})_{v\in Val(F)}$. Ce sont des \'el\'ements de $\tilde{G}'_{1}({\mathbb A}_{F})$, resp. $\tilde{G}'({\mathbb A}_{F})$, $\tilde{G}({\mathbb A}_{F})$. Nous allons d\'efinir un facteur de transfert global $\Delta(\delta_{1},\gamma)$.

   On fixe   un sous-groupe de Borel $B'$ de $G'$, d\'efini sur $\bar{F}$ et contenant $T'$. Identifions $\underline{la}$ paire de Borel \'epingl\'ee de $G$ \`a une telle paire ${\cal E}^*=(B^*,T^*,(E^*_{\alpha})_{\alpha\in \Delta})$ d\'efinie sur $\bar{F}$.  On fixe une application $\sigma\mapsto u_{{\cal E}^*}(\sigma)$ \`a valeurs dans $G_{SC}(\bar{F})$,  de sorte que $ad_{u_{{\cal E}^*}(\sigma)}\circ \sigma$ conserve ${\cal E}^*$ pour tout $\sigma\in \Gamma_{F}$. On peut supposer que cette application est continue, qu'elle se factorise par un quotient fini de $\Gamma_{F}$ et que $u_{{\cal E}^*}(1)=1$. On construit  une action quasi-d\'eploy\'ee de $\Gamma_{F}$ sur $G$ not\'ee $\sigma\mapsto \sigma_{G^*}=ad_{u_{{\cal E}^*}(\sigma)}\circ\sigma $ (on note simplement $\sigma\mapsto \sigma$ l'action naturelle de $\Gamma_{F}$ sur $G$). Des deux paires de Borel se d\'eduit un homomorphisme $\xi_{T^*,T'}:T^*\to T'$. Il existe un cocycle $\omega_{T}:\Gamma_{F}\to W^{\theta}$ de sorte que $\xi_{T^*,T'}\circ \omega_{T}(\sigma)\circ\sigma_{G^*}=\sigma_{G'}\circ\xi_{T^*,T'}$ sur $T^*$ pour tout $\sigma\in \Gamma_{F}$.   Fixons de plus $e\in Z(\tilde{G},{\cal E}^*;\bar{F})$ et posons $\theta^*=ad_{e}$. D'apr\`es [K1] corollaire  2.2, on peut trouver $x\in G_{SC}^{\theta^*}(\bar{F})$ tel que, pour tout $\sigma\in \Gamma_{F}$, $x\sigma_{G^*}(x)^{-1}$ normalise $T^*$ et ait pour image $\omega_{T}(\sigma)$ dans $W$. On note $T$ le tore $T^*$ muni de l'action galoisienne $\sigma\mapsto \sigma_{T}=\omega_{T}(\sigma)\circ \sigma_{G^*}$.

Fixons une extension galoisienne finie $E$ de $F$ telle que

- $B'$ et $T'$ soient d\'efinis sur $E$ et $T'$ soit d\'eploy\'e sur $E$;

- ${\cal E}^*$ soit d\'efinie sur $E$ (pour l'action galoisienne naturelle) et $T^*$ soit d\'eploy\'e sur $E$;

- $e\in Z(\tilde{G},{\cal E}^*;E)$, $x\in G_{SC}^{\theta^*}(E)$;

- l'application $\sigma\mapsto u_{{\cal E}^*}(\sigma)$  se factorise par $Gal(E/F)$ et, pour tout $\sigma\in Gal(E/F)$, on a $u_{{\cal E}^*}(\sigma)\in G_{SC}(E)$.

Il en r\'esulte que l'action galoisienne quasi-d\'eploy\'ee et l'action naturelle co\"{\i}ncident sur $\Gamma_{E}$. On  a $u_{{\cal E}^*}(\sigma)=1$ pour $\sigma\in \Gamma_{E}$. La paire de Borel $(B^*,T^*)$ est aussi d\'efinie sur $E$ et les tores $T^*$ et $T$ sont d\'eploy\'es sur $E$. Rappelons que, pour tout $\sigma\in \Gamma_{F}$, il existe un unique $z(\sigma)\in Z(G;\bar{F})$ tel que $ \sigma_{G^*}(e)=z(\sigma)^{-1}e$.  L'application $\sigma\mapsto z(\sigma)$ se factorise par $Gal(E/F)$ et prend ses valeurs dans $Z(G;E)$.

On note ${\mathbb A}_{E}$ l'anneau des ad\`eles de $E$. Le groupe $Gal(E/F)$ agit sur ${\mathbb A}_{E}$. Pour tout groupe alg\'ebrique $H$ d\'efini sur $F$, le groupe $Gal(E/F)$ agit sur $H({\mathbb A}_{E})$. Soit $v\in Val(F)$, posons $E_{v}=\prod_{w\vert  v}E_{w}$ o\`u $w$ parcourt les places de $E$ divisant $v$. Alors le groupe $Gal(E/F)$ agit sur $E_{v}$. Pour tout groupe alg\'ebrique $H_{v}$ d\'efini sur $F_{v}$, le groupe $Gal(E/F)$ agit sur $H_{v}(E_{v})$.

Soit $v\in Val(F)$. Comme on l'a dit en 1.1, notre notion de localisation en  $v$ d\'epend du choix d'un prolongement $\bar{v}$ de $v$ \`a $\bar{F}$. Le corps $\bar{F}_{v}$ est la cl\^oture alg\'ebrique de $F_{v}$ dans le compl\'et\'e de $\bar{F}$ en $\bar{v}$. Par abus de notation, notons-le plus pr\'ecis\'ement $\bar{F}_{\bar{v}}$. Le groupe $\Gamma_{F_{v}}$ est le fixateur de $\bar{v}$ dans $\Gamma_{F}$, notons-le plus pr\'ecis\'ement $\Gamma_{\bar{v}}$.
 Ici, parce que l'on va travailler avec le corps $E$, on va devoir faire varier $\bar{v}$. Notons $w$ la restriction de $\bar{v}$ \`a $E$.
En reprenant la preuve du lemme 1.10 de [I], on voit que  pour tout $v\in Val(F)$, on peut fixer un diagramme $(\delta_{v},B',T',B_{w},T_{v},\gamma_{v})$ o\`u $(B',T')$ est la paire d\'ej\`a fix\'ee. 
Le groupe $B_{w}$ est d\'efini sur $\bar{F}_{\bar{v}}$ et la condition d'\'equivariance du diagramme est relative \`a $\Gamma_{\bar{v}}$. 

{\bf Remarques.} (3) Ce diagramme est unique. En effet, comme on l'a d\'ej\`a dit, $T_{v}$ est uniquement d\'etermin\'e et $B_{w}$ est en tout cas bien d\'etermin\'e modulo l'action de $W^{\theta}$ (en identifiant $W$ au groupe de Weyl de $G$ relatif \`a $T_{v}$). Du diagramme se d\'eduit un homomorphisme $\xi_{T_{v},T'}:T_{v}\to T'$, puis une application 
$$\tilde{\xi}_{T_{v},T'}:(T_{v}/(1-\theta)(T_{v}))\times_{{\cal Z}(G)}{\cal Z}(\tilde{G})\to T'\times_{{\cal Z}(G')}{\cal Z}(\tilde{G}').$$
Cette application doit envoyer l'image de $\gamma_{v}$ dans l'espace de d\'epart sur l'image de $\delta_{v}$ dans l'espace d'arriv\'ee. Si l'on remplace $B_{w}$ par $\omega(B_{w})$, avec $\omega\in W^{\theta}$, $\tilde{\xi}_{T_{v},T'}$ est remplac\'e par $\tilde{\xi}_{T_{v},T'}\circ \omega^{-1}$. La forte r\'egularit\'e de $\gamma_{v}$ et la propri\'et\'e [I] 1.3(5) entra\^{\i}nent que cette nouvelle application ne v\'erifie plus la propri\'et\'e pr\'ec\'edente, sauf si $\omega=1$.

(4) Fixons $g\in G$ tel que $ad_{g}(B_{w},T_{v})=(B^*,T^*)$. Alors $ad_{g}$ identifie $T_{v}$ au tore $T$, plus exactement \`a son localis\'e en la place $v$. 

 \bigskip
 Soit $w'$ une autre place de $E$ au-dessus de $v$. Fixons $\tau\in \Gamma_{F}$ tel que $\tau(w)=w'$. Notons $\bar{v}'$ l'image de $\bar{v}$ par $\tau$.  L'\'el\'ement $\tau$ d\'efinit naturellement un isomorphisme de $\bar{F}_{\bar{v}}$ sur $\bar{F}_{\bar{v}'}$ et, pour tout groupe alg\'ebrique $H_{v}$ d\'efini sur $F_{v}$, un isomorphisme de $H_{v}(\bar{F}_{\bar{v}})$ sur $H_{v}(\bar{F}_{\bar{v}'})$. On note encore $\tau$ ces isomorphismes. Puisque $T_{v}$ et $\gamma_{v}$ sont d\'efinis sur $F_{v}$ donc  invariants par $\tau$, le couple $(\tau(B_{w}),T_{v})$ est une paire de Borel de $G$ d\'efinie sur $\bar{F}_{\bar{v}'}$ et conserv\'ee par $ad_{\gamma_{v}}$. En identifiant gr\^ace \`a cette paire le groupe de Weyl $W$ au groupe de Weyl de $T_{v}$, on d\'efinit le sous-groupe de Borel $B_{w'}=\omega_{T}(\tau)^{-1}\tau(B_{w})$.  Puisque $\omega_{T}(\tau)\in W^{\theta}$, la paire $(B_{w'},T_{v})$ est encore conserv\'ee par $ad_{\gamma_{v}}$, cf. [I] 1.3(2). Montrons que

(5) le sextuplet  $(\delta_{v},B',T',B_{w'},T_{v},\gamma_{v})$ est un diagramme, le corps $\bar{F}_{v}$ \'etant identifi\'e \`a $\bar{F}_{\bar{v}'}$.

Preuve. 
Les tores $T'$ et $T_{v}$ sont   d\'eploy\'es sur $E_{w}$. Cela implique que le groupe de Borel $B_{w}$ est d\'efini sur $E_{w}$ ($B'$ aussi, mais c'est d\'ej\`a dans l'hypoth\`ese sur $E$). Les deux paires de Borel $(B_{w},T_{v})$ et $(B^*,T^*)$ \'etant toutes deux d\'efinies sur $E_{w}$, on peut fixer $g_{w}\in G_{SC}(E_{w})$ tel que $ad_{g_{w}}(B_{w},T_{v})=(B^*,T^*)$.   Posons $g_{w'}=x\tau_{G^*}(x)^{-1}u_{{\cal E}^*}(\tau)\tau(g_{w})$. C'est un \'el\'ement de $G_{SC}(E_{w'})$. Montrons que

(6) on a l'\'egalit\'e $ad_{g_{w'}}(B_{w'},T_{v})=(B^*,T^*)$. 

  Puisque $ad_{g_{w}}(B_{w},T_{v})=(B^*,T^*)$, on a $ad_{\tau(g_{w})}(\tau(B_{w}),T_{v})=\tau(B^*,T^*)$. Donc 
$$ad_{u_{{\cal E}^*}(\tau)\tau(g_{w})}(\tau(B_{w}),T_{v})=ad_{u_{{\cal E}^*}(\tau)}\circ\tau (B^*,T^*)=\tau_{G^*}(B^*,T^*)=(B^*,T^*).$$
C'est donc l'isomorphisme $ad_{u_{{\cal E}^*}(\tau)\tau(g_{w})}$ qui identifie comme plus haut le groupe de Weyl $W$ au groupe de Weyl de $T_{v}$. C'est-\`a-dire que $B_{w'}$ est le groupe de Borel tel que
$$ad_{u_{{\cal E}^*}(\tau)\tau(g_{w})}(B_{w'},T_{v})=\omega_{T}(\tau)^{-1}ad_{u_{{\cal E}^*}(\tau)\tau(g_{w})}(\tau(B_{w}),T_{v})=\omega_{T}(\tau)^{-1}(B^*,T^*).$$
D'o\`u
$$\omega_{T}(\tau)ad_{u_{{\cal E}^*}(\tau)\tau(g_{w})}(B_{w'},T_{v})=(B^*,T^*).$$
Puisque $x\tau_{G^*}(x)^{-1}$ a pour image $\omega_{T}(\tau)$ dans $W$, on peut remplacer dans l'\'egalit\'e ci-dessus $\omega_{T}(\tau)$ par $ad_{x\tau_{G^*}(x)^{-1}}$ et on obtient (6).

Aux choix des paires de Borel $(B',T')$ et $(B_{w},T_{v})$, resp. $(B_{w'},T_{v})$, sont associ\'es des homomorphismes de $T_{v}$ dans $T'$ not\'es  pr\'ecis\'ement $\xi_{B_{w},T_{v},B',T'} $ et  $\xi_{B_{w'},T_{v},B',T'} $. Montrons qu'ils v\'erifient la relation

(7)  $\xi_{B_{w'},T_{v},B',T'} \circ \tau= \tau\circ \xi_{B_{w},T_{v},B',T'} $.

Les deux homomorphismes sont les restrictions \`a $T_{v}$ de $\xi_{T^*,T'}\circ ad_{g_{w'}}$, resp. $\xi_{T^*,T'}\circ ad_{g_{w}}$. Par d\'efinition de $g_{w'}$, on a 
$$ad_{g_{w'}}\circ \tau=\omega_{T}(\tau)\circ \tau_{G^*}\circ ad_{g_{w}}.$$
Par d\'efinition de $\omega_{T}(\tau)$, on a aussi
$$\xi_{T^*,T'}\circ \omega_{T}(\tau)\circ \tau_{G^*}=\tau\circ \xi_{T^*,T'}.$$
La relation (7) en r\'esulte. 

Pour prouver (5), on doit d'abord montrer que l'homomorphisme $\xi_{B_{w'},T_{v},B',T'}:T_{v}\to T'$ est \'equivariant pour l'action de $\Gamma_{\bar{v}'}$.  Pour $\sigma\in \Gamma_{\bar{v}'}$, on a
$$\xi_{B_{w'},T_{v},B',T'}\circ \sigma= \xi_{B_{w'},T_{v},B',T'}\circ \tau\circ (\tau^{-1} \sigma\tau)\circ\tau^{-1}=\tau\circ\xi_{B_{w},T_{v},B',T'}\circ (\tau^{-1}\sigma\tau)\circ \tau^{-1} .$$
Mais $\tau^{-1}\sigma\tau$ appartient \`a $\Gamma_{\bar{v}}$. En utilisant l'\'equivariance de $\xi_{B_{w},T_{v},B',T'}$, on obtient
$$\xi_{B_{w'},T_{v},B',T'}\circ \sigma=\tau\circ(\tau^{-1}\sigma\tau)\circ\xi_{B_{w},T_{v},B',T'}\circ\tau^{-1}=\sigma\circ\tau\circ\xi_{B_{w},T_{v},B',T'}\circ\tau^{-1}=\sigma\circ\xi_{B_{w'},T_{v},B',T'}$$
comme on le voulait.

On doit aussi prouver que l'application d\'eduite
$$\tilde{\xi}_{B_{w'},T_{v},B',T'}:(T_{v}/(1-\theta)(T_{v}))\times_{{\cal Z}(G)}{\cal Z}(\tilde{G})\to T'\times_{{\cal Z}(G')}{\cal Z}(\tilde{G}')$$
 envoie l'image de $\gamma_{v}$ dans l'espace de d\'epart sur l'image de $\delta_{v}$ dans l'espace d'arriv\'ee. La d\'efinition de l'action galoisienne sur ${\cal Z}(\tilde{G})$ permet d'\'etendre la relation (7) \`a ces applications, c'est-\`a-dire que l'on a la relation
$$\tilde{\xi}_{B_{w'},T_{v},B',T'} \circ \tau= \tau\circ \tilde{\xi}_{B_{w},T_{v},B',T'}. $$
Puisque $\tilde{\xi}_{B_{w},T_{v},B',T'}$ envoie l'image de $\gamma_{v}$ sur celle de $\delta_{v}$ et puisque les \'el\'ements $\gamma_{v}$ et $\delta_{v}$ sont tous deux invariants par $\tau$, la relation ci-dessus implique l'assertion cherch\'ee. Cela prouve (5). $\square$

  D\'efinissons un homomorphisme
$$\xi_{T_{v},T'}:T_{v}(E_{v})=\prod_{w'\vert v}T_{v}(E_{w'})\to T'(E_{v})=\prod_{w'\vert v}T'(E_{w'})$$
comme le produit sur les $w'$ divisant $v$ des homomorphismes $\xi_{B_{w'},T_{v},B',T'}$.  On a

(8) $\xi_{T_{v},T'}$ est \'equivariant pour les actions de $Gal(E/F)$.

 Cela  r\'esulte de ce  que, pour tout $w'$,  l'homomorphisme $\xi_{B_{w'},T_{v},B',T'}$ est \'equivariant pour l'action de $\Gamma_{\bar{v}'}$ et que l'on a la relation (7) ci-dessus.

Soit $v\in Val(F)$ telle que $v\not\in V_{ram}({\bf G}')$ et $v$ soit non-ramifi\'ee dans $E$. Rappelons que le sous-groupe compact hypersp\'ecial $K_{v}$ est attach\'e \`a un sch\'ema en groupes ${\cal K}_{v}$ d\'efini sur l'anneau des entiers $\mathfrak{o}_{v}$ de $F_{v}$. Si $w'$ est une place de $E$ au-dessus de $v$, posons $K_{w'}={\cal K}_{v}(\mathfrak{o}_{w'})$ et $\tilde{K}_{w'}=K_{w'}\tilde{K}_{v}$. Alors $K_{w'}$ est un sous-groupe compact hypersp\'ecial de $G(E_{w'})$ et $\tilde{K}_{w'}$ est un sous-espace hypersp\'ecial de $\tilde{G}(E_{w'})$.   On note aussi $K_{v}^{nr}={\cal K}_{v}(\mathfrak{o}_{v}^{nr})$, $F_{v}^{nr}$ \'etant identifi\'e \`a une extension de $E_{w'}$. De $K_{v}$ se d\'eduisent des sous-groupes hypersp\'eciaux $K_{v,sc}$ de $G_{SC}(F_{v})$ et $K_{v,ad}$ de $G_{AD}(F_{v})$ et on utilise pour ces groupes des notations similaires.   On peut fixer un ensemble fini $V'$ de places de $F$, contenant $V$ et les places ramifi\'ees dans $E$, de sorte que, pour $v\not\in V'$ et pour toute place $w'$ de $E$ au-dessus de $v$, on ait:

- la condition (2) est v\'erifi\'ee pour $v$;

(9) $K_{w'}$ est le sous-groupe compact hypersp\'ecial issu de la paire de Borel \'epingl\'ee ${\cal E}^*$;

(10) $e\in \tilde{K}_{w'}$, $x\in K_{w'}$ et, pour tout $\sigma\in \Gamma_{F}$, $u_{{\cal E}^*}(\sigma)\in K_{w'}$ et $z(\sigma)\in K_{w'}$.

Montrons que

(11) pour $v\not\in V'$ et pour toute place $w'$ de $E$ au-dessus de $v$, il existe $g_{w'}\in K_{w',sc}$ tel que $ad_{g_{w'}}(B_{w'},T_{v})=(B^*,T^*)$. 

Preuve. On ne perd rien ici \`a supposer  $w'=w$ (qui est la restriction de $\bar{v}$ \`a $E$).   Fixons un entier $N\geq1$ tel que $(\theta^*)^N=1$. L'hypoth\`ese $v\not\in V_{ram}$ assure que l'on peut prendre $N$ premier \`a la caract\'eristique r\'esiduelle $p$ de $F_{v}$. Introduisons le groupe non connexe $G^+=G\rtimes\{1,\theta^*,...,(\theta^*)^{N-1}\}$ d\'efini sur $E_{w}$. L'espace $\tilde{G}$ s'identifie \`a la composante $G\theta^*$: pour $g\in G$,  $ge$ s'identifie \`a $g\theta^*$.   Le sous-espace $\tilde{K}_{w}$ s'identifie \`a $K_{w}\theta^*$. Dans cette situation, on a d\'efini en [W] 5.2 la notion d'\'el\'ement compact de $\tilde{G}(E_{w})$: un \'el\'ement est compact si et seulement si le sous-groupe qu'il engendre dans $G^+(E_{w})$ est d'adh\'erence compacte. Puisque $\gamma_{v}\in \tilde{K}_{v}\subset K_{w}\theta^*$ et que $K_{w}\rtimes\{1,\theta^*,...,(\theta^*)^{N-1}\}$ est un groupe compact, l'\'el\'ement $\gamma_{v}$ est compact. D'apr\`es [W] 5.2, on peut d\'ecomposer $\gamma_{v}$ en $u\gamma_{v,p'}$, o\`u $\gamma_{v,p'}$ est d'ordre fini premier \`a $p$ et $u$ est topologiquement unipotent. Ces \'el\'ements appartiennent \`a l'adh\'erence du groupe engendr\'e par $\gamma_{v}$ et sont d\'efinis sur $F_{v}$. Comme on vient de le voir, l'intersection de cette adh\'erence avec $\tilde{G}$, resp. $G$, est contenue dans $\tilde{K}_{w}$, resp. $K_{w}$. Donc $\gamma_{v,p'}\in \tilde{K}_{v}$ et $u\in K_{v}$.    Les \'el\'ements $u$ et $\gamma_{v,p'}$  commutent entre eux et commutent donc \`a $\gamma_{v}$. Cela entra\^{\i}ne que $u\in Z_{G}(\gamma_{v};F_{v})=T_{v}^{\theta}(F_{v})$. Ecrivons $\gamma_{v}=t_{v}e_{v}$ comme dans la relation (2). Alors $\gamma_{v,p'}=u^{-1}t_{v}e_{v}$. Puisque $u$ est topologiquement unipotent, les valeurs de $(N\alpha)(u)$ sont des  \'el\'ements de $\mathfrak{o}_{w}$ de r\'eduction $1$ dans le corps r\'esiduel ${\mathbb E}_{w}$ de $\mathfrak{o}_{w}$. Alors la relation (2) est encore v\'erifi\'ee par $\gamma_{v,p'}$, a fortiori $\gamma_{v,p'}$ est r\'egulier.   Le lemme [W] 5.4 implique l'existence de $k\in K_{v}^{nr}$ tel que $ad_{k}(\gamma_{v,p'})\in T^*e$. Puisqu'il s'agit d'\'el\'ements r\'eguliers, on a automatiquement l'\'egalit\'e $ad_{k}(T_{v})=T^*$. Donc aussi $ad_{k}(\gamma_{v})\in T^*e$.  L'automorphisme $ad_{k}$ envoie $B_{w}$  sur un groupe de Borel contenant $T^*$ et conserv\'e par $\theta^*$. Un tel groupe est de la forme $\omega(B^*)$,  o\`u $\omega\in W^{\theta}$. Relevons $\omega$ en un \'el\'ement  $h\in K_{w}\cap G_{e}$. En rempla\c{c}ant $k$ par $h^{-1}k$, on obtient l'\'egalit\'e $ad_{k}(B_{w},T_{v})=(B^*,T^*)$.  Les th\'eor\`emes de structure de Bruhat-Tits entra\^{\i}nent que tout \'el\'ement de $K_{v}^{nr}$ est produit d'un \'el\'ement de $T^*(\mathfrak{o}_{v}^{nr})$ et d'un \'el\'ement de $K_{v,sc}^{nr}$. Quitte \`a multiplier $k$ \`a gauche par un \'el\'ement de $T^*(\mathfrak{o}_{v}^{nr})$, on peut supposer $k\in K_{v,sc}^{nr}$. La relation (4) entra\^{\i}ne que $ad_{k}:T_{v}\to T^*$ est \'equivariant pour l'action de $\Gamma_{E_{w}}$. Pour $\sigma\in Gal(F_{v}^{nr}/E_{w})$,  on a donc $k\sigma(k)^{-1}\in T^*\cap K_{v,sc}^{nr}=T^*(\mathfrak{o}_{v}^{nr})$. On obtient un cocycle de $Gal(F_{v}^{nr}/E_{w})$ dans $T^*(\mathfrak{o}_{v}^{nr})$. Un tel cocycle est un cobord.  Cela signifie que, quite \`a multiplier  $k$ \`a gauche par un \'el\'ement de $T^*(\mathfrak{o}_{v}^{nr})$, on peut supposer $k\sigma(k)^{-1}=1$ pour tout $\sigma\in Gal(F_{v}^{nr}/E_{w})$. Autrement dit $k\in K_{w,sc}$.   Cela prouve (11). $\square$

 Pour toute place $w'$ de $E$, on introduit le diagramme $(\delta_{v},B',T',B_{w'},T_{v},\gamma_{v})$ qui est unique d'apr\`es (3). On introduit un \'el\'ement $g_{w'}\in G_{SC}(E_{w'})$ tel que $ad_{g_{w'}}(B_{w'},T_{v})=(B^*,T^*)$. On suppose $g_{w'}\in K_{w',sc}$ pour toute place $w'$ au-dessus d'une place $v\not\in V'$. On pose $g=(g_{w'})_{w'\in Val(E)}$. C'est un \'el\'ement de $G_{SC}({\mathbb A}_{E})$. Pour $\sigma\in \Gamma_{F}$, posons
 $$V_{T}(\sigma)=x\sigma_{G^*}(x)^{-1}u_{{\cal E}^*}(\sigma) \sigma(g)g^{-1}.$$
 En utilisant la relation (8), on voit facilement que $V_{T}(\sigma)$ appartient \`a $T_{sc}({\mathbb A}_{E})$ (c'est-\`a-dire $T^*_{sc}({\mathbb A}_{E})$ muni de l'action galoisienne $\sigma\mapsto \sigma_{T}$). Le cobord de la cocha\^{\i}ne $V_{T}$ est \'egal \`a
   celui de $\sigma\mapsto u_{{\cal E}^*}(\sigma)$ donc prend ses valeurs dans $T_{sc}(E)$. En poussant $V_{T}$ en une cocha\^{\i}ne \`a valeurs dans $T_{sc}({\mathbb A}_{E})/T_{sc}(E)$, $V_{T}$ devient
     un cocycle.
On peut le consid\'erer comme un cocycle de $\Gamma_{F}$ dans $T_{sc}({\mathbb A}_{\bar{F}})/T_{sc}(\bar{F})$, o\`u ${\mathbb A}_{\bar{F}}$ est la limite inductive des ${\mathbb A}_{E'}$ sur les extensions finies $E'$ de $F$.
Comme en [I] 2.2, notons ${\cal T}_{1}$ le produit fibr\'e de $T'_{1}$ et $T$ au-dessus de $T'$. Il est muni de l'action galoisienne produit de l'action naturelle sur $T'_{1}$ et de l'action $\sigma\mapsto \sigma_{T}$ sur $T$. Notons $e'$ l'image naturelle de $e$ dans ${\cal Z}(\tilde{G}')\subset \tilde{G}'$ et fixons un rel\`evement $e'_{1}$ de $e'$ dans ${\cal Z}(\tilde{G}'_{1},E)$. Ecrivons $g\gamma g^{-1}=\nu e$, $\delta_{1}=\mu_{1}e'_{1}$. On a $\nu\in T({\mathbb A}_{E})$, $\mu_{1}\in T'_{1}({\mathbb A}_{E})$. On note $\nu_{1}$ l'image de $(\mu_{1},\nu)$ dans ${\cal T}_{1}({\mathbb A}_{\bar{F}})/{\cal T}_{1}(\bar{F})$. Pour deux tores $S_{1}$ et $S_{2}$ d\'efinis sur $F$ et pour un homomorphisme $f:S_{1}\to S_{2}$ d\'efini sur $F$, on note selon l'usage $H^{1,0}({\mathbb A}_{F}/F;S_{1}\stackrel{f}{\to}S_{2})$ le groupe $H^{1,0}(\Gamma_{F};S_{1}({\mathbb A}_{\bar{F}})/S_{1}(\bar{F})\stackrel{f}{\to} S_{2}({\mathbb A}_{\bar{F}})/S_{1}(\bar{F}))$, c'est-\`a-dire la limite inductive des $H^{1,0}(Gal(E'/F);S_{1}({\mathbb A}_{E'})/S_{1}(E')\stackrel{f}{\to}S_{2}({\mathbb A}_{E'})/S_{2}(E'))$ sur les extensions galoisiennes finies $E'$ de $F$. On note aussi $Z^{1,0}({\mathbb A}_{F}/F;S_{1}\stackrel{f}{\to}S_{2})$ l'ensemble de cocycles correspondant. On v\'erifie que le couple $(V_{T},\nu_{1})$ appartient \`a $Z^{1,0}({\mathbb A}_{F}/F; T_{sc}\stackrel{1-\theta}{\to}{\cal T}_{1})$. Sa classe dans $H^{1,0}({\mathbb A}_{F}/F; T_{sc}\stackrel{1-\theta}{\to}{\cal T}_{1})$ ne d\'epend pas du choix de l'\'el\'ement $g$: on ne peut changer $g$ qu'en le multipliant \`a gauche par un \'el\'ement de $T_{sc}({\mathbb A}_{E})$, ce qui multiplie $V_{T}$ par un cobord.

Dans le cas local, on a construit en [I] 2.2 (en suivant Kottwitz et Shelstad) un cocycle $\hat{V}_{{\cal T}_{1}}$ de $W_{F}$ dans le tore dual $\hat{\cal T}_{1}$. La m\^eme construction vaut dans le cas global, \`a condition bien s\^ur  d'utiliser des $\chi$-data globales (c'est-\`a-dire que les $\chi_{\alpha}$ sont des caract\`eres automorphes d'extensions $F_{\alpha}$ de $F$). Comme dans cette r\'ef\'erence, on \'ecrit $\tilde{s}=s\hat{\theta}$, avec $s\in \hat{T}$. On note $s_{ad}$ l'image de $s$ dans $\hat{T}_{ad}$. On v\'erifie que le couple $(\hat{V}_{{\cal T}_{1}},s_{ad})$ appartient \`a $Z^{1,0}(W_{F};\hat{\cal T}_{1}\stackrel{1-\hat{\theta}}{\to}\hat{T}_{ad})$. 

D'apr\`es [KS] (C.2.3), on dispose d'un accouplement
$$H^{1,0}({\mathbb A}_{F}/F; T_{sc}\stackrel{1-\theta}{\to}{\cal T}_{1})\times H^{1,0}(W_{F};\hat{\cal T}_{1}\stackrel{1-\hat{\theta}}{\to}\hat{T}_{ad})\to {\mathbb C}^{\times},$$
que l'on note $<.,.>$. On pose
$$\Delta_{imp}(\delta_{1},\gamma)=<(V_{T},\nu_{1});(\hat{V}_{{\cal T}_{1}},s_{ad})>^{-1}.$$
On a d\'ej\`a fix\'e des $\chi$-data globales pour $T$. On fixe aussi des $a$-data globales, c'est-\`a-dire que les $a_{\alpha}$ appartiennent \`a $\bar{F}^{\times}$. On peut alors d\'efinir un facteur $\Delta_{II}(\delta,\gamma)$ comme en [I] 2.2. Le point est ici qu'une expression comme $(N\alpha)(\nu)-1$ est un \'el\'ement du groupe d'id\`eles de l'extension $F_{\alpha}$ parce que (2) entra\^{\i}ne que c'est une unit\'e pour presque tout $v$. On pose
$$\Delta(\delta_{1},\gamma)=\Delta_{II}(\delta,\gamma)\Delta_{imp}(\delta_{1},\gamma).$$
On a utilis\'e de nombreuses donn\'ees auxiliaires mais on va montrer que $\Delta(\delta_{1},\gamma)$ ne d\'epend que de $\delta_{1}$ et $\gamma$. Plus g\'en\'eralement, reprenons la construction \`a partir d'un autre tore $\underline{T}'$ v\'erifiant aussi l'hypoth\`ese (Hyp). On souligne les donn\'ees relatives \`a ce nouveau tore. On introduit de nouveaux \'el\'ements $\underline{\delta}_{1}$, $\underline{\gamma}$. Pour tout $v\in Val(F)$, les couples $(\delta_{1,v},\gamma_{v})$ et $(\underline{\delta}_{1,v},\underline{\gamma}_{v})$ appartiennent \`a ${\cal D}_{1,v}$, le facteur $\boldsymbol{\Delta}_{1,v}(\delta_{1,v},\gamma_{v};\underline{\delta}_{1,v},\underline{\gamma}_{v})$ est donc bien d\'efini.

\ass{Proposition}{(i) Pour presque tout $v$, on a $\boldsymbol{\Delta}_{1,v}(\delta_{1,v},\gamma_{v};\underline{\delta}_{1,v},\underline{\gamma}_{v})=1$.

(ii) On a l'\'egalit\'e
$$\Delta(\delta_{1},\gamma)\Delta(\underline{\delta}_{1},\underline{\gamma})^{-1}=\prod_{v\in Val(F)}\boldsymbol{\Delta}_{1,v}(\delta_{1,v},\gamma_{v};\underline{\delta}_{1,v},\underline{\gamma}_{v}).$$

(iii) Le terme $\Delta(\delta_{1},\gamma)$ ne d\'epend pas des donn\'ees auxiliaires utilis\'ees pour le d\'efinir.}

Preuve.  Dans la construction de $\Delta(\delta_{1},\gamma)$, on a utilis\'e une paire de Borel \'epingl\'ee ${\cal E}^*$ d\'efinie sur $\bar{F}$, une cocha\^{\i}ne $u_{{\cal E}^*}$ \`a valeurs dans $G_{SC}(\bar{F})$, un \'el\'ement $e\in Z(\tilde{G},{\cal E}^*;\bar{F})$ et un \'el\'ement $e'_{1}\in {\cal Z}(\tilde{G}'_{1};\bar{F})$. Les choix des termes $u_{{\cal E}^*}$, $e$ et $e'_{1}$ n'influent pas sur $\Delta(\delta_{1},\gamma)$. En effet, le choix de $u_{{\cal E}^*}$ ne change pas l'action galoisienne $\sigma\mapsto \sigma_{G^*}$. On ne peut modifier $u_{{\cal E}^*}$ que par des \'el\'ements qui appartiennent  \`a $Z(G_{SC};\bar{F})$, donc \`a $T_{sc}(\bar{F})$, et de tels termes ne changent pas l'image de $V_{T}$ dans $T_{sc}({\mathbb A}_{\bar{F}})/T_{sc}(\bar{F})$. De m\^eme, on ne peut modifier le couple $(e,e'_{1})$  que par un  \'el\'ement du produit fibr\'e de  $Z(G;\bar{F})$ et $Z(G'_{1};\bar{F})$ au-dessus de $Z(G';\bar{F})$, ce   qui ne modifie pas l'image de $\nu_{1}$ modulo ${\cal T}_{1}(\bar{F})$.  Dans la construction de $\Delta(\underline{\delta}_{1},\underline{\gamma}_{1})$, on utilise d'autres donn\'ees $\underline{{\cal E}}^*$, $u_{\underline{{\cal E}}^*}$, $\underline{e}$, $\underline{e}'_{1}$. On peut fixer $r\in G_{SC}(\bar{F})$ qui conjugue ${\cal E}^*$ en $\underline{{\cal E}}^*$. D'apr\`es ce que l'on vient de dire, on peut supposer que $u_{\underline{{\cal E}}^*}(\sigma)=ru_{{\cal E}^*}(\sigma)\sigma(r)^{-1}$ pour tout $\sigma\in \Gamma_{F}$ et que $\underline{e}=ad_{r}(e)$. On peut alors supposer que $\underline{e}'_{1}=e'_{1}$, puisque $e$ et $\underline{e}$ ont alors la m\^eme image $e'$ dans ${\cal Z}(\tilde{G}';\bar{F})$. 

On peut aussi modifier la d\'efinition de $V_{T}$ (et aussi de $V_{\underline{T}}$) de la fa\c{c}on suivante. A l'aide des $a$-data globales que l'on a fix\'ees, on peut d\'efinir une cocha\^{\i}ne $r_{T}$ comme dans le cas local, cf. [I] 2.2. On peut alors  d\'efinir $V_{T}$ par
$$V_{T}(\sigma)=r_{T}(\sigma)n_{{\cal E}^*}(\omega_{T}(\sigma))u_{{\cal E}^*}(\sigma)\sigma(g)g^{-1}.$$
En effet,  on passe de la d\'efinition pr\'ec\'edente \`a celle-ci en multipliant \`a gauche par $r_{T}(\sigma)n_{{\cal E}^*}(\omega_{T}(\sigma)) \sigma_{G^*}(x)x^{-1}$. Par d\'efinition de $x$, c'est un \'el\'ement de $T_{sc}(\bar{F})$.  Il ne modifie donc pas l'image de $V_{T}$ modulo ce groupe.
 
  Ces modifications \'etant faites, on a l'\'egalit\'e
  $$\Delta_{imp}(\delta_{1},\gamma)\Delta_{imp}(\underline{\delta}_{1},\underline{\gamma})^{-1}=<((V_{T},V_{\underline{T}}^{-1}),(\nu_{1},\underline{\nu}_{1}^{-1})),((\hat{V}_{{\cal T}_{1}},\hat{V}_{\underline{{\cal T}}_{1}}),(s_{ad},s_{ad}))>^{-1},$$
  le produit \'etant celui sur
$$H^{1,0}({\mathbb A}_{F}/F; (T_{sc}\times \underline{T}_{sc})\stackrel{1-\theta}{\to}({\cal T}_{1}\times \underline{{\cal T}}_{1}))\times H^{1,0}(W_{F};(\hat{\cal T}_{1}\times \hat{\underline{{\cal T}}}_{1})\stackrel{1-\hat{\theta}}{\to}(\hat{T}_{ad}\times \hat{\underline{T}}_{ad})).$$
On introduit les tores $U$ et $S_{1}$ de [I] 2.2. Dans cette r\'ef\'erence, le tore $U$ est \'egal \`a $ad_{g^{-1}}(T)\times ad_{\underline{g}^{-1}}(\underline{T}))/diag_{-}(Z(G_{SC}))$, mais on peut l'identifier \`a $(T\times \underline{T})/diag_{-}(Z(G_{SC}))$.  De m\^eme pour $S_{1}$. Alors les deux tores sont d\'efinis sur $F$. Rappelons que $\hat{S}_{1}$ est le tore des $(t,\underline{t},t_{sc})\in \hat{\cal T}_{1}\times \hat{\underline{{\cal T}}}_{1}\times \hat{T}_{sc}$ tels que $j(t_{sc})=t\underline{t}^{-1}$, o\`u on a identifi\'e $\hat{T}$ et $\hat{\underline{T}}$ \`a un tore commun (muni de deux actions galoisiennes en g\'en\'eral distinctes) et o\`u $j:\hat{T}_{sc}\to \hat{T}$ est l'homomorphisme naturel. La structure galoisienne sur $\hat{S}_{1}$ est un peu compliqu\'ee, les formules sont les m\^emes que dans le cas local. On a aussi $\hat{U}=(\hat{T}_{sc}\times \hat{\underline{T}}_{sc})/diag(Z(\hat{G}_{SC}))$. On a un diagramme commutatif \'evident
$$\begin{array}{ccc}\hat{S}_{1}&\stackrel{1-\hat{\theta}}{\to}& \hat{U}\\ \downarrow&&\downarrow\\ 
 \hat{\cal T}_{1}\times \hat{\underline{{\cal T}}}_{1}&\stackrel{1-\hat{\theta}}{\to}&\hat{T}_{ad}\times \hat{\underline{T}}_{ad}\\ \end{array}$$
 d'o\`u un homomorphisme
 $$H^{1,0}(W_{F};\hat{S}_{1}\stackrel{1-\hat{\theta}}{\to} \hat{U})\to H^{1,0}(W_{F};(\hat{\cal T}_{1}\times \hat{\underline{{\cal T}}}_{1})\stackrel{1-\hat{\theta}}{\to}(\hat{T}_{ad}\times \hat{\underline{T}}_{ad})).$$
 En copiant les d\'efinitions du cas local, on d\'efinit un \'el\'ement $(\hat{V}_{1},{\bf s})\in H^{1,0}(W_{F};\hat{S}_{1}\stackrel{1-\hat{\theta}}{\to} \hat{U})$. Il r\'esulte des d\'efinitions que son  image dans $H^{1,0}(W_{F};(\hat{\cal T}_{1}\times \hat{\underline{{\cal T}}}_{1})\stackrel{1-\hat{\theta}}{\to}(\hat{T}_{ad}\times \hat{\underline{T}}_{ad}))$ est \'egale \`a
 $$((\hat{V}_{{\cal T}_{1}},\hat{V}_{\underline{{\cal T}}_{1}}),(s_{ad},s_{ad})).$$
 Notons $(V,\boldsymbol{\nu}_{1})$ l'image de 
$$((V_{T},V_{\underline{T}}^{-1}),(\nu_{1},\underline{\nu}_{1}^{-1}))$$
dans $H^{1,0}({\mathbb A}_{F}/F;U\stackrel{1-\theta}{\to}S_{1})$ par l'homomorphisme dual du pr\'ec\'edent. Par compatibilit\'e des produits, on a
  $$\Delta_{imp}(\delta_{1},\gamma)\Delta_{imp}(\underline{\delta}_{1},\underline{\gamma})^{-1}=<(V,\boldsymbol{\nu}_{1}),(\hat{V}_{1},{\bf s})>^{-1}.$$
  Le terme $(V,\boldsymbol{\nu}_{1})$ se rel\`eve en un \'el\'ement de $H^{1,0}(\Gamma_{F};U({\mathbb A}_{\bar{F}})\stackrel{1-\theta}{\to}S_{1}({\mathbb A}_{\bar{F}}))$, d\'efini exactement par les m\^emes formules (apr\`es les modifications apport\'ees ci-dessus). Cet \'el\'ement appartient en fait \`a $H^{1,0}(Gal(E/F); U({\mathbb A}_{E})\stackrel{1-\theta}{\to}S_{1}({\mathbb A}_{E}))$, si $E$ d\'esigne maintenant une extension finie de $F$ v\'erifiant les m\^emes propri\'et\'es que plus haut mais pour nos deux ensembles de donn\'ees.  Pour toute place $v\in Val(F)$, cet \'el\'ement relev\'e d\'efinit un \'el\'ement $(V_{v},\boldsymbol{\nu}_{1,v})\in H^{1,0}(\Gamma_{F_{v}};U\stackrel{1-\theta}{\to}S_{1})$. De m\^eme, $(\hat{V}_{1},{\bf s})$ se restreint en un \'el\'ement $(\hat{V}_{1,v},{\bf s}_{v})\in H^{1,0}(W_{F_{v}};\hat{S}_{1}\stackrel{1-\hat{\theta}}{\to}\hat{U})$. La compatibilit\'e des produits et le lemme C.1.B de [KS] assure que
  
  (12) on a l'\'egalit\'e $<  (V_{v},\boldsymbol{\nu}_{1,v}),(\hat{V}_{1,v},{\bf s}_{v})>=1$ pour presque tout $v$;
  
  $$(13) \qquad <(V,\boldsymbol{\nu}_{1}),(\hat{V}_{1},{\bf s})>=\prod_{v\in Val(F)}<  (V_{v},\boldsymbol{\nu}_{1,v}),(\hat{V}_{1,v},{\bf s}_{v})>.$$
  
  On a aussi
  $$\Delta_{II}(\delta,\gamma)=\prod_{v\in Val(F)}\Delta_{II}(\delta_{v},\gamma_{v})$$
 et  les termes du produit sont  presque tous \'egaux \`a $1$. On en d\'eduit
  $$\Delta_{imp}(\delta_{1},\gamma)\Delta_{imp}(\underline{\delta}_{1},\underline{\gamma})^{-1}=\prod_{v\in Val_{F}}\Delta_{II}(\delta_{v},\gamma_{v})\Delta_{II}(\underline{\delta}_{v},\underline{\gamma}_{v})^{-1}<  (V_{v},\boldsymbol{\nu}_{1,v}),(\hat{V}_{1,v},{\bf s}_{v})>^{-1}.$$
  Pour achever la preuve des deux premi\`eres assertions de l'\'enonc\'e, il suffit de prouver que, pour tout $v$, on a l'\'egalit\'e
  $$\Delta_{II}(\delta_{v},\gamma_{v})\Delta_{II}(\underline{\delta}_{v},\underline{\gamma}_{v})^{-1}<  (V_{v},\boldsymbol{\nu}_{1,v}),(\hat{V}_{1,v},{\bf s}_{v})>^{-1}=\boldsymbol{\Delta}_{1,v}(\delta_{1,v},\gamma_{v};\underline{\delta}_{1,v},\underline{\gamma}_{v}).$$
    Pour d\'efinir le membre de droite, on utilise les paires de Borel \'epingl\'ees ${\cal E}_{w}=ad_{g_{w}^{-1}}({\cal E}^*)$ et $\underline{{\cal E}}_{w}=ad_{\underline{g}_{w}^{-1}}(\underline{{\cal E}}^*)$, o\`u $w$ est encore la place de $E$ fix\'ee au-dessus de $v$. On choisit $u_{{\cal E}_{w}}(\sigma)=g_{w}^{-1}u_{{\cal E}^*}(\sigma)\sigma(g_{w})$ et $u_{\underline{{\cal E}}_{w}}(\sigma)=\underline{g}_{w}^{-1}u_{\underline{{\cal E}}^*}(\sigma)\sigma(\underline{g}_{w})$ pour tout $\sigma\in \Gamma_{F_{v}}$. On constate alors que les  deux membres de l'\'egalit\'e ci-dessus sont d\'efinis de la m\^eme fa\c{c}on, apr\`es l'identification que l'on a faite des tores $ad_{g^{-1}}(T)$ et $ad_{\underline{g}^{-1}}(\underline{T})$ \`a $T$ et $\underline{T}$.   Cela prouve les deux premi\`eres assertions de l'\'enonc\'e.
  
  L'assertion (iii) s'en d\'eduit: puisque les termes $\Delta_{1}(\underline{\delta}_{1},\underline{\gamma})$ et $\boldsymbol{\Delta}_{1,v}(\delta_{1,v},\gamma_{v};\underline{\delta}_{1,v},\underline{\gamma}_{v})$ ne d\'ependent pas des donn\'ees auxiliaires utilis\'ees pour d\'efinir $\Delta_{1}(\delta_{1},\gamma)$, ce dernier terme n'en d\'epend pas non plus. $\square$
  
  Revenons \`a notre tore $T'$. Il reste \`a v\'erifier que l'on peut choisir des \'el\'ements $\delta$ et $\gamma$ v\'erifiant les conditions (1) et  (2). On introduit les m\^emes donn\'ees qu'au d\'ebut du paragraphe, en particulier le corps $E$. On modifie la d\'efinition de l'ensemble $V'$. On note maintenant $V'$ un ensemble fini de places de $F$, contenant $V$ et les places ramifi\'ees dans $E$, de sorte que, pour tout $v\not\in V'$ et pour toute place $w'$ de $E$ au-dessus de $v$, les conditions (9) et (10) soient v\'erifi\'ees, ainsi que les conditions (14) et (15) ci-dessous.    
   Pour $v\not\in V$ telle que $v$ soit non ramifi\'ee dans $E$, les tores $T$  et $T'$ sont non ramifi\'es en $v$ et ont  donc une structure naturelle sur $\mathfrak{o}_{v}$. On note ${\mathbb T}_{v}=T\times_{\mathfrak{o}_{v}}{\mathbb F}_{v}$ la fibre de $T$  sur ${\mathbb F}_{v}$. Un \'el\'ement  de $\Sigma(T)$ est aussi un caract\`ere de ce tore. On impose
   
   (14) l'image naturelle $e'\in {\cal Z}(\tilde{G}')\subset \tilde{G}'$ de $e$ appartient \`a $\tilde{K}'_{w'}$ et $T'(\mathfrak{o}_{w'})$ est inclus dans $K'_{w'}$;
  
  (15) soit $t_{0}\in T(\mathfrak{o}_{w'})$; alors il existe $t\in T(\mathfrak{o}_{v})t_{0}$ dont la r\'eduction $\underline{t}\in {\mathbb T}_{v}({\mathbb E}_{w'})$ v\'erifie $N\alpha(\underline{t})\not=\pm 1$ pour tout $\alpha\in \Sigma(T)$.
  
  La condition (14) est satisfaite presque partout.  Il faut montrer qu'il en est de m\^eme de (15). Il s'agit de montrer que, pour presque tout $v$, ${\mathbb T}_{v}({\mathbb F}_{v})$ n'est pas contenu dans la r\'eunion sur $\alpha\in \Sigma(T)$ et $\epsilon=\pm 1$ des sous-ensembles $\{\underline{t}\in {\mathbb T}_{v}({\mathbb F}_{v}); N\alpha(\underline{t}\underline{t}_{0})=\epsilon\}$, o\`u $\underline{t}_{0}$ est la r\'eduction de $t_{0}$. Notons $d$ la dimension de $T$ et $q_{v}$ le nombre d'\'el\'ements de ${\mathbb F}_{v}$. Il existe $c>0$ ind\'ependant de $v$ tel que le nombre d'\'el\'ements de ${\mathbb T}_{v}({\mathbb F}_{v})$ soit au moins \'egal \`a $cq_{v}^d$. Il suffit de d\'emontrer qu'il existe $c'>0$ ind\'ependant de $v$ tel que chacun des sous-ensembles ci-dessus ait un nombre d'\'el\'ements au plus \'egal \`a $c'q_{v}^{d-1}$. Consid\'erons le sous-ensemble $\{\underline{t}\in {\mathbb T}_{v}({\mathbb F}_{v}); N\alpha(\underline{t}\underline{t}_{0})=\epsilon\}$. Il peut \^etre vide.
Sinon, il est en bijection avec $\{\underline{t}\in {\mathbb T}_{v}({\mathbb F}_{v}); N\alpha(\underline{t})=1\}$, ou encore avec
$$\{\underline{t}\in {\mathbb T}_{v}({\mathbb F}_{v}); \forall \sigma\in Gal(E_{w'}/F_{v}), (N\sigma\alpha)(\underline{t})=1\}.$$
 Introduisons le tore $S$ sur ${\mathbb F}_{v}$ qui est la restriction des scalaires du tore multiplicatif ${\mathbb G}_{m}$ sur ${\mathbb E}_{w'}$. L'homomorphisme
$$\begin{array}{cccc}\bar{\tau}:&{\mathbb T}_{v}(\bar{{\mathbb F}}_{v})&\to& S(\bar{{\mathbb F}}_{v})=\prod_{\sigma\in Gal(E_{w'}/F_{v})}\bar{{\mathbb F}}_{v}^{\times}\\& \underline{t}&\mapsto&((N\sigma\alpha)(\underline{t}))_{\sigma\in Gal(E_{w'}/F_{v})}\\ \end{array}$$
est d\'efini sur ${\mathbb F}_{v}$. L'ensemble pr\'ec\'edent est le noyau de l'homomorphisme
$$\tau:{\mathbb T}_{v}({\mathbb F}_{v})\to S({\mathbb F}_{v}).$$
On montre ais\'ement que le nombre de composantes connexes du noyau de $\bar{\tau}$
  est born\'e par un nombre qui ne d\'epend que de l'homomorphisme $X_{*}({\mathbb T}_{v})\to X_{*}(S)$ d\'eduit de $\bar{\tau}$ et de la structure de ces $Gal(E_{w'}/F_{v})$-modules. De m\^eme, la composante neutre de ce noyau est un tore d\'efini sur ${\mathbb F}_{v}$ et de dimension au plus $d-1$, dont la structure ne d\'epend que des m\^emes donn\'ees. Or ces donn\'ees ne varient que dans un ensemble fini, car le groupe $Gal(E_{w'}/F_{v})$ est lui-m\^eme toujours un sous-groupe de $Gal(E/F)$.  Il en r\'esulte que le nombre d'\'el\'ements du noyau est bien born\'e par $c'q_{v}^{d-1}$ pour un $c'>0$ ind\'ependant de $v$. Cela prouve que (15) est v\'erifi\'e pour presque tout $v$. 

Soit  $v\not\in V'$, notons $w$ la restriction de $\bar{v}$ \`a  $E$. Notons $g_{ad}$ l'image dans $G_{AD}$ d'un \'el\'ement $g\in G$. D'apr\`es (10), l'application $\sigma\mapsto x_{ad}\sigma_{G^*}(x_{ad})^{-1}u_{{\cal E}^*,ad}(\sigma)=x_{ad}u_{{\cal E}^*,ad}(\sigma)\sigma(x_{ad})^{-1}$ est un cocycle de $Gal(E_{w}/F_{v})$ dans $K_{w,ad}$.   Par le th\'eor\`eme de Lang, un tel cocycle est trivial (cf. [W2] lemme 4.2(ii)). On peut fixer $y_{ad}\in K_{w,ad}$ tel que $x_{ad}\sigma_{G^*}(x_{ad})^{-1}u_{{\cal E}^*,ad}(\sigma)=y_{ad}\sigma(y_{ad})^{-1}$ pour tout $\sigma\in Gal(E_{w}/F_{v})$. Puisque $K_{w}$ est le groupe hypersp\'ecial associ\'e \`a ${\cal E}^*$ d'apr\`es (9),  on a  $T^*_{ad}(\mathfrak{o}_{w})=T_{ad}(\mathfrak{o}_{w})\subset K_{w,ad}$ et  la d\'ecomposition d'Iwasawa montre que l'application 
 $$\begin{array}{ccc}T_{ad}(\mathfrak{o}_{w})\times K_{w,sc}&\to& K_{w,ad}\\ (t,x)&\mapsto &t\pi_{ad}(x)\\ \end{array}$$
 est surjective. On rel\`eve $y_{ad}$ en $(t,g_{w})\in T_{ad}(\mathfrak{o}_{w})\times K_{w,sc}$. On a alors
 $$x_{ad}\sigma_{G^*}(x_{ad})^{-1}u_{{\cal E}^*,ad}(\sigma)=t\pi_{ad}(g_{w}\sigma(g_{w})^{-1})\sigma(t)^{-1},$$
 d'o\`u
 $$t^{-1}x_{ad}\sigma_{G^*}(x_{ad})^{-1}u_{{\cal E}^*,ad}(\sigma)\sigma(t)=\pi_{ad}(g_{w}\sigma(g_{w})^{-1}).$$
 On a l'\'egalit\'e
 $$x_{ad}\sigma_{G^*}(x_{ad})^{-1}u_{{\cal E}^*,ad}(\sigma)\sigma(t)=\sigma_{T}(t)x_{ad}\sigma_{G^*}(x_{ad})^{-1}u_{{\cal E}^*,ad}(\sigma),$$
 d'o\`u
 $$t^{-1}\sigma_{T}(t)x_{ad}\sigma_{G^*}(x_{ad})^{-1}u_{{\cal E}^*,ad}(\sigma)\pi_{ad}(\sigma(g_{w})g_{w}^{-1})=1.$$
 A fortiori,
 $$x_{ad}\sigma_{G^*}(x_{ad})^{-1}u_{{\cal E}^*,ad}(\sigma)\pi_{ad}(\sigma(g_{w})g_{w}^{-1})\in T_{ad}(\mathfrak{o}_{w}),$$
 d'o\`u il r\'esulte que $x\sigma_{G^*}(x)^{-1}u_{{\cal E}^*}(\sigma)\sigma(g_{w})g_{w}^{-1}\in T_{sc}(\mathfrak{o}_{w})$. Notons $t_{sc}(\sigma)$ cet \'el\'ement. Posons ${\cal E}_{w}=ad_{g_{w}^{-1}}({\cal E}^*)$ et notons $(B_{w},T_{v})$ la paire de Borel sous-jacente \`a ${\cal E}_{w}$. Par le m\^eme calcul que dans la preuve de (6), la relation pr\'ec\'edente entra\^{\i}ne que $T_{v}$ est d\'efini sur $F_{v}$ et que l'homomorphisme $\xi_{B_{w},T_{v},B',T'}$ d\'eduit des paires $(B_{w},T_{v})$ et $(B',T')$ est \'equivariant pour les actions de $\Gamma_{\bar{v}}$. Posons $e_{v}=ad_{g_{w}^{-1}}(e)$. C'est un \'el\'ement de $Z(\tilde{G},{\cal E}_{w};E_{w})$. D'apr\`es (10), c'est aussi un \'el\'ement de $\tilde{K}_{w}$. Pour $\sigma\in \Gamma_{F_{v}}$, on a les \'egalit\'es
 $$\sigma(e_{v})=ad_{\sigma(g_{w}^{-1})}\circ\sigma(e)=ad_{\sigma(g_{w})^{-1}u_{{\cal E}^*}(\sigma)^{-1}\sigma_{G^*}(x)x^{-1}}\circ ad_{x\sigma_{G^*}(x)^{-1}u_{{\cal E}^*}(\sigma)}\circ \sigma(e).$$
 On a $ad_{u_{{\cal E}^*}(\sigma)}\circ\sigma(e)=z(\sigma)^{-1}e$. Puisque $x\in G_{SC}^{\theta^*}(E_{w})$, l'\'el\'ement $x\sigma_{G^*}(x)^{-1}$ commute \`a $e$ et aussi, bien s\^ur, \`a $z(\sigma)$. Donc 
 $$ad_{x\sigma_{G^*}(x)^{-1}u_{{\cal E}^*}(\sigma)}\circ \sigma(e)=z(\sigma)^{-1}e.$$
 De plus 
 $$\sigma(g_{w})^{-1}u_{{\cal E}^*}(\sigma)^{-1}\sigma_{G^*}(x)x^{-1}=g_{w}^{-1}t_{sc}(\sigma)^{-1}=t_{sc,w}(\sigma)^{-1}g_{w}^{-1},$$
 o\`u $t_{sc,w}=ad_{g_{w}^{-1}}(t_{sc}(\sigma))\in T_{v}(\mathfrak{o}_{w})$. D'o\`u
 $$\sigma(e_{v})=z(\sigma)^{-1} ad_{t_{sc,w}^{-1}g_{w}^{-1}}(e)=z(\sigma)^{-1}t_{sc,w}^{-1}e_{w}t_{sc,w}=z(\sigma)^{-1}(\theta-1)(t_{sc,w})e_{v}.$$
 L'application $\sigma\mapsto z(\sigma)^{-1}(\theta-1)(t_{sc,w}(\sigma))$  prend ses valeurs dans $T_{v}(\mathfrak{o}_{w})$ et la relation ci-dessus implique que c'est un cocycle de $Gal(E_{w}/F_{v})$ \`a valeurs dans ce groupe.  Un tel cocycle est forc\'ement trivial. On peut donc fixer $t_{0}\in T_{v}(\mathfrak{o}_{w})$ tel que $\sigma(t_{0})=t_{0} z(\sigma)(1-\theta)(t_{sc}(\sigma))$ pour tout $\sigma$. Alors $t_{0}e_{v}\in \tilde{G}(F_{v})$. 
   On peut multiplier $t_{0}$ par un \'el\'ement de $T_{v}(\mathfrak{o}_{v})$ de sorte que le produit $t$ v\'erifie la conclusion de (15) transport\'ee \`a $T_{v}$ par l'isomorphisme $ad_{g_{w}^{-1}}$. En multipliant encore par un \'el\'ement de $T_{v}(\mathfrak{o}_{v})$ assez voisin de l'origine, on peut assurer que $te_{v}$ est fortement r\'egulier. Posons $\gamma_{v}=te_{v}$ et $\delta_{v}=\xi_{B_{w},T_{v},B',T'}(t)e'$, o\`u $e'\in {\cal Z}(\tilde{G}')$ est l'image naturelle de $e$ (ou $e_{v}$, c'est pareil). Les constructions impliquent que $\gamma_{v}\in \tilde{K}_{v}$, $\delta_{v}\in \tilde{K}'_{v}$ et $(\delta_{v},B',T',B_{w},T_{v},\gamma_{v})$ est un diagramme.  Le choix de $t$ implique que la condition (2) est satisfaite.  Pour $v\not\in V'$, on a donc construit des \'el\'ements v\'erifiant (1) et (2).

   \bigskip
  
  \subsection{Utilisation du facteur de transfert global, cas particulier}

Soit ${\bf G}'=(G',{\cal G}',\tilde{s})$ une donn\'ee endoscopique relevante de $(G,\tilde{G},{\bf a})$. On suppose qu'elle v\'erifie l'hypoth\`ese (Hyp) du paragraphe pr\'ec\'edent. Soit $V$ un ensemble fini de places de $F$ contenant $V_{ram}({\bf G}')$. On a d\'efini l'espace $C_{c}^{\infty}({\bf G}'_{V})$ en 3.3. Pour $v\in V$, on a d\'efini en [I] 2.5 un espace $C_{c}^{\infty}({\bf G}'_{v})$ par une tout autre m\'ethode. 

\ass{Proposition}{Il existe un isomorphisme canonique
$$C_{c}^{\infty}({\bf G}'_{V})\simeq \otimes_{v\in V}C_{c}^{\infty}({\bf G}'_{v}).$$}

Preuve. Consid\'erons des donn\'ees auxiliaires $G'_{1}$, $\tilde{G}'_{1}$, $C_{1}$, $\hat{\xi}_{1}$, $(\tilde{K}'_{1,v})_{v\not\in V}$  non ramifi\'ees hors de $V$. Pour $v\in Val(F)-V$, le choix des espaces hypersp\'eciaux d\'etermine un facteur de transfert $\Delta_{1,v}$, cf. [I] 6.3. Pour $v\in V$, on ne sait pas normaliser le facteur de transfert. Mais on peut normaliser le produit sur $v\in V$ de ces facteurs. En effet,  construisons des \'el\'ements comme dans le paragraphe pr\'ec\'edent, et il est plus simple ici de les souligner. On a donc des \'el\'ements $\underline{\delta}_{1}=(\underline{\delta}_{1,v})_{v\in Val(F)}\in \tilde{G}_{1}'({\mathbb A}_{F})$, $\underline{\gamma}=(\underline{\gamma}_{v})_{v\in Val(F)}\in \tilde{G}({\mathbb A}_{F})$ et le facteur $\Delta_{1}(\underline{\delta}_{1},\underline{\gamma})$.   Soient $\delta_{1,V}=(\delta_{1,v})_{v\in V}\in \tilde{G}'_{1}(F_{V})$, $\gamma_{V}=(\gamma_{v})_{v\in V}\in \tilde{G}(F_{V})$. Supposons $(\delta_{1,V},\gamma_{V})\in {\cal D}_{1,V}$, on entend par l\`a que $(\delta_{1,v} ,\gamma_{v})\in {\cal D}_{1,v}$ pour tout $v\in V$. On pose
$$\Delta_{1,V}(\delta_{1,V},\gamma_{V})=\Delta_{1}(\underline{\delta}_{1},\underline{\gamma})\left(\prod_{v\not\in V}\Delta_{1,v}(\underline{\delta}_{1,v},\underline{\gamma}_{v})^{-1}\right)\left(\prod_{v\in V}\boldsymbol{\Delta}_{1,v}(\delta_{1,v},\gamma_{v};\underline{\delta}_{1,v},\underline{\gamma}_{v})\right) .$$
Il r\'esulte des calculs du paragraphe pr\'ec\'edent que les termes du premier produit sont presque tous \'egaux \`a $1$. Le terme ainsi d\'efini est un facteur de transfert. La proposition du paragraphe pr\'ec\'edent montre qu'il ne d\'epend pas des donn\'ees auxiliaires $\underline{\delta}_{1}$ et $\underline{\gamma}$.

D'apr\`es 3.3, le choix des $(\tilde{K}'_{1,v})_{v\not\in V}$ permet d'identifier $C_{c}^{\infty}({\bf G}'_{V})$ \`a $C_{c,\lambda_{1}}^{\infty}(\tilde{G}'_{1}(F_{V}))=\otimes_{v\in V}C_{c,\lambda_{1}}^{\infty}(\tilde{G}'_{1}(F_{v}))$. D'apr\`es [I] 2.5, le choix de $\Delta_{1,V}$ permet d'identifier ce dernier espace \`a  $\otimes_{v\in V}C_{c}^{\infty}({\bf G}'_{v})$. D'o\`u l'isomorphisme de l'\'enonc\'e. Pour qu'il soit "canonique", il suffit qu'il ne d\'epende pas des donn\'ees auxiliaires. Consid\'erons une autre famille $G'_{2}$, $\tilde{G}'_{2}$, $C_{2}$, $\hat{\xi}_{2}$, $(\tilde{K}'_{2,v})_{v\not\in V}$ de donn\'ees auxiliaires  non ramifi\'ees hors de $V$. Il y a deux isomorphismes de recollement entre les espaces $C_{c,\lambda_{1}}^{\infty}(\tilde{G}'_{1}(F_{V}))$ et $C_{c,\lambda_{2}}^{\infty}(\tilde{G}'_{2}(F_{V}))$: celle de 3.3 utilisant les espaces hypersp\'eciaux $(\tilde{K}'_{1,v})_{v\not\in V}$ et $(\tilde{K}'_{2,v})_{v\not\in V}$; celle de [I] 2.5 utilisant les facteurs de transfert $\Delta_{1,V}$ et $\Delta_{2,V}$. On doit prouver que ce sont les m\^emes. Les deux isomorphismes $f_{1}\mapsto f_{2}$ sont d\'efinis pas une formule $f_{2}(\delta_{2,V})=\tilde{\lambda}_{12,V}(\delta_{1,V},\delta_{2,V})f_{1}(\delta_{1,V})$ o\`u $\delta_{1,V}$ est un \'el\'ement quelconque tel que $(\delta_{1,V},\delta_{2,V})\in \tilde{G}'_{12}(F_{V})$, mais la fonction $\tilde{\lambda}_{12,V}$ n'est pas a priori la m\^eme pour les deux recollements. Notons $\tilde{\lambda}_{12,K,V}$ celle pour le premier recollement et $\tilde{\lambda}_{12,\Delta,V}$ celle pour le second. Soient $(\delta_{1,V},\delta_{2,V})\in \tilde{G}'_{12}(F_{V})$ et $(x_{1},x_{2})\in G'_{12}(F_{V})$. On a en tout cas
$$\tilde{\lambda}_{12,K,V}(x_{1}\delta_{1,V},x_{2}\delta_{2,V})=\lambda_{12,V}(x_{1},x_{2})\tilde{\lambda}_{12,K,V}(\delta_{1,V},\delta_{2,V}),$$ 
$$\tilde{\lambda}_{12,\Delta,V}(x_{1}\delta_{1,V},x_{2}\delta_{2,V})=\lambda_{12,V}(x_{1},x_{2})\tilde{\lambda}_{12,\Delta,V}(\delta_{1,V},\delta_{2,V}),$$
pour un m\^eme caract\`ere $\lambda_{12,V}$ de $G'_{12}(F_{V})$. Il suffit donc de prouver l'\'egalit\'e   $\tilde{\lambda}_{12,K,V}(\delta_{1,V},\delta_{2,V})=\tilde{\lambda}_{12,\Delta,V}(\delta_{1,V},\delta_{2,V})$ pour un seul couple $(\delta_{1,V},\delta_{2,V})$. On choisit ce couple ainsi: on construit des \'el\'ements $\delta\in \tilde{G}'({\mathbb A}_{F})$ et $\gamma\in \tilde{G}({\mathbb A}_{F})$ comme en 3.6, on rel\`eve $\delta$ en $\delta_{1}\in \tilde{G}'_{1}({\mathbb A}_{F})$ et $\delta_{2}\in \tilde{G}'_{2}({\mathbb A}_{F})$; on prend pour $\delta_{1,V}$ et $\delta_{2,V}$ les produits sur $v\in V$ des composantes locales de $\delta_{1}$ et $\delta_{2}$. Par d\'efinition de $\tilde{\lambda}_{12,\Delta,V}$, on a l'\'egalit\'e
$$(1) \qquad \Delta_{2,V}(\delta_{2,V},\gamma_{V})=\tilde{\lambda}_{12,\Delta,V}(\delta_{1,V},\delta_{2,V})\Delta_{1,V}(\delta_{1,V},\gamma_{V}).$$
En 1.15, on a normalis\'e une fonction $\tilde{\lambda}_{12}$ sur $\tilde{G}'_{12}({\mathbb A})$ de sorte qu'elle vaille $1$ sur $\tilde{G}'_{12}(F)$ et des fonctions $\tilde{\lambda}_{12,v}$ pour $v\not\in V$. Par d\'efinition,
$$\tilde{\lambda}_{12,K,V}(\delta_{1,V},\delta_{2,V})=\tilde{\lambda}_{12}(\delta_{1},\delta_{2})\prod_{v\not\in V}\tilde{\lambda}_{12,v}(\delta_{1,v},\delta_{2,v})^{-1}.$$
On a

(2) $\Delta_{2,v}(\delta_{2,v},\gamma_{v})=\tilde{\lambda}_{12,v}(\delta_{1,v},\delta_{2,v})\Delta_{1,v}(\delta_{2,v},\gamma_{v})$ pour tout $v\not\in V$.

En effet,  avec les notations de [I] 6.3, on a l'\'egalit\'e 
$$\Delta_{2,v}(\delta_{2,v},\gamma_{v})=\tilde{\lambda}_{\zeta_{1}}(\delta_{1,v})\tilde{\lambda}_{\zeta_{2}}(\delta_{2,v})^{-1}\Delta_{1,v}(\delta_{1,v},\gamma_{v}).$$
 Il suffit de comparer les d\'efinitions pour constater que
 $$\tilde{\lambda}_{\zeta_{1}}(\delta_{1,v})\tilde{\lambda}_{\zeta_{2}}(\delta_{2,v})^{-1}= \tilde{\lambda}_{12,v}(\delta_{1,v},\delta_{2,v}).$$
 D'o\`u (2).
 
 Alors
$$\tilde{\lambda}_{12,K,V}(\delta_{1,V},\delta_{2,V})=\tilde{\lambda}_{12}(\delta_{1},\delta_{2})\prod_{v\not\in V}\Delta_{2,v}(\delta_{2,v},\gamma_{v})^{-1}\Delta_{1,v}(\delta_{1,v},\gamma_{v})$$
$$=\tilde{\lambda}_{12}(\delta_{1},\delta_{2})\Delta_{2}(\delta_{2},\gamma)^{-1}\Delta_{1}(\delta_{1},\gamma)\Delta_{2,V}(\delta_{2,V},\gamma_{V})\Delta_{1,V}(\delta_{1,V},\gamma_{V})^{-1}.$$
En comparant avec (1), il reste \`a montrer l'\'egalit\'e
$$(3) \qquad \Delta_{2}(\delta_{2},\gamma)=\tilde{\lambda}_{12}(\delta_{1},\delta_{2})\Delta_{1}(\delta_{1},\gamma).$$
La d\'emonstration est similaire \`a celle du lemme [I] 2.5, nous n'en donnons que le squelette. De fa\c{c}on g\'en\'erale, pour un groupe r\'eductif connexe $H$ d\'efini sur $F$, un \'el\'ement de $H^1(W_{F},Z(\hat{H}))$ d\'etermine non seulement un caract\`ere de $H({\mathbb A}_{F})$ trivial sur $H(F)$, mais plus g\'en\'eralement un caract\`ere du groupe $( H({\mathbb A}_{\bar{F}})/Z(H;\bar{F}))^{\Gamma_{F}}$, trivial sur $(  H(\bar{F})/Z(H;\bar{F}))^{\Gamma_{F}}=H_{AD}(F)$. Avec les notations de [I] 2.5, le cocycle $w\mapsto (\zeta_{1}(w),\zeta_{2}(w)^{-1})$ de $W_{F}$ dans $Z(\hat{G}'_{12})$ d\'etermine donc un caract\`ere de $(  G'_{12}({\mathbb A}_{\bar{F}})/Z(G'_{12};\bar{F}))^{\Gamma_{F}}$, trivial sur $G'_{12,AD}(F)$. Notons ce caract\`ere $\underline{\tilde{\lambda}}_{12}$. L'ensemble $\tilde{G}'_{12}({\mathbb A}_{F})$ s'envoie naturellement dans  $(G'_{12}({\mathbb A}_{\bar{F}})/Z(G'_{12};\bar{F}))^{\Gamma_{F}}$. En effet, pour $(\underline{\delta}_{1},\underline{\delta}_{2})\in \tilde{G}'_{12}({\mathbb A}_{F})$, on choisit $(e'_{1},e'_{2})\in {\cal Z}(\tilde{G}'_{12};\bar{F})$, on \'ecrit $\underline{\delta}_{1}=x_{1}e'_{1}$, $\underline{\delta}_{2}=x_{2}e'_{2}$ avec $(x_{1},x_{2})\in G'_{12}({\mathbb A}_{\bar{F}})$. L'image de $(x_{1},x_{2})$ dans $  G'_{12}({\mathbb A}_{\bar{F}})/Z(G'_{12};\bar{F})$ ne d\'epend pas des choix de $e'_{1}$ et $e'_{2}$ et est invariante par $\Gamma_{F}$. L'application cherch\'ee est $(\underline{\delta}_{1},\underline{\delta}_{2})\mapsto (x_{1},x_{2})$. Par cette application, $\underline{\tilde{\lambda}}_{12}$ devient une fonction sur $\tilde{G}'_{12}({\mathbb A}_{F})$. Les m\^emes calculs qu'en [I] 2.5 conduisent \`a l'\'egalit\'e
$$(4) \qquad \Delta_{2}(\delta_{2},\gamma)=\underline{\tilde{\lambda}}_{12}(\delta_{1},\delta_{2})\Delta_{1}(\delta_{1},\gamma).$$
Or il r\'esulte des constructions que les fonctions $\tilde{\lambda}_{12}$ et $\underline{\tilde{\lambda}}_{12}$ se transforment selon le m\^eme caract\`ere $\lambda_{12}$  de $G'_{12}({\mathbb A}_{F})$. Par d\'efinition, $\tilde{\lambda}_{12}$ vaut $1$ sur $\tilde{G}'_{12}(F)$ et la construction ci-dessus montre que $\underline{\tilde{\lambda}}_{12}$ v\'erifie la m\^eme propri\'et\'e. Les deux fonctions sont donc \'egales et l'\'egalit\'e (4) \'equivaut \`a (3). $\square$

 \bigskip
 
 \subsection{Une construction auxiliaire}
Soit ${\bf G}'=(G',{\cal G}',\tilde{s})$ une donn\'ee endoscopique relevante de $(G,\tilde{G},{\bf a})$. Soit $(H,\tilde{H},{\bf b})$ un triplet similaire \`a $(G,\tilde{G},{\bf a})$ et soit ${\bf H}'=(H',{\cal H}',\tilde{t})$  une donn\'ee endoscopique  pour ce triplet. Consid\'erons les hypoth\`eses (1) \`a (6) suivantes.

(1) Il y a une suite exacte
$$1\to G\stackrel{\iota}{\to} H\to D\to 1$$
d'homomorphismes de groupes d\'efinis sur $F$, o\`u $D$ est un tore; il y a un plongement $\tilde{G}\stackrel{\tilde{\iota}}{\to} \tilde{H}$ d\'efini sur $F$ compatible avec $\iota$.

Pour $\tilde{h}\in \tilde{H}$, l'automorphisme $ad_{\tilde{h}}$ se quotiente en un automorphisme de $D$ qui ne d\'epend pas de $\tilde{h}$. On le note $\theta_{D}$. On ne demande pas qu'il soit l'identit\'e. De la suite (1) se d\'eduit une suite duale
$$1\to \hat{D}\to \hat{H}\stackrel {\hat{\iota}}{\to}\hat{G}\to 1$$
et une projection $\hat{\tilde{\iota}}:{^L\tilde{H}}\to {^L\tilde{G}}$ compatible avec $\hat{\iota}$ (on rappelle que $^L\tilde{G}={^LG}\hat{\boldsymbol{\theta}}$, cf. [I] 1.4). Notons $\hat{T}^H$ le tore d'une paire de Borel de $\hat{H}$ comme en [I] 1.4.

(2) $\hat{D}\cap \hat{T}^{H,\hat{\boldsymbol{\theta}},0}=\hat{D}^{\hat{\boldsymbol{\theta}},0}$.

(3) On a l'\'egalit\'e $\tilde{s}=\hat{\tilde{\iota}}(\tilde{t})$.

Puisque $\hat{H}'$ et $\hat{G}'$ sont les composantes neutres des commutants de $\tilde{t}$ et $\tilde{s}$, on a une suite exacte

$$1\to \hat{D}^{\hat{\boldsymbol{\theta}},0}\to \hat{H}'\to \hat{G}'\to 1.$$

(4) On a l'\'egalit\'e ${\cal G}'=\hat{\tilde{\iota}}({\cal H}')$.

Il en r\'esulte que ${\bf a}$ est le compos\'e de ${\bf b}$ et de la projection $Z(\hat{H})\to Z(\hat{G})$.

(5) On a l'\'egalit\'e $V_{ram}({\bf H}')=V_{ram}({\bf H})$.

(6) La donn\'ee ${\bf H}'$ est relevante et  v\'erifie l'hypoth\`ese (Hyp) de  [I] 6.4. 

 Consid\'erons pour $i=1,2$ des familles de donn\'ees $(H_{i},\tilde{H}_{i},{\bf b}_{i})$ et ${\bf H}'_{i}=(H'_{i},{\cal H}'_{i},\tilde{t}_{i})$  v\'erifiant  les hypoth\`eses (1) \`a (6).  On peut dire que la famille index\'ee par $2$ domine la famille index\'ee par $1$ s'il existe un homomorphisme injectif
 $$\kappa:H_{1}\to H_{2}$$
 et une application compatible $\tilde{\kappa}:\tilde{H}_{1}\to \tilde{H}_{2}$ de sorte que les hypoth\`eses suivantes soient v\'erifi\'ees:
 
 - les diagrammes
 $$\begin{array}{ccc}&&H_{1}\\ &\nearrow \iota_{1}&\\ G&&\downarrow \kappa\\ &\searrow \iota_{2}& \\ &&H_{2}\\ \end{array}$$
 $$\begin{array}{ccc}&&\tilde{H}_{1}\\ &\nearrow \tilde{\iota}_{1}&\\ \tilde{G}&&\downarrow \tilde{\kappa}\\ &\searrow \tilde{\iota}_{2}& \\ &&\tilde{H}_{2}\\ \end{array}$$
 sont commutatifs;
 
 en notant $\hat{\kappa}:{^LH_{2}}\to {^LH_{1}}$ et $\hat{\tilde{\kappa}}:{^L\tilde{H}_{2}}\to {^L\tilde{H}_{1}}$ les applications d\'eduites de $\kappa$,
 
 - ${^L\tilde{{\kappa}}}(\tilde{t}_{2})=\tilde{t}_{1}$ et $\hat{\kappa}({\cal H}'_{2})={\cal H}'_{1}$.

\ass{Lemme}{(i) Il existe des donn\'ees v\'erifiant les hypoth\`eses (1) \`a (6).

(ii) Pour deux familles de donn\'ees $(H_{i},\tilde{H}_{i},{\bf b}_{i})$ et ${\bf H}'_{i}=(H'_{i},{\cal H}'_{i},\tilde{t}_{i})$ pour $i=1,2$ v\'erifiant toutes deux les hypoth\`eses (1) \`a (6), il existe une troisi\`eme famille v\'erifiant les m\^emes hypoth\`eses et les dominant toutes deux.}

Preuve. Notons $T^*$ $\underline{le}$ tore maximal de $G$ muni de l'action galoisienne quasi-d\'eploy\'ee. Il est aussi muni de l'automorphisme $\theta^*$. Posons $H=(G\times T^*)/diag_{-}(Z(G))$ o\`u $diag_{-}$ est le plongement anti-diagonal et notons $\tilde{H}$ le quotient de $\tilde{G}\times T^*$ par $Z(G)$ agissant anti-diagonalement par multiplication \`a gauche. On d\'efinit deux actions de $G\times T^*$ sur $\tilde{H}$ par
$$(g,t)(\gamma,\tau)(g',t')=(g\gamma g',t\tau \theta^*(t')).$$
Ces actions se descendent en des actions de $H$ sur $\tilde{H}$ qui font de $\tilde{H}$ un espace tordu sous $H$. On a une suite exacte
$$1\to G\stackrel{\iota}{\to} H\to D=T^*_{ad}\to 1$$
et un plongement compatible $\tilde{\iota}:\tilde{G}\to \tilde{H}$ qui \`a $\gamma$ associe l'image de $(\gamma,1)$ dans $\tilde{H}$. Notons que le centre de $H$ est $(Z(G)\times T^*)/diag_{-}(Z(G))\simeq T^*$. Donc

(7) le centre $Z(H)$ est connexe.

On choisit une paire de Borel \'epingl\'ee $\hat{{\cal E}}=(\hat{B},\hat{T},(\hat{E}_{\alpha})_{\alpha\in \Delta})$ de $\hat{G}$  adapt\'ee \`a $\tilde{s}$ et on \'ecrit $\tilde{s}=s\hat{\theta}$ avec $s\in\hat{T}$, cf. [I] 1.5. Cette paire se rel\`eve en une paire de Borel \'epingl\'ee de $\hat{H}$. On note $\hat{T}^{H}$ le tore de cette paire et encore $\hat{\theta}$ l'automorphisme de $\hat{H}$ associ\'e \`a cette paire. 
Prouvons que l'\'egalit\'e (2) est v\'erifi\'ee. Il s'agit de prouver que la suite 
$$1\to \hat{D}^{\hat{\boldsymbol{\theta}},0}\to \hat{T}^{H,\hat{\boldsymbol{\theta}},0}\to \hat{T}^{\hat{\boldsymbol{\theta}},0}\to 1$$
est exacte. Il revient au m\^eme de prouver que la suite
$$0\to X_{*}(\hat{D})^{\hat{\boldsymbol{\theta}}}\to X_{*}(\hat{T}^H)^{\hat{\boldsymbol{\theta}}}\to X_{*}(\hat{T})^{\hat{\boldsymbol{\theta}}}\to 0$$
est exacte. Seule la surjectivit\'e finale pose probl\`eme. Les actions galoisiennes n'interviennent pas ici. On peut travailler sur $\bar{F}$ et identifier $T^*$ \`a un sous-tore de $G$.  Posons $T^{*H}=(T^*\times T^*)/diag_{-}(Z(G))$, qui est un sous-tore maximal de $H$. La surjectivit\'e \`a prouver \'equivaut \`a celle de l'homomorphisme
$$X^*(T^{*H})^{\theta}\to X^*(T^*)^{\theta}$$
issue du plongement $t\mapsto (t,1)$ de $ T^*$ dans $T^{*H}$. Mais on a aussi un homomorphisme $T^{*H}\to T^*$ d\'efini par $(t_{1},t_{2})\mapsto t_{1}t_{2}$ dont le compos\'e avec le pr\'ec\'edent est l'identit\'e de $T^*$. Ainsi l'homomorphisme ci-dessus s'inscrit dans une suite
$$X^*(T^*)^{\theta}\to X^*(T^{*H})^{\theta}\to X^*(T^*)^{\theta}$$
dont le compos\'e est l'identit\'e. La deuxi\`eme fl\`eche est donc surjective, comme on le voulait.

Pour $v\in Val(F)$, $v\not\in V_{ram}(\tilde{G},{\bf a})$, le groupe $H$ est non ramifi\'e sur $F_{v}$.   Puisque $G_{AD}=H_{AD}$, le sous-groupe compact hypersp\'ecial $K_{v}$ de $G(F_{v})$ d\'etermine un tel sous-groupe $K^H_{v}$ de $H(F_{v})$. L'espace $\tilde{K}^H_{v}=K^H_{v}\tilde{\iota}(\tilde{K}_{v})$ est un espace hypersp\'ecial pour ce groupe.

 On choisit $t\in \hat{T}^{H}$ d'image $s$ dans $\hat{T}$ et on pose $\tilde{t}=t\hat{\theta}$. La relation (3) est v\'erifi\'ee. On peut identifier $\hat{T}^{\hat{\theta},0}$ au tore d'une paire de Borel \'epingl\'ee de $\hat{G}'$. Il s'en d\'eduit une structure galoisienne sur ce tore, de la forme $\sigma\mapsto \sigma_{G'}=\omega_{G'}(\sigma)\circ \sigma_{G}$, o\`u $\omega_{G'}$ est un cocycle \`a valeurs dans $W^{\theta}$. Ce groupe $W^{\theta}$ est le m\^eme pour $G$ ou $H$. On peut donc relever l'action pr\'ec\'edente en une action $\sigma\mapsto \sigma_{H'}=\omega_{G'}(\sigma)\circ \sigma_{H}$ de $\Gamma_{F}$ sur $\hat{T}^{H,\hat{\theta},0}$. Ces actions se prolongent en des actions sur $\hat{T}$ et $\hat{T}^{H}$. Remarquons que ces actions sont non ramifi\'ees en $v$ pour $v\in Val(F)-V_{ram}({\bf G}')$. Notons $\hat{T}'$, resp. $\hat{T}^{G'}$, $\hat{T}^{'H}$, $\hat{T}^{H'}$, les tores $\hat{T}$, resp. $\hat{T}^{\hat{\theta},0}$, $\hat{T}^{H}$,  $\hat{T}^{H,\hat{\theta},0}$, munis de ces structures.   On  note $T^{G'}$ et $T^{H'}$ les tores duaux de $\hat{T}^{G'}$ et $\hat{T}^{H'}$ d\'efinis sur $F$. D'apr\`es le lemme 3.2 et la relation (3) de sa preuve, on peut prolonger le plongement $\hat{T}^{H'}\to \hat{H}^{\hat{\theta},0}$, resp.  $\hat{T}^{G'}\to {\cal G}'$ en des plongements
$$\begin{array}{ccc}{^LT}^{H'}&\to& \hat{H}^{\hat{\theta},0}\rtimes W_{F}\\ (x,w)&\mapsto &(xh^1(w),w),\\ \end{array}$$
resp. 
$$\begin{array}{ccc} {^LT}^{G'}&\to&  {\cal G}'\subset {^LG}\\ (x,w)&\mapsto &(xg'(w),w),\\ \end{array}$$
tels que, pour $v\not\in V_{ram}({\bf G}')$ et $w\in I_{v}$, on ait $h^1(w)=1$ et $g'(w)=1$.
Quotientons le premier par $\hat{D}\cap \hat{T}^{H'}=\hat{D}\cap \hat{T}^{H,\hat{\theta},0}$ (on rappelle que $\hat{D}$ est le tore dual de $D\simeq T^*_{ad}$, d'o\`u $\hat{D}\simeq \hat{T}_{sc}$). On obtient un plongement
$$\begin{array}{ccc}{^LT}^{G'}&\to& \hat{G}^{\hat{\theta},0}\rtimes W_{F}\\ (x,w)&\mapsto &(x\hat{\iota}(h^1(w)),w)\\ \end{array}$$
 Les deux plongements pr\'ec\'edents ne peuvent diff\'erer que par un cocycle. C'est-\`a-dire qu'il existe un cocycle $u:W_{F}\to \hat{T}$ tel que $g'(w)=u(w)\hat{\iota}(h^1(w))$ pour tout $w$. Pour $v\in Val(F)-V_{ram}({\bf G}')$ et $w\in I_{v}$, on a $u(w)=1$. En appliquant 3.2(2) \`a la suite exacte 
$$1\to \hat{D}\to \hat{T}^{H}\to \hat{T}\to 1,$$
on voit que l'on peut relever $u$ en un cocycle $u^{H}:W_{F}\to \hat{T}^{H}$ tel que $\hat{\iota}\circ u^{H}=u$ de sorte que, pour tous  $v\in Val(F)-V_{ram}({\bf G}')$ et $w\in I_{v}$, on ait $u^{H}(w)=1$ . On pose $h'(w)=u^{H}(w)h^1(w)$.    L'application
$$\begin{array}{ccc}{^LT}^{H'}&\to& {^LH}\\ (x,w)&\mapsto &(x h'(w),w)\\ \end{array}$$
est alors un homomorphisme. Pour $w\in W_{F}$, soit $b(w)\in \hat{H}$ tel que
$$(8) \qquad ad_{\tilde{t}}(h'(w),w)=(b(w)h'(w),w).$$
En projetant dans $^LG$ on voit que $b(w)$ se projette sur l'\'el\'ement $a(w)$ de $Z(\hat{G})$ tel que $ad_{\tilde{s}}(g'(w),w)=(a(w)g'(w),w)$. Donc $b(w)\in Z(\hat{H})$. L'\'equation (8) oblige $b$ \`a \^etre un cocycle. On note ${\bf b}$ sa classe modulo $ker^1(W_{F};Z(\hat{H}))$. Elle se projette sur ${\bf a}$. On note $\hat{H}'=Z_{\hat{H}}(\tilde{t})^0$. L'\'equation (8) oblige $(h'(w),w)$ \`a normaliser $\hat{H}'$. On note ${\cal H}'$ le groupe engendr\'e par $\hat{H}'$ et les $(h'(w),w)$ pour $w\in W_{F}$.  C'est une extension de $\hat{H}'$. Elle d\'etermine une action galoisienne sur ce groupe qui en conserve une paire de Borel \'epingl\'ee: par exemple le rel\`evement dans $\hat{H}'$ d'une  telle paire de $\hat{G}'$ conserv\'ee par l'action galoisienne. On note $H'$ le groupe quasi-d\'eploy\'e sur $F$ dont le $L$-groupe est $\hat{H}'$ muni de cette structure galoisienne. Alors $(H',{\cal H}',\tilde{t})$ est une donn\'ee endoscopique pour $(H,\tilde{H},{\bf b})$ et la relation (3) est v\'erifi\'ee. Cette donn\'ee est non ramifi\'ee en tout $v\in Val(F)-V_{ram}({\bf G}') $ car, pour une telle place, on a $h'(w)=1$ pour $w\in I_{v}$.
Montrons que cette donn\'ee est relevante. Dualement \`a l'homomorphisme $\hat{H}'\stackrel{\hat{\tilde{\iota}}}{\to }\hat{G}'$, on a un homomorphisme $\iota':G'\to H'$ qui est d\'efini sur $F$. On a aussi des plongements compatibles ${\cal Z}(G)\to {\cal Z}(H)$ et ${\cal Z}(\tilde{G})\to {\cal Z}(\tilde{H})$. L'homomorphisme $\iota'$ se prolonge en 
$$\tilde{\iota}':\tilde{G}'=G'\times_{{\cal Z}(G)}{\cal Z}(\tilde{G})\to\tilde{H}'= H'\times_{{\cal Z}(H)}{\cal Z}(\tilde{H}).$$
Puisque $\tilde{G}'(F)$ n'est pas vide, $\tilde{H}'(F)$ ne l'est pas non  plus.
Remarquons que, par construction, on a $G_{SC}=H_{SC}$ et $G'_{SC}=H'_{SC}$. Il y a donc une bijection entre paires de Borel pour $G$ et pour $H$, et entre paires de Borel pour $G'$ et pour $H'$. Soit $v\in Val(F)$. Puisque ${\bf G}'$ est relevante, on peut fixer un diagramme $(\delta,B',T',B,T,\gamma)$ avec $\delta\in \tilde{G}'(F_{v})$, $\gamma\in \tilde{G}(F_{v})$ et $\gamma$ fortement r\'egulier. Notons $\delta^{H'}$ et $\gamma^{H}$ les images de $\delta$ dans $\tilde{H}'(F_{v})$ et de $\gamma$ dans $\tilde{H}(F_{v})$. Notons $(B^{H'},T^{H'})$ la paire de Borel de $H'$ correspondant \`a $(B',T')$ et $(B^{H},T^{H})$ la paire de Borel de $H$ correspondant \`a $(B,T)$. Alors $(\delta^{H'},B^{H'},T^{H'},B^{H},T^{H},\gamma^{H})$ est un diagramme et $\gamma^{H}$ est fortement r\'egulier. Donc ${\bf H}'_{v}$ est relevante.

Pour achever la preuve de (i), il reste \`a prouver que ${\bf H}'$ v\'erifie l'hypoth\`ese (Hyp). Elle va \^etre assur\'ee par (7). On peut aussi bien revenir aux donn\'ees initiales et prouver

(9) si $Z(G)$est connexe et ${\bf G}'$ est relevante, alors ${\bf G}'$ v\'erifie l'hypoth\`ese (Hyp).

Pour $v\in V_{ram}({\bf G}')$, l'hypoth\`ese que ${\bf G}'_{v}$ est relevante permet de fixer un sous-tore maximal  $T'_{v}$ de $G'$, d\'efini sur $F_{v}$, tel qu'il existe $(\delta_{v},\gamma_{v})\in {\cal D}({\bf G}'_{v})$ de sorte que $\delta_{v}\in \tilde{T}'_{v}(F_{v})$. Fixons un \'el\'ement  $Y_{v}\in \mathfrak{t}'_{v}(F_{v})$ r\'egulier dans $\mathfrak{g}'(F_{v})$. On peut fixer un \'el\'ement $Y\in \mathfrak{g}'(F)$ dont la composante en $v$ soit aussi proche que l'on veut de $Y_{v}$ pour tout $v\in V_{ram}({\bf G}')$. Notons $T'$ le commutant de $Y$. C'est un sous-tore maximal de $G'$, d\'efini sur $F$. Si la composante de $Y$ en $v$ est assez proche de $Y_{v}$, ce tore est conjugu\'e \`a $T'_{v}$ par un \'el\'ement de $G'(F_{v})$. Il v\'erifie donc la m\^eme condition que $T'_{v}$. Il faut montrer qu'il v\'erifie aussi  cette condition pour $v\not\in V_{ram}({\bf G}')$.    Pour une telle place, on peut identifier $\underline{la}$ paire de Borel \'epingl\'ee de $G$ \`a une paire ${\cal E}^*=(B^*,T^*,(E^*_{\alpha})_{\alpha\in \Delta})$ d\'efinie sur $F_{v}$ et dont $K_{v}$ soit le groupe hypersp\'ecial associ\'e. D'apr\`es [I] 6.2, on peut fixer $e\in Z(\tilde{G},{\cal E}^*)(F_{v}^{nr})\cap T^*(\mathfrak{o}^{nr}_{v})\tilde{K}_{v}$, avec les notations de cette r\'ef\'erence. Soit $z:\Gamma_{F_{v}}\to Z(G;\bar{F}_{v})$ l'application telle que $\sigma(e)=z(\sigma)^{-1}e$. Alors $z$ est un cocycle non ramifi\'e \`a valeurs dans $Z(G;\bar{F}_{v})\cap T^*(\mathfrak{o}^{nr}_{v})=Z(G;\mathfrak{o}^{nr}_{v})$. Puisque $Z(G)$ est connexe, un tel cocycle est trivial. Quitte \`a multiplier $e$ par un \'el\'ement de $Z(G;\mathfrak{o}^{nr}_{v})$, on a donc $e\in Z(\tilde{G},{\cal E}^*)(F_{v})\cap \tilde{K}_{v}$. Posons $\theta=ad_{e}$ et fixons un sous-groupe de Borel $B'$ de $G'$ contenant $T'$. Gr\^ace \`a [K1]  corollaire 2.2, on peut fixer $x\in G_{SC}^{\theta}(\bar{F}_{v})$ de sorte qu'en posant $ad_{x^{-1}}(B^*,T^*)=(B_{v},T_{v})$, le tore $T_{v}$ soit d\'efini sur $F_{v}$ et l'homomorphisme $\xi_{T_{v},T'}$ d\'eduit de $(B_{v},T_{v})$ et de $(B',T')$ soit \'equivariant pour les actions galoisiennes. Posons ${\cal E}_{v}=ad_{x^{-1}}({\cal E}^*)$. Puisque $x$ commute \`a $\theta$, on a encore $e\in Z(\tilde{G},{\cal E}_{v})(F_{v})$.
On fixe $\nu\in T_{v}(F_{v})$ en position g\'en\'erale, on pose $\mu=\xi_{T_{v},T'}(\nu)$, $\gamma_{v}=\nu e$ et $\delta_{v}=\mu e'$, o\`u $e'\in {\cal Z}(\tilde{G}';F_{v})$ est l'image de $e$. Alors $(\delta_{v},B',T',B,T,\gamma_{v})$ est un diagramme avec $\delta_{v}\in \tilde{T}'(F_{v})$ et $\gamma_{v}$ fortement r\'egulier. Cela d\'emontre (9) et le (i) de l'\'enonc\'e. 

Soient maintenant deux familles comme dans le (ii) de l'\'enonc\'e. Pour $i=1,2$, on peut \'ecrire $H_{i}=(G\times Z(H_{i}))/diag_{-}(Z(G))$  et identifier $\tilde{H}_{i}$ au quotient de $(\tilde{G}\times Z(H_{i}))$ par $Z(G)$ agissant anti-diagonalement par multiplication \`a gauche. Posons $Z_{12}(G)=\{(z,z_{1},z_{2})\in Z(G); zz_{1}z_{2}=1\}$. Posons $H=(G\times Z(H_{1})\times Z(H_{2}))/Z_{12}(G)$ et notons $\tilde{H}$ le quotient de $\tilde{G}\times Z(H_{1})\times Z(H_{2})$ par l'action de $Z_{12}(G)$ par multiplication \`a gauche. On munit $\tilde{H}$ d'une structure d'espace tordu sur $H$ comme au d\'ebut de la preuve de (i). Il y a un diagramme naturel d'homomorphismes
$$\begin{array}{ccccc}&&H_{1}&&\\ &\nearrow \iota_{1}&&\searrow \kappa_{1}&\\ G&&&& H\\ &\searrow \iota_{2}&&\nearrow \kappa_{2}&\\ && H_{2}&&\\ \end{array}$$
On v\'erifie qu'ils ont tous injectifs. L'homomorphisme compos\'e $\iota$ s'ins\`ere dans une suite exacte
$$1\to G\stackrel{\iota}{\to }H \to D_{1 }\times D_{2}\to 1$$
 Tous les homomorphismes se prolongent en des applications compatibles entre les espaces tordus correspondants. Du c\^ot\'e dual, $\hat{H}$ est le produit fibr\'e de $\hat{H}_{1}$ et $\hat{H}_{2} $ au-dessus de $\hat{G}$. Comme dans la preuve de (i), on fixe une paire de Borel \'epingl\'ee de $\hat{G}$ de tore $\hat{T}$ de sorte que $\tilde{s}=s\hat{\theta}$ avec $s\in \hat{T}$. Elle se rel\`eve en des paires pour $\hat{H}_{1}$ et $\hat{H}_{2}$ de tores $\hat{T}_{1}$ et $\hat{T}_{2}$. Pour $i=1,2$, on a $\tilde{t}_{i}=t_{i}\hat{\theta}$  avec $t_{i}\in \hat{T}_{i}$. On  pose $t=(t_{1},t_{2})$ et $\tilde{t}=t\hat{\theta}$. On d\'efinit ${\cal H}'$ comme l'ensemble des \'el\'ements $(x_{1},x_{2},w)\in {^LH}$ tels que $(x_{1},w)\in {\cal H}'_{1}$ et $(x_{2},w)\in {\cal H}'_{2}$. Comme dans la preuve de (i), on associe \`a ces donn\'ees un groupe $H'$ d\'efini et quasi-d\'eploy\'e sur $F$, ainsi qu'une classe de cocycle ${\bf b}$. On laisse le lecteur v\'erifier que les donn\'ees $(H,\tilde{H},{\bf b})$ et ${\bf H}'=(H',{\cal H}',\tilde{t})$ satisfont les conditions requises. $\square$

 {\bf Remarque.} La preuve fournit un groupe $H$ et un espace $\tilde{H}$ qui sont ind\'ependants de la donn\'ee endoscopique ${\bf G}'$.

\bigskip

\subsection{Facteur de transfert global, cas g\'en\'eral}

Soit ${\bf G}'=(G',{\cal G}',\tilde{s})$ une donn\'ee endoscopique relevante de $(G,\tilde{G},{\bf a})$.  
 Soit $V$ un ensemble fini de places de $F$ contenant $V_{ram}({\bf G}')$. 
 Consid\'erons des donn\'ees $(H,\tilde{H},{\bf b})$ et ${\bf H}'=(H',{\cal H}',\tilde{t})$ v\'erifiant les hypoth\`eses (1) \`a (6) du paragraphe pr\'ec\'edent.    Consid\'erons aussi des donn\'ees auxiliaires $H'_{1}$, $\tilde{H}'_{1}$, $C_{1,H}$, $\hat{\xi}_{1,H}$ pour ${\bf H}'$, non ramifi\'ees hors de $V$. On a la projection $H'_{1}\to H'$. On a aussi un homomorphisme $\iota':G'\to H'$, d\'efini sur $F$, et la projection duale $\hat{\iota}':\hat{H}'\to \hat{G}'$. Notons $G'_{1}$ le produit fibr\'e de $H'_{1}$ et $G'$ au-dessus de $H'$ et posons $C_{1}=C_{1,H}$. On a la suite exacte
 $$1\to C_{1}\to G'_{1}\to G'\to 1$$
  Comme on l'a vu dans la preuve de 3.8, de $\iota'$ se d\'eduit  une application $\tilde{\iota}':\tilde{G}'=G'\times_{ {\cal Z}(G)}{\cal Z}(\tilde{G})\to \tilde{H}'=H'\times_{{\cal Z}(H)}{\cal Z}(\tilde{H})$. On d\'efinit $\tilde{G}'_{1}$ comme le produit fibr\'e de $\tilde{H}'_{1}$ et $\tilde{G}'$ au-dessus de $\tilde{H}'$. La projection naturelle $\tilde{G}'_{1}\to \tilde{G}'$ est compatible avec la suite exacte ci-dessus. Du c\^ot\'e des groupes duaux, on v\'erifie que l'on a un diagramme commutatif
 $$\begin{array}{ccccccccc}&&1&&1&&&&\\&&\downarrow&&\downarrow&&&&\\&&\hat{D}^{\hat{\theta},0}&=&\hat{D}^{\hat{\theta},0}&&&&\\&&\downarrow&&\downarrow&&&&\\ 1&\to&\hat{H}'&\to&\hat{H}'_{1}&\to&\hat{C}_{1}&\to&1\\&&\downarrow&&\downarrow&&\parallel&&\\ 1&\to&\hat{G}'&\to&\hat{G}'_{1}&\to&\hat{C}_{1}&\to&1\\ &&\downarrow&&\downarrow&&&&\\ &&1&&1&&&&\\ \end{array}$$
 dont toutes les suites sont exactes. On a aussi une suite exacte
 $$1\to \hat{D}^{\hat{\theta},0}\to {\cal H}'\to {\cal G}'\to 1$$
 En quotientant par $\hat{D}^{\hat{\theta},0}$, il se d\'eduit du plongement $\hat{\xi}_{1,H}:{\cal H}'\to {^LH}'_{1}$ un plongement
 $$\hat{\xi}_{1}:{\cal G}'\to {^LG}'_{1}$$
 Les donn\'ees $G'_{1}$, $\tilde{G}'_{1}$, $C_{1}$, $\hat{\xi}_{1}$ sont des donn\'ees auxiliaires pour ${\bf G}'$ qui sont non ramifi\'ees hors de $V$.  On les compl\`ete par une famille d'espaces hypersp\'eciaux $(\tilde{K}'_{1,v})_{v\not\in V}$ v\'erifiant les conditions usuelles, cf. 1.1.  Les groupes $G'$ et $H'$ ont m\^eme groupe adjoint. Pour $v\not\in V$, le sous-groupe compact hypersp\'ecial $K'_{v}$ d\'etermine donc de tels sous-groupes $K'_{H,v}$ de $H'(F_{v})$ puis $K'_{1,H,v}$ de $H'_{1}(F_{v})$.     Alors l'ensemble $K'_{1,H,v}\tilde{\iota}_{1}'(\tilde{K}'_{1,v})$ est un espace hypersp\'ecial de $\tilde{H}'_{1,v}(F_{v})$. On le note $\tilde{K}'_{1,H,v}$. La famille $(\tilde{K}'_{1,H,v})_{v\not\in V}$ v\'erifie la condition de compatibilit\'e globale de 1.1. On compl\`ete les donn\'ees auxiliaires $H'_{1}$ etc... par cette famille. 
 
Pour toute place $v$, on introduit les ensembles ${\cal D}_{v}$ et ${\cal D}_{1,v}$ relatifs \`a $\tilde{G}'$ et $\tilde{G}'_{1}$ sur $F_{v}$ et les ensembles similaires ${\cal D}_{H,v}$ et ${\cal D}_{1,H,v}$ relatifs \`a $\tilde{H}'$ et $\tilde{H}'_{1}$. On a vu dans la preuve de 3.8 que pour $(\delta,\gamma)\in {\cal D}_{v}$, on a $(\tilde{\iota}'(\delta),\tilde{\iota}(\gamma))\in {\cal D}_{H,v}$. Il en r\'esulte que, pour $(\delta_{1},\gamma)\in {\cal D}_{1,v}$, on a $(\tilde{\iota}_{1}'(\delta_{1}),\tilde{\iota}(\gamma))\in {\cal D}_{1,H,v}$, o\`u $\tilde{\iota}'_{1}:\tilde{G}'_{1}\to \tilde{H}'_{1}$ est l'application naturelle. Puisque ${\bf H}'$ v\'erifie l'hypoth\`ese (Hyp), on peut lui appliquer les constructions de la preuve de la proposition 3.7: de nos choix d'espaces hypersp\'eciaux se d\'eduit un facteur de transfert normalis\'e, notons-le $\Delta_{1,H,V}$ sur ${\cal D}_{1,H,V}=\prod_{v\in V}{\cal D}_{1,H,v}$. On d\'efinit une fonction $\Delta_{1,V}$ sur ${\cal D}_{1,V}=\prod_{v\in V}{\cal D}_{1,v}$ par $\Delta_{1,V}(\delta_{1},\gamma)=\Delta_{1,H,V}(\tilde{\iota}'_{1}(\delta_{1}),\tilde{\iota}(\gamma))$.  On a

(1) $\Delta_{1,V}$ est un facteur de transfert.

Preuve. Puisque $\Delta_{1,H,V}$ en est un, il suffit de prouver que

  (2) pour  toute place $v$ et tous  $(\delta_{1},\gamma), (\underline{\delta}_{1},\underline{\gamma})\in {\cal D}_{1,v}$, on a l'\'egalit\'e
 $$\boldsymbol{\Delta}_{1,H,v}(\tilde{\iota}'_{1}(\delta_{1}),\tilde{\iota}(\gamma);\tilde{\iota}'_{1}(\underline{\delta}_{1}),\tilde{\iota}(\underline{\gamma}))=\boldsymbol{\Delta}_{1,v}(\delta_{1},\gamma;\underline{\delta}_{1},\underline{\gamma}).$$
 
 On reprend les d\'efinitions de [I] 2.2 en ajoutant judicieusement des indices $H$ pour les termes relatifs \`a $\tilde{H}$. Les facteurs $\Delta_{II}$ intervenant sont les m\^emes des deux c\^ot\'es car ces facteurs sont insensibles aux centres et on a $G_{AD}=H_{AD}$, $G'_{AD}=H'_{AD}$. Il faut comparer les facteurs $\boldsymbol{\Delta}_{imp,H,v}$ et $\boldsymbol{\Delta}_{imp,v}$. On a des \'egalit\'es
 $$\boldsymbol{\Delta}_{imp,v}(\delta_{1},\gamma;\underline{\delta}_{1},\underline{\gamma})=<(V,\boldsymbol{\nu}_{1}),(\hat{V}_{1},{\bf s})>^{-1},$$
 $$\boldsymbol{\Delta}_{imp,H,v}(\tilde{\iota}_{1}'(\delta_{1}),\tilde{\iota}(\gamma);\tilde{\iota}'_{1}(\underline{\delta}_{1}),\tilde{\iota}(\underline{\gamma}))=<(V_{H},\boldsymbol{\nu}_{1,H}),(\hat{V}_{1,H},{\bf t})>^{-1},$$
 les produits \'etant respectivement ceux sur
 $$H^{1,0}(\Gamma_{F_{v}};U\stackrel{1-\theta}{\to} S_{1})\times H^{1,0}(W_{F_{v}}; \hat{S}_{1}\stackrel{1-\hat{\theta}}{\to}\hat{U})$$
 et
  $$H^{1,0}(\Gamma_{F_{v}};U_{H}\stackrel{1-\theta}{\to} S_{1,H})\times H^{1,0}(W_{F_{v}}; \hat{S}_{1,H}\stackrel{1-\hat{\theta}}{\to}\hat{U}_{H}).$$
  En fait, on a $U_{H}=U$ et des homomorphismes duaux $S_{1}\to S_{1,H}$, $\hat{S}_{1,H}\to \hat{S}_{1}$. En choisissant convenablement les donn\'ees auxiliaires intervenant, on v\'erifie que $(\hat{V}_{1},{\bf s})$ est l'image naturelle de $(\hat{V}_{1,H},{\bf t})$ par l'homorphisme
  $$H^{1,0}(W_{F_{v}}; \hat{S}_{1,H}\stackrel{1-\hat{\theta}}{\to}\hat{U}_{H})\to H^{1,0}(W_{F_{v}}; \hat{S}_{1}\stackrel{1-\hat{\theta}}{\to}\hat{U}),$$
  tandis que $(V_{H},\boldsymbol{\nu}_{1,H})$ est l'image naturelle de $(V,\boldsymbol{\nu}_{1})$ par l'homomorphisme dual
  $$H^{1,0}(\Gamma_{F_{v}};U\stackrel{1-\theta}{\to} S_{1})\to H^{1,0}(\Gamma_{F_{v}};U_{H}\stackrel{1-\theta}{\to} S_{1,H}).$$
  L'\'egalit\'e (2) r\'esulte alors simplement de la compatibilit\'e des produits. Cela prouve (2) et (1). $\square$
  
   Pour $v\not\in V$, on a deux facteurs de transfert normalis\'es $\Delta_{1,v}$ sur ${\cal D}_{1,v}$ et $\Delta_{1,H,v}$ sur ${\cal D}_{1,H,v}$. Pour $(\delta_{1,v},\gamma_{v})\in {\cal D}_{1,v}$, on a l'\'egalit\'e
 
 (3) $\Delta_{1,v}(\delta_{1,v},\gamma_{v})=\Delta_{1,H,v}(\tilde{\iota}'_{1}(\delta_{1,v}),\tilde{\iota}(\gamma_{v}))$.
 
 La preuve est similaire \`a celle de (2).

    Comme en 3.7, de l'existence du facteur de transfert  $\Delta_{1,V}$ va r\'esulter la proposition suivante.
  
    \ass{Proposition}{Il existe un isomorphisme canonique
 $$C_{c}^{\infty}({\bf G}'_{V})\simeq \otimes_{v\in V}C_{c}^{\infty}({\bf G}'_{v}).$$}
 
 Preuve. On choisit des donn\'ees $(H,\tilde{H},{\bf b})$ etc... et  on construit le facteur de transfert $\Delta_{1,V}$ sur ${\cal D}_{1,V}$.  
  Comme en 3.7, on a alors les isomorphismes
  $$C_{c}^{\infty}({\bf G}'_{V})\simeq C_{c,\lambda_{1}}^{\infty}(\tilde{G}'_{1}(F_{V}))\simeq \otimes_{v\in V}C_{c}^{\infty}({\bf G}'_{v}),$$
  le premier \'etant relatif aux donn\'ees   $(\tilde{K}'_{1,v})_{v\not\in V}$ et le second au facteur de transfert $\Delta_{1,V}$. Le compos\'e de ces isomorphismes fournit celui de l'\'enonc\'e.

  On doit montrer qu'il est "canonique", c'est-\`a-dire qu'il ne d\'epend pas des donn\'ees auxiliaires.

 Conservons les donn\'ees $(H,\tilde{H},{\bf b})$, ${\bf H}'=(H',{\cal H}',\tilde{t})$, mais rempla\c{c}ons  $H'_{1}$, $\tilde{H}'_{1}$, $C_{1,H}$, $\hat{\xi}_{1,H}$ par d'autres donn\'ees auxiliaires  $H'_{2}$, $\tilde{H}'_{2}$, $C_{2,H}$, $\hat{\xi}_{2,H}$. On en d\'eduit de nouvelles donn\'ees auxiliaires $G'_{2}$, $\tilde{G}'_{2}$, $C_{2}$, $\hat{\xi}_{2}$ pour ${\bf G}'$. On fixe des familles  $(\tilde{K}'_{1,v})_{v\not\in V}$ et $(\tilde{K}'_{2,v})_{v\not\in V}$ o\`u $\tilde{K}'_{1,v}$, resp. $\tilde{K}'_{2,v}$ est un sous-espace hypersp\'ecial de $\tilde{G}'_{1}(F_{v})$, resp. $\tilde{G}'_{2}(F_{v})$. Il s'en d\'eduit des familles  $(\tilde{K}'_{1,H,v})_{v\not\in V}$ et $(\tilde{K}'_{2,H,v})_{v\not\in V}$ o\`u $\tilde{K}'_{1,H,v}$, resp. $\tilde{K}'_{2,H,v}$ est un sous-espace hypersp\'ecial de $\tilde{H}'_{1}(F_{v})$, resp. $\tilde{H}'_{2}(F_{v})$. Comme dans la preuve de 3.7, on a deux isomorphismes
 $$C_{c,\lambda_{1}}^{\infty}(\tilde{G}'_{1}(F_{V}))\simeq C_{c,\lambda_{2}}^{\infty}(\tilde{G}'_{2}(F_{V}))$$
 dont nous voulons prouver qu'ils sont \'egaux. Ils sont donn\'es par des fonctions de recollement $\tilde{\lambda}_{12,K,V}$ et $\tilde{\lambda}_{12,\Delta,V}$ sur $\tilde{G}'_{12}(F_{V})$. Comme en 3.7, il s'agit de prouver que ces fonctions sont \'egales. Mais pour ${\bf H}'$, on a aussi des fonctions de recollement $\tilde{\lambda}_{12,H,K,V}$ et $\tilde{\lambda}_{12,H,\Delta,V}$ sur $\tilde{H}'_{12}(F_{V})$.  Des applications $\tilde{\iota}'_{1}$ et $\tilde{\iota}'_{2}$ (avec des notations \'evidentes) se d\'eduit une application $\tilde{\iota}'_{12}:\tilde{G}'_{12}\to \tilde{H}'_{12}$. Il r\'esulte des d\'efinitions et de (2) et (3) que l'on a les \'egalit\'es
 $$\tilde{\lambda}_{12,K,V}=\tilde{\lambda}_{12,H,K,V}\circ\tilde{\iota}'_{12},\,\,\tilde{\lambda}_{12,\Delta,V}=\tilde{\lambda}_{12,H,\Delta,V}\circ\tilde{\iota}'_{12}.$$
 Parce que ${\bf H}$ v\'erifie l'hypoth\`ese (Hyp), on peut lui appliquer la preuve de 3.7 qui montre l'\'egalit\'e $\tilde{\lambda}_{12,H,K,V}=\tilde{\lambda}_{12,H,\Delta,V}$. L'\'egalit\'e cherch\'ee $\tilde{\lambda}_{12,K,V}=\tilde{\lambda}_{12,\Delta,V}$ s'ensuit.

  Rempla\c{c}ons maintenant  les donn\'ees auxiliaires $(H,\tilde{H},{\bf b})$, ${\bf H}'=(H',{\cal H}',\tilde{t})$, $H'_{1}$, $\tilde{H}'_{1}$, $C_{1,H}$, $\hat{\xi}_{1,H}$ par d'autres $(\bar{H},\tilde{\bar{H}},\overline{{\bf b}})$ etc... v\'erifiant les m\^emes conditions. En utilisant le lemme 3.8 (ii), on introduit des donn\'ees $(\underline{H},\underline{\tilde{H}},\underline{{\bf b}})$ et $\underline{{\bf H}}'=(\underline{H}',\underline{{\cal H}}',\underline{\tilde{t}}')$ qui dominent   $(H,\tilde{H},{\bf b})$ et ${\bf H}'$ ainsi que  $(\bar{H},\tilde{\bar{H}},\overline{{\bf b}})$ et $\overline{{\bf H}}'$. On peut fixer des donn\'ees auxiliaires $\underline{H}'_{1}$, $\underline{\tilde{H}}'_{1}$, $\underline{C}_{1,H}$, $\underline{\hat{\xi}}_{1,H}$. On peut d\'ecomposer la preuve en deux: prouver que l'isomorphisme ne change pas quand on remplace les donn\'ees  $H$ etc.. par les donn\'ees $\underline{H}$ etc... puis qu'il ne change pas quand on remplace les donn\'ees $\underline{H}$ etc... par $\bar{H}$ etc... Les deux assertions sont similaires. On peut ne d\'emontrer que la premi\`ere partie et oublier les donn\'ees $\bar{H}$ etc... On s'est ainsi ramen\'e au cas o\`u les donn\'ees $(\underline{H},\underline{\tilde{H}},\underline{{\bf b}})$ et $\underline{{\bf H}}'=(\underline{H}',\underline{{\cal H}}',\underline{\tilde{t}}')$  dominent   $(H,\tilde{H},{\bf b})$ et ${\bf H}'$. En particulier, on a des plongements compatibles $\kappa:H\to \underline{H}$, $\tilde{\kappa}:\tilde{H}\to \underline{\tilde{H}}$ et, dualement, des plongements $\hat{\kappa}$ et $\hat{\tilde{\kappa}}$. On a aussi des homomorphismes compatibles $\kappa':H'\to \underline{H}'$ et $\tilde{\kappa}':\tilde{H}'\to \underline{\tilde{H}}'$.  Le m\^eme proc\'ed\'e qui nous a permis de d\'eduire de $H'_{1}$ etc... des donn\'ees auxiliaires $G'_{1}$ etc... nous permet maintenant de d\'eduire des donn\'ees $\underline{H}'_{1}$ etc... des donn\'ees auxiliaires pour ${\bf H}'$, que l'on note $H'_{2}$, $\tilde{H}'_{2}$ etc... Par exemple, $H'_{2}$ est le produit fibr\'e de $H'$ et $\underline{H}'_{1}$ au-dessus de $\underline{H}'$. D'apr\`es ce que l'on a d\'ej\`a d\'emontr\'e, notre isomorphisme est insensible au changement de ces donn\'ees auxiliaires, on peut donc supposer que $H'_{1}=H'_{2}$ etc...  Il r\'esulte des constructions que les donn\'ees auxiliaires $G'_{1}$, $\tilde{G}'_{1}$, $C_{1}$, $\hat{\xi}_{1}$ pour ${\bf G}'$ d\'eduites de $H'_{1}$ etc... sont les m\^emes que celles d\'eduites de $\underline{H}'_{1}$ etc... On a d'ailleurs $C_{1}=C_{1,H}=C_{1,\underline{H}}$. On fixe une famille $(\tilde{K}'_{1,v})_{v\not\in V}$ d'espaces hypersp\'eciaux, dont on d\'eduit des familles $(\tilde{K}'_{1,H,v})_{v\not\in V}$ et $(\tilde{K}'_{1,\underline{H}, v})_{v\not\in V}$. Pour $v\not\in V$, ces familles d\'eterminent des facteurs de transfert normalis\'es $\Delta_{1,H,v}$ et $\Delta_{1,\underline{H},v}$.  On d\'emontre les assertions similaires \`a (2) et (3):
  
  (4) pour $v\in Val(F)$ et  $(\delta_{1,v},\gamma_{v}), (\underline{\delta}_{1,v},\underline{\gamma}_{v})\in {\cal D}_{1,H,v}$, on a l'\'egalit\'e
 $$\boldsymbol{\Delta}_{1,\underline{H},v}(\tilde{\kappa}'_{1}(\delta_{1,v}),\tilde{\kappa}(\gamma_{v});\tilde{\kappa}'_{1}(\underline{\delta}_{1,v}),\tilde{\kappa}(\underline{\gamma}_{v}))=\boldsymbol{\Delta}_{1,H,v}(\delta_{1,v},\gamma_{v};\underline{\delta}_{1,v},\underline{\gamma}_{v});$$
 
 (5) pour $v\not\in V$ et $(\delta_{1,v},\gamma_{v})\in {\cal D}_{1,H,v}$, on a l'\'egalit\'e
 $$\Delta_{1,\underline{H},v}(\tilde{\kappa}'_{1}(\delta_{1,v}),\tilde{\kappa}(\gamma_{v}))=\Delta_{1,H,v}(\delta_{1,v},\gamma_{v}).$$
 
 Fixons un tore maximal $T'_{H}$ de $H'$ v\'erifiant l'hypoth\`ese (Hyp) et construisons des \'el\'ements $(\delta_{1},\gamma)\in \tilde{H}'_{1}({\mathbb A}_{F})\times \tilde{H}({\mathbb A}_{F})$ comme en 3.6. Alors le commutant $T'_{\underline{H}}$ de $\kappa'(T')$ dans $\underline{H}'$  v\'erifie l'hypoth\`ese (Hyp) et les \'el\'ements $(\tilde{\kappa}'(\delta_{1}),\tilde{\kappa}(\gamma))\in \underline{\tilde{H}}'_{1}({\mathbb A}_{F})\times \underline{\tilde{H}}({\mathbb A}_{F})$ v\'erifient les hypoth\`eses de ce paragraphe. De plus, on a
 
 (6) $\Delta_{H}(\delta_{1},\gamma)=\Delta_{\underline{H}}(\tilde{\kappa}'(\delta_{1}),\tilde{\kappa}(\gamma))$.
 
 La preuve est similaire \`a celle de (2).
 
 Il  r\'esulte de (4), (5) et (6) que, pour $(\delta_{1,V},\gamma_{V})\in {\cal D}_{1,H,V}$, on a l'\'egalit\'e
 $$\Delta_{1,\underline{H},V}(\tilde{\kappa}'_{1}(\delta_{1,V}),\kappa(\gamma_{V}))=\Delta_{1,H,V}(\delta_{1,V},\gamma_{V}).$$
 Le plongement $\tilde{\iota}_{\underline{H}}:\tilde{G}\to \underline{\tilde{H}}$ est le compos\'e du plongement similaire $\tilde{\iota}_{H}:\tilde{G}\to \tilde{H}$ et de $\tilde{\kappa}$. De m\^eme, l'application $\tilde{\iota}'_{\underline{H}}$ est la compos\'ee de $\tilde{\iota}'_{H}$ et de $\tilde{\kappa}'$. L'\'egalit\'e pr\'ec\'edente montre que le facteur de transfert $\Delta_{1,V}$ pour ${\bf G}'$ d\'eduit de $\Delta_{1,\underline{H},V}$ est le m\^eme que celui d\'eduit de $\Delta_{1,H,V}$. Donc l'isomorphisme
  $$  \otimes_{v\in V}C_{c}^{\infty}({\bf G}'_{v})\simeq  C_{c,\lambda_{1}}^{\infty}(\tilde{G}'_{1}(F_{V}))$$
  est inchang\'e quand on remplace les donn\'ees $H$ etc... par $\underline{H}$ etc... L'isomorphisme
 $$C_{c}^{\infty}({\bf G'}_{V})\simeq C_{c,\lambda_{1}}^{\infty}(\tilde{G}'_{1}(F_{V}))$$
 ne change pas non plus puisqu'il ne d\'epend que de la famille $(\tilde{K}'_{1,v})_{v\not\in V}$.   Cela ach\`eve la preuve. $\square$
 
 {\bf Remarque.} En 3.3, on a d\'efini l'espace $C_{c}^{\infty}({\bf G'}_{V})$ en supposant ${\bf G}'$ relevante et $V\supset V_{ram}({\bf G}')$. Supposant toujours ${\bf G}'$ relevante et soit $V$ un ensemble fini quelconque de places de $F$. On peut poser $C_{c}^{\infty}({\bf G'}_{V})= \otimes_{v\in V}C_{c}^{\infty}({\bf G}'_{v})$. Il n'y a pas d'ambigu\"{\i}t\'e puisque la proposition pr\'ec\'edente affirme que, dans le domaine commun des deux d\'efinitions (c'est-\`a-dire quand $V$ contient $V_{ram}({\bf G}')$), les deux espaces ainsi d\'efinis sont canoniquement isomorphes. 
 
  \bigskip

  \subsection{Adaptation aux $K$-espaces}
 Consid\'erons un $K$-espace comme en 1.16. Soit ${\bf G}'=(G',{\cal G}',\tilde{s})$ une donn\'ee endoscopique relevante de ce $K$-espace. On fixe un ensemble fini $V$ de places contenant $V_{ram}$ et des donn\'ees auxiliaires non ramifi\'ees hors de $V$.
 
 Soit $v$ une place de $F$. Le bifacteur $\boldsymbol{\Delta}_{1,v}$ s'\'etend \`a tout le $K$-espace, c'est-\`a-dire que l'on peut d\'efinir des termes $\boldsymbol{\Delta}_{1,v}(\delta_{1},\gamma;\underline{\delta}_{1},\underline{\gamma})$ pour deux couples $(\delta_{1},\gamma)$, $(\underline{\delta}_{1},\underline{\gamma})$ se correspondant, avec $\gamma\in \tilde{G}_{p}(F_{v})$ et $\underline{\gamma}\in \tilde{G}_{\underline{p}}(F_{v})$. La d\'efinition est la m\^eme qu'en [I] 2.3.  Avec les d\'efinitions de 1.16, on a les \'egalit\'es

$\boldsymbol{\Delta}_{1,v}(\delta_{1},\gamma;\underline{\delta}_{1},\underline{\gamma})=\boldsymbol{\Delta}_{1,v}(\delta_{1},\tilde{\phi}'_{p',p}(\gamma);\underline{\delta}_{1},\tilde{\phi}'_{p',\underline{p}}(\underline{\gamma}))$, si $v$ est finie ou complexe;

 $\boldsymbol{\Delta}_{1,v}(\delta_{1},\gamma;\underline{\delta}_{1},\underline{\gamma})=\omega(h_{p}h_{\underline{p}}^{-1})\boldsymbol{\Delta}_{1,v}(\delta_{1},\tilde{\phi}'_{p',p}(\gamma);\underline{\delta}_{1},\tilde{\phi}'_{\underline{p}',\underline{p}}(\underline{\gamma}))$ si $v$ est r\'eelle.

   Consid\'erons les deux hypoth\`eses
 
 ${\bf Hyp'}$ il existe un sous-tore maximal $T'$ de $F$, d\'efini sur $F$, de sorte que, pour toute place $v$ de $F$, il existe $p_{v}\in \Pi$ et un couple $(\delta_{v},\gamma_{v})\in {\cal D}({\bf G}'_{v})$ avec $\delta_{v}\in \tilde{T}'(F_{v})$ et $\gamma_{v}\in \tilde{G}_{p_{v}}(F_{v})$; 
 
  ${\bf Hyp''}$ il existe un sous-tore maximal $T'$ de $F$, d\'efini sur $F$, et un \'el\'ement $p\in \Pi$ de sorte que, pour toute place $v$ de $F$, il existe un couple $(\delta_{v},\gamma_{v})\in {\cal D}({\bf G}'_{v})$ avec $\delta_{v}\in \tilde{T}'(F_{v})$ et $\gamma_{v}\in \tilde{G}_{p}(F_{v})$.
  
  Elles sont \'equivalentes. En effet, la seconde implique \'evidemment la premi\`ere. Supposons ${\bf Hyp'}$ v\'erifi\'ee. Pour toute place r\'eelle, fixons $p_{v}$ v\'erifiant cette hypoth\`ese. D'apr\`es le lemme 1.16, il existe  $p\in \Pi$, d'ailleurs unique, tel que $p$ soit $v$-\'equivalent \`a $p_{v}$ pour toute place $v$ r\'eelle. On voit alors que cet $p$ v\'erifie ${\bf Hyp''}$.  
 
 Supposons ces hypoth\`eses v\'erifi\'ees.   On fixe $T'$ et $p$ v\'erifiant ${\bf Hyp''}$. Appliquant la construction  de 3.6 \`a la composante $\tilde{G}_{p}$, on construit un couple $(\delta_{1},\gamma)$ avec $\gamma\in \tilde{G}_{p}({\mathbb A}_{F})$  et un terme $\Delta(\delta_{1},\gamma)$.  Si on remplace les donn\'ees $T'$, $p$, $\delta_{1}$ et $\gamma$ par d'autres donn\'ees $\underline{T}'$, $\underline{p}$, $\underline{\delta}_{1}$ et $\underline{\gamma}$, la proposition 3.6(ii) est encore v\'erifi\'ee, les bifacteurs \'etant \'etendus comme ci-dessus au $K$-espace. La seule modification \`a faire \`a la d\'emonstration est la suivante: en 3.6, on a fix\'e $r\in G_{SC}(\bar{F})$ tel que $ad_{r}({\cal E}^*)=\underline{{\cal E}}^*$ et on a pos\'e $u_{\underline{{\cal E}}^*}(\sigma)=ru_{{\cal E}^*}(\sigma)\sigma(r)^{-1}$; il faut maintenant fixer $r\in G_{\underline{p},SC}(\bar{F})$ tel que $ad_{r}\circ \phi_{\underline{p},p}({\cal E}^*)=\underline{{\cal E}}^*$ et poser 
 $$u_{\underline{{\cal E}}^*}(\sigma)=r\phi_{\underline{p},p}(u_{{\cal E}^*}(\sigma))\nabla_{\underline{p},p}(\sigma)\sigma(r)^{-1}.$$
 
 Cette construction une fois faite, le reste des paragraphes 3.6 et 3.7 s'adapte sans changement aux $K$-espaces. 
 
 Les paragraphes 3.8 et 3.9 s'adaptent aussi de la fa\c{c}on suivante. Fixons $p\in \Pi$ et effectuons les constructions de 3.8 pour la composante $G_{p}$. On construit donc un groupe $H_{p}$ et un espace $\tilde{H}_{p}$ comme en 3.8 relatif \`a $G_{p}$ et $\tilde{G}_{p}$. Pour $q\in \Pi$, on introduit un groupe $H_{q}$ et un espace $\tilde{H}_{q}$: sur $\bar{F}$, ils sont \'egaux \`a $H_{p}$ et $\tilde{H}_{p}$; on notant les identit\'es $\phi_{p,q}^H:H_{q}\to H_{p}$ et $\tilde{\phi}_{p,q}^H:\tilde{H}_{q}\to \tilde{H}_{p}$, on d\'efinit les actions galoisiennes sur $H_{q}$ et $\tilde{H}_{q}$ de sorte que l'on ait les \'egalit\'es $\phi_{p,q}^H\circ \sigma(\phi_{p,q}^H)^{-1}=ad_{\nabla_{p,q}(\sigma)}$ et $\tilde{\phi}_{p,q}^H\circ \sigma(\tilde{\phi}^H_{p,q})^{-1}=ad_{\nabla_{p,q}(\sigma)}$ pour tout $\sigma\in \Gamma_{F}$. On a alors des plongements $G_{q}\to H_{q}$ et $\tilde{G}_{q}\to \tilde{H}_{q}$ d\'efinis sur $F$ de sorte que les diagrammes suivants soient commutatifs
 $$\begin{array}{ccccccc}G_{p}&\to& H_{p}&\,\,&\tilde{G}_{p}&\to&\tilde{H}_{p}\\ \,\,\uparrow \phi_{p,q}&&\,\,\uparrow\phi_{p,q}^H&\,\,&\,\,\uparrow\tilde{\phi}_{p,q}&&\,\,\uparrow\tilde{\phi}_{p,q}^H\\ G_{q}&\to& H_{q}&\,\,&\tilde{G}_{q}&\to&\tilde{H}_{q}\\ \end{array}$$
 La collection $(\tilde{H}_{q})_{q\in \Pi}$ n'a pas de raison d'\^etre un $K$-espace au sens de 1.16 car les applications $H^1(F;G_{q})\to H^1(F;H_{q})$ ne sont pas bijectives en g\'en\'eral. Mais les conditions pr\'ecises impos\'ees en 1.16 ne servent qu'aux formules d'inversion de [I] 4.9. Elles ne sont pas n\'ecessaires pour d\'efinir les facteurs de transfert de la m\^eme fa\c{c}on que ci-dessus.  Par abus de notations, on notera encore $K\tilde{H}$ la collection $(\tilde{H}_{q})_{q\in \Pi}$. Alors, comme en 3.9, on d\'efinit les facteurs de transfert pour cette collection puis on les restreint au $K$-espace de d\'epart $K\tilde{G}$. On obtient les m\^emes cons\'equences qu'en 3.9.

 \bigskip
 
 \section{Int\'egrales orbitales pond\'er\'ees et endoscopie}

 \bigskip
 
 \subsection{Int\'egrales orbitales pond\'er\'ees invariantes stables}
 On suppose $(G,\tilde{G},{\bf a})$ quasi-d\'eploy\'e et \`a torsion int\'erieure. Soit $\tilde{M}\in {\cal L}(\tilde{M}_{0})$ et soit $V$ un ensemble fini de places.  Posons
 $$D_{\diamond}^{st}(\tilde{M}(F_{V}))=D_{g\acute{e}om,\tilde{G}-\acute{e}qui}^{st}(\tilde{M}(F_{V}))+ D_{tr-orb}^{st}(\tilde{M}(F_{V}))\subset D_{g\acute{e}om}^{st}(\tilde{M}(F_{V})).$$
 On va d\'efinir une forme bilin\'eaire
 $$(\boldsymbol{\delta},{\bf f})\mapsto S_{\tilde{M}}^{\tilde{G}}(\boldsymbol{\delta},{\bf f})$$
 sur le produit
 $$\left( D^{st}_{\diamond}(\tilde{M}(F_{V}))\otimes Mes(M(F_{V}))^*\right)\times \left(C_{c}^{\infty}(\tilde{G}(F_{V}))\otimes Mes(G(F_{V}))\right).$$
 
 Le th\'eor\`eme qui suit \'enonce une de ses propri\'et\'es cl\'es. On le prouvera au paragraphe suivant.
 
 \ass{Th\'eor\`eme}{Pour tout $\boldsymbol{\delta}\in D^{st}_{\diamond}(\tilde{M}(F_{V}))\otimes Mes(M(F_{V}))^*$, la distribution ${\bf f}\mapsto  S_{\tilde{M}}^{\tilde{G}}(\boldsymbol{\delta},{\bf f})$ est stable, c'est-\`a-dire se factorise par la projection de $C_{c}^{\infty}(\tilde{G}(F_{V}))\otimes Mes(G(F_{V}))$ dans 
  $SI(\tilde{G}(F_{V}))\otimes Mes(G(F_{V}))$.}
 
  Conform\'ement \`a nos hypoth\`eses de r\'ecurrence, on admet cette propri\'et\'e  pour les couples $(G',\tilde{G}')$ v\'erifiant les m\^emes hypoth\`eses que $(G,\tilde{G})$ et tels que $dim(G'_{SC})< dim(G_{SC})$.  
 
 Supposons d\'efinie notre forme bilin\'eaire. On doit \'etendre la d\'efinition \`a la situation "avec caract\`ere central" de 1.15. On doit aussi montrer que, pour deux extensions $\tilde{G}_{1}$, $\tilde{G}_{2}$ de $\tilde{G}$   et pour une fonction de recollement $\tilde{\lambda}_{12,V}$ comme dans ce paragraphe, les deux d\'efinitions de la forme bilin\'eaire se recollent. Admettons cela par r\'ecurrence pour les couples $(G',\tilde{G}')$ comme ci-dessus. Soit $s\in Z(\hat{M})^{\Gamma_{F}}$. On d\'efinit la donn\'ee endoscopique ${\bf G}'(s)$ comme dans le cas local (cf. [I] 3.3), dont la donn\'ee endoscopique maximale ${\bf M}$ est une "donn\'ee de Levi". Pour $s\not\in Z(\hat{G})^{\Gamma_{F}}$, on a $dim(G'(s)_{SC})<dim(G_{SC})$ et les propri\'et\'es formelles ci-dessus permettent de d\'efinir $S_{{\bf M}}^{{\bf G}'(s)}(\boldsymbol{\delta},{\bf f})$ pour $\boldsymbol{\delta}\in D^{st}_{\diamond}({\bf M}_{V})\otimes Mes(M(F_{V}))^*$ et ${\bf f}\in SI({\bf G}'(s)_{V})\otimes Mes(G'(F_{V}))$. On d\'efinit d'autre part
 $$i_{\tilde{M}}(\tilde{G},\tilde{G}'(s))=\left\lbrace\begin{array}{cc} [Z(\hat{G}'(s))^{\Gamma_{F}}:Z(\hat{G})^{\Gamma_{F}}]^{-1},&\text{ si }{\bf G}'(s)\text{ est elliptique,}\\ 0,&\text{ sinon.}\\ \end{array}\right.$$
 
 On peut alors poser la d\'efinition
 $$S_{\tilde{M}}^{\tilde{G}}(\boldsymbol{\delta},{\bf f})=I_{\tilde{M}}^{\tilde{G}}(\boldsymbol{\delta},{\bf f})-\sum_{s\in Z(\hat{M})^{\Gamma_{F}}/Z(\hat{G})^{\Gamma_{F}},s\not=1}i_{\tilde{M}}(\tilde{G},\tilde{G}'(s))S_{{\bf M}}^{{\bf G}'(s)}(\boldsymbol{\delta},{\bf f}^{{\bf G}'(s)}).$$
 
 Il est clair que cette forme est invariante en ${\bf f}$, c'est-\`a-dire ne d\'epend que de l'image de ${\bf f}$ dans $I(\tilde{G}(F_{V}))\otimes Mes(G(F_{V}))$.
 
 \bigskip
 
 \subsection{Formules de d\'ecomposition}
 
 La situation est la m\^eme que dans le paragraphe pr\'ec\'edent. Le lien entre la forme bilin\'eaire qu'on y a d\'efinie et ses avatars locaux d\'efinis en [II] 1.10  et [V] 1.4  et 2.4 est donn\'e par la proposition suivante. Soient deux espaces de Levi $\tilde{L}_{1},\tilde{L}_{2}\in {\cal L}(\tilde{M})$ tels que $d_{\tilde{M}}^{\tilde{G}}(\tilde{L}_{1},\tilde{L}_{2})\not=0$. Alors ${\cal A}_{L_{1}}^G\cap {\cal A}_{L_{2}}^G=\{0\}$. Cette condition entra\^{\i}ne dualement que l'homomorphisme
 $$Z(\hat{M})^{\Gamma_{F}}/Z(\hat{G})^{\Gamma_{F}}\to Z(\hat{M})^{\Gamma_{F}}/Z(\hat{L}_{1})^{\Gamma_{F}}\oplus Z(\hat{M})^{\Gamma_{F}}/Z(\hat{L}_{2})^{\Gamma_{F}}$$
 est surjectif et de noyau fini. On note $k_{\tilde{M}}^{\tilde{G}}(\tilde{L}_{1},\tilde{L}_{2})$ le nombre d'\'el\'ements de ce noyau et on pose 
  $e_{\tilde{M}}^{\tilde{G}}(\tilde{L}_{1},\tilde{L}_{2})=d_{\tilde{M}}^{\tilde{G}}(\tilde{L}_{1},\tilde{L}_{2})k_{\tilde{M}}^{\tilde{G}}(\tilde{L}_{1},\tilde{L}_{2})^{-1}$.    Si au contraire $d_{\tilde{M}}^{\tilde{G}}(\tilde{L}_{1},\tilde{L}_{2})=0$, on pose $e_{\tilde{M}}^{\tilde{G}}(\tilde{L}_{1},\tilde{L}_{2})=0$. Soit $\tilde{L}^V=(\tilde{L}^v)_{v\in V}\in {\cal L}(\tilde{M}_{V})$ tel que $d_{\tilde{M}_{V}}^{\tilde{G}}(\tilde{M},\tilde{L}^V)\not=0$. Cette condition entra\^{\i}ne encore que l'homomorphisme
  $$Z(\hat{M})^{\Gamma_{F}}/Z(\hat{G})^{\Gamma_{F}}\to \oplus_{v\in V}Z(\hat{M}_{v})^{\Gamma_{F_{v}}}/Z(\hat{L}^v)^{\Gamma_{F_{v}}}$$
  est surjectif et de noyau fini. On note $k_{\tilde{M}_{V}}^{\tilde{G}}(\tilde{M},\tilde{L}^V)$ le nombre d'\'el\'ements de ce noyau et on pose 
   $e_{\tilde{M}_{V}}^{\tilde{G}}(\tilde{M},\tilde{L}^V)=d_{\tilde{M}_{V}}^{\tilde{G}}(\tilde{M},\tilde{L}^V)k_{\tilde{M}_{V}}^{\tilde{G}}(\tilde{M},\tilde{L}^V)^{-1}$. Si au contraire $d_{\tilde{M}_{V}}^{\tilde{G}}(\tilde{M},\tilde{L}^V)=0$, on pose $e_{\tilde{M}_{V}}^{\tilde{G}}(\tilde{M},\tilde{L}^V)=0$.
 
 \ass{Proposition}{Soient $\boldsymbol{\delta}=\otimes_{v\in V}\boldsymbol{\delta}_{v}\in D^{st}_{\diamond}(\tilde{M}(F_{V}))\otimes Mes(M(F_{V}))^*$ et ${\bf f}=\otimes_{v\in V}{\bf f}_{v}\in C_{c}^{\infty}(\tilde{G}(F_{V}))\otimes Mes(G(F_{V}))$. 
 
 (i) On a l'\'egalit\'e
 $$S_{\tilde{M}}^{\tilde{G}}(\boldsymbol{\delta},{\bf f})=\sum_{\tilde{L}^V\in {\cal L}(\tilde{M}_{V})}e_{\tilde{M}_{V}}^{\tilde{G}}(\tilde{M},\tilde{L}^V)\prod_{v\in V}S_{\tilde{M}_{v}}^{\tilde{L}^v}(\boldsymbol{\delta}_{v},{\bf f}_{v,\tilde{L}^v}).$$
 
 (ii) Supposons $V$ r\'eunion disjointe de $V_{1}$ et $V_{2}$. Pour $i=1,2$, posons  $\boldsymbol{\delta}_{i}=\otimes_{v\in V_{i}}\boldsymbol{\delta}_{v}$ et ${\bf f}_{i}=\otimes_{v\in V_{i}}{\bf f}_{v}$. Alors
 $$S_{\tilde{M}}^{\tilde{G}}(\boldsymbol{\delta},{\bf f})=\sum_{\tilde{L}_{1},\tilde{L}_{2}\in {\cal L}(\tilde{M})}  e_{\tilde{M}}^{\tilde{G}}(\tilde{L}_{1},\tilde{L}_{2})S_{\tilde{M}}^{\tilde{L}_{1}}(\boldsymbol{\delta}_{1},{\bf f}_{1,\tilde{L}_{1}})S_{\tilde{M}}^{\tilde{L}_{2}}(\boldsymbol{\delta}_{2},{\bf f}_{2,\tilde{L}_{2}}).$$
 
 (iii) Supposons $V$ r\'eduit \`a une place $v$. Alors
 $$S_{\tilde{M}}^{\tilde{G}}(\boldsymbol{\delta},{\bf f})=\sum_{\tilde{L}^v\in {\cal L}(\tilde{M}_{v})}e_{\tilde{M}_{v}}^{\tilde{G}}(\tilde{M},\tilde{L}^v)S_{\tilde{M}_{v}}^{\tilde{L}^v}(\boldsymbol{\delta},{\bf f}_{\tilde{L}^v}).$$}
 
 Preuve. On ne d\'emontre pas cette proposition car nous ferons une d\'emonstration analogue au paragraphe 4.4, dans une situation plus g\'en\'erale. On va simplement montrer ici que (i) est \'equivalente \`a la r\'eunion des  deux assertions  (ii) et (iii).
 
 On va commencer par l'implication (ii)+(iii) implique (i).On raisonne par r\'ecurrence sur le nombre d'\'el\'ements de $V$. Si $V$ est r\'eduit \`a une place $v$, l'assertion (i) est identique \`a (iii). Si $V$ a au moins deux \'el\'ements, on d\'ecompose $V$ en union disjointe $V_{1}\sqcup V_{2}$ de deux sous-ensembles non vides, donc de nombres d'\'el\'ements strictement inf\'erieurs \`a celui de $V$. On applique (ii). Pour $i=1,2$, on applique (i)   par r\'ecurrence aux termes $S_{\tilde{M}}^{\tilde{L}_{i}}(\boldsymbol{\delta}_{i},{\bf f}_{i,\tilde{L}_{i}})$  qui interviennent.  Il appara\^{\i}t des ensembles de sommation ${\cal L}^{\tilde{L}_{i}}(\tilde{M}_{V_{i}})\subset {\cal L}(\tilde{M}_{V_{i}})$. L'ensemble ${\cal L}(\tilde{M}_{V})$ est le produit de ${\cal L}(\tilde{M}_{V_{1}})$ et de ${\cal L}(\tilde{M}_{V_{2}})$. On pose
  $${\cal L}^{\tilde{L}_{1},\tilde{L}_{2}}(\tilde{M}_{V})={\cal L}^{\tilde{L}_{1}}(\tilde{M}_{V_{1}})\times {\cal L}^{\tilde{L}_{2}}(\tilde{M}_{V_{2}}).$$
 Un produit pour $i=1,2$ de sommes sur ${\cal L}^{\tilde{L}_{i}}(\tilde{M}_{V_{i}})$ est donc une somme sur ${\cal L}^{\tilde{L}_{1},\tilde{L}_{2}}(\tilde{M}_{V})$. D'o\`u
 $$S_{\tilde{M}}^{\tilde{L}_{1}}(\boldsymbol{\delta}_{1},{\bf f}_{1,\tilde{L}_{1}})S_{\tilde{M}}^{\tilde{L}_{2}}(\boldsymbol{\delta}_{2},{\bf f}_{2,\tilde{L}_{2}})=\sum_{\tilde{L}^V\in {\cal L}^{\tilde{L}_{1},\tilde{L}_{2}}(\tilde{M}_{V})}e_{\tilde{M}_{V_{1}}}^{\tilde{L}_{1}}(\tilde{M},\tilde{L}^{V_{1}})e_{\tilde{M}_{V_{2}}}^{\tilde{L}_{2}}(\tilde{M},\tilde{L}^{V_{2}})\prod_{v\in V}S_{\tilde{M}_{v}}^{\tilde{L}^v}(\boldsymbol{\delta}_{v},{\bf f}_{v,\tilde{L}^v})$$
 (on a not\'e $\tilde{L}^{V_{i}}$ pour $i=1,2$ les deux composantes de $\tilde{L}^{V}$). En inversant les sommes en $\tilde{L}_{1}, \tilde{L}_{2}$ et en $\tilde{L}^V$, on obtient
  $$S_{\tilde{M}}^{\tilde{G}}(\boldsymbol{\delta},{\bf f})=\sum_{\tilde{L}^V\in {\cal L}(\tilde{M}_{V})}E_{\tilde{M}_{V}}^{\tilde{G}}(\tilde{M},\tilde{L}^V)\prod_{v\in V}S_{\tilde{M}_{v}}^{\tilde{L}^v}(\boldsymbol{\delta}_{v},{\bf f}_{v,\tilde{L}^v}),$$
  o\`u
  $$E_{\tilde{M}_{V}}^{\tilde{G}}(\tilde{M},\tilde{L}^V)=\sum_{\tilde{L}_{1},\tilde{L}_{2}\in {\cal L}(\tilde{M}),\tilde{L}^V\in {\cal L}^{\tilde{L}_{1},\tilde{L}_{2}}(\tilde{M}_{V})}  e_{\tilde{M}}^{\tilde{G}}(\tilde{L}_{1},\tilde{L}_{2})e_{\tilde{M}_{V_{1}}}^{\tilde{L}_{1}}(\tilde{M},\tilde{L}^{V_{1}})e_{\tilde{M}_{V_{2}}}^{\tilde{L}_{2}}(\tilde{M},\tilde{L}^{V_{2}}).$$
Pour obtenir (i), il reste \`a prouver que, pour $\tilde{L}^V\in {\cal L}(\tilde{M}_{V})$, on a  
 $$(1) \qquad E_{\tilde{M}_{V}}^{\tilde{G}}(\tilde{M},\tilde{L}^V)=e_{\tilde{M}_{V}}^{\tilde{G}}(\tilde{M},\tilde{L}^V).$$
 On montre que
 
 (2) s'il existe $\tilde{L}_{1},\tilde{L}_{2}\in {\cal L}(\tilde{M})$ tels que $\tilde{L}^V\in {\cal L}^{\tilde{L}_{1},\tilde{L}_{2}}(\tilde{M}_{V})$ et 
 $$ d_{\tilde{M}}^{\tilde{G}}(\tilde{L}_{1},\tilde{L}_{2})d_{\tilde{M}_{V_{1}}}^{\tilde{L}_{1}}(\tilde{M},\tilde{L}^{V_{1}})d_{\tilde{M}_{V_{2}}}^{\tilde{L}_{2}}(\tilde{M},\tilde{L}^{V_{2}})\not=0,$$
  alors $d_{\tilde{M}_{V}}^{\tilde{G}}(\tilde{M},\tilde{L}^{V})\not=0$.
 
 L'hypoth\`ese $d_{\tilde{M}_{V_{1}}}^{\tilde{L}_{1}}(\tilde{M},\tilde{L}^{V_{1}})\not=0$ signifie que
 $$\Delta_{V_{1}}({\cal A}_{\tilde{M}}^{\tilde{L}_{1}})\oplus {\cal A}_{\tilde{L}^{V_{1}}}^{\tilde{L}_{1}}={\cal A}_{\tilde{M}_{V_{1}}}^{\tilde{L}_{1}},$$
 o\`u on utilise les notations de 1.4.
 En ajoutant l'espace $\oplus_{v\in V_{1}}{\cal A}_{\tilde{L}_{1}}^{\tilde{G}}$, qui est en somme directe avec les pr\'ec\'edents, on obtient
 $$\Delta_{V_{1}}({\cal A}_{\tilde{M}}^{\tilde{L}_{1}})\oplus {\cal A}_{\tilde{L}^{V_{1}}}^{\tilde{G}}={\cal A}_{\tilde{M}_{V_{1}}}^{\tilde{G}}.$$
 De m\^eme en rempla\c{c}ant les indices $1$ par $2$. En sommant les \'egalit\'es obtenues, on obtient
 $$(3) \qquad \Delta_{V_{1}}({\cal A}_{\tilde{M}}^{\tilde{L}_{1}})\oplus \Delta_{V_{2}}({\cal A}_{\tilde{M}}^{\tilde{L}_{2}})\oplus {\cal A}_{\tilde{L}^V}^{\tilde{G}}={\cal A}_{\tilde{M}_{V}}^{\tilde{G}}.$$
 Montrons que l'on a l'inclusion
 $$(4) \qquad \Delta_{V_{1}}({\cal A}_{\tilde{M}}^{\tilde{L}_{1}})\oplus \Delta_{V_{2}}({\cal A}_{\tilde{M}}^{\tilde{L}_{2}})\subset \Delta_{V}({\cal A}_{\tilde{M}}^{\tilde{G}})+( \Delta_{V_{1}}({\cal A}_{\tilde{L}_{1}}^{\tilde{G}})\oplus \Delta_{V_{2}}({\cal A}_{\tilde{L}_{2}}^{\tilde{G}})).$$
 L'hypoth\`ese $ d_{\tilde{M}}^{\tilde{G}}(\tilde{L}_{1},\tilde{L}_{2})\not=0$ entra\^{\i}ne que l'application lin\'eaire somme des projections orthogonales
 $$\begin{array}{ccc}{\cal A}_{\tilde{M}}^{\tilde{G}}&\to&{\cal A}_{\tilde{M}}^{\tilde{L}_{1}}\oplus {\cal A}_{\tilde{M}}^{\tilde{L}_{2}}\\ H&\mapsto&(H^{\tilde{L}_{1}},H^{\tilde{L}_{2}})\\ \end{array}$$
 est bijective. Pour $i=1,2$, soit $H_{i}\in {\cal A}_{\tilde{M}}^{\tilde{L}_{i}}$. On introduit l'image r\'eciproque $H$ de $(H_{1},H_{2})$ par l'application pr\'ec\'edente. Alors
 $$\Delta_{V_{1}}(H_{1})+\Delta_{V_{2}}(H_{2})=\Delta_{V}(H)-\Delta_{V_{1}}(H_{\tilde{L}_{1}})-\Delta_{V_{2}}(H_{\tilde{L}_{2}}).$$
 Le membre de droite appartient au membre de droite de (4), ce qui prouve cette assertion.

 Le dernier espace du membre de gauche de (3) contient   $ \Delta_{V_{1}}({\cal A}_{\tilde{L}_{1}}^{\tilde{G}})\oplus \Delta_{V_{2}}({\cal A}_{\tilde{L}_{2}}^{\tilde{G}})$. Gr\^ace \`a (4), on obtient l'inclusion
 $$(5) \qquad  \Delta_{V_{1}}({\cal A}_{\tilde{M}}^{\tilde{L}_{1}})\oplus \Delta_{V_{2}}({\cal A}_{\tilde{M}}^{\tilde{L}_{2}})\oplus {\cal A}_{\tilde{L}^V}^{\tilde{G}}\subset \Delta_{V}({\cal A}_{\tilde{M}}^{\tilde{G}})+{\cal A}_{\tilde{L}^V}^{\tilde{G}}.$$
 La somme des dimensions est la m\^eme des deux c\^ot\'es. En effet, $\Delta_{V}({\cal A}_{\tilde{M}}^{\tilde{G}})$ a m\^eme dimension que ${\cal A}_{\tilde{M}}^{\tilde{G}}$ et, pour $i=1,2$, $\Delta_{V_{i}}({\cal A}_{\tilde{M}}^{\tilde{L}_{i}})$ a m\^eme dimension que ${\cal A}_{\tilde{M}}^{\tilde{L}_{i}}$. Or  l'hypoth\`ese $ d_{\tilde{M}}^{\tilde{G}}(\tilde{L}_{1},\tilde{L}_{2})\not=0$ entra\^{\i}ne que ${\cal A}_{\tilde{M}}^{\tilde{G}}={\cal A}_{\tilde{M}}^{\tilde{L}_{1}}\oplus {\cal A}_{\tilde{M}}^{\tilde{L}_{2}}$. Cela implique l'\'egalit\'e voulue des dimensions. Mais celle-ci implique que les espaces du membre de droite de (5) sont en somme directe et que l'inclusion est une \'egalit\'e. En utilisant (3), on obtient l'\'egalit\'e
 $$(6) \qquad \Delta_{V}({\cal A}_{\tilde{M}}^{\tilde{G}})\oplus {\cal A}_{\tilde{L}^V}^{\tilde{G}}={\cal A}_{\tilde{M}_{V}}^{\tilde{G}}.$$
 Cela signifie que la conclusion de (2) est v\'erifi\'ee.
 
 Inversement, on a
  
  (7) si $d_{\tilde{M}_{V}}^{\tilde{G}}(\tilde{M},\tilde{L}^{V})\not=0$, alors il existe un unique couple $(\tilde{L}_{1},\tilde{L}_{2})$ v\'erifiant les hypoth\`eses de (2).
  
 Montrons d'abord l'unicit\'e. La condition $\tilde{L}^{V_{1}}\in {\cal L}^{\tilde{L}_{1}}(\tilde{M}_{V_{1}})$ signifie que ${\cal A}_{\tilde{L}_{1}}\subset {\cal A}_{\tilde{L}_{v}}$ pour tout $v\in V_{1}$. La condition  $d_{\tilde{M}_{V_{1}}}^{\tilde{L}_{1}}(\tilde{M},\tilde{L}^{V_{1}})\not=0$ entra\^{\i}ne  
 $${\cal A}_{\tilde{M}}^{\tilde{L}_{1}}\cap(\cap_{v\in V_{1}}{\cal A}_{\tilde{L}_{v}}^{\tilde{L}_{1}})=\{0\}.$$
 Leur conjonction conduit \`a l'\'egalit\'e
 $$(8) \qquad {\cal A}_{\tilde{L}_{1}}={\cal A}_{\tilde{M}}\cap(\cap_{v\in V_{1}}{\cal A}_{\tilde{L}_{v}}).$$
 D'o\`u l'unicit\'e de $\tilde{L}_{1}$ et de m\^eme celle de $\tilde{L}_{2}$. Inversement, notons ${\cal B}$ l'intersection de droite ci-dessus, soit $T$ le  sous-tore d\'efini sur $F$ de $A_{\tilde{M}}$  tel que $X_{*}(T)=X_{*}(A_{\tilde{M}})\cap {\cal B}$.  Soit $\tilde{L}_{1}$ le commutant de $T$ dans $\tilde{G}$. Il contient $\tilde{M}$ et les $\tilde{L}_{v}$ pour $v\in V_{1}$. A fortiori, il est non vide et c'est donc un espace de Levi. Les inclusions ci-dessus entra\^{\i}nent que ${\cal A}_{\tilde{L}_{1}}$ est contenu dans ${\cal B}$. Mais il contient par d\'efinition $X_{*}(T)$ qui engendre ${\cal B}$. Il est donc \'egal \`a ${\cal B}$.  Cela prouve l'existence d'un espace de Levi $\tilde{L}_{1}$ v\'erifiant (8). On d\'efinit de m\^eme l'espace $\tilde{L}_{2}$. On va montrer que le couple $(\tilde{L}_{1},\tilde{L}_{2})$ v\'erifie les hypoth\`eses de (2). La d\'efinition de ces espaces n'entra\^{\i}ne pas l'\'egalit\'e (3), mais elle entra\^{\i}ne que les espaces du membre de gauche de cette \'egalit\'e sont en somme directe. L'hypoth\`ese $d_{\tilde{M}_{V}}^{\tilde{G}}(\tilde{M},\tilde{L}^{V})\not=0$ entra\^{\i}ne l'\'egalit\'e (6). On d\'eduit de ces deux faits l'in\'egalit\'e
  $$dim({\cal A}_{\tilde{M}}^{\tilde{G}})\geq dim({\cal A}_{\tilde{M}}^{\tilde{L}_{1}})+ dim({\cal A}_{\tilde{M}}^{\tilde{L}_{2}}).$$
  L'\'egalit\'e (8) et son analogue pour $\tilde{L}_{2}$ entra\^{\i}nent que 
  $${\cal A}_{\tilde{L}_{1}}\cap {\cal A}_{\tilde{L}_{2}}={\cal A}_{\tilde{M}}\cap(\cap_{v\in V}{\cal A}_{\tilde{L}_{v}}).$$
 Mais cette intersection est r\'eduite \`a ${\cal A}_{\tilde{G}}$ en vertu de l'hypoth\`ese $d_{\tilde{M}_{V}}^{\tilde{G}}(\tilde{M},\tilde{L}^{V})\not=0$. Donc les espaces ${\cal A}_{\tilde{L}_{1}}^{\tilde{G}}$ et ${\cal A}_{\tilde{L}_{2}}^{\tilde{G}}$ sont en somme directe. En prenant les orthogonaux, on obtient l'\'egalit\'e
 $$ {\cal A}_{\tilde{M}}^{\tilde{G}}={\cal A}_{\tilde{M}}^{\tilde{L}_{1}}+ {\cal A}_{\tilde{M}}^{\tilde{L}_{2}}.$$
 L'in\'egalit\'e de dimensions prouv\'ee ci-dessus entra\^{\i}ne que la somme de droite est directe. D'o\`u $d_{\tilde{M}}^{\tilde{G}}(\tilde{L}_{1},\tilde{L}_{2})\not=0$. D'o\`u aussi: cette in\'egalit\'e de dimensions est une \'egalit\'e. En reprenant le raisonnement conduisant \`a cette in\'egalit\'e, on voit que l'\'egalit\'e (3)  est v\'erifi\'ee. En la projetant sur $\oplus_{v\in V_{i}}{\cal A}_{\tilde{M}}^{\tilde{L}_{i}}$, on obtient 
  $$\Delta_{V_{i}}({\cal A}_{\tilde{M}}^{\tilde{L}_{i}})\oplus {\cal A}_{\tilde{L}^{V_{i}}}^{\tilde{L}_{i}}={\cal A}_{\tilde{M}_{V_{i}}}^{\tilde{L}_{i}}.$$
 Cela signifie que $d_{\tilde{M}_{V_{i}}}^{\tilde{L}_{i}}(\tilde{M},\tilde{L}^{V_{i}})\not=0$ pour $i=1,2$. Cela prouve (7).

 Gr\^ace \`a (2) et (7), on obtient
 $$E_{\tilde{M}_{V}}^{\tilde{G}}(\tilde{M},\tilde{L}^V)=\left\lbrace\begin{array}{cc}0,&\text{ si }d_{\tilde{M}_{V}}^{\tilde{G}}(\tilde{M},\tilde{L}^{V})=0,\\ e_{\tilde{M}}^{\tilde{G}}(\tilde{L}_{1},\tilde{L}_{2})e_{\tilde{M}_{V_{1}}}^{\tilde{L}_{1}}(\tilde{M},\tilde{L}^{V_{1}})e_{\tilde{M}_{V_{2}}}^{\tilde{L}_{2}}(\tilde{M},\tilde{L}^{V_{2}}),&\text{ si }d_{\tilde{M}_{V}}^{\tilde{G}}(\tilde{M},\tilde{L}^{V})\not=0,\\ \end{array}\right.$$
 o\`u $(\tilde{L}_{1},\tilde{L}_{2})$ est le couple d\'etermin\'e par (7). Pour prouver (1), on peut supposer $d_{\tilde{M}_{V}}^{\tilde{G}}(\tilde{M},\tilde{L}^{V})\not=0$ et il faut prouver l'\'egalit\'e
 $$e_{\tilde{M}_{V}}^{\tilde{G}}(\tilde{M},\tilde{L}^V)=e_{\tilde{M}}^{\tilde{G}}(\tilde{L}_{1},\tilde{L}_{2})e_{\tilde{M}_{V_{1}}}^{\tilde{L}_{1}}(\tilde{M},\tilde{L}^{V_{1}})e_{\tilde{M}_{V_{2}}}^{\tilde{L}_{2}}(\tilde{M},\tilde{L}^{V_{2}}).$$
 Elle se d\'ecompose en les deux \'egalit\'es
  $$d_{\tilde{M}_{V}}^{\tilde{G}}(\tilde{M},\tilde{L}^V)=d_{\tilde{M}}^{\tilde{G}}(\tilde{L}_{1},\tilde{L}_{2})d_{\tilde{M}_{V_{1}}}^{\tilde{L}_{1}}(\tilde{M},\tilde{L}^{V_{1}})d_{\tilde{M}_{V_{2}}}^{\tilde{L}_{2}}(\tilde{M},\tilde{L}^{V_{2}}),$$
  et
$$(9) \qquad k_{\tilde{M}_{V}}^{\tilde{G}}(\tilde{M},\tilde{L}^V)=k_{\tilde{M}}^{\tilde{G}}(\tilde{L}_{1},\tilde{L}_{2})k_{\tilde{M}_{V_{1}}}^{\tilde{L}_{1}}(\tilde{M},\tilde{L}^{V_{1}})k_{\tilde{M}_{V_{2}}}^{\tilde{L}_{2}}(\tilde{M},\tilde{L}^{V_{2}}).$$
La premi\`ere se prouve en reprenant les d\'emonstrations de (2) et (7) et en pr\'ecisant comment se comportent les mesures selon les diff\'erentes d\'ecompositions en somme directe. On laisse les d\'etails fastidieux au lecteur. Prouvons (9). L'homomorphisme
$$\frac{Z(\hat{M})^{\Gamma_{F}}}{Z(\hat{G})^{\Gamma_{F}}}\to \oplus_{v\in V}\frac{Z(\hat{M}_{v})^{\Gamma_{F_{v}}}}{Z(\hat{L}^v)^{\Gamma_{F_{v}}}}$$
se d\'ecompose en le produit de 
$$\frac{Z(\hat{M})^{\Gamma_{F}}}{Z(\hat{G})^{\Gamma_{F}}}\to\frac{Z(\hat{M})^{\Gamma_{F}}}{Z(\hat{L}_{1})^{\Gamma_{F}}}\oplus \frac{Z(\hat{M})^{\Gamma_{F}}}{Z(\hat{L}_{2})^{\Gamma_{F}}}$$
et de
$$\frac{Z(\hat{M})^{\Gamma_{F}}}{Z(\hat{L}_{1})^{\Gamma_{F}}}\oplus \frac{Z(\hat{M})^{\Gamma_{F}}}{Z(\hat{L}_{2})^{\Gamma_{F}}}\to 
\left( \oplus_{v\in V_{1}}\frac{Z(\hat{M}_{v})^{\Gamma_{F_{v}}}}{Z(\hat{L}^v)^{\Gamma_{F_{v}}}}\right)\oplus \left( \oplus_{v\in V_{2}}\frac{Z(\hat{M}_{v})^{\Gamma_{F_{v}}}}{Z(\hat{L}^v)^{\Gamma_{F_{v}}}}\right).$$
Tous ces homomorphismes sont surjectifs. Donc le nombre d'\'el\'ements $k_{\tilde{M}_{V}}^{\tilde{G}}(\tilde{M},\tilde{L}^V)$  du noyau du compos\'e  est le produit des nombres d'\'el\'ements des noyaux des deux homomorphismes ci-dessus. Ceux-ci sont respectivement $k_{\tilde{M}}^{\tilde{G}}(\tilde{L}_{1},\tilde{L}_{2})$ et $k_{\tilde{M}_{V_{1}}}^{\tilde{L}_{1}}(\tilde{M},\tilde{L}^{V_{1}})k_{\tilde{M}_{V_{2}}}^{\tilde{L}_{2}}(\tilde{M},\tilde{L}^{V_{2}})$. Cela prouve (9) et ach\`eve la preuve de l'implication (ii)+(iii) implique (i).
 
  En fait, on a prouv\'e que, si on admettait (i) pour les facteurs du membre de droite de (ii), alors ce membre de droite \'etait \'egal \`a celui de (i). Mais cela d\'emontre que (i) implique (ii). De plus, (iii) n'est que (i) dans le cas particulier o\`u $V$ n'a qu'un \'el\'ement. Donc (i) implique (ii)+(iii). $\square$
 
Le (i) de cette proposition ram\`ene la preuve des propri\'et\'es requises de la forme bilin\'eaire $(\boldsymbol{\delta},{\bf f})\mapsto S_{\tilde{M}}^{\tilde{G}}(\boldsymbol{\delta},{\bf f})$ \`a celle des m\^emes propri\'et\'es pour ses avatars locaux. En ce qui concerne les propri\'et\'es formelles, on a fait cette preuve en [II] 1.10.  La propri\'et\'e de stabilit\'e a \'et\'e prouv\'ee en [III] 2.8 dans le cas non-archim\'edien, en [V] 1.5 et section 4 dans le cas archim\'edien. Cela prouve le th\'eor\`eme 4.1.

{\bf Variante.} Supposons donn\'e un syst\`eme de fonctions $B$ comme en 1.10. En rempla\c{c}ant les int\'egrales $I_{\tilde{M}}^{\tilde{G}}(\boldsymbol{\gamma},{\bf f})$ par leurs variantes $I_{\tilde{M}}^{\tilde{G}}(\boldsymbol{\gamma},B,{\bf f})$ dans les constructions pr\'ec\'edentes, on d\'efinit les variantes  $S_{\tilde{M}}^{\tilde{G}}(\boldsymbol{\delta},B,{\bf f})$ des int\'egrales orbitales pond\'er\'ees stables. Elles v\'erifient des propri\'et\'es analogues aux pr\'ec\'edentes.

\bigskip

\subsection{Une propri\'et\'e de support}
Les hypoth\`eses sont les m\^emes que dans le paragraphe pr\'ec\'edent mais on suppose que $V$ contient $V_{ram}$.

\ass{Lemme}{Soit $\Xi\subset {\cal A}_{\tilde{M}}$ un ensemble compact et soit ${\bf f}\in C_{c}^{\infty}(\tilde{G}(F_{V}))\otimes Mes(G(F_{V}))$. Alors il existe un sous-ensemble compact $\tilde{C}_{V}$ de $\tilde{M}(F_{V})$ tel que, pour tout $\boldsymbol{\delta}\in D^{st}_{\diamond}(\tilde{M}(F_{V}))\otimes Mes(M(F_{V}))^*$ v\'erifiant les deux conditions:

- l'image par $\tilde{H}_{\tilde{M}_{V}}$ du support de $\boldsymbol{\delta}$ est contenu dans $ \Xi$,

- $S_{\tilde{M}}^{\tilde{G}}(\boldsymbol{\delta},{\bf f})\not=0$,

\noindent il existe un \'el\'ement du support de $\boldsymbol{\delta}$ qui soit conjugu\'e \`a un \'el\'ement de $\tilde{C}_{V}$ par un \'el\'ement de $M(F_{V})$.}

Preuve. On utilise la d\'efinition. Pour que $S_{\tilde{M}}^{\tilde{G}}(\boldsymbol{\delta},{\bf f}) $ soit non nul, il faut que $I_{\tilde{M}}^{\tilde{G}}(\boldsymbol{\delta},{\bf f})$ soit non nul ou qu'il existe $\tilde{s}\in Z(\hat{M})^{\Gamma_{F}}/Z(\hat{G})^{\Gamma_{F}}$ avec $\tilde{s}\not=1$ de sorte que $S_{{\bf M}}^{{\bf G}'(\tilde{s})}(\boldsymbol{\delta},{\bf f}^{{\bf G}'(\tilde{s})})$ soit non nul. Dans le premier cas, on conclut par le lemme 1.12. Dans le deuxi\`eme, l'application $\tilde{H}_{\tilde{M}_{V}}$ ne d\'ependant pas des espaces ambiants $\tilde{G}$ ou $\tilde{G}'(\tilde{s})$, les hypoth\`eses restent v\'erifi\'ees si l'on remplace $\tilde{G}$ par $\tilde{G}'(\tilde{s})$. Puisque $\tilde{s}\not=1$, on peut appliquer le lemme par r\'ecurrence, ce qui conduit encore \`a la conclusion. $\square$

\bigskip

\subsection{Le syst\`eme de fonctions $B^{\tilde{G}}$}
Revenons au cas g\'en\'eral. Soit ${\bf G}'=(G',{\cal G}',\tilde{s})$ une donn\'ee endoscopique de $(G,\tilde{G},{\bf a})$. Nous allons munir $\tilde{G}'$ d'un syst\`eme de fonctions comme en 1.10 que nous noterons $B^{\tilde{G}}$. On fixe comme toujours une paire de Borel \'epingl\'ee $\hat{{\cal E}}$ de $\hat{G}$ pour laquelle on utilise les notations usuelles, cf. [I] 1.5. En particulier, on suppose $\tilde{s}=s\hat{\theta}$, avec $s\in \hat{T}$. 

 Fixons des paires de Borel \'epingl\'ees ${\cal E}=(B,T(,E_{\alpha})_{\alpha\in \Delta})$ de $G$ et ${\cal E}'=(B',T',(E'_{\alpha'})_{\alpha'\in \Delta'})$ de $G'$.   On suppose ${\cal E}'$ d\'efinie sur $F$. On note $\tilde{T}$ l'ensemble des $\gamma\in \tilde{G}$ tels que $ad_{\gamma}$ conserve $(B,T)$ et $\tilde{T}'$ l'ensemble des $\delta\in \tilde{G}'$ tels que $ad_{\delta}$ conserve $(B',T')$.  On a un homomorphisme $\xi:T\to T'$, qui se prolonge en une application $\tilde{\xi}:\tilde{T}\to \tilde{T}'$. Soit $\epsilon\in \tilde{T}'$, fixons $\eta\in \tilde{T}$ tel que $\tilde{\xi}(\eta)=\epsilon$, \'ecrivons $\eta=\nu e$, avec $\nu\in T$ et $e\in Z(\tilde{G},{\cal E})$.   Notons $\Sigma(T)$ l'ensemble des racines de $T$ dans $\mathfrak{g}$ et $\Sigma^{G'}_{\epsilon}(T')$ celui des racines de $T'$ dans $\mathfrak{g}'_{\epsilon}$. On a d\'ecrit maintes fois ce dernier ensemble.   C'est la r\'eunion des ensembles

(a) les $N\alpha$ pour $\alpha\in \Sigma(T)$ de type 1 tels que $N\alpha(\nu)=1$ et $N\hat{\alpha}(s)=1$; 

(b) les $2N\alpha$ pour $\alpha\in \Sigma(T)$ de type 2 tels que $N\alpha(\nu)=1$ et $N\hat{\alpha}(s)=1$; 

(c) les $2N\alpha$ pour $\alpha\in \Sigma(T)$ de type 2 tels que $N\alpha(\nu)=-1$ et $N\hat{\alpha}(s)=1$; 

(d) les $N\alpha$ pour $\alpha\in \Sigma(T)$ de type 3 tels que $N\alpha(\nu)=1$ et $N\hat{\alpha}(s)=-1$.

On a introduit en 1.10  une d\'ecomposition $\tilde{T}'=\sqcup_{\Omega\in \underline{\Omega}}\Omega$. Soit $\Omega\in \underline{\Omega}$ tel que $\epsilon\in \Omega$.  Soit $\epsilon'$ un autre \'el\'ement de $\Omega$, que l'on rel\`eve en $\eta'=\nu'e\in\tilde{T}$.
Les ensembles $\Sigma^{G'_{\epsilon}}(T')$ et $\Sigma^{G'_{\epsilon'}}(T')$ sont \'egaux par d\'efinition de $\underline{\Omega}$. Une racine $N\alpha$ avec $\alpha$ de type 1 ne saurait \^etre \'egale \`a une racine $2N\beta$ avec $\beta$ de type 2, ni \`a $N\beta$ avec $\beta$ du type 3. Donc les racines de type (a) pour $\epsilon$ sont aussi de type (a) pour $\epsilon'$.  De m\^eme, les racines de type (d) pour $\epsilon$ sont aussi de type (d) pour $\epsilon'$. Par contre, une racine de type (b), resp. (c), pour $\epsilon$ pourrait \^etre de type (b) ou (c) pour $\epsilon'$. Mais, en tout cas, pour une racine $2N\alpha$ de type (b) ou (c) pour $\epsilon$, on a forc\'ement $N\alpha(\nu')=\pm 1$. Or $\Omega$ est connexe par d\'efinition. Son image r\'eciproque dans $\tilde{T}$ l'est aussi. Cela entra\^{\i}ne que $N\alpha(\nu')$ est constant quand $\epsilon'$ parcourt $\Omega$. Alors les racines de type (b), resp. (c), pour $\epsilon$ sont aussi de type (b), resp. (c), pour tout $\epsilon'\in \Omega$.  
On d\'efinit alors une fonction $B^{\tilde{G}}_{\Omega}$ sur $ \Sigma(\Omega)=\Sigma^{G'_{\epsilon}}(T')$ de la fa\c{c}on suivante.  Dans le cas (a), $B^{\tilde{G}}_{\Omega}(N\alpha)=n_{\alpha}$; dans le cas (b), $B^{\tilde{G}}_{\Omega}(2N\alpha)=2n_{\alpha}$; dans le cas (c), $B^{\tilde{G}}_{\Omega}(2N\alpha)=n_{\alpha}$; dans le cas (d), $B^{\tilde{G}}_{\Omega}(N\alpha)=2n_{\alpha}$. La m\^eme preuve qu'en [II] 1.11 montre que la famille de fonctions $(B^{\tilde{G}}_{\Omega})_{\Omega\in \underline{\Omega}}$ ainsi d\'efinie  v\'erifie les conditions de 1.10. 

La m\^eme construction vaut si l'on travaille avec un $K$-triplet $(KG,K\tilde{G},{\bf a})$ puisque seule intervient $\underline{la}$ paire de Borel \'epingl\'ee associ\'ee \`a ce $K$-triplet.

\bigskip

\subsection{Int\'egrales orbitales pond\'er\'ees $\omega$-\'equivariantes endoscopiques} 
 Commen\c{c}ons par quelques rappels locaux. Pour quelques instants, supposons que $F$ soit un corps local non-archim\'edien de caract\'eristique nulle. Consid\'erons un triplet $(G,\tilde{G},{\bf a})$ d\'efini sur $F$, un espace de Levi $\tilde{M}$ de $\tilde{G}$ et une donn\'ee endoscopique ${\bf M}'$ de $(M,\tilde{M},{\bf a}_{M})$ elliptique et relevante.  Pour $\boldsymbol{\delta}\in D_{g\acute{e}om}^{st}({\bf M}')\otimes Mes(M'(F))^*$ et ${\bf f}\in I(\tilde{G}(F),\omega)\otimes Mes(G(F))$, on a d\'efini en [II] 1.12 un terme $I_{\tilde{M}}^{\tilde{G},{\cal E}}({\bf M}',\boldsymbol{\delta},{\bf f})$. On en a d\'eduit en [II] 1.15 une forme bilin\'eaire 
 $$(\boldsymbol{\gamma},{\bf f})\mapsto I_{\tilde{M}}^{\tilde{G},{\cal E}}(\boldsymbol{\gamma},{\bf f})$$
 sur
 $$(D_{g\acute{e}om}(\tilde{M}(F),\omega)\otimes Mes(M(F))^*)\times( I(\tilde{G}(F),\omega)\otimes Mes(G(F))).$$
 Le proc\'ed\'e consistait \`a \'ecrire $\boldsymbol{\gamma}$ comme somme $\sum_{i=1,...,n}transfert(\boldsymbol{\delta}_{i})$, o\`u $\boldsymbol{\delta}_{i}$ est une distribution g\'eom\'etrique stable  dans une donn\'ee endoscopiques ${\bf M}'_{i}$, et \`a poser
$$I_{\tilde{M}}^{\tilde{G},{\cal E}}(\boldsymbol{\gamma},{\bf f})=\sum_{i=1,...,n}I_{\tilde{M}}^{\tilde{G},{\cal E}}({\bf M}'_{i},\boldsymbol{\delta}_{i},{\bf f}).$$
Pour la suite, on doit g\'en\'eraliser la d\'efinition de $I_{\tilde{M}}^{\tilde{G},{\cal E}}({\bf M}',\boldsymbol{\delta},{\bf f})$ au cas o\`u ${\bf M}'$ est une donn\'ee endoscopique de $(M,\tilde{M},{\bf a})$ qui est relevante mais pas forc\'ement elliptique. Dans ce cas, il correspond \`a $\tilde{M}'$ un Levi $\tilde{R}$ de $\tilde{M}$ et ${\bf M}'$ est une donn\'ee endoscopique elliptique et relevante de $(R,\tilde{R},{\bf a})$. On d\'efinit
$$I_{\tilde{M}}^{\tilde{G},{\cal E}}({\bf M}',\boldsymbol{\delta},{\bf f})=\sum_{\tilde{L}\in {\cal L}(\tilde{R})}d_{\tilde{R}}^{\tilde{G}}(\tilde{M},\tilde{L})I_{\tilde{R}}^{\tilde{L},{\cal E}}({\bf M}',\boldsymbol{\delta},{\bf f}_{\tilde{L},\omega}).$$
En vertu de la relation [II] 1.15(1), on a encore l'\'egalit\'e
$$(1) \qquad I_{\tilde{M}}^{\tilde{G},{\cal E}}(transfert(\boldsymbol{\delta}),{\bf f})=I_{\tilde{M}}^{\tilde{G},{\cal E}}({\bf M}',\boldsymbol{\delta},{\bf f}).$$
Supposons maintenant que $F={\mathbb R}$ ou $F={\mathbb C}$. Des constructions analogues valent d'apr\`es [V] 1.7, 1.8 et 2.4. Il y a quelques modifications. On doit travailler avec un $K$-espace $(KG,K\tilde{G},{\bf a})$ et un $K$-espace de Levi $K\tilde{M}$. Pour une donn\'ee endoscopique ${\bf M}'$ de $(KM,K\tilde{M},{\bf a}_{M})$ elliptique et relevante, on a d\'efini en [V] 1.7 l'espace 
$D_{g\acute{e}om,\tilde{G}-\acute{e}qui}^{st}({\bf M}')$. On note $D_{\diamond}^{st}({\bf M}')$ la somme de cet espace et de $D^{st}_{tr-orb}({\bf M}')$. Les termes $ I_{K\tilde{M}}^{K\tilde{G},{\cal E}}({\bf M}',\boldsymbol{\delta},{\bf f})$ sont d\'efinis pour $\boldsymbol{\delta}\in D_{\diamond}^{st}({\bf M}')\otimes Mes(M'(F))^*$. Comme ci-dessus, les constructions se g\'en\'eralisent au cas o\`u ${\bf M}'$ est relevante mais pas elliptique. On a une \'egalit\'e analogue \`a (1). Les termes $I_{K\tilde{M}}^{K\tilde{G},{\cal E}}(\boldsymbol{\gamma},{\bf f})$ sont d\'efinis pour $\boldsymbol{\gamma}\in D_{g\acute{e}om,\tilde{G}-\acute{e}qui}(\tilde{M}(F),\omega)\otimes Mes(M(F))^*$.

 Revenons \`a notre corps de nombres $F$ et consid\'erons un $K$-triplet $(KG,K\tilde{G},{\bf a})$  comme en 1.16.  
 Soient $V$ un ensemble fini de places de $F$, $K\tilde{M} $ un \'el\'ement de ${\cal L}(K\tilde{M}_{0})$  et ${\bf M}'=(M',{\cal M}',\tilde{\zeta})$ une donn\'ee endoscopique elliptique et relevante de $(KM,K\tilde{M},{\bf a}_{M})$. Pour $\tilde{s}\in \tilde{\zeta}Z(\hat{M})^{\Gamma_{F},\hat{\theta}}$, on d\'efinit comme dans le cas local la donn\'ee endoscopique ${\bf G}'(\tilde{s})$ qui contient ${\bf M}'$ comme donn\'ee de Levi. On pose $i_{\tilde{M}'}(\tilde{G},\tilde{G}'(\tilde{s}))=0$ si ${\bf G}'(\tilde{s})$ n'est pas elliptique (on utilise ici et dans la suite la notation symbolique $\tilde{G}$ au lieu de $K\tilde{G}$ chaque fois que cela peut se faire sans ambigu\"{\i}t\'e). Si cette donn\'ee est elliptique, on pose
$$i_{\tilde{M}'}(\tilde{G},\tilde{G}'(\tilde{s}))=[Z(\hat{M}')^{\Gamma_{F}}:(Z(\hat{M}')^{\Gamma_{F}}\cap Z(\hat{M}))][Z(\hat{G}'(\tilde{s}))^{\Gamma_{F}}:(Z(\hat{G}'(\tilde{s}))^{\Gamma_{F}}\cap Z(\hat{G}))]^{-1}.$$
Comme en [II] 1.12, il y a un homomorphisme naturel
$$Z(\hat{M})^{\Gamma_{F},\hat{\theta}}/Z(\hat{G})^{\Gamma_{F},\hat{\theta}}\to Z(\hat{M}')^{\Gamma_{F}}/Z(\hat{G}'(\tilde{s}))^{\Gamma_{F}}.$$
Il est surjectif et de noyau fini. Alors $i_{\tilde{M}'}(\tilde{G},\tilde{G}'(\tilde{s}))$ est l'inverse du nombre d'\'el\'ements de ce noyau.

 On a d\'efini en [V] 1.7 et ci-dessus les espaces  $ D^{st}_{g\acute{e}om,\tilde{G}-\acute{e}qui}({\bf M}'_{v})$ et  $D^{st}_{\diamond}({\bf M}'_{v})$ pour une place $v$ archim\'edienne.   Il s'en d\'eduit comme en 1.8 des espaces $ D^{st}_{g\acute{e}om,\tilde{G}-\acute{e}qui}({\bf M}'_{V})$ et $D^{st}_{\diamond}({\bf M}'_{v})$. Pour $\boldsymbol{\delta}\in D^{st}_{\diamond}({\bf M}'_{V})\otimes Mes(M'(F_{V}))^*$ et ${\bf f}\in I(K\tilde{G}(F_{V}),\omega)\otimes Mes(G(F_{V}))$, on peut alors d\'efinir
$$(2) \qquad I_{K\tilde{M}}^{K\tilde{G},{\cal E}}({\bf M}',\boldsymbol{\delta},{\bf f})=\sum_{\tilde{s}\in \tilde{\zeta}Z(\hat{M})^{\Gamma_{F},\hat{\theta}}/Z(\hat{G})^{\Gamma_{F},\hat{\theta}}}i_{\tilde{M}'}(\tilde{G},\tilde{G}'(\tilde{s}))S_{{\bf M}'}^{{\bf G}'(\tilde{s})}(\boldsymbol{\delta}, B^{\tilde{G}},{\bf f}^{{\bf G}'(\tilde{s})}).$$
Pour d\'efinir des termes $ I_{K\tilde{M}}^{K\tilde{G},{\cal E}}(\boldsymbol{\gamma},{\bf f})$, on ne peut pas appliquer le m\^eme proc\'ed\'e  que dans le cas local car, dans le cas global, il n'y a pas en g\'en\'eral suffisamment de donn\'ees endoscopiques d\'efinies sur $F$  de $(KM,K\tilde{M},{\bf a}_{M})$ pour \'ecrire $\boldsymbol{\gamma}$ comme somme de transfert \`a partir de telles donn\'ees. Mais  soient $\boldsymbol{\gamma}=\otimes \boldsymbol{\gamma}_{v}\in D_{g\acute{e}om,\tilde{G}-\acute{e}qui}(K\tilde{M}(F_{V}),\omega)\otimes Mes(F_{V})^*$ et ${\bf f}=\otimes {\bf f}_{v}\in I(K\tilde{G}(F_{V}),\omega)\otimes Mes(G(F_{V}))$. Puisqu'on a d\'ej\`a d\'efini les formes bilin\'eaires locales, on peut poser
$$(3) \qquad I_{K\tilde{M}}^{K\tilde{G},{\cal E}}(\boldsymbol{\gamma},{\bf f})=\sum_{K\tilde{L}^V\in {\cal L}(K\tilde{M}_{V})}d_{\tilde{M}_{V}}^{\tilde{G}}(\tilde{M},\tilde{L}^V)\prod_{v\in V}I_{K\tilde{M}_{v}}^{K\tilde{L}^v,{\cal E}}(\boldsymbol{\gamma}_{v},{\bf f}_{v,\tilde{L}^v,\omega}).$$
Cette d\'efinition se prolonge par multilin\'earit\'e \`a tout $\boldsymbol{\gamma}$ et tout ${\bf f}$.

\ass{Proposition}{Soient ${\bf M}'$ une donn\'ee endoscopique elliptique et relevante de $(KM,K\tilde{M},{\bf a}_{M})$   et ${\bf f}\in I(K\tilde{G}(F_{V}),\omega)\otimes Mes(G(F_{V}))$.

(i)  Soit $\boldsymbol{\delta}=\otimes_{v\in V} \boldsymbol{\delta}_{v}\in D^{st}_{\diamond}({\bf M}'_{V})\otimes Mes(M'(F_{V}))^*$. On a l'\'egalit\'e
 $$I_{K\tilde{M}}^{K\tilde{G},{\cal E}}({\bf M}',\boldsymbol{\delta},{\bf f})=\sum_{\tilde{L}^V\in {\cal L}(\tilde{M}_{V})}d_{\tilde{M}_{V}}^{\tilde{G}}(\tilde{M},\tilde{L}^V)\prod_{v\in V}I_{K\tilde{M}_{v}}^{K\tilde{L}^v,{\cal E}}({\bf M}'_{v},\boldsymbol{\delta}_{v},{\bf f}_{v,K\tilde{L}^v,\omega}).$$
 
 (ii)  Soit $\boldsymbol{\delta} \in D^{st}_{g\acute{e}om,\tilde{G}-\acute{e}qui}({\bf M}'_{V})\otimes Mes(M'(F_{V}))^*$. On a l'\'egalit\'e
$$I_{K\tilde{M}}^{K\tilde{G},{\cal E}}(transfert(\boldsymbol{\delta}),{\bf f})=I_{K\tilde{M}}^{K\tilde{G},{\cal E}}({\bf M}',\boldsymbol{\delta},{\bf f}).$$}

Preuve. Dans le cas local, la relation (ii) est vraie d'apr\`es (1). Dans le cas pr\'esent, cette relation r\'esulte donc de (i) et de la d\'efinition (3). On pourrait prouver (i) directement. Mais, pour \^etre plus clair, on la d\'ecompose en deux assertions:

(4) si $V=V_{1}\sqcup V_{2}$, on a l'\'egalit\'e
$$I_{K\tilde{M}}^{K\tilde{G},{\cal E}}({\bf M}',\boldsymbol{\delta},{\bf f})=\sum_{K\tilde{L}_{1},K\tilde{L}_{2}\in {\cal L}(K\tilde{M})}d_{\tilde{M}}^{\tilde{G}}(\tilde{L}_{1},\tilde{L}_{2})I_{K\tilde{M}}^{K\tilde{L}_{1},{\cal E}}({\bf M}',\boldsymbol{\delta}_{1},{\bf f}_{1,K\tilde{L}_{1},\omega})I_{K\tilde{M}}^{K\tilde{L}_{2},{\cal E}}({\bf M}',\boldsymbol{\delta}_{2},{\bf f}_{2,K\tilde{L}_{2},\omega})$$
(avec des notations \'evidentes);

(5) si $V$ est r\'eduit \`a une place $v$, on a l'\'egalit\'e
$$I_{K\tilde{M}}^{K\tilde{G},{\cal E}}({\bf M}',\boldsymbol{\delta},{\bf f})=\sum_{K\tilde{L}^v\in {\cal L}(K\tilde{M}_{v})}d_{\tilde{M}_{v}}^{\tilde{G}}(\tilde{M},\tilde{L}^v)I_{K\tilde{M}_{v}}^{K\tilde{L}^v,{\cal E}}({\bf M}'_{v},\boldsymbol{\delta}_{v},{\bf f}_{v,K\tilde{L}^v,\omega}).$$

Comme dans la d\'emonstration de 4.1, la r\'eunion de ces deux assertions \'equivaut \`a (i).  

Prouvons (4). Dans le membre de droite de la d\'efinition (2), on applique la formule (ii) de la proposition 4.2. Notons que les fonctions $B^{\tilde{L}_{i}}$ pour $i=1,2$ qui interviennent sont les restrictions de la fonction $B^{\tilde{G}}$. Pour simplifier, on les note encore $B^{\tilde{G}}$. 
On obtient
$$ I_{K\tilde{M}}^{K\tilde{G},{\cal E}}({\bf M}',\boldsymbol{\delta},{\bf f})=\sum_{\tilde{s}\in \tilde{\zeta}Z(\hat{M})^{\Gamma_{F},\hat{\theta}}/Z(\hat{G})^{\Gamma_{F},\hat{\theta}}}i_{\tilde{M}'}(\tilde{G},\tilde{G}'(\tilde{s}))$$
$$\sum_{\tilde{L}'_{\tilde{s},1},\tilde{L}'_{\tilde{s},2}\in {\cal L}^{\tilde{G}'(\tilde{s})}(\tilde{M})}e_{\tilde{M}}^{\tilde{G}'(\tilde{s})}(\tilde{L}'_{\tilde{s},1},\tilde{L}'_{\tilde{s},2})S_{{\bf M}'}^{{\bf L}'_{1}(\tilde{s})}(\boldsymbol{\delta}_{1},B^{\tilde{G}},({\bf f}_{1}^{{\bf G}'(\tilde{s})})_{{\bf L}'_{1}(\tilde{s})})S_{{\bf M}'}^{{\bf L}'_{2}(\tilde{s})}(\boldsymbol{\delta}_{2},B^{\tilde{G}},({\bf f}_{2}^{{\bf G}'(\tilde{s})})_{{\bf L}'_{2}(\tilde{s})}).$$
Expliquons la notation. 
Un espace de Levi $\tilde{L}'_{\tilde{s},1}$ intervenant dans cette formule d\'etermine un  $K$-espace de Levi   $K\tilde{L}_{1} \in {\cal L}(K\tilde{M})$ par l'\'egalit\'e ${\cal A}_{\tilde{L}_{1}}={\cal A}_{\tilde{L}'_{\tilde{s},1}}$. Alors $\tilde{L}'_{\tilde{s},1}$ s'identifie \`a l'espace $\tilde{L}'_{1}(\tilde{s})$ associ\'e \`a la donn\'ee endoscopique elliptique ${\bf L}'_{1}(\tilde{s})$ de $(KL_{1}, K\tilde{L}_{1},{\bf a}_{L_{1}})$.
On a l'\'egalit\'e $({\bf f}^{{\bf G}'(\tilde{s})})_{{\bf L}'_{\tilde{s},1}}=({\bf f}_{\tilde{L}_{1},\omega})^{{\bf L}'_{1}(\tilde{s})}$. Regroupons les termes selon les couples $(K\tilde{L}_{1},K\tilde{L}_{2})$ d'espaces de Levi qui apparaissent. Pour un tel couple, notons  $S(K\tilde{L}_{1},K\tilde{L}_{2})$  l'ensemble des $\tilde{s}\in \tilde{\zeta}Z(\hat{M})^{\Gamma_{F},\hat{\theta}}/Z(\hat{G})^{\Gamma_{F},\hat{\theta}}$ tels que ${\bf L}'_{i}(\tilde{s})$ soit elliptique pour $i=1,2$. On obtient
$$(6) \qquad  I_{K\tilde{M}}^{K\tilde{G},{\cal E}}({\bf M}',\boldsymbol{\delta},{\bf f})=\sum_{K\tilde{L}_{1},K\tilde{L}_{2}\in {\cal L}(K\tilde{M})}X(K\tilde{L}_{1},K\tilde{L}_{2}),$$
o\`u
$$X(K\tilde{L}_{1},K\tilde{L}_{2})=\sum_{\tilde{s}\in S(K\tilde{L}_{1},K\tilde{L}_{2})}i_{\tilde{M}'}(\tilde{G},\tilde{G}'(\tilde{s}))e_{\tilde{M}}^{\tilde{G}'(\tilde{s})}(\tilde{L}'_{1}(\tilde{s}),\tilde{L}'_{2}(\tilde{s}))$$
$$S_{\tilde{M}'}^{{\bf L}'_{1}(\tilde{s})}(\boldsymbol{\delta}_{1},B^{\tilde{G}},({\bf f}_{1,\tilde{L}_{1},\omega})^{{\bf L}'_{1}(\tilde{s})})S_{\tilde{M}'}^{{\bf L}'_{2}(\tilde{s})}(\boldsymbol{\delta}_{2},B^{\tilde{G}},({\bf f}_{2,\tilde{L}_{2},\omega})^{{\bf L}'_{2}(\tilde{s})}).$$
  Fixons $K\tilde{L}_{1},K\tilde{L}_{2}\in {\cal L}(K\tilde{M})$. 
   Soit $\tilde{s}$ contribuant de fa\c{c}on non nulle \`a la somme $X(\tilde{L}_{1},\tilde{L}_{2})$. La donn\'ee ${\bf G}'(\tilde{s})$ doit \^etre elliptique (sinon 
$i_{\tilde{M}'}(\tilde{G},\tilde{G}'(\tilde{s}))=0$). Les donn\'ees ${\bf L}'_{i}(\tilde{s})$ aussi, pour $i=1,2$. Le coefficient $e_{\tilde{M}}^{\tilde{G}'(\tilde{s})}(\tilde{L}'_{1}(\tilde{s}),\tilde{L}'_{2}(\tilde{s}))$ contient $d_{\tilde{M}}^{\tilde{G}'(\tilde{s})}(\tilde{L}'_{1}(\tilde{s}),\tilde{L}'_{2}(\tilde{s}))$. Celui-ci ne d\'epend que des espaces ${\cal A}_{\tilde{M}}$, ${\cal A}_{\tilde{G}'(\tilde{s})}$ et  ${\cal A}_{\tilde{L}'_{i}(\tilde{s})}$ pour $i=1,2$. Mais, par ellipticit\'e, ces derniers sont  respectivement \'egaux \`a ${\cal A}_{\tilde{G}}$ et ${\cal A}_{\tilde{L}_{i}}$ pour $i=1,2$. On obtient l'\'egalit\'e
$$d_{\tilde{M}}^{\tilde{G}'(\tilde{s})}(\tilde{L}'_{1}(\tilde{s}),\tilde{L}'_{2}(\tilde{s}))=d_{\tilde{M}}^{\tilde{G}}(\tilde{L}_{1},\tilde{L}_{2}).$$
Si ce dernier terme est nul, on a donc $X(K\tilde{L}_{1},K\tilde{L}_{2})=0$. Supposons $d_{\tilde{M}}^{\tilde{G}}(\tilde{L}_{1},\tilde{L}_{2})\not=0$. Le calcul ci-dessus permet de r\'ecrire 
$$(7) \qquad X(K\tilde{L}_{1},K\tilde{L}_{2})=d_{\tilde{M}}^{\tilde{G}}(\tilde{L}_{1},\tilde{L}_{2})\sum_{\tilde{s}\in S(K\tilde{L}_{1},K\tilde{L}_{2})}i_{\tilde{M}'}(\tilde{G},\tilde{G}'(\tilde{s}))k_{\tilde{M}}^{\tilde{G}'(\tilde{s})}(\tilde{L}'_{1}(\tilde{s}),\tilde{L}'_{2}(\tilde{s}))^{-1}$$
$$S_{\tilde{M}'}^{{\bf L}'_{1}(\tilde{s})}(\boldsymbol{\delta}_{1},B^{\tilde{G}},({\bf f}_{1,\tilde{L}_{1},\omega})^{{\bf L}'_{1}(\tilde{s})})S_{\tilde{M}'}^{{\bf L}'_{2}(\tilde{s})}(\boldsymbol{\delta}_{2},B^{\tilde{G}},({\bf f}_{2,\tilde{L}_{2},\omega})^{{\bf L}'_{2}(\tilde{s})}).$$
Introduisons l'homomorphisme
$$q:Z(\hat{M})^{\Gamma_{F},\hat{\theta}}/Z(\hat{G})^{\Gamma_{F},\hat{\theta}}\to Z(\hat{M})^{\Gamma_{F},\hat{\theta}}/Z(\hat{L}_{1})^{\Gamma_{F},\hat{\theta}}\oplus Z(\hat{M})^{\Gamma_{F},\hat{\theta}}/Z(\hat{L}_{2})^{\Gamma_{F},\hat{\theta}}.$$
L'hypoth\`ese $d_{\tilde{M}}^{\tilde{G}}(\tilde{L}_{1},\tilde{L}_{2})\not=0$ implique qu'il est surjectif et de noyau fini. Pour $\tilde{s}=s\tilde{\zeta}$ avec $s\in Z(\hat{M})^{\Gamma_{F},\hat{\theta}}/Z(\hat{G})^{\Gamma_{F},\hat{\theta}}$, posons $q(s)=(s_{1},s_{2})$ et $\tilde{s}_{i}=s_{i}\tilde{\zeta}$
 pour $i=1,2$. On a $\tilde{L}_{i}(\tilde{s})=\tilde{L}_{i}(\tilde{s}_{i})$. Montrons que
 
 (8) pour $\tilde{s}\in S(K\tilde{L}_{1},K\tilde{L}_{2})$, on a l'\'egalit\'e
 $$i_{\tilde{M}'}(\tilde{G},\tilde{G}'(\tilde{s}))k_{\tilde{M}}^{\tilde{G}'(\tilde{s})}(\tilde{L}'_{1}(\tilde{s}),\tilde{L}'_{2}(\tilde{s}))^{-1}=k_{\tilde{M}}^{\tilde{G}}(\tilde{L}_{1},\tilde{L}_{2})^{-1}i_{\tilde{M}'}(\tilde{L}_{1},\tilde{L}'_{1}(\tilde{s}_{1}))i_{\tilde{M}'}(\tilde{L}_{2},\tilde{L}'_{2}(\tilde{s}_{2})).$$

 On a 
$${\cal A}_{\tilde{G}'(\tilde{s})}\subset {\cal A}_{\tilde{L}'_{1}(\tilde{s})}\cap {\cal A}_{\tilde{L}'_{1}(\tilde{s})}={\cal A}_{\tilde{L}_{1}}\cap {\cal A}_{\tilde{L}_{2}}={\cal A}_{\tilde{G}}.$$
Donc ${\bf G}'(\tilde{s})$ est elliptique et $i_{\tilde{M}'}(\tilde{G},\tilde{G}'(\tilde{s}))$ est l'inverse du nombre d'\'el\'ements du noyau de l'homomorphisme
$$r:Z(\hat{M})^{\Gamma_{F},\hat{\theta}}/Z(\hat{G})^{\Gamma_{F},\hat{\theta}}\to Z(\hat{M}')^{\Gamma_{F}}/Z(\hat{G}'(\tilde{s}))^{\Gamma_{F}}.$$
Introduisons les homomorphismes
$$\begin{array}{cc}r_{12}:Z(\hat{M})^{\Gamma_{F},\hat{\theta}}/Z(\hat{L}_{1})^{\Gamma_{F},\hat{\theta}}\oplus Z(\hat{M})^{\Gamma_{F},\hat{\theta}}/Z(\hat{L}_{2})^{\Gamma_{F},\hat{\theta}}\to &Z(\hat{M}')^{\Gamma_{F}}/Z(\hat{L}'_{1}(\tilde{s}_{1}))^{\Gamma_{F}}\\ &\oplus Z(\hat{M}')^{\Gamma_{F}}/Z(\hat{L}'_{2}(\tilde{s}_{2}))^{\Gamma_{F}}\\ \end{array}$$
et
$$q_{12}:Z(\hat{M}')^{\Gamma_{F}}/Z(\hat{G}'(\tilde{s}))^{\Gamma_{F}}\to Z(\hat{M}')^{\Gamma_{F}}/Z(\hat{L}'_{1}(\tilde{s}_{1}))^{\Gamma_{F}}\oplus Z(\hat{M}')^{\Gamma_{F}}/Z(\hat{L}'_{2}(\tilde{s}_{2}))^{\Gamma_{F}}.$$
On a l'\'egalit\'e $r_{12}\circ q=q_{12}\circ r$. Tous ces homomorphismes sont surjectifs et de noyaux finis. Le nombre d'\'el\'ements du noyau du compos\'e est donc \'egal au produit des nombres d'\'el\'ements des noyaux de $r_{12}$ et de $q$, ou de $q_{12}$ et de $r$. Ces  noyaux ont pour nombre d'\'el\'ements $i_{\tilde{M}'}(\tilde{L}_{1},\tilde{L}'_{1}(\tilde{s}_{1}))^{-1}i_{\tilde{M}'}(\tilde{L}_{2},\tilde{L}'_{2}(\tilde{s}_{2}))^{-1}$ pour $r_{12}$, $k_{\tilde{M}}^{\tilde{G}}(\tilde{L}_{1},\tilde{L}_{2})$ pour $q$, $k_{\tilde{M}}^{\tilde{G}'(\tilde{s})}(\tilde{L}'_{1}(\tilde{s}),\tilde{L}'_{2}(\tilde{s}))$ pour $q_{12}$ et $i_{\tilde{M}'}(\tilde{G},\tilde{G}'(\tilde{s}))^{-1}$ pour $r$. L'assertion (8) s'ensuit. 

Il r\'esulte de (8) que le terme que l'on somme dans (7) ne d\'epend que du couple $(\tilde{s}_{1},\tilde{s}_{2})$. On peut r\'ecrire (7) comme une somme sur ces couples $(\tilde{s}_{1},\tilde{s}_{2})$, le terme que l'on somme \'etant multipli\'e par le nombre des $\tilde{s}\in S(K\tilde{L}_{1},K\tilde{L}_{2})$ qui se projettent sur le couple en question. Par d\'efinition de $S(K\tilde{L}_{1},K\tilde{L}_{2})$, ce nombre est nul si l'un des $\tilde{L}'_{i}(\tilde{s}_{i})$ n'est pas elliptique. Sinon, c'est le nombre d'\'el\'ements du noyau de $q$, c'est-\`a-dire $k_{\tilde{M}}^{\tilde{G}}(\tilde{L}_{1},\tilde{L}_{2})$. Ce terme compense le m\^eme facteur intervenant dans le membre de droite de (8). D'autre part, on se rappelle que la condition $i_{\tilde{M}'}(\tilde{L}_{i},\tilde{L}'_{i}(\tilde{s}_{i}))\not=0$ \'equivaut \`a l'ellipticit\'e de $\tilde{L}'_{i}(\tilde{s}_{i})$. On obtient l'\'egalit\'e
$$X(K\tilde{L}_{1},K\tilde{L}_{2})=d_{\tilde{M}}^{\tilde{G}}(\tilde{L}_{1},\tilde{L}_{2})\prod_{i=1,2}
\sum_{\tilde{s_{i}}\in  \tilde{\zeta}Z(\hat{M})^{\Gamma_{F},\hat{\theta}}/Z(\hat{L}_{i})^{\Gamma_{F},\hat{\theta}}}i_{\tilde{M}'}(\tilde{L}_{i},\tilde{L}_{i}'(\tilde{s}_{i})) $$
$$
S_{\tilde{M}'}^{{\bf L}'_{i}(\tilde{s}_{i})}(\boldsymbol{\delta}_{i},B^{\tilde{G}},({\bf f}_{i,\tilde{L}_{i},\omega})^{{\bf L}'_{i}(\tilde{s}_{i})}).$$
 La somme en $\tilde{s}_{i}$ est  \'egale par d\'efinition \`a $I_{K\tilde{M}}^{K\tilde{L}_{i},{\cal E}}({\bf M}',\boldsymbol{\delta}_{i},{\bf f}_{i,\tilde{L}_{i},\omega})$. Donc
 $$X(K\tilde{L}_{1},K\tilde{L}_{2})=d_{\tilde{M}}^{\tilde{G}}(\tilde{L}_{1},\tilde{L}_{2})I_{K\tilde{M}}^{K\tilde{L}_{1},{\cal E}}({\bf M}',\boldsymbol{\delta}_{1},{\bf f}_{1,\tilde{L}_{1},\omega})I_{K\tilde{M}}^{K\tilde{L}_{2},{\cal E}}({\bf M}',\boldsymbol{\delta}_{2},{\bf f}_{2,\tilde{L}_{2},\omega}).$$
 Les membres de droite de (4) et (6) sont alors \'egaux, ce qui d\'emontre (4). 
 
 Esquissons la preuve de (5). On a ici $V=\{v\}$.  La donn\'ee ${\bf M}'_{v}$ de $(KM_{v},K\tilde{M}_{v},{\bf a}_{M_{v}})$ n'est pas forc\'ement elliptique mais elle est relevante. Il lui correspond un $K$-espace de Levi $K\tilde{R}$ de $K\tilde{M}$. D'apr\`es nos d\'efinitions, le membre de droite de (5) est \'egal \`a
 $$\sum_{K\tilde{L}^v\in {\cal L}(K\tilde{M}_{v})}d_{\tilde{M}_{v}}^{\tilde{G}}(\tilde{M},\tilde{L}^v)\sum_{K\tilde{S}^v\in {\cal L}(K\tilde{R}_{v})}d_{\tilde{R}_{v}}^{\tilde{L}^v}(\tilde{M}_{v},\tilde{S}^v)I_{K\tilde{R}_{v}}^{K\tilde{S}^v,{\cal E}}({\bf M}'_{v},\boldsymbol{\delta}_{v},{\bf f}_{v,K\tilde{S}^v,\omega}).$$
 En \'echangeant les lettres $L$ et $S$, on obtient
$$ \sum_{K\tilde{L}^v\in {\cal L}(K\tilde{R}_{v})}X_{\tilde{R}_{v}}^{\tilde{G}}(\tilde{M},\tilde{L}^v)I_{K\tilde{R}_{v}}^{K\tilde{L}^v,{\cal E}}({\bf M}'_{v},\boldsymbol{\delta}_{v},{\bf f}_{v,K\tilde{L}^v,\omega})$$
o\`u
$$X_{\tilde{R}_{v}}^{\tilde{G}}(\tilde{M},\tilde{L}^v)=\sum_{K\tilde{S}^v\in {\cal L}(K\tilde{L}^v)}d_{\tilde{M}_{v}}^{\tilde{G}}(\tilde{M},\tilde{S}^v)d_{\tilde{R}_{v}}^{\tilde{S}^v}(\tilde{M}_{v},\tilde{L}^v).$$
Un calcul d\'ej\`a fait plusieurs fois prouve que $X_{\tilde{R}_{v}}^{\tilde{G}}(\tilde{M},\tilde{L}^v)=d_{\tilde{R}_{v}}^{\tilde{G}}(\tilde{M},\tilde{L}^v)$. La relation (5) \'equivaut donc \`a
$$I_{K\tilde{M}}^{K\tilde{G},{\cal E}}({\bf M}',\boldsymbol{\delta},{\bf f})= \sum_{K\tilde{L}^v\in {\cal L}(K\tilde{R}_{v})}d_{\tilde{R}_{v}}^{\tilde{G}}(\tilde{M},\tilde{L}^v)I_{K\tilde{R}_{v}}^{K\tilde{L}^v,{\cal E}}({\bf M}'_{v},\boldsymbol{\delta}_{v},{\bf f}_{v,K\tilde{L}^v,\omega}).$$
La preuve de cette \'egalit\'e est similaire \`a celle de (4)  et d'ailleurs presque identique \`a celle de la proposition [II] 1.14(i). On la laisse au lecteur. $\square$

{\bf Remarque.} La preuve de (4) est r\'eversible en ce sens que, si l'on suppose v\'erifi\'ee la relation (i) de l'\'enonc\'e ainsi que la formule (ii) de la proposition 4.2 pour tous les termes sauf un du membre de droite de la d\'efinition (2), on en d\'eduit cette formule (ii) pour le terme restant. Dans cette direction, cela prouve par r\'ecurrence cette formule (ii) de la proposition 4.2, puisque dans la situation de ce paragraphe, l'analogue de la relation (i) est connue pour le terme $I_{\tilde{M}}^{\tilde{G}}(\boldsymbol{\delta},{\bf f})$.

\bigskip

Pourquoi avoir travaill\'e ici avec des $K$-espaces? Parce que, dans le cas local et pour un unique triplet $(G,\tilde{G},{\bf a})$, on ne conna\^{\i}t pas (pas encore, plut\^ot) l'analogue de la relation (ii) de l'\'enonc\'e. Expliquons cela. Supposons qu'un triplet $(G,\tilde{G},{\bf a})$ soit une composante connexe d'un $K$-triplet $(KG,K\tilde{G},{\bf a})$, disons que c'est la composante index\'ee par $p\in \Pi$. On peut appliquer la relation (ii) de l'\'enonc\'e \`a une fonction ${\bf f}\in I(\tilde{G}(F_{V}),\omega)\otimes Mes(G(F_{V}))$ identifi\'ee \`a un \'el\'ement de $I(K\tilde{G}(F_{V}),\omega)\otimes Mes(G(F_{V}))$ nul sur les autres composantes connexes. Mais $transfert(\boldsymbol{\delta})$ vit sur toutes les composantes connexes, il est de la forme $\oplus_{q\in \Pi}\boldsymbol{\gamma}_{q}$, o\`u $\boldsymbol{\gamma}_{q}\in D_{g\acute{e}om}(\tilde{M}_{q}(F_{V}),\omega)\otimes Mes(M(F_{V}))^*$. Le membre de gauche de (ii) est \'egal \`a
$$\sum_{q\in \Pi}I_{K\tilde{M}}^{K\tilde{G},{\cal E}}(\boldsymbol{\gamma}_{q},{\bf f}).$$
Bien que ${\bf f}$ soit concentr\'ee sur la composante $\tilde{G}_{p}$, on ne sait pas a priori que $I_{K\tilde{M}}^{K\tilde{G},{\cal E}}(\boldsymbol{\gamma}_{q},{\bf f})=0$ pour $q\not=p$. C'est l'unique raison, nous semble-t-il, pour laquelle nous devrons travailler avec des $K$-espaces. Notons toutefois que, pour un seul triplet $(G,\tilde{G},{\bf a})$, on peut parfaitement d\'efinir $I_{\tilde{M}}^{\tilde{G}}({\bf M}',\boldsymbol{\delta},{\bf f})$ comme le membre de droite de la formule (1). Si on inclut $(G,\tilde{G},{\bf a})$ comme composante connexe d'un triplet $(KG,K\tilde{G},{\bf a})$, ce terme co\"{\i}ncide avec $I_{K\tilde{M}}^{K\tilde{G}}({\bf M}',\boldsymbol{\delta},{\bf f})$, o\`u ${\bf f}$ est identifi\'e \`a une fonction sur $K\tilde{G}(F_{V})$ nulle sur les autres composantes. 

Dans le cas o\`u $(G,\tilde{G},{\bf a})$ est quasi-d\'eploy\'e et \`a torsion int\'erieure, si l'on fixe un syst\`eme de fonctions $B$ comme en 1.10, on peut aussi d\'efinir le terme $I_{\tilde{M}}^{\tilde{G}}({\bf M}',\boldsymbol{\delta},B,{\bf f})$.

\bigskip

\subsection{Le r\'esultat de comparaison des int\'egrales orbitales pond\'er\'ees $\omega$-\'equivariantes}
Les hypoth\`eses sont les m\^emes que dans le paragraphe pr\'ec\'edent.
\ass{Proposition (\`a prouver)}{Pour tout $\boldsymbol{\gamma}\in D_{g\acute{e}om,\tilde{G}-\acute{e}qui}(K\tilde{M}(F_{V}),\omega)\otimes Mes(M(F_{V}))^*$ et tout ${\bf f}\in I(K\tilde{G}(F_{V}),\omega)\otimes Mes(G(F_{V}))$, on a l'\'egalit\'e $I^{K\tilde{G},{\cal E}}_{K\tilde{M}}(\boldsymbol{\gamma},{\bf f})=I^{K\tilde{G}}_{K\tilde{M}}(\boldsymbol{\gamma},{\bf f})$.  }

D'apr\`es la d\'efinition du membre de gauche et la formule de descente 1.11(1), la proposition r\'esulte de ses analogues locales, c'est-\`a-dire des th\'eor\`emes 1.16 de [II] et 1.10 de [V], qui restent \`a prouver.

\bigskip

\subsection{Une autre forme du r\'esultat de comparaison}
On conserve la m\^eme situation. 
\ass{Proposition}{On admet la validit\'e des th\'eor\`emes 1.16 de [II] et  1.10 de [V]. Soit ${\bf M}'$ une donn\'ee endoscopique de $(KM,K\tilde{M},{\bf a}_{M})$ elliptique et relevante. Soit $\boldsymbol{\delta}\in D_{tr-orb}^{st}({\bf M}'_{V})\otimes Mes(M'(F_{V}))^*$.  Alors, pour tout ${\bf f}\in I(K\tilde{G}(F_{V}),\omega)\otimes Mes(G(F_{V}))^*$, on a l'\'egalit\'e
$$I_{K\tilde{M}}^{K\tilde{G}}(transfert(\boldsymbol{\delta}),{\bf f})=I_{K\tilde{M}}^{K\tilde{G},{\cal E}}({\bf M}',\boldsymbol{\delta},{\bf f}).$$}

Remarquons que $transfert(\boldsymbol{\delta})$ appartient \`a $D_{tr-orb}(K\tilde{M}(F_{V}),\omega)\otimes Mes(M(F_{V}))^*$. Comme on l'a vu dans la section 5 de  [V], la validit\'e du th\'eor\`eme [V] 1.10 permet de d\'efinir   le membre de gauche de l'\'egalit\'e de l'\'enonc\'e. 

Preuve. Le membre de gauche de l'\'egalit\'e \`a prouver v\'erifie la formule de descente 1.11(1). Celui de droite v\'erifie la formule parall\`ele de la proposition 4.5(i). Cela nous ram\`ene \`a prouver l'analogue locale de l'\'egalit\'e. En une place non-archim\'edienne, cette \'egalit\'e r\'esulte directement du th\'eor\`eme [II] 1.16. En une place r\'eelle, on a vu dans la section 5 de [V] qu'elle r\'esultait (moins directement) du th\'eor\`eme [V] 1.10. $\square$

\bigskip

\subsection{Le cas quasi-d\'eploy\'e et \`a torsion int\'erieure}
Soient $(G,\tilde{G},{\bf a})$ un triplet quasi-d\'eploy\'e et \`a torsion int\'erieure, $B$ un syst\`eme de fonctions comme en 1.10 et  $V$ un ensemble fini de places de $F$. Soient $\tilde{M}$ un espace de Levi de $\tilde{G}$ et ${\bf M}'$ une donn\'ee endoscopique elliptique et relevante de $(M,\tilde{M})$.  

\ass{Corollaire}{Pour tout $\boldsymbol{\delta}\in D_{tr-orb}({\bf M}')\otimes Mes(M'(F_{V}))^*$ et tout ${\bf f}\in C_{c}^{\infty}(\tilde{G}(F_{V}))\otimes Mes(G(F_{V}))$, on a l'\'egalit\'e 
 $$I_{\tilde{M}}^{\tilde{G}}(transfert(\boldsymbol{\delta}),B,{\bf f})=I_{\tilde{M}}^{\tilde{G}}({\bf M}',\boldsymbol{\delta},B,{\bf f}).$$}
 
 L'argument est le m\^eme que dans le paragraphe pr\'ec\'edent. Parce que $(G,\tilde{G},{\bf a})$ est quasi-d\'eploy\'e et \`a torsion int\'erieure, le th\'eor\`eme [II] 1.16 est prouv\'e (cf. [III] proposition 2.9) et un substitut du th\'eor\`eme [V] 1.10 aussi: c'est la proposition [V] 1.13 dont on a vu dans la section 4 de [V] qu'elle suffisait \`a notre propos.

\bigskip

\section{La formule des traces stable}

\bigskip

 \subsection{Quelques d\'efinitions}
Soit $V$ un ensemble fini de places de $F$ contenant $V_{ram}$. On fixe un ensemble de repr\'esentants ${\cal E}(\tilde{G},{\bf a},V)$ des classes d'\'equivalence de donn\'ees endoscopiques de $(G,\tilde{G},{\bf a})$ qui sont elliptiques, relevantes et non ramifi\'ees hors de $V$.

Soit ${\bf G}'=(G',{\cal G}',\tilde{s})\in {\cal E}(\tilde{G},{\bf a},V)$. On munit $G'$ des objets auxiliaires de 3.1. En particulier, pour $v\not\in V_{ram}({\bf G}')$, on fixe un sous-espace hypersp\'ecial $\tilde{K}^{G'}_{v}$ de $\tilde{G}'(F_{v})$ correspondant \`a $\tilde{K}_{v}$. De l'homomorphisme naturel $A_{\tilde{G}}\to A_{G'}$ se d\'eduit un homomorphisme ${\cal A}_{\tilde{G}}\to {\cal A}_{G'}$. C'est un isomorphisme par l'hypoth\`ese d'ellipticit\'e. Il pr\'eserve les mesures par d\'efinition de celles-ci.

Consid\'erons $\underline{les}$ paires de Borel \'epingl\'ees de $G$ et $G'$, dont on note les tores $T^*$ et $T^{'*}$. Les classes de conjugaison semi-simples dans $\tilde{G}(\bar{F})$, resp. $\tilde{G}'(\bar{F})$, sont classifi\'ees par 
$$((T^*/(1-\theta)(T^*))/W^{\theta})\times_{ {\cal Z}(G)}{\cal Z}(\tilde{G}),$$
resp.
$$(T^{'*}/W^{G'})\times_{Z(G')}{\cal Z}(\tilde{G}').$$
Il y a une application naturelle du second ensemble dans le premier, autrement dit une application qui, \`a une classe de conjugaison semi-simple dans $\tilde{G}'(\bar{F})$, associe une telle classe dans $\tilde{G}(\bar{F})$. Cette application est \'equivariante pour les actions galoisiennes.  Les classes de conjugaison g\'eom\'etriques (c'est-\`a-dire par $G(\bar{F})$) semi-simples dans $\tilde{G}(F)$ sont classifi\'ees par un sous-ensemble de
$$\left(((T^*/(1-\theta)(T^*))/W^{\theta})\times_{ {\cal Z}(G)}{\cal Z}(\tilde{G})\right)^{\Gamma_{F}}.$$
La description exacte de ce sous-ensemble est compliqu\'ee.

On d\'efinit la conjugaison stable entre \'el\'ements semi-simples de $\tilde{G}(F)$ comme dans le cas local. Deux \'el\'ements $\eta,\eta'\in \tilde{G}_{ss}(F)$ sont stablement conjugu\'es si et seulement s'il existe $g\in G(\bar{F})$ tel que $g^{-1}\eta g=\eta'$ et $g\sigma(g)^{-1}\in I_{\eta}(\bar{F})$ pour tout $\sigma\in \Gamma_{F}$, o\`u $I_{\eta}=Z(G)^{\theta}G_{\eta}$. La correspondance ci-dessus se raffine en une correspondance entre classes de conjugaison stable dans $\tilde{G}_{ss}'(F)$ et classes de conjugaison stable dans $\tilde{G}_{ss}(F)$. 

De m\^eme, si on se place sur  l'anneau $F_{V}$, la correspondance ci-dessus se raffine en une correspondance entre classes de conjugaison stable dans $\tilde{G}_{ss}'(F_{V})$ et classes de conjugaison stable dans $\tilde{G}_{ss}(F_{V})$. Si ${\cal O}$ est une classe de conjugaison stable dans $\tilde{G}_{ss}(F_{V})$, on note ${\cal O}_{\tilde{G}'}$ la r\'eunion (finie, \'eventuellement vide) des classes de conjugaison stable dans $\tilde{G}'_{ss}(F_{V})$ qui correspondent \`a ${\cal O}$.

On note $Aut({\bf G}')$ le groupe d'automorphismes de ${\bf G}'$. En fixant une paire de Borel \'epingl\'ee convenable de $\hat{G}$, on a une suite exacte
$$1\to (Z(\hat{G})/(Z(\hat{G})\cap \hat{T}^{\hat{\theta},0}))^{\Gamma_{F}}\to Aut({\bf G}')/\hat{G}'\to Out({\bf G}')\to 1,$$
o\`u
$Out({\bf G}')$ est un sous-groupe du groupe des automorphismes ext\'erieurs de $G'$. On note selon l'usage $\pi_{0}(X)$ le groupe des composantes connnexes d'un groupe alg\'ebrique complexe $X$. On pose
$$i(\tilde{G},\tilde{G}')=\vert Out({\bf G}')\vert ^{-1}\vert det((1-\theta)_{\vert \mathfrak{A}_{G}/\mathfrak{A}_{\tilde{G}}})\vert ^{-1}\vert \pi_{0}(Z(\hat{G})^{\Gamma_{F}})\vert \vert ker^1(F,Z(\hat{G}))\vert ^{-1}$$
$$\vert \pi_{0}((Z(\hat{G})/(Z(\hat{G})\cap \hat{T}^{\hat{\theta},0}))^{\Gamma_{F}})\vert ^{-1}\vert \pi_{0}(Z(\hat{G})^{\Gamma_{F},0}\cap \hat{T}^{\hat{\theta},0})\vert \vert \pi_{0}(Z(\hat{G}')^{\Gamma_{F}})\vert^{-1} \vert ker^1(F,Z(\hat{G}'))\vert  .$$
Remarquons que, par exemple, le produit $ \vert \pi_{0}(Z(\hat{G})^{\Gamma_{F}})\vert \vert ker^1(F,Z(\hat{G}))\vert ^{-1}$ peut s'interpr\'eter comme le nombre de Tamagawa de $G'$. Notons-le $\tau(G')$. La formule ci-dessus se r\'ecrit
$$i(\tilde{G},\tilde{G}')=\vert Aut({\bf G}')\vert ^{-1}det((1-\theta)_{\vert \mathfrak{A}_{G}/\mathfrak{A}_{\tilde{G}}})\vert ^{-1}\vert \pi_{0}(Z(\hat{G})^{\Gamma_{F},0}\cap \hat{T}^{\hat{\theta},0})\vert \tau(G)\tau(G')^{-1}.$$

On a d\'efini la distribution $A^{\tilde{G}}(V,{\cal O},\omega)\in D_{orb}(\tilde{G}(F_{V}),\omega)\otimes Mes(G(F_{V}))^*$ pour une classe de conjugaison par $G(F_{V})$ dans $\tilde{G}_{ss}(F_{V})$. Si ${\cal O}$ est maintenant une r\'eunion finie de telles classes, on note $A^{\tilde{G}}(V,{\cal O},\omega)$ la somme des distributions associ\'ees \`a chacune de ces classes. En particulier, $A^{\tilde{G}}(V,{\cal O},\omega)$ est d\'efini pour une classe de conjugaison stable ${\cal O}$ dans $\tilde{G}_{ss}(F_{V})$.

\bigskip

\subsection{Les distributions $SA^{\tilde{G}}(V,{\cal O})$}

La situation est la m\^eme que dans le paragraphe pr\'ec\'edent. On suppose de plus $(G,\tilde{G},{\bf a})$ quasi-d\'eploy\'e et \`a torsion int\'erieure. Pour toute classe de conjugaison stable ${\cal O}$ dans $\tilde{G}_{ss}(F_{V})$, on va d\'efinir une distribution $SA^{\tilde{G}}(V,{\cal O})\in D_{tr-orb}(\tilde{G}(F_{V}))\otimes Mes(G(F_{V}))^*$. Comme toujours, on a besoin de supposer par r\'ecurrence que cette distribution v\'erifie certaines  propri\'et\'es. Il y a les propri\'et\'es formelles qui permettent de "recoller" ces distributions dans la situation de 1.15. Elles sont faciles \`a v\'erifier par r\'ecurrence et on les abandonne au lecteur. Il y a une autre propri\'et\'e plus subtile.  La distribution $A^{\tilde{G}}(V,{\cal O})$ d\'epend des $\tilde{K}_{v}$ pour $v\not\in V$, ou plus exactement des classes de conjugaison par $G(F_{v})$ des $\tilde{K}_{v}$. La d\'efinition ci-dessous fournit une distribution $SA^{\tilde{G}}(V,{\cal O})$ qui en d\'epend aussi. De fait, elle en d\'epend. Mais on a besoin de savoir que

(1) elle ne  d\'epend que des classes de conjugaison  par $G_{AD}(F_{v})$ des $\tilde{K}_{v}$, pour $v\not\in V$.

On a surtout besoin de supposer par r\'ecurrence que cette distribution $SA^{\tilde{G}}(V,{\cal O})$ est stable. Modulo ces propri\'et\'es, si ${\bf G}'=(G',{\cal G}',s)$ est une donn\'ee endoscopique de $(G,\tilde{G},{\bf a})$ non ramifi\'ee hors de $V$,  avec $dim(G'_{SC})<dim(G_{SC})$, et si ${\cal O}'$ est une classe de conjugaison stable dans $\tilde{G}'_{ss}(F_{V})$, on peut d\'efinir la distribution $SA^{{\bf G}'}(V,{\cal O}')\in D^{st}_{g\acute{e}om}({\bf G}'_{V})\otimes Mes(G'(F_{V}))^*$. 

{\bf Remarque.} Quant \`a la d\'ependance des espaces hypersp\'eciaux, une extension formelle de la propri\'et\'e (1) montre que cette distribution ne d\'epend que des 
  classes de conjugaison par $G_{AD}'(F_{v})$ des $\tilde{K}'_{v}$, pour $v\not\in V$. Mais ces classes sont   bien d\'etermin\'ees par les $\tilde{K}_{v}$. Donc  $SA^{{\bf G}'}(V,{\cal O}')$ ne d\'epend que des $\tilde{K}_{v}$. 
  \bigskip
  
  On peut alors poser la d\'efinition
$$SA^{\tilde{G}}(V,{\cal O})=A^{\tilde{G}}(V,{\cal O})-\sum_{{\bf G}'\in {\cal E}(\tilde{G},V), G'\not=G}i(\tilde{G},\tilde{G}')transfert(SA^{{\bf G}'}(V,{\cal O}_{\tilde{G}'})).$$
 Remarquons que la d\'efinition entra\^{\i}ne par r\'ecurrence que $SA^{\tilde{G}}(V,{\cal O})$ est \`a support dans l'ensemble des \'el\'ements de $\tilde{G}(F_{V})$ dont la partie semi-simple appartient \`a ${\cal O}$. 
 
   Enon\c{c}ons les propri\'et\'es de notre distribution sous la forme d'un th\'eor\`eme \`a prouver.
  
  \ass{Th\'eor\`eme (\`a prouver)}{On suppose $(G,\tilde{G},{\bf a})$ quasi-d\'eploy\'e et \`a torsion int\'erieure. Soit ${\cal O}$ une classe de conjugaison stable semi-simple dans $\tilde{G}(F_{V})$. Alors $SA^{\tilde{G}}(V,{\cal O})$ est stable et v\'erifie (1).}
  
  \bigskip

\subsection{Propri\'et\'es des distributions $SA^{\tilde{G}}(V,{\cal O})$}
On suppose $(G,\tilde{G},{\bf a})$ quasi-d\'eploy\'e et \`a torsion int\'erieure. Soient $V$ un ensemble fini de places de $F$ contenant $V_{ram}$ et ${\cal O}$ une classe de conjugaison stable dans $\tilde{G}_{ss}(F_{V})$. On a

(1) si $SA^{\tilde{G}}(V,{\cal O})\not=0$ alors il existe un \'el\'ement semi-simple $\gamma\in \tilde{G}_{ss}(F)$ dont la projection  dans $\tilde{G}(F_{V})$ appartienne \`a ${\cal O}$ et tel que $\tilde{H}_{\tilde{G}_{V}}(\gamma)=0$.

Preuve. D'apr\`es la d\'efinition de 5.1, la condition $SA^{\tilde{G}}(V,{\cal O})\not=0$ entra\^{\i}ne que $A^{\tilde{G}}(V,{\cal O})\not=0$ ou qu'il existe ${\bf G}'\in {\cal E}(\tilde{G},V)$, avec $G'\not=G$, tel que $SA^{{\bf G}'}(V,{\cal O}_{\tilde{G}'})\not=0$. Si $A^{\tilde{G}}(V,{\cal O})\not=0$,  la d\'efinition de 2.7 entra\^{\i}ne qu'il existe $\gamma\in \tilde{G}_{ss}(F)$ dont la projection dans $\tilde{G}(F_{V})$ appartienne \`a ${\cal O}$ et tel que, pour $v\not\in V$, $\gamma$ soit conjugu\'e \`a un \'el\'ement de $\tilde{K}_{v}$ par un \'el\'ement de $G(F_{v})$. Cette derni\`ere condition et la formule de produit entra\^{\i}ne que $\tilde{H}_{\tilde{G}_{V}}(\gamma)=0$.  Soit ${\bf G}'\in {\cal E}(\tilde{G},V)$, avec $G'\not=G$, supposons  $SA^{{\bf G}'}(V,{\cal O}_{\tilde{G}'})\not=0$. En raisonnant par r\'ecurrence, on peut supposer qu'il existe $\gamma'\in \tilde{G}'_{ss}(F)$ dont la projection  dans $\tilde{G}'(F_{V})$ appartienne \`a ${\cal O}_{\tilde{G}'}$ et tel que $\tilde{H}_{\tilde{G}'_{V}}(\gamma')=0$.  Puisque ${\cal O}_{\tilde{G}'}$ est l'ensemble des classes de conjugaison stable correspondant \`a ${\cal O}$, il existe $\gamma\in {\cal O}$ et un diagramme $(\gamma',B',T',B,T,\gamma)$.   Puisque ${\bf G}'$ est elliptique, les applications $\tilde{H}_{\tilde{G}_{V}}$ et $\tilde{H}_{\tilde{G}'_{V}}$ sont compatibles, c'est-\`a-dire que $\tilde{H}_{\tilde{G}_{V}}(\gamma)=\tilde{H}_{\tilde{G}'_{V}}(\gamma')=0$. $\square$

Soit $\tilde{M}$ un espace de Levi de $\tilde{G}$ et soit ${\cal O}$ une classe de conjugaison stable dans $\tilde{M}_{ss}(F_{V})$. Fixons un syst\`eme de fonctions $B$ comme en 1.10. Rappelons que l'on impose que les valeurs des fonctions $B_{\eta}$ pour $\eta\in \tilde{G}(\bar{F})$ sont des unit\'es hors de $V_{ram}$, a fortiori hors de $V$.
On sait  d\'efinir $S_{\tilde{M}}^{\tilde{G}}(SA^{\tilde{M}}(V,{\cal O}),B,{\bf f})$ pour tout ${\bf f}\in I(\tilde{G}(F_{V}))\otimes Mes(G(F_{V}))$. Si le support de ${\cal O}$ est form\'e d'\'el\'ements $\gamma\in \tilde{M}(F_{V})$  qui sont $\tilde{G}$-\'equisinguliers, c'est-\`a-dire tels  que $M_{\gamma}=G_{\gamma}$, on a
$$S_{\tilde{M}}^{\tilde{G}}(SA^{\tilde{M}}(V,{\cal O}),B,{\bf f})=S_{\tilde{M}}^{\tilde{G}}(SA^{\tilde{M}}(V,{\cal O}),{\bf f}).$$
En effet, cette \'egalit\'e est v\'erifi\'ee par d\'efinition si l'on remplace $SA^{\tilde{M}}(V,{\cal O})$ par n'importe quel \'el\'ement de $D_{g\acute{e}om}(\tilde{M}(F_{V}))\otimes Mes(M(F_{V}))^*$ v\'erifiant la condition de support indiqu\'ee.

\ass{Lemme}{   L'\'egalit\'e
$$S_{\tilde{M}}^{\tilde{G}}(SA^{\tilde{M}}(V,{\cal O}),B,{\bf f})=S_{\tilde{M}}^{\tilde{G}}(SA^{\tilde{M}}(V,{\cal O}),{\bf f})$$
est v\'erifi\'ee pour tout ${\cal O}$.}

Preuve. En utilisant la d\'efinition de 4.1 et en raisonnant par r\'ecurrence, il suffit de prouver l'\'egalit\'e
$$I_{\tilde{M}}^{\tilde{G}}(SA^{\tilde{M}}(V,{\cal O}),B,{\bf f})=I_{\tilde{M}}^{\tilde{G}}(SA^{\tilde{M}}(V,{\cal O}),{\bf f}).$$
En utilisant la d\'efinition de 5.1, cette \'egalit\'e r\'esulte de l'\'egalit\'e
$$(2) \qquad I_{\tilde{M}}^{\tilde{G}}(A^{\tilde{M}}(V,{\cal O}),B,{\bf f})=I_{\tilde{M}}^{\tilde{G}}(A^{\tilde{M}}(V,{\cal O}),{\bf f}),$$
et, pour tout ${\bf M}'=(M',{\cal M}',\zeta)\in {\cal E}(\tilde{M},V)$, avec $M'\not=M$, de l'\'egalit\'e
$$(3) \qquad I_{\tilde{M}}^{\tilde{G}}(transfert(SA^{{\bf M}'}(V,{\cal O}_{\tilde{M}'})),B,{\bf f})=I_{\tilde{M}}^{\tilde{G}}(transfert(SA^{{\bf M}'}(V,{\cal O}_{\tilde{M}'})),{\bf f}).$$
La distribution $A^{\tilde{M}}(V,{\cal O})$ est combinaison lin\'eaire d'int\'egrales orbitales v\'erifiant l'hypoth\`ese du lemme 1.14 et l'\'egalit\'e (2) r\'esulte de ce lemme.  Les membres de gauche et de droite de (3) sont respectivement \'egaux \`a 
 $I_{\tilde{M}}^{\tilde{G},{\cal E}}({\bf M}',SA^{{\bf M}'}(V,{\cal O}_{\tilde{M}'}),B,{\bf f})$ et
 
 \noindent $ I_{\tilde{M}}^{\tilde{G},{\cal E}}({\bf M}',SA^{{\bf M}'}(V,{\cal O}_{\tilde{M}'}),{\bf f})$ d'apr\`es le corollaire 4.8. En utilisant la d\'efinition 4.5(1) de ce dernier terme et la variante ce cette d\'efinition pour le premier, on voit que, pour prouver qu'ils sont \'egaux, il suffit de fixer $s\in \zeta Z(\hat{M})^{\Gamma_{F}}/Z(\hat{G})^{\Gamma_{F}}$ et de prouver l'\'egalit\'e
$$S_{{\bf M}'}^{{\bf G}'(s)}(SA^{{\bf M}'}(V,{\cal O}_{\tilde{M}'}),B,{\bf f}^{{\bf G}'(s)})=S_{{\bf M}'}^{{\bf G}'(s)}(SA^{{\bf M}'}(V,{\cal O}_{\tilde{M}'}),{\bf f}^{{\bf G}'(s)}).$$
C'est l'\'egalit\'e de l'\'enonc\'e que l'on peut appliquer par r\'ecurrence puisque $M'\not=M$ donc $G'(s)\not=G$. $\square$

  \bigskip
  
  \subsection{Les distributions $A^{\tilde{G},{\cal E}}(V,{\cal O},\omega)$ }

Revenons au cas g\'en\'eral.   Pour une classe de conjugaison stable semi-simple ${\cal O}$ dans $\tilde{G}(F_{V})$, on pose
$$A^{\tilde{G},{\cal E}}(V,{\cal O},\omega)=\sum_{{\bf G}'\in {\cal E}(\tilde{G},{\bf a},V)}i(\tilde{G},\tilde{G}')transfert(SA^{{\bf G}'}(V,{\cal O}_{\tilde{G}'})).$$
Cette distribution $A^{\tilde{G},{\cal E}}(V,{\cal O},\omega)$ est un \'el\'ement de $D_{tr-orb}(\tilde{G}(F_{V}),\omega)\otimes Mes(G(F_{V}))^*$. Son support est contenu  dans l'ensemble des \'el\'ements de $\tilde{G}(F_{V})$ dont la partie semi-simple appartient \`a ${\cal O}$. 

{\bf Remarque.} Le cas o\`u $(G,\tilde{G},{\bf a})$ est quasi-d\'eploy\'e et \`a torsion int\'erieure est un cas particulier. Dans ce cas, nos hypoth\`eses de r\'ecurrence ne s'appliquent pas \`a la donn\'ee endoscopique maximale ${\bf G}$. Il convient de remplacer le terme $transfert(SA^{{\bf G}}(V,{\cal O}))$ intervenant dans la somme par $SA^{\tilde{G}}(V,{\cal O})$. On a alors $A^{\tilde{G},{\cal E}}(V,{\cal O})=A^{\tilde{G}}(V,{\cal O})$ par d\'efinition de $SA^{\tilde{G}}(V,{\cal O})$.

\ass{Th\'eor\`eme (\`a prouver)}{Soit ${\cal O}$ une classe de conjugaison stable semi-simple dans $\tilde{G}(F_{V})$. Alors, on a l'\'egalit\'e $A^{\tilde{G},{\cal E}}(V,{\cal O},\omega)=A^{\tilde{G}}(V,{\cal O},\omega)$.}

La d\'efinition ci-dessus s'adapte imm\'ediatement aux $K$-espaces. Pour un tel $K$-espace, une classe de conjugaison stable semi-simple ${\cal O}$ est r\'eunion disjointe de telles classes ${\cal O}_{p}$ pour $p\in \Pi$. On a simplement $A^{K\tilde{G},{\cal E}}(V,{\cal O},\omega)=\oplus_{p\in \Pi}A^{\tilde{G}_{p},{\cal E}}(V,{\cal O}_{p},\omega)$. 

\bigskip

\subsection{Le th\'eor\`eme d'Arthur}
Supposons ici $G=\tilde{G}$, ${\bf a}=1$, $\tilde{K}_{v}=K_{v}$ pour tout $v\not\in V$.

\ass{Th\'eor\`eme}{Sous ces hypoth\`eses, les th\'eor\`emes 5.2 et 5.4 sont v\'erifi\'es.}

C'est  le Global Theorem 1' de [A1].  La preuve que nous donnerons  des th\'eor\`emes 5.2 et 5.4 \'etant directement inspir\'ee de celle d'Arthur, nous pourrions aussi bien  red\'emontrer l'\'enonc\'e ci-dessus. Mais cela n'aurait aucun int\'er\^et. Nous pr\'ef\'erons simplifier un peu la n\^otre en utilisant le r\'esultat d'Arthur.  La propri\'et\'e 5.2(1) des distributions $SA^G(V,{\cal O})$ n'est pas clairement \'enonc\'ee par Arthur, mais est incluse dans sa d\'emonstration. Elle r\'esulte en tout cas de la d\'emonstration plus g\'en\'erale de cette propri\'et\'e qui sera donn\'ee  le moment venu.

 \bigskip

\subsection{Un th\'eor\`eme compl\'ementaire concernant l'endoscopie non standard}
Dans ce paragraphe, on consid\`ere un triplet endoscopique non standard $(G_{1},G_{2},j_{*})$ d\'efini sur $F$. La d\'efinition est la m\^eme que dans le cas local ([W2] 1.7), rappelons-la. Les termes   $G_{1}$ et $G_{2}$  sont des groupes r\'eductifs connexes d\'efinis sur $F$, quasi-d\'eploy\'es et simplement connexes. On consid\`ere $\underline{leurs}$ paires de Borel \'epingl\'ees, dont on note les tores $T_{1}$ et $T_{2}$. Ils sont munis d'actions de $\Gamma_{F}$. Pour $i=1,2$, on note $\Sigma_{i}$ l'ensemble de racines de $T_{i}$ dans $G_{i}$. Pour $\alpha\in \Sigma_{i}$, on note $\check{\alpha}$ la coracine associ\'ee. On pose $X_{i,*,{\mathbb Q}}=X_{*}(T_{i})\otimes_{{\mathbb Z}}{\mathbb Q}$ et $X^*_{i,{\mathbb Q}}=X^*(T_{i})\otimes_{{\mathbb Z}}{\mathbb Q}$.  Le terme $j_{*}$ est  un isomorphisme $j_{*}:X_{1,*,{\mathbb Q}}\to X_{2,*,{\mathbb Q}}$. On note  $j^*:X_{2,{\mathbb Q}}^*\to X_{1,{\mathbb Q}}^*$ l'isomorphisme dual. On suppose

(1) $j_{*}$ est \'equivariant pour les actions de $\Gamma_{F}$.

On suppose qu'il existe une bijection $\tau:\Sigma_{2}\to \Sigma_{1}$ et une fonction $b:\Sigma_{2}\to {\mathbb Q}_{>0}$ telles que

(2) $j^*(\alpha_{2})=b(\alpha_{2})\tau(\alpha_{2})$ pour tout $\alpha_{2}\in \Sigma_{2}$;

(3) $j_{*}(\check{\alpha}_{1})=b(\alpha_{2})\check{\alpha}_{2}$ pour tout $\alpha_{1}\in \Sigma_{1}$, o\`u $\alpha_{2}=\tau^{-1}(\alpha_{1})$. 

Remarquons que $\tau$ et $b$ sont uniquement d\'etermin\'es par ces conditions.
Les groupes de Weyl $W_{1}$ et $W_{2}$ de $G_{1}$ et $G_{2}$ s'identifient, la sym\'etrie relative \`a une racine $\alpha_{2}\in \Sigma_{2}$ s'identifiant \`a celle relative \`a $\tau(\alpha_{2})$. Soit $v$ une place de $F$. Alors $j_{*}$ d\'efinit un isomorphisme 
$$((\mathfrak{t}_{1}(\bar{F}_{v})\cap \mathfrak{g}_{1,reg}(\bar{F}_{v}))/W_{1})^{\Gamma_{F_{v}}}\to ((\mathfrak{t}_{1}(\bar{F}_{v})\cap \mathfrak{g}_{2,reg}(\bar{F}_{v}))/W_{1})^{\Gamma_{F_{v}}}.$$
 Autrement dit une bijection entre classes de conjugaison stable d'\'el\'ements semi-simples r\'eguliers dans $\mathfrak{g}_{1}(F_{v})$ et $\mathfrak{g}_{2}(F_{v})$.  Pour $X_{1}\in \mathfrak{g}_{2,reg}(F_{v})$ et $X_{2}\in \mathfrak{g}_{2,reg}(F_{v})$ dont les classes de conjugaison stable se correspondent,   notons $S_{1}$ et $S_{2}$ leurs commutants dans $G_{1}$ et $G_{2}$. L'isomorphisme $j_{*}$ d\'efinit un isomorphisme $\mathfrak{s}_{1}(F_{v})\to \mathfrak{s}_{2}(F_{v})$. On munit ces alg\`ebres de mesures de Haar se correspondant par cet isomorphisme. Pour $i=1,2$, on munit $S_{i}(F_{v})$ de la mesure de Haar telle que le jacobien de l'exponentielle  vaille $1$ au point $0 \in \mathfrak{s}_{i}(F_{v})$.  On dispose alors de l'int\'egrale orbitale 
$$f_{i}\otimes dg_{i}\mapsto I^{G_{i}}(X_{i},f_{i}\otimes dg_{i})=D^{G_{i}}(X_{i})^{1/2}\int_{S_{i}(F_{v})\backslash G_{i}(F_{v})}f_{i}(ad_{g_{i}^{-1}}(X_{i}))\,dg_{i}$$
sur $C_{c}^{\infty}(\mathfrak{g}_{i}(F_{v}))\otimes Mes(G_{i}(F_{v}))$, puis de l'int\'egrale orbitale stable
$$f_{i}\otimes dg_{i}\mapsto S^{G_{i}}(X_{i},f_{i}\otimes dg_{i}).$$
 Disons que $f_{1}\otimes dg_{1}$ et $f_{2}\otimes dg_{2}$ se correspondent si et seulement si
$$S^{G_{1}}(X_{1},f_{1}\otimes dg_{1})=S^{G_{2}}(X_{2},f_{2}\otimes dg_{2})$$
pour tout couple $(X_{1},X_{2})$ comme ci-dessus. On a

(4) cette correspondance se quotiente en un isomorphisme
$$SI(\mathfrak{g}_{1}(F_{v}))\otimes Mes(G_{1}(F_{v}))\simeq SI(\mathfrak{g}_{2}(F_{v}))\otimes Mes(G_{2}(F_{v})).$$

 Si $v$ est finie, c'est la proposition 1.8(ii) de [W2]. Le cas complexe se ram\`ene au cas r\'eel en rempla\c{c}ant les groupes complexes par les groupes sur ${\mathbb R}$ obtenus par restriction des scalaires. Le cas r\'eel se d\'eduit facilement des r\'esultats de Shelstad comme on l'a vu en [V] 5.1.  
 
 Dualement, on en d\'eduit un isomorphisme
 $$D_{g\acute{e}om}^{st}(\mathfrak{g}_{1}(F_{v}))\otimes Mes(G_{1}(F_{v}))^*\simeq D_{g\acute{e}om}^{st}(\mathfrak{g}_{2}(F_{v}))\otimes Mes(G_{2}(F_{v}))^*.$$
 Pour $i=1,2$, notons $D_{nil}^{st}(\mathfrak{g}_{i}(F_{v}))$ le sous-espace des \'el\'ements de  $D_{g\acute{e}om}^{st}(\mathfrak{g}_{i}(F_{v}))$ \`a support nilpotent. L'isomorphisme ci-dessus se restreint en un isomorphisme
 $$(5) \qquad D_{nil}^{st}(\mathfrak{g}_{1}(F_{v}))\otimes Mes(G_{1}(F_{v}))^*\simeq D_{nil}^{st}(\mathfrak{g}_{2}(F_{v}))\otimes Mes(G_{2}(F_{v}))^*.$$

 Soit $V$ un ensemble fini de places de $F$ tel que
 
 - $V$ contient les places archim\'ediennes de $F$;
 
 - $G_{1}$ et $G_{2}$ sont non ramifi\'es hors de $V$;
 
 - pour $v\not\in V$, notons $p$ la caract\'eristique r\'esiduelle de $F_{v}$ et $e_{v}=[F_{v}:{\mathbb Q}_{p}]$; alors $p>e_{v}N(G_{i})+1$ pour $i=1,2$, o\`u $N(G_{i})$ est l'entier d\'efini en [W2] 4.3;
 
 - les  valeurs de la fonction $b$ sont des unit\'es hors de $V$.
 
 Pour $i=1,2$, on a d\'efini en 5.2 la distribution $SA^{G_{i}}(V,{\cal O})$ associ\'ee \`a une classe de conjugaison stable semi-simple  ${\cal O}\subset G_{i}(F_{V})$. On consid\`ere ici le cas ${\cal O}=\{1\}$ et on note plut\^ot $SA^{G_{i}}_{unip}(V)$ la distribution associ\'ee \`a cette classe. On est ici dans le cas non tordu, cette distribution est relative \`a des sous-groupes compacts hypersp\'eciaux $K_{i,v}$ de $G_{i}(F_{v})$ pour $v\not\in V$. D'apr\`es le th\'eor\`eme d'Arthur cit\'e en 5.5, la condition 5.2(1) est v\'erifi\'ee. Puisque deux sous-groupes compacts hypersp\'eciaux de $G_{i}(F_{v})$ sont toujours conjugu\'es par $G_{i,AD}(F_{v})$, la distribution $SA^{G_{i}}_{unip}(V)$ ne d\'epend pas de ces choix. D'apr\`es le m\^eme th\'eor\`eme, elle est stable. On sait qu'elle est \`a support dans l'ensemble des \'el\'ements unipotents de $G_{i}(F_{V})$. Par l'exponentielle, on la descend en une distribution \`a support nilpotent sur $\mathfrak{g}_{i}(F_{V})$, qui est encore stable. Cela l'identifie \`a un \'el\'ement de $D_{g\acute{e}om}^{st}(\mathfrak{g}_{1}(F_{v}))\otimes Mes(G_{1}(F_{v}))^*$.
 
 \ass{Th\'eor\`eme (\`a prouver)}{Sous ces hypoth\`eses et modulo les identifications ci-dessus, les distributions $SA^{G_{1}}_{unip}(V)$ et $SA^{G_{2}}_{unip}(V)$ se correspondent par le produit tensoriel sur les $v\in V$ des isomorphismes (5).}
 
 Consid\'erons deux triplets endoscopiques non standard $(G_{1},G_{2},j_{*})$ et $(G'_{1},G'_{2},j'_{*})$. On dit qu'ils sont \'equivalents si et seulement s'il existe des isomorphismes $\iota_{i}:G_{i}\to G'_{i}$ d\'efinis sur $F$ et un rationnel $b$ non nul de sorte que le diagramme suivant soit commutatif 
 $$\begin{array}{ccc}X_{1,*,{\mathbb Q}}&\stackrel{j_{*}}{\to}&X_{2,*,{\mathbb Q}}\\ \iota_{1}\downarrow&&\downarrow \iota_{2}\\ X'_{1,*,{\mathbb Q}}&\stackrel{bj'_{*}}{\to}&X'_{2,*,{\mathbb Q}}\\ \end{array}$$
 Comme dans le cas local, on peut classifier tous les triplets possibles.  Appelons triplet \'el\'ementaire un triplet de l'un des types suivants

 (6) $G_{1}=G_{2}$ et $j_{*}$ est l'identit\'e;
 
 (7) quitte \`a \'echanger $G_{1}$ et $G_{2}$, $G_{1}$ est de type $B_{n}$ et $G_{2}$ est de type $C_{n}$, avec $n\geq2$; $j_{*}$ envoie une coracine courte sur  une coracine longue et une coracine longue sur $2$ fois une coracine courte;
 
 (8) $G_{1}$ et $G_{2}$ sont de type $F_{4}$ et $j_{*}$ envoie une coracine courte  sur une coracine longue et une coracine longue sur $2$ fois une coracine courte;
 
 (9) $G_{1}$ et $G_{2}$ sont de type $G_{2}$ (sic!) et $j_{*}$ envoie une coracine courte  sur une coracine longue et une coracine longue sur $3$ fois une coracine courte.
 
 Appelons triplet quasi-\'el\'ementaire un triplet qui se d\'eduit par restriction des scalaires d'un triplet \'el\'ementaire d\'efini sur une extension finie de $F$. Alors tout  triplet endoscopique non standard $(G_{1},G_{2},j_{*})$  est produit de triplets  dont chacun est \'equivalent \`a un triplet quasi-\'el\'ementaire. 
 Evidemment, si un triplet est produit de deux triplets, le th\`eor\`eme pour chacun de ces deux triplets implique celui pour le produit.  On verra au paragraphe suivant que le th\'eor\`eme est insensible \`a une \'equivalence. Il suffit donc de le d\'emontrer pour un triplet quasi-\'el\'ementaire. Remarquons que le cas (6) est insensible \`a une restriction des scalaires et que le th\'eor\`eme est tautologique dans ce cas. Il suffit donc de traiter le cas des triplets d\'eduits par restriction des scalaires de triplets \'el\'ementaires des cas (7), (8) ou (9).

 \bigskip
 
 \subsection{R\'eduction du th\'eor\`eme 5.6}
 On consid\`ere un triplet endoscopique non standard $(G_{1},G_{2},j_{*})$  et un rationnel $b\not=0$. Soit $V$ un ensemble fini de places de $F$ v\'erifiant les conditions du paragraphe pr\'ec\'edent pour $(G_{1},G_{2},j_{*})$. On suppose que
    $b$ est une unit\'e hors de $V$. Alors $V$ v\'erifie aussi ces conditions pour le triplet $(G_{1},G_{2},bj_{*})$. 
 
 \ass{Lemme}{Sous ces hypoth\`eses,  si le th\'eor\`eme 5.5 est v\'erifi\'e pour le triplet $(G_{1},G_{2},j_{*})$, il l'est pour le triplet $(G_{1},G_{2},bj_{*})$.}
 
 Preuve. Il y a deux correspondances entre classes de conjugaison stable semi-simples dans $\mathfrak{g}_{1}(F_{V})$ et $\mathfrak{g}_{2}(F_{V})$, qui sont d\'eduites de $j_{*}$ et de $bj_{*}$. Il est clair que la seconde est la compos\'ee de la premi\`ere et de la correspondance entre classes de conjugaison stable semi-simples dans $\mathfrak{g}_{2}(F_{V})$ d\'eduite de l'homoth\'etie $X\mapsto bX$. Il y a deux isomorphismes 
 $$SI(\mathfrak{g}_{1}(F_{V}))\otimes Mes(G_{1}(F_{V}))\simeq SI(\mathfrak{g}_{2}(F_{V}))\otimes Mes(G_{2}(F_{V})),$$
 qui sont d\'eduits de $j_{*}$ et de $bj_{*}$. Il r\'esulte de ce qui pr\'ec\`ede que le second est le compos\'e du premier et de l'automorphisme de $SI(\mathfrak{g}_{2}(F_{V}))\otimes Mes(G_{2}(F_{V}))$ d\'eduit de l'homoth\'etie $X\mapsto bX$. 
  Puisque $b$ est une unit\'e hors de $V$, on a $\prod_{v\in V}\vert b\vert _{F_{v}}=1$, d'o\`u $D^{\tilde{G}_{2}}(bX)=D^{\tilde{G}_{2}}(X)$ pour tout $X\in \mathfrak{g}_{2}(F_{V})$. De m\^eme, si $X$ est fortement r\'egulier, la multiplication par $b$ pr\'eserve les mesures sur $\mathfrak{s}(F_{V})$, o\`u $S$ est le commutant de $X$. On voit alors que  l'automorphisme ci-dessus de $SI(\mathfrak{g}_{2}(F_{V}))\otimes Mes(G_{2}(F_{V}))$ est d\'eduit de l'automorphisme 
  $f\mapsto f^{b^{-1}}$ de $C_{c}^{\infty}(\mathfrak{g}_{2}(F_{V}))$, o\`u $f^{b^{-1}}(X)=f(b^{-1}X)$.   
 On note encore ${\bf f}\mapsto {\bf f}^{b^{-1}}$  cet automorphisme de $SI(\mathfrak{g}_{2}(F_{V}))\otimes Mes(G_{2}(F_{V}))$. Soient ${\bf f}_{1}\in SI(\mathfrak{g}_{1}(F_{V}))\otimes Mes(G_{1}(F_{V}))$ et ${\bf f}_{2}\in SI(\mathfrak{g}_{2}(F_{V}))\otimes Mes(G_{2}(F_{V}))$ se correspondant par l'isomorphisme d\'eduit de $j_{*}$. Le th\'eor\`eme pour  $(G_{1},G_{2},j_{*})$ affirme que $S^{G_{1}}(SA^{G_{1}}_{unip}(V),{\bf f}_{1})=S^{G_{2}}(SA^{G_{2}}_{unip}(V),{\bf f}_{2})$. Le th\'eor\`eme pour  $(G_{1},G_{2},bj_{*})$ affirme que $S^{G_{1}}(SA^{G_{1}}_{unip}(V),{\bf f}_{1})=S^{G_{2}}(SA^{G_{2}}_{unip}(V),{\bf f}_{2}^{b^{-1}})$. Pour d\'emontrer que ces assertions ont \'equivalentes, il suffit de prouver que $S^{G_{2}}(SA^{G_{2}}_{unip}(V),{\bf f}_{2})=S^{G_{2}}(SA^{G_{2}}_{unip}(V),{\bf f}_{2}^{b^{-1}})$. Seul le groupe $G_{2}$ intervient ici. On peut simplifier la notation en posant $G=G_{2}$ et en supprimant les indices $2$. Toujours pour simplifier, on peut remplacer $b$ par $b^{-1}$ et fixer des mesures de Haar sur les groupes intervenant, ce qui nous d\'ebarrasse des espaces de mesures. Puisque l'automorphisme $ f\mapsto f^{b}$ se rel\`eve \'evidemment en un automorphisme de $I(\mathfrak{g}(F_{V}))\otimes Mes(G(F_{V}))$, on peut d\'emontrer la relation
 
 (1) $I^G(SA^{G}_{unip}(V), f)=I^G(SA^{G}_{unip}(V), f^b)$ pour tout $ f\in I(\mathfrak{g}(F_{V}))$. 
 
 En utilisant la d\'efinition de 5.2, on est ramen\'e \`a prouver

(2) $I^G(A^G_{unip}(V),f)=I^G(A^G_{unip}(V), f^b)$,

\noindent et, pour tout ${\bf G}'\in {\cal E}(G,V)$ avec $G'\not=G$, 

(3)  $I^{{\bf G}'}(SA^{{\bf G}'}_{unip}(V),f^{{\bf G}'})=I^{{\bf G}'}(SA^{{\bf G}'}_{unip}(V),( f^b)^{{\bf G}'})$.
 
Commen\c{c}ons par (3). Soit ${\bf G}'\in {\cal E}(G,V)$ avec $G'\not=G$. Puisqu'on travaille ici avec des alg\`ebres de Lie, les donn\'ees auxiliaires ne jouent gu\`ere de r\^ole et on peut identifier $f^{{\bf G}'}$ \`a une fonction $f^{G'}$ sur $\mathfrak{g}'(F_{V})$. D'apr\`es [F] lemme 3.2.1, il existe un caract\`ere automorphe $\chi$ de ${\mathbb A}_{F}^{\times}/F^{\times}$ tel que $(f^{\lambda})^{G'}=\chi(\lambda)(f^{G'})^{\lambda}$ pour tout $\lambda\in  F_{V}^{\times}$. Le fait que ${\bf G}'$ soit non ramifi\'e hors de $V$ et la d\'efinition de ce caract\`ere $\chi$ entra\^{\i}nent que $\chi$ est non ramifi\'e hors de $V$. Donc, en notant $b_{V}$ et $b^V$ les projections de $b$ dans $F_{V}^{\times}$ et ${\mathbb A}_{F}^{V,\times}$, on a $\chi(b_{V})=\chi(b^V)^{-1}=1$. D'o\`u $(f^b)^{G'}=(f^{G'})^b$. Puisque $G'\not=G$, on peut appliquer (1) par r\'ecurrence en y rempla\c{c}ant $G$ par $G'$. On obtient alors la relation (3).

Pour prouver (2), introduisons  sur $D_{orb}(\mathfrak{g}(F_{V}))$ l'action duale de l'homoth\'etie de rapport $b$. C'est-\`a-dire que, pour $\boldsymbol{\gamma}\in D_{orb}(\mathfrak{g}(F_{V}))$, on note $\boldsymbol{\gamma}^b$ la distribution telle que $I^G(\boldsymbol{\gamma},f^b)=I^G(\boldsymbol{\gamma}^b,f)$ pour tout $f$. Si $\boldsymbol{\gamma}$ est l'int\'egrale orbitale associ\'ee \`a un \'el\'ement $X\in \mathfrak{g}(F_{V})$ et \`a une mesure sur $G_{X}(F_{V})$, $\boldsymbol{\gamma}^b$ est l'int\'egrale orbitale associ\'ee \`a l'\'el\'ement $bX$ et \`a la m\^eme mesure sur $G_{bX}(F_{V})=G_{X}(F_{V})$. Introduisons le syst\`eme de fonctions $B$ sur $G$ ainsi d\'efini: pour tout \'el\'ement semi-simple $\eta\in G$ et toute \'el\'ement $\alpha$ du syst\`eme de racines de $G_{\eta}$, $B_{\eta}(\alpha)=b^{-1}$.  
Soient $M$ un Levi standard de $G$ et $\boldsymbol{\gamma}\in D_{orb}( \mathfrak{m}(F_{V}))$. Montrons que l'on a l'\'egalit\'e

(4) $J_{M}^G(\boldsymbol{\gamma},B,f^b)=J_{M}^G(\boldsymbol{\gamma}^b,f)$ pour tout $f$.

Les formules de descente habituelles nous ram\`enent au cas local. C'est alors l'assertion [III] 6.8(5). Le corps local \'etait non-archim\'edien dans cette r\'ef\'erence, mais cela n'importait pas pour cette assertion.

 Revenons \`a la d\'efinition de $A^G_{unip}(V)$.  Dans l'\'egalit\'e (2) \`a prouver, $f$ est une fonction sur $\mathfrak{g}(F_{V})$. Notons-la plut\^ot $f_{V}$ et notons comme ci-dessus $b_{V}$ la projection de $b$ dans $F_{V}$. La fonction not\'ee pr\'ec\'edemment $f^b$ s'\'ecrit maintenant $(f_{V})^{b_{V}}$. Par l'exponentielle, on rel\`eve ces deux fonctions en des fonctions d\'efinies sur un voisinage de $1$ dans $G(F_{V})$ invariant par conjugaison. On continue de noter ces fonctions $f_{V}$ et $(f_{V})^{b_{V}}$. On les compl\`ete globalement en des fonctions $\dot{f}={\bf 1}_{K^V}\otimes  f_{V}$, $\dot{f}^{b_{V}}={\bf 1}_{K^V}\otimes (f_{V})^{b_{V}}$. Montrons  que l'on a l'\'egalit\'e
 
 (5) $J_{unip}^G(\dot{f}^{b_{V}})=J_{unip}^G(\dot{f})$.
 
 D'apr\`es les d\'efinitions de 2.1, il suffit de fixer un sous-groupe parabolique standard $P=MU_{P}$ de $G$ et de prouver l'\'egalit\'e
 $$K_{P,unip}(\dot{f}^{b_{V}},g)=K_{P,unip}(\dot{f},g)$$
 pour tout $g\in G({\mathbb A}_{F})$. Rappelons que
 $$K_{P,unip}(\dot{f}^{b_{V}},g)=\int_{U_{P}(F)\backslash U_{P}({\mathbb A}_{F})}\sum_{\gamma\in P_{unip}(F)}\dot{f}^{b_{V}}(g^{-1}\gamma ug)\,du=\int_{U_{P}({\mathbb A}_{F})}\sum_{\gamma\in M_{unip}(F)}\dot{f}^{b_{V}}(g^{-1}\gamma ug)\,du,$$
 o\`u $P_{unip}$ et $M_{unip}$ sont les ensembles d'\'el\'ements unipotents dans $P$ et $M$. En introduisant l'ensemble $\mathfrak{m}_{nil}$ des \'el\'ements nilpotents de $\mathfrak{m}$, on a 
 $$K_{P,unip}(\dot{f}^{b_{V}},g)=\int_{U_{P}({\mathbb A}_{F})}\sum_{X\in \mathfrak{m}_{nil}(F)}\dot{f}^{b_{V}}(g^{-1}exp(X)ug)\, du.$$
 On descend ais\'ement l'int\'egrale \`a l'alg\`ebre de Lie et on obtient
$$K_{P,unip}(\dot{f}^{b_{V}},g)=\int_{\mathfrak{u}_{P}({\mathbb A}_{F})}\sum_{X\in \mathfrak{m}_{nil}(F)}\dot{f}^{b_{V}}(g^{-1}exp(X+N)g)\, dN.$$
Puisque $b\in {\mathbb Q}^{\times}$, on peut remplacer $X$ par $b^{-1}X$ et $N$ par $b^{-1}N$. En d\'ecomposant les int\'egrales selon les places dans $v$ et celles hors de $V$, on obtient
$$K_{P,unip}(\dot{f}^{b_{V}},g)=\sum_{X\in \mathfrak{m}_{nil}(F)}C_{V}((f_{V})^{b_{V}},b^{-1}X)C^V({\bf 1}_{K^V},b^{-1}X),$$
o\`u
$$C_{V}((f_{V})^{b_{V}},b^{-1}X)=\int_{u_{P}(F_{V})}(f_{V})^{b_{V}}(g^{-1}(b^{-1}X+b^{-1}N)g)\, dN$$
et
$$C^V({\bf 1}_{K^V},b^{-1}X)=\int_{u_{P}({\mathbb A}_{F}^V)}{\bf 1}_{K^V}(g^{-1}exp(b^{-1}X+b^{-1}N)g)\,dN.$$
Il r\'esulte de la d\'efinition de $(f_{V})^{b_{V}}$ que l'on a l'\'egalit\'e
$$C_{V}((f_{V})^{b_{V}},b^{-1}X)=C_{V}(f_{V},X).$$
Si $v$ est une place hors de $V$, on a $v\not\in V_{ram}$ et la caract\'eristique r\'esiduelle est "grande". Cela assure la propri\'et\'e suivante. Notons $\mathfrak{o}_{v}$ l'anneau des entiers de $F_{v}$. Au groupe hypersp\'ecial $K_{v}$ de $G(F_{v})$ est associ\'e une $\mathfrak{o}_{v}$-alg\`ebre de Lie $\mathfrak{k}_{v}$. Alors, pour tout \'el\'ement nilpotent $N\in \mathfrak{g}(F_{v})$, on a $exp(N)\in K_{v}$ si et seulement si $N\in \mathfrak{k}_{v}$. Puisque $b$ est une unit\'e en $v$, on a alors l'\'egalit\'e
 $${\bf 1}_{K_{v}}(g^{-1}exp(b^{-1}X+b^{-1}N)g)={\bf 1}_{K_{v}}(g^{-1}exp(X+N)g)$$
 pour tous $X$, $N$, $g$ intervenant ci-dessus. D'o\`u l'\'egalit\'e
 $$C^V({\bf 1}_{K^V},b^{-1}X)=C^V({\bf 1}_{K^V},X),$$
 puis 
$$K_{P,unip}(\dot{f}^{b_{V}},g)=\sum_{X\in \mathfrak{m}_{nil}(F)}C_{V}(f_{V},X)C^V({\bf 1}_{K^V},X).$$
A partir de l\`a, le m\^eme calcul que ci-dessus, en sens inverse, conduit \`a l'\'egalit\'e 
$$K_{P,unip}(\dot{f}^{b_{V}},g)=K_{P,unip}(\dot{f},g)$$
cherch\'ee. D'o\`u (5).

 Par d\'efinition, on a
 $$(6) \qquad I^G(A^G_{unip}(V),(f_{V})^{b_{V}})=J^G_{unip}(\dot{f}^{b_{V}})-\sum_{M\in {\cal L}(M_{0}),M\not=G}\vert W^M\vert \vert W^G\vert ^{-1}J_{M}^G(A^M_{unip}(V),(f_{V})^{b_{V}}).$$
 En appliquant (5), le premier terme devient $J^G_{unip}(\dot{f})$. Soit $M$ un Levi tel que $M\not=G$.  On est dans la situation o\`u le lemme 1.14 s'applique. Ce lemme est \'enonc\'e pour les int\'egrales orbitales pond\'er\'ees invariantes, mais sa preuve s'applique aussi bien \`a leurs  versions non invariantes. D'o\`u
$$J_{M}^G(A^M_{unip}(V),(f_{V})^{b_{V}})=J_{M}^G(A^M_{unip}(V),B,(f_{V})^{b_{V}}).$$ 
En appliquant (4), on obtient
$$J_{M}^G(A^M_{unip}(V),(f_{V})^{b_{V}})=J_{M}^G(A^M_{unip}(V)^{b_{V}},f_{V}).$$
  En raisonnant comme toujours par r\'ecurrence, on peut supposer l'analogue de (2) connu si l'on remplace $G$ par $M$. Cette relation \'equivaut \`a l'\'egalit\'e $A^M_{unip}(V)^{b_{V}}=A^M_{unip}(V)$. D'o\`u encore 
 $$J_{M}^G(A^M_{unip}(V),(f_{V})^{b_{V}})=J_{M}^G(A^M_{unip}(V),f_{V}).$$
 Mais alors le membre de droite de (6) est \'egal \`a la m\^eme expression o\`u $(f_{V})^{b_{V}}$ est remplac\'ee par $f_{V}$. D'o\`u l'\'egalit\'e
 $$ I^G(A^G_{unip}(V),(f_{V})^{b_{V}})= I^G(A^G_{unip}(V),f_{V}),$$
 c'est-\`a-dire (2). Cela ach\`eve la d\'emonstration. $\square$
 
   $\square$
   
   \bigskip
   
   \subsection{Insertion du th\'eor\`eme 5.6 dans les hypoth\`eses de r\'ecurrence}
   
   En [III] 6.1, on a associ\'e un entier $N(G_{1},G_{2},j_{*})$ \`a tout triplet endoscopique non standard $(G_{1},G_{2},j_{*})$. On a aussi d\'efini des triplets particuliers $(G,\tilde{G},{\bf a})$ en [III] 6.3. Les constructions de cette r\'ef\'erence valent sur notre corps de nombres $F$. 
   
   Les hypoth\`eses de r\'ecurrence pos\'ees en 1.1 doivent \^etre compl\'et\'ees de la fa\c{c}on suivante. Pour d\'emontrer une assertion concernant un triplet $(G,\tilde{G},{\bf a})$ quasi-d\'eploy\'e et \`a torsion int\'erieure, on ne se soucie pas du th\'eor\`eme 5.6 qu'on n'utilisera pas dans ce cas. Pour d\'emontrer une assertion concernant l'un des triplets particuliers $(G,\tilde{G},{\bf a})$ de [III] 6.3, on suppose connu le th\'eor\`eme 5.6 pour tout triplet $(G_{1},G_{2},j_{*})$ tel que $N(G_{1},G_{2},j_{*})< dim(G_{SC})$. Pour d\'emontrer une assertion concernant un triplet $(G,\tilde{G},{\bf a})$  qui n'est pas l'un de ces triplets particuliers et qui n'est pas quasi-d\'eploy\'e et \`a torsion int\'erieure, on suppose connu le th\'eor\`eme 5.6 pour tout triplet $(G_{1},G_{2},j_{*})$ tel que $N(G_{1},G_{2},j_{*})\leq dim(G_{SC})$. Pour d\'emontrer une assertion concernant un triplet 
$(G_{1},G_{2},j_{*})$, on suppose connu le th\'eor\`eme 5.6 pour tout triplet $(G'_{1},G'_{2},j'_{*})$ tel que $N(G'_{1},G'_{2},j'_{*})< N(G_{1},G_{2},j_{*})$. On suppose connus tous les r\'esultats concernant les triplets $(G,\tilde{G},{\bf a})$ quelconques tels que $dim(G_{SC})< N(G_{1},G_{2},j_{*})$. On suppose connus tous les r\'esultats concernant les triplets $(G,\tilde{G},{\bf a})$ qui v\'erifient les deux conditions suivantes:

- ils sont quasi-d\'eploy\'es et \`a torsion int\'erieure ou ils font partie des triplets particuliers de [III] 6.3;

- on a $dim(G_{SC})=N(G_{1},G_{2},j_{*})$.

\bigskip

\subsection{La formule stable}

   Supposons  $(G,\tilde{G},{\bf a})$ quasi-d\'eploy\'e et \`a torsion int\'erieure. Soit $V$ un ensemble fini de places contenant $V_{ram}$. On note $\tilde{G}_{ss}(F_{V})/st-conj$ l'ensemble des classes de conjugaison stable semi-simples dans $\tilde{G}(F_{V})$.  Pour ${\bf f}\in SI(\tilde{G}(F_{V}))\otimes Mes(G(F_{V}))$, on pose
$$S^{\tilde{G}}_{g\acute{e}om}({\bf f})=\sum_{\tilde{M}\in {\cal L}(\tilde{M}_{0})}\vert W^M\vert \vert W^G\vert ^{-1}\sum_{{\cal O}\in \tilde{M}_{ss}(F_{V})/st-conj} S_{\tilde{M}}^{\tilde{G}}(SA^{\tilde{M}}({\cal O},V),{\bf f}).$$

\ass{Lemme}{Cette somme est finie. }

Preuve. On peut \'evidemment fixer $\tilde{M}$. En utilisant 5.3(1) et le lemme 4.3, on voit qu'il existe un sous-ensemble compact $\tilde{C}_{V}$ de $\tilde{M}(F_{V})$ tel que, si  $S_{\tilde{M}}^{\tilde{G}}(SA^{\tilde{M}}({\cal O},V),{\bf f}')\not=0$, alors ${\cal O}$ coupe $\tilde{C}_{V}$. Il reste \`a prouver

(1) il n'y a qu'un nombre fini de ${\cal O}\in \tilde{M}_{ss}(F_{V})/st-conj$ tels que $SA^{\tilde{M}}({\cal O},V)\not=0$ et ${\cal O}$ coupe $\tilde{C}_{V}$.

On peut aussi bien supposer $\tilde{M}=\tilde{G}$. On utilise la d\'efinition de 5.2. Si $SA^{\tilde{G}}({\cal O},V)\not=0$, alors $A^{\tilde{G}}({\cal O},V)\not=0$ ou il existe ${\bf G}'\in {\cal E}(\tilde{G},V)$, avec $G'\not=G$, tel que $SA^{{\bf G}'}({\cal O}_{\tilde{G}'},V)\not=0$. Dans le premier cas, il existe $\gamma\in \tilde{G}_{ss}(F)$ dont la projection dans $\tilde{G}(F_{V})$ appartient \`a ${\cal O}$ et tel que, pour tout $v\not\in V$, $\gamma$ soit conjugu\'e \`a un \'el\'ement de $\tilde{K}_{v}$ par un \'el\'ement de $G(F_{v})$. En imposant la condition que ${\cal O}$ coupe $\tilde{C}_{V}$, l'ensemble de ces $\gamma$ forme un nombre fini de classes de conjugaison par $G(F)$ (lemme 2.1). A fortiori l'ensemble des classes de conjugaison stable ${\cal O}$ contenant un tel \'el\'ement est fini. Dans le deuxi\`eme cas, le compact $\tilde{C}_{V}$ d\'etermine un compact $\tilde{C}_{V,\tilde{G}'}$ de $\tilde{G}'(F_{V})$ tel que la condition que ${\cal O}$ coupe $\tilde{C}_{V}$ entra\^{\i}ne que ${\cal O}_{\tilde{G}'}$ coupe $\tilde{C}_{V,\tilde{G}'}$. Puisque $G'\not=G$, on peut appliquer (1) par r\'ecurrence: l'ensemble des ${\cal O}_{\tilde{G}'}\in \tilde{G}'_{ss}(F_{V})/st-conj$ tels que $SA^{{\bf G}'}({\cal O}_{\tilde{G}'},V)\not=0$ et ${\cal O}_{\tilde{G}'}$ coupe $\tilde{C}_{V,\tilde{G}'}$ est fini. L'ensemble des classes ${\cal O}$ qui leur correspondent l'est aussi. $\square$

\bigskip

\subsection{Le th\'eor\`eme principal}

Levons l'hypoth\`ese que $(G,\tilde{G},{\bf a})$ est quasi-d\'eploy\'e et \`a torsion int\'erieure et travaillons plut\^ot avec un $K$-triplet $(KG,K\tilde{G},{\bf a})$. Soit $V$ un ensemble fini de places de $F$ contenant $V_{ram}$. Par un proc\'ed\'e formel familier, de la d\'efinition du paragraphe pr\'ec\'edent se d\'eduit celle des termes $S_{g\acute{e}om}^{{\bf G}'}({\bf f}^{{\bf G}'})$ apparaissant dans l'\'enonc\'e suivant.

\ass{Th\'eor\`eme (\`a prouver)}{ Pour tout ${\bf f}\in C_{c}^{\infty}(K\tilde{G}(F_{V}),\omega)\otimes Mes(G(F_{V}))$, on a l'\'egalit\'e
$$I^{K\tilde{G}}_{g\acute{e}om}({\bf f},\omega)=\sum_{{\bf G}'\in {\cal E}(\tilde{G},{\bf a},V)} i(\tilde{G},\tilde{G}')S^{{\bf G}'}_{g\acute{e}om}({\bf f}^{{\bf G}'}).$$}

On montrera dans la section 6 que ce th\'eor\`eme r\'esulte assez facilement des autres th\'eor\`emes pr\'ec\'edemment \'enonc\'es. 

\bigskip

 \section{Preuve conditionnelle du th\'eor\`eme 5.10}

 \bigskip
\subsection{Rappel}
On consid\`ere un triplet $(G,\tilde{G},{\bf a})$ et une donn\'ee endoscopique ${\bf G}'$. 
On rappelle le lemme suivant qui se trouve d\'ej\`a dans [A11].

\ass{Lemme}{Les fl\`eches naturelles, $ker^1(F,Z(\hat{G}))\rightarrow ker^1(F,\hat{T})$ et 
$$ker^1(F,Z(\hat{G}'))\rightarrow ker^1(F,\hat{T}^{\hat{\theta},0})$$
sont des isomorphismes.}

La preuve suit [K2] preuve du lemme 4.3.2 et [A11] lemme 2: ces deux fl\`eches sont \'evidemment analogues puisque $\hat{T}^{\hat{\theta},0}$ est un tore maximal de $\hat{G'}$. On d\'emontre donc la premi\`ere.

On a la suite exacte:
$1\rightarrow Z(\hat{G}) \rightarrow \hat{T}\rightarrow \hat{T}/Z(\hat{G}) \rightarrow 1
$ dont se d\'eduit une suite exacte de groupes de cohomologie. Le tore $\hat{T}/Z(\hat{G})$ est induit c'est-\`a-dire que son groupe des caract\`eres a une base sur laquelle $\Gamma_{F}$ agit; ainsi  $(\hat{T}/Z(\hat{G}))^{\Gamma_{F}}$ est connexe et  la fl\`eche de l'\'enonc\'e est injective. Pour la surjectivit\'e, il suffit de remarquer que $ker^1(F, (\hat{T}/Z(\hat{G})))=0$: en effet  on se ram\`ene au cas o\`u $\Gamma_{F}$ agit transitivement sur une base des caract\`eres de $\hat{T}/Z(\hat{G})$. Dans ce cas $H^1(W_{F}, (\hat{T}/Z(\hat{G})))$ s'identifie \`a $H^1(W_{F'}, {\mathbb C}^*)$ o\`u $F'$ est une extension galoisienne de $F$ d\'etermin\'e par le sous-groupe de $\Gamma_{F}$ stabilisant un \'el\'ement de la base des caract\`eres avec l'action triviale de $\Gamma_{F'}$ sur ${\mathbb C}^*$; ainsi ce groupe de cohomologie n'est autre que le groupe des caract\`eres de $W_{F'}$ (continus \`a valeurs complexes). Un tel caract\`ere correspond \`a un \'el\'ement du  sous-groupe $ker^1(F, (\hat{T}/Z(\hat{G})))$ si le localis\'e du caract\`ere en toute place $v$ est trivial; le caract\`ere est alors n\'ecessairement trivial.

\bigskip
\subsection{Au sujet des constantes}
Pour $\tilde{G}$ un espace tordu, on note $j(\tilde{G}):= \vert det_{{\mathcal A}_{G}/{\mathcal A}_{\tilde{G}}} (1-{\theta})\vert$ et, pour une donn\'ee endoscopique elliptique ${\bf G}'$, 
on a pos\'e en  5.1,
$i({\tilde{G},\tilde{G}'})=$
$$j(\tilde{G})^{-1}\vert \pi_{0} (Aut_{\hat{G}}({\bf {G}'})/\hat{G}')\vert^{-1}\vert ker^1(F,Z(\hat{G}'))\vert \vert ker^1(F,Z(\hat{G}))\vert^{-1} \delta(\tilde{G},\tilde{G}'),$$ o\`u
$
\delta(\tilde{G},\tilde{G}')=\vert \pi_{0}(Z(\hat{G})^{\Gamma_{F}})\vert\vert \pi_{0}(Z(\hat{G}')^{\Gamma_{F}})\vert^{-1} \vert \pi_{0}(Z(\hat{G})^{\Gamma_{F},0}\cap Z(\hat{G}'))\vert.
$
 On a pr\'ecis\'e ici la notation $Aut({\bf G}')$ de 5.1 en $Aut_{\hat{G}}({\bf {G}'})$. Avec la description de  5.1, la composante neutre de $Aut_{\hat{G}}({\bf {G}'})/\hat{G}'$ est exactement l'image de $Z(\hat{G})^{\Gamma_{F},0}$ dans ce groupe d'automorphismes. Ainsi $$\vert \pi_{0}(Aut_{\hat{G}}({\bf {G}'}))\vert= \vert Aut_{\hat{G}}({\bf {G}'})/Z(\hat{G})^{\Gamma_{F},0} \hat{G}'\vert.$$
 On pose $\overline{Aut}_{\hat{G}}({\bf {G}'})= {Aut}_{\hat{G}}({\bf {G}'})/\hat{G}'Z(\hat{G})^{\Gamma_{F}}$ qui est donc un groupe fini. On fixe aussi des espaces de Levi de $\tilde{G}$ et $\tilde{G}'$, not\'es $\tilde{M}$ et $\tilde{M}'$.   On suppose que $\tilde{M}'$ est un espace endoscopique de $\tilde{M}$.   On a alors d\'efini en  4.4, la constante
$$
i_{\tilde{M}'}(\tilde{G},\tilde{G}'):= \vert Z(\hat{G}')^{\Gamma_{F}}/Z(\hat{G})^{\Gamma_{F}}\cap Z(\hat{G}')^{\Gamma_{F}}\vert^{-1} \vert Z(\hat{M}')^{\Gamma_{F}}/Z(\hat{M})^{\Gamma_{F}}\cap Z(\hat{M}')^{\Gamma_{F}}\vert.
$$ Ici on modifie cette constante car au lieu de sommer \`a l'int\'erieur d'une classe sous $Z(\hat{M})^{\Gamma_{F},\hat{\theta}}/Z(\hat{G})^{\Gamma_{F},\hat{\theta}}$ on va sommer sur une classe sous $
Z(\hat{M})^{\Gamma_{F}}Z(\hat{G})/Z(\hat{G})(1-\hat{\theta})(Z(\hat{M})^{\Gamma_{F}})$ (cf. [I] 3.3 (2)) et on pose
$
i'_{\tilde{M}'}(\tilde{G},\tilde{G}'):=$
$$j(\tilde{G})^{-1}j(\tilde{M}) \vert Z(\hat{G}')^{\Gamma_{F}}/Z(\hat{G})^{\Gamma_{F}}\cap Z(\hat{G}')^{\Gamma_{F}}\vert^{-1} \vert Z(\hat{M}')^{\Gamma_{F}}/Z(\hat{M})^{\Gamma_{F}}\cap Z(\hat{M}')^{\Gamma_{F}}\vert.
$$
  et le but de ce paragraphe est de montrer la proposition suivante:
  \ass{Proposition}{ Soient $\tilde{G},{\bf {G}'},\tilde{M}, {\bf {M}'}$; on a
$$
i(\tilde{G},\tilde{G}')i'_{\tilde{M}'}(\tilde{G},\tilde{G}')^{-1}i(\tilde{M},\tilde{M}')^{-1}=\vert \overline{Aut}_{\hat{G}}({\bf {G}'})\vert^{-1}\vert \overline{Aut}_{\hat{M}}({\bf {M}'})\vert.
$$}
 
 Il est clair que les $j(?)$ se compensent.
 On r\'ecrit diff\'eremment $\delta(\tilde{G},\tilde{G}')$: par ellipticit\'e, $Z(\hat{G}')^{\Gamma_{F},0}= Z(\hat{G})^{\Gamma_{F},\hat{\theta},0}$, d'o\`u l'inclusion $Z(\hat{G}')^{\Gamma_{F},0}\subset Z(\hat{G})^{\Gamma_{F},0}$. Ainsi $$\vert \pi_{0}(Z(\hat{G}')^{\Gamma_{F}})\vert^{-1} \vert 
 \pi_{0}(Z(\hat{G})^{\Gamma_{F},0}\cap Z(\hat{G}'))\vert=
 \vert Z(\hat{G}')/Z(\hat{G})^{\Gamma_{F},0}\cap Z(\hat{G}')\vert^{-1}.
 $$Et
$
\delta(\tilde{G},\tilde{G}')= \vert Z(\hat{G})^{\Gamma_{F}}/Z(\hat{G})^{\Gamma_{F},0}\vert \vert 
Z(\hat{G}')^{\Gamma_{F}}/Z(\hat{G})^{\Gamma_{F},0}\cap Z(\hat{G}')^{\Gamma_{F}}\vert^{-1}.
$ On consid\`ere la suite exacte:
$$
1\rightarrow Z(\hat{G})^{\Gamma_{F}}\cap Z(\hat{G}')^{\Gamma_{F}}/Z(\hat{G})^{\Gamma_{F},0}\cap Z(\hat{G}')^{\Gamma_{F}}\rightarrow Z(\hat{G})^{\Gamma_{F}}/Z(\hat{G})^{\Gamma_{F},0} $$
$$\rightarrow Z(\hat{G})^{\Gamma_{F}}/Z(\hat{G})^{\Gamma_{F},0} (Z(\hat{G})^{\Gamma_{F}}\cap Z(\hat{G}')^{\Gamma_{F}}).$$
Et la suite exacte:
$$
1\rightarrow Z(\hat{G})^{\Gamma_{F}}\cap Z(\hat{G}')^{\Gamma_{F}}/Z(\hat{G})^{\Gamma_{F},0}\cap Z(\hat{G}')^{\Gamma_{F}}\rightarrow Z(\hat{G}')^{\Gamma_{F}}/Z(\hat{G})^{\Gamma_{F},0}\cap Z(\hat{G}')^{\Gamma_{F}}$$
$$\rightarrow Z(\hat{G}')^{\Gamma_{F}}/Z(\hat{G})^{\Gamma_{F}}\cap Z(\hat{G}')^{\Gamma_{F}}.
$$
Et on obtient $\delta(\tilde{G},\tilde{G}')=$
$$ \vert Z(\hat{G})^{\Gamma_{F}}/Z(\hat{G})^{\Gamma_{F},0}(Z(\hat{G}')^{\Gamma_{F}}\cap Z(\hat{G})^{\Gamma_{F}})\vert
\vert Z(\hat{G}')^{\Gamma_{F}}/Z(\hat{G})^{\Gamma_{F}}\cap Z(\hat{G}')^{\Gamma_{F}}\vert^{-1}.
$$
On remarque que $Z(\hat{G})^{\Gamma_{F},0}(Z(\hat{G})^{\Gamma_{F}}\cap Z(\hat{G}')^{\Gamma_{F}})=
Z(\hat{G})^{\Gamma_{F}}\cap Z(\hat{G})^{\Gamma_{F},0} \hat{G}'
$. Ainsi
$$
 \vert Aut_{\hat{G}}({\bf {G}'})/Z(\hat{G})^{\Gamma_{F},0} \hat{G}'\vert  \vert Z(\hat{G})^{\Gamma_{F}}/Z(\hat{G})^{\Gamma_{F},0}(Z(\hat{G}')^{\Gamma_{F}}\cap Z(\hat{G})^{\Gamma_{F}})\vert^{-1}=
 $$
 $$
 \vert \overline{Aut}_{\hat{G}}({\bf {G}'})\vert.
 $$
Et
 $$
 i(\tilde{G},\tilde{G}')=\vert \overline{Aut}_{\hat{G}}({\bf {G}'})\vert^{-1} \vert Z(\hat{G}')^{\Gamma_{F}}/Z(\hat{G})^{\Gamma_{F}}\cap Z(\hat{G}')^{\Gamma_{F}}\vert^{-1} $$
 $$\times \vert ker^1(F,Z(\hat{G}'))\vert \vert ker^1(F,Z(\hat{G}))\vert^{-1}
 $$
 $$
 =\vert \overline{Aut}_{\hat{G}}({\bf {G}'})\vert^{-1} i_{\tilde{M}'}(G,G')^{-1}  \vert Z(\hat{M}')^{\Gamma_{F}}/
 Z(\hat{M}')^{\Gamma_{F}}\cap Z(\hat{M}')^{\Gamma_{F}}\vert^{-1}
 $$
 $$
 \times \vert ker^1(F,Z(\hat{M}'))\vert \vert ker^1(F,Z(\hat{M}))\vert^{-1}=
 $$
 $$
 \vert \overline{Aut}_{\hat{G}}({\bf {G}'})\vert^{-1} i_{\tilde{M}'}(G,G')^{-1}  \vert \overline{Aut}_{\hat{M}}({\bf {M}'})\vert i(\tilde{M},\tilde{M}'),
 $$
 ce qui est l'assertion cherch\'ee.
 
 \bigskip
\subsection{Combinatoire des sommes}
On donnera en 6.5 un analogue du lemme 10.2 de [A1]. Auparavant, il faut rappeler que, si ${\bf G}'$ est une donn\'ee endoscopique elliptique de $(G,\tilde{G},{\bf a})$ et si ${\bf M}'$ est un espace de Levi de ${\bf G}'$ (en un sens compr\'ehensible), il ne correspond pas forc\'ement \`a $\tilde{M}'$ un espace de Levi $\tilde{M}$ de $\tilde{G}$. Mais il lui correspond un groupe de Levi $\hat{M}$ de $\hat{G}$  tel qu'il existe un sous-groupe parabolique $\hat{P}\in {\cal P}(\hat{M})$ qui soit stable par l'action galoisienne et par $\hat{\theta}$, cf. [I] 3.4.  Dans la suite, $\hat{M}$ est suppos\'e v\'erifier cette propri\'et\'e. On dira que ${\bf M}'$ est une donn\'ee endoscopique elliptique de $\hat{M}$. Les constantes d\'efinies dans le paragraphe pr\'ec\'edent sont encore d\'efinies dans cette situation: $\tilde{M}$ n'y intervenait que via $\hat{M}$. On les utilise en rempla\c{c}ant $\tilde{M}$ par $\hat{M}$ dans les notations.

On fixe une fonction sur l'ensemble des triplets ${\bf {G}'}, {\bf {M}'}, \hat{M}$ form\'e d'un espace endoscopique elliptique de $\tilde{G}$ et d'un espace de Levi de cet espace endoscopique,  not\'ee $S({\bf {G}'}, {\bf {M}'})$ et d'un espace de Levi de $\hat{G}$ (v\'erifiant la condition ci-dessus) tel que ${\bf {M}'}$ en soit un espace endoscopique elliptique et on suppose que cette fonction est invariante sous l'action par conjugaison de $\hat{G}$. On va sommer de deux fa\c{c}ons diff\'erentes cette fonction (avec des coefficients) sur l'ensemble des triplets ${\bf {G}'}, {\bf {M}'}, \hat{M}$ modulo conjugaison sous $\hat{G}$; il faut pr\'eciser quelques notations.

Les triplets consid\'er\'es sont form\'es d'un espace de Levi $\hat{M}$ de $\hat{G}$ et d'un couple $\bf{G}'$,  $\bf{M}'$, o\`u $\bf{M}'$ est un espace de Levi de la donn\'ee endoscopique $\bf{G}'$  et o\`u $\bf{G}'$ est une donn\'ee endoscopique elliptique de $\tilde{G}$ tandis que $\bf{M}'$ est une donn\'ee endoscopique elliptique de $\hat{M}$; donc en particulier, dans la donn\'ee endoscopique ${\bf {G}'}$, on a un \'el\'ement $\tilde{s}_{\bf{G}'}\in \hat{G}\hat{\theta}/Z(\hat{G})$ et dans la donn\'ee endoscopique ${\bf M}'$, on a un \'el\'ement $\tilde{s}_{\bf{M}'}\in \hat{M}\hat{\theta}/Z(\hat{M})$. Et on a n\'ecessairement $\tilde{s}_{\bf{G}'}= \tilde{s}_{\bf{M}'} Z(\hat{M})$.  En tenant compte de [I] 3.2 (1), on impose (ce qui est loisible) \`a $\tilde{s}_{\bf{M}'}$ d'\^etre dans la classe canonique sous $Z(\hat{M})^{\Gamma_{F}}Z(\hat{G})$ d\'efinie en loc. cite, \`a l'int\'erieur de sa classe sous $Z(\hat{M})$. Ainsi $\tilde{s}_{\bf{G}'}\in \tilde{s}_{\bf{M}'}Z(\hat{M})^{\Gamma_{F}}Z(\hat{G})/Z(\hat{G})$.

On a besoin de remarquer que $Aut_{\hat{M}}(\bf{M}')$ agit par conjugaison sur $\tilde{s}_{\bf{M}'}$ en laissant stable sa classe modulo $Z(\hat{M})^{\Gamma_{F}}Z(\hat{G})$. En fait cela r\'esulte d'un r\'esultat g\'en\'eral prouv\'e dans le paragraphe suivant, que l'on applique avec $\tilde{G},{\bf{G}'}$ remplac\'e par $\hat{M}, {\bf M}'$.

\bigskip
\subsection{Remarque sur l'action des groupes d'automorphismes de donn\'ees endoscopiques}

\ass{Remarque}{Soit $\bf{G}'$ une donn\'ee endoscopique (non n\'ecessairement elliptique) de $\tilde{G}$ d'o\`u en particulier un \'el\'ement $\tilde{s}_{\bf{G}'}\in \hat{G}\hat{\theta}$. Alors pour tout $x\in Aut_{\hat{G}}(\bf{G}')$, $x \tilde{s}_{\bf{G}'} x^{-1}\in \tilde{s}_{\bf{G}'} Z(\hat{G})^{\Gamma_{F}}$.}

On note $z$ l'\'el\'ement de $Z(\hat{G})$ tel que $x\tilde{s}_{\bf{G}'}x^{-1}=z\tilde{s}_{\bf{G}'}$; pour faire agir $\Gamma_{F}$ sur $z$, on utilise les \'el\'ements de ${\mathcal G}'$: pour tout $w\in \Gamma_{F}$, on fixe $h_{w}\in {\mathcal G}'$ dont l'image dans la projection de ${\mathcal G}'$ sur $\Gamma_{F}$ est $w$. On a $h_{w}zh_{w}^{-1}=w z w^{-1}$. On conjugue cette \'egalit\'e par $x$:
$$
x h_{w}x^{-1} z\,  xh_{w}^{-1}x^{-1}= x (w z w^{-1})x^{-1}=w zw^{-1},
$$
car $wzw^{-1}\in Z(\hat{G})$ et $x\in \hat{G}$.  On utilise le fait que $h_{w}$ commute \`a $\tilde{s}_{{\bf G}'}$ \`a un cocycle pr\`es d'apr\`es les d\'efinitions de 3.1, cocycle not\'e $a$ comme en loc.cite. En agissant par conjugaison, $x$ laisse stable ${\mathcal G}'$ et donc $xh_{w} x^{-1}$ commute aussi \`a $\tilde{s}_{\bf{G}'}$ au m\^eme cocycle pr\`es et pour tout $w\in \Gamma_{F}$, $a(w)\in Z(\hat{G})$. D'o\`u:
$$z \tilde{s}_{\bf{G}'}=
x \tilde{s}_{\bf{G}'} x^{-1}= x (a(w)^{-1}h_{w}\tilde{s}_{\bf{G}'} h_{w}^{-1})x^{-1}$$
$$=a(w)^{-1}xh_{w}x^{-1}\, x \tilde{s}_{\bf{G}'} x^{-1}\, xh_{w}^{-1}x^{-1}
$$
$$=
a(w)^{-1}xh_{w}x^{-1} z xh_{w}^{-1}x^{-1}\, xh_{w}x^{-1}\,  \tilde{s}_{\bf{G}'} \, xh_{w}^{-1}x^{-1}=w zw^{-1} \tilde{s}_{\bf{G}'}.
$$
Cela donne l'\'egalit\'e cherch\'ee $z= wzw^{-1}$.

\bigskip

\subsection{La combinatoire}
On a pr\'ecis\'e les notations ${\bf{G}}',\hat{M},{\bf {M}}'$. Dans l'\'enonc\'e ci-dessous, on \'ecrit $\sim {H}$ pour indiquer que l'on prend l'\'el\'ement consid\'er\'e \`a conjugaison pr\`es sous le groupe $H$. On d\'efinit aussi $W(\hat{M}):=Norm_{\hat{G}}\hat{M}/\hat{M}$, qui est muni d'une action galoisienne et d'une action de $\hat{\theta}$. On d\'efinit de fa\c{c}on identique $W(\hat{M}')=Norm_{\hat{G}'}\hat{M'}/\hat{M}'$, qui est muni d'une action galoisienne provenant de ${\cal M}'$.

\ass{Proposition}{ Dans les deux sommes suivantes, les ${\bf {G}'}, {\bf {M}'}, \hat{M}$ sont des triplets comme ci-dessus. Et on a:
$$
\sum_{{\bf {G}'}/\sim \hat{G}}i(\tilde{G},\tilde{G}') \sum_{{\bf {M}'}/\sim \hat{G}'} \vert W(\hat{M}')^{\Gamma_F}\vert^{-1} S({\bf {G}'}, {\bf {M}'},\hat{M}); \eqno(1)$$
$$
=
$$
 $$
\sum_{\hat{M}/\sim \hat{G}}\vert W(\hat{M})^{\Gamma_F,\hat{\theta}}\vert^{-1} \sum_{ {\bf {M}'}/\sim \hat {M}}i(\hat{M},\tilde{M}') $$
$$\sum_{{\bf G}'={\bf G}'(\tilde{s}); \tilde{s}\in \tilde{s}_{{\bf M}'} Z(\hat{M})^{\Gamma_{F}}/Z(\hat{G})^{\Gamma_{F}}(1-\hat{\theta})(Z(\hat{M})^{\Gamma_{F}})}i'_{\tilde{M}'}(\tilde{G},\tilde{G}')S({\bf {G}'}, {\bf {M}'},\hat{M}) \eqno(2)$$
$$
=
$$
$$
\sum_{\hat{M}/\sim \hat{G}}\vert W(\hat{M})^{\Gamma_F,\hat{\theta}}\vert^{-1} \sum_{ {\bf {M}'}/\sim \hat {M}}i(\hat{M},\tilde{M}') $$
$$\sum_{{\bf G}'={\bf G}'(\tilde{s}); \tilde{s}\in \tilde{s}_{{\bf M}'} Z(\hat{M})^{\Gamma_{F},\hat{\theta}}/Z(\hat{G})^{\Gamma_{F},\hat{\theta}} }i_{\tilde{M}'}(\tilde{G},\tilde{G}')S({\bf {G}'}, {\bf {M}'},\hat{M}).\eqno(3)$$}

L'\'egalit\'e de (2) et (3) r\'esulte de [I] 3.3(2). On prouve l'\'egalit\'e de (1) et (2). 
Chaque somme est une somme sur les triplets ${\bf {G}'}, {\bf {M}'},\hat{M}$ comme expliqu\'e pr\'ec\'edemment. Pour $i=1,2$ et pour ${\bf {G}'}, {\bf {M}'},\hat{M}$  on note $n_{i}({\bf {G}'}, {\bf {M}'},\hat{M})$ le nombre de repr\'esentants de la classe de conjugaison sous $\hat{G}$ de ce triplet qui apparaissent dans la somme (i) (on v\'erifiera que ces nombres sont finis). Et on doit montrer pour tout tel triplet que:
$$
\vert W(\hat{M})^{\Gamma_F,\hat{\theta}}\vert i(\tilde{G},\tilde{G}')n_{1}({\bf {G}'}, {\bf {M}'},\hat{M})\times( \vert W(\hat{M}')^{\Gamma_F}\vert i'_{\tilde{M}'}(\tilde{G},\tilde{G}')i(\hat{M},\tilde{M}')n_{2}({\bf {G}'}, {\bf {M}'},\hat{M}))^{-1}$$vaut 1. En tenant compte de la proposition 6.2, cela revient au m\^eme que de d\'emontrer:
$$
\vert W(\hat{M})^{\Gamma_F,\hat{\theta}}\vert n_{1}({\bf {G}'}, {\bf {M}'},\hat{M})\vert W(\hat{M}')^{\Gamma_F}\vert^{-1}n_{2}({\bf {G}'}, {\bf {M}'},\hat{M})^{-1} \eqno(4)
$$
$$
=\vert \overline{Aut}_{\hat{G}}({\bf {G}'})\vert\vert \overline{Aut}_{\hat{M}}({\bf {M}'})\vert^{-1}
$$
Dans (1), le groupe  $Aut_{\hat{G}}({\bf {G}'})/\hat{G}'$ op\`ere sur les classes de $\hat{G}'$-conjugaison form\'ees d'\'el\'ements ${\bf {M}'}$. Ainsi $n_{1}({\bf {G}'}, {\bf {M}'},\hat{M})= \vert Aut_{\hat{G}}({\bf {G}'})/\hat{G}' Aut_{\hat{G}}({\bf {G}'}, {\bf {M}'})\vert.
$
Le sous-groupe $Z(\hat{G})^{\Gamma_{F}}$ de $Aut_{G}({\bf {G}'})$ op\`ere trivialement sur tout $ {\bf {M}'}$ espace de Levi de ${\bf {G}'}$ et on peut donc voir $Aut_{\hat{G}}({\bf {G}'})/\hat{G}' Aut_{\hat{G}}({\bf {G}'}, {\bf {M}'})$ comme un espace quotient de $\overline{Aut}_{G}({\bf {G}'})$ et on a:
$$
n_{1}({\bf {G}'}, {\bf {M}'},\hat{M})= \vert \overline{Aut}_{G}({\bf {G}'})\vert \vert  Aut_{\hat{G}}({\bf {G}'}, {\bf {M}'})/ Aut_{\hat{G}}({\bf {G}'}, {\bf {M}'})\cap \hat{G}'Z(\hat{G})^{\Gamma_{F}}\vert^{-1}.$$
Or $Aut_{\hat{G}}({\bf {G}'}, {\bf {M}'})\cap \hat{G}'Z(\hat{G})^{\Gamma_{F}}=Norm_{\hat{G}'}({\bf M}')Z(\hat{G})^{\Gamma_{F}}$; d'o\`u
$$
\vert  Aut_{\hat{G}}({\bf {G}'}, {\bf {M}'})/ Aut_{\hat{G}}({\bf {G}'}, {\bf {M}'})\cap \hat{G}'Z(\hat{G})^{\Gamma_{F}}\vert=$$
$$ \vert 
Aut_{\hat{G}}({\bf {G}'}, {\bf {M}'})/ \hat{M}' Z(\hat{G})^{\Gamma_{F}}\vert \vert
Norm_{\hat{G}'}({\bf M}')Z(\hat{G})^{\Gamma_{F}}/\hat{M}' Z(\hat{G})^{\Gamma_{F}}\vert^{-1}$$
$$
=\vert  Aut_{\hat{G}}({\bf {G}'}, {\bf {M}'})/ \hat{M}' Z(\hat{G})^{\Gamma_{F}}\vert \vert W(\hat{M}')^{\Gamma_F}\vert^{-1}\vert
Z(\hat{G})^{\Gamma_{F}}\cap \hat{G}'/\hat{M}' \cap Z(\hat{G})^{\Gamma_{F}}\vert.
$$Or $Z(\hat{G})^{\Gamma_{F}}\cap \hat{G}' \subset Z(\hat{G}')\subset Z(\hat{M}')$,
et on trouve donc que 
$
n_{1}({\bf {G}'}, {\bf {M}'},\hat{M})=$
$$ \vert \overline{Aut}_{G}({\bf {G}'})\vert \vert W(\hat{M}')^{\Gamma_F}\vert \vert 
Aut_{\hat{G}}({\bf {G}'}, {\bf {M}'})/ \hat{M}' Z(\hat{G})^{\Gamma_{F}}\vert^{-1}.
$$

Dans (2), il y a d'abord l'image r\'eciproque $Norm_{\hat{G}}(\hat{M})^*$ de $W(\hat{M})^{\Gamma_{F},\hat{\theta}}$ dans $Norm_{\hat{G}}(\hat{M})$ qui op\`ere sur les ${\bf {M}'}$ alors que dans la somme on n'a pris en compte que  l'action de $\hat{M}$; ensuite, sur les classes de conjugaison de ${\bf {G}'}$ modulo $Z(\hat{M})^{\Gamma_{F}}$ op\`ere $Aut_{\hat{M}}({\bf {M}'})/\hat{M}'Z(\hat{M})^{\Gamma_{F}}$. Ainsi
$$
n_{2}({\bf {G}'}, {\bf {M}'},\hat{M})=
\vert Norm_{\hat{G}}(\hat{M})^*/\hat{M} Aut_{\hat{G}}({\bf {G}'}, {\bf {M}'},\hat{M})\vert$$
$$ \times \vert Aut_{\hat{M}}({\bf {M}'})/ Aut_{\hat{M}}({\bf {G}'}, {\bf {M}'},\hat{M})Z(\hat{M})^{\Gamma_{F}}\vert. 
$$
Ce nombre est \'evidemment fini et on peut le r\'ecrire sous la forme:
$$
\vert W(\hat{M})^{\Gamma_F,\hat{\theta}}\vert \vert Aut_{\hat{G}}({\bf {G}'}, {\bf {M}'},\hat{M})/Aut_{\hat{M}}({\bf {G}'}, {\bf {M}'},\hat{M})\vert^{-1}$$
$$\times \vert \overline{Aut}_{\hat{M}}({\bf {M}'})\vert
\vert Aut_{\hat{M}}({\bf {G}'}, {\bf {M}'},\hat{M})/\hat{M}'(Z(\hat{M})^{\Gamma_{F}}\cap Aut_{\hat{M}}({\bf {G}'}, {\bf {M}'},\hat{M}) )\vert^{-1}
$$
$$
= \vert W(\hat{M})^{\Gamma_F,\hat{\theta}}\vert \vert \overline{Aut}_{\hat{M}}({\bf {M}'})\vert \vert Aut_{\hat{G}}({\bf {G}'}, {\bf {M}'},\hat{M})/\hat{M}'(Z(\hat{M})^{\Gamma_{F}}\cap Aut_{\hat{M}}({\bf {G}'}, {\bf {M}'},\hat{M}) )\vert^{-1}.
$$
Ainsi d\'emontrer (4) est \'equivalent \`a montrer que
$$
\vert 
Aut_{\hat{G}}({\bf {G}'}, {\bf {M}'})/ \hat{M}' Z(\hat{G})^{\Gamma_{F}}\vert =$$
$$
\vert Aut_{\hat{G}}({\bf {G}'}, {\bf {M}'},\hat{M})/\hat{M}'(Z(\hat{M})^{\Gamma_{F}}\cap Aut_{\hat{M}}({\bf {G}'}, {\bf {M}'},\hat{M}) )\vert.
$$
Comme $\hat{M}$ est uniquement d\'etermin\'e par ${\bf {G}'}$ et son espace de Levi $\tilde{M}'$, $Aut_{\hat{G}}({\bf {G}'}, {\bf {M}'})=Aut_{\hat{G}}({\bf {G}'}, {\bf {M}'},\hat{M})$ et 
il suffit de montrer  l'inclusion
$$
Z(\hat{M})^{\Gamma_{F}}\cap Aut_{\hat{G}}({\bf {G}'})\hookrightarrow \hat{M}'Z(\hat{G})^{\Gamma_{F}}.
$$
Pour montrer cette inclusion, on consid\`ere l'image de $Z(\hat{M})^{\Gamma_{F}}\cap Aut_{\hat{G}}({\bf {G}'})$ dans $\hat{T}/Z(\hat{G})$; l'image est incluse dans $(\hat{T}/Z(\hat{G}))^{\Gamma_{F},\hat{\theta}}$; or $\hat{T}/Z(\hat{G})$ est un tore induit au sens suivant: le groupe de ces caract\`eres admet une base sur laquelle $\hat{\theta}$ et $\Gamma_{F}$ op\`ere par permutations. Ainsi $(\hat{T}/Z(\hat{G}))^{\Gamma_{F},\hat{\theta}}$ est un tore connexe et est n\'ecessairement l'image de $\hat{T}^{\Gamma_{F},\hat{\theta},0}$; or $\hat{T}^{\Gamma_{F},\hat{\theta},0}$ est un sous-groupe de $({\hat{M}'})^{\Gamma_{F}}$ et ainsi $Z(\hat{M})^{\Gamma_{F}}\cap Aut_{\hat{G}}({\bf {G}'})$ est un sous-groupe de $(\hat{M}')^{\Gamma_{F}}Z(\hat{G})$ et donc de $(\hat{M}')^{\Gamma_{F}}Z(\hat{G})^{\Gamma_{F}}$.

  \bigskip

\subsection{Un r\'esultat d'annulation}
On travaille ici avec des $K$-espaces. On fixe un ensemble fini $V$  de places contenant $V_{ram}$. On fixe comme toujours une paire de Borel \'epingl\'ee $\hat{{\cal E}}=(\hat{B},\hat{T},(\hat{E}_{\alpha})_{\alpha\in \Delta})$ de $\hat{G}$, conserv\'ee par l'action galoisienne. 
Si $K\tilde{M}$ est un $K$-espace de Levi de $K\tilde{G}$, on sait d\'efinir la notion de donn\'ee endoscopique  de $(KM,K\tilde{M},{\bf a}_{M})$. En fait, $K\tilde{M}$ n'intervient dans cette d\'efinition que via le Levi $\hat{M}$ de $\hat{G}$ qui lui est associ\'e. On sait que l'on peut supposer ce Levi standard et invariant par $\hat{\theta}$ et par l'action galoisienne. Pour se donner un peu plus de libert\'e, on peut imposer la condition plus faible

(1) $\hat{T}\subset \hat{M}$ et il existe $\hat{P}\in {\cal P}(\hat{M})$ de sorte que $(\hat{P},\hat{M})$ soit conserv\'e par $\hat{\theta}$ et par l'action galoisienne.

Consid\'erons non plus un $K$-espace de Levi $K\tilde{M}$ de $K\tilde{G}$ mais un Levi $\hat{M}$ de $\hat{G}$ qui v\'erifie la condition (1). Le cocycle ${\bf  a}$ se pousse en un cocycle ${\bf a}_{\hat{M}}$ \`a valeurs dans $Z(\hat{M})$. On d\'efinit comme en 3.1 la notion de donn\'ee endoscopique disons de $(\hat{M},{\bf a}_{\hat{M}})$. Si $\hat{M}$ correspond \`a un $K$-espace de Levi $K\tilde{M}$, cette notion co\"{\i}ncide avec celle de donn\'ee endoscopique de $(KM,K\tilde{M},{\bf a}_{M})$. Mais la pr\'esente notion est un peu plus g\'en\'erale puisque $\hat{M}$ ne corresponde pas toujours \`a un  tel $K$-espace de Levi. Consid\'erons donc un Levi $\hat{M}$ comme ci-dessus et une donn\'ee endoscopique ${\bf M}'=(M',{\cal M}',\tilde{\zeta})$ de $(\hat{M},{\bf a}_{\hat{M}})$. Un tore maximal de $\hat{M}'$ est isomorphe \`a $\hat{T}^{\hat{\theta},0}$. Si on introduit des sous-tores maximaux $T$ de $G$ et $T'$ de $M'$, on a par dualit\'e un homomorphisme $\xi:T\to T'$. Il n'est d\'efini qu'\`a conjugaison pr\`es mais sa restriction \`a $Z(G)$ est bien d\'efinie et envoie $Z(G)$ dans $Z(M')$. On peut donc d\'efinir l'espace tordu $\tilde{M}'=M'\times_{Z(G)}{\cal Z}(\tilde{G})$ comme en [I] 1.7. Pour une place $v$  hors de $V$, la situation est non ramifi\'ee. Il existe donc un espace de Levi $\tilde{M}_{v}$ de $\tilde{G}$ d\'efini sur $F_{v}$ qui correspond \`a $\hat{M}$. La donn\'ee localis\'ee ${\bf M}'_{v}$ est une donn\'ee endoscopique de $(M_{v},\tilde{M}_{v},{\bf a}_{M_{v}})$.  Le sous-espace hypersp\'ecial $\tilde{K}_{v}\cap \tilde{M}_{v}$ d\'etermine alors un sous-espace hypersp\'ecial $\tilde{K}^{M'}_{v}$ de $\tilde{M}'(F_{v})$, bien d\'efini modulo conjugaison par $M'_{AD}(F_{v})$. On fixe de tels sous-espaces, soumis \`a la condition de compatibilit\'e globale de 1.1.
La notion de donn\'ees auxiliaires $M'_{1}$, $\tilde{M}'_{1}$, $C_{1}$, $\hat{\xi}_{1}$ d\'efinies sur $F$ et non ramifi\'ees hors de $V$ se d\'efinit comme en 3.3 et la preuve du lemme de ce paragraphe montre que de telles donn\'ees existent.  On adjoint \`a ces donn\'ees une familles d'espaces hypersp\'eciaux $(\tilde{K}^{M'}_{1,v})_{v\not\in V}$ relevant les $\tilde{K}^{M'}_{v}$, soumise \`a la m\^eme condition de compatibilit\'e globale. Consid\'erons maintenant une autre s\'erie de donn\'ees auxilaires $M'_{2}$,... $(\tilde{K}^{M'}_{2,v})_{v\not\in V}$. La m\^eme construction qu'en 1.15 d\'efinit une fonction  $\tilde{\lambda}_{12,V}$ qui permet de recoller $C_{c,\lambda_{1}}^{\infty}(\tilde{M}'_{1}(F_{V}))$ \`a  $C_{c,\lambda_{2}}^{\infty}(\tilde{M}'_{1}(F_{V}))$, du moins si $\tilde{M}'(F)\not=\emptyset$. On d\'efinit alors l'espace $C_{c}^{\infty}({\bf M}'_{V})$ comme la limite inductive de ces espaces. On a de m\^eme des espaces $I({\bf M}'_{V})$, $SI({\bf M}'_{V})$ etc... Toutes les constructions formelles que l'on a faites dans le cas o\`u $\hat{M}$ correspondait \`a un $K$-espace de Levi de $K\tilde{G}$ s'\'etendent dans notre situation.

 Cela \'etant, supposons $\tilde{M}'(F)\not=\emptyset$ et ${\bf M}'$ elliptique, c'est-\`a-dire que  $Z(\hat{M}')^{\Gamma_{F},0}=Z(\hat{M})^{\Gamma_{F},\hat{\theta},0}$. Soient $\boldsymbol{\delta}\in D_{\diamond}^{st}({\bf M}'_{V})\otimes Mes(M'(F_{V}))^*$ et ${\bf f}\in I(K\tilde{G}(F_{V}),\omega)\otimes Mes(G(F_{V}))$. On peut poser
$$(2) \qquad I_{*}^{K\tilde{G},{\cal E}}({\bf M}',\boldsymbol{\delta},{\bf f})=\sum_{\tilde{s}\in \tilde{\zeta}Z(\hat{M})^{\Gamma_{F},\hat{\theta}}/Z(\hat{G})^{\Gamma_{F},\hat{\theta}}}i_{\tilde{M}'}(\tilde{G},\tilde{G}'(\tilde{s}))S_{{\bf M}'}^{{\bf G}'(\tilde{s})}(\boldsymbol{\delta}, B^{\tilde{G}},{\bf f}^{{\bf G}'(\tilde{s})}).$$
Par convention, ${\bf f}^{{\bf G}'(\tilde{s})}=0$ si ${\bf G}'(\tilde{s})$ n'est pas relevante.
S'il existe un espace de Levi $K\tilde{M}$ correspondant \`a $\hat{M}$ et si ${\bf M}'$ est relevante pour cet espace de Levi, ce  n'est autre que le terme $I_{K\tilde{M}}^{K\tilde{G},{\cal E}}({\bf M}',\boldsymbol{\delta},{\bf f})$ du paragraphe 4.4.

\ass{Proposition}{Supposons soit qu'il n'existe aucun espace de Levi $K\tilde{M}$ de $K\tilde{G}$ correspondant \`a $\tilde{M}$, soit qu'un tel espace $K\tilde{M}$  existe mais que ${\bf M}'$ soit une donn\'ee non relevante de $(KM,K\tilde{M},{\bf a}_{M})$. Alors $I_{*}^{K\tilde{G},{\cal E}}({\bf M}',\boldsymbol{\delta},{\bf f})=0$ pour tout $\boldsymbol{\delta}\in D_{\diamond}^{st}({\bf M}'_{V})\otimes Mes(M'(F_{V}))^*$ et tout ${\bf f}\in I(K\tilde{G}(F_{V}),\omega)\otimes Mes(G(F_{V}))$.}

 Les places archim\'ediennes compliquent grandement la d\'emonstration, cf.  la remarque (30) de 6.10. On va \'enoncer deux propositions auxiliaires. On montrera  en 6.9 que la seconde entra\^{\i}ne la premi\`ere et que celle-ci entra\^{\i}ne la proposition ci-dessus. On prouvera la seconde proposition auxiliaire en 6.10. Dans  les quatre paragraphes suivants, on conserve la pr\'esente situation et on impose les hypoth\`eses de la proposition. Pour simplifier, on fixe des mesures de Haar sur tous les groupes intervenant, ce qui nous d\'ebarrasse des espaces de mesures.
 
 \bigskip
 
 \subsection{Une premi\`ere proposition auxiliaire}
Soient $v\in V$ et $\hat{\mathfrak{L}}_{v}$ un Levi de $\hat{G}$. Consid\'erons la condition

(1) $\hat{T}\subset \hat{\mathfrak{L}}_{v}$ et il existe $\hat{\mathfrak{P}}_{v}\in {\cal P}(\hat{\mathfrak{L}}_{v})$ de sorte que $(\hat{\mathfrak{P}}_{v},\hat{\mathfrak{L}}_{v})$ soit conserv\'e par $\hat{\theta}$ et par l'action de $\Gamma_{F_{v}}$,

\noindent qui est l'analogue locale de 6.6(1).  
Notons $\hat{\mathfrak{M}}_{v}$ le commutant de $Z(\hat{M}')^{\Gamma_{F_{v}},0}$ dans $\hat{G}$. C'est un Levi de $\hat{G}$ inclus dans $\hat{M}$. L'inclusion peut \^etre stricte. Ce Levi v\'erifie  (1).  
Posons $\hat{\mathfrak{M}}_{V}=(\hat{\mathfrak{M}}_{v})_{v\in V}$. Notons ${\cal L}(\hat{\mathfrak{M}}_{V})$ l'ensemble des familles $\hat{\mathfrak{L}}_{V}=(\hat{\mathfrak{L}}_{v})_{v\in V}$ telles que, pour tout $v\in V$, $\hat{\mathfrak{L}}_{v}$ soit un Levi de $\hat{G}$ v\'erifiant (1) et contenant $\hat{\mathfrak{M}}_{v}$. Consid\'erons une telle famille. On pose
$${\cal Z}(\hat{\mathfrak{L}}_{V})=Z(\hat{M})^{\Gamma_{F},\hat{\theta}}\cap (\cap_{v\in V}Z(\hat{\mathfrak{L}}_{v}))=Z(\hat{M})^{\Gamma_{F},\hat{\theta}}\cap (\cap_{v\in V}Z(\hat{\mathfrak{L}}_{v})^{\Gamma_{F_{v}},\hat{\theta}}).$$
Pour une place $v\in V$, les m\^emes consid\'erations que dans le paragraphe pr\'ec\'edent s'appliquent: on d\'efinit la notion de donn\'ee endoscopique (locale) de $(\hat{\mathfrak{L}}_{v},{\bf a}_{\hat{\mathfrak{L}}_{v}})$. En particulier, puisque ${\bf M}'_{v}$ est une donn\'ee endoscopique elliptique de $(\hat{\mathfrak{M}}_{v},{\bf a}_{\hat{\mathfrak{M}}_{v}})$ et que $\hat{\mathfrak{M}}_{v}$ est un Levi de $\hat{\mathfrak{L}}_{v}$, on peut d\'efinir la donn\'ee endoscopique $\boldsymbol{ \mathfrak{L}}'_{v}(\tilde{s})$ de  $(\hat{\mathfrak{L}}_{v},{\bf a}_{\hat{\mathfrak{L}}_{v}})$ pour $\tilde{s}\in \tilde{\zeta}Z(\hat{\mathfrak{M}})^{\Gamma_{F_{v}},\hat{\theta}}/Z(\hat{\mathfrak{L}})^{\Gamma_{F_{v}},\hat{\theta}}$. Consid\'erons un \'el\'ement  $s\in \tilde{\zeta}Z(\hat{M})^{\Gamma_{F},\hat{\theta}}/{\cal Z}(\hat{\mathfrak{L}}_{V})$. Alors $\boldsymbol{ \mathfrak{L}}'_{v}(\tilde{s})$ est d\'efini pour tout $v\in V$ et on pose $\boldsymbol{\mathfrak{L}}'_{V}(\tilde{s})=(\boldsymbol {\mathfrak{L}}'_{v}(\tilde{s}))_{v\in V}$. Relevons $\tilde{s}$ en un \'el\'ement de $\tilde{\zeta} Z(\hat{M})^{\Gamma_{F},\hat{\theta}}/Z(\hat{G})^{\Gamma_{F},\hat{\theta}}$. Alors la donn\'ee globale ${\bf G}'(\tilde{s})$ est bien d\'efinie. Pour tout $v\in V$, $\boldsymbol {\mathfrak{L}}'_{v}(\tilde{s})$ est une "donn\'ee de Levi" de la donn\'ee locale ${\bf G}'_{v}(\tilde{s})$. Fixons des donn\'ees auxiliaires globales $G'_{1}(\tilde{s})$,... , $(\tilde{K}'_{1,v})_{v\not\in V}$    pour ${\bf G}'(\tilde{s})$. Il s'en d\'eduit par restriction des donn\'ees auxilaires locales pour les $\boldsymbol {\mathfrak{L}}'_{v}(\tilde{s})$, ainsi que des donn\'ees auxiliaires globales pour ${\bf M}'$.  Conform\'ement \`a notre habitude, on note par exemple $\tilde{\mathfrak{L}}'_{1,v}(\tilde{s})$ ou $M'_{1}(\tilde{s})$  l'image r\'eciproque de  $\tilde{\mathfrak{L}}'_{v}(\tilde{s})$ ou $M'$ dans $\tilde{G}'_{1}(\tilde{s})$.  En posant $\tilde{\mathfrak{L}}'_{1,V}(\tilde{s};F_{V})=\prod_{v\in V}\mathfrak{\tilde{L}}'_{1,v}(\tilde{s},F_{v})$, on a un espace de fonctions  $SI_{\lambda_{1}}^{\infty}(\tilde{\mathfrak{L}}'_{1,V}(\tilde{s},F_{V}))$. Faisons varier le rel\`evement de $\tilde{s}$ et les donn\'ees auxiliaires. On remplace l'indice $1$ par $2$ pour ces nouvelles donn\'ees. On a une fonction de recollement d\'efinie sur le produit fibr\'e 
$\tilde{\mathfrak{L}}'_{12,V}(\tilde{s};F_{V})$ de $\tilde{\mathfrak{L}}'_{1,V}(\tilde{s};F_{V})$ et $\tilde{\mathfrak{L}}'_{2,V}(\tilde{s};F_{V})$ au-dessus de $\tilde{\mathfrak{L}}_{V}'(\tilde{s};F_{V})$. Elle n'est d\'efinie qu'\`a multiplication pr\`es par un scalaire. 
 Mais on a remarqu\'e au paragraphe pr\'ec\'edent que l'on pouvait normaliser canoniquement  les restrictions de ces fonctions au   produit fibr\'e similaire $\tilde{M}'_{12}(F_{V})$ qui est inclus dans $\tilde{\mathfrak{L}}'_{12,V}(\tilde{s};F_{V})$. Cette normalisation fixe la fonction de recollement sur tout $\tilde{\mathfrak{L}}'_{12,V}(\tilde{s};F_{V})$. On peut alors d\'efinir par limite inductive un espace que l'on note $SI(\boldsymbol{\mathfrak{L}}_{V}'(\tilde{s}))$. Dualement, on a un espace $D^{st}_{g\acute{e}om}(\boldsymbol{\mathfrak{L}}'_{V}(\tilde{s}))$.

  Soient $v\in V$ et $\delta\in \tilde{\mathfrak{L}}_{v}'(\tilde{s};F_{v})$. Sur la cl\^oture alg\'ebrique $\bar{F}_{v}$, il existe un espace de Levi $\tilde{\mathfrak{L}}_{v}$ de $\tilde{G}$ qui correspond \`a $\hat{\mathfrak{L}}_{v}$. La classe de conjugaison de la partie semi-simple de $\delta$ correspond \`a une classe de conjugaison d'un \'el\'ement semi-simple $\gamma$ de $\tilde{\mathfrak{L}}_{v}$. On dit que $\delta$ est $\tilde{G}$-r\'egulier, resp. $\tilde{G}$-\'equisingulier, si $\gamma$ est fortement $\tilde{G}$-r\'egulier, resp.  $\tilde{G}$-\'equisingulier. 
   Un \'el\'ement $\delta=(\delta_{v})_{v\in V}\in \tilde{\mathfrak{L}}'_{V}(\tilde{s};F_{V})$ est dit $\tilde{G}$-r\'egulier, resp. $\tilde{G}$-\'equisingulier,  aux places archim\'ediennes si
    $\delta_{v}$ est  $\tilde{G}$-r\'egulier, resp. $\tilde{G}$-\'equisingulier,  pour tout $v$ archim\'edienne.  
     
   Soit $\hat{\mathfrak{R}}_{V}\in {\cal L}(\hat{\mathfrak{M}}_{V})$. On note ${\cal L}(\hat{\mathfrak{R}}_{V})$ l'ensemble des $\hat{\mathfrak{L}}_{V}\in {\cal L}(\hat{\mathfrak{M}}_{V})$ telles que $\hat{\mathfrak{R}}_{v}\subset \hat{\mathfrak{L}}_{v}$ pour tout $v\in V$. Soient $s\in \tilde{\zeta}Z(\hat{M})^{\Gamma_{F},\hat{\theta}}/{\cal Z}(\hat{\mathfrak{L}}_{V})$, ${\bf f}\in SI(\boldsymbol{\mathfrak{L}}'_{V}(\tilde{s})) $ et $\boldsymbol{\delta}\in D^{st}_{ g\acute{e}om}(\boldsymbol{\mathfrak{R}}'_{V}(\tilde{\zeta}))= D^{st}_{ g\acute{e}om}(\boldsymbol{\mathfrak{R}}'_{V}(\tilde{s})) $.  Supposons 
   
   (2)  le support de $\boldsymbol{\delta}$ soit form\'e d'\'el\'ements $\tilde{G}$-\'equisinguliers aux places archim\'ediennes.  
   
   Fixons des donn\'ees auxiliaires comme ci-dessus. Pour simplifier, on peut supposer que $\boldsymbol{\delta}$ et ${\bf f}$  s'identifient respectivement \`a $\otimes_{v\in V}\boldsymbol{\delta}_{1,v}$ et $\otimes_{v\in V}{\bf f}_{1,v}$. Les termes du produit
 $$\prod_{v\in V}S_{\tilde{\mathfrak{R}}'_{1,v}(\tilde{s}),\lambda_{1}}^{\tilde{\mathfrak{L}}'_{1,v}(\tilde{s})}(\boldsymbol{\delta}_{1,v},B^{\tilde{G}},{\bf f}_{1,v})$$
 sont bien d\'efinis.  Notons que l'hypoth\`ese (2) nous permet si l'on veut de supprimer la mention du syst\`eme de fonctions $B^{\tilde{G}}$ aux places archim\'ediennes. Quand on fait varier les donn\'ees auxiliaires, le produit ne change pas. On le note
 $S_{\boldsymbol{\mathfrak{R}}'_{V}(\tilde{\zeta})}^{\boldsymbol{\mathfrak{L}}'_{V}(\tilde{s})}(\boldsymbol{\delta},B^{\tilde{G}},{\bf f})$. 
 
  Soient $\hat{\mathfrak{R}}_{V}\in {\cal L}(\hat{\mathfrak{M}}_{V})$ et $\tilde{s}\in \tilde{\zeta}{\cal Z}(\hat{\mathfrak{R}}_{V}) /Z(\hat{G})^{\Gamma_{F},\hat{\theta}}$. Supposons ${\bf G}'(\tilde{s})$ et $\boldsymbol{\mathfrak{R}}'_{V}(\tilde{\zeta})$ elliptiques. Notons ${\cal L}^{\tilde{G}'(\tilde{s})}(\tilde{\mathfrak{R}}'_{V}(\tilde{\zeta}))$ l'ensemble des familles $\tilde{\mathfrak{L}}'_{V}=(\tilde{\mathfrak{L}}'_{v})_{v\in V}$ telles que, pour tout $v\in V$, $\tilde{\mathfrak{L}}'_{v}$ soit un espace de Levi de $\tilde{G}'(\tilde{s})$ d\'efini sur $F_{v}$ et contenant $\tilde{\mathfrak{R}}'_{v}(\tilde{\zeta})$. Pour une telle famille et pour $v\in V$, notons $\hat{\mathfrak{L}}_{v}$ le commutant de $Z(\hat{\mathfrak{L}}'_{v})^{\Gamma_{F_{v}},0}$ dans $\hat{G}$. Alors la famille $\hat{\mathfrak{L}}_{V}=(\hat{\mathfrak{L}}_{v})_{v\in V}$ appartient \`a ${\cal L}(\hat{\mathfrak{R}}_{V})$. La famille $\tilde{\mathfrak{L}}'_{V}$ appara\^{\i}t comme la famille d'espaces de Levi associ\'ee \`a la donn\'ee endoscopique elliptique $\boldsymbol{\mathfrak{L}}'_{V}(\tilde{s})$ de $(\hat{\mathfrak{L}}_{V},{\bf a}_{\hat{\mathfrak{L}}_{V}})$. Pour simplifier les notations, nous noterons directement $\tilde{\mathfrak{L}}'_{V}(\tilde{s})$ les \'el\'ements de ${\cal L}^{\tilde{G}'(\tilde{s})}(\tilde{\mathfrak{R}}'_{V}(\tilde{\zeta}))$. Remarquons que l'application "terme constant"
  $$\begin{array}{ccc}SI({\bf G}'(\tilde{s})) &\to&SI(\boldsymbol{\mathfrak{L}}'(\tilde{s})) \\ {\bf f}&\mapsto &{\bf f}_{\boldsymbol{\mathfrak{L}}'(\tilde{s})}\\ \end{array}$$
  est bien d\'efinie.
  
  Soient $\hat{\mathfrak{R}}_{V}\in {\cal L}(\hat{\mathfrak{M}}_{V})$, $\boldsymbol{\delta}\in D^{st}_{g\acute{e}om}(\boldsymbol{\mathfrak{R}}'_{V}(\tilde{\zeta}))$ et $f\in I(K\tilde{G}(F_{V}),\omega)$. On suppose (2) v\'erifi\'ee. On suppose aussi
  
  (3) $\boldsymbol{\mathfrak{R}}'(\tilde{\zeta})$ est elliptique.
  
  On d\'efinit
  $$J(\hat{\mathfrak{R}}_{V},\boldsymbol{\delta},f)=\sum_{\tilde{s}\in \tilde{\zeta}{\cal Z}(\hat{\mathfrak{R}}_{V})/Z(\hat{G})^{\Gamma_{F},\hat{\theta}}} i_{\tilde{M}'}(\tilde{G},\tilde{G}'(\tilde{s}))$$
  $$\sum_{\tilde{\mathfrak{L}}'_{V}(\tilde{s})\in {\cal L}^{\tilde{G}'(\tilde{s})}(\tilde{\mathfrak{R}}'_{V}(\tilde{\zeta}))}e_{\tilde{M}'_{V}}^{\tilde{G}'(\tilde{s})}(\tilde{M}',\tilde{\mathfrak{L}}'_{V}(\tilde{s}))S_{\boldsymbol{\mathfrak{R}}'_{V}(\tilde{\zeta})}^{\boldsymbol{\mathfrak{L}}'_{V}(\tilde{s})}(\boldsymbol{\delta},B^{\tilde{G}},(f^{{\bf G}'(\tilde{s})})_{\boldsymbol{\mathfrak{L}}'_{V}(\tilde{s})}).$$
  
  \ass{Proposition}{Soient $\hat{\mathfrak{R}}_{V}\in {\cal L}(\hat{\mathfrak{M}}_{V})$, $\boldsymbol{\delta}\in D^{st}_{g\acute{e}om}(\boldsymbol{\mathfrak{R}}'_{V}(\tilde{\zeta}))$ et $f\in I(K\tilde{G}(F_{V}),\omega)$.   On suppose que (2) et (3) sont v\'erifi\'ees et que $\boldsymbol{\delta}$  est l'image par induction d'un \'el\'ement de $D^{st}_{g\acute{e}om}({\bf M}'_{V})$. Alors $J(\hat{\mathfrak{R}}_{V},\boldsymbol{\delta}, f)=0$.}
  
  \bigskip
  
  \subsection{Une deuxi\`eme proposition auxiliaire}
  
  \ass{Proposition}{Soient $\hat{\mathfrak{R}}_{V}\in {\cal L}(\hat{\mathfrak{M}}_{V})$, $\boldsymbol{\delta}\in D^{st}_{g\acute{e}om}(\boldsymbol{\mathfrak{R}}'_{V}(\tilde{\zeta}))$ et $ f\in I(K\tilde{G}(F_{V}),\omega)$.  On suppose  que $\boldsymbol{\mathfrak{R}}'(\tilde{\zeta})$ est elliptique et que $\boldsymbol{\delta}$  est  l'image par induction d'un \'el\'ement de $D^{st}_{orb}({\bf M}'_{V})$ dont le support est form\'e d'\'el\'ements $\tilde{G}$-r\'eguliers aux places archim\'ediennes. Alors $J(\hat{\mathfrak{R}}_{V},\boldsymbol{\delta},f)=0$.}

   \bigskip
   
   \subsection{R\'eduction de la proposition 6.6}
   Evidemment, la proposition 6.8 est un cas particulier de la proposition 6.7. Mais nous allons prouver qu'inversement, elle implique cette proposition. Consid\'erons la situation de cette proposition 6.7. On peut supposer $f=\otimes_{v\in V}f_{v}$. On fait jouer  aux donn\'ees ${\bf G}'(\tilde{\zeta})$ et $\boldsymbol{\mathfrak{R}}'_{V}(\tilde{\zeta})$  un r\^ole de r\'ef\'erence. Pour simplifier, on supprime le terme $\tilde{\zeta}$ des notations, en posant ${\bf G}'={\bf G}'(\tilde{\zeta})$, $\boldsymbol{\mathfrak{R}}'_{V}=\boldsymbol{\mathfrak{R}}'_{V}(\tilde{\zeta})$.  Fixons des donn\'ees suppl\'ementaires $G'_{1}$,...,$(\tilde{K}'_{1,v})_{v\not\in V}$ pour ${\bf G}'$. On peut supposer que $\boldsymbol{\delta}$ s'identifie \`a $\otimes_{v\in V}\boldsymbol{\delta}_{1,v}$, avec $\boldsymbol{\delta}_{1,v}\in D_{g\acute{e}om,\lambda_{1,v}}^{st}(\tilde{\mathfrak{R}}'_{1,v}(F_{v}))$. On peut fixer pour tout $v\in V$  une classe de conjugaison stable semi-simple $ {\cal O}'_{v}$ dans $\tilde{M}'(F_{v})$ de sorte que $\boldsymbol{\delta}_{1,v}$ soit induite d'un \'el\'ement de $D^{st}_{g\acute{e}om,\lambda_{1,v}}(\tilde{M}'_{1}(F_{v}))$ dont le support est form\'e d'\'el\'ements de partie semi-simple dans ${\cal O}'_{v}$. Notons  ${\cal O}_{v}^{\tilde{\mathfrak{R}}'_{v}}$ la classe de conjugaison stable dans $\tilde{\mathfrak{R}}'_{v}(F_{v})$ contenant ${\cal O}'$. L'hypoth\`ese sur $\boldsymbol{\delta}$ est que ${\cal O}_{v}^{\tilde{\mathfrak{R}}'_{v}}$ est $\tilde{G}$-\'equisinguli\`ere pour toute place $v$ archim\'edienne.  Nous consid\'erons comme fix\'es $f$ et les composantes $\delta_{1,v}$ pour $v$ non archim\'edienne. On va \'etudier comment $J(\hat{\mathfrak{R}}_{V},\boldsymbol{\delta},f)$ d\'epend des $\delta_{1,v}$ pour $v$ archim\'edienne. Pour cela, consid\'erons pour toute place archim\'edienne un \'el\'ement $\boldsymbol{\tau}_{1,v}\in D_{g\acute{e}om,\lambda_{1,v}}^{st}(\tilde{\mathfrak{R}}'_{1,v}(F_{v}))$.  On suppose soit que $\boldsymbol{\tau}_{1,v}=\boldsymbol{\delta}_{1,v}$, soit que $\boldsymbol{\tau}_{1,v}$ appartient \`a $D_{orb,\lambda_{1,v}}^{st}(\tilde{\mathfrak{R}}'_{1,v}(F_{v}))$ et que son support est form\'e d'\'el\'ements $\tilde{G}$-r\'eguliers proches de ${\cal O}_{v}^{\tilde{\mathfrak{R}}'_{v}}$. Pour unifier les notations, on pose $\boldsymbol{\tau}_{1,v}=\boldsymbol{\delta}_{1,v}$ pour toute $v\in V$ non-archim\'edienne.
   On note $\boldsymbol{\tau}$ l'\'el\'ement de $D^{st}_{g\acute{e}om}(\boldsymbol{\mathfrak{R}}'_{V})$ auquel s'identifie $\otimes_{v\in V}\boldsymbol{\tau}_{1,v}$. 
   
   Pour tout $\tilde{s}\in \tilde{\zeta}Z(\hat{\mathfrak{R}}_{V})/Z(\hat{G})^{\Gamma_{F},\hat{\theta}}$, fixons des donn\'ees suppl\'ementaires $G'_{1}(\tilde{s})$,...,$(\tilde{K}'_{1,v}(\tilde{s}))_{v\not\in V}$ pour ${\bf G}'(\tilde{s})$. On a deux s\'eries de donn\'ees auxiliaires pour la donn\'ee $\boldsymbol{\mathfrak{R}}'_{V}=\boldsymbol{\mathfrak{R}}'_{V}(\tilde{s})$. 
  Comme on l'a dit,  les espaces $SI_{\lambda_{1}(\tilde{s})}(\tilde{\mathfrak{R}}'_{1,V}(\tilde{s};F_{V}))$ et $SI_{\lambda_{1}}(\tilde{\mathfrak{R}}'_{1,V}(F_{V}))$ se recollent canoniquement. On peut d\'ecomposer   cet isomorphisme canonique en produit d'isomorphismes sur toutes les places $v\in V$. On a alors des isomorphismes duaux entre espaces de distributions. 
  Pour tout $v\in V$,  $\boldsymbol{\tau}_{1,v}$ s'identifie ainsi \`a un \'el\'ement $\boldsymbol{\tau}_{1,v}(\tilde{s})\in D^{st}_{g\acute{e}om,\lambda_{1,v}(\tilde{s})}(\tilde{\mathfrak{R}}'_{1,v}(\tilde{s};F_{v}))$. D'autre part, $f^{{\bf G}'(\tilde{s})}$ s'identifie \`a un \'el\'ement $\otimes_{v\in V}f_{1,v}(\tilde{s})$.
   On a l'\'egalit\'e
   $$(1) \qquad  S_{\boldsymbol{\mathfrak{R}}'_{V}(\tilde{s})}^{\boldsymbol{\mathfrak{L}}'_{V}(\tilde{s})}(\boldsymbol{\tau},B^{\tilde{G}},(f^{{\bf G}'(\tilde{s})})_{\boldsymbol{\mathfrak{L}}'_{V}(\tilde{s})})=\prod_{v\in V}S_{\tilde{\mathfrak{R}}'_{1,v}(\tilde{s}),\lambda_{1,v}(\tilde{s})}^{\tilde{\mathfrak{L}}'_{1,v}(\tilde{s})}(\boldsymbol{\tau}_{1,v}(\tilde{s}),B^{\tilde{G}},f_{1,v}(\tilde{s})_{\tilde{\mathfrak{L}}'_{1,v}(\tilde{s})}).$$
   Si $v$ est archim\'edienne, on a dit que l'on pouvait oublier le syst\`eme de fonctions $B^{\tilde{G}}$ en vertu de l'hypoth\`ese sur le support de $\boldsymbol{\tau}_{1,v}$. On a vu en [V] 1.4 (2) et (3) qu'il existait $\varphi_{1,v}\in SI_{\lambda_{1,v}(\tilde{s})}(\tilde{\mathfrak{R}}'_{1,v}(\tilde{s};F_{v}))$ de sorte que
   
 (2)  pour tout $\boldsymbol{\tau}_{1,v} $ comme ci-dessus, on a 
 $$S_{\tilde{\mathfrak{R}}'_{1,v}(\tilde{s}),\lambda_{1,v}(\tilde{s})}^{\tilde{\mathfrak{L}}'_{1,v}(\tilde{s})}(\boldsymbol{\tau}_{1,v},f_{1,v}(\tilde{s})_{\tilde{\mathfrak{L}}'_{1,v}(\tilde{s})})=S^{\tilde{\mathfrak{R}}'_{1,v}(\tilde{s})}_{\lambda_{1,v}(\tilde{s})}(\boldsymbol{\tau}_{1,v}(\tilde{s}),\varphi_{1,v}).$$

A l'aide des recollements fix\'es, on peut identifier $\varphi_{1,v}$ \`a un \'el\'ement de $SI_{\lambda_{1,v}}(\tilde{\mathfrak{R}}'_{1,v}(F_{v}))$. Comme cet \'el\'ement d\'epend de $\tilde{s}$ et de $\boldsymbol{\mathfrak{L}}'_{V}(\tilde{s})$, notons-le $\phi_{1,v}[\tilde{s},\boldsymbol{\mathfrak{L}}'_{V}(\tilde{s})]$. Le membre de droite de (2)  devient 
$$S^{\tilde{\mathfrak{R}}'_{1,v}}_{\lambda_{1,v}}(\boldsymbol{\tau}_{1,v},\phi_{1,v}[\tilde{s},\boldsymbol{\mathfrak{L}}'_{V}(\tilde{s})]).$$
Notons $V_{\infty}$ l'ensemble des places archim\'ediennes de $F$ et indiquons par un indice $\infty$ les produits ou produits tensoriels sur les places $v\in V_{\infty}$. Par exemple $\boldsymbol{\tau}_{1,\infty}=\otimes_{v\in V_{\infty}}\boldsymbol{\tau}_{1,v}$.  
L'\'egalit\'e (1) devient
$$S_{\boldsymbol{\mathfrak{R}}'_{V}(\tilde{s})}^{\boldsymbol{\mathfrak{L}}'_{V}(\tilde{s})}(\boldsymbol{\tau},B^{\tilde{G}},(f^{{\bf G}'(\tilde{s})})_{\boldsymbol{\mathfrak{L}}'_{V}(\tilde{s})})=c[\tilde{s},\boldsymbol{\mathfrak{L}}'_{V}(\tilde{s})] S^{\tilde{\mathfrak{R}}'_{1,\infty}}_{\lambda_{1,\infty}}(\boldsymbol{\tau}_{1,\infty},\phi_{1,\infty}[\tilde{s},\boldsymbol{\mathfrak{L}}'_{V}(\tilde{s})]),$$
o\`u   $c[\tilde{s},\boldsymbol{\mathfrak{L}}'_{V}(\tilde{s})]$ est ind\'ependant des $\boldsymbol{\tau}_{1,v}$ pour $v\in Val_{\infty}(F)$. Posons
$$\phi_{1,\infty}=\sum_{\tilde{s}\in \tilde{\zeta}{\cal Z}(\hat{\mathfrak{R}}_{V})/Z(\hat{G})^{\Gamma_{F},\hat{\theta}}} i_{\tilde{M}'}(\tilde{G},\tilde{G}'(\tilde{s}))$$
$$\sum_{\tilde{\mathfrak{L}}'_{V}(\tilde{s})\in {\cal L}^{\tilde{G}'(\tilde{s})}(\tilde{\mathfrak{R}}'_{V}(\tilde{\zeta}))}e_{\tilde{M}'_{V}}^{\tilde{G}'(\tilde{s})}(\tilde{M}',\tilde{\mathfrak{L}}'_{V}(\tilde{s}))c[\tilde{s},\boldsymbol{\mathfrak{L}}'_{V}(\tilde{s})]\phi_{1,\infty}[\tilde{s},\boldsymbol{\mathfrak{L}}'_{V}(\tilde{s})].$$
Alors

$$(3) \qquad J(\hat{\mathfrak{R}}_{V},\boldsymbol{\tau},f)=S^{\tilde{\mathfrak{R}}'_{1,\infty}}_{\lambda_{1,\infty}}(\boldsymbol{\tau}_{1,\infty},\phi_{1,\infty}).$$
Soit $\boldsymbol{\mu}_{1,\infty}\in D^{st}_{g\acute{e}om,\lambda_{1,\infty}}(\tilde{M}'_{1}(F_{\infty}))$, \`a support dans les \'el\'ements de $\tilde{M}'(F_{\infty})$ de partie semi-simple dans ${\cal O}'_{\infty}$, tel que $\boldsymbol{\delta}_{1,\infty}$ soit l'induite de $\boldsymbol{\mu}_{1,\infty}$. En appliquant (3) \`a $\boldsymbol{\tau}=\boldsymbol{\delta}$, on obtient
$$J(\hat{\mathfrak{R}}_{V},\boldsymbol{\delta},f)=S^{\tilde{M}'_{1,\infty}}_{\lambda_{1,\infty}}(\boldsymbol{\mu}_{1,\infty},\phi_{1,\infty,\tilde{M}'_{1,\infty}}).$$
On veut prouver que le membre de gauche est nul. Il suffit de prouver que $\phi_{1,\infty,\tilde{M}'_{1,\infty}}$ est nul au voisinage de ${\cal O}'_{\infty}$. Pr\'ecis\'ement, il suffit de prouver que, pour $\boldsymbol{\nu}_{1,\infty}\in D^{st}_{orb,\lambda_{1,\infty}}(\tilde{M}'_{1}(F_{\infty}))$, \`a support $\tilde{G}$-r\'egulier et proche de ${\cal O}'_{\infty}$, on a $S^{\tilde{M}'_{1,\infty}}_{\lambda_{1,\infty}}(\boldsymbol{\nu}_{1,\infty},\phi_{1,\infty,\tilde{M}'_{1,\infty}})=0$. Fixons un tel $\nu_{1,\infty}$, notons $\boldsymbol{\tau}_{1,\infty}$ l'induite de $\boldsymbol{\nu}_{1,\infty}$ \`a $\tilde{\mathfrak{R}}'_{1,\infty}(F_{\infty})$. Compl\'etons $\boldsymbol{\tau}_{1,\infty}$ en un \'el\'ement $\boldsymbol{\tau}$ de composantes $\boldsymbol{\delta}_{1,v}$ aux places de $V$ non-archim\'ediennes. Le m\^eme calcul que ci-dessus montre que
 $$S^{\tilde{M}'_{1,\infty}}_{\lambda_{1,\infty}}(\boldsymbol{\nu}_{1,\infty},\phi_{1,\infty,\tilde{M}'_{1,\infty}})=J(\hat{\mathfrak{R}}_{V},\boldsymbol{\tau},f).$$
 Mais $\boldsymbol{\tau}$ v\'erifie les hypoth\`eses de 6.8. Donc le membre de droite ci-dessus est nul. D'o\`u l'assertion cherch\'ee, ce qui prouve la proposition 6.7.
 
 Nous allons maintenant prouver que cette proposition 6.7 entra\^{\i}ne la proposition 6.6.  Consid\'erons la d\'efinition 6.6(2). On utilise la proposition 4.2(i) pour d\'evelopper chaque terme $S_{{\bf M}'}^{{\bf G}'(\tilde{s})}(\boldsymbol{\delta},B^{\tilde{G}},{\bf f}^{{\bf G}'(\tilde{s})})$. Avec les notations que l'on a introduites, on obtient
 $$S_{{\bf M}'}^{{\bf G}'(\tilde{s})}(\boldsymbol{\delta},B^{\tilde{G}},{\bf f}^{{\bf G}'(\tilde{s})})=\sum_{\tilde{\mathfrak{L}}'_{V}(\tilde{s})\in {\cal L}^{\tilde{ G}'(\tilde{s})}(\tilde{M}'_{V})}e_{\tilde{M}'_{V}}^{\tilde{G}'(\tilde{s})}(\tilde{M}',\tilde{\mathfrak{L}}'_{V}(\tilde{s}))S_{{\bf M}'_{V}}^{\boldsymbol{\mathfrak{L}}'_{V}(\tilde{s})}(\boldsymbol{\delta},B^{\tilde{G}},(f^{{\bf G}'(\tilde{s})})_{\boldsymbol{\mathfrak{L}}'_{V}(\tilde{s})}).$$
 Fixons des donn\'ees auxiliaires comme plus haut. Les int\'egrales du membre de droite se d\'ecomposent alors en produit sur $v\in V$ d'int\'egrales locales. Consid\'erons une place $v$ archim\'edienne. Le terme local est
 $$(4) \qquad S_{\tilde{M}'_{1,v}(\tilde{s}),\lambda_{1,v}(\tilde{s})}^{\tilde{\mathfrak{L}}'_{1,v}(\tilde{s})}(\boldsymbol{\delta}_{1,v}(\tilde{s}),B^{\tilde{G}},(f_{v}^{\tilde{G}'_{1,v}(\tilde{s})})_{\tilde{\mathfrak{L}}'_{1,v}(\tilde{s})}).$$
 On peut supposer comme plus haut que les \'el\'ements du support de $\boldsymbol{\delta}_{1,v}$ ont leur partie semi-simple dans une classe de conjugaison stable ${\cal O}'_{v}$. 
  Fixons $H_{v}\in {\cal A}_{\tilde{M}'_{v}}$ en position g\'en\'erale. Relevons-le en un \'el\'ement $H_{1,v}\in {\cal A}_{\tilde{M}'_{1,v}(\tilde{s})}$. Consid\'erons  un Levi de $\tilde{\mathfrak{L}}'_{v}(\tilde{s})$ contenant $\tilde{M}'_{v}$. Conform\'ement aux notations introduites en 6.7, notons-le $\tilde{\mathfrak{R}}'_{v}(\tilde{s})$.
 Les d\'efinitions de [V] 6.3 s'\'etendent au cas des distributions se transformant selon le caract\`ere $\lambda_{1,v}(\tilde{s})$ de $C_{1}(\tilde{s};F_{v})$. On a d\'efini dans cette r\'ef\'erence un \'el\'ement
 $$ \xi^{\tilde{\mathfrak{R}}'_{1,v}(\tilde{s}),st}(\boldsymbol{\delta}_{1,v}(\tilde{s}),B^{\tilde{G}},H_{1,v}^{\tilde{\mathfrak{R}}'_{1,v}(\tilde{s})}).$$
 C'est une distribution induite \`a $\tilde{\mathfrak{R}}'_{1,v}(\tilde{s};F_{v})$ d'un \'el\'ement de $D^{st}_{g\acute{e}om}(\tilde{M}'_{1,v}(\tilde{s};F_{v}))$ dont le support est form\'e d'\'el\'ements de partie semi-simple dans ${\cal O}'_{v}$.   
  La proposition  [V] 6.3 entra\^{\i}ne que le  (3)  est faiblement \'equivalent \`a la fonction qui, \`a $H_{1,v} $ associe
 $$\sum_{\tilde{\mathfrak{R}}'_{v}(\tilde{s})\in {\cal L}^{\tilde{\mathfrak{L}}'_{v}(\tilde{s})}(\tilde{M}'_{v})}S_{\tilde{\mathfrak{R}}'_{1,v}(\tilde{s}),\lambda_{1,v}(\tilde{s})}^{\tilde{\mathfrak{L}}'_{1,v}(\tilde{s})}(exp(H_{1,v,\tilde{\mathfrak{R}}'_{1,v}(\tilde{s})})\xi^{\tilde{\mathfrak{R}}'_{1,v}(\tilde{s}),st}(\boldsymbol{\delta}_{1,v}(\tilde{s}),B^{\tilde{G}},H_{1,v}^{\tilde{\mathfrak{R}}'_{1,v}(\tilde{s})}),(f_{v}^{\tilde{G}'_{1,v}(\tilde{s})})_{\tilde{\mathfrak{L}}'_{1,v}(\tilde{s})}).$$
 Cela implique que, si l'on remplace dans l'expression ci-dessus $H_{1,v}$ par $H_{1,v}/n$, pour un entier $n\geq1$, la limite de cette expression quand $n$ tend vers l'infini est \'egale \`a (4). Posons 
 $$\boldsymbol{\tau}_{1,v}^{\tilde{\mathfrak{R}}'_{1,v}(\tilde{s})}(n)=exp(H_{1,v,\tilde{\mathfrak{R}}'_{1,v}(\tilde{s})}/n)\xi^{\tilde{\mathfrak{R}}'_{1,v}(\tilde{s}),st}(\boldsymbol{\delta}_{1,v}(\tilde{s}),B^{\tilde{G}},H_{1,v}^{\tilde{\mathfrak{R}}'_{1,v}(\tilde{s})}/n).$$
 Pour unifier les notations, pour une place $v\in V$ non archim\'edienne, posons 
  $$\boldsymbol{\tau}_{1,v}^{\tilde{M}'_{1,v}(\tilde{s})}(n)=\delta_{1,v}(\tilde{s}).$$
  Enfin, on a d\'efini en 6.7 l'ensemble ${\cal L}^{\tilde{ G}'(\tilde{s})}(\tilde{ M}'_{V})$. Notons ${\cal L}_{\infty}^{\tilde{G}'(\tilde{s})}(\tilde{M}'_{V})$ le sous-ensemble des $\tilde{\mathfrak{R}}'_{V}(\tilde{s})\in {\cal L}^{\tilde{ G}'(\tilde{s})}(\tilde{ M}'_{V})$ tels que $\tilde{\mathfrak{R}}'_{v}(\tilde{s})=\tilde{M}'_{v}$ pour toute place archim\'edienne. Avec ces d\'efinitions, on obtient que $S_{{\bf M}'}^{{\bf G}'(\tilde{s})}(\boldsymbol{\delta},B^{\tilde{G}},{\bf f}^{{\bf G}'(\tilde{s})})$ est \'egale \`a la limite quand $n$ tend vers l'infini de
  $$(5) \qquad \sum_{\tilde{\mathfrak{R}}'_{V}(\tilde{s})\in {\cal L}_{\infty}^{\tilde{ G}'(\tilde{s})}(\tilde{M}'_{V})}\sum_{\tilde{\mathfrak{L}}'_{V}(\tilde{s})\in {\cal L}^{\tilde{G}'(\tilde{s})}(\tilde{\mathfrak{R}}'_{V}(\tilde{s}))}e_{\tilde{M}'_{V}}^{\tilde{G}'(\tilde{s})}(\tilde{M}',\tilde{\mathfrak{L}}'_{V}(\tilde{s}))$$
  $$\prod_{v\in V}S_{\tilde{\mathfrak{R}}'_{1,v}(\tilde{s}),\lambda_{1,v}(\tilde{s})}^{\tilde{\mathfrak{L}}'_{1,v}(\tilde{s})}(\boldsymbol{\tau}_{1,v}^{\tilde{M}'_{1,v}(\tilde{s})}(n),B^{\tilde{G}},(f_{v}^{\tilde{G}'_{1,v}(\tilde{s})})_{\tilde{\mathfrak{L}}'_{1,v}(\tilde{s})}).$$
   Les distributions  $\boldsymbol{\tau}_{1,v}^{\tilde{\mathfrak{R}}'_{1,v}(\tilde{s})}(n)$ d\'ependent  des donn\'ees auxiliaires mais on v\'erifie sans peine qu'elles se recollent convenablement quand on change de donn\'ees auxiliaires. On doit toutefois prendre garde au fait que la translation par $exp(H_{1,v,\tilde{\mathfrak{R}}'_{1,v}(\tilde{s})})$ pour $v$ archim\'edienne n'est compatible au recollement qu'\`a un caract\`ere pr\`es. Plus pr\'ecis\'ement, pour une telle place, on a introduit en [IV] 2.1 un caract\`ere $\lambda_{\mathfrak{A}_{\tilde{G}'_{1,v}(\tilde{s})}}$ de ${\cal A}_{\tilde{G}'_{1}(\tilde{s})}$. Alors les distributions 
   $$(\prod_{v\in Val_{\infty}(F)}\lambda_{\mathfrak{A}_{\tilde{G}'_{1,v}(\tilde{s})}}(H_{1,v,\tilde{G}'_{1,v}(\tilde{s})}/n)\otimes_{v\in V}\boldsymbol{\tau}_{1,v}^{\tilde{\mathfrak{R}}'_{1,v}(\tilde{s})}(n)$$
   se recollent en une distribution que l'on note $\boldsymbol{\tau}^{\boldsymbol{\mathfrak{R}}'_{V}(\tilde{s})}(n)$. 
    Cette multiplication par le produit des $\lambda_{\mathfrak{A}_{\tilde{G}'_{1,v}(\tilde{s})}}(H_{1,v,\tilde{G}'_{1,v}(\tilde{s})}/n)$ ne nous g\^ene pas: on peut multiplier (5) par ce terme sans changer les propri\'et\'es de cette expression.  Alors (5) se r\'ecrit
    $$\sum_{\tilde{\mathfrak{R}}'_{V}(\tilde{s})\in {\cal L}_{\infty}^{\tilde{ G}'(\tilde{s})}(\tilde{M}'_{V})}\sum_{\tilde{\mathfrak{L}}'_{V}(\tilde{s})\in {\cal L}^{\tilde{G}'(\tilde{s})}(\tilde{\mathfrak{R}}'_{V}(\tilde{s}))}e_{\tilde{M}'_{V}}^{\tilde{G}'(\tilde{s})}(\tilde{M}',\tilde{\mathfrak{L}}'_{V}(\tilde{s}))S_{\boldsymbol{\mathfrak{R}}'_{V}(\tilde{s})}^{\boldsymbol{\mathfrak{L}}'_{V}(\tilde{s})}(\boldsymbol{\tau}^{\boldsymbol{\mathfrak{R}}'_{V}(\tilde{s})}(n),B^{\tilde{G}},(f^{{\bf G}'(\tilde{s})})_{\boldsymbol{\mathfrak{L}}'_{V}(\tilde{s})}).$$
 En revenant \`a 6.6(2), on voit que $I_{*}^{K\tilde{G},{\cal E}}({\bf M}',\boldsymbol{\delta},f)$ est \'egale \`a la limite quand $n$ tend vers l'infini de
 $$(6) \qquad \sum_{\tilde{s}\in \tilde{\zeta}Z(\hat{M})^{\Gamma_{F},\hat{\theta}}/Z(\hat{G})^{\Gamma_{F},\hat{\theta}}}i_{\tilde{M}'}(\tilde{G},\tilde{G}'(\tilde{s}))  \sum_{\tilde{\mathfrak{R}}'_{V}(\tilde{s})\in {\cal L}_{\infty}^{\tilde{ G}'(\tilde{s})}(\tilde{M}'_{V})}$$
 $$\sum_{\tilde{\mathfrak{L}}'_{V}(\tilde{s})\in {\cal L}^{\tilde{G}'(\tilde{s})}(\tilde{\mathfrak{R}}'_{V}(\tilde{s}))}e_{\tilde{M}'_{V}}^{\tilde{G}'(\tilde{s})}(\tilde{M}',\tilde{\mathfrak{L}}'_{V}(\tilde{s}))S_{\boldsymbol{\mathfrak{R}}'_{V}(\tilde{s})}^{\boldsymbol{\mathfrak{L}}'_{V}(\tilde{s})}(\boldsymbol{\tau}^{\boldsymbol{\mathfrak{R}}'_{V}(\tilde{s})}(n),B^{\tilde{G}},(f^{{\bf G}'(\tilde{s})})_{\boldsymbol{\mathfrak{L}}'_{V}(\tilde{s})}).$$ 
 A tout $\boldsymbol{\mathfrak{R}}'_{V}(\tilde{s})$ intervenant est associ\'e un \'el\'ement $\hat{\mathfrak{R}}_{V}\in {\cal L}(\hat{\mathfrak{M}}_{V})$. Plus pr\'ecis\'ement, celui-ci  appartient au sous-ensemble ${\cal L}_{\infty}(\hat{\mathfrak{M}}_{V})$ d\'efini par les conditions $\hat{\mathfrak{R}}_{v}=\hat{\mathfrak{M}}_{v}$ pour $v\in V$ non-archim\'edienne. On peut remplacer la somme en $\boldsymbol{\mathfrak{R}}'_{V}(\tilde{s})$ par une somme en $\hat{\mathfrak{R}}_{V}$, que l'on permute avec la somme en $\tilde{s}$. On peut ensuite d\'ecomposer cette derni\`ere somme en une somme sur $\tilde{t}\in \tilde{\zeta}Z(\hat{M})^{\Gamma_{F},\hat{\theta}}/{\cal Z}(\hat{\mathfrak{R}}_{V})$ et une somme en $\tilde{s}\in \tilde{t}{\cal Z}(\hat{\mathfrak{R}}_{V})/Z(\hat{G})^{\Gamma_{F},\hat{\theta}}$. L'expression (6) devient la somme sur $\hat{\mathfrak{R}}_{V}\in {\cal L}_{\infty}(\hat{\mathfrak{M}}_{V})$ de
 $$\sum_{\tilde{t}\in \tilde{\zeta}Z(\hat{M})^{\Gamma_{F},\hat{\theta}}/{\cal Z}(\hat{\mathfrak{R}}_{V}); \boldsymbol{\mathfrak{R} }'_{V}(\tilde{t})\text{ elliptique}}J(\hat{\mathfrak{R}}_{V},\tilde{t},\boldsymbol{\tau}^{\boldsymbol{\mathfrak{R}}'_{V}(\tilde{t})}(n),f),$$
 o\`u
 $$J(\hat{\mathfrak{R}}_{V},\tilde{t},\boldsymbol{\tau}^{\boldsymbol{\mathfrak{R}}'_{V}(\tilde{t})}(n),f)=\sum_{\tilde{s}\in \tilde{t}{\cal Z}(\hat{\mathfrak{R}}_{V})/Z(\hat{G})^{\Gamma_{F},\hat{\theta}}}i_{\tilde{M}'}(\tilde{G},\tilde{G}'(\tilde{s}))$$
 $$\sum_{\tilde{\mathfrak{L}}'_{V}(\tilde{s})\in {\cal L}^{\tilde{G}'(\tilde{s})}(\tilde{\mathfrak{R}}'_{V}(\tilde{s}))}e_{\tilde{M}'_{V}}^{\tilde{G}'(\tilde{s})}(\tilde{M}',\tilde{\mathfrak{L}}'_{V}(\tilde{s}))S_{\boldsymbol{\mathfrak{R}}'_{V}(\tilde{s})}^{\boldsymbol{\mathfrak{L}}'_{V}(\tilde{s})}(\boldsymbol{\tau}^{\boldsymbol{\mathfrak{R}}'_{V}(\tilde{s})}(n),B^{\tilde{G}},(f^{{\bf G}'(\tilde{s})})_{\boldsymbol{\mathfrak{L}}'_{V}(\tilde{s})}).$$
 Dans le cas o\`u $\tilde{t}=\tilde{\zeta}$, cette expression est \'egale \`a l'expression $J(\hat{\mathfrak{R}}_{V},\boldsymbol{\tau}^{\boldsymbol{\mathfrak{R}}'_{V}(\tilde{t})}(n),f)$ de 6.7. Pour $\tilde{t}$ quelconque, elle est \'egale \`a  l'analogue de cette expression quand on remplace la donn\'ee de d\'epart ${\bf M}'=(M',{\cal M}',\tilde{\zeta})$ par la donn\'ee \'equivalente $(M',{\cal M}',\tilde{t})$. On peut donc lui appliquer cette proposition 6.7. Les distributions $\boldsymbol{\tau}^{\boldsymbol{\mathfrak{R}}'_{V}(\tilde{t})}(n)$ v\'erifient par construction les hypoth\`eses de cette proposition: en une place archim\'edienne $v$, la translation par $exp(H_{1,v,\tilde{\mathfrak{R}}'_{1,v}(\tilde{s})}/n)$ assure que le support  de la distribution est $\tilde{G}$-\'equisingulier puisque $H_{v}$ est en position g\'en\'erale. La proposition 6.7 implique donc que $J(\hat{\mathfrak{R}}_{V},\tilde{t},\boldsymbol{\tau}^{\boldsymbol{\mathfrak{R}}'_{V}(\tilde{t})}(n),f)=0$. Alors l'expression  (6) est nulle. Sa limite $I_{*}^{K\tilde{G},{\cal E}}({\bf M}',\boldsymbol{\delta},f)$ est nulle elle aussi, ce qui prouve la proposition 6.6. 
 
 \subsection{Preuve de la proposition 6.8}
 On fixe des donn\'ees $\hat{\mathfrak{R}}_{V} $, $\boldsymbol{\delta}$ et $f$ comme dans l'\'enonc\'e de la proposition 6.8. On suppose $f=\otimes_{v\in V}f_{v}$. On a besoin de facteurs de transfert globaux. Pour cela, on fixe les extensions
 $$1\to G\to H\to D\to 1,\,\, \tilde{G}\to \tilde{H}$$
 construites dans la preuve de 3.8. Pour tout $\tilde{s}\in \tilde{\zeta} Z(\hat{M})^{\Gamma_{F},\hat{\theta}}/Z(\hat{G})^{\Gamma_{F},\hat{\theta}}$, on \'etend comme dans ce paragraphe la donn\'ee ${\bf G}'(\tilde{s})$ en une donn\'ee ${\bf H}'(\tilde{t})$ et on fixe des donn\'ees auxiliaires $H'_{1}(\tilde{t})$,...,$\hat{\xi}_{1}(\tilde{t})$. On en d\'eduit comme en 3.9 des donn\'ees auxiliaires $G'_{1}(\tilde{s})$,...,$\hat{\xi}_{1}(\tilde{s})$ pour ${\bf G}'(\tilde{s})$, que l'on compl\`ete par le choix d'espaces hypersp\'eciaux $(\tilde{K}'_{1,v}(\tilde{s}))_{v\not\in V}$. Il s'en d\'eduit un facteur de transfert canonique comme en 3.9, que l'on d\'ecompose en produit sur $v\in V$ de facteurs locaux. On note $f_{v}^{\tilde{G}'_{1}(\tilde{s})}$ le transfert de $f_{v}$ calcul\'e \`a l'aide de ce facteur.  Comme en  6.9, on consid\`ere les donn\'ees auxiliaires relatives \`a $\tilde{\zeta}$ comme des donn\'ees de r\'ef\'erence et, pour celles-ci, on supprime $\tilde{\zeta}$ de la notation: ${\bf G}'={\bf G}'(\tilde{\zeta})$, $G'_{1}=G'_{1}(\tilde{\zeta})$ etc... On peut supposer que $\boldsymbol{\delta}$ s'identifie \`a $\otimes_{v\in V}\boldsymbol{\delta}_{1,v}$, avec $\boldsymbol{\delta}_{1,v}\in D_{g\acute{e}om,\lambda_{1,v}}^{st}(\tilde{\mathfrak{R}}'_{1,v}(F_{v}))$ pour tout $v\in V$. 
  Comme on l'a vu en 6.7, on dispose pour tout $\tilde{s}\in \tilde{\zeta}{\cal Z}(\hat{\mathfrak{R}}_{V})/Z(\hat{G})^{\Gamma_{F},\hat{\theta}}$ de recollements canoniques entre $SI_{\lambda_{1}}(\tilde{M}'_{1}(F_{V}))$ et $SI_{\lambda_{1}(\tilde{s})}(\tilde{M}'_{1}(\tilde{s};F_{V}))$ et aussi entre 
   $SI_{\lambda_{1}}(\tilde{\mathfrak{R}}'_{1,V}(F_{V}))$ et $SI_{\lambda_{1}(\tilde{s})}(\tilde{\mathfrak{R}}'_{1,V}(\tilde{s};F_{V}))$. On d\'ecompose ceux-ci en produit tensoriel d'isomorphismes locaux. On fait de m\^eme pour les espaces de distributions.  Alors chaque $\boldsymbol{\delta}_{1,v}$ s'identifie \`a $\boldsymbol{\delta}_{1,v}(\tilde{s})\in D_{g\acute{e}om,\lambda_{1,v}(\tilde{s})}^{st}(\tilde{\mathfrak{R}}'_{1,v}(\tilde{s};F_{v}))$. La d\'efinition de $J(\hat{\mathfrak{R}}_{V},\boldsymbol{\delta},f)$ se r\'ecrit
   $$(1) \qquad J(\hat{\mathfrak{R}}_{V},\boldsymbol{\delta},f)=\sum_{\tilde{s}\in \tilde{\zeta}{\cal Z}(\hat{\mathfrak{R}}_{V})/Z(\hat{G})^{\Gamma_{F},\hat{\theta}}}i_{\tilde{M}'}(\tilde{G},\tilde{G}'(\tilde{s}))$$
   $$\sum_{\tilde{\mathfrak{L}}'_{V}(\tilde{s})\in {\cal L}^{\tilde{G}'(\tilde{s})}(\tilde{\mathfrak{R}}'_{V}(\tilde{s}))}e_{\tilde{M}'_{V}}^{\tilde{G}'(\tilde{s})}(\tilde{M}',\tilde{\mathfrak{L}}'_{V}(\tilde{s}))\prod_{v\in V}S_{\tilde{\mathfrak{R}}'_{1,v}(\tilde{s}),\lambda_{1,v}(\tilde{s})}^{\tilde{\mathfrak{L}}'_{1,v}(\tilde{s})}(\boldsymbol{\delta}_{1,v}(\tilde{s}),B^{\tilde{G}},(f_{v}^{\tilde{G}'_{1,v}(\tilde{s})})_{\tilde{\mathfrak{L}}'_{1,v}(\tilde{s})}).$$

   On peut supposer que, pour tout $v\in V$, il y a une classe de conjugaison stable semi-simple ${\cal O}'_{v}\subset \tilde{M}'(F_{v})$ de sorte que
 
 - $\boldsymbol{\delta}_{1,v}$ est induite d'un \'el\'ement de $D^{st}_{orb,\lambda_{1,v}}(\tilde{M}'_{1,v}(F_{v}))$ dont le support est form\'e d'\'el\'ements de partie semi-simple dans ${\cal O}'_{v}$;
 
 - si $v$ est archim\'edienne, ${\cal O}'_{v}$ est form\'e d'\'el\'ements $\tilde{G}$-r\'eguliers. 
 
 Fixons $\epsilon_{v}\in {\cal O}'_{v}$ tel que $M'_{\epsilon_{v}}$ soit quasi-d\'eploy\'e sur $F_{v}$. On note $\tilde{R}'_{v}$ le commutant de $A_{M'_{\epsilon_{v}}}$ dans $\tilde{M}'$. On a $A_{\tilde{R}'_{v}}=A_{M'_{\epsilon_{v}}}$.  On peut fixer une distribution ${\bf d}_{v}\in D^{st}_{g\acute{e}om,\lambda_{1}}(\tilde{R}'_{1,v}(F_{v}))$ de sorte que
 
 - $\boldsymbol{\delta}_{1,v}$ soit l'induite de ${\bf d}_{v}$ \`a $\tilde{\mathfrak{R}}'_{1,v}(F_{v})$;
 
 - si $v$ est non-archim\'edienne, les projections dans $\tilde{R}'_{v}(F_{v})$ des \'el\'ements du support de ${\bf d}_{v}$ ont leur partie semi-simple dans la classe de conjugaison stable ${\cal O}_{\tilde{R}'_{v}}$ de $\epsilon_{v}$ dans $\tilde{R}'_{v}(F_{v})$;
 
 - si $v$ est archim\'edienne, ${\bf d}_{v}$ est une int\'egrale orbitale stable associ\'ee \`a un rel\`evement de $\epsilon_{v}$ dans $\tilde{R}'_{1,v}(F_{v})$. 
 
 On a
 
 (2)  $\epsilon_{v}$ appartient \`a un sous-tore tordu maximal elliptique de $\tilde{R}'_{v}$.
 
 En effet, soit $T_{v}^{\flat}$ un sous-tore maximal elliptique de $M'_{\epsilon_{v}}$ et $T_{v}$ son commutant dans $R'_{v}$. Alors l'ensemble $T_{v}\epsilon_{v}$ r\'epond \`a la question. Pour que cette construction sois-correcte, il faut \'evidemment que $M'_{\epsilon_{v}}$ poss\`ede un sous-tore maximal elliptique. C'est toujours vrai si $v$ est non-archim\'edienne. Si $v$ est archim\'edienne, $\epsilon_{v}$ est $\tilde{G}$-r\'egulier par hypoth\`ese donc $M'_{\epsilon_{v}}$ est lui-m\^eme un tore et l'assertion s'ensuit.

  On note $\hat{R}_{v}$ le commutant de $Z(\hat{R}'_{v})^{\Gamma_{F_{v}},0}$ dans $\hat{G}$. C'est un Levi de $\hat{\mathfrak{M}}_{v}$. Quitte \`a remplacer la donn\'ee locale ${\bf M}'_{v}$ par une donn\'ee \'equivalente, on peut supposer que tous ces Levi sont standard et que  la donn\'ee locale ${\bf M}'_{v}$ provient d'une donn\'ee ${\bf R}'_{v}$ pour $\hat{R}_{v}$, cf. [I] 3.4. C'est-\`a-dire qu'en posant ${\cal R}'_{v}={\cal M}'_{v}\cap {^LR}_{v}$, le triplet ${\bf R}'_{v}=(R'_{v},{\cal R}'_{v},\tilde{\zeta})$ est une donn\'ee endoscopique de $\hat{R}_{v}$ et que ${\cal M}'_{v}=\hat{M}'{\cal R}'_{v}$. Notons que ${\bf R}'_{v}$ est une donn\'ee elliptique par construction. 
 
 {\bf Remarque.} On pourrait poser des d\'efinitions plus sophistiqu\'ees \'evitant de remplacer ${\bf M}'_{v}$ par une donn\'ee \'equivalente. En tout cas, ce remplacement ne perturbera pas la suite de la d\'emonstration.
 
Le Levi $\hat{R}_{v}$ v\'erifie la condition 6.7(1).  On pose $\hat{R}_{V}=(\hat{R}_{v})_{v\in V}$. On d\'efinit comme en 6.7 l'ensemble ${\cal L}(\hat{R}_{V})$ des familles $\hat{L}_{V}=(\hat{L}_{v})_{v\in V}$ telles que, pour tout $v\in V$, $\hat{L}_{v}$ soit un Levi de $\hat{G}$ v\'erifiant 6.7(1) et contenant $\hat{R}_{v}$. Consid\'erons une telle famille. Pour $\tilde{s}\in \tilde{\zeta}Z(\hat{M})^{\Gamma_{F},\hat{\theta}}/Z(\hat{G})^{\Gamma_{F},\hat{\theta}}$ et pour une place $v\in V$, on d\'efinit la donn\'ee endoscopique ${\bf L}'_{v}(\tilde{s})$ de $\hat{L}_{v}$ associ\'ee \`a $\tilde{s}$ et \`a la donn\'ee endoscopique ${\bf R}'_{v}$ du Levi $\hat{R}_{v}$ de $\hat{L}_{v}$. En posant 
$${\cal Z}(\hat{L}_{V})=Z(\hat{M})^{\Gamma_{F},\hat{\theta}}\cap(\cap_{v\in V}Z(\hat{L}_{v})),$$
cette donn\'ee ne d\'epend que de l'image de $\tilde{s}$ modulo ${\cal Z}(\hat{L}_{V})$. Consid\'erons deux \'el\'ements $\tilde{s}_{1}$ et $\tilde{s}_{2}$ de $\tilde{\zeta}{\cal Z}(\hat{\mathfrak{R}}_{V})/Z(\hat{G})^{\Gamma_{F},\hat{\theta}}$ ayant m\^eme image modulo ce groupe. On a choisi ci-dessus des donn\'ees auxiliaires pour les donn\'ees ${\bf G}'(\tilde{s}_{1})$ et ${\bf G}'(\tilde{s}_{2})$ qui se restreignent pour toute place $v\in V$ en des donn\'ees auxiliaires pour ${\bf L}'_{v}(\tilde{s}_{1})={\bf L}'_{v}(\tilde{s}_{2})$. On dispose donc d'espaces $SI_{\lambda_{1,v}(\tilde{s}_{1})}(\tilde{L}'_{1,v}(\tilde{s}_{1};F_{v}))$ et $SI_{\lambda_{1,v}(\tilde{s}_{2})}(\tilde{L}'_{1,v}(\tilde{s}_{2};F_{v}))$. Ces espaces sont canoniquement isomorphes. En effet, on sait que ces espaces sont isomorphes, l'isomorphisme n'\'etant en g\'en\'eral d\'efini qu'\`a un scalaire pr\`es. Il s'agit de normaliser celui-ci. On a d\'ej\`a fix\'e les isomorphismes
$$SI_{\lambda_{1,v}(\tilde{s}_{1})}(\tilde{\mathfrak{R}}'_{1,v}(\tilde{s}_{1};F_{v}))\simeq SI_{\lambda_{1,v}}(\tilde{\mathfrak{R}}'_{1,v}(F_{v}))\simeq SI_{\lambda_{1,v}(\tilde{s}_{2})}(\tilde{\mathfrak{R}}'_{1,v}(\tilde{s}_{2};F_{v})).$$
On normalise nos isomorphismes de sorte que le diagramme suivant soit commutatif
$$\begin{array}{ccc}SI_{\lambda_{1,v}(\tilde{s}_{1})}(\tilde{\mathfrak{R}}'_{1,v}(\tilde{s}_{1};F_{v}))&\to&SI_{\lambda_{1,v}(\tilde{s}_{2})}(\tilde{\mathfrak{R}}'_{1,v}(\tilde{s}_{2};F_{v}))\\ \downarrow&&\downarrow\\ SI_{\lambda_{1,v}(\tilde{s}_{1})}(\tilde{R}'_{1,v}(\tilde{s}_{1};F_{v}))&\to&SI_{\lambda_{1,v}(\tilde{s}_{2})}(\tilde{R}'_{1,v}(\tilde{s}_{1};F_{v}))\\ \uparrow&&\uparrow\\ SI_{\lambda_{1,v}(\tilde{s}_{1})}(\tilde{L}'_{1,v}(\tilde{s}_{1};F_{v}))&\to&SI_{\lambda_{1,v}(\tilde{s}_{2})}(\tilde{L}'_{1,v}(\tilde{s}_{2};F_{v}))\\ \end{array}$$
o\`u les applications verticales sont les applications "termes constants". De m\^emes consid\'erations valent pour les espaces de distributions. 

Posons $\tilde{R}'_{V}=\prod_{v\in V}\tilde{R}'_{v}$. Soit  $\tilde{\zeta}{\cal Z}(\hat{\mathfrak{R}}_{V})/Z(\hat{G})^{\Gamma_{F},\hat{\theta}}$, supposons ${\bf G}'(\tilde{s})$ elliptique. Notons ${\cal L}^{\tilde{G}'(\tilde{s})}(\tilde{R}'_{V})$ l'ensemble des familles $\tilde{L}'_{V}=(\tilde{L}'_{v})_{v\in V}$ telles que, pour tout $v\in V$, $\tilde{L}'_{v}$ soit un espace de Levi de $\tilde{G}'(\tilde{s})$ d\'efini sur $F_{v}$ et contenant $\tilde{R}'_{v}$. Pour une telle famille et pour $v\in V$, notons $\hat{L}_{v}$ le commutant de $Z(\hat{L}'_{v})^{\Gamma_{F_{v}},0}$ dans $\hat{G}$. Alors la famille $\hat{L}_{V}=(\hat{L}_{v})_{v\in V}$ appartient \`a ${\cal L}(\hat{R}_{V})$. La famille $\tilde{L}'_{V}$ appara\^{\i}t comme la famille d'espaces de Levi associ\'ee \`a la donn\'ee endoscopique elliptique  ${\bf L}'(\tilde{s})$ de $(\hat{L}_{V},{\bf a}_{\hat{L}_{V}})$. Comme en 6.7, on notera directement $\tilde{L}'_{V}(\tilde{s})$ les \'el\'ements de ${\cal L}^{\tilde{G}'(\tilde{s})}(\tilde{R}'_{V})$.

 Consid\'erons la formule (1). Fixons $\tilde{s}$, $\tilde{\mathfrak{L}}'_{V}(\tilde{s})$ y apparaissant et une place $v\in V$. Par les isomorphismes canoniques, la distribution ${\bf d}_{v}$ introduite ci-dessus s'identifie \`a ${\bf d}_{v}(\tilde{s})\in D^{st}_{g\acute{e}om,\lambda_{1}(\tilde{s})}(\tilde{R}'_{1,v}(\tilde{s};F_{v}))$. La distribution $\boldsymbol{\delta}_{1,v}(\tilde{s})$ est l'induite de ${\bf d}_{v}(\tilde{s})$ \`a $\tilde{\mathfrak{R}}'_{1,v}(\tilde{s};F_{v})$. On applique les propositions [II] 1.14(ii) ou [V] 1.6(ii). On obtient
 $$S_{\tilde{\mathfrak{R}}'_{1,v}(\tilde{s}),\lambda_{1,v}(\tilde{s})}^{\tilde{\mathfrak{L}}'_{1,v}(\tilde{s})}(\boldsymbol{\delta}_{1,v}(\tilde{s}),B^{\tilde{G}},(f_{v}^{\tilde{G}'_{1,v}(\tilde{s})})_{\tilde{\mathfrak{L}}'_{1,v}(\tilde{s})})=\sum_{\tilde{L}'_{v}(\tilde{s})\in {\cal L}^{\tilde{\mathfrak{L}}'_{v}(\tilde{s})}(\tilde{R}'_{v})}e_{\tilde{R}'_{v}}^{\tilde{\mathfrak{L}}'_{v}(\tilde{s})}(\tilde{\mathfrak{R}}'_{v},\tilde{L}'_{v}(\tilde{s}))$$
 $$S_{\tilde{R}'_{1,v}(\tilde{s}),\lambda_{1,v}(\tilde{s})}^{\tilde{L}'_{1,v}(\tilde{s})}({\bf d}_{v}(\tilde{s}),B^{\tilde{G}},(f_{v}^{\tilde{G}'_{1,v}(\tilde{s})})_{\tilde{L}'_{1,v}(\tilde{s})}).$$
 Pour tout $\hat{L}_{V}\in {\cal L}(\hat{R}_{V})$, posons 
$$(3) \qquad J(\hat{\mathfrak{R}}_{V},\hat{L}_{V},\boldsymbol{\delta},f)=  \sum_{\tilde{s}\in \tilde{\zeta}Z(\hat{\mathfrak{R}}_{V})/Z(\hat{G})^{\Gamma_{F},\hat{\theta}}, {\bf L}'_{V}(\tilde{s})\text{ elliptique}} E(\hat{L}_{V},\tilde{s})$$
$$\prod_{v\in V}S_{\tilde{R}'_{1,v}(\tilde{s}),\lambda_{1,v}(\tilde{s})}^{\tilde{L}'_{1,v}(\tilde{s})}({\bf d}_{v}(\tilde{s}),B^{\tilde{G}},(f_{v}^{\tilde{G}'_{1,v}(\tilde{s})})_{\tilde{L}'_{1,v}(\tilde{s})}),$$
 o\`u
 $$(4) \qquad E(\hat{L}_{V},\tilde{s})=i_{\tilde{M}'}(\tilde{G},\tilde{G}'(\tilde{s}))\sum_{\tilde{\mathfrak{L}}'_{V}(\tilde{s})\in {\cal L}^{\tilde{G}'(\tilde{s})}(\tilde{\mathfrak{R}}'_{V}(\tilde{s}))\cap {\cal L}^{\tilde{G}'(\tilde{s})}(\tilde{L}'_{V}(\tilde{s}))}e_{\tilde{M}'_{V}}^{\tilde{G}'(\tilde{s})}(\tilde{M}',\tilde{\mathfrak{L}}'_{V}(\tilde{s}))$$
 $$\prod_{v\in V}e_{\tilde{R}'_{v}}^{\tilde{\mathfrak{L}}'_{v}(\tilde{s})}(\tilde{\mathfrak{R}}'_{v},\tilde{L}'_{v}(\tilde{s})).$$
 Alors la formule (1) se r\'ecrit
 $$J(\hat{\mathfrak{R}}_{V},\boldsymbol{\delta},f)=\sum_{\hat{L}_{V}\in {\cal L}(\hat{R}_{V})}J(\hat{\mathfrak{R}}_{V},\hat{L}_{V},\boldsymbol{\delta},f).$$
 Pour prouver que le membre de gauche est nul, il nous suffit de fixer $\hat{L}_{V}$ et de prouver que $J(\hat{\mathfrak{R}}_{V},\hat{L}_{V},\boldsymbol{\delta},f)=0$. 
 
 Fixons d\'esormais $\hat{L}_{V}\in {\cal L}(\hat{R}_{V})$.   
 Pour $v\in V$, aux Levi $\hat{\mathfrak{M}}_{v}$, $\hat{\mathfrak{R}}_{v}$, $\hat{R}_{v}$ et $\hat{L}_{v}$ sont associ\'es des espaces ${\cal A}_{\tilde{\mathfrak{M}}_{v}}$, ${\cal A}_{\tilde{\mathfrak{R}}_{v}}$, ${\cal A}_{\tilde{R}_{v}}$ et ${\cal A}_{\tilde{L}_{v}}$. Par exemple, ${\cal A}_{\tilde{\mathfrak{M}}_{v}}=X^*(Z(\hat{\mathfrak{M}}_{v})^{\Gamma_{F_{v}},\hat{\theta},0})\otimes_{{\mathbb Z}}{\mathbb R}$.  On a
 $$\left.\begin{array}{ccc}\left.\begin{array}{c}{\cal A}_{\tilde{\mathfrak{R}}_{v}}\\ {\cal A}_{\tilde{M}_{v}}\\ \end{array}\right\rbrace&\subset& {\cal A}_{\tilde{\mathfrak{M}}_{v}}\\ &&{\cal A}_{\tilde{L}_{v}}\\ \end{array}\right\rbrace\subset {\cal A}_{\tilde{R}_{v}} .$$
 Posons par exemple ${\cal A}_{\tilde{\mathfrak{M}}_{V}}=\oplus_{v\in V}{\cal A}_{\tilde{\mathfrak{M}}_{v}}$. Rappelons que l'on a un plongement diagonal $\Delta:{\cal A}_{\tilde{M}}^{\tilde{G}}\to {\cal A}_{\tilde{\mathfrak{M}}_{V}}$. Par composition, on obtient des homomorphismes
  $${\cal A}_{\tilde{M}}^{\tilde{G}}\to {\cal A}_{\tilde{\mathfrak{M}}_{V}}/{\cal A}_{\tilde{\mathfrak{R}}_{V}} $$
  et 
   $${\cal A}_{\tilde{M}}^{\tilde{G}}\to {\cal A}_{\tilde{\mathfrak{M}}_{V}}\to {\cal A}_{\tilde{R}_{V}}\to {\cal A}_{\tilde{R}_{V}}/{\cal A}_{\tilde{L}_{V}}.$$
   On note 
    $$D:{\cal A}_{\tilde{M}}^{\tilde{G}}\to {\cal A}_{\tilde{\mathfrak{M}}_{V}}/{\cal A}_{\tilde{\mathfrak{R}}_{V}}\oplus {\cal A}_{\tilde{R}_{V}}/{\cal A}_{\tilde{L}_{V}}$$
leur somme directe.  On obtient dualement un homomorphisme
 $$\hat{D}:Z(\hat{M})^{\Gamma_{F},\hat{\theta}}/Z(\hat{G})^{\Gamma_{F},\hat{\theta}}\to \oplus_{v\in V}(Z(\hat{\mathfrak{M}}_{v})^{\Gamma_{F_{v}},\hat{\theta}}/Z(\hat{\mathfrak{R}}_{v})^{\Gamma_{F_{v}},\hat{\theta}}\oplus Z(\hat{R}_{v})^{\Gamma_{F_{v}},\hat{\theta}}/Z(\hat{L}_{v})^{\Gamma_{F_{v}},\hat{\theta}}).$$
 Supposons que $D$ soit un isomorphisme. Alors $\hat{D}$ est surjectif et de noyau fini. On note $k(D)$ le nombre d'\'el\'ements  de ce noyau. On note $d(D)$ le nombre tel que $D$ identifie la mesure sur son ensemble de d\'epart avec $d(\Delta)$ fois celle sur son ensemble d'arriv\'ee. On pose $e(D)=d(D)k(D)^{-1}$. Pour $v\in V$, l'homomorphisme
 $$ Z(\hat{\mathfrak{M}}_{v})^{\Gamma_{F_{v}},\hat{\theta}}/Z(\hat{\mathfrak{R}}_{v})^{\Gamma_{F_{v}},\hat{\theta}}\to Z(\hat{\mathfrak{M}}'_{v})^{\Gamma_{F_{v}}}/Z(\hat{\mathfrak{R}}'_{v})^{\Gamma_{F_{v}}}$$
 est surjectif et de noyau fini (car les donn\'ees $\boldsymbol{\mathfrak{M}}'_{v}$ et $\boldsymbol{\mathfrak{R}}'_{v}$ sont elliptiques par d\'efinition). On note $i_{\tilde{\mathfrak{M}}'_{v}}(\hat{\mathfrak{R}}_{v},\tilde{\mathfrak{R}}'_{v})$ l'inverse du nombre d'\'el\'ements de son noyau. De m\^eme, soit $\tilde{s}\in \tilde{\zeta}{\cal Z}(\hat{\mathfrak{R}}_{V})/Z(\hat{G})^{\Gamma_{F},\hat{\theta}}$ tel que ${\bf L}'(\tilde{s})$ soit elliptique. Alors l'homomorphisme
 $$Z(\hat{R}_{v})^{\Gamma_{F_{v}},\hat{\theta}}/Z(\hat{L}_{v})^{\Gamma_{F_{v}},\hat{\theta}}\to Z(\hat{R}'_{v})^{\Gamma_{F_{v}}}/Z(\hat{L}'_{v}(\tilde{s}))^{\Gamma_{F_{v}}}$$
 est surjectif et de noyau fini. On note   $i_{\tilde{R}'_{v}}(\hat{L}_{v},\tilde{L}'_{v}(\tilde{s}))$ l'inverse du nombre d'\'el\'ements de son noyau. Soit $\tilde{s}$ comme ci-dessus. On va montrer
 
 (5) si $D$ n'est pas un isomorphisme, $E(\hat{L}_{V},\tilde{s})=0$;
  
 (6) supposons que $D$ soit un isomorphisme; alors
 $$E(\hat{L}_{V},\tilde{s})=e(D)\prod_{v\in V}i_{\tilde{\mathfrak{M}}'_{v}}(\hat{\mathfrak{R}}_{v},\tilde{\mathfrak{R}}'_{v})i_{\tilde{R}'_{v}}(\hat{L}_{v},\tilde{L}'_{v}(\tilde{s})).$$
 
Supposons $E(\hat{L}_{V},\tilde{s})\not=0$. Alors $i_{\tilde{M}'}(\tilde{G},\tilde{G}'(\tilde{s}))\not=0$ donc ${\bf G}'(\tilde{s})$ est elliptique. On peut fixer $\hat{\mathfrak{L}}_{V}\in {\cal L}(\hat{\mathfrak{R}}_{V})\cap {\cal L}(\hat{L}_{V})$ de sorte que $\boldsymbol{\mathfrak{L}}'_{V}(\tilde{s})$ soit elliptique et que les constantes $d_{\tilde{M}'_{V}}^{\tilde{G}'(\tilde{s})}(\tilde{M}',\tilde{\mathfrak{L}}'_{V}(\tilde{s}))$ et $d_{\tilde{R}'_{v}}^{\tilde{\mathfrak{L}}'_{v}(\tilde{s})}(\tilde{\mathfrak{R}}'_{v},\tilde{L}'_{v}(\tilde{s}))$, pour $v\in V$, soient non nulles. Notons par exemple ${\cal A}_{\tilde{R}_{v}}^{\tilde{\mathfrak{L}}_{v}}$ l'orthogonal de ${\cal A}_{\tilde{\mathfrak{L}}_{v}}$ dans ${\cal A}_{\tilde{R}_{v}}$ et ${\cal A}_{\tilde{R}_{V}}^{\tilde{\mathfrak{L}}_{V}}=\oplus_{v\in V}{\cal A}_{\tilde{R}_{v}}^{\tilde{\mathfrak{L}}_{v}}$. 
Toutes les donn\'ees endoscopiques intervenant \'etant elliptiques (${\bf M}'_{v}$ \'etant consid\'er\'ee comme une donn\'ee de $\hat{\mathfrak{M}}_{v}$), ces non-nullit\'es signifient que l'on a les \'egalit\'es
$$(7) \qquad {\cal A}_{\tilde{\mathfrak{M}}_{V}}^{\tilde{G}}=\Delta({\cal A}_{\tilde{M}}^{\tilde{G}})\oplus {\cal A}_{\tilde{\mathfrak{L}}_{V}}^{\tilde{G}}$$
  et
$$(8) \qquad {\cal A}_{\tilde{R}_{V}}^{\tilde{\mathfrak{L}}_{V}}={\cal A}_{\tilde{\mathfrak{R}}_{V}}^{\tilde{\mathfrak{L}}_{V}}\oplus {\cal A}_{\tilde{L}_{V}}^{\tilde{\mathfrak{L}}_{V}},$$
  L'application $D$ est la compos\'ee des deux applications
 $${\cal A}_{\tilde{M}}^{\tilde{G}}\stackrel{\Delta}{\to }{\cal A}_{\tilde{\mathfrak{M}}_{V}}^{\tilde{G}}\to 
 {\cal A}_{\tilde{\mathfrak{M}}_{V}}/{\cal A}_{\tilde{\mathfrak{L}}_{V}}$$
 et
 $$(9) \qquad {\cal A}_{\tilde{\mathfrak{M}}_{V}}/{\cal A}_{\tilde{\mathfrak{L}}_{V}}\to  {\cal A}_{\tilde{\mathfrak{M}}_{V}}/{\cal A}_{\tilde{\mathfrak{R}}_{V}}\oplus {\cal A}_{\tilde{R}_{V}}/{\cal A}_{\tilde{L}_{V}}.$$
 La premi\`ere est un isomorphisme d'apr\`es (7). On peut d\'ecomposer l'espace de d\'epart de (9) en ${\cal A}_{\tilde{\mathfrak{M}}_{V}}^{\tilde{\mathfrak{R}}_{V}}\oplus {\cal A}_{\tilde{\mathfrak{R}}_{V}}/{\cal A}_{\tilde{\mathfrak{L}}_{V}}$. Alors (9) devient une application triangulaire. Les termes diagonaux sont les applications
 $${\cal A}_{\tilde{\mathfrak{M}}_{V}}^{\tilde{\mathfrak{R}}_{V}}\to {\cal A}_{\tilde{\mathfrak{M}}_{V}}/{\cal A}_{\tilde{\mathfrak{R}}_{V}}$$ 
 et
 $${\cal A}_{\tilde{\mathfrak{R}}_{V}}/{\cal A}_{\tilde{\mathfrak{L}}_{V}}\to {\cal A}_{\tilde{R}_{V}}/{\cal A}_{\tilde{L}_{V}}.$$
 La premi\`ere est \'evidemment un isomorphisme et la seconde l'est d'apr\`es (8). Donc (9) est un isomorphisme et $D$ aussi. Cela d\'emontre (5). Remarquons qu'en pr\'ecisant ces calculs,   on obtient l'\'egalit\'e
 $$(10) \qquad d(D)=d_{\tilde{M}'_{V}}^{\tilde{G}'(\tilde{s})}(\tilde{M}',\tilde{\mathfrak{L}}'_{V}(\tilde{s}))\prod_{v\in V}d_{\tilde{R}'_{v}}^{\tilde{\mathfrak{L}}'_{v}(\tilde{s})}(\tilde{\mathfrak{R}}'_{v},\tilde{L}'_{v}(\tilde{s})).$$
 Inversement, supposons que $D$ soit un isomorphisme. Il y a au plus un espace $\hat{\mathfrak{L}}_{V}$ qui peut contribuer \`a la somme $E(\hat{L}_{V},\tilde{s})$. En effet, l'\'egalit\'e (8) doit \^etre v\'erifi\'ee par cet espace, ce qui d\'etermine 
 $$(11) \qquad {\cal A}_{\tilde{\mathfrak{L}}_{V}}={\cal A}_{\tilde{\mathfrak{R}}_{V}}\cap {\cal A}_{\tilde{L}_{V}}.$$
 D\'efinissons ainsi cet espace. Pour qu'il intervienne vraiment, la donn\'ee $\boldsymbol{\mathfrak{L}}'_{V}(\tilde{s})$ doit \^etre elliptique. Par hypoth\`ese, les donn\'ees $\boldsymbol{\mathfrak{R}}'_{V}$ et ${\bf L}'(\tilde{s})$ sont elliptiques. Le membre de droite de (11) est donc \'egal \`a
 $${\cal A}_{\tilde{\mathfrak{R}}'_{V}}\cap {\cal A}_{\tilde{L}'_{V}(\tilde{s})}.$$
 L'espace ${\cal A}_{\tilde{\mathfrak{L}}'_{V}(\tilde{s})}$ est inclus dans cette intersection. Il est donc inclus dans le membre de gauche de (11), ce qui implique que $\boldsymbol{\mathfrak{L}}'_{V}(\tilde{s})$ est elliptique. En inversant le calcul fait ci-dessus, on voit que l'\'egalit\'e (8) et l'hypoth\`ese que $D$ est un isomorphisme impliquent (7). Puisque $\boldsymbol{\mathfrak{L}}'_{V}(\tilde{s})$ est elliptique, le dernier terme de cette \'egalit\'e est \'egal \`a ${\cal A}_{\tilde{\mathfrak{L}}'_{V}(\tilde{s})}^{\tilde{G}}$. Cet espace contient $\Delta({\cal A}_{\tilde{G}'(\tilde{s})}^{\tilde{G}})$. Puisque les deux termes du membre de droite de (7) sont en somme directe, l'espace ${\cal A}_{\tilde{G}'(\tilde{s})}^{\tilde{G}}$ est nul. Donc la donn\'ee ${\bf G}'(\tilde{s})$ est elliptique. Consid\'erons la d\'efinition (4), o\`u n'intervient plus que notre Levi $\tilde{\mathfrak{L}}'_{V}(\tilde{s})$. Une constante telle que $e_{\tilde{M}'_{V}}^{\tilde{G}'(\tilde{s})}(\tilde{M}',\tilde{\mathfrak{L}}'_{V}(\tilde{s}))$ est le produit de  $d_{\tilde{M}'_{V}}^{\tilde{G}'(\tilde{s})}(\tilde{M}',\tilde{\mathfrak{L}}'_{V}(\tilde{s}))$ et de l'inverse de  $k_{\tilde{M}'_{V}}^{\tilde{G}'(\tilde{s})}(\tilde{M}',\tilde{\mathfrak{L}}'_{V}(\tilde{s}))$, ce terme \'etant le nombre d'\'el\'ements d'un certain noyau. L'\'egalit\'e (10) montre que le produit des constantes $d$ est \'egal \`a $d(D)$. Pour prouver (6), il suffit de prouver l'\'egalit\'e
 $$(12) \qquad k(D)^{-1}\prod_{v\in V}i_{\tilde{\mathfrak{M}}'_{v}}(\hat{\mathfrak{R}}_{v},\tilde{\mathfrak{R}}'_{v})i_{\tilde{R}'_{v}}(\hat{L}_{v},\tilde{L}'_{v}(\tilde{s}))=i_{\tilde{M}'}(\tilde{G},\tilde{G}'(\tilde{s})) k_{\tilde{M}'_{V}}^{\tilde{G}'(\tilde{s})}(\tilde{M}',\tilde{\mathfrak{L}}'_{V}(\tilde{s}))^{-1}$$
 $$\prod_{v\in V}k_{\tilde{R}'_{v}}^{\tilde{\mathfrak{L}}'_{v}(\tilde{s})}(\tilde{\mathfrak{R}}'_{v},\tilde{L}'_{v}(\tilde{s}))^{-1}.$$
 On a un diagramme commutatif
 $$\begin{array}{ccc} Z(\hat{M})^{\Gamma_{F},\hat{\theta}}/Z(\hat{G})^{\Gamma_{F},\hat{\theta}}&\stackrel{\hat{D}}{\to}&\oplus_{v\in V}(Z(\hat{\mathfrak{M}}_{v})^{\Gamma_{F_{v}},\hat{\theta}}/Z(\hat{\mathfrak{R}}_{v})^{\Gamma_{F_{v}},\hat{\theta}}\oplus Z(\hat{R}_{v})^{\Gamma_{F_{v}},\hat{\theta}}/Z(\hat{L}_{v})^{\Gamma_{F_{v}},\hat{\theta}})\\ \downarrow&&\downarrow\\ 
Z(\hat{M}')^{\Gamma_{F}}/Z(\hat{G}'(\tilde{s}))^{\Gamma_{F}}&\stackrel{\hat{D}'}{\to}&
\oplus_{v\in V}(Z(\hat{M}'_{v})^{\Gamma_{F_{v}}}/Z(\hat{\mathfrak{R}}'_{v})^{\Gamma_{F_{v}}}\oplus Z(\hat{R}'_{v})^{\Gamma_{F_{v}}}/Z(\hat{L}'_{v}(\tilde{s}))^{\Gamma_{F_{v}}})\\ \end{array}$$
  Tous les homomorphismes sont surjectifs et de noyaux finis. Calculons l'inverse du nombre d'\'el\'ements de l'homomorphisme compos\'e. Si on utilise le chemin nord-est, on trouve le membre de gauche de (12). Si on utilise le chemin sud-ouest, on trouve $i_{\tilde{M}'}(\tilde{G},\tilde{G}'(\tilde{s}))\vert ker(\hat{D}')\vert ^{-1}$. L'homomorphisme $\hat{D}'$ se d\'ecompose en
  $$Z(\hat{M}')^{\Gamma_{F}}/Z(\hat{G}'(\tilde{s}))^{\Gamma_{F}}\to \oplus_{v\in V}Z(\hat{M}'_{v})^{\Gamma_{F_{v}}}/Z(\hat{\mathfrak{L}}'_{v}(\tilde{s}))^{\Gamma_{F_{v}}}$$
  $$\stackrel{\oplus_{v\in V}\hat{\iota}_{v}}{\to}
\oplus_{v\in V}(Z(\hat{M}'_{v})^{\Gamma_{F_{v}}}/Z(\hat{\mathfrak{R}}'_{v})^{\Gamma_{F_{v}}}\oplus Z(\hat{R}'_{v})^{\Gamma_{F_{v}}}/Z(\hat{L}'_{v}(\tilde{s}))^{\Gamma_{F_{v}}}).$$
De nouveau, les homomophismes sont surjectifs et de noyaux finis. Le nombre d'\'el\'ements du noyau du premier homomorphisme est \'egal \`a $k_{\tilde{M}'_{V}}^{\tilde{G}'(\tilde{s})}(\tilde{M}',\tilde{\mathfrak{L}}'_{V}(\tilde{s}))$. Pour obtenir (12), il reste \`a prouver que, pour tout $v\in V$, on a
$$(13) \qquad \vert ker(\hat{\iota}_{v})\vert =k_{\tilde{R}'_{v}}^{\tilde{\mathfrak{L}}'_{v}(\tilde{s})}(\tilde{\mathfrak{R}}'_{v},\tilde{L}'_{v}(\tilde{s})).$$
  On  a un diagramme commutatif
  $$\begin{array}{ccccccc}1&&&1&&&\\ \downarrow&&&\downarrow&&&\\ Z(\hat{M}'_{v})^{\Gamma_{F_{v}}}/Z(\hat{\mathfrak{L}}'_{v}(\tilde{s}))^{\Gamma_{F_{v}}}&\stackrel{\hat{\iota}_{v}}{\to}&(Z(\hat{M}'_{v})^{\Gamma_{F_{v}}}/Z(\hat{\mathfrak{R}}'_{v})^{\Gamma_{F_{v}}}&\oplus &Z(\hat{R}'_{v})^{\Gamma_{F_{v}}}/Z(\hat{L}'_{v}(\tilde{s}))^{\Gamma_{F_{v}}})&\to&1\\ \downarrow&&&\downarrow&&&\\ Z(\hat{R}'_{v})^{\Gamma_{F_{v}}}/Z(\hat{\mathfrak{L}}'_{v}(\tilde{s}))^{\Gamma_{F_{v}}}&\stackrel{\hat{\kappa}_{v}}\to&(Z(\hat{R}'_{v})^{\Gamma_{F_{v}}}/Z(\hat{\mathfrak{R}}'_{v})^{\Gamma_{F_{v}}}&\oplus& Z(\hat{R}'_{v})^{\Gamma_{F_{v}}}/Z(\hat{L}'_{v}(\tilde{s}))^{\Gamma_{F_{v}}})&\to&1\\ \downarrow&&&\downarrow&&&\\ Z(\hat{R}'_{v})^{\Gamma_{F_{v}}}/Z(\hat{M}'_{v})^{\Gamma_{F_{v}}}&=& (Z(\hat{R}'_{v})^{\Gamma_{F_{v}}}/Z(\hat{M}'_{v})^{\Gamma_{F_{v}}}&\oplus&1)&&\\ \downarrow&&&\downarrow&&&\\ 1&&&1&&&\\ \end{array}$$
  Les lignes et colonnes sont exactes. Il en resulte que  $\vert ker(\hat{\iota}_{v})\vert =\vert ker(\hat{\kappa}_{v})\vert $. Mais ce dernier nombre est \'egal au membre de droite de (13). Cela prouve (13), (12) et ach\`eve la preuve de (6). 
  
  Si $D$ n'est pas un isomorphisme, (5) entra\^{\i}ne que $J(\hat{\mathfrak{R}}_{V},\hat{L}_{V},\boldsymbol{\delta},f)=0$ et on a fini. On suppose maintenant que $D$ est un isomorphisme. On peut supposer qu'il existe $\tilde{s}$ intervenant dans la d\'efinition (3) tel que $(f_{v}^{\tilde{G}'_{1,v}(\tilde{s})})_{\tilde{L}'_{1,v}(\tilde{s})}$ soit non nul pour tout $v\in V$. Cela entra\^{\i}ne que le Levi $\tilde{L}'_{v}(\tilde{s})$ est relevant. A fortiori, $\hat{L}_{v}$ correspond \`a un $K$-espace de Levi $K\tilde{L}_{v}$ de $K\tilde{G}_{v}$. On peut fixer une collection $K\tilde{L}_{V}=(K\tilde{L}_{v})_{v\in V}$ de tels $K$-espaces de Levi et identifier ${\bf L}_{v}(\tilde{s})$ \`a une donn\'ee endoscopique elliptique de $(KL_{V},K\tilde{L}_{V},{\bf a}_{L})$. On a pour tout $v$ l'\'egalit\'e $(f_{v}^{\tilde{G}'_{1,v}(\tilde{s})})_{\tilde{L}'_{1,v}(\tilde{s})}=(f_{v,K\tilde{L}_{v},\omega_{v}})^{\tilde{L}'_{1,v}(\tilde{s})}$. On se sert ici de la restriction du facteur de transfert fix\'e sur $\tilde{G}'_{1}(\tilde{s},F_{v})\times \tilde{G}(F_{v})$ pour d\'efinir le transfert.   Rappelons que le terme $i_{\tilde{R}'_{v}}(\tilde{L}_{v},\tilde{L}'_{v}(\tilde{s}))$   est nul si ${\bf L}'_{V}(\tilde{s})$ n'est pas elliptique. Compte tenu de l'\'egalit\'e (6), la d\'efinition (3) se r\'ecrit
  $$(14) \qquad J(\hat{\mathfrak{R}}_{V},\hat{L}_{V},\boldsymbol{\delta},f)=c\sum_{\tilde{s}\in Z(\hat{\mathfrak{R}}_{V})/Z(\hat{G})^{\Gamma_{F},\hat{\theta}}}X(\tilde{s})\prod_{v\in V} S_{\tilde{R}'_{1,v}(\tilde{s}),\lambda_{1,v}(\tilde{s})}^{\tilde{L}'_{1,v}(\tilde{s})}({\bf d}_{v}(\tilde{s}),B^{\tilde{G}}(f_{v,K\tilde{L}_{v},\omega_{v}})^{\tilde{L}'_{1,v}(\tilde{s})}),$$
  o\`u
  $$c=e(D)\prod_{v\in V}i_{\tilde{\mathfrak{M}}'_{v}}(\hat{\mathfrak{R}}_{v},\tilde{\mathfrak{R}}'_{v})$$
  et 
  $$X(\tilde{s})=\prod_{v\in V}i_{\tilde{R}'_{v}}(\tilde{L}_{v},\tilde{L}'_{v}(\tilde{s})).$$
  Pour toute place $v\in V$, introduisons le groupe $Z(\hat{R}_{v})_{*}$ image r\'eciproque dans $Z(\hat{R})$ de $(Z(\hat{R}_{v})/Z(\hat{R}_{v})\cap  Z(\hat{R}'_{v}))^{\Gamma_{F_{v}}}$. Le sous-groupe des \'el\'ements invariants par $\hat{\theta}$ dans 
  $$\cap_{v\in V}Z(\hat{L}_{v})^{\Gamma_{F_{v}}}(1-\hat{\theta})(Z(\hat{R}_{v})_{*})$$ contient $\cap_{v\in V}Z(\hat{L}_{v})^{\Gamma_{F_{v}},\hat{\theta}}$ comme sous-groupe d'indice fini. Le groupe
 $$({\cal Z}(\hat{\mathfrak{R}}_{V})\cap\cap_{v\in V}Z(\hat{L}_{v})^{\Gamma_{F_{v}},\hat{\theta}}) /Z(\hat{G})^{\Gamma_{F},\hat{\theta}}$$
 est le noyau de $\hat{D}$ donc est fini. Il en r\'esulte que le groupe
$$\left({\cal Z}(\hat{\mathfrak{R}}_{V})\cap \cap_{v\in V}Z(\hat{L}_{v})^{\Gamma_{F_{v}}}(1-\hat{\theta})(Z(\hat{R}_{v})_{*})\right)/Z(\hat{G})^{\Gamma_{F},\hat{\theta}}$$
est fini. Il nous suffit d'en trouver un sous-groupe ${\cal Z}$ v\'erifiant la propri\'et\'e suivante.  Fixons $\tilde{s}_{0}\in \tilde{\zeta}{\cal Z}(\hat{\mathfrak{R}}_{V})/Z(\hat{G})^{\Gamma_{F},\hat{\theta}}$. Alors, dans l'expression (14),  la sous-somme sur $\tilde{s}\in {\cal Z}\tilde{s}_{0}$ est nulle. Consid\'erons donc un sous-groupe ${\cal Z}$ que l'on pr\'ecisera plus tard. La propri\'et\'e requise d\'epend d'un \'el\'ement $\tilde{s}_{0}$. Pour la simplicit\'e de l'\'ecriture, on peut supposer $\tilde{s}_{0}=\tilde{\zeta}$. On a

(15) pour $\tilde{s}\in {\cal Z}\tilde{\zeta}$, les donn\'ees endoscopiques ${\bf L}'_{V}(\tilde{s})$ et ${\bf L}'_{V}={\bf L}'_{V}(\tilde{\zeta})$ sont \'equivalentes; si ces donn\'ees sont elliptiques, on a $X(\tilde{s})=X(\tilde{\zeta})$.

 Soit $z$ un repr\'esentant dans $Z(\hat{M})$ d'un \'el\'ement de ${\cal Z}$. Pour toute place $v\in V$, on \'ecrit $z=\tau_{v}(1-\hat{\theta})(\rho_{v})$, avec $\tau_{v}\in Z(\hat{L}_{v})^{\Gamma_{F_{v}}}$ et $\rho_{v}\in Z(\hat{R}_{v})_{*}$. Puisque $\rho_{v}\in \hat{R}_{v}\subset \hat{L}_{v}$,  $ad_{\rho_{v}}$ est un automorphisme int\'erieur de $\hat{L}_{v}$. On a l'\'egalit\'e $z\tilde{\zeta}=\rho_{v}\tau_{v}\tilde{\zeta}\rho_{v}^{-1}$, donc $ad_{\rho_{v}}$ envoie $\hat{L}'_{v}(\tau_{v}\tilde{\zeta})$ sur $\hat{L}'_{v}(z\tilde{\zeta})$. Puisque $\tau_{v}\in Z(\hat{L}_{v})$, on a $\hat{L}'_{v}(\tau_{v}\tilde{\zeta})=\hat{L}'_{v}$.   Les groupes ${\cal L}'_{v}$ et ${\cal L}'_{v}(z\tilde{\zeta})$ sont engendr\'es par ${\cal R}'_{v}={\cal M}_{v}'\cap {^LR}_{v}$ et $\hat{L}'_{v}$, respectivement $\hat{L}'_{v}(z\tilde{\zeta})$. Puisque $\rho_{v}\in Z(\hat{R}_{v})_{*}$, $ad_{\rho_{v}}$ conserve ${\cal R}'_{v}$. Donc $ad_{\rho_{v}}$ envoie ${\cal L}'_{v}$ sur ${\cal L}'_{v}(z\tilde{\zeta})$. Autrement dit, $\rho_{v}$ d\'efinit une \'equivalence entre les donn\'ees ${\bf L}'_{v}$ et ${\bf L}'_{v}(z\tilde{\zeta})$. L'\'egalit\'e $X(\tilde{s})=X(\tilde{\zeta})$  r\'esulte \'evidemment de la d\'efinition de ces termes. Cela prouve (15).

On peut supposer ${\bf L}'_{V}$ elliptique, sinon $X(\tilde{\zeta})=0$.  La somme que l'on consid\`ere est proportionnelle \`a
 $$\sum_{z\in {\cal Z}}\prod_{v\in V}S_{\tilde{R}'_{1,v}(z\tilde{\zeta}),\lambda_{1,v}(z\tilde{\zeta})}^{\tilde{L}'_{1,v}(z\tilde{\zeta})}({\bf d}_{v}(z\tilde{\zeta}),B^{\tilde{G}}, (f_{\tilde{L}_{v},\omega})^{\tilde{L}'_{1,v}(z\tilde{\zeta})} ).$$
 Fixons $z\in Z(\hat{M})^{\Gamma_{F},\hat{\theta}}$ relevant un \'el\'ement de ${\cal Z}$. Puisque $Z(\hat{M})^{\Gamma_{F},\hat{\theta}}=Z(\hat{M})^{\Gamma_{F},\hat{\theta},0}Z(\hat{G})^{\Gamma_{F},\hat{\theta}}$, on peut supposer $z\in Z(\hat{M})^{\Gamma_{F},\hat{\theta},0}$. Comme dans la preuve de (15), on \'ecrit $z=\tau_{v}(1-\hat{\theta})(\rho_{v})$ pour toute place $v\in V$, o\`u $\tau_{v}\in Z(\hat{L}_{v})^{\Gamma_{F_{v}}}$ et $\rho_{v}\in Z(\hat{R}_{v})_{*}$. Les deux s\'eries de donn\'ees auxiliaires $M'_{1}=M'_{1}(\tilde{\zeta})$,... et $M'_{1}(z\tilde{\zeta})$ fournissent une fonction de recollement $\tilde{\lambda}(z)^M$ sur le produit fibr\'e $\tilde{M}'_{1,V}(F_{V})\times \tilde{M}'_{1,V}(z\tilde{\zeta};F_{V})$.  Par restriction puis dualit\'e, on en d\'eduit un isomorphisme
 $$\iota(z)^{M,*}:\prod_{v\in V}D_{g\acute{e}om,\lambda_{1,v}(z\tilde{\zeta})}^{st}(\tilde{R}'_{1,v}(z\tilde{\zeta};F_{v}))\simeq \prod_{v\in V}D_{g\acute{e}om,\lambda_{1,v}}^{st}(\tilde{R}'_{1,v}(F_{v})).$$
En posant ${\bf d}=\otimes_{v\in V}{\bf d}_{v}$, on a par d\'efinition, ${\bf d}=\iota(z)^{M,*}({\bf d}(z\tilde{\zeta}))$.
 
 Pour tout $v\in V$, on compl\`ete la paire de Borel $(\hat{B}\cap \hat{L}'_{v},\hat{T}^{\hat{\theta},0})$ de $\hat{L}'_{v}$ en une paire de Borel \'epingl\'ee invariante par $\Gamma_{F_{v}}$. De m\^eme pour $\hat{L}'_{v}(z\tilde{\zeta})$. On peut supposer que les restrictions \`a $\hat{R}'_{v}$ de ces deux \'epinglages co\"{\i}ncident. Quitte \`a multiplier $\rho_{v}$ par un \'el\'ement de $Z(\hat{R}_{v})\cap \hat{T}^{\hat{\theta},0}$, ce qui ne change pas $(1-\hat{\theta})(\rho_{v})$, on peut supposer que $ad_{\rho_{v}}$ \'echange ces deux \'epinglages. Alors $ad_{\rho_{v}}$ devient \'equivariant pour les actions galoisiennes. On peut identifier les groupes $L'_{v}$ et $L'_{v}(z\tilde{\zeta})$ ainsi que les espaces $\tilde{L}'_{v}$ et $\tilde{L}'_{v}(z\tilde{\zeta})$. Les donn\'ees auxiliaires $L'_{1,v}(z\tilde{\zeta})$ etc... se transportent en des donn\'ees auxiliaires pour ${\bf L}'_{v}$. Le plongement de ${\cal L}'_{v}$ dans $^LL'_{1,v}(z\tilde{\zeta})$ est le compos\'e de $ad_{\rho_{v}}$ et du plongement $\hat{\xi}_{1}(z\tilde{\zeta}):{\cal L}'_{v}(z\tilde{\zeta})\to {^LL}_{1,v}(z\tilde{\zeta})$. Le facteur de transfert n'est bien d\'efini que sur le produit de ces donn\'ees auxiliaires sur toutes les places $v\in V$. C'est alors la restriction du facteur canonique $\Delta_{1}(z\tilde{\zeta})$ issu de la donn\'ee ${\bf G}'(z\tilde{\zeta})$. On se retrouve avec deux s\'eries de donn\'ees auxiliaires pour ${\bf L}'_{V}$, d'o\`u une fonction de recollement $\tilde{\lambda}(z)^L$. Il s'en d\'eduit encore un isomorphisme
 $$\iota(z)^{L,*}:\prod_{v\in V}D_{g\acute{e}om,\lambda_{1,v}(z\tilde{\zeta})}^{st}(\tilde{R}'_{1,v}(z\tilde{\zeta};F_{v}))\simeq \prod_{v\in V}D_{g\acute{e}om,\lambda_{1,v}}^{st}(\tilde{R}'_{1,v}(F_{v})).$$
 Le transfert commute au recollement donc celui-ci envoie $\otimes_{v\in V}(f_{\tilde{L}_{v},\omega})^{\tilde{L}'_{1,v}(z\tilde{\zeta})}$ sur $\otimes_{v\in V}(f_{\tilde{L}_{v},\omega})^{\tilde{L}'_{1,v}}$. Il envoie ${\bf d}(z\tilde{\zeta})$ sur $\iota(z)^{L,*}({\bf d}(z\tilde{\zeta}))=\iota(z)^{L,*}\circ (\iota(z)^{M,*})^{-1}({\bf d})$. En d\'ecomposant nos isomorphismes de recollement en produits tensoriels sur les $v\in V$, on obtient l'\'egalit\'e
$$\prod_{v\in V}S_{\tilde{R}'_{1,v}(z\tilde{\zeta}),\lambda_{1,v}(z\tilde{\zeta})}^{\tilde{L}'_{1,v}(z\tilde{\zeta})}({\bf d}_{v}(z\tilde{\zeta}),B^{\tilde{G}},(f_{\tilde{L}_{v},\omega})^{\tilde{L}'_{1,v}(z\tilde{\zeta})} )=$$
$$\prod_{v\in V}S_{\tilde{R}'_{1,v},\lambda_{1,v}}^{\tilde{L}'_{1,v}}((\iota(z)_{v}^{L,*}\circ (\iota(z)_{v}^{M,*})^{-1}({\bf d}_{v}),B^{\tilde{G}},(f_{\tilde{L}_{v},\omega})^{\tilde{L}'_{1,v}} ).$$
Il nous suffit de prouver que
$$\sum_{z\in {\cal Z}}\iota(z)^{L,*}\circ (\iota(z)^{M,*})^{-1}({\bf d})=0.$$
L'automorphisme $(\iota(z)^{M})^{-1}\circ \iota(z)^{L}$ de $\prod_{v\in V}C_{c,\lambda_{1,v}}^{\infty}(\tilde{R}'_{1,v}(F_{v}))$ est de la forme $\varphi\mapsto \tilde{\lambda}_{z}\varphi$, o\`u $\tilde{\lambda}_{z}$ est une fonction lisse sur $\prod_{v\in V}\tilde{R}'_{1,v}(F_{v})$. Il nous suffit encore de prouver que

(16) $\sum_{z\in {\cal Z}}\tilde{\lambda}_{z}(\gamma)=0$ pour tout $\gamma\in \prod_{v\in V}\tilde{R}'_{1,v}(F_{v})$ dans un voisinage invariant par conjugaison stable de l'\'el\'ement $\epsilon=\prod_{v\in V}\epsilon_{v}$ fix\'e plus haut. 

On fixe de nouveau $z$ que l'on \'ecrit comme ci-dessus. On va calculer $\tilde{\lambda}_{z}$.  On simplifie les notations en posant $\tilde{\lambda}=\tilde{\lambda}_{z}$,  $\tilde{\lambda}^M=\tilde{\lambda}(z)^M$ etc... On pose
$\tilde{\zeta}_{1}=\tilde{\zeta}$ et $\tilde{\zeta}_{2}=z\tilde{\zeta}$. On  supprime autant que c'est possible ces termes de la notation. Les donn\'ees relatives \`a $\tilde{\zeta}_{1}$ seront affect\'ees d'un indice $1$ et celles relatives \`a $\tilde{\zeta}_{2}$ d'un indice $2$ (par exemple, on note $\tilde{R}'_{2,v}$ l'espace not\'e pr\'ec\'edemment $\tilde{R}'_{1,v}(z\tilde{\zeta})$). Soit $r'\in\prod_{v\in V} \tilde{R}'_{v}(F_{v})$ et, pour $i=1,2$, soit $r'_{i}\in \prod_{v\in V}\tilde{R}'_{i,v}(F_{v})$ un \'el\'ement au-dessus de $r'$.  Par d\'efinition,
$$\tilde{\lambda}(r'_{1})=\tilde{\lambda}^M(r'_{1},r'_{2})^{-1}\tilde{\lambda}^L(r'_{1},r'_{2}),$$
o\`u $\tilde{\lambda}^M$ et $\tilde{\lambda}^L$ sont les fonctions de recollement introduites ci-dessus. Fixons $m'\in \tilde{M}'(F)$ et, pour $i=1,2$, un \'el\'ement $m'_{i}\in \tilde{M}'_{i}(F)$ au-dessus de $m'$. Soit $(b_{1},b_{2})$ l'\'el\'ement du produit fibr\'e $M'_{1}(F_{V})\times_{M'(F_{V})}M'_{2}(F_{V})$ tel que $(r'_{1},r'_{2})=(b_{1}m'_{1},b_{2}m'_{2})$. On a l'\'egalit\'e
$$\tilde{\lambda}^M(r'_{1},r'_{2})=\lambda^M(b_{1},b_{2})\tilde{\lambda}^M(m'_{1},m'_{2}).$$
On se rappelle la d\'efinition de $\tilde{\lambda}^M$. Parce que les \'el\'ements $m'_{1}$ et $m'_{2}$ sont d\'efinis sur $F$, on a
$$\tilde{\lambda}^M(m'_{1},m'_{2})=\prod_{v\not\in V}\tilde{\lambda}^M_{v}(m'_{1},m'_{2})^{-1},$$
o\`u les $\tilde{\lambda}^M_{v}$ pour $v\not\in V$ sont les fonctions de recollement associ\'ees aux espaces hypersp\'eciaux $\tilde{K}_{i,v}\cap \tilde{M}'_{i}(F_{v})$ pour $i=1,2$. 
Puisqu'on a suppos\'e que ${\bf L}'_{v}={\bf L}'_{v}(\tilde{\zeta}_{1})\simeq {\bf L}'_{v}(\tilde{\zeta}_{2})$ \'etait relevant pour tout $v\in V$, on peut fixer $l\in \prod_{v\in V}K\tilde{L}_{v}(F_{v})$ assez r\'egulier et $l'\in \prod_{v\in V}\tilde{L}'_{v}(F_{v})$ dont les classes de conjugaison stable se correspondent. On note $\tilde{L}_{v}$ la composante de $K\tilde{L}_{v}$ telle que $l$ appartienne \`a  $\prod_{v\in V}\tilde{L}_{v}(F_{v})$. On note $\tilde{G}$ la composante de $K\tilde{G}$ telle que $\tilde{L}_{v}\subset \tilde{G}_{v}$ pour tout $v\in V$ et on note $\tilde{H}$ la composante de $K\tilde{H}$ qui contient $\tilde{G}$. Pour $i=1,2$ on fixe $l'_{i}\in \prod_{v\in V}\tilde{L}'_{i,v}(F_{v})$ au-dessus de $l'$. Notons $(a_{1},a_{2})$ l'\'el\'ement du produit fibr\'e $\prod_{v\in V}L'_{1,v}(F_{v})\times_{L'(F_{v})}L'_{2,v}(F_{v})$ tel que $(r'_{1},r'_{2})=(a_{1}l'_{1},a_{2}l'_{2})$. On a l'\'egalit\'e
$$\tilde{\lambda}^L(r'_{1},r'_{2})=\lambda^L(a_{1},a_{2})\tilde{\lambda}^L(l'_{1},l'_{2})=\lambda^L(a_{1},a_{2})\Delta_{2}(l'_{2},l)\Delta_{1}(l'_{1},l)^{-1}.$$
Pour calculer les facteurs de transfert, on utilise la d\'efinition de 3.9.  Pour $i=1,2$, on a rel\`eve la donn\'ee endoscopique ${\bf G}'(\tilde{\zeta}_{i})$  en une donn\'ee de $(KH,K\tilde{H},{\bf b})$, que l'on note pour simplifier ${\bf H}'(\tilde{\zeta}_{i})$, munie de donn\'ees auxiliaires.  On choisit un  couple comme en 3.6, qui est not\'e $(\delta_{1},\gamma)$ dans ce paragraphe, et que nous noterons ici $(h'_{i},h(\tilde{\zeta}_{i}))$. On dispose du facteur global que l'on note $\Delta_{i,glob}(h'_{i},h(\tilde{\zeta}_{i}))$ pour $i=1,2$. Par d\'efinition
$$\Delta_{i}(l'_{i},l)=\Delta_{i,glob}(h'_{i},h(\tilde{\zeta}_{i}))\left(\prod_{v\not\in V}\Delta_{i,v}(h'_{i},h(\tilde{\zeta}_{i}))^{-1}\right)\left(\prod_{v\in V}\boldsymbol{\Delta}_{i,v}(l'_{i},l;h'_{i},h(\tilde{\zeta}_{i}))\right).$$
A ce point, on obtient l'\'egalit\'e
$$(17) \qquad \tilde{\lambda}(r'_{1})=\lambda^M(b_{1},b_{2})^{-1}\lambda^L(a_{1},a_{2})\Delta_{2,glob}(h'_{2},h(\tilde{\zeta}_{2}))\Delta_{1,glob}(h'_{1},h(\tilde{\zeta}_{1}))^{-1}$$
$$\left(\prod_{v\not\in V}\tilde{\lambda}_{v}^M(m'_{1},m'_{2})\Delta_{1,v}(h'_{1},h(\tilde{\zeta}_{1}))\Delta_{2,v}(h'_{2},h(\tilde{\zeta}_{2}))^{-1}\right)$$
$$\left(\prod_{v\in V}\boldsymbol{\Delta}_{2,v}(l'_{2},l;h'_{2},h(\tilde{\zeta}_{2}))\boldsymbol{\Delta}_{1,v}(l'_{1},l';h'_{1},h(\tilde{\zeta}_{1}))^{-1}\right).$$

Pour tout $v\in V$, on fixe un diagramme $(l',B^{L'}_{v},T^{L'}_{v},B^L_{v},T^L_{v},l)$ et  des $a$-data et des $\chi$-data relatives \`a ce diagramme. On suppose les $\chi$-data triviales sur les orbites asym\'etriques. On a la complication que l'on doit plonger $G$ et $G'(\tilde{\zeta}_{i})$ pour $i=1,2$ dans les groupes plus gros $H$ et $H'(\tilde{\zeta}_{i})$. Il est clair que notre diagramme se prolonge en des diagrammes  relatifs \`a ces groupes plus gros, que l'on note  $(l',B^{L'}_{v}(\tilde{\zeta}_{i}),T^{L'}_{v}(\tilde{\zeta}_{i}),B^{L,H}_{v},T^{L,H}_{v},l)$. On utilise les diagrammes prolong\'es et les m\^emes $a$-data et $\chi$-data pour calculer les facteurs de transfert intervenant ci-dessus. On doit aussi fixer des diagrammes 
$$(h'(\tilde{\zeta}_{i}),B'(\tilde{\zeta}_{i}),T'(\tilde{\zeta}_{i}), B(\tilde{\zeta}_{i}),T(\tilde{\zeta}_{i}),h(\tilde{\zeta}_{i}))$$
 pour $i=1,2$, o\`u $h'(\tilde{\zeta}_{i})$ est l'image de $h'_{i}$ dans $\tilde{H}'(\tilde{\zeta}_{i})$. On peut supposer et on suppose que le tore $T'(\tilde{\zeta}_{i})$ est d\'efini sur $F$ et que le groupe de Borel $B'(\tilde{\zeta}_{i})$ est d\'efini sur $\bar{F}$. Dans chaque facteur de transfert interviennent des facteurs $\Delta_{II}$. Ceux relatifs aux couples $(h'_{i},h(\tilde{\zeta}_{i}))$ disparaissent car leur intervention dans le facteur global compense leurs interventions dans les facteurs locaux. Ceux relatifs aux couples $(l'_{i},l)$ disparaissent aussi. En effet, pour chaque $v\in V$, la contribution aux facteurs $\boldsymbol{\Delta}_{i,v}$ des orbites galoisiennes dans $\tilde{L}_{v}$ est la m\^eme pour les deux facteurs. Les orbites hors de $\tilde{L}_{v}$ sont asym\'etriques et leur contribution est triviale d'apr\`es le choix de nos $\chi$-data. On peut donc remplacer chaque facteur de transfert par le facteur $\Delta_{imp}$ correspondant. 

Fixons provisoirement $v\in V$ et abandonnons les indices $v$ pour simplifier. Soit $i=1,2$. On introduit les tores
$$U_{i} =(T_{sc}^{L}\times T_{sc}(\tilde{\zeta}_{i}))/diag_{-}(Z(G_{SC}));$$

${\cal T}^{L}(\tilde{\zeta}_{i})$ le produit fibr\'e de $T^{L'}_{i}(\tilde{\zeta}_{i})$ et de $T^{L,H}$ au-dessus de $T^{L'}(\tilde{\zeta}_{i})$, o\`u $T^{L'}_{i}(\tilde{\zeta}_{i})$ est l'image r\'eciproque de $T^{L'}(\tilde{\zeta}_{i})$ dans la donn\'ee auxiliaire $H'_{i}(\tilde{\zeta}_{i})$;

 ${\cal T}(\tilde{\zeta}_{i})$ le produit fibr\'e de $T'_{i}(\tilde{\zeta}_{i})$ et de $T(\tilde{\zeta}_{i})$ au-dessus de $T'(\tilde{\zeta}_{i})$, o\`u $T'_{i}(\tilde{\zeta}_{i})$ est l'image r\'eciproque de $T'(\tilde{\zeta}_{i})$ dans la donn\'ee auxiliaire $H'_{i}(\tilde{\zeta}_{i})$;

 $S_{i}=({\cal T}^{L}(\tilde{\zeta}_{i})\times {\cal T}(\tilde{\zeta}_{i}))/diag_{-}({\cal Z}_{i}^H)$, o\`u ${\cal Z}^H_{i}$ est le produit fibr\'e de $Z(H'_{i}(\tilde{\zeta}_{i}))$ et de $Z(H)$ au dessus de $Z(H'(\tilde{\zeta}_{i}))$;

 le tore dual $\hat{U}_{i}=(\hat{T}_{sc}^{L}\times \hat{T}_{sc}(\tilde{\zeta}_{i}))/diag(Z(\hat{G}_{SC}))$;

 le tore dual $\hat{{\cal T}}^{L}(\tilde{\zeta}_{i})$ qui est le quotient de $\hat{T}^{L'}_{i}(\tilde{\zeta}_{i})$ et de $\hat{T}^{L,H}$ par $\hat{T}^{L'}(\tilde{\zeta}_{i})$ plong\'e par $t'\mapsto (\hat{\xi}^H_{i}(t')^{-1},t')$ (on note $\hat{\xi}_{i}^H:{\cal H}'(\tilde{\zeta}_{i})\to {^LH}'_{i}(\tilde{\zeta}_{i})$ le plongement fix\'e);

 le tore dual $\hat{{\cal T}}(\tilde{\zeta}_{i})$ qui se d\'ecrit de la m\^eme fa\c{c}on;

 le tore dual $\hat{S}_{i}$  qui se d\'ecrit comme le sous-groupe des $(t^L,t,t_{sc})\in \hat{{\cal T}}^{L}(\tilde{\zeta}_{i})\times \hat{{\cal T}}(\tilde{\zeta}_{i})\times \hat{T}_{sc}$ tels que $j(t_{sc})=t^Lt^{-1}$; on a ici identifi\'e tous les tores  \`a un tore commun, en oubliant leurs actions galoisiennes; $j$ est l'application naturelle $j:\hat{T}_{sc}\to \hat{{\cal T}}^{L}(\tilde{\zeta}_{i})\simeq \hat{{\cal T}}(\tilde{\zeta}_{i})$; on renvoie \`a [I] 2.2 pour la description de l'action galoisienne sur $\hat{S}_{i}$.

{\bf Remarque.} Dans la description de $U_{i}$ et $\hat{U}_{i}$ apparaissent a priori les groupes $H$ et $\hat{H}$. Mais ils n'apparaissent que via leurs rev\^etements simplement connexes, qui s'identifient \`a $G_{SC}$ et $\hat{G}_{SC}$.

\bigskip
On construit comme en [II] 2.2, 2.3:

- une cocha\^{\i}ne $(V_{i}^{L}, V_{i}^{-1}):\Gamma_{F_{v}}\to T_{sc}^{L}\times T_{sc}(\tilde{\zeta}_{i})$;

- un \'el\'ement $(\nu^{L}(\tilde{\zeta}_{i}),\nu_{i}^{-1})\in {\cal T}^{L}(\tilde{\zeta}_{i})\times {\cal T}(\tilde{\zeta}_{i})$ (la dissym\'etrie de cette notation et de la suivante s'expliquera plus loin);

- un cocycle $(\hat{V}^{L}(\tilde{\zeta}_{i}),\hat{V}_{i}, \hat{V}_{i,sc}):W_{F_{v}}\to \hat{S}_{i}$;

- un \'el\'ement  $(\zeta_{i,sc},\zeta_{i,sc})\in \hat{T}_{sc}^{L}\times \hat{T}_{sc}^{i}$;  pour cela, on  note  $\zeta$ l'\'el\'ement de $\hat{T}$ tel que $\tilde{\zeta}=\zeta\hat{\theta}$ et on  choisit des rel\`evements $\zeta_{sc}$ et $z_{sc}$ dans $\hat{T}_{sc}$ des images de $\zeta$ et $z$ dans $\hat{T}_{ad}$; on pose $\zeta_{1,sc}=\zeta_{sc}$ et $\zeta_{2,sc}=z_{sc}\zeta_{sc}$.

On note encore $(V_{i}^{L}, V_{i}^{-1})$, resp. $(\nu^{L}(\tilde{\zeta}_{i}),\nu_{i}^{-1})$, $(\zeta_{i,sc},\zeta_{i,sc})$,   les images de ces termes   dans $U_{i}$, resp. $S_{i}$, $\hat{U}_{i}$. Pour d\'efinir les \'el\'ements $(V_{i}^{L}, V_{i}^{-1})$ et $(\nu^{L}(\tilde{\zeta}_{i}),\nu_{i}^{-1})$, on doit compl\'eter $(B^{L,H}_{v},T^{L,H}_{v})$ en une paire de Borel \'epingl\'ee ${\cal E}^{H}$ et choisir un \'el\'ement $e^{H}\in Z(\tilde{H},{\cal E}^{H})$ ainsi qu'une cocha\^{\i}ne $u_{{\cal E}^{H}}:\Gamma_{F_{v}}\to G_{SC}(\bar{F}_{v})$ comme en [I] 1.2. Pour cela, on fixe une fois pour toutes une paire de Borel \'epingl\'ee ${\cal E}^*$ de $G$ d\'efinie sur $\bar{F}$, un \'el\'ement $e^*\in Z(\tilde{G},{\cal E}^*)$ et une cochaine $u_{{\cal E}^*}:\Gamma_{F}\to G_{SC}(\bar{F})$. On fixe un \'el\'ement $g_{sc,v}\in G_{SC}(\bar{F}_{v})$ tel que $ad_{g_{sc,v}}(B^*,T^*)=(B^{L,H}_{v}, T^{L,H}_{v})$. L'application $ad_{g_{sc,v}}$ envoie l'\'epinglage de ${\cal E}^*$ sur un \'epinglage qui compl\`ete la paire $(B^{L,H}_{v},T^{L,H}_{v})$. On suppose que ${\cal E}^{H}$ est cette paire de Borel compl\'et\'ee par cet \'epinglage. On pose $e^H=ad_{g_{sc,v}}(e^*)$ et $u_{{\cal E}^{H}}(\sigma)=g_{sc,v}u_{{\cal E}^*}(\sigma)\sigma(g_{sc,v})^{-1}$ pour tout $\sigma\in \Gamma_{F_{v}}$. On notera simplement $e=e^H$ dans la suite. 

Par d\'efinition, on a l'\'egalit\'e
$$\boldsymbol{\Delta}_{i,imp,v}(l'_{i},l';h'_{i},h(\tilde{\zeta}_{i}))=<\left((V_{i}^{L}, V_{i}^{-1}),(\nu^{L}(\tilde{\zeta}_{i}),\nu_{i}^{-1})\right),\left( (\hat{V}^{L}(\tilde{\zeta}_{i}),\hat{V}_{i}, \hat{V}_{i,sc}), (\zeta_{i,sc},\zeta_{i,sc})\right)>^{-1},$$
o\`u il s'agit du produit sur
$$H^{1,0}(\Gamma_{F_{v}};U_{i}\stackrel{1-\theta}{\to}S_{i})\times H^{1,0}(W_{F_{v}}; \hat{S}_{i}\stackrel{1-\hat{\theta}}{\to}\hat{U}_{i})\to {\mathbb C}^{\times}.$$

 On a des inclusions $G'(\tilde{\zeta}_{i})\subset H'(\tilde{\zeta}_{i})$. On note $G'_{i}(\tilde{\zeta}_{i})$ l'image r\'eciproque de $G'(\tilde{\zeta}_{i})$ dans la donn\'ee auxiliaire $H'_{i}(\tilde{\zeta}_{i})$ et $T^{L'}_{i}$ l'image r\'eciproque de $T^{L'}$ dans $G'_{i}(\tilde{\zeta}_{i})$. 
Notons ${\cal T}^{L}_{i}$ le produit fibr\'e de $T^{L'}_{i}$ et de $T^{L}$ au-dessus de $T^{L'}$.   La diff\'erence avec ${\cal T}^{L}(\tilde{\zeta}_{i})$ est qu'ici, les tores restent dans des groupes issus de $G$ et non pas de $H$. Posons  $\underline{S}_{i}=({\cal T}^{L}_{i}\times {\cal T}(\tilde{\zeta}_{i}))/diag_{-}({\cal Z}_{i})$, o\`u ${\cal Z}_{i}$ est le produit fibr\'e de $Z(G'_{i}(\tilde{\zeta}_{i}))$ et de $Z(G)$ au dessus de $Z(G'(\tilde{\zeta}_{i}))$. On a une application naturelle $\underline{S}_{i}\to S_{i}$.  Parce que les \'el\'ements $l'_{i}$ et $l$ sont dans $\tilde{G}'_{i}(\tilde{\zeta}_{i})$ et $\tilde{G}$ (et pas seulement dans $\tilde{H}'_{i}(\tilde{\zeta}_{i})$ et $\tilde{H}$) et parce que l'on a choisi $e\in \tilde{G}$, on  v\'erifie que le couple $(\nu^{L}(\tilde{\zeta}_{i}),\nu_{i}^{-1})$ appartient \`a ${\cal T}^{L}_{i}\times {\cal T}^{i}$ et d\'efinit donc un \'el\'ement de $\underline{S}_{i}$, que l'on note plut\^ot $(\nu^L_{i},\nu_{i}^{-1})$. Toujours d'apr\`es le choix de $e$, on v\'erifie que le couple $\left((V_{i}^{L}, V_{i}^{-1}),(\nu_{i}^{L},\nu_{i}^{-1})\right)$ est un cocycle appartenant \`a $Z^{1,0}(\Gamma_{F_{v}};U_{i}\stackrel{1-\theta}{\to}\underline{S}_{i})$. Le tore dual $\hat{{\cal T}}^L_{i}$ est le quotient de $\hat{T}^{L'}_{i}\times \hat{T}^L$ par $\hat{T}^{L'}$ plong\'e par $t'\mapsto(\hat{\xi}_{i}(t')^{-1},t')$. Le tore dual $\hat{\underline{S}}_{i}$ est le groupe des $(t^L,t,t_{sc})\in \hat{{\cal T}}^{L}_{i}\times \hat{{\cal T}}(\tilde{\zeta}_{i})\times \hat{T}_{sc}$ tels que $t^L$ est l'image naturelle dans $\hat{{\cal T}}^L_{i}$ de l'\'el\'ement $j(t_{sc})t$. On note $\hat{V}_{i}^L$ l'image naturelle de $\hat{V}^{L}(\tilde{\zeta}_{i})$ dans $\hat{\cal T}^L_{i}$. Par compatibilit\'e des produits, on obtient l'\'egalit\'e
$$\boldsymbol{\Delta}_{i,imp,v}(l'_{i},l';h'_{i},h(\tilde{\zeta}_{i}))=<\left((V_{i}^{L}, V_{i}^{-1}),(\nu_{i}^{L},\nu_{i}^{-1})\right),\left( (\hat{V}_{i}^{L},\hat{V}_{i}, \hat{V}_{i,sc}), (\zeta_{i,sc},\zeta_{i,sc})\right)>^{-1},$$
o\`u il s'agit du produit sur
$$H^{1,0}(\Gamma_{F_{v}};U_{i}\stackrel{1-\theta}{\to}\underline{S}_{i})\times H^{1,0}(W_{F_{v}}; \hat{\underline{S}}_{i}\stackrel{1-\hat{\theta}}{\to}\hat{U}_{i}).$$
Notons $U_{12}$ le produit $T_{sc}^L\times T^1_{sc}\times T^2_{sc}$ quotient\'e par $Z(G_{SC})$ plong\'e par $z\mapsto (z,z^{-1},z^{-1})$.  Parce que l'on a choisi les m\^emes diagrammes et les m\^emes objets auxiliaires pour construire  $V_{1}^L$ et $V_{2}^L$, ces deux cocha\^{\i}nes sont \'egales. Notons-les simplement $V^L$. Alors $(V^L,V_{1}^{-1},V_{2}^{-1})$ est une cocha\^{\i}ne \`a valeurs dans $T_{sc}^L\times T^1_{sc}\times T^2_{sc}$, qui se descend en un cocycle \`a valeurs dans $U_{12}$. Rappelons la construction des \'el\'ements $\nu_{i}^L$. On note $e'_{i}$ l'image de $e$ dans ${\cal Z}(\tilde{G}'(\tilde{\zeta}_{i}))$. On \'ecrit $l'_{i}=\mu_{i}e'_{i}$ et $l=\nu e$. Alors $\nu_{i}^L$ est le couple $(\mu_{i},\nu)$.  Notons ${\cal T}^L_{12}$ le produit fibr\'e de $T^{L'}_{1}$, $T^{L'}_{2}$ et $T^L$ au-dessus de $T^{L'}$. On d\'efinit l'\'el\'ement $\nu_{12}^L=(\mu_{1},\mu_{2},\nu)\in {\cal T}^L_{12}$. Notons ${\cal Z}_{12}$ le groupe des triplets $(z_{1},z_{2},z)\in Z(G'_{1}(\tilde{\zeta}_{1}))\times Z(G'_{2}(\tilde{\zeta}_{2}))\times Z(G)$ tels que, pour $i=1,2$, $z_{i}$ a m\^eme image que $z$ dans $Z(G'(\tilde{\zeta}_{i}))$.   Notons $\underline{S}_{12}$ le quotient de ${\cal T}_{12}^L\times {\cal T}(\tilde{\zeta}_{1})\times {\cal T}(\tilde{\zeta}_{2})$ par le groupe ${\cal Z}_{12}$ plong\'e par $(z_{1},z_{2},z)\mapsto ((z_{1},z_{2},z),(z_{1},z)^{-1},(z_{2},z)^{-1})$. Le triplet $(\nu_{12}^L,\nu_{1}^{-1},\nu_{2}^{-1})$ d\'efinit un \'el\'ement de ce quotient et on voit que la paire $\left((V^L,V_{1}^{-1},V_{2}^{-1}),(\nu_{12}^L,\nu_{1}^{-1},\nu_{2}^{-1})\right)$ est un cocycle appartenant \`a $Z^{1,0}(\Gamma_{F_{v}};U_{12}\stackrel{1-\theta}{\to}\underline{S}_{12})$. Il y a des homomorphismes d'oubli d'une s\'erie de variables
$$\begin{array}{ccccc}&&U_{12}\stackrel{1-\theta}{\to}\underline{S}_{12}&&\\ &\,\,\swarrow p_{1}&&\,\,\searrow p_{2}&\\ U_{1}\stackrel{1-\theta}{\to}\underline{S}_{1}&&&&U_{2}\stackrel{1-\theta}{\to}\underline{S}_{2}\\ \end{array}$$
En notant encore $p_{i}:H^{1,0}(\Gamma_{F_{v}};U_{12}\stackrel{1-\theta}{\to}\underline{S}_{12})\to H^{1,0}(\Gamma_{F_{v}};U_{i}\stackrel{1-\theta}{\to}\underline{S}_{i})$ l'homomorphisme d\'eduit par fonctorialit\'e, on a $p_{i}((V^L,V_{1}^{-1},V_{2}^{-1}),(\nu_{12}^L,\nu_{1}^{-1},\nu_{2}^{-1}))=((V^{L}_{i},V_{i}^{-1}),(\nu_{i}^L,\nu_{i}^{-1})$.  On note $\hat{p}_{i}:H^{1,0}(W_{F_{v}}; \hat{\underline{S}}_{i}\stackrel{1-\hat{\theta}}{\to}\hat{U}_{i})\to H^{1,0}(W_{F_{v}}; \hat{\underline{S}}_{12}\stackrel{1-\hat{\theta}}{\to}\hat{U}_{12})$ l'homomorphisme dual de $p_{i}$. Par compatibilit\'e des produits, on en d\'eduit
que $\boldsymbol{\Delta}_{i,imp,v}(l'_{i},l';h'_{i},h(\tilde{\zeta}_{i}))$ est \'egal \`a 
$$<\left((V^L,V_{1}^{-1},V_{2}^{-1}),(\nu_{12}^L,\nu_{1}^{-1},\nu_{2}^{-1})\right),\hat{p}_{i}\left((\hat{V}_{i}^{L},\hat{V}_{i}, \hat{V}_{i,sc}), (\zeta_{i,sc},\zeta_{i,sc})\right)>^{-1},$$
o\`u il s'agit du produit sur
$$H^{1,0}(\Gamma_{F_{v}};U_{12}\stackrel{1-\theta}{\to}\underline{S}_{12})\times H^{1,0}(W_{F_{v}}; \hat{\underline{S}}_{12}\stackrel{1-\hat{\theta}}{\to}\hat{U}_{12}).$$
Le tore dual $\hat{\underline{S}}_{12}$ est le groupe des $(t^L,t_{1},t_{2},t_{sc})\in \hat{{\cal T}}^L_{12}\times \hat{{\cal T}}(\tilde{\zeta}_{1})\times \hat{{\cal T}}(\tilde{\zeta}_{2})\times \hat{T}_{sc}$ tels que $t^L$ soit le produit des images naturelles de $t_{sc}$, $t_{1}$ et $t_{2}$ dans $\hat{{\cal T}}_{12}^L$. Le tore dual $\hat{U}_{12}$ est le quotient de $\hat{T}_{sc}^L\times \hat{T}_{sc}(\tilde{\zeta}_{1})\times \hat{T}_{sc}(\tilde{\zeta}_{2})$ par le groupe des triplets $(z,z_{1},z_{2})\in Z(\hat{G}_{SC})^3$ tels que $z=z_{1}z_{2}$. L'\'el\'ement 
$$\hat{p}_{1}\left((\hat{V}_{1}^{L},\hat{V}_{1}, \hat{V}_{1,sc}), (\zeta_{1,sc},\zeta_{1,sc})\right)\hat{p}_{2}\left((\hat{V}_{2}^{L},\hat{V}_{2}, \hat{V}_{2,sc}), (\zeta_{2,sc},\zeta_{2,sc})\right)^{-1}$$
est de la forme
$$\left((\hat{V}_{12}^L,\hat{V}_{1},\hat{V}_{2}^{-1},\hat{V}_{12,sc}),(z_{sc}^{-1},\zeta_{sc},z_{sc}^{-1}\zeta_{sc}^{-1})\right).$$
On obtient
$$\boldsymbol{\Delta}_{2,imp,v}(l'_{2},l';h'_{2},h(\tilde{\zeta}_{2}))\boldsymbol{\Delta}_{1,imp,v}(l'_{1},l';h'_{1},h(\tilde{\zeta}_{1}))^{-1}=$$
$$<\left((V^L,V_{1}^{-1},V_{2}^{-1}),(\nu_{12}^L,\nu_{1}^{-1},\nu_{2}^{-1})\right),\left((\hat{V}_{12}^L,\hat{V}_{1},\hat{V}_{2}^{-1},\hat{V}_{12,sc}),(z_{sc}^{-1},\zeta_{sc},z_{sc}^{-1}\zeta_{sc}^{-1})\right)>.$$
Le calcul de $\hat{V}_{12}^L$ et $\hat{V}_{12,sc}$ est long mais on l'a fait dans la preuve de la proposition 1.14(iii) de [II] et on va d\'ecrire le r\'esultat. On a implicitement fix\'e des paires de Borel \'epingl\'ees de $\hat{G}'(\tilde{\zeta}_{1})$ et $\hat{G}'(\tilde{\zeta}_{2})$, les paires de Borel sous-jacentes \'etant \'evidemment les intersections de $(\hat{B},\hat{T})$ avec chacun des groupes. Puisque ces deux groupes ont en commun le Levi $\hat{M}'$, on peut supposer que ces paires prolongent en un sens plus ou moins clair une paire de Borel \'epingl\'ee de de Levi. Pour $i=1,2$, le groupe $\hat{H}(\tilde{\zeta}_{i})$ est en r\'ealit\'e d\'etermin\'e par le choix d'un rel\`evement $\tilde{t}_{i}$ de $\tilde{\zeta}_{i}$ dans $\hat{H}$. Notons $\hat{M}^H$ l'image r\'eciproque de $\hat{M}$ dans $\hat{H}$. Puisque $z\in Z(\hat{M})^{\Gamma_{F},\hat{\theta},0}$, on peut le relever en un \'el\'ement $z^H\in Z(\hat{M}^H)^{\Gamma_{F},\hat{\theta},0}$ et supposer $\tilde{t}_{2}=z^H\tilde{t}_{1}$. Alors $\hat{H}'(\tilde{\zeta}_{1})$ et $\hat{H}'(\tilde{\zeta}_{2})$ ont en commun un Levi $\hat{M}^{H'}$ qui rel\`eve $\hat{M}'$. Les paires de Borel \'epingl\'ees de $\hat{G}'(\tilde{\zeta}_{1})$, $\hat{G}'(\tilde{\zeta}_{2})$ et $\hat{M}'$ se rel\`event en de telles paires de $\hat{H}'(\tilde{\zeta}_{1})$, $\hat{H}'(\tilde{\zeta}_{2})$ et $\hat{M}^{H'}$. Ces paires de Borel \'epingl\'ees d\'eterminent pr\'ecis\'ement les actions galoisiennes sur chaque groupe. Ainsi, les actions galoisiennes sur $\hat{H}'(\tilde{\zeta}_{1})$ et $\hat{H}'(\tilde{\zeta}_{2})$ se restreignent en une m\^eme action sur $\hat{M}^{H'}$. Pour $w\in W_{F}$ et $i=1,2$, on fixe des \'el\'ements $(h_{i}(w),w)\in {\cal H}'(\tilde{\zeta}_{i})$ tels que $ad_{h_{i}(w)}\circ w_{H}=w_{H'(\tilde{\zeta}_{i})}$.  La propri\'et\'e pr\'ec\'edente assure que $h_{2}(w)=m^{H'}(w)h_{1}(w)$, avec $m^{H'}(w)\in Z(\hat{M}^{H'})$. Puisque ce groupe est produit de $Z(\hat{M}^{H'})^0$ et de $Z(\hat{H}'(\tilde{\zeta}_{2}))$, on peut supposer $m^{H'}(w)\in Z(\hat{M}^{H'})^0$. On note $g_{1}(w)$, $g_{2}(w)$ et $m'(w)$ les projections de $h_{1}(w)$, $h_{2}(w)$ et $m^{H'}(w)$ dans $\hat{G}$. Ces \'el\'ements v\'erifient des propri\'et\'es analogues aux pr\'ec\'edents. Pour $i=1,2$, on note $\hat{\xi}_{i}:{\cal G}'(\tilde{\zeta}_{i})\to {^LG}'_{i}(\tilde{\zeta}_{i})$ le plongement fix\'e et, pour $w\in W_{F}$, on pose $\hat{\xi}_{i}(g_{i}(w),w)=(\varphi_{i}(w),w)$. L'\'el\'ement $\varphi_{i}(w)$ appartient \`a $Z(\hat{G}'_{i}(\tilde{\zeta}_{i}))$. On pose $\varphi'_{2}(w)= \hat{\xi}_{2}(m'(w))^{-1}\varphi_{2}(w)$, autrement dit $\hat{\xi}_{2}(g_{1}(w),w)=(\varphi'_{2}(w),w)$. On fixe un rel\`evement $\rho_{v,sc}\in \hat{G}_{SC}$ de l'image de $\rho_{v}$ dans $\hat{G}_{AD}$ et, pour tout $w\in W_{F}$, un rel\`evement $m'_{sc}(w)$ dans $\hat{G}_{SC}$ de l'image de $m'(w)$ dans $\hat{G}_{AD}$. On suppose ainsi qu'il est loisible que $m'_{sc}(w)\in Z(\hat{M}'_{sc})^0$. Pour $w\in W_{F_{v}}$, on a alors
$$\hat{V}_{12}^L(w)=(\varphi_{1}(w),\varphi'_{2}(w)^{-1},w_{G}(\rho_{v})w_{T^L}(\rho_{v})^{-1}),$$
$$\hat{V}_{12,sc}(w)=(   t_{sc}(\tilde{\zeta}_{1})(w)^{-1}t_{sc}(\tilde{\zeta}_{2})(w)m'_{sc}(w)w_{G}(\rho_{v,sc})w_{T^L}(\rho_{v,sc})^{-1}),$$
o\`u, pour $i=1,2$,  $t_{sc}(\tilde{\zeta}_{i})$ est une cocha\^{\i}ne ne d\'ependant que des objets issus de $(h'_{i},h(\tilde{\zeta}_{i}))$. 

{\bf Remarque.} Il y a un changement de signe par rapport \`a la r\'ef\'erence [II] car on y calculait un rapport $\boldsymbol{\Delta}_{1}/\boldsymbol{\Delta}_{2}$ alors que l'on calcule ici un rapport inverse. Par ailleurs, le r\'esultat n'est pas tout-\`a-fait exact. Il faudrait multiplier les formules ci-dessus par un cobord qui dispara\^{\i}t imm\'ediatement dans la suite du calcul.

\bigskip

On a fix\'e plus haut des \'el\'ements $m',m'_{1},m'_{2}$. On peut les supposer assez r\'eguliers. Notons $T^{M'},T^{M'}_{1}, T^{M'}_{2}$ leurs commutants et notons $T_{12}^{M'}$ le produit fibr\'e de $T^{M'}_{1}$ et $T^{M'}_{2}$ au-dessus de $T^{M'}$. On \'ecrit $m'_{i}=\mu^{M'}_{i}e'_{i}$ pour $i=1,2$. Alors le couple $\mu^{M'}_{12}=(\mu^{M'}_{1},\mu^{M'}_{2})$ appartient \`a $T^{M'}_{12}$. 
Notons $\Sigma$ le quotient de ${\cal T}^L_{12}\times T^{M'}_{12}\times {\cal T}(\tilde{\zeta}_{1})\times {\cal T}(\tilde{\zeta}_{2})$ par le groupe ${\cal Z}_{12}$ plong\'e par $(z_{1},z_{2},z)\mapsto ((z_{1},z_{2},z),(z_{1},z_{2})^{-1},(z_{1},z)^{-1},(z_{2},z)^{-1})$. Le quadruplet $(\nu_{12}^L,(\mu^{M'}_{12})^{-1},\nu_{1}^{-1},\nu_{2}^{-1})$ d\'efinit un \'el\'ement de $\Sigma$. On a un homomorphisme d'oubli
$$\begin{array}{c}U_{12}\stackrel{1-\theta}{\to}\Sigma\\ \downarrow\\ U_{12}\stackrel{1-\theta}{\to}\underline{S}_{12}\\ \end{array}$$
Il est clair que 
$$\left((V^L,V_{1}^{-1},V_{2}^{-1}),(\nu_{12}^L,\nu_{1}^{-1},\nu_{2}^{-1})\right)$$
est l'image de 
$$\left((V^L,V_{1}^{-1},V_{2}^{-1}),(\nu_{12}^L,(\mu^{M'}_{12})^{-1}\nu_{1}^{-1},\nu_{2}^{-1})\right)$$
par l'homorphisme
$$H^{1,0}(\Gamma_{F_{v}}; U_{12}\stackrel{1-\theta}{\to}\Sigma)\to H^{1,0}(\Gamma_{F_{v}};U_{12}\stackrel{1-\theta}{\to}\underline{S}_{12})$$
d\'eduit par fonctorialit\'e du pr\'ec\'edent. Le tore dual $\hat{\Sigma}$ est le groupe des $(t^L,t^{M'},t_{1},t_{2},t_{sc})\in \hat{{\cal T}}^L_{12}\times \hat{T}^{M'}_{12}\times \hat{{\cal T}}(\tilde{\zeta}_{1})\times \hat{{\cal T}}(\tilde{\zeta}_{2})\times \hat{T}_{sc}$ tels que $t^L$ soit le produit des images naturelles de $t^{M'}$, $t_{1}$, $t_{2}$ et $t_{sc}$. Par l'homomorphisme dual du pr\'ec\'edent, l'\'el\'ement
$$\left((\hat{V}_{12}^L,\hat{V}_{1},\hat{V}_{2}^{-1},\hat{V}_{12,sc}),(z_{sc}^{-1},\zeta_{sc},z_{sc}^{-1}\zeta_{sc}^{-1})\right)$$
s'envoie sur
$$(18) \qquad \left((\hat{V}_{12}^L,1, \hat{V}_{1},\hat{V}_{2}^{-1},\hat{V}_{12,sc}),(z_{sc}^{-1},\zeta_{sc},z_{sc}^{-1}\zeta_{sc}^{-1})\right).$$
Par compatibilit\'e des produits, on obtient
$$\boldsymbol{\Delta}_{2,imp,v}(l'_{2},l';h'_{2},h(\tilde{\zeta}_{2}))\boldsymbol{\Delta}_{1,imp,v}(l'_{1},l';h'_{1},h(\tilde{\zeta}_{1}))^{-1}=$$
$$<\left((V^L,V_{1}^{-1},V_{2}^{-1}),(\nu_{12}^L,(\mu^{M'})^{-1},\nu_{1}^{-1},\nu_{2}^{-1})\right),\left((\hat{V}_{12}^L,1,\hat{V}_{1},\hat{V}_{2}^{-1},\hat{V}_{12,sc}),(z_{sc}^{-1},\zeta_{sc},z_{sc}^{-1}\zeta_{sc}^{-1})\right)>,$$
o\`u il s'agit du produit sur
$$H^{1,0}(\Gamma_{F_{v}};U_{12}\stackrel{1-\theta}{\to}\Sigma)\times H^{1,0}(W_{F_{v}};\hat{\Sigma}\stackrel{1-\hat{\theta}}{\to}\hat{U}_{12}).$$
Notons $\Sigma_{ML}$ le quotient de ${\cal T}_{12}^L\times T^{M'}_{12}$ par ${\cal Z}_{12}$ plong\'e par $(z_{1},z_{2},z)\mapsto ((z_{1},z_{2},z),(z_{1},z_{2})^{-1})$. Son dual $\hat{\Sigma}_{ML}$ est le groupe des $(t^L,t^{M'},t_{sc})$ tels que $t^L(t^{M'})^{-1}=j(t_{sc})$. Pour $w\in W_{F_{v}}$, posons
$$X_{ML}^L(w)=\hat{V}_{12}^L(w)\in \hat{{\cal T}}_{12}^L,\,\,X_{ML}^M(w)=(\varphi_{1}(w),\varphi'_{2}(w)^{-1})\in \hat{T}^{M'}_{12}$$
$$X_{ML,sc}(w)=( w_{G}(\rho_{v,sc})w_{T^L}(\rho_{v,sc})^{-1})\in \hat{T}_{sc}.$$
Posons $X_{ML}(w)=(X_{ML}^L(w),X_{ML}^M(w),X_{ML,sc}(w))$. Ce terme appartient \`a $\hat{\Sigma}_{ML}$. On a fix\'e arbitrairement l'\'el\'ement $z_{sc}$. Mais on se rappelle que, par d\'efinition, $z$ appartient \`a $Z(\hat{M})^{\Gamma_{F},\hat{\theta}}$. L'image de ce groupe dans $\hat{G}_{AD}$ \'etant connexe, on peut supposer et on suppose que $z_{sc}\in Z(\hat{M}_{sc})^{\Gamma_{F},\hat{\theta},0}$. On v\'erifie alors que le couple $(X_{ML},z_{sc}^{-1})$ est un cocycle qui d\'efinit un \'el\'ement de $H^{1,0}(W_{F_{v}};\hat{\Sigma}_{ML}\stackrel{1-\hat{\theta}}{\to}\hat{T}^L_{sc})$. On a un homomorphisme plus ou moins \'evident
$$\begin{array}{c}\hat{\Sigma}_{ML}\stackrel{1-\hat{\theta}}{\to}\hat{T}^L_{sc}\\ \downarrow\\ \hat{\Sigma}\stackrel{1-\hat{\theta}}{\to}\hat{U}_{12}\\ \end{array}$$
Par l'homomorphisme qui s'en d\'eduit par fonctorialit\'e, $(X_{ML},z_{sc}^{-1})$ s'envoie sur le cocycle
$$\left((\hat{V}_{12}^L,X_{ML}^M,1,1,X_{ML,sc}),(z_{sc}^{-1},1,1)\right).$$
Le cocycle (18) est le produit de celui-ci avec le cocycle
$$(19) \qquad \left((1,(X_{ML}^M)^{-1},\hat{V}_{1},\hat{V}_{2}^{-1},\hat{\underline{V}}_{12,sc}),(1,\zeta_{sc},z_{sc}^{-1}\zeta_{sc}^{-1})\right),$$
o\`u
$$\hat{\underline{V}}_{12,sc}(w)=( t_{sc}(\tilde{\zeta}_{1})(w)^{-1}t_{sc}(\tilde{\zeta}_{2})(w)m'_{sc}(w)).$$
Il est clair que ce dernier cocycle vit dans des groupes plus petits, o\`u l'on supprime la premi\`ere composante. C'est-\`a-dire, notons $\Sigma_{\star}$ le quotient de $T^{M'}_{12}\times {\cal T}(\tilde{\zeta}_{1})\times {\cal T}(\tilde{\zeta}_{2})$ par le groupe ${\cal Z}_{12}$ plong\'e par
$(z_{1},z_{2},z)\mapsto ((z_{1},z_{2}),(z_{1},z),(z_{2},z))$. Notons $U_{\star}=(T_{sc}(\tilde{\zeta}_{1})\times T_{sc}(\tilde{\zeta}_{2}))/diag(Z(G_{SC}))$.  On a un homomorphisme
$$\begin{array}{c}\hat{\Sigma}_{\star}\stackrel{1-\hat{\theta}}{\to}\hat{U}_{\star}\\ \downarrow \\ \hat{\Sigma}\stackrel{1-\hat{\theta}}{\to}\hat{U}_{12}\\ \end{array}$$
Par l'homomorphisme qui s'en d\'eduit par fonctorialit\'e, le cocycle (19) est l'image du cocycle
$$ \left(((X_{ML}^M)^{-1},\hat{V}_{1},\hat{V}_{2}^{-1},\hat{\underline{V}}_{12,sc}),(\zeta_{sc},z_{sc}^{-1}\zeta_{sc}^{-1})\right)\in H^{1,0}(W_{F_{v}};\hat{\Sigma}_{\star}\stackrel{1-\hat{\theta}}{\to}\hat{U}_{\star}).$$
On a des homomorphismes duaux aux pr\'ec\'edents
$$\begin{array}{ccc}&&H^{1,0}(\Gamma_{F_{v}};T^L_{ad}\stackrel{1-\theta}{\to}\Sigma_{ML})\\ &\nearrow&\\ H^{1,0}(\Gamma_{F_{v}};U_{12}\stackrel{1-\theta}{\to }\Sigma)&&\\ &\searrow&\\ &&H^{1,0}(\Gamma_{F_{v}}; U_{\star}\stackrel{1-\theta}{\to}\Sigma_{\star})\\ \end{array}$$
Par ces homomorphismes, le cocycle 
$$\left((V^L,V_{1}^{-1},V_{2}^{-1}),(\nu_{12}^L,(\mu^{M'})^{-1},\nu_{1}^{-1},\nu_{2}^{-1})\right)$$
s'envoie respectivement sur $(V_{L,ad},(\nu_{12}^L,(\mu^{M'})^{-1}))$ et l'inverse de $((V_{1},V_{2}), (\mu^{M'},\nu_{1},\nu_{2}))$. La d\'ecomposition ci-dessus et la compatibilit\'e des produits conduit \`a l'\'egalit\'e
$$(20) \qquad \boldsymbol{\Delta}_{2,imp,v}(l'_{2},l';h'_{2},h(\tilde{\zeta}_{2}))\boldsymbol{\Delta}_{1,imp,v}(l'_{1},l';h'_{1},h(\tilde{\zeta}_{1}))^{-1}=A_{v}B_{v}^{-1},$$
o\`u
$$A_{v}=<(V_{L,ad},(\nu_{12}^L,(\mu^{M'})^{-1})), (X_{ML},z_{sc}^{-1})>,$$
$$B_{v}=<((V_{1},V_{2}), (\mu^{M'},\nu_{1},\nu_{2})), \left(((X_{ML}^M)^{-1},\hat{V}_{1},\hat{V}_{2}^{-1},\hat{\underline{V}}_{12,sc}),(\zeta_{sc},z_{sc}^{-1}\zeta_{sc}^{-1})\right)>.$$
On va d'abord se pr\'eoccuper du terme $B_{v}$. Evidemment, ce qui nous int\'eresse est le produit de ces termes sur toutes les places $v\in V$. Compte tenu du choix des \'el\'ements $m'_{i}$, $h'_{i}$ et $h(\tilde{\zeta}_{i})$, le terme $B_{v}$ peut aussi bien \^etre d\'efini pour une place $v\not\in V$. De plus, bien que les \'el\'ements $h'_{i}$ et $h(\tilde{\zeta}_{i})$ soient seulement ad\'eliques, les tores qui interviennent dans les d\'efinitions sont les localis\'es de tores d\'efinis sur $F$. Les cocha\^{\i}nes intervenant sont aussi "ad\'eliques". Par exemple, le terme not\'e $V_{1}$ est la localis\'ee d'une cocha\^{\i}ne encore not\'ee $V_{1}:\Gamma_{F}\to T_{sc}(\tilde{\zeta}_{1})$. On peut donc d\'efinir un terme
$$B=<((V_{1},V_{2}), (\mu^{M'},\nu_{1},\nu_{2})), \left(((X_{ML}^M)^{-1},\hat{V}_{1},\hat{V}_{2}^{-1},\hat{\underline{V}}_{12,sc}),(\zeta_{sc},z_{sc}^{-1}\zeta_{sc}^{-1})\right)>,$$
o\`u il s'agit cette fois du produit dans
$$H^{1,0}({\mathbb A}_{F}/F;U_{\star}\stackrel{1-\theta}{\to}\Sigma_{\star})\times H^{1,0}(W_{F}; \hat{\Sigma}_{\star}\stackrel{1-\hat{\theta}}{\to}\hat{U}_{\star}).$$
D'apr\`es les propri\'et\'es g\'en\'erales de ce produit, on a
$$B=\prod_{v\in Val(F)}B_{v},$$
les termes du produit \'etant presque tous \'egaux \`a $1$. On a un homomorphisme naturel
$$\begin{array}{c}T_{sc}(\tilde{\zeta}_{1})\times T_{sc}(\tilde{\zeta}_{2})\stackrel{1-\theta}{\to}{\cal T}(\tilde{\zeta}_{1})\times {\cal T}(\tilde{\zeta}_{2})\\ \downarrow\\ U_{\star}\stackrel{1-\theta}{\to}\Sigma_{\star}.\\ \end{array}$$
Le quadruplet $((V_{1},V_{2}),(\nu_{1},\nu_{2}))$ d\'efinit naturellement un \'el\'ement de 
$$H^{1,0}({\mathbb A}_{F}/F;T_{sc}(\tilde{\zeta}_{1})\times T_{sc}(\tilde{\zeta}_{2})\stackrel{1-\theta}{\to}{\cal T}(\tilde{\zeta}_{1})\times {\cal T}(\tilde{\zeta}_{2})).$$
Le cocycle $((V_{1},V_{2}), (\mu^{M'},\nu_{1},\nu_{2}))$ en est l'image dans $H^{1,0}({\mathbb A}_{F}/F;U_{\star}\stackrel{1-\theta}{\to}\Sigma_{\star})$ par l'homomorphisme d\'eduit par fonctorialit\'e du pr\'ec\'edent. En effet, parce que $\mu'_{1}$ et $\mu'_{2}$ sont d\'efinis sur $F$ et parce que les \'el\'ements $e'_{1}$ et $e'_{2}$ sont d\'efinis sur $\bar{F}$, le terme $\mu^{M'}$ appartient \`a $T^{M'}_{12}(\bar{F})$. Il dispara\^{\i}t par d\'efinition des groupes de cohomologie "globaux". Par compatibilit\'e des produits, on obtient
$$B=<((V_{1},V_{2}),(\nu_{1},\nu_{2})),((\hat{V}_{1},\hat{V}_{2}^{-1}),(\zeta_{ad},z_{ad}^{-1}\zeta_{ad}^{-1}))>,$$
o\`u il s'agit du produit sur 
$$H^{1,0}({\mathbb A}_{F}/F;T_{sc}(\tilde{\zeta}_{1})\times T_{sc}(\tilde{\zeta}_{2})\stackrel{1-\theta}{\to}{\cal T}(\tilde{\zeta}_{1})\times {\cal T}(\tilde{\zeta}_{2}))\times H^{1,0}(W_{F};\hat{{\cal T}}(\tilde{\zeta}_{1})\times \hat{{\cal T}}(\tilde{\zeta}_{2})\stackrel{1-\hat{\theta}}{\to}\hat{T}_{ad}(\tilde{\zeta}_{1})\times \hat{T}_{ad}(\tilde{\zeta}_{2})).$$
Ce produit se d\'ecompose selon les composantes index\'ees par $1$ et $2$. On n'a pas de mal \`a reconna\^{\i}tre ces composantes comme les facteurs $\Delta_{i,glob}(h'_{i},h(\tilde{\zeta}_{i}))$ priv\'es de leurs facteurs $\Delta_{II}$. Notons ces facteurs $\Delta_{i,imp,glob}(h'_{i},h(\tilde{\zeta}_{i}))$. Il y a toutefois une inversion de signe sur le facteur d'indice $1$ et on obtient
$$(21) \qquad B=\Delta_{2,imp,glob}(h'_{2},h(\tilde{\zeta}_{2}))\Delta_{1,imp,glob}(h'_{1},h(\tilde{\zeta}_{1}))^{-1}.$$
Pour $v\in V$, la relation (20) et le fait que le terme $A_{v}$ ne d\'epend pas des couples $(h'_{i},h(\tilde{\zeta}_{i}))$ entra\^{\i}ne que, si l'on remplace dans les constructions ces couples par d'autres $(\underline{h}'_{i},\underline{h}(\tilde{\zeta}_{i}))$, et si l'on note $\underline{B}_{v}$ le terme obtenu, on a l'\'egalit\'e
$$B_{v}=\underline{B}_{v}\boldsymbol{\Delta}_{1,imp,v}(\underline{h}'_{1},\underline{h}(\tilde{\zeta}_{1});h'_{1},h(\tilde{\zeta}_{1}))\boldsymbol{\Delta}_{2,imp,v}(\underline{h}'_{2},\underline{h}(\tilde{\zeta}_{2});h'_{2},h(\tilde{\zeta}_{2}))^{-1}.$$
On a

(22) cette propri\'et\'e perdure pour tout $v\in Val(F)$. 

 Dans le calcul conduisant \`a l'\'egalit\'e (20), on est parti d'un produit d\'ependant de trois donn\'ees $(l'_{1},l'_{2},l)$, $(h'_{1},h(\tilde{\zeta}_{1}))$ et $(h'_{2},h(\tilde{\zeta}_{2}))$. On a ins\'er\'e une quatri\`eme donn\'ee $(m'_{1},m'_{2})$, puis on a d\'ecompos\'e le produit obtenu en deux produits, l'un relatif aux donn\'ees $(l'_{1},l'_{2},l)$ et $(m'_{1},m'_{2})$, l'autre aux donn\'ees $(m'_{1},m'_{2})$,  $(h'_{1},h(\tilde{\zeta}_{1}))$ et $(h'_{2},h(\tilde{\zeta}_{2}))$. Le m\^eme proc\'ed\'e permet d'ins\'erer dans $B_{v}$ de nouvelles donn\'ees, disons $(\underline{h}'_{1},\underline{h}(\tilde{\zeta}_{1}))$ puis de d\'ecomposer le produit obtenu en deux produits, l'un relatif aux donn\'ees $(m'_{1},m'_{2})$, $(\underline{h}'_{1},\underline{h}(\tilde{\zeta}_{1}))$  et $(h'_{2},h(\tilde{\zeta}_{2}))$, l'autre relatif aux donn\'ees $(\underline{h}'_{1},\underline{h}(\tilde{\zeta}_{1}))$  et $(h'_{1},h(\tilde{\zeta}_{1}))$. On reconna\^{\i}t ces produits comme \'etant $\underline{B}_{v}$ et $\boldsymbol{\Delta}_{1,imp,v}(\underline{h}'_{1},\underline{h}(\tilde{\zeta}_{1});h'_{1},h(\tilde{\zeta}_{1}))$. On laisse les d\'etails au lecteur. Cela prouve (22).

Fixons une place $v\not\in V$. La situation \'etant non ramifi\'ee, ${\bf M}_{v}'$ est relevant. On peut fixer un \'el\'ement $y\in \tilde{M}_{v}(F_{v})$ assez r\'egulier et $y'\in M'(F_{v})$ de sorte que leurs classes de conjugaison stable se correspondent. On fixe des rel\`evements $y'_{1}$ et $y'_{2}$ de $y'$ dans $\tilde{G}'_{1}(\tilde{\zeta}_{1})$, resp. $ \tilde{G}'_{2}(\tilde{\zeta}_{2})$. On construit $\underline{B}_{v}$ comme ci-dessus, relatif aux couples $(\underline{h}'_{i},\underline{h}(\tilde{\zeta}_{i}))=(y'_{i},y)$ pour $i=1,2$. La situation \'etant non ramifi\'ee, on dispose de facteurs de transfert canoniques et la relation pr\'ec\'edant  (22) devient
$$(23) \qquad B_{v}=\underline{B}_{v}\Delta_{1,imp,v}(y'_{1},y)\Delta_{2,imp,v}(y'_{2},y)^{-1}\Delta_{1,imp,v}(h'_{1},h(\tilde{\zeta}_{1}))^{-1}\Delta_{2,imp,v}(h'_{2},h(\tilde{\zeta}_{2})).$$
On va calculer $\underline{B}_{v}$. Par rapport \`a la situation ant\'erieure, les tores $T(\tilde{\zeta}_{1})$ et $T(\tilde{\zeta}_{2})$ se confondent en un unique tore $T_{y}^H$. Les cocycles $V_{1}$ et $V_{2}$ sont \'egaux \`a un unique cocycle $V_{y}$. Les \'el\'ements $\nu_{1}$ et $\nu_{2}$ ne sont pas \'egaux, mais sont de la forme $(\mu_{y,1},\nu_{y})$, $(\mu_{y,2},\nu_{y})$, o\`u les $\mu_{y,i}$ appartiennent aux sous-tores $T^{'H}_{y,i}$ de $H'_{i}(\tilde{\zeta}_{i})$ associ\'es \`a $T_{y}^H$. En fait, on peut simplifier puisqu'on a choisi $y\in \tilde{M}_{v}(F_{v})$ et non pas seulement $y\in \tilde{M}_{v}^H(F_{v})$.  On note $T_{y}=G\cap T_{y}^H$ et $T'_{y,i}=T^{'H}_{y,i}\cap G'_{i}(\tilde{\zeta}_{i})$ pour $i=1,2$. On introduit le tore ${\cal T}_{y}$ produit fibr\'e de $T'_{y,1}$, $T'_{y,2}$ et $T_{y}$ au-dessus du tore $T'_{y}$ de $M'$ (qui est un Levi commun de $G'(\tilde{\zeta}_{1})$ et $G'(\tilde{\zeta}_{2})$). On note $\nu_{y,12}$ l'\'el\'ement $(\mu_{y,1},\mu_{y,2},\nu_{y})$ de ${\cal T}_{y}$. On note $\Sigma_{y}$ le quotient de $T^{M'}\times {\cal T}_{y}$ par ${\cal Z}_{12}$ plong\'e par $(z_{1},z_{2},z)\mapsto ((z_{1},z_{2}),(z_{1},z_{2},z))$. Le couple $(V_{y},(\mu^{M'},\nu_{y,12}))$ d\'efinit un \'el\'ement de $H^{1,0}(\Gamma_{F_{v}};T_{y,sc}\stackrel{1-\theta}{\to}\Sigma_{y})$. On a un homomorphisme \'evident
$$\begin{array}{c}T_{y,sc}\stackrel{1-\theta}{\to}\Sigma_{y}\\ \downarrow\\ U_{\star}\stackrel{1-\theta}{\to}\Sigma_{\star}\\ \end{array}$$
Par l'homomorphisme fonctoriellement associ\'e, le cocycle pr\'ec\'edent s'envoie sur celui intervenant dans la d\'efinition de $\underline{B}_{v}$. Pour utiliser la compatibilit\'e des produits,  
on doit calculer l'image de $(((X_{ML}^M)^{-1},\hat{V}_{1},\hat{V}_{2},\hat{\underline{V}}_{12,sc}), (\zeta_{sc},z_{sc}^{-1}\zeta_{sc}^{-1}))$ par l'homomorphisme dual
$$H^{1,0}(W_{F_{v}};\hat{\Sigma}_{\star}\stackrel{1-\hat{\theta}}{\to}\hat{U}_{\star})\to H^{1,0}(W_{F_{v}}; \hat{\Sigma}_{y}\stackrel{1-\hat{\theta}}{\to}\hat{T}_{y,ad}).$$
Le premier terme $X_{ML}^M$ se conserve tel quel. On peut remplacer $\hat{V}_{i}$ par son image dans l'analogue de $\hat{{\cal T}}(\tilde{\zeta}_{i})$ o\`u le groupe $\hat{G}$ remplace $\hat{H}$. Une fois fait ce remplacement, reportons-nous aux d\'efinitions de [I] 2.2. On a une \'egalit\'e $\hat{V}_{i}(w)=(\varphi_{i}(w),t(\tilde{\zeta}_{i})(w))$ pour tout $w\in W_{F_{v}}$. Parce que $y$ appartient \`a l'espace de Levi commun $\tilde{M}_{v}$, on voit que $t(\tilde{\zeta}_{i})$ est de la forme $t(\tilde{\zeta}_{i})(w)=N(w)g_{i}(w)^{-1}N'(w)$, o\`u $N(w)\in \hat{M}_{v,sc}$ et $N'(w)\in \hat{M}'_{sc}$ sont les m\^emes pour $i=1,2$. L'image de $(\hat{V}_{1},\hat{V}_{2})$ dans $\hat{{\cal T}}_{y}$ est le triplet $(\varphi_{1},\varphi_{2}^{-1},t(\tilde{\zeta}_{1})t(\tilde{\zeta}_{2})^{-1})$, que l'on calcule gr\^ace aux formules ci-dessus. C'est $(\varphi_{1},\varphi_{2}^{-1},m')$. Ce terme vit en fait dans un quotient par le groupe $(\hat{T}'_{y})^2$ plong\'e par $(t_{1},t_{2})\mapsto (\hat{\xi}_{1}(t_{1}),\hat{\xi}_{2}(t_{2}),t_{1}^{-1}t_{2}^{-1})$. Cela permet de remplacer le triplet pr\'ec\'edent par $(\varphi_{1},(\varphi_{2}')^{-1},1)$.
Le terme $\underline{V}_{12,sc}$ se simplifie. Avec les m\^emes notations que ci-dessus et en fixant un rel\`evement $g_{1,sc}(w)$ de l'image de $h_{1}(w)$ dans $\hat{G}_{AD}$, on a 
$$\underline{V}_{12,sc}(w)=t_{sc}(\tilde{\zeta}_{1})(w)^{-1}t_{sc}(\tilde{\zeta}_{2})(w)m'_{sc}(w)$$
$$= N'(w)^{-1}g_{1,sc}(w)N(w)^{-1}N(w)g_{1,sc}(w)^{-1}m'_{sc}(w)^{-1}N'(w)m'_{sc}(w)=1$$
(on se rappelle que $m'_{sc}(w)$ appartient \`a $Z(\hat{M}'_{sc})^0$). 
Enfin, le couple $(\zeta_{sc},z_{sc}^{-1}\zeta_{sc}^{-1})$ s'envoie \'evidemment sur $z_{ad}^{-1}$. On obtient
$$\underline{B}_{v}=<(V_{y},((\mu^{M'},\nu_{y,12})),((\varphi_{1}^{-1},\varphi'_{2}),(\varphi_{1},(\varphi'_{2})^{-1},1),1), z_{ad}^{-1})>.$$
D'apr\`es nos choix, on a $z_{ad}\in Z(\hat{M}_{ad})^{\Gamma_{F_{v}},0}$. Ce groupe est contenu dans $\hat{T}_{y,ad}^{\Gamma_{F_{v}},0}$ puisque $\hat{T}_{y}\subset \hat{M}$. Or le groupe $\hat{T}_{y,ad}^{\Gamma_{F_{v}},0}$ est le noyau de l'accouplement sur 
$$H^{1,0}(\Gamma_{F_{v}}; T_{y,sc}\stackrel{1-\theta}{\to}\Sigma_{y})\times H^{1,0}(W_{F_{v}}; \hat{\Sigma}_{y}\stackrel{1-\hat{\theta}}{\to}\hat{T}_{y,ad}).$$
On peut donc aussi bien supprimer le terme $z_{ad}$ de la formule ci-dessus.  Le couple $(\varphi_{1},(\varphi_{2}')^{-1})$, modulo le tore $\hat{T}'_{y}$ plong\'e par $t\mapsto (\hat{\xi}_{1}(t),\hat{\xi}_{2}(t)^{-1})$, ne d\'epend pas du choix de $g_{1}(w)$.  On peut modifier cette cocha\^{\i}ne $g_{1}$. La situation \'etant non ramifi\'ee, on peut supposer que c'est un cocycle non ramifi\'e. Il en est alors de m\^eme de $\varphi_{1}$ et $\varphi'_{2}$.  Pour $i=1,2$, notons $Z(G_{i}'(\tilde{\zeta}_{i});G)$ la projection dans $Z(G'_{i}(\tilde{\zeta}_{i}))$ du produit fibr\'e de ce groupe avec $Z(G)$, au-dessus de $Z(G'(\tilde{\zeta}_{i}))$.  Introduisons le tore $Y_{i}$ quotient de $T^{M'}_{i}\times T'_{y,i}$ par  $Z(G'_{i}(\tilde{\zeta}_{i});G)$ agissant diagonalement. Son dual $\hat{Y}_{i}$ est l'ensemble des $(t^{M'},t_{y},t_{sc})\in \hat{T}^{M'}\times \hat{T}'_{y,i}\times \hat{T}_{sc}^{\hat{\theta}}$  tels que $t^{M'}t_{y}t_{sc}=1$. On a un homomorphisme
$$\begin{array}{ccc}\hat{Y}_{i}&&\\ \downarrow&&\\ \hat{\Sigma}_{y}&\stackrel{1-\hat{\theta}}{\to}&\hat{T}_{y,ad}\\ \end{array}$$
dont on d\'eduit un homomorphisme
$$H^1(W_{F_{v}}; \hat{Y}_{i})\to H^{1,0}(W_{F_{v}}; \hat{\Sigma}_{y}\stackrel{1-\hat{\theta}}{\to}\hat{T}_{y,ad}).$$
Le triplet $(\varphi_{1}^{-1},\varphi_{1},1)$ d\'efinit un \'el\'ement de $H^1(W_{F_{v}};\hat{Y}_{1})$ tandis que le triplet $(\varphi'_{2},(\varphi'_{2})^{-1},1)$ d\'efinit un \'el\'ement de $H^1(W_{F_{v}};\hat{Y}_{2})$. Le cocycle $((\varphi_{1}^{-1},\varphi'_{2}),(\varphi_{1},(\varphi'_{2})^{-1},1),1),1)$ est le produit des images de ces deux cocycles. Par compatibilit\'e des produits ([KS2] 4.3), on obtient
$$\underline{B}_{v}=<(\mu'_{1},\mu_{y,1}),(\varphi_{1}^{-1},\varphi_{1},1)><(\mu'_{2},\mu_{y,2}),(\varphi'_{2},(\varphi'_{2})^{-1},1)>,$$
o\`u il s'agit des produits sur $H^0(\Gamma_{F_{v}}; Y_{i})\times H^1(W_{F_{v}};\hat{Y}_{i})$.  Posons 
$$^2M'_{i}=(M'_{i}(\tilde{\zeta}_{i})\times M'_{i}(\tilde{\zeta}_{i}))/diag(Z(G'_{i}(\tilde{\zeta}_{i});G)).$$Le tore $Y_{i}$ est un sous-tore maximal de ce groupe. On a un homomorphisme
$$\begin{array}{ccc}H^1(W_{F_{v}}; Z(\hat{M}'_{1}(\tilde{\zeta}_{1})))&\to  &H^1(W_{F_{v}}; Z(^2\hat{M}'_{1}))\\ \varphi_{1}&\mapsto &(\varphi_{1},\varphi_{1}^{-1},1)\\ \end{array}$$
et un homomorphisme analogue concernant $\varphi'_{2}$.
Notons $\lambda_{\varphi_{1}}$ le caract\`ere de $M'_{1}(\tilde{\zeta}_{1})$ d\'etermin\'e par $\varphi_{1}$. On d\'efinit de m\^eme $\lambda_{\varphi'_{2}}$. On calcule alors
$$\underline{B}_{v}=\lambda_{\varphi_{1}}((\mu'_{1})^{-1}\mu_{y,1})\lambda_{\varphi'_{2}}(\mu'_{2}\mu_{y,2}^{-1}).$$
On note aussi $\tilde{\lambda}_{\varphi_{1}}$ l'unique fonction sur $\tilde{M}'_{1}(\tilde{\zeta}_{1})$ qui se transforme selon le caract\`ere $\lambda_{\varphi_{1}}$ et qui vaut $1$ sur l'espace  hypersp\'ecial fix\'e par les donn\'ees auxiliaires. On peut aussi bien r\'ecrire
$$\underline{B}_{v}=\tilde{\lambda}_{\varphi_{1}}(m'_{1})^{-1}\tilde{\lambda}_{\varphi_{1}}(y'_{1})\tilde{\lambda}_{\varphi'_{2}}(m'_{2})\tilde{\lambda}_{\varphi'_{2}}(y'_{2})^{-1}.$$
Il r\'esulte des d\'efinitions que
$$\tilde{\lambda}_{\varphi_{1}}(y'_{1})\tilde{\lambda}_{\varphi'_{2}}(y'_{2})^{-1}=\Delta_{1,imp,v}(y'_{1},y)^{-1}\Delta_{2,imp,v}(y'_{2},y).$$
Cela assure en m\^eme temps que 
$$\tilde{\lambda}_{v}^M(m'_{1},m'_{2})=\tilde{\lambda}_{\varphi_{1}}(m'_{1})\tilde{\lambda}_{\varphi'_{2}}(m'_{2})^{-1}.$$
D'o\`u l'\'egalit\'e
$$\underline{B}_{v}=\tilde{\lambda}_{v}^M(m'_{1},m'_{2})^{-1}\Delta_{1,imp,v}(y'_{1},y)^{-1}\Delta_{2,imp,v}(y'_{2},y).$$
Gr\^ace \`a (23), on obtient
$$B_{v}=\tilde{\lambda}_{v}^M(m'_{1},m'_{2})^{-1}\Delta_{1,imp,v}(h'_{1},h(\tilde{\zeta}_{1}))^{-1}\Delta_{2,imp,v}(h'_{2},h(\tilde{\zeta}_{2})).$$
Rassemblons cette \'egalit\'e avec les \'egalit\'es (20) et (21). On obtient
$$\prod_{v\in V}\boldsymbol{\Delta}_{2,imp,v}(l'_{2},l';h'_{2},h(\tilde{\zeta}_{2}))\boldsymbol{\Delta}_{1,imp,v}(l'_{1},l';h'_{1},h(\tilde{\zeta}_{1}))^{-1}$$
$$=(\prod_{v\in V}A_{v})B^{-1}\prod_{v\not\in V}B_{v}$$
$$=(\prod_{v\in V}A_{v})\Delta_{1,imp,glob}(h'_{1},h(\tilde{\zeta}_{1}))\Delta_{2,imp,glob}(h'_{2},h(\tilde{\zeta}_{2}))^{-1}$$
$$\prod_{v\not\in V}\tilde{\lambda}_{v}^M(m'_{1},m'_{2})^{-1}\Delta_{1,imp,v}(h'_{1},h(\tilde{\zeta}_{1}))^{-1}\Delta_{2,imp,v}(h'_{2},h(\tilde{\zeta}_{2})).$$
Revenons maintenant \`a la formule (18), en se rappelant que l'on a d\'ej\`a \'elimin\'e les facteurs $\Delta_{II}$ de cette formule. On obtient
$$(24) \qquad \tilde{\lambda}_{z}(r'_{1})=\lambda^M(b_{1},b_{2})^{-1}\lambda^L(a_{1},a_{2})\prod_{v\in V}A_{v},$$
o\`u on a r\'etabli l'indice $z$ pour plus de pr\'ecision.

Il faut revenir au calcul du terme $A_{v}$. Il s'av\`ere que c'est exactement le m\^eme que celui qui intervenait dans la preuve de la proposition 1.14(iii) de [II]. On peut utiliser les calculs de cette preuve. On d\'ecompose le probl\`eme en deux. Supposons d'abord 

(25) il existe une place $v\in V$ telle que $\hat{R}_{v}$ ne corresponde \`a aucun $K$-espace de Levi de $K\tilde{G}$. 

Dans ce cas, on choisit pour ${\cal Z}$ le groupe
$$(26) \qquad \left({\cal Z}(\hat{\mathfrak{R}}_{V})\cap\bigcap_{v\in V}Z(\hat{L}_{v})^{\Gamma_{F_{v}}}(1-\hat{\theta})\circ\pi(Z(\hat{R}_{sc})^{\Gamma_{F_{v}}})\right)/Z(\hat{G})^{\Gamma_{F},\hat{\theta}}$$
(on note comme toujours $\pi:\hat{G}_{SC}\to \hat{G}$ l'homomorphisme naturel). Dans les constructions pr\'ec\'edentes, on peut supposer que, pour $v\in V$, $\rho_{v}$ est l'image d'un \'el\'ement $\rho_{v,sc}\in Z(\hat{R}_{sc})^{\Gamma_{F_{v}}}$. On rappelle que l'on a choisi un rel\`evement $z_{sc}\in Z(\hat{M}_{sc})^{\Gamma_{F},\hat{\theta},0}$. On pose $\tau_{v,sc}=z_{sc}(1-\hat{\theta})(\rho_{v,sc})^{-1}$. On voit qu'il s'agit d'un \'el\'ement de $Z(\hat{L}_{v,sc})^{\Gamma_{F_{v}}}$. L'homomorphisme 
$$Z(\hat{G}_{SC})^{\Gamma_{F_{v}}}/Z(\hat{G}_{SC})^{\Gamma_{F_{v}},0}\to Z(\hat{L}_{v,sc})^{\Gamma_{F_{v}}}/Z(\hat{L}_{v,sc})^{\Gamma_{F_{v}},0}$$
 est surjectif. Fixons $x_{v}\in Z(\hat{G}_{SC})^{\Gamma_{F_{v}}}/Z(\hat{G}_{SC})^{\Gamma_{F_{v}},0}$ qui s'envoie sur l'image de $\tau_{v,sc}$ dans le groupe de droite. Rappelons qu'au groupe $G$ est associ\'e un \'el\'ement de $H^1(F;G_{AD})$: la classe du cocycle $\sigma\mapsto u_{{\cal E}}(\sigma)_{ad}$ pour une paire de Borel \'epingl\'ee ${\cal E}$ quelconque. On la note $u_{G}$.  Pour tout $v\in V$, le groupe $H^1(F;G_{AD})$ s'envoie dans  son  analogue $H^1(F_v;G_{AD})$. On note $u_{G,v}$ l'image de $u_{G}$. Pour tout $v\in V$, on a un produit naturel sur
 $$H^1(F_v;G_{AD})\times Z(\hat{G}_{SC})^{\Gamma_{F_{v}}}/Z(\hat{G}_{SC})^{\Gamma_{F_{v}},0}.$$
 A partir de l'\'egalit\'e (24), les calculs de [II] prouvent que, pour $z$ dans le groupe ${\cal Z}$ d\'efini par (26), la fonction $\tilde{\lambda}_{z}$ est constante, de valeur
 $$\prod_{v\in V}<u_{G,v},x_{v}>.$$
 La famille $(u_{G,v})_{v\in V}$ d\'efinit un caract\`ere de $\prod_{v\in V} Z(\hat{G}_{SC})^{\Gamma_{F_{v}}}/Z(\hat{G}_{SC})^{\Gamma_{F_{v}},0}$, notons-le $\chi_{G}$. La formule ci-dessus est l'\'evaluation de ce caract\`ere $\chi_{G}$ au point $(x_{v})_{v\in V}$. On a effectu\'e divers choix pour d\'efinir les $x_{v}$. La formule montre que ces choix n'affectent la famille $(x_{v})_{v\in V}$ qu'en la multipliant par un \'el\'ement du noyau de $\chi_{G}$. On a donc une application
 $${\cal Z}\to \left(\prod_{v\in V} Z(\hat{G}_{SC})^{\Gamma_{F_{v}}}/Z(\hat{G}_{SC})^{\Gamma_{F_{v}},0}\right)/Ker(\chi_{G}),$$
 qui est clairement un homomorphisme. Pour obtenir (17), il suffit de prouver que cet homomorphisme est non trivial. L'hypoth\`ese (25) et le lemme [I] 3.5 (qui reprend le lemme 2.1 de [A8]) entra\^{\i}nent que le sous-groupe 
 $$\prod_{v\in V}(Z(\hat{G}_{SC})^{\Gamma_{F_{v}}}\cap Z(\hat{R}_{v,sc})^{\Gamma_{F_{v}},0})/Z(\hat{G}_{SC})^{\Gamma_{F_{v}},0}\subset \prod_{v\in V}Z(\hat{G}_{SC})^{\Gamma_{F_{v}}}/Z(\hat{G}_{SC})^{\Gamma_{F_{v}},0}$$
 n'est pas contenu dans $Ker(\chi_{G})$.
 
   Il suffit de prouver que tout \'el\'ement de ce sous-groupe peut \^etre choisi comme famille $(x_{v})_{v\in V}$ associ\'ee \`a un \'el\'ement de ${\cal Z}$. Soit donc pour tout $v\in V$ un \'el\'ement $x_{v}\in Z(\hat{G}_{SC})^{\Gamma_{F_{v}}}\cap Z(\hat{R}_{v,sc})^{\Gamma_{F_{v}},0}$. Gr\^ace \`a l'\'egalit\'e 
 $$Z(\hat{R}_{v,sc})^{\Gamma_{F_{v}},0}=Z(\hat{R}_{v,sc})^{\Gamma_{F_{v}},\hat{\theta},0}(1-\hat{\theta})(Z(\hat{R}_{v,sc})^{\Gamma_{F_{v}},0},$$
 on peut \'ecrire $x_{v}=y_{v}(1-\hat{\theta})(\rho_{v,sc})^{-1}$, avec $y_{v}\in Z(\hat{R}_{v,sc})^{\Gamma_{F_{v}},\hat{\theta},0}$ et $\rho_{v,sc}\in Z(\hat{R}_{v,sc})^{\Gamma_{F_{v}},0}$. Notons ${\cal Z}(\hat{\mathfrak{R}}_{V})_{sc}$ l'image r\'eciproque de ${\cal Z}(\hat{\mathfrak{R}}_{V})$ dans $\hat{G}_{SC}$ et ${\cal Z}(\hat{\mathfrak{R}}_{V})_{sc}^0$ sa composante neutre. Puisque ${\cal Z}(\hat{\mathfrak{R}}_{V})\subset Z(\hat{M})^{\Gamma_{F},\hat{\theta}}$, on a ${\cal Z}(\hat{\mathfrak{R}}_{V})_{sc}^0\subset Z(\hat{M}_{sc})^{\Gamma_{F},\hat{\theta},0}$.  On d\'efinit un homomorphisme
$$(27)\qquad   {\cal Z}(\hat{\mathfrak{R}}_{V})_{sc}^0\times \prod_{v\in V}Z(\hat{L}_{v,sc})^{\Gamma_{F_{v}}, \hat{\theta},0}\to \prod_{v\in V}Z(\hat{R}_{v,sc})^{\Gamma_{F_{v}},\hat{\theta},0}.$$
Sur ${\cal Z}(\hat{\mathfrak{R}}_{V})_{sc}^0$, c'est le plongement diagonal. Sur chaque $Z(\hat{L}_{v,sc})^{\Gamma_{F_{v}}, \hat{\theta},0}$, c'est le plongement naturel $Z(\hat{L}_{v,sc})^{\Gamma_{F_{v}}, \hat{\theta},0}\to Z(\hat{R}_{v,sc})^{\Gamma_{F_{v}}, \hat{\theta},0}$. La surjectivit\'e de $\hat{D}$ entra\^{\i}ne que l'homomorphisme ainsi d\'efini est
 surjectif.  On peut donc choisir $z_{sc}\in {\cal Z}(\hat{\mathfrak{R}}_{V})_{sc}^0\subset  Z(\hat{M}_{sc})^{\Gamma_{F},\hat{\theta},0}$ et, pour tout $v\in V$, un \'el\'ement $\tau'_{v,sc}\in Z(\hat{L}_{v,sc})^{\Gamma_{F_{v}}, \hat{\theta},0}$ de sorte que $y_{v}=z_{sc}(\tau'_{v,sc})^{-1}$. Pour tout $v\in V$, on pose $\tau_{v,sc}=x_{v}\tau'_{v,sc}$. On a alors $z_{sc}=\tau_{v,sc}(1-\hat{\theta})(\rho_{v,sc})$. On voit que l'\'el\'ement $z=\pi(z_{sc})$ appartient \`a ${\cal Z}$ et que, pour cet \'el\'ement, on peut choisir la famille $(x_{v})_{v\in V}$ comme famille associ\'ee. Cela ach\`eve la preuve sous l'hypoth\`ese (25).

Supposons maintenant que pour tout $v\in V$, le Levi $\hat{R}_{v}$ est associ\'e \`a un $K$-espace de Levi de $\tilde{G}$ (sur $F_{v}$). On fixe un tel espace $K\tilde{R}_{v}$. On introduit le groupe $R_{v,0}$ quasi-d\'eploy\'e sur $F_{v}$ et dual de $\hat{R}_{v}^{\hat{\theta},0}$. Il lui est associ\'e un espace $\tilde{R}_{v,0}$. D'apr\`es [I] 1.12, on a des homomorphismes
$$\begin{array}{ccc}KR_{v,ab}(F_{v})&&\\ &N^{R_{v}}\searrow&\\ && R_{v,0,ab}(F_{v})\\&N^{R'_{v},R_{v}} \nearrow&\\ R'_{v,ab}(F_{v})&&\\ \end{array}$$
et des applications compatibles
$$\begin{array}{ccc}K\tilde{R}_{v,ab}(F_{v})&&\\ &N^{\tilde{R}_{v}}\searrow&\\ && \tilde{R}_{v,0,ab}(F_{v})\\&N^{\tilde{R}'_{v},\tilde{R}_{v}} \nearrow&\\ \tilde{R}'_{v,ab}(F_{v})&&\\ \end{array}$$
Tout \'el\'ement $\rho_{v}\in Z(\hat{R}_{v})_{*}$ d\'efinit un caract\`ere $\chi_{\rho_{v}}$ de $R_{v,0,ab}(F_{v})/N^{R_{v}}(KR_{v,ab}(F_{v}))$: c'est le caract\`ere associ\'e au cocycle $w\mapsto w_{G}(\rho_{v})\rho_{v}^{-1}$ de $W_{F_{v}}$ dans $Z(\hat{R}_{v,0})$.   L'application $\rho_{v}\mapsto \chi_{\rho_{v}}$ se quotiente en une surjection  $$Z(\hat{R}_{v})_{*}/(Z(\hat{R}_{v})\cap \hat{T}^{\hat{\theta},0})Z(\hat{R}_{v})^{\Gamma_{F_{v}}}\to (R_{v,0,ab}(F_{v})/N^{R_{v}}(KR_{v,ab}(F_{v})))^{\vee},$$
l'exposant $^{\vee}$ d\'esignant le groupe des caract\`eres. On note $\tilde{\chi}_{\rho_{v}}$ la fonction sur $\tilde{R}_{v,0,ab}(F_{v})$ qui se transforme selon le caract\`ere $\chi_{\rho_{v}}$ et qui vaut $1$ sur l'image de $N^{\tilde{R}_{v}}$.

 On choisit pour ${\cal Z}$ le groupe
$$\left({\cal Z}(\hat{\mathfrak{R}}_{V})\cap \bigcap_{v\in V}Z(\hat{L}_{v})^{\Gamma_{F_{v}}}(1-\hat{\theta})(Z(\hat{R}_{v})_{*})\right)/Z(\hat{G})^{\Gamma_{F},\hat{\theta}}$$
tout entier.  Comme pr\'ec\'edemment, pour $z\in {\cal Z}$ et $v\in V$, on \'ecrit $z=\tau_{v}(1-\hat{\theta})(\rho_{v})$, avec $\tau_{v}\in Z(\hat{L}_{v})^{\Gamma_{F_{v}}}$ et $\rho_{v}\in Z(\hat{R}_{v})_{*}$. A partir de (24), le calcul de [II] conduit cette fois au r\'esultat suivant: on a l'\'egalit\'e
$$\tilde{\lambda}_{z}(r'_{1})=\prod_{v\in V}\tilde{\chi}_{\rho_{v}}(N^{\tilde{R}'_{v},\tilde{R}_{v}}(r'_{v})),$$
o\`u on a \'ecrit $r'=(r'_{v})_{v\in V}$ l'image de $r'_{1}$ dans $\prod_{v\in V}\tilde{R}'_{v}(F_{v})$. Ici encore, le r\'esultat montre que le membre de droite ne d\'epend que de $z$ et pas des choix des $\rho_{v}$. On peut donc \'ecrire
$$(28) \qquad \sum_{z\in {\cal Z}}\tilde{\lambda}_{z}(r'_{1})=\sum_{(\rho_{v})_{v\in V}\in J}\prod_{v\in V}\tilde{\chi}_{\rho_{v}}(N^{\tilde{R}'_{v},\tilde{R}_{v}}(r'_{v})),$$
o\`u  $J$ est le sous-groupe de 
$$(29)\qquad \prod_{v\in V}Z(\hat{R}_{v})_{*}/(Z(\hat{R}_{v})\cap \hat{T}^{\hat{\theta},0})Z(\hat{R}_{v})^{\Gamma_{F_{v}}}$$
form\'e des images naturelles des familles $(\rho_{v})_{v\in V}$ que l'on peut associer comme ci-dessus \`a un \'el\'ement $z\in {\cal Z}$. On se rappelle la propri\'et\'e (2):  pour tout $v\in V$,  l'\'el\'ement  $\epsilon_{v}$ appartient \`a un sous-tore tordu elliptique de $\tilde{R}'_{v}$. Il y a donc pour tout $v\in V$ des \'el\'ements elliptiques r\'eguliers dans $\tilde{R}'_{v}(F_{v})$ qui sont aussi voisins qu'on le veut de $\epsilon_{v}$. Mais, puisque ${\bf M}'$ n'est pas relevante, il y a a fortiori au moins une place $v\in V$ o\`u $\tilde{R}'_{v}$ n'est pas relevant. La proposition [I] 1.14 implique alors que, pour cette place $v$, on a $N^{\tilde{R}'_{v},\tilde{R}_{v}}(\epsilon_{v})\not\in N^{\tilde{R}_{v}}(K\tilde{R}_{v,ab}(F_{v}))$.  Pour prouver (16), il suffit de prouver que le membre de droite de (28) est nul pour $r'$ assez voisin de $\epsilon$. D'apr\`es la propri\'et\'e que l'on vient de voir  de cet \'el\'ement, il suffit que $J$ soit \'egal au groupe (29) tout entier. Soit donc $(\rho_{v})_{v\in V}\in \prod_{v\in V}Z(\hat{R}_{v})_{*}$. Pour tout $v\in V$, on a
$$(1-\hat{\theta})(\rho_{v})\in Z(\hat{R}_{v})^{\Gamma_{F_{v}}}=Z(\hat{G})^{\Gamma_{F_{v}}}Z(\hat{R}_{v})^{\Gamma_{F_{v}},0}=Z(\hat{G})^{\Gamma_{F_{v}}}Z(\hat{R}_{v})^{\Gamma_{F_{v}},\hat{\theta},0}(1-\hat{\theta})(Z(\hat{R}_{v})^{\Gamma_{F_{v}},0}).$$
Quitte \`a multiplier $\rho_{v}$ par un \'el\'ement de $Z(\hat{R}_{v})^{\Gamma_{F_{v}},0}$, ce qui ne change pas l'image dans (29) de la famille $(\rho_{v})_{v\in V}$, on peut supposer $(1-\hat{\theta})(\rho_{v})=\xi_{v}x_{v}$, avec $\xi_{v}\in Z(\hat{G})^{\Gamma_{F_{v}}}$ et $x_{v}\in Z(\hat{R}_{v})^{\Gamma_{F_{v}},\hat{\theta},0}$. On a un homomorphisme 
$${\cal Z}(\hat{\mathfrak{R}}_{V})^0\times\prod_{v\in V}Z(\hat{L}_{v})^{\Gamma_{F_{v}},\hat{\theta},0}\to \prod_{v\in V}Z(\hat{R}_{v})^{\Gamma_{F_{v}},\hat{\theta},0}$$
similaire \`a (27), qui est surjectif car $\hat{D}$ l'est. 
On peut donc  trouver $z\in  {\cal Z}(\hat{\mathfrak{R}}_{V})^0$ et, pour tout $v\in V$, un \'el\'ement $\tau'_{v}\in Z(\hat{L}_{v})^{\Gamma_{F_{v}},\hat{\theta},0}$, de sorte que $z\tau'_{v}=x_{v}$. 
Posons $\tau_{v}=(\xi_{v}\tau'_{v})^{-1}$. C'est un \'el\'ement de $Z(\hat{L}_{v})^{\Gamma_{F_{v}}}$ et on a $z=\tau_{v}(1-\hat{\theta})(\rho_{v})$. On voit alors que $z$ appartient \`a ${\cal Z}$ et que $(\rho_{v})_{v\in V}$ est une famille associ\'ee \`a $z$ comme plus haut. Cela ach\`eve cette br\`eve d\'emonstration. $\square$

{\bf Remarques.} (30) Les complications du d\'ebut de la preuve, \`a savoir le passage par les propositions auxiliaires 6.7 et 6.8, sont dues aux places archim\'ediennes. Cela parce que, sur un corps de base $F_{v}$ archim\'edien, on n'a pas d\'efini les int\'egrales orbitales pond\'er\'ees stables $S_{\tilde{M}}^{\tilde{G}}(\boldsymbol{\delta},{\bf f})$ pour tout $\boldsymbol{\delta}\in D_{g\acute{e}om}^{st}(\tilde{M}( F_{v}))\otimes Mes(M(F_{v}))^*$. Cela nous oblige \`a des contorsions pour assurer que les distributions ${\bf d}_{v}$ de la d\'emonstration ci-dessus sont dans le domaine de d\'efinition de ces int\'egrales orbitales pond\'er\'ees stables. Si on remplace ces int\'egrales  par leurs avatars spectraux, cette difficult\'e dispara\^{\i}t et une m\^eme d\'emonstration s'applique tout en se simplifiant. 

(31) Si on oublie cette difficult\'e aux places archim\'ediennes, on voit que l'on d\'emontre en fait un r\'esultat plus fin que la proposition 6.6, qui est le suivant. Ecrivons $\boldsymbol{\delta}=(\boldsymbol{\delta}_{v})_{v\in V}$ et supposons que, pour tout $v\in V$, $\boldsymbol{\delta}_{v}$ soit induite \`a partir d'un espace de Levi $\tilde{R}'_{v}$ de $\tilde{M}'_{v}$. Supposons qu'il existe $v$ tel que $\tilde{R}'_{v}$ ne corresponde \`a aucun $K$-espace de Levi de $K\tilde{G}_{v}$ ou qu'il existe un tel $K$-espace de Levi mais que $\tilde{R}'_{v}$ ne soit pas relevant. Alors $I_{*}^{K\tilde{G},{\cal E}}({\bf M}',\boldsymbol{\delta},{\bf f})=0$.

 \bigskip
 
 \subsection{Le th\'eor\`eme 5.10}
 Dans ce paragraphe, on suppose d\'emontr\'es les th\'eor\`emes   [II] 1.16, [V] 1.10 et 5.2 et 5.4 ci-dessus.  On va alors prouver le th\'eor\`eme 5.10. Notons
 $$X=\sum_{{\bf G}'\in {\cal E}( \tilde{G},{\bf a},V)}i(\tilde{G},\tilde{G}')S^{{\bf G}'}_{g\acute{e}om}({\bf f}^{{\bf G}'})$$
  le membre de droite de l'\'egalit\'e de ce th\'eor\`eme. D\'eveloppons $X$ en appliquant les d\'efinitions. On obtient
 $$ X=\sum_{{\bf G}'\in {\cal E}( \tilde{G},{\bf a},V)}i(\tilde{G},\tilde{G}')\sum_{\tilde{M}'\in {\cal L}(\tilde{M}'_{0})}\vert W^{M'}\vert \vert W^{G'}\vert ^{-1}\sum_{{\cal O}'\in \tilde{M}'_{ss}(F_{V})/st-conj}S_{{\bf M}'}^{{\bf G}'}(SA^{{\bf M}'}({\cal O}',V),{\bf f}^{{\bf G}'}).$$
 Le terme ${\bf M}'$ est le triplet $(M',{\cal M}',\tilde{\zeta})$ associ\'e naturellement aux donn\'ees ${\bf G}'$ et \`a l'espace de Levi $\tilde{M}'\in {\cal L}(\tilde{M}'_{0})$.  Notons que $\tilde{M}'(F)\not=\emptyset$ pour tous les $\tilde{M}'$ intervenant.
 Le lemme 5.3 nous autorise \`a remplacer les termes  $S_{{\bf M}'}^{{\bf G}'}(SA^{{\bf M}'}({\cal O}',V),{\bf f}^{{\bf G}'})$ par leurs variantes $S_{{\bf M}'}^{{\bf G}'}(SA^{{\bf M}'}({\cal O}',V),B^{\tilde{G}},{\bf f}^{{\bf G}'})$.
 Il n'est pas difficile d'adapter la proposition 6.5 \`a nos pr\'esentes notations. Elle entra\^{\i}ne que l'on peut r\'ecrire cette somme sous la forme
 $$X=\sum_{\hat{M}}\vert W^{\tilde{M}}\vert \vert W^{\tilde{G}}\vert ^{-1}\sum_{{\bf M}'\in {\cal E}_{*}(\tilde{M},{\bf a}_{M},V)}i(\tilde{M},\tilde{M}')\sum_{\tilde{s}\in \tilde{\zeta}Z(\hat{M})^{\Gamma_{F},\hat{\theta}}/Z(\hat{G})^{\Gamma_{F},\hat{\theta}}}i_{\tilde{M}'}(\tilde{G},\tilde{G}'(\tilde{s}))$$
 $$\sum_{{\cal O}'\in \tilde{M}'_{ss}(F_{V})/st-conj}S_{{\bf M}'}^{{\bf G}'(\tilde{s})}(SA^{{\bf M}'}({\cal O}',B,V),B^{\tilde{G}},{\bf f}^{{\bf G}'(\tilde{s})}).$$
 Ici, $\hat{M}$ parcourt les Levi de $\hat{G}$ contenant $\hat{M}_{0}$ qui sont des composantes de Levi de sous-groupes paraboliques invariants par $\hat{\theta}$ et par $\Gamma_{F}$. Un tel Levi ne correspond pas toujours \`a un $K$-espace de Levi $K\tilde{M}$ de $K\tilde{G}$ mais on peut lui associer divers objets que,  par anticipation, nous notons comme si un tel $K$-espace existait. Par exemple $W^{\tilde{M}}$ est le sous-groupe des \'el\'ements invariants par $\hat{\theta}$ et $\Gamma_{F}$ dans $Norm_{\hat{G}}(\hat{M})/\hat{M}$. Si $K\tilde{M}$ existe, l'ensemble ${\cal E}_{*}(\tilde{M},{\bf a}_{M},V)$ est celui des classes d'\'equivalence de donn\'ees endoscopiques elliptiques $(M',{\cal M}',\tilde{\zeta})$ de $(KM,K\tilde{M},{\bf a}_{M})$ qui sont non ramifi\'ees hors de $V$ et telles que $\tilde{M}'(F)\not=\emptyset$. La d\'efinition de ces notions ne faisant intervenir que le groupe $\hat{M}$, la d\'efinition s'\'etend au cas o\`u $K\tilde{M}$ n'existe pas. L'indice $*$ signifie que l'on n'impose pas que la donn\'ee soit relevante, contrairement \`a nos ensembles habituels   ${\cal E}(\tilde{M},{\bf a}_{M},V)$.    Avec la d\'efinition de 6.6, on obtient
 $$X=\sum_{\hat{M}}\vert W^{\tilde{M}}\vert \vert W^{\tilde{G}}\vert ^{-1}\sum_{{\bf M}'\in {\cal E}_{*}(\tilde{M},{\bf a}_{M},V)}i(\tilde{M},\tilde{M}')\sum_{{\cal O}'\in \tilde{M}'_{ss}(F_{V})/st-conj}I_{*}^{K\tilde{G},{\cal E}}({\bf M}',SA^{{\bf M}'}({\cal O}',V),{\bf f}).$$
 La proposition 6.6 nous dit que, pour que le terme que l'on somme soit non nul, il faut que $\hat{M}$  corresponde \`a un $K$-espace de Levi $K\tilde{M}$ de $\tilde{G}$ et que ${\bf M}'$ soit relevant pour $(KM,K\tilde{M},{\bf a}_{M})$. Cela permet de r\'ecrire
 $$X=\sum_{K\tilde{M}\in {\cal L}(\tilde{M}_{0})}\vert W^{\tilde{M}}\vert \vert W^{\tilde{G}}\vert ^{-1}\sum_{{\bf M}'\in {\cal E}(\tilde{M},{\bf a}_{M},V)}i(\tilde{M},\tilde{M}')$$
 $$\sum_{{\cal O}'\in \tilde{M}'_{ss}(F_{V})/st-conj}I_{K\tilde{M}}^{K\tilde{G},{\cal E}}({\bf M}',SA^{{\bf M}'}({\cal O}',V),{\bf f}).$$
 Les hypoth\`eses de la proposition 4.6 sont v\'erifi\'ees. Cela nous permet de remplacer les termes $I_{K\tilde{M}}^{K\tilde{G},{\cal E}}({\bf M}',SA^{{\bf M}'}({\cal O}',V),{\bf f})$ par $I_{K\tilde{M}}^{K\tilde{G}}(transfert(SA^{{\bf M}'}({\cal O}',V)),{\bf f})$. 
  On peut maintenant se limiter aux classes de conjugaison stable ${\cal O}'$ qui correspondent \`a une telle classe dans le $K$-espace de Levi correspondant $K\tilde{M}(F_{V})$: pour les autres, $transfert(SA^{{\bf M}'}({\cal O}',V))=0$.
 On peut regrouper ces classes selon la classe qui leur correspond dans $K\tilde{M}(F_{V})$. On obtient
  $$X=\sum_{K\tilde{M}\in {\cal L}(K\tilde{M}_{0})}\vert W^{\tilde{M}}\vert \vert W^{\tilde{G}}\vert ^{-1}\sum_{{\cal O}\in K\tilde{M}_{ss}(F_{V})/st-conj}\sum_{{\bf M}'\in {\cal E}(\tilde{M},{\bf a}_{M},V)}i(\tilde{M},\tilde{M}')$$
 $$\sum_{{\cal O}'\in \tilde{M}'_{ss}(F_{V})/st-conj; {\cal O}'\mapsto {\cal O}}I_{K\tilde{M}}^{K\tilde{G}}(transfert(SA^{{\bf M}'}({\cal O}',V)),{\bf f}),$$
 o\`u ${\cal O}'\mapsto {\cal O}$ d\'esigne la correspondance entre classes de conjugaison stable. Pour tous $K\tilde{M}$ et ${\cal O}$ intervenant ci-dessus, on a l'\'egalit\'e 
$$ \sum_{{\bf M}'\in {\cal E}(\tilde{M},{\bf a}_{M},V)}i(\tilde{M},\tilde{M}')\sum_{{\cal O}'\in \tilde{M}'_{ss}(F_{V})/st-conj; {\cal O}'\mapsto {\cal O}}transfert(SA^{{\bf M}'}({\cal O}',V))$$
$$=A^{K\tilde{M},{\cal E}}({\cal O},V,{\bf a})$$
d'apr\`es la d\'efinition de 5.4. C'est encore \'egal \`a $A^{K\tilde{M}}({\cal O},V,\omega)$ d'apr\`es le th\'eor\`eme 5.4. Les hypoth\`eses de la proposition 4.6 sont v\'erifi\'ees et on obtient
$$X=\sum_{K\tilde{M}\in {\cal L}(K\tilde{M}_{0})}\vert W^{\tilde{M}}\vert \vert W^{\tilde{G}}\vert ^{-1}\sum_{{\cal O}\in \tilde{M}_{ss}(F_{V})/st-conj}I_{K\tilde{M}}^{K\tilde{G}}(A^{K\tilde{M}}({\cal O},V,\omega),{\bf f}).$$
Ceci n'est autre que $I^{K\tilde{G}}_{g\acute{e}om}({\bf f},\omega)$, cf. 2.9. Cela prouve le th\'eor\`eme 5.10.

\bigskip

 {\bf Bibliographie}
 
  [A1] J. Arthur: {\it A stable trace formula I. General expansions}, JIMJ 1 (2002), p. 175-277
 
 [A2] -----------: {\it The invariant trace formula I. Local theory}, J. AMS 1 (1988), p. 323-383
 
 [A3]-----------: {\it The invariant trace formula II. Global theory}, J. AMS 1 (1988), p. 501- 554

 [A4] -----------: {\it Canonical normalization of weighted characters and a transfer conjecture}, C. R. Math. Acad. Sci. Canada 20 (1998), p. 35-52
 
 [A5] -----------: {\it The local behaviour of weighted orbital integrals}, Duke Math. J. 56 (1988), p. 223-293
 
 [A6] -----------: {\it On a family of distributions obtained from Eisenstein series II: explicit formulas}, Amer. J. of Math. 104 (1982), p. 1289-1336
 
  [A7]-----------: {\it A measure on the unipotent variety}, Can. J. Math. 37 (1985), p. 1237-1274
  
   [A8] -----------: {\it On the transfer of distributions: weighted orbital integrals}, Duke Math. J. 99 (1999), p.209-283
 
  [A9]------------: {\it The trace formula in invariant form}, Annals of Math. 114 (1981), p. 1-74

 [A10] -----------: {\it On a family of distributions obtained from orbits}, Can. J. Math. 38 (1986), p. 179-214
 
 [A11] -----------: {\it Endoscopic $L$-functions and a combinatorial identity}, Can. J. of Math. 51 (1999), p. 1135-1148

 [F] A. Ferrari: {\it Th\'eor\`eme de l'indice et formule des traces}, manuscripta math. 124 (2007), p. 363-390
 
 [K1] R. Kottwitz: {\it Rational conjugacy classes in reductive groups}, Duke Math. J. 49 (1982), p. 785-806
 
 [K2]--------------: {\it Stable trace formula: cuspidal tempered terms}, Duke Math. J. 51 (1984), p. 611-650
 
  [K3]--------------: {\it Stable trace formula: elliptic singular terms}, Math. Annalen 275 (1986), p. 365-399
 
 [KS] ------------, D. Shelstad: {\it Foundations of twisted endoscopy}, Ast\'erisque 255 (1999)
 
 [KS2] ----------------------------: {\it On splitting invariants and sign conventions in endoscopic transfer}, pr\'epublication 2012
 
 [Lab] J.-P. Labesse: { \it Cohomologie, stabilisation et changement de base}, Ast\'erisque 257 (1999)
 
 [LabW] ---------------, J.-L. Waldspurger: {\it La formule des traces tordue d'apr\`es le Friday Morning Seminar}, pr\'epublication 2012
 
 [Lan] R. P. Langlands: {\it Stable conjugacy: definitions and lemmas}, Can. J. of Math 31 (1979), p. 700-725
 
 [S] J.-L. Sansuc: {\it Groupe de Brauer et arithm\'etique des groupes alg\'ebriques lin\'eaires sur un corps de nombres}, Journal f\"{u}r die reine u. ang. Math. 327 (1981), p.12-80
 
 [T] J. Tate: {\it Number theoretic background}, in {\it Automorphic forms, representations and $L$-functions}, A. Borel et W. Casselman ed., Proc. of Symp. in Pure Math. XXXIII, part 2, AMS 1979
 
  [W1] J.-L. Waldspurger: {\it La formule des traces locale tordue}, pr\'epublication 2012
 
 [W2] ---------------------: {\it L'endoscopie tordue n'est pas si tordue}, Memoirs AMS 908 (2008)
 
 [I] ---------------------: {\it Stabilisation de la formule des traces tordue I: endoscopie tordue sur un corps local}, pr\'epublication 2014
 
 [II] ----------------------: {\it Stabilisation de la formule des traces tordue II: int\'egrales orbitales et endoscopie sur un corps local non-archim\'edien; d\'efinitions et \'enonc\'es des r\'esultats}, pr\'epublication 2012
 
 [III] ---------------------:  {\it Stabilisation de la formule des traces tordue III: int\'egrales orbitales et endoscopie sur un corps local non-archim\'edien; r\'eductions et preuves}, pr\'epublication 2014
 
 [IV]---------------------:  {\it Stabilisation de la formule des traces tordue IV: transfert spectral archim\'edien}, pr\'epublication 2014
 
 [V]---------------------: {\it Stabilisation de la formule des traces tordue V: int\'egrales orbitales et endoscopie sur le corps r\'eel}, pr\'epublication 2014
 
\bigskip
 
 Institut de Math\'ematiques de Jussieu -CNRS
 
 2 place Jussieu, 75005 Paris 
 
 e-mail: colette.moeglin@imj-prg.fr, jean-loup.waldspurger@imj-prg.fr

\end{document}